\documentclass[11pt, reqno]{amsart}
\allowdisplaybreaks

\usepackage{comment}
\usepackage[rgb,dvipsnames]{xcolor}
\usepackage{tikz-cd}
\usepackage{faktor}
\usepackage{stmaryrd}
\usepackage{amsfonts}
\usepackage{amssymb}
\usepackage{amsmath,amsthm}
\usepackage{mathtools}
\usepackage{latexsym}
\usepackage{amscd}
\usepackage{esint}
\usepackage{faktor}
\usepackage{mathrsfs}
\usepackage{dsfont}
\usepackage{setspace}
\usepackage{enumitem}
\usepackage[margin=1.2in]{geometry}
\usepackage{bm}
\usepackage{hyperref}

\newtheorem{prop}{Proposition}[section]
\newtheorem{thm}[prop]{Theorem}
\newtheorem{lemm}[prop]{Lemma}
\newtheorem{coro}[prop]{Corollary}

\newtheorem{claim}{Claim}[section]
\newtheorem*{claim*}{Claim}
\newtheorem*{lemm*}{Lemma}
\newtheorem*{thm*}{Theorem}

\theoremstyle{definition}
\newtheorem{defi}[prop]{Definition}
\newtheorem{rmk}[prop]{Remark}

\newcommand{\CC}{\mathbb{C}}

\newcommand{\HH}{\mathbb{H}}

\newcommand{\NN}{\mathbb{N}}

\newcommand{\RR}{\mathbb{R}}

\newcommand{\ZZ}{\mathbb{Z}}

\newcommand{\cA}{\mathcal A}
\newcommand{\cB}{\mathcal B}
\newcommand{\cC}{\mathcal C}
\newcommand{\cD}{\mathcal D}
\newcommand{\cE}{\mathcal E}
\newcommand{\cF}{\mathcal F}

\newcommand{\cH}{\mathcal H}

\newcommand{\cK}{\mathcal K}

\newcommand{\cM}{\mathcal M}
\newcommand{\cN}{\mathcal N}
\newcommand{\cO}{\mathcal O}
\newcommand{\cP}{\mathcal P}

\newcommand{\cR}{\mathcal R}
\newcommand{\cS}{\mathcal S}

\newcommand{\cU}{\mathcal U}
\newcommand{\cV}{\mathcal V}
\newcommand{\cW}{\mathcal W}
\newcommand{\cX}{\mathcal X}
\newcommand{\cY}{\mathcal Y}
\newcommand{\cZ}{\mathcal Z}

\newcommand{\sD}{\mathscr{D}}
\newcommand{\sE}{\mathscr{E}}

\def\fB{\mathfrak{B}}

\def\fc{\mathfrak{c}}

\def\bB{\mathbf{B}}
\def\bC{\mathbf{C}}

\def\bP{\mathbf{P}}

\def\bt{\mathbf{t}}

\def\bh{\mathbf{h}}

\def\ba{\mathbf{a}}
\def\bb{\mathbf{b}}
\def\bc{\mathbf{c}}
\def\be{\mathbf{e}}
\def\bg{\mathbf{g}}

\def\bu{\mathbf{u}}
\def\bv{\mathbf{v}}
\def\bw{\mathbf{w}}
\def\bx{\mathbf{x}}
\def\bsig{\boldsymbol{\sigma}}

\def\bphi{\boldsymbol{\varphi}}

\def\0{{\bf 0}}

\def\bF{\mathbf{F}}
\def\bM{\mathbf{M}}
\def\bI{\mathbf{I}}

\DeclareMathOperator{\Aut}{Aut}
\DeclareMathOperator{\Int}{Int}
\DeclareMathOperator{\tr}{tr}

\DeclareMathOperator{\im}{Im}
\DeclareMathOperator{\Span}{span}

\DeclareMathOperator{\supp}{supp}

\DeclareMathOperator{\loc}{loc}
\DeclareMathOperator{\osc}{osc}
\DeclareMathOperator{\re}{Re}

\DeclareMathOperator{\diam}{diam}
\DeclareMathOperator{\dist}{dist}

\DeclareMathOperator{\Inte}{Int}

\DeclareMathOperator{\vol}{vol}
\DeclareMathOperator{\Vol}{Vol}
\DeclareMathOperator{\Area}{Area}

\DeclareMathOperator{\inj}{inj}
\DeclareMathOperator{\Met}{Met}
\DeclareMathOperator{\hyp}{hyp}
\DeclareMathOperator{\euc}{euc}
\DeclareMathOperator{\pro}{prod}
\DeclareMathOperator{\Diff}{Diff}

\newcommand{\ep}{\varepsilon}

\newcommand{\reg}{\text{reg}}

\newcommand{\id}{\text{id}}

\newcommand{\bangle}[1]{\big\langle #1 \big\rangle}
\newcommand{\pa}[2]{\frac{\partial #1}{\partial #2}}

\newcommand{\paop}[1]{\pa{}{#1}}

\newcommand{\rom}[1]{\expandafter\romannumeral #1}
\newcommand{\Rom}[1]{\uppercase\expandafter{\romannumeral #1}}

\newcommand{\current}[1]{\llbracket #1 \rrbracket}
\setlist[enumerate]{leftmargin = 2em}
\setlist[itemize]{leftmargin = 2em}
\allowdisplaybreaks

\numberwithin{equation}{subsection}

\makeatletter
\def\@tocline#1#2#3#4#5#6#7{\relax
  \ifnum #1>\c@tocdepth 
  \else
    \par \addpenalty\@secpenalty\addvspace{#2}%
    \begingroup \hyphenpenalty\@M
    \@ifempty{#4}{%
      \@tempdima\csname r@tocindent\number#1\endcsname\relax
    }{%
      \@tempdima#4\relax
    }%
    \parindent\z@ \leftskip#3\relax \advance\leftskip\@tempdima\relax
    \rightskip\@pnumwidth plus4em \parfillskip-\@pnumwidth
    #5\leavevmode\hskip-\@tempdima
      \ifcase #1
       \or\or \hskip 1.5em \or \hskip 3em \else \hskip 4.5em \fi%
      #6\nobreak\relax
    \hfill\hbox to\@pnumwidth{\@tocpagenum{#7}}\par
    \nobreak
    \endgroup
  \fi}
\makeatother

\title{Existence of classical minimal surfaces in $4$ and $5$-manifolds}
\author{Da Rong Cheng}
\address{Department of Mathematics, University of Miami, Coral Gables, FL 33146}
\email{darong.cheng@miami.edu}

\author{Xin Zhou}
\address{Department of Mathematics, Cornell University, Ithaca, NY 14853}
\email{xinzhou@cornell.edu}

\begin{document}

\begin{abstract} 
We prove that every closed Riemannian $4$ or $5$-manifold $M$ contains a branched immersed closed minimal surface. That is, there exists a non-constant weakly conformal harmonic map from some closed Riemann surface into $M$. We rely on the existence of multisections in dimensions $4$ and $5$ to generate a non-trivial class of sweepouts of $M$ by mappings from a closed surface $S$ of genus at least $2$. To each sweepout in a minimizing sequence within the class, through the intermediary of quasiconformal maps of the upper half-plane, we associate a family of hyperbolic metrics on $S$ with respect to which the mappings in the sweepout have nearly equal energy and area. The harmonic replacement method of Colding and Minicozzi is then applied to obtain a min-max sequence that converges to a bubble tree of branched minimal immersions. 
\end{abstract}

\maketitle 

\tableofcontents

\section{Introduction}\label{sec:intro}

\subsection{Statements}\label{subsec:statements}
Suppose $S$ is a closed surface and $M$ is a closed Riemannian $n$-manifold, both assumed to be oriented and connected, with the latter isometrically embedded into some Euclidean space $\RR^q$. For $v \in W^{1, 2}(S; M)$, we denote by $A(v)$ its area, and by $E(\gamma, v)$ its Dirichlet energy with respect to a given Riemannian metric $\gamma$ on $S$. (We recall their definitions in Section~\ref{subsec:notation}.) For convenience, we sometimes refer to $E(\gamma, v)$ as the $\gamma$-energy of $v$.  Next, introducing $I := [0, 1]$ and $I_{\delta} := [\delta, 1 - \delta]$, the latter for $\delta \in (0, \frac{1}{4})$, and fixing some $m \in \NN$, we define the relevant $m$-parameter families of maps in $(C^0 \cap W^{1, 2})(S; M)$ to be used in our min-max construction.
\begin{defi}\label{defi:class-cP}
We say that a continuous map $\bv$ from $I^m$ into the space $C^0(S; M)$ belongs to the collection $\cP = \cP(m, S)$ if
\begin{itemize}
\item $\bv(t) \in (C^0 \cap W^{1, 2})(S; M)$ for all $t \in \Int(I^m)$.
\vskip 1mm
\item $\bv|_{\Int(I^m)}: \Int(I^m) \to (C^0 \cap W^{1, 2})(S; M)$ is continuous with respect to both the $C^0$ and $W^{1, 2}$-norms. That is, as $t \to t_0  \in \Int(I^m)$, we have that $\bv(t) \to \bv(t_0)$ uniformly on $S$, and that $d[\bv(t)] \to d[\bv(t_0)]$ in $L^2(S)$.
\vskip 1mm
\item The composition $A \circ \bv$, defined on $\Int(I^m)$, satisfies the boundary condition
\begin{equation}\label{eq:no-area-at-boundary}
A(\bv(t)) \to 0,\quad \text{as } \dist(t, \partial I^m) \to 0.
\end{equation}
\end{itemize}
Depending on the context, we sometimes denote a given $\bv \in \cP$ by $\{\bv(t)\}_{t \in I^m}$.
\end{defi}
Next, we partition $\cP$ into suitable ``homotopy classes''.
\begin{defi}\label{defi:relation-on-P}
For $\bv, \bw \in \cP$, we write ``$\bv \sim \bw$'' if there exist $\delta \in (0, \frac{1}{4})$ and a continuous map 
\[
\bh: [0, 1] \times \Int(I^m) \to (C^0 \cap W^{1, 2})(S; M) 
\]
such that $\bv(t) = \bw(t)$ for all $t \in I^m \setminus I_{\delta}^m$, that
\begin{equation}\label{eq:cP-homotopy-1}
\bh(0, t) = \bv(t) \quad \text{and} \quad \bh(1, t) = \bw(t),\ \text{ for all } t \in \Int(I^m),
\end{equation}
and that
\begin{equation}\label{eq:cP-homotopy-2}
\bh(s, t) = \bv(t),\quad \text{for all }(s, t) \in [0, 1] \times (\Int(I^m)\setminus I_{\delta}^m).
\end{equation}
\end{defi}
It is straightforward to verify that ``$\sim$'' is an equivalence relation. The equivalence class of a given $\bv_{0} \in \cP$ is denoted $[\bv_{0}]$, to which we associate the \emph{width}, defined by
\begin{equation}\label{eq:width-defi}
\cW([\bv_{0}]) : = \inf_{\bv \in [\bv_{0}]}\sup_{t \in \Int(I^m)}A(\bv(t)).
\end{equation}
Note that, by the continuity of $A$ with respect to strong $W^{1,2}$-convergence, and the boundary condition~\eqref{eq:no-area-at-boundary}, we have $\sup_{t \in \Int(I^m)}A(\bv(t)) < \infty$ for all $\bv \in \cP$. Following standard terminology, a sequence $(\bv_{j}) \subset [\bv_{0}]$ satisfying $\lim_{j \to \infty}\sup_{t \in \Int(I^m)}A(\bv_{j}(t))=\cW([\bv_{0}])$ is called a \emph{minimizing sequence} in $[\bv_{0}]$. Given a minimizing sequence, if $(t_j) \subset \Int(I^m)$ is such that $A(\bv_{j}(t_j)) \to \cW([\bv_{0}])$ as $j \to \infty$, then we say that $(\bv_{j}(t_j))$ is a \emph{min-max sequence}.

The two main results of this paper are the following.
\begin{thm}\label{thm:min-max-existence}
Let the surface $S$ and Riemannian manifold $M$ be as above. Suppose that $S$ has genus $g > 1$, and that there exists some $m \in \NN$ and $\bv_0 \in \cP(m, S)$ such that $\cW([\bv_0]) > 0$. Then $M$ contains a closed, branched immersed minimal surface of genus at most $g$. In other words, there exists a closed Riemann surface $\Sigma$ of genus at most $g$ which admits a non-constant, weakly conformal, harmonic mapping into $M$.
\end{thm}
\begin{thm}\label{thm:existence-of-sweepout}
Let $M$ be a closed, connected, oriented Riemannian manifold. Suppose for some $m \in \{2, 3\}$ that the homology group $H_{m+2}(M, \ZZ)$ is nontrivial. Then there exists a closed, oriented surface $S$ of genus $g > 1$ and an element $\bv_{0} \in \cP(m, S)$ such that $\cW([\bv_0]) > 0$.
\end{thm}

\begin{rmk}
In the case $m + 2 = \dim M$, that is, when $M$ has dimension $4$ or $5$, the assumption on homology in Theorem~\ref{thm:existence-of-sweepout} is fulfilled because of the fundamental class. Applying Theorem~\ref{thm:min-max-existence} to the resulting element $\bv_0 \in \cP(m, S)$, we conclude that $M$ contains a branched immersed closed minimal surface.
\end{rmk}

We next outline the proof of the two theorems. The proof of Theorem~\ref{thm:existence-of-sweepout} consists of three parts, corresponding respectively to Sections~\ref{subsec:mappings-and-currents},~\ref{subsec:Almgren isomorphism}, and~\ref{subsec:multisections}. For simplicity, here we only describe the argument in the case $m + 2 = \dim M$.
\vskip 1mm
\begin{itemize}
\item First, denoting by $\cZ_k(M)$ the space of integral $k$-cycles in $M$, with the help of the compactness theorem of integral currents and some standard constructions, we associate to each $\bv \in \cP(m, S)$ an element $[f_{\bv}] \in \pi_m(\cZ_2(M), \{0\})$ which depends only on the equivalence class of $\bv$, and deduce from a result of Almgren~\cite[Theorem 8.2]{Almgren62} that $\cW([\bv])$ would be positive provided $[f_{\bv}] \neq 0$ (Proposition \ref{prop:classes-of-maps-to-currents}). 
\vskip 1mm
\item We then recall Almgren's isomorphism 
\[
F_A:\pi_m(\cZ_k(M), \{0\}) \to H_{m + k}(M),
\]
which reduces the requirement on $\bv$ further to $F_{A}([f_{\bv}]) \neq 0$. Our strategy is to obtain such a $\bv$ by using suitable Lipschitz maps as a bridge. Specifically, we prove that to each Lipschitz map $\Phi:I^{m} \times S \to M$ satisfying certain conditions on its regularity, injectivity, and boundary behavior, we can assign a continuous map $\widetilde{f}_{\Phi}:(I^{m}, \partial I^{m}) \to (\cZ_{2}(M), \{0\})$ on the one hand, and an element $\bv \in \cP(m, S)$ on the other, in such a way that $[f_{\bv}] = [\widetilde{f}_{\Phi}]$, that $F_{A}([\widetilde{f}_{\Phi}])$ is represented by $\Phi_{\#}(\current{I^{m}} \times \current{S})$, and that the latter is homologically nontrivial. (See Corollary~\ref{coro:non-trivial-sweepout-in-P} and Proposition~\ref{prop:induced-sweepout}.)
\vskip 1mm
\item To find a surface $S$ and a Lipschitz map $\Phi: I^m \times S \to M$ which meets the conditions mentioned in the previous step, we invoke the work of Gay and Kirby~\cite{Gay-Kirby2016}, and subsequent generalizations by Ben Aribi, Courte, Golla, and Moussard~\cite{BACGM2023}, on the existence of \emph{multisections} of closed manifolds in dimensions $4$ and $5$. Taking the case $\dim M = 4$ for example, it was proved in~\cite{Gay-Kirby2016} that every closed, connected, oriented (smooth) $4$-manifold can be decomposed into three $4$-dimensional handlebodies in such a way that each pairwise intersection is a $3$-dimensional handlebody, and the triple intersection is a closed surface $S$, which can be assumed to have genus larger than $1$ without loss of generality. In analogy with how Heegaard splittings of closed $3$-manifolds give rise to $1$-parameter sweepouts, from the above decomposition, known as a \textit{trisection}, we obtain a Lipschitz map $I^{2} \times S \to M$ which feeds into the results from the previous step (Proposition~\ref{prop:multisection-to-sweepout}) to yield $\bv \in \cP(2, S)$ such that $[f_{\bv}] \neq 0$, implying the positivity of width we want.
\end{itemize}
\vskip 2mm

Turning to Theorem~\ref{thm:min-max-existence}, the strategy is to follow the earlier work of the second named author~\cite{Zhou10,Zhou17b} and apply min-max methods to the Dirichlet energy, thought of as depending on both the map $S \to M$ and the conformal class of the domain metric. Some further notation is needed to be more precise about this last point, so let $\Met_{-1}$ be the set of all metrics on $S$ with constant curvature $-1$, and fix a reference element $\gamma_0 \in \Met_{-1}$. The identity component $\Diff_0$ of the group of orientation-preserving diffeomorphisms of $S$ acts by pullback on $\Met_{-1}$, and it is known classically that $\Met_{-1}/ \Diff_0$ is one of several equivalent models for the Teichm\"uller space of closed genus-$g$ surfaces. Moreover, by standard results about harmonic maps (see Section~\ref{subsec:hyperbolic}), each $\Diff_0$ orbit contains exactly one representative from the set below:
\begin{equation}\label{eq:hyperbolic-harmonic-id}
\Met^*_{-1} := \{\gamma \in \Met_{-1}\ |\  \id: (S, \gamma_0) \to (S, \gamma) \text{ is a harmonic map}\}.
\end{equation}
In view of these facts, we consider the Dirichlet energy as a functional defined on the space
\begin{equation}\label{eq:cM'-definition}
\cM' := \Met^*_{-1} \times (C^0 \cap W^{1, 2})(S; M),
\end{equation}
equipped with the product metric coming from smooth convergence of metrics and $C^0 \cap W^{1,2}$-convergence of mappings.

Having fixed the domain of $E$, we proceed to discuss the two most important analytical ingredients in the proof of Theorem \ref{thm:min-max-existence}, namely \emph{conformal reparametrization} and \emph{harmonic replacement}. Given $\ep, \Lambda > 0$, what we call conformal reparametrization, discussed in Section \ref{sec:conformal-reparametrization}, assigns a pair $(\sigma, \widetilde{v}) \in \cM'$ to each $C^2$-map $v:S \to M$ satisfying $E(\gamma_0, v) \leq \Lambda$, such that 
\begin{equation}\label{eq:almost-conformal-intro}
A(\widetilde{v}) \leq E(\sigma, \widetilde{v}) < A(\widetilde{v}) + \ep.
\end{equation}
Loosely speaking, the construction of the map $v \mapsto (\sigma, \widetilde{v})$, carried out in Section \ref{subsec:area-energy}, proceeds as follows. Fix a covering map $p: \HH \to S$ and let $\Gamma$ be the group of deck transformations. Let $g_M$ denote the metric on $M$, and take some $\eta> 0$ to be determined later. Given $v \in C^2(S; M)$, we can always find a complex-valued $C^1$-function $\mu$ on $\HH$ such that $\|\mu\|_{\infty;\HH} < 1$, and that
\[
p^*(v^*g_M + \eta \gamma_0)\quad\text{is conformal to}\quad |dz + \mu d\overline{z}|^2.
\]
By the work of Ahlfors--Bers~\cite{Ahlfors-Bers} on quasiconformal maps, part of which we recall in Section~\ref{subsec:quasi-conformal}, we get a $C^1$-diffeomorphism $w: \HH \to \HH$ which satisfies 
\begin{equation}\label{eq:Beltrami-intro}
w_{\overline{z}} = \mu w_{z},
\end{equation}
and is thus a conformal map from $(\HH, |dz + \mu d\overline{z}|^2)$ to $(\HH, |dz|^2)$. Further, the invariance of $p^*(v^*g_M + \eta\gamma_0)$ with respect to the action of $\Gamma$ gives $\mu$ a certain symmetry that forces the conjugated group $w \Gamma w^{-1}$ to remain in $\Aut(\HH)$. In particular, the standard hyperbolic metric on $\HH$, denoted $g_{\hyp}$, descends to a metric $\gamma$ with constant curvature $-1$ on the quotient $S' : = \HH /w \Gamma w^{-1}$, and $w$ descends to a conformal map 
\[
G:(S, v^*g_M + \eta\gamma_0) \to (S', \gamma).
\]
By standard results about harmonic mappings into negatively curved targets, together with the work of Schoen--Yau~\cite{Schoen-Yau1978} on the case of mappings between surfaces, within the homotopy class of $G$, there is a unique harmonic map 
\[
\varphi:(S, \gamma_0) \to (S', \gamma),
\]
and moreover this $\varphi$ is an orientation-preserving diffeomorphism. We then let
\begin{equation}\label{eq:Upsilon-intro}
\Upsilon_{\eta}(v) = (\sigma, \widetilde{v}): = (\varphi^*\gamma, v \circ G^{-1} \circ \varphi),
\end{equation}
and observe that $(\sigma, \widetilde{v})$ lies in $\cM'$. Tracing the definitions, we find $G^{-1} \circ \varphi$ to be a conformal map from $(S, \sigma)$ to $(S, v^*g_M + \eta \gamma_0)$, and a straightforward computation gives~\eqref{eq:almost-conformal-intro} provided $\eta$ is small enough.

The idea of harmonic replacement, which we recall in Section \ref{sec:harmonic-replacement}, originated with Colding and Minicozzi \cite{Colding-Minicozzi08b}. As a direct consequence of their work, there exists a threshold $\ep_0 > 0$, depending only on $M$, with the following property. Given $(\sigma, v) \in \cM'$, if $\fB = \{B_{j}\}_{j \in J}$ is a finite, disjoint collection of geodesic disks in $(S, \sigma)$ satisfying 
\[
\sum_{j \in J} \int_{B_j}|\nabla v|_{\sigma}^2 \vol_{\sigma} < \ep_0,
\]
then for each $j \in J$, among maps into $M$ that agree with $v$ on $\partial B_j$, there is a unique $\sigma$-energy minimizer. Moreover, if we replace $v$ with the said energy minimizer on each $B_j$, and leave it unchanged elsewhere, then the resulting map $\cR(\sigma, v, \fB)$ still lies in $C^0 \cap W^{1, 2}$, and there holds the following \emph{convexity estimate}:
\begin{equation}\label{eq:convexity-intro}
\frac{1}{4}\int_{S} |\nabla v - \nabla \cR(\sigma, v, \fB)|_{\sigma}^2 \vol_{\sigma} \leq E(\sigma, v) - E(\sigma, \cR(\sigma, v, \fB)).
\end{equation}
We refer to the right-hand side of~\eqref{eq:convexity-intro} as the \textit{energy drop}, with its dependence on $\sigma$, $v$ and $\fB$ emphasized as needed.

Similar to \cite{Colding-Minicozzi08b}, to extend harmonic replacement of individual pairs in $\cM'$ to an \emph{energy decreasing process} that would serve as a kind of gradient flow when applied to a continuous family of pairs in $\cM'$, say of the form $\{(\sigma_{t}, v_{t})\}_{t \in I^m}$, the key step is to establish that the maximal amount of energy drop when performing replacement on $(\sigma_t, v_t)$ is an upper semi-continuous function of $t$ (Proposition \ref{prop:e-semicontinuous}). Compared to the situation in \cite{Colding-Minicozzi08b}, an added difficulty we have is that the metric on $S$ varies with $t$, and we need a well-defined procedure to transplant collections of geodesic disks in $(S, \sigma_t)$ across different values of $t$. The approach in \cite{Zhou10,Zhou17b}, which we adopt here, is to go through the universal cover $\HH$. Specifically, using again the results of Ahlfors and Bers on the equation~\eqref{eq:Beltrami-intro}, we single out for each $t$ a covering map $p_t:\HH \to S$, and show that lifts of geodesic disks in $(S, \sigma_t)$ along $p_{t}$ also get projected to geodesic disks for nearby values of $t$. (See especially Definition~\ref{defi:canonical-lift}, Remark~\ref{rmk:upstairs-and-downstairs}, and Corollary~\ref{coro:replacement-downstairs-continuity}.) The required upper semi-continuity is then proved via a contradiction argument.

With this at hand, we follow the steps in \cite{Colding-Minicozzi08b}, constructing the energy decreasing process locally and then patching the pieces together by suitably scaling the disks on which replacement occurs. The outcome is a new family $\{(\sigma_t, \widehat{v}_t)\}_{t \in I^m}$ which is homotopic to the original one in a certain sense, and has the property that, at any point $t \in I^m$ where the area of the original map $v_{t}$ is large,
\begin{equation}\label{eq:further-drop-controlled-outline}
\text{(energy drop of any further replacement on $(\sigma_t, \widehat{v}_t)$)} \leq \Theta(E(\sigma_t, v_t) - E(\sigma_t, \widehat{v}_t)),
\end{equation}
where $\Theta$ is an increasing continuous function with $\Theta(0) = 0$ (Proposition \ref{prop:energy-decreasing}). It turns out that, in the context of min-max constructions, this last property leads to compactness results for min-max sequences when combined with the convexity estimate \eqref{eq:convexity-intro} and the regularity theory of harmonic maps.

We next explain the steps in the actual proof of Theorem \ref{thm:min-max-existence}. The main line of reasoning is presented in Section \ref{sec:proof-of-main}. 
\vskip 1mm
\begin{itemize}
\item We begin by taking a minimizing sequence $(\bw_n)$ in the given class $[\bv_0]$, say with 
\[
\sup_{t \in \Int(I^{m})} A(\bw_n(t)) < \cW([\bv_0]) + \frac{1}{8n}.
\]
Similar to~\cite{Colding-Minicozzi08b}, to upgrade from $C^0 \cap W^{1, 2}$ to the $C^2$-regularity needed for conformal reparametrization, and to create room for harmonic replacement, we apply to each family $\{\bw_n(t)\}_{t \in I^m}$ a suitable mollification followed by a cut-and-paste procedure (Proposition \ref{prop:area-min-seq}). In fact, to remain in the class $[\bv_0]$, we leave $\bw_n(t)$ unchanged when $t$ is near $\partial I^m$ and only carry out the above operations away from $\partial I^m$. Some care is then required to transition continuously between these two regions in such a way that the area stays controlled, so that the resulting families still constitute a minimizing sequence. 
\vskip 1mm
\item Next we apply the map $\Upsilon_{\eta}$ (with $\eta$ appropriately chosen depending on $n$), and again a key point is to continuously ``turn off'' its effect near $\partial I^m$ in a way that keeps area under control. For this purpose, alongside $\Upsilon_{\eta}$ we also construct a map
\begin{equation}\label{eq:Xi-intro}
\Xi_{\eta}:[0, 1] \times C^2(S; M) \to C^1(S; M),
\end{equation}
with the property that $\Xi_{\eta}(\cdot, v)$ is, in the notation of equation \eqref{eq:Upsilon-intro}, a path connecting $v$ to $\widetilde{v}$ along which the area remains constant (Proposition \ref{prop:conformal-reparametrization}). The construction is inspired by an argument in~\cite{Earle-McMullen1988}. Still in the notation from the paragraph containing~\eqref{eq:Upsilon-intro}, we first lift the identity map $\id: (S, \gamma_0) \to (S, \sigma)$ to a quasiconformal map $\widetilde{f}: \HH \to \HH$, and note the standard fact, based primarily on the homotopy between $\varphi^{-1} \circ G$ and $\id$, that 
\[
\widetilde{f} \circ \gamma \circ \widetilde{f}^{-1} = w\circ \gamma \circ w^{-1},\quad\text{for all }\gamma \in \Gamma.
\]
In other words, $\widetilde{f}$ intertwines the actions of $\Gamma$ and $w\cdot \Gamma \cdot w^{-1}$. Then, by solving~\eqref{eq:Beltrami-intro} with a suitable $1$-parameter family of $\mu$, and adjusting the resulting maps using the famous Douady--Earle extension \cite{Douady-Earle1986} (see also Section \ref{subsec:DE-extension}), we obtain a path of quasiconformal $C^{1}$-diffeomorphisms leading from $w$ to $\widetilde{f}$, such that each map along the way intertwines the actions of $\Gamma$ and $w\Gamma w^{-1}$. Consequently, the path descends to give a family of $C^{1}$-diffeomorphisms $\{H_s\}_{s \in [0, 1]}$ leading from the map $\varphi^{-1} \circ G$ to the identity map of $S$. We then set 
\[
\Xi_{\eta}(s, v) = v \circ (G^{-1}\circ \varphi) \circ H_s.
\]
The map $\Xi_{\eta}$, in turn, is what is needed to achieve the ``turning off'' described at the beginning of this paragraph (see the proof of Proposition~\ref{prop:good-min-max-seq}, especially~\eqref{eq:energy-decreasing-homotopy-defi}). We note also that, throughout this construction, we rely on several \textit{a priori} estimates on quasiconformal maps (see for instance Propositions \ref{prop:qc-smooth-dependence} and \ref{prop:metric-inverse-convergence}) to guarantee the continuity of $\Xi_{\eta}$ and $\Upsilon_{\eta}$.
\vskip 1mm
\item From the previous two steps, we obtain a minimizing sequence $(\bw_{1, n})$ in $[\bv_0]$ and, for each $n$, a family of metrics $\{\bsig_n(t)\}_{t \in I^m}$ in $\Met_{-1}^{*}$, so that, roughly speaking, when $t$ is away from $\partial I^m$ there holds
\begin{equation}\label{eq:almost-conformal-outline}
E(\bsig_n(t), \bw_{1, n}(t)) < A(\bw_{1, n}(t)) + \frac{1}{4n} < \cW([\bv_0]) + \frac{3}{4n}.
\end{equation}
Next we run the energy decreasing process on each family $\{(\bsig_{n}(t), \bw_{1, n}(t))\}_{t \in I^m}$ in the minimizing sequence, and denote the result by $\{(\bsig_n(t), \widehat{\bw}_{1,n}(t))\}_{t \in I^m}$. By construction, the sequence $(\widehat{\bw}_{1, n})$ still lies in $[\bv_0]$. In particular we may choose $(t_n)$ so that, say,
\[
A(\widehat{\bw}_{1, n}(t_n)) \geq \cW([\bv_0]) - \frac{1}{4n}.
\]
For each $n$, we can in fact arrange for $t_n$ to be sufficiently far from the boundary, so that \eqref{eq:almost-conformal-outline} holds, and we obtain the following string of inequalities:
\[
\begin{split}
\cW([\bv_0]) - \frac{1}{4n} \leq\ & A(\widehat{\bw}_{1,n}(t_n)) \quad \text{(choice of $t_n$)} \\
\leq\ & E(\bsig_n(t_n), \widehat{\bw}_{1, n}(t_n))\quad  \text{(energy bounds area)}\\
\leq\ & E(\bsig_n(t_n), \bw_{1, n}(t_n)) \quad \text{(replacement decreases energy)}\\
\leq\ & A(\bw_{1, n}(t_n)) + \frac{1}{4n} < \cW([\bv_0]) + \frac{3}{4n} \quad \text{(inequality \eqref{eq:almost-conformal-outline}).}
\end{split}
\]
From this we infer that $\lim_{n \to \infty}A(\widehat{\bw}_{1, n}(t_n))= \cW([\bv_0])$, that 
\begin{equation}\label{eq:limit-conformal-outline}
\lim_{n \to \infty}\big[E(\bsig_n(t_n), \widehat{\bw}_{1, n}(t_n)) - A(\widehat{\bw}_{1, n}(t_n))\big]= 0,
\end{equation}
and that
\begin{equation}\label{eq:drop-squeezed-outline}
\lim_{n \to \infty}\big[E(\bsig_n(t_n), \bw_{1, n}(t_n)) - E(\bsig_n(t_n), \widehat{\bw}_{1, n}(t_n)) \big] = 0.
\end{equation}
\vskip 1mm
\item Let $(\sigma_n, v_n): = (\bsig_n(t_n), \widehat{\bw}_{1, n}(t_n))$. The final step is to establish a bubble tree convergence result (Theorem \ref{thm:main-existence}). To start, it follows from \eqref{eq:drop-squeezed-outline}, \eqref{eq:further-drop-controlled-outline}, and the convexity estimate \eqref{eq:convexity-intro} that $v_n$ is ``almost harmonic'' in the sense that on any disjoint collection of geodesic disks in $(S, \sigma_n)$ where $v_n$ has small energy, it is $W^{1, 2}$-close to a collection of $\sigma_n$-energy minimizers. Now, by standard results on how hyperbolic metrics on surfaces can degenerate~\cite[Chapter IV]{Hummel}, we split $(S, \sigma_n)$ into the union of a nice region where the metrics $\sigma_n$ are converging smoothly up to diffeomorphisms, and possibly several collar regions around shrinking geodesic loops. Within both types of regions, using the almost harmonic property just noted, it is fairly standard to identify subregions where bubbles develop, as well as neck regions connecting them. The property \eqref{eq:limit-conformal-outline}, and a subtle estimate for almost harmonic maps on long cylinders due to Colding and Minicozzi \cite[Proposition B.19]{Colding-Minicozzi08b}, imply that neck regions carry no energy in the limit (Proposition \ref{prop:neck}), and we get that $\cW([\bv_0])$ is accounted for by a harmonic map from a closed Riemann surface of genus at most $g$, together with at most finitely many harmonic $2$-spheres. Moreover, both the base map and the bubbles are weakly conformal by \eqref{eq:limit-conformal-outline}. (For the bubbles this also follows from the well-known Hopf differential argument.) This concludes the proof of Theorem \ref{thm:min-max-existence}.
\end{itemize}

\subsection{Context}

The construction of closed geodesics or minimal submanifolds by variational methods is an endeavor with a long history. In general, direct minimization could lead to trivial solutions, and one seeks instead non-minimizing critical points of the length and area functionals. Results of this type can be traced back to the work of G. D. Birkhoff~\cite{Birkhoff17} in the 1910s, who established a min-max principle to produce immersed closed geodesics in any Riemannian manifold diffeomorphic to $S^2$. The idea was extended by Lusternik and Schnirelmann~\cite{Lyusternik-Snirelman47} to provide an approach to finding, as conjectured by Poincar\'e~\cite{Poincare05}, at least three simple closed geodesics in any Riemannian $2$-sphere, a program subsequently completed by Grayson \cite{Grayson89}. On a related front, the min-max principle was also an inspiration for the development of Morse theory; see~\cite[page 921]{Bott80}. 

Motivated by these works, Almgren \cite{Almgren62, Almgren65} in the 1960s developed a framework for applying min-max methods to the area functional based on geometric measure theory, and used it to prove the existence of closed minimal subvarieties of any dimension and codimension in a given closed Riemannian manifold. These minimal varieties belong to the class of integer multiplicity rectifiable varifolds, which are roughly speaking locally Lipschitz submanifolds, and the next step was to improve their regularity. Breakthrough in the codimension-one case came in 1981, when Pitts \cite{Pitts81} showed that these minimal varieties are in fact smoothly embedded minimal hypersurfaces for ambient dimensions between 3 and 6. The approach had two essential components: the first was the partial regularity of area-minimizing hypersurfaces (see~\cite[Chapter 5]{Federer}), which combined the work of many, Almgren included, and the second were a priori estimates for stable minimal hypersurfaces in these ambient dimensions due to Schoen, Simon, and Yau~\cite{Schoen-Simon-Yau75}. 
Schoen and Simon~\cite{Schoen-Simon81} then extended the a priori estimates to higher dimensions by allowing the stable minimal hypersurfaces to have mild singularities, thereby establishing in all dimensions that min-max minimal hypersurfaces have the same regularity as area-minimizing hypersurfaces. We note that the Almgren--Pitts min-max theory in codimension one has since led to many remarkable results on the existence of minimal hypersurfaces, as well as their application to other problems in geometry. See for instance~\cite{Marques-Neves14, Marques-Neves17, Song18, Wang-Zhou23-four-spheres}, and the surveys~\cite{Marques-Neves21,Zhou-ICM22}.

Regarding the regularity question in higher codimensions, a celebrated result of Almgren~\cite{Almgren00} established that the singular set of an area-minimizing integral current has codimension at least two. (See also the recent new and simpler proof by De Lellis and Spadaro~\cite{DeLellis-Spadaro14, DeLellis-Spadaro16a, DeLellis-Spadaro16b, DeLellis16}.) Moreover, when the minimizing current has dimension two, the singular set was shown to be discrete by Chang \cite{Chang-XD88} and De Lellis--Spadaro--Spolaor~\cite{DeLellis-Spadaro-Spolaor17, DeLellis-Spadaro-Spolaor18, DeLellis-Spadaro-Spolaor20}. It is conjectured that, as in the codimension-one case, Almgren's min-max varifolds also have the same regularity. However, a major challenge for generalizing the approach in~\cite{Pitts81} to higher codimensions is that there is as yet no analogue of the a priori estimates of Schoen--Simon--Yau \cite{Schoen-Simon-Yau75} or Schoen--Simon \cite{Schoen-Simon81}, even in the case where the submanifolds in question are two-dimensional. 

In the two-dimensional case specifically, there exists another effective method for constructing minimal surfaces, which originated in the solution by Douglas, and independently Rad\'o, to the classical Plateau problem. In this approach, one works with Sobolev mappings from a given surface into the ambient manifold, and seeks to produce non-constant weakly conformal harmonic maps, that is, non-trivial critical points of the Dirichlet integral with equal area and energy, the key fact being that such maps parametrize immersed minimal surfaces away from isolated singularities known as branch points. Min-max constructions in this context was first carried out in the seminal work of Sacks and Uhlenbeck~\cite{Sacks-Uhlenbeck81}. Under the assumption that ambient manifold has nontrivial $k$-th homotopy group for some $k \geq 2$, they proved that there must exist a branched immersed minimal $S^2$. The techniques of Sacks--Uhlenbeck greatly influenced the subsequent development of geometric analysis, and are still widely in use. On the other hand, for the purpose of producing minimal surfaces of higher genus, it appears that their method of perturbing the Dirichlet energy will have to be complemented by new ideas. The main issue is that harmonic maps from higher-genus surfaces are not automatically weakly conformal, which in turn reflects the fact that there is a positive-dimensional space of conformal structures on such surfaces. We mention that, in the case where the ambient manifold contains an incompressible surface, this difficulty was overcome by Schoen--Yau~\cite{Schoen-Yau1979}, and independently Sacks--Uhlenbeck~\cite{Sacks-Uhlenbeck82}, who obtained a branched minimal immersion by minimizing the Dirichlet energy first in the space of mappings, and then in the moduli space of conformal structures. It is however not immediately clear how this method can be extended to min-max constructions.

In our approach, as outlined in the previous section, we work with the mapping and the conformal structure at the same time, rather than in succession. The two main ingredients mentioned earlier, namely the harmonic replacement due to Colding--Minicozzi~\cite{Colding-Minicozzi08b}, and the conformal reparametrization coming from the work of the second named author~\cite{Zhou10,Zhou17b}, serve respectively to provide the compactness necessary for extracting convergent sequences, and to ensure that the limiting object is weakly conformal. In addition, our method also incorporates elements of geometric measure theory, as the link between $C^0 \cap W^{1, 2}$-maps and integral currents, also mentioned in the previous section, allows us to apply the deep results of Almgren~\cite{Almgren62} to ascertain that the family of sweepouts obtained from multisections of closed $4$ and $5$-manifolds has positive width. At this point, we also note the viscosity method developed by Rivi\'ere and Pigati~\cite{Riviere17, Riviere17b, Pigati-Riviere20a, Pigati-Riviere20b}, which returns to the perturbation idea, but applies it directly to the area functional by adding a small multiple of the $L^{2p}$-integral of the second fundamental form. Performing min-max constructions in this setting yields mappings which parametrize stationary, integer rectifiable varifolds in a suitable sense, and their regularity result~\cite{Riviere17b,Pigati-Riviere20a} is a delicate combination of mapping theory and geometric measure theory.

\subsection{Notation}\label{subsec:notation}
Let $S$ be a closed, oriented surface, and $M$ a closed, connected, oriented Riemannian $n$-manifold. We assume that $M$ is isometrically embedded in some Euclidean space $\RR^q$, let $\cV$ be a tubular neighborhood, and write $\Pi: \cV \to M$ for the nearest-point projection. Without loss of generality, we may assume that $\Pi$ has bounded derivatives of all orders on $\cV$. As usual, the Sobolev space $W^{1, 2}(S; M)$ is defined to be
\[
W^{1, 2}(S; M) = \{ u\in W^{1, 2}(S; \RR^q)\ |\ u(x) \in M \text{ for almost every }x \in S\}.
\]
A norm on $W^{1, 2}(S; \RR^q)$ compatible with its topology is given by 
\[
\|u\|_{2}^2 + \|du\|_{2}^2 := \int_{S} |u|^2 \vol_{\gamma_0} +\int_{S} |du|_{\gamma_0}^2 \vol_{\gamma_0},
\]
where $\gamma_0$ is a choice of background metric on $S$. Since $S$ is compact, any two metrics yield equivalent norms. Next, given $u \in W^{1, 2}(S; M)$, recall that its \textit{mapping area} is by definition
\begin{equation}\label{eq:mapping-area-defi}
A(u) := \int_{S} |du(e_1) \wedge du(e_2)| \vol_{\gamma_0} = \int_{S}\sqrt{|du(e_1)|^2|du(e_2)|^2 - \bangle{du(e_1), du(e_2)}^2} \vol_{\gamma_0},
\end{equation}
while its \textit{Dirichlet energy} with respect to a given metric $\gamma$ on $S$ is
\begin{equation}\label{eq:Dirichlet-energy-defi}
E(\gamma, u) = \frac{1}{2}\int_{S} |d u|_{\gamma}^2 \vol_{\gamma}= \frac{1}{2}\int_{S} |du(e_1)|^2 + |du(e_2)|^2 \vol_{\gamma}.
\end{equation}
In both~\eqref{eq:mapping-area-defi} and~\eqref{eq:Dirichlet-energy-defi}, $e_1, e_2$ is any local orthonormal frame with respect to the metric used to define the volume form. Also, the Dirichlet energy depends only on the conformal class of the domain metric, whereas the mapping area does not depend on the choice of domain metric at all. Given $u$ and $\gamma$ as above, it is an elementary fact that 
\begin{equation}\label{eq:A-E-integrand}
|du(e_1) \wedge du(e_2)|  \leq \frac{|du(e_1)|^2 + |du(e_2)|^2}{2},
\end{equation}
where $e_1, e_2$ is any local $\gamma$-orthonormal frame, and both sides are well-defined $L^{1}$-functions on $S$. As a result, for all measurable subset $K \subset S$, there holds
\begin{equation}\label{eq:A-E-integrated}
A(u; K) \leq E(\gamma, u; K),
\end{equation}
with equality holding if and only if we have equality in~\eqref{eq:A-E-integrand} almost everywhere on $K$. Also, the difference $E(\gamma, u; K) - A(u; K)$ depends in a non-decreasing manner on the set $K$. Finally, we note that both $A(\cdot)$ and $E(\gamma, \cdot)$ are continuous with respect to strong convergence in $W^{1, 2}(S; M)$. For the Dirichlet energy this is because $S$ is compact, while for the mapping area this follows, for instance, from~\cite[page 2565, estimate (A.4)]{Colding-Minicozzi08b}.

We next collect some standard notations in geometric measure theory, and refer to \cite{Simon1983}, \cite{Almgren62} and \cite[Section 2.1]{Pitts81} for further material. We denote by $\bI_{k}(M)$ the space of $k$-dimensional integral currents in $\RR^q$ with support in $M$. By $\cZ_{k}(M)$ we mean the space of integral currents $T\in\bI_{k}(M)$ with $\partial T=0$. The flat norm \cite[Section 31]{Simon1983} and mass norm \cite[Section 26]{Simon1983} on $\bI_k(M)$ are denoted, respectively, by $\cF$ and $\bM$. Recall in particular that the flat norm of an integral current $T \in \bI_{k}(M)$ is given by
\begin{equation}\label{eq:flat-norm-defi}
\cF(T) = \inf\{\bM(R) + \bM(Q)\ |\ R \in \bI_{k+1}(M),\ Q \in \bI_{k}(M),\ T = \partial R + Q\}.
\end{equation}
Unless otherwise stated, we equip $\bI_k(M)$ and $\cZ_k(M)$ with the topology induced by the flat norm. Also, $\bI_*(M) = \bigoplus_{k}\bI_k(M)$ denotes the chain complex of integral currents in $M$.

Turning to notation pertaining to Riemann surfaces, by $\CC$, $\HH$, and $\bB$ we mean, respectively, the complex plane, the (open) upper half-plane, and the (open) unit disk. The group of biholomorphic maps from $\HH$ to itself is denoted $\Aut(\HH)$, and similarly for $\bB$. Recall that $\Aut(\HH)$ consists precisely of the real fractional linear transformations, which have the form
\[
\gamma(z) = \frac{az + b}{cz + d},
\]
for some $a, b, c, d\in \RR$ with $ad - bc = 1$. Furthermore, letting $g_{\euc}$ be the Euclidean metric $|dz|^2$ on $\CC$, and defining on $\HH$ the hyperbolic metric
\begin{equation}\label{eq:g-hyp-definition}
g_{\hyp} := \frac{|dz|^2}{(\im z)^2} = \frac{g_{\euc}}{(\im z)^2},
\end{equation}
it is a classical fact that $\Aut(\HH)$ coincides with the group of orientation-preserving isometries of $(\HH, g_{\hyp})$. A closed Riemann surface is a pair $\Sigma = (S, [\gamma])$ consisting of an underlying closed, oriented, $2$-dimensional smooth manifold $S$, and a conformal class of Riemannian metrics $[\gamma]$ on $S$. The orientation and conformal class uniquely determines a complex structure $J:TS \to TS$ by the relation
\begin{equation}\label{eq:J-from-c}
\gamma(JV, W) = \vol_{\gamma}(V, W), \quad \text{ for all }V, W \in \Gamma(TS).
\end{equation}
Notice that $J$ indeed depends only on the conformal class of $\gamma$. Finally, assuming that the surface $S$ has genus at least $2$, then the uniformization theorem gives a holomorphic covering map $p: \HH \to \Sigma$ with deck transformation group contained in $\Aut(\HH)$, and $g_{\hyp}$ descends along $p$ to the unique metric in $[\gamma]$ with constant curvature $-1$.

\subsection{Organization}\label{subsec:organization}
Section~\ref{sec:set-up} is devoted to the proof of Theorem~\ref{thm:existence-of-sweepout}, while the remaining sections center around Theorem \ref{thm:min-max-existence}. Section \ref{sec:Teichmuller} recalls three equivalent models of the Teichm\"uller space of surfaces with genus $g > 1$, involving respectively quasiconformal maps, Fuchsian models, and hyperbolic metrics. With the help of this material, in Section \ref{sec:conformal-reparametrization} we construct the conformal reparametrization map \eqref{eq:Upsilon-intro} and the homotopy \eqref{eq:Xi-intro}, while in Section \ref{sec:harmonic-replacement} we recall the notion of harmonic replacement and define the energy decreasing process, which satisfies~\eqref{eq:further-drop-controlled-outline} among other properties. Finally, in Section \ref{sec:proof-of-main} we give the proof Theorem \ref{thm:min-max-existence}. 

Appendix~\ref{appendix:patching-adapted} fleshes out the construction of sweepouts sketched in the part of Section~\ref{subsec:multisections} after Theorem~\ref{thm:existence of multisection}. In Appendix~\ref{appendix:quasiconformal} we derive, using only standard methods, some convergence results for quasiconformal maps. Appendices \ref{appendix:estimates-iterated-replacement} and \ref{appendix:semi-continuity} concern estimates that are crucial for applying harmonic replacement to families of maps. In Appendix \ref{appendix:operations} we construct the mollification and cut-and-paste procedures that, in the proof of Theorem \ref{thm:min-max-existence}, serve to prepare for the application of conformal reparametrization and harmonic replacement. Appendix~\ref{appendix:long-cylinders} contains two technical results needed for the neck analysis in the proof of Theorem~\ref{thm:min-max-existence}.

\subsection*{Acknowledgments} 
X.Z. would like to than Ian Agol for explaining to him the work of Gay and Kirby \cite{Gay-Kirby2016}, as well as Tara Holm and Inna Zakharevich for helpful discussions. X.Z. acknowledges the support by NSF grants DMS-1945178, DMS-2506717 and a grant from the Simons Foundation. 

\section{Setting up the min-max problem}\label{sec:set-up}
An overview of each of the three sections to follow is already given in Section~\ref{subsec:statements} after the statement of Theorem~\ref{thm:existence-of-sweepout}, the proof of which we complete at the end of Section~\ref{subsec:multisections}.

\subsection{Mappings and currents}\label{subsec:mappings-and-currents}
We first observe that the current $\current{S}$ represented by the oriented surface $S$ can be pushed forward by maps $S \to M$ of class $C^0 \cap W^{1, 2}$ to yield integral $2$-cycles in $M$. Given a smooth, compactly supported $2$-form $\omega$ on $\RR^{q}$ and a map $u = (u^{1},\cdots, u^{q}) \in  C^0 \cap W^{1, 2}(S; M)$, regarded as going into $\RR^q$, we define
\begin{equation}\label{eq:u-push-defi}
\bangle{u_{\#}\current{S}, \omega} := \int_{S} \omega_u(du(e_1), du(e_2)) \vol_{\gamma_0} = \int_{S} \sum_{i,j}(\omega_{ij})_{u} du^{i}(e_1) du^{j}(e_2) \vol_{\gamma_0},
\end{equation}
where $\omega_u$ stands for $\omega \circ u$, and $e_1, e_2$ is any oriented local $\gamma_0$-orthonormal frame, its choice being irrelevant. In particular, the function being integrated is a well-defined object on $S$ of class $L^1$. Note also that changing the choice of metric $\gamma_0$ does not affect the above definition. 

\begin{lemm}\label{lemm:mapping-integral}
In the above notation, $u_{\#}\current{S}$ is an integral $2$-cycle, and we have
\begin{equation}\label{eq:mass-area}
\bM(u_{\#}\current{S}) \leq A(u).
\end{equation}
\end{lemm}
\begin{proof}
This is most likely a known result, but we give a proof anyway. Rewriting the last expression in~\eqref{eq:u-push-defi} as
\[
\begin{split}
\omega_u(du(e_1), du(e_2)) =\ &   \sum_{i < j}(\omega_{ij})_u\big(du^i(e_1) du^j(e_2) - du^j(e_1) du^i(e_2) \big),
\end{split}
\]
we obtain, by the Cauchy-Schwarz inequality, the standard fact that
\begin{equation}\label{eq:duality}
\big| \omega_u(du(e_1), du(e_2))\big| \leq |\omega_u| \big| du(e_1)\wedge du(e_2) \big|,
\end{equation}
and thus, in view of~\eqref{eq:mapping-area-defi}, we have
\[
|\bangle{u_{\#}\current{S}, \omega}| \leq \|\omega_u\|_{\infty; S} \cdot A(u) \leq \|\omega\|_{\infty; \RR^q} \cdot A(u),
\]
which implies that $u_\#\current{S}$ is a continuous linear functional on smooth, compactly supported $2$-forms on $\RR^q$, and so defines a $2$-current. The estimate~\eqref{eq:mass-area} also follows.

To show that $u_{\#}\current{S}$ is an integral cycle, take a sequence $(u_n)$ in $C^\infty(S; M)$ such that
\begin{equation}\label{eq:approximation-for-current}
\|du_n - du\|_{2} + \|u_n - u\|_{\infty} \to 0\quad \text{as }n \to \infty.
\end{equation}
A straightforward computation with the help of~\eqref{eq:duality} shows that
\[
\begin{split}
&\big| \omega_{u_n}\big( du_n(e_1), du_n(e_2) \big) - \omega_{u}\big( du(e_1), du(e_2) \big) \big| \\
&\leq\ |\omega_{u_n} - \omega_u| |du_n|_{\gamma_0}^2 + |\omega_u| |du_n - du|_{\gamma_0} |du_n|_{\gamma_0} +  |\omega_u| |du|_{\gamma_0} |du_n - du|_{\gamma_0}.
\end{split}
\]
Integrating this over $S$ with respect to $\vol_{\gamma_0}$ and using H\"older's inequality gives
\[
\Big|\bangle{(u_n)_{\#}\current{S} - u_{\#}\current{S}, \omega}\Big| \leq \|\omega_{u_n} - \omega_{u}\|_{\infty} \cdot \|du_n\|_2^2 + \|\omega_{u}\|_{\infty}\cdot \|du_n - du\|_2 (\|du_n\|_2 + \|du\|_2).
\]
Consequently
\begin{equation}\label{eq:current-convergence}
(u_n)_{\#}\current{S} \to u_{\#}\current{S}\quad \text{as currents}.
\end{equation} 
On the other hand, since $u_n$ is smooth, each $(u_n)_{\#}\current{S}$ is an integer-multiplicity rectifiable current (see~\cite[Remark 27.2(3)]{Simon1983}, or~\cite[4.1.28 and 4.1.30]{Federer}). Observing furthermore that $\partial \big((u_n)_{\#}\current{S}\big) = 0$ and that, by~\eqref{eq:mass-area},~\eqref{eq:approximation-for-current} and the estimate (A.4) from~\cite[page 2565]{Colding-Minicozzi08b}, 
\[
\bM((u_n)_{\#}\current{S}) \leq A(u_n) \to A(u)\quad \text{as }n \to \infty,
\]
we may invoke the compactness theorem (\cite[Sections 27 and 32]{Simon1983}) to conclude that a subsequence of $(u_n)_{\#}\current{S}$ converges to an integral $2$-cycle. This together with~\eqref{eq:current-convergence} shows that $u_{\#}\current{S}$ is an integral $2$-cycle.
\end{proof}
For the next result we choose any $\delta > 0$ such that $\{y \in \RR^q\ |\ \dist(y, M) \leq 3\delta \}$ is contained in the tubular neighborhood $\cV$ of $M$ fixed in Section \ref{subsec:notation}. Recall also that $\Pi:\cV \to M$ denotes the nearest-point projection.
\begin{lemm}\label{lemm:segment}
Suppose $u, v \in C^{0} \cap W^{1, 2}(S; M)$ are such that $\|u - v\|_{\infty} < \delta$, and define $F = F_{u, v} : [0, 1] \times S \to M$ by
\[
F(t, \cdot) = \Pi(tv(\cdot) + (1 - t) u(\cdot)) , \text{ for }t \in [0, 1].
\]
Given a smooth, compactly supported $3$-form $\alpha$ on $\RR^{q}$, define
\begin{equation}\label{eq:path-current-defi}
\begin{split}
&\bangle{F_{\#}\current{[0, 1] \times S}, \alpha} =\int_{[0, 1] \times S} \alpha_{F}(dF(\partial_t), dF(e_1), dF(e_2)) \, dt \wedge \vol_{\gamma_0},
\end{split}
\end{equation}
where as before $\alpha_F$ means $\alpha \circ F$, and $e_1, e_2$ is any choice of oriented local orthonormal frame on $S$ with respect to $\gamma_0$. Then $F_{\#}\current{[0, 1] \times S}$ is an integral $3$-current, with boundary given by
\begin{equation}\label{eq:path-current-boundary}
\partial \big( F_{\#}\current{[0, 1] \times S}\big) = v_{\#}\current{S} - u_{\#}\current{S},
\end{equation}
and with mass satisfying
\begin{equation}\label{eq:path-mass}
\bM(F_{\#}\current{[0, 1] \times S}) \leq L^3 \cdot \|u - v\|_{\infty} \cdot \big( \|du\|_2^2 + \|dv\|_2^2 \big),
\end{equation}
where $L = \sup_{y \in \cV_{2\delta}}|(d\Pi)_y|$.
\end{lemm}
\begin{proof}
Again the $3$-form being integrated in~\eqref{eq:path-current-defi} is well-defined object of class $L^{1}$ on $[0, 1] \times S$ irrespective of the choice of $\gamma_0$ and $e_1, e_2$. Notice also that
\begin{equation}\label{eq:F-derivatives}
dF(\partial_t) = (d\Pi)_{tv + (1 - t)u}(v - u),\quad dF(e_i) = (d\Pi)_{tv + (1 - t)u}(tdv(e_i) + (1 - t)du(e_i)),
\end{equation}
and recall the following elementary fact:
\[
\begin{split}
\big| \alpha_{F}(dF(\partial_t), dF(e_1), dF(e_2)) \big| \leq \ & |\alpha_F| \big|dF(\partial_t) \wedge dF(e_1) \wedge dF(e_2)\big| \\ \leq\ & |\alpha_F| \big|dF(\partial_t) \big| \big| dF(e_1) \big| \big| dF(e_2)\big|,
\end{split}
\]
which combines with~\eqref{eq:F-derivatives} to give
\[
\begin{split}
\big| \alpha_{F}(dF(\partial_t), dF(e_1), dF(e_2)) \big|
\leq\ &|\alpha_F| |(d\Pi)_{tv + (1 - t)u}|^3  \cdot |v - u| \cdot \frac{|t dv + (1 - t)du|_{\gamma_0}^2}{2}\\
\leq\ & \frac{L^3}{2} |\alpha_F| |v - u|\cdot (|dv|_{\gamma_0} + |du|_{\gamma_0})^2.
\end{split}
\]
Integrating over $[0, 1] \times S$ leads to
\[
\big| \bangle{F_{\#}\current{[0, 1] \times S}, \alpha} \big| \leq L^3 \|\alpha\|_{\infty; \RR^q} \cdot \|v - u\|_{\infty; S}\cdot (\|dv\|_2^2 + \|du\|_2^2),
\]
which shows that $F_{\#}\current{[0, 1] \times S}$ is a $3$-current and also gives~\eqref{eq:path-mass}. 

Next we let $(u_n)$ and $(v_n)$ be sequences in $C^{\infty}(S; M)$ so that 
\[
\|du_n - du\|_{2} + \|u_n - u\|_{\infty} + \|dv_n - dv\|_{2} + \|v_n - v\|_{\infty} \to 0 \text{ as }n \to \infty.
\]
Then eventually $\|u_n - v_n\|_{\infty} < \delta$, so we may define
\[
F_n(t, \cdot) =  \Pi(tv_n(\cdot) + (1 - t) u_n(\cdot)),\ \ t \in [0, 1],
\]
and observe that $F_n$ converges uniformly on $[0, 1] \times S$ to $F$. Furthermore, by~\eqref{eq:F-derivatives} and its counterpart for $F_{n}$, we see that
\[
\sup_{[0, 1] \times S} |dF_n(\partial_t) - dF(\partial_t)|
 + \int_{[0, 1] \times S} \sum_{i = 1}^2 |dF_n(e_i) - dF(e_i)|^2 \, dt \wedge \vol_{\gamma_0} \to 0 \text{ as }n \to \infty.
\]
Recalling~\eqref{eq:path-current-defi}, and arguing as in the lines between~\eqref{eq:approximation-for-current} and~\eqref{eq:current-convergence}, we deduce that
\[
(F_n)_{\#}\current{[0, 1] \times S} \to F_{\#}\current{[0, 1] \times S} \quad  \text{as currents,}
\]
from which we get~\eqref{eq:path-current-boundary} upon noting
\begin{equation}\label{eq:path-current-boundary-n}
\partial \big( (F_n)_{\#}\current{[0, 1] \times S}  \big) = (v_n)_{\#}\current{S} - (u_n)_{\#}\current{S},
\end{equation}
and recalling~\eqref{eq:current-convergence}. Now, by~\eqref{eq:path-mass}, the masses of the integral currents $(F_n)_{\#}\current{[0, 1] \times {S}} $ are uniformly bounded. By~\eqref{eq:path-current-boundary-n} and Lemma~\ref{lemm:mapping-integral}, so are the masses of $\partial \big( (F_n)_{\#}\current{[0, 1] \times {S}}  \big)$. The compactness theorem and the convergence of $(F_n)_{\#}\current{[0, 1] \times {S}}$ to $F_{\#}\current{[0, 1] \times {S}}$ then implies that the latter is an integer-multiplicity rectifiable $3$-current. Lemma~\ref{lemm:mapping-integral} together with~\eqref{eq:path-current-boundary} shows that $\partial \big( F_{\#}\current{[0, 1] \times {S}}  \big)$ is integer-multiplicity rectifiable as well.
\end{proof}
Recalling from Section~\ref{subsec:notation} the definition of the flat norm $\cF$, we have the following consequence of Lemmas~\ref{lemm:mapping-integral} and~\ref{lemm:segment}.
\begin{coro}\label{coro:C0-W12-flat}
Suppose $(u_n)$ is a sequence in $(C^0 \cap W^{1, 2})(S; M)$.
\vskip 1mm
\begin{enumerate}
\item[(a)] If $u_n$ converges in $C^0 \cap W^{1, 2}$ to some $u: {S} \to M$, then we have
\[
\lim_{n \to \infty}\cF((u_n)_{\#}\current{{S}}- u_{\#}\current{{S}} ) = 0.
\]
\vskip 1mm
\item[(b)] Assume instead that $A(u_n) \to 0$ as $n \to \infty$. Then
\[
\lim_{n \to \infty}\cF((u_n)_{\#}\current{{S}}) = 0.
\]
We emphasize that for part (b), the maps $u_n$ themselves are not assumed to converge.
\end{enumerate}
\end{coro}
\begin{proof}
For part (a), we note that Lemma~\ref{lemm:segment} is applicable for sufficiently large $n$, giving us
\[
\cF((u_n)_{\#}\current{{S}}- u_{\#}\current{{S}} ) \leq C \|u_n - u\|_{\infty} \cdot (\|du_n\|_2^2 + \|du\|_{2}^{2}) \to 0\quad \text{as }n \to \infty.
\]
For part (b), by Lemma~\ref{lemm:mapping-integral} we see that $(u_n)_{\#}\current{S}$ is a sequence of integral cycles with mass tending to $0$. The result then follows from the isoperimetric theorem (\cite[Section 30]{Simon1983}).
\end{proof}
 
Proposition~\ref{prop:classes-of-maps-to-currents} below is the main result of this section. Recall that we equip the space $\cZ_k(M)$ of integral $k$-cycles in $M$ with the topology induced by the flat norm. Fixing some $m \in \NN$, for all $\bw \in \cP = \cP(m, S)$, we define $f_{\bw}: (I^m, \partial I^m) \to (\cZ_2(M), \{0\})$ by 
\begin{equation}\label{eq:fw-defi}
f_{\bw}(t) =\left\{
\begin{array}{ll}
{\bw}(t)_{\#}\current{{S}},& \text{ for }t \in \Int(I^m),\\
0, & \text{ for } t \in \partial I^m.
\end{array}
\right.
\end{equation}
\begin{prop}\label{prop:classes-of-maps-to-currents}
In the above notation, we have the following.
\vskip 1mm
\begin{enumerate}
\item[(a)] $f_{\bw}: (I^m, \partial I^m) \to (\cZ_2(M), \{0\})$ is continuous.
\vskip 1mm
\item[(b)] If $\bw \sim \bw'$ in $\cP$ in the sense of Definition~\ref{defi:relation-on-P}, then $f_{\bw}$ and $f_{\bw'}$ represent the same class in $\pi_m(\cZ_2(M), \{0\})$. 
\vskip 1mm
\item[(c)] There exists $\mu_{m, M} > 0$ such that if $f_{\bw}$ represents a non-trivial class in $\pi_m(\cZ_2(M), \{0\})$, then $\sup_{t \in \Int(I^m)}A(\bw(t)) \geq \mu_{m, M}$.
\end{enumerate}
\end{prop}
\begin{proof}
For part (a), by the continuity of ${\bw}|_{\Int(I^m)}$ as a map into $(C^0 \cap W^{1, 2})(S; M)$ along with Corollary~\ref{coro:C0-W12-flat}(a) we see that $f_{\bw}$ is continuous on $\Int(I^m)$. On the other hand, by the boundary condition~\eqref{eq:no-area-at-boundary} along with Corollary~\ref{coro:C0-W12-flat}(b) we see that $f_{\bw}$ is continuous at points on $\partial I^m$.

For part (b), suppose $\bw \sim \bw'$ in $\cP$, so that there exists a continuous map 
\[
\bh: [0, 1] \times \Int(I^m) \to (C^0 \cap W^{1, 2})(S; M)
\]
satisfying both~\eqref{eq:cP-homotopy-1} and~\eqref{eq:cP-homotopy-2}, with $\bv$ and $\bw$ replaced respectively by $\bw$ and $\bw'$. Define $f_{\bh}:([0, 1] \times I^m, [0, 1] \times \partial I^m) \to (\cZ_2(M), \{0\})$ by 
\[
f_{\bh}(s, t) = \left\{\begin{array}{ll}
\bh(s, t)_{\#}\current{S}, & \text{ if }t \in \Int(I^m),  \\
0, & \text{ if }t \in \partial I^m.
\end{array}
\right.
\]
Then, with the help of Corollary~\ref{coro:C0-W12-flat}, we see that $f_{\bh}$ is a continuous map. Since $f_{\bh}(0, \cdot) = f_{\bw}$ and $f_{\bh}(1, \cdot) = f_{\bw'}$, we get part (b). Finally, part (c) follows from Lemma~\ref{lemm:mapping-integral}, the isoperimetric theorem~\cite[Section 30]{Simon1983}, and the work of Almgren~\cite[Theorem 8.2]{Almgren62}.
\end{proof}

\subsection{Almgren's isomorphism}\label{subsec:Almgren isomorphism}
We begin by recalling the isomorphism established by Almgren between certain homotopy groups of the flat chain space and the homology groups of $M$ (\cite{Almgren62}; see also \cite{Almgren65, Pitts81, Zhou15, Marques-Neves17}). In addition to the notation in Section~\ref{subsec:notation}, we also need the cubical complexes utilized in the discretization scheme in \cite{Almgren62, Pitts81}. For each $j \in\NN$, let $I(1, j)$ be the cubical complex on the unit interval $I=[0, 1]$ whose 1-cells and 0-cells (the latter are also referred to as vertices) are, respectively,
\[ [0, 3^{-j}], [3^{-j}, 2\cdot 3^{-j}], \cdots, [1-3^{-j}, 1]\, \text{ and }\, [0], [3^{-j}], \cdots, [1-3^{-j}], [1]. \]
We then denote by $I(m, j)$ the induced cell complex on $I^m$:
\begin{equation}\label{eq:m-dimensiona-complex}
I(m, j)=  \underbrace{I(1, j)\otimes\cdots\otimes I(1, j)}_{\text{$m$ times}}.
\end{equation}
In particular, $\alpha=\alpha_1\otimes\cdots\otimes\alpha_m$ is a $d$-cell of $I(m, j)$ if and only if each $\alpha_i$ is a cell of $I(1, j)$, and $\sum_{i=1}^m\dim(\alpha_i)=d$. We often identify a $d$-cell $\alpha$ with its support $\alpha_1\times\cdots\times \alpha_m\subset I^m$, and write the set of all $d$-cells as $I(m, j)_d$. Also, we use $I_0(m, j)$ to denote the subcomplex of $I(m, j)$ consisting of all cells lying in $\partial I^m$, and let 
\[
I_0(m, j)_d  = I_0(m, j)\cap I(m, j)_d.
\]
The boundary homomorphism $\partial: I(m, j) \to I(m, j)$ is defined in the usual way~\cite[page 266]{Almgren62}. In particular, the support of each cell of dimension $d > 0$ can be oriented in such a way that the associated currents satisfy
\[
\partial\current{\alpha} = \left\{
\begin{array}{ll}
\delta_{t} - \delta_{s}, & \text{ if }d = 1 \text{ and }\partial\alpha = t - s \text{ with }t, s \in I(m, j)_0,\\
\current{\partial\alpha}, & \text{ if }d > 1,
\end{array}
\right.
\]
where by $\delta_{t}$ we mean the delta measure supported at $t$. Also, to be precise in the second case, we first express $\partial \alpha$ as $\sum_{i = 1}^{d}(-1)^{i-1}(\beta_{i}^+ - \beta_{i}^{-})$ according to the formula above point (4) in~\cite[page 266]{Almgren62}, with each $\beta_{i}^{+}, \beta_i^{-} \in I(m, j)_{d-1}$, and $\current{\partial\alpha}$ is then understood to mean $\sum_{i = 1}^{d}(-1)^{i-1}(\current{\beta_{i}^+} - \current{\beta_{i}^{-}})$.

Given integers $0\leq k\leq n = \dim M$ and $m>0$, Almgren constructed (see especially~\cite[Section 3]{Almgren62}) an isomorphism
\[
F_A: \pi_m\big(\cZ_k(M), \{0\}\big) \to H_{k+m}(M, \ZZ),
\]
which we now describe. By the isoperimetric theorem of Federer and Fleming, there exists a small threshold $\nu = \nu(M, m) > 0$ such that for any continuous map 
\[
f: (I^m, \partial I^m) \to \big(\cZ_k(M), \{0\}\big)
\]
representing a homotopy class, as long as $j$ is large enough so that
\begin{equation}\label{eq:fineness}
\cF\big(f(t), f(t')\big)\leq \nu \text{ for any two adjacent vertices $t, t'\in I(m, j)_0$,}
\end{equation}
one may define a degree-$k$ chain map $\phi: I(m, j) \to \bI_*(M)$ by inductively filling in with \emph{isoperimetric choices}, as follows: 
\begin{enumerate}
    \item[(A1)] For each vertex $t\in I(m, j)_0$, set $\phi(t) = f(t)$;
    \item[(A2)] For each 1-cell $\alpha\in I(m, j)_1$, set $\phi(\alpha)\in \bI_{k+1}(M)$ to be an $\cF$-isoperimetric choice for $\phi(\partial\alpha)$~\cite[Definition 1.15]{Almgren62}, where by definition $\phi(\partial\alpha)=\phi(t)-\phi(s)$ if $\partial\alpha$ is expressed as $t - s$ for some $t, s \in I(m, j)_0$;
    \item[(A3)] For $2\leq d\leq m$ and each $d$-cell $\alpha\in I(m, j)_d$, set $\phi(\alpha)\in \bI_{k+d}(M)$ to be an $\bM$-isoperimetric choice for $\phi(\partial\alpha)$~\cite[Definition 1.12]{Almgren62}. Here $\phi(\partial\alpha)$ is defined in the same way as in (A2).
\end{enumerate}
Since $f(\partial I^m) = \{0\}$, the construction guarantees that $\phi(\alpha) = 0$ for all $\alpha \in I_0(m, j)$. Almgren \cite[Section 3 and Section 6]{Almgren62} then proved that the following homology class 
\[
\left[\sum_{\text{all $m$-cells }\alpha \text{ of }I(m, j)} \phi(\alpha)\right] \in  H_{k+m}(M, \ZZ)
\]
is independent of the choice of $f$ representing the homotopy class, and the choice of $j$ for which~\eqref{eq:fineness} holds. Consequently, the following map from $\pi_m\big(\cZ_k(M), \{0\}\big)$ to $H_{k+m}(M, \ZZ)$ is well-defined, and in fact Almgren showed that it is an isomorphism:
\begin{equation}\label{eq:Almgren-isomorphism} 
F_A: [f] \longmapsto \left[\sum_{\text{all $m$-cells }\alpha \text{ of }I(m, j)} \phi(\alpha)\right]. 
\end{equation}

We are interested primarily in the case when the initial representative $f: (I^m, \partial I^m) \to (\cZ_k(M), \{0\})$ arises from a Lipschitz map $I^m \times \Sigma_0 \to M$, where $\Sigma_0$ is some oriented $k$-manifold. The next two propositions together yield conditions on the Lipschitz map under which we can relate its image to $F_{A}([f])$.

\begin{prop}\label{prop:homology-class-from-Lipschitz}
Given $k, m \in \NN$ such that $k \geq 2$ and $m + k \leq n = \dim M$, let $\Sigma_0$ be a closed, oriented $k$-manifold, and suppose $\Phi: I^m \times \Sigma_0 \to M$ is a Lipschitz map satisfying
\begin{equation}\label{eq:boundary-lower-dimension}
\cH^{k}(\Phi(\{ t \} \times \Sigma_0)) = 0 , \quad \text{for all }t \in \partial I^m.
\end{equation}
Then we have the following.
\begin{enumerate}
\item[(a)] $\Phi(t, \cdot)_{\#}\current{\Sigma_0}=0$ in $\cZ_k(M)$, for all $t \in \partial I^m$. 
\vskip 1mm
\item[(b)] $\partial\Phi_{\#}(\current{I^m} \times \current{\Sigma_0}) = 0$ in $\cZ_{m-1 + k}(M)$. Also, $\cH^{m-1 + k}(\Phi(\partial I^m \times \Sigma_0)) = 0$.
\vskip 1mm
\item[(c)] The map $\widetilde{f}_{\Phi}:(I^m, \partial I^m) \to (\cZ_k(M), \{0\})$ defined by
\begin{equation}\label{eq:f-Phi-defi} 
\widetilde{f}_{\Phi}: t\mapsto \Phi(t,\cdot)_{\#}\current{\Sigma_0},
\end{equation}
is continuous, and the image of $[\widetilde{f}_{\Phi}]$ under Almgren's isomorphism $F_A$ coincides with the homology class represented by $\Phi_{\#}(\current{I^m} \times \current{\Sigma_0})$. 
\end{enumerate}
\end{prop}
\begin{proof}
Without loss of generality we assume that $\Sigma_0$ is an embedded submanifold in some Euclidean space, and equip it with the induced metric. Given $t \in \partial I^m$ and writing $\Phi_{t}$ for $\Phi(t, \cdot)$, by the assuption~\eqref{eq:boundary-lower-dimension} and the area formula, we have
\begin{equation}\label{eq:vanishing-Jacobian}
|d\Phi_t(v_1) \wedge \cdots \wedge d\Phi_t(v_k)| = 0,\quad \text{$\cH^{k}$-a.e. on $\Sigma_0$},
\end{equation}
where $v_1, \cdots, v_k$ is any oriented orthonormal basis for the tangent space to $\Sigma_0$ at the point of interest. Using this in the integral formula for $\Phi(t, \cdot)_{\#}\current{\Sigma_0}$ (see~\cite[Remark 27.2(3)]{Simon1983} or~\cite[4.1.30]{Federer}) shows that the latter vanishes, and we are done with (a). 

Moving on to parts (b) and (c), by the Lipschitz assumption on $\Phi$ and the definition of the flat norm, we see that~\eqref{eq:f-Phi-defi} defines a continuous map into $\cZ_{k}(M)$, which satisfies $\widetilde{f}_{\Phi}(\partial I^m) = \{0\}$ by part (a). To prove that $\Phi_{\#}(\current{I^m} \times \current{\Sigma_0})$ has no boundary and hence determines a homology class, and that the latter coincides with $F_A([\widetilde{f}_{\Phi}])$, we let $\nu_1, \nu_2 > 0$ be the constants from~\cite[Proposition 1.11]{Almgren62}, which depend only on $M$ in the present context since $\cZ_{k}(M)$ is non-trivial only when $k \leq n$, and choose $\eta > 0$ so that
\begin{equation}\label{eq:eta-choice}
\eta < \frac{1}{8m}\min\{\nu_1, 1\}, \quad \text{and} \quad \nu_2 \cdot (8m\eta)^{\frac{1}{n + m}} < \frac{1}{8m}.
\end{equation}
By the continuity of $\widetilde{f}_{\Phi}$ into $\cZ_{k}(M)$, together with the Lipschitz continuity of $\Phi$, for all sufficiently large $j$ we have that~\eqref{eq:fineness} holds with $f = \widetilde{f}_{\Phi}$, and that
    \begin{equation}\label{eq:mass threshold}
        \bM\big(\Phi_{\#}(\current{\alpha} \times \current{\Sigma_0}) \big)< \eta, 
    \end{equation}
for each $d\in \{1, \cdots, m\}$ and $\alpha\in I(m, j)_d$. Fixing such a $j$, we let 
\[
\phi:I(m, j) \to \bI_{*}(M)
\]
denote the degree-$k$ chain map obtained from $\widetilde{f}_{\Phi}$ via the procedure summarized after~\eqref{eq:fineness}. On the other hand, since each $\Phi_{\#}(\current{\alpha} \times \current{\Sigma_0})$ lies in $\bI_{\dim(\alpha) + k}(M)$ and has boundary equal to $\Phi_{\#}(\current{\partial\alpha} \times \current{\Sigma_0})$, we get another degree-$k$ chain map $\widetilde{\phi}: I(m, j)\to \bI_*(M)$ upon setting
    \begin{itemize}
        \item $\widetilde\phi(t) = \Phi(t, \cdot)_{\#}\current{\Sigma_0}$ for each vertex $t\in I(m, j)_0$. 
        \vskip 1mm
        \item $\widetilde\phi(\alpha) = \Phi_{\#}(\current{\alpha} \times \current{\Sigma_0})$ for each $d$-cell $\alpha\in I(m, j)_d$, for $d = 1, \cdots, m$. 
    \end{itemize}
Observe from part (a) that 
\[
\widetilde{\phi}(t) = 0,\quad \text{for all } t \in I_0(m, j)_0.
\]
On the other hand, given $\alpha \in I_0(m, j)_{d}$ with $d > 0$, note that for $\cH^{d}$-almost every $t \in \alpha$, the derivative $d(\Phi|_{\alpha \times \Sigma_0})_{(t, p)}$ exists for $\cH^{k}$-almost every $p \in \Sigma_0$, and that each such $p$ is also a point of differentiability for $\Phi_{t}$, with 
\[
d(\Phi|_{\alpha \times \Sigma_0})_{(t, p)}(v) = d\Phi_{t}(v)\quad\text{for all }v \in T_{p}\Sigma_0.
\]
Since~\eqref{eq:vanishing-Jacobian} is true for all $t \in \partial I^m$, we deduce that
\[
d(\Phi|_{\alpha \times \Sigma_0})_{(t, p)}(v_1) \wedge \cdots\wedge d(\Phi|_{\alpha \times \Sigma_0})_{(t, p)}(v_k) = 0,\quad \text{for $\cH^{d + k}$-a.e. $(t, p) \in \alpha \times \Sigma_0$},
\]
where $v_1,\cdots, v_k$ is any orthonormal basis of $T_{p}\Sigma_0$. Combining this with the definition of $\Phi_{\#}(\current{\alpha} \times \current{\Sigma_0})$ and, respectively, the area formula, we deduce that 
\[
\widetilde{\phi}(\alpha) = 0 \text{ in }\bI_{d + k}(M),\quad \text{and}\quad \cH^{d + k}(\Phi(\alpha \times \Sigma_0)) = 0.
\]
Applying this to the $(m-1)$-cells on $\partial I^m$ gives conclusion (b). In particular $\Phi_{\#}(\current{I^m} \times \current{\Sigma_0})$ determines a homology class in $H_{m + k}(M)$, which is captured by $\widetilde{\phi}$ via the relation:
\begin{equation}\label{eq:phi-tilde-image}
\left[\sum_{\text{all $m$-cells }\alpha} \widetilde\phi(\alpha)\right] = \big[\Phi_{\#}(\current{I^m} \times \current{\Sigma_0})\big].
\end{equation}

It remains to prove that this class coincides with $F_A([\widetilde{f}_{\Phi}])$, which will be accomplished by constructing a chain homotopy~\cite[Definition 2.3]{Almgren62} between $\phi$ and $\widetilde\phi$; that is, a graded homomorphism $\psi: I(m, j)\to \bI_*(M)$ of degree $k + 1$ such that
    \begin{equation}\label{eq:chain homotopy relation}
        \partial\psi(\alpha) + \psi(\partial\alpha) = \widetilde\phi(\alpha) - \phi(\alpha), \text{ for each }\alpha \in I(m, j),
    \end{equation}
    \begin{equation}\label{eq:chain homotopy boundary}
        \psi(\alpha) = 0, \text{ for each }\alpha \in I_0(m, j).
    \end{equation}
   The construction proceeds inductively:
    \begin{enumerate}[label=(\roman*)]
        \vskip 1mm
        \item For each vertex $t\in I(m, j)_0$, we let $\psi(t)=0 \in \bI_{k + 1}(M)$.
        \vskip 1mm
        \item For each 1-cell $\alpha\in I(m, j)_1$, we let $\psi(\alpha)\in\bI_{k+2}(M)$ to be an $\bM$-isoperimetric choice for $\widetilde\phi(\alpha)-\phi(\alpha)$. Such a choice always exists because, first of all, by definition $\widetilde\phi$ and $\phi$ are identical when restricted to all vertices in $I(m, j)_0$, so that
        \[\partial(\widetilde\phi(\alpha)-\phi(\alpha)) = \widetilde\phi(\partial\alpha)-\phi(\partial\alpha) =0;\] 
        secondly, since $\phi(\alpha)$ is an $\cF$-isoperimetric choice for $\phi(\partial\alpha)=\widetilde\phi(\partial\alpha)$, we have $\bM(\phi(\alpha))\leq \bM(\widetilde\phi(\alpha))$, and hence by~\eqref{eq:mass threshold} and our choice of $\eta$, 
        \[
        \bM\big(\widetilde\phi(\alpha)-\phi(\alpha)\big) \leq 2\bM\big( \widetilde{\phi}(\alpha) \big) < 2\eta < \nu_1.
        \]
        Thus the definition of $\psi(\alpha)$ makes sense. The isoperimetric inequality in~\cite[Proposition 1.11]{Almgren62} along with~\eqref{eq:eta-choice} then gives
        \[
        \begin{split}
        \bM\big(\psi(\alpha) \big) \leq\ & \nu_2 \cdot \big[ \bM\big(\widetilde\phi(\alpha)-\phi(\alpha)\big) \big]^{\frac{k + 2}{k + 1}}\\
        \leq\ & \nu_2 \cdot (2\eta)^{\frac{k + 2}{k + 1}} \leq \nu_2 \cdot (8m\eta)^{\frac{1}{n + m}} \cdot (2\eta) < \eta,
        \end{split}
        \]
        where the penultimate inequality uses $8m\eta < 1$. Also, the first inequality implies that $\psi(\alpha) = 0$ if $\alpha \in I_0(m, j)_1$.
        \vskip 1mm
        \item Suppose for some $2 \leq d \leq m$ that $\psi$ has been defined for all cells of dimension less than $d$ in $I(m, j)$ so that both~\eqref{eq:chain homotopy relation},~\eqref{eq:chain homotopy boundary} are satisfied, and that the images of these cells under $\psi$ each has mass smaller than $\eta$. Given $\alpha\in I(m, j)_d$, we let $\psi(\alpha)\in\bI_{k+d+1}(M)$ be an $\bM$-isoperimetric choice for $\widetilde\phi(\alpha)-\phi(\alpha)-\psi(\partial\alpha)$. Similar to (ii), this is possible because, first, by~\eqref{eq:chain homotopy relation} applied to $\partial\alpha$, there holds
        \[
        \begin{split}
            \partial\big(\widetilde\phi(\alpha)-\phi(\alpha)-\psi(\partial\alpha)\big) 
            & = \widetilde\phi(\partial\alpha)-\phi(\partial\alpha)-\partial\psi(\partial\alpha)\\
            & = \psi(\partial\partial\alpha) + \partial\psi(\partial\alpha) - \partial\psi(\partial\alpha) = 0;
        \end{split}
        \]
        secondly, since $\phi(\alpha)$ is an $\bM$-isoperimetric choice, and since the number of $(d-1)$-cells making up $\partial \alpha$ is at most $2m$, we have by~\eqref{eq:mass threshold} and the mass bound on $\psi(\partial\alpha)$ in the induction hypothesis that
        \[
        \begin{split}
        \bM\big( \widetilde\phi(\alpha)-\phi(\alpha)-\psi(\partial\alpha)\big) \leq 2\bM\big( \widetilde\phi(\alpha)- \psi(\partial\alpha) \big) \leq 2 ( 1 + 2m)\eta < 8m\eta < \nu_1.
        \end{split}
        \]
        Using the isoperimetric inequality and~\eqref{eq:eta-choice} as in the end of step (ii), we see that
        \[
        \begin{split}
        \bM\big( \psi(\alpha) \big) \leq\ & \nu_2 \cdot \big[\bM\big( \widetilde\phi(\alpha)-\phi(\alpha)-\psi(\partial\alpha) \big)\big]^{\frac{k + d + 1}{k + d}} \\
        \leq\ & \nu_2 \cdot (8m\eta)^{\frac{k + d + 1}{k + d}} \leq \nu_2 \cdot (8m\eta)^{\frac{1}{n+m}} \cdot (8m\eta) < \eta.
        \end{split}
        \]
        Also, for $\alpha \in I_0(m, j)_d$ we must have $\psi(\alpha) = 0$, since in this case $\widetilde\phi(\alpha)-\phi(\alpha)-\psi(\partial\alpha) = 0$. Thus, having already built~\eqref{eq:chain homotopy relation} into the definition of $\psi(\alpha)$, the inductive construction can continue, and we obtain the desired chain homotopy between $\phi$ and $\widetilde{\phi}$.
    \end{enumerate}
    \vskip 1mm
    To finish, by our choice of $j$ and the definition of $F_A$, together with the properties~\eqref{eq:chain homotopy relation} and~\eqref{eq:chain homotopy boundary} of $\psi$, we have 
    \begin{equation}\label{eq:result-from-chain-homotopy}
    F_A([\widetilde{f}_{\Phi}]) = \left[\sum_{\text{all $m$-cells }\alpha}\phi(\alpha)\right] = \left[\sum_{\text{all $m$-cells }\alpha}\widetilde\phi(\alpha)\right].  
    \end{equation}
    Combining this with~\eqref{eq:phi-tilde-image} gives $F_A([\widetilde{f}_{\Phi}]) =  \big[\Phi_{\#}(\current{I^m}\times \current{\Sigma_0})\big]$. The proof is complete.
\end{proof}

\begin{prop}\label{prop:induced-sweepout}
Let $k$, $m$, $n$, and $\Sigma_0$ be as in Proposition~\ref{prop:homology-class-from-Lipschitz}, and assume in addition that $m + k = n$. Suppose ${X} \subset \RR^{m}$ is a bounded domain such that there is a Lipschitz map $\tau: (I^{m}, \partial I^{m}) \to (\overline{X}, \partial X)$ with a Lipschitz inverse $\tau^{-1}:(\overline{X}, \partial X) \to (I^{m}, \partial I^{m})$. Furthermore, let $\Psi: \overline{{X}} \times \Sigma_0 \to M$ be a Lipschitz map that satisfies
\vskip 1mm
\begin{enumerate}
\item[(i)] $\cH^{k}(\Psi(\{x\} \times \Sigma_0)) = 0$ for all $x \in \partial {X}$. 
\vskip 1mm
\item[(ii)] $\Psi$ restricts to a $C^{1}$-diffeomorphism near some $(x_0, p_0) \in {X} \times \Sigma_0$, and $\Psi(x_0, p_0) \not \in \Psi((\overline{{X}} \times \Sigma_0)\setminus \{(x_0, p_0)\})$.
\end{enumerate}
Then, letting $\Phi = \Psi \circ (\tau \times \id_{\Sigma_0})$, we have
\[
\Phi_{\#}(\current{I^m} \times \current{\Sigma_0}) = \pm\current{M}.
\]
\end{prop}
\begin{proof}
To reduce notation, we let 
\[
T = \Phi_{\#}(\current{I^m} \times \current{\Sigma_0}).
\]
By assumption (i), we can apply Proposition~\ref{prop:homology-class-from-Lipschitz} to see that $\partial T = 0$. Since $\supp(T) \subset M$, the constancy theorem~\cite[Theorem 26.27]{Simon1983} then yields some integer $l \in \ZZ$ such that
\begin{equation}\label{eq:constancy}
T = l\current{M}.
\end{equation}
Next, by condition (ii), there are neighborhoods $U$ of $(x_0, p_0)$ in ${X} \times \Sigma_0$ and $V$ of $\Psi(x_0, p_0)$ in $M$ such that $\Psi|_{U}$ is a $C^1$-diffeomorphism from $U$ onto $V$. Defining
\[
W = V \setminus \Psi((\overline{{X}} \times \Sigma_0) \setminus U),\quad\quad \widetilde{U} = (\tau \times \id)^{-1}(U),
\]
we infer by the second part of (ii) that $W$ is a neighborhood of $\Psi(x_0, p_0)$ in $M$. Note also that $W \subset \Phi(\widetilde{U})$, that $\Phi^{-1}(W) \subset \widetilde{U} \subset \Int(I^{m}) \times \Sigma_0$, and that $\Phi$ is injective on $\widetilde{U}$.

To continue, define
\[
\widetilde{U}_{+} = \{(t, p) \in \widetilde{U}\ |\ d\Phi_{(t, p)} \text{ exists and has full rank}\},
\]
and note by the Lipschitz continuity of $\Phi$ and the area formula that
\begin{equation}\label{eq:bad-set-has-null-image}
\cH^{m+k}(\Phi(\widetilde{U} \setminus \widetilde{U}_{+})) = 0.
\end{equation}
Then, for $y \in W \cap \Phi(\widetilde{U}_{+})$, we let
\[
\theta(y) = \bangle{ \vol_M|_{y}, \sum_{(t, p) \in \Phi^{-1}(\{y\}) \cap \widetilde{U}_{+}} \frac{d\Phi_{(t, p)}(e_1\wedge \cdots \wedge e_m \wedge v_1 \wedge \cdots \wedge v_k)}{|d\Phi_{(t, p)}(e_1\wedge \cdots \wedge e_m \wedge v_1 \wedge \cdots \wedge v_k)|}},
\]
where $\bangle{\cdot, \cdot}$ denotes the pairing between $(m+k)$-vectors and $(m+k)$-covectors, and $v_1, \cdots, v_k$ is any oriented orthonormal basis of $T_{p}\Sigma_0$. 
Since $d\Phi_{(t, p)}:\RR^{m} \times T_{p}\Sigma_0 \to T_{y}M$ is a linear isomorphism for each $(t, p)$ as in the summation, and since $\Phi$ is injective on $\widetilde{U}$, we see that 
\[
\theta(y) \in \{-1, 1\},\quad\text{for all }y \in W \cap \Phi(\widetilde{U}_{+}).
\]
Moreover, the area formula implies that $\theta$ is $\cH^{m+k}$-measurable, and also relates it to $T$ as follows: given $\zeta \in C^{\infty}_{c}(W)$, letting $\alpha$ be an arbitrary smooth extension of $\zeta\vol_{M}$ to a compactly supported $(m + k)$-form on $\RR^{q}$, we have 
\[
\begin{split}
T(\alpha) =\ &\int_{\Phi^{-1}(W)} \bangle{\alpha \circ \Phi, d\Phi(e_1\wedge \cdots \wedge e_m \wedge v_1 \wedge \cdots \wedge v_k)}\,d\cH^{m + k}\\
=\ & \int_{\Phi^{-1}(W) \cap \widetilde{U}_{+}} \bangle{\alpha \circ \Phi, d\Phi(e_1\wedge \cdots \wedge e_m \wedge v_1 \wedge \cdots \wedge v_k)}\,d\cH^{m + k}\\
=\ & \int_{W \cap \Phi(\widetilde{U}_{+})} \zeta(y)\theta(y)\, d\cH^{m + k}(y) = \int_{W} \zeta(y)\theta(y)\, d\cH^{m + k}(y),
\end{split}
\]
where the last equality uses~\eqref{eq:bad-set-has-null-image}. Since $\current{M}(\alpha) = \int_{W}\zeta\, d\cH^{m + k}$, we deduce from~\eqref{eq:constancy} and the arbitrariness of $\zeta \in C^{\infty}_{c}(W)$ that $l = \pm 1$. The proof is complete.
\end{proof}

We end this section by specializing to the case $k = 2$, so that $\Sigma_0$ is a closed oriented surface, which we write as $S$. The following is a corollary of the previous two propositions along with Proposition~\ref{prop:classes-of-maps-to-currents}. Recall that, given $\bv_{0}$ in the collection $\cP(m, S)$ defined in Section~\ref{subsec:statements}, the width of the equivalence class $[\bv_{0}]$ is defined by~\eqref{eq:width-defi}.

\begin{coro}\label{coro:non-trivial-sweepout-in-P}
Suppose $m + 2 \leq n$. Let $N$ be a closed, oriented, connected $(m+2)$-manifold embedded in some Euclidean space, and $h: N \to M$ a smooth map satisfying
\[
\big[h_{\#}\current{N}\big] \neq 0 \text{ in }H_{m+2}(M, \ZZ).
\]
Suppose also that the domain ${X} \subset \RR^{m}$, along with the maps $\tau: (I^{m}, \partial I^{m}) \to (\overline{X}, \partial X)$ and  $\tau^{-1}:(\overline{X}, \partial X) \to (I^{m}, \partial I^{m})$, are as in Proposition~\ref{prop:induced-sweepout}, and that $\Psi:\overline{{X}} \times S \to N$ is a Lipschitz map with the following properties:
\vskip 1mm
\begin{enumerate}
\item[(i)] $x \mapsto \Psi(x, \cdot) $ defines a continuous mapping from ${X}$ into $C^{1}(S; N)$.
\vskip 1mm
\item[(ii)] $\Psi$ is a $C^{1}$-diffeomorphism near some $(x_0, p_0) \in {X} \times S$, and moreover 
\[
\Psi(x_0, p_0) \not\in \Psi((\overline{{X}} \times S) \setminus \{(x_0, p_0)\}).
\]
\vskip 1mm
\item[(iii)] For all $x \in \partial {X}$, there holds $\cH^{2}(\Psi(\{x\} \times S)) = 0$, and that 
\[
A(\Psi(x', \cdot)) \to 0\quad\text{as }x' \to x \text{ from within }{X},
\]
where we use the induced metric on $N$ to compute the mapping area.
\end{enumerate} 
\vskip 1mm
Then, there exists $\bw \in \cP(m, S)$ such that
\[
\cW([\bw]) \geq \mu_{m, M},
\]
where $\mu_{m, M}$ is given by Proposition~\ref{prop:classes-of-maps-to-currents}.
\end{coro}
\begin{proof}
Define
\[
\Phi = h \circ \Psi \circ (\tau \times \id): I^m \times S \to M.
\]
Then by assumption (ii) and the first part of (iii), we may invoke Proposition~\ref{prop:induced-sweepout} to get
\begin{equation}\label{eq:class-of-Phi}
\Phi_{\#}(\current{I^m}\times \current{S}) = \pm h_{\#}\current{N}.
\end{equation}
Next, since $h$ is smooth and thus Lipschitz, we have by (iii) that
\[
\cH^{2}(\Phi(\{t\} \times S)) = \cH^{2}((h \circ \Psi)(\{\tau(t)\} \times S)) = 0\quad\text{for all }t \in \partial I^m.
\]
Proposition~\ref{prop:homology-class-from-Lipschitz} is then applicable to $\Phi$, which together with~\eqref{eq:class-of-Phi} gives
\begin{equation}\label{eq:class-of-Psi-nontrivial}
F_A([\widetilde{f}_{\Phi}]) = \pm\big[h_{\#}\current{N}\big]  \neq 0 \quad \text{in }H_{m+2}(M, \ZZ).
\end{equation}
On the other hand, letting 
\[
\bw(t) := \Phi(t, \cdot) = (h\circ \Psi)(\tau(t), \cdot) \quad\text{for }t \in I^m,
\]
and noting that $\tau(\Int(I^{m})) \subset X$, we see from the continuity of $\Phi$ and assumption (i) that $\bw$ and  $\bw|_{\Int(I^m)}$ are continuous as maps into $C^0(S; M)$ and $(C^0 \cap W^{1, 2})(S; M)$, respectively. Also, the second part of (iii) implies that, for all $t_0 \in \partial I^m$ and sequence $(t_i)$ in $\Int(I^{m})$ converging to $t_0$, we have
\[
A(\bw(t_i)) \to 0 \quad\text{as }i \to \infty,
\]
which gives~\eqref{eq:no-area-at-boundary} since $\partial I^m$ is compact. Thus $\bw$ belongs to $\cP(m, S)$. Comparing~\eqref{eq:f-Phi-defi} with~\eqref{eq:fw-defi}, and using Proposition~\ref{prop:homology-class-from-Lipschitz}(a), we have
\[
\widetilde{f}_{\Phi}(t) = f_{\bw}(t)\quad \text{for all }t \in I^m.
\]
Recalling that $F_A$ is an isomorphism, we obtain from~\eqref{eq:class-of-Psi-nontrivial} that
\[
[f_{\bw}] \neq 0\quad \text{in }\pi_m(\cZ_2(M), \{0\}).
\]
Thus, given $\bv \in [\bw]$, we have by Proposition~\ref{prop:classes-of-maps-to-currents}(b)(c) that 
\[
\sup_{t \in \Int(I^m)}A(\bv(t)) \geq \mu_{m, M}.
\]
This proves the asserted lower bound on the width of $[\bw]$.
\end{proof}
\subsection{Sweepouts from multisections}\label{subsec:multisections}
A \emph{multisection} of a closed manifold is a decomposition in which each collection of pieces has a specific type of intersection. By the seminal work of Gay--Kirby~\cite{Gay-Kirby2016}, where the study of multisections was initiated, as well as subsequent developments due to Lambert-Cole--Miller~\cite{Lambert-Cole-Miller19} and Ben Aribi--Courte--Golla--Moussard~\cite{BACGM2023}, it is known that every closed, connected, orientable smooth manifold of dimension $4$ or $5$ admits a multisection. The purpose of this section is to show that, just as Heegard splittings yield sweepouts of closed $3$-manifolds by surfaces, we can likewise use multisections to obtain non-trivial sweepouts in $\cP(m, S)$ for some $m \in \NN$ and closed surface $S$. For convenience, in this section, as well as the related sections of Appendix~\ref{appendix:patching-adapted}, we take the dimension of our closed Riemannian manifold $M$ to be $n + 1$, contrary to the convention used elsewhere in this paper. That said, we continue to assume that $M$ is isometrically embedded in some Euclidean space.

The following definition of multisections, or more precisely $n$-sections, is taken from~\cite{BACGM2023}. As preparation, to each integer $1 \leq k \leq n$ we associate a $(k-1)$-dimensional simplex $\Delta^{k-1}_{0} \subset \RR^{n}$ with barycenter at the origin, whose vertices are labeled $v_{1; k}, \cdots, v_{k; k}$. By $E^{k-1}$ we mean the $(k-1)$-plane containing $\Delta^{k-1}_{0}$, and for each $i \in \{1, \cdots, k\}$, we write
\[
E^{k-1}_{i} := \big\{\sum_{j \in \{1, \cdots, k\}\setminus \{i\}} \lambda_{j}v_{j; k}\ |\ \text{ each }\lambda_{j} \geq 0\big\}.
\]
Also, following the usage in~\cite{BACGM2023}, by a \textit{$1$-handlebody}, we mean specifically an orientable smooth manifold which has dimension at least $3$ and is obtained by attaching a finite number of $1$-handles to a $0$-handle. The number of $1$-handles is called the \textit{genus} of the $1$-handlebody, which in the present situation determines its diffeomorphism type. (See~\cite[Corollary VI.11.4]{Kosinski1993}, or~\cite[Example 4.1.4(b)]{Gompf-Stipsicz99}.)

\begin{defi}[\cite{BACGM2023}, pages 2 to 3]\label{defi:multisections}
Suppose $\dim M = n+ 1 \geq 3$. A decomposition of $M = \cup_{i = 1}^{n}M_{i}$ into a collection of $n$ compact subsets with pairwise disjoint interiors is said to be an \textit{$n$-section of genus $g$} if it satisfies the following additional requirements.
\vskip 1mm
\begin{enumerate}
\item[(1)] $M_{1} \cap \cdots \cap M_{n}$ is an embedded, closed, orientable surface of genus $g$, called the \textit{central surface},
\vskip 1mm
\item[(2)] Given $p \in M$, letting $I = I_{p} = \{i \in \{1, \cdots, n\}\ |\ p \in M_{i}\}$ and writing $k$ for its length $|I|$, there exist
\vskip 1mm
\begin{itemize}
\item a neighborhood $U \subset M$ of $p$,
\vskip 1mm
\item a bijection $\sigma:I \to \{1, \cdots, k\}$,
\vskip 1mm
\item a diffeomorphism $\varphi$ from a neighborhood $\widetilde{U} \subset \RR^{n + 2 - k} \times E^{k-1}$ of $(0, 0)$ onto $U$,
\end{itemize}
\vskip 1mm
with the property that
\begin{equation}\label{eq:adapted-chart-property}
\varphi(\widetilde{U} \cap (\RR^{n+ 2 - k} \times E^{k-1}_{\sigma(i)})) = U \cap M_{i}\quad\text{for all }i \in I.
\end{equation}
We refer to the map $\varphi$ as an \textit{adapted chart} at $p$. Note that there is no loss of generality in assuming further that 
\begin{equation}\label{eq:adapted-chart-disjoint}
\overline{U} \cap M_{i} = \emptyset\quad\text{whenever }i \not\in I.
\end{equation}
\vskip 1mm
\item[(3)] Given a non-empty subset $I \subset \{1, \cdots, n\}$, define 
\[
M_{I} = \cap_{i \in I}M_{i}, \quad  \mathring{M}_{I} = M_{I} \setminus \big( \cup_{i \not\in I}M_{i} \big), \quad \partial M_{I} = M_{I} \setminus \mathring{M}_{I}.
\]
Noting that $p \in \mathring{M}_{I}$ if and only if $\{1 \leq i \leq n\ |\ p\in M_{i}\} = I$, we deduce from~\eqref{eq:adapted-chart-property} and~\eqref{eq:adapted-chart-disjoint} that, in the case $|I| < n$, the set $\mathring{M}_{I} \cup \big( \cup_{i \not\in I}\mathring{M}_{I \cup \{i\}} \big)$ is an embedded $(n+2-k)$-submanifold with boundary, and $\mathring{M}_{I}$ is the set of interior points. With the above observations in mind, whenever $1 \leq |I| < n$, we further demand that $\mathring{M}_{I}$ be diffeomorphic to the interior of a $1$-handlebody.
\end{enumerate}
\end{defi}

As already mentioned, the existence of $n$-sections in low dimensions was established in the fundamental works~\cite{Gay-Kirby2016} and~\cite{BACGM2023}. We combine the relevant results into the following statement.
\begin{thm}[{\cite{Gay-Kirby2016, BACGM2023}}]\label{thm:existence of multisection}
    Suppose $n=3$ or $4$. Every closed, connected, orientable smooth $(n+1)$-dimensional manifold admits an $n$-section of some genus $g>1$. 
\end{thm}
\begin{proof}
    The existence of an $n$-section is contained in~\cite[Theorem 4]{Gay-Kirby2016} for $n = 3$, and in \cite[Theorem 7.3]{BACGM2023} for $n = 4$. The requirement that $g > 1$ can be fulfilled by performing stabilizations if necessary, as in~\cite[Lemma 10]{Gay-Kirby2016} and~\cite[Proposition 5.7]{BACGM2023}.\\
\end{proof}

Continuing to take $n =3$ or $4$, and denoting the central surface by $S_0$, in what follows we explain how, from the existence of an $n$-section, one can construct a domain $X$ and a Lipschitz map $\Psi:\overline{X} \times S_0 \to M$ to which the results of the previous section, particularly Corollary~\ref{coro:non-trivial-sweepout-in-P}, can be applied. Technical details omitted from the outline below can often be found in Appendix~\ref{appendix:patching-adapted}, and we include references to the latter where appropriate.

For convenience, we work with a different collection of simplices than those mentioned prior to Definition~\ref{defi:multisections}. Let $\be_{1}, \cdots, \be_{n}$ denote the standard basis of $\RR^{n}$. Given a non-empty index set $I \subset \{1, \cdots, n\}$, and writing $k = |I|$ as above, we define
\[
\bc_{I} = \frac{1}{k} \sum_{i \in I}\be_{i}, \quad\quad \ba_{i; I} = \be_{i} - \bc_{I} \quad\text{for }i \in I.
\]
Then $\{\ba_{i; I}\}_{i \in I}$ is the vertex set of a $(k-1)$-simplex with barycenter at the origin, denoted $\Delta^{k-1, I}$, which is contained in the $(k-1)$-plane given by
\[
V^{k-1, I} = \Span\{\ba_{i; I}\ |\ i \in I\}.
\]
With $B_{r} \subset \RR^{n}$ denoting open balls centered at the origin, we write 
\[
B_{r}^{k-1, I} = B_{r} \cap V^{k-1, I}.
\]
Also, for each $i \in I$, we set
\[
V^{k-1, I}_{i} = \Big\{ \sum_{j \in I \setminus \{i\}} c_{j}\ba_{j; I}  \ \big| \ \text{each }c_{j} \geq 0  \Big\}.
\]
When $I$ is equal to $\{1, \cdots, n\}$, we drop it from the notation just introduced. Thus, for example, we have
\[
\ba_{i} = \ba_{i; \{1, \cdots, n\}}, \quad V^{n-1} = V^{n-1; \{1, \cdots, n\}},\quad B^{n-1} = B^{n-1; \{1, \cdots, n\}},
\]
and so forth. Among the convenient properties of these objects, we mention that given another non-empty index set $I' \subset I$, with $l: = |I'| < k$, there holds the following orthogonal decomposition (see~\eqref{eq:decomposition-wrt-lower-simplex} and~\eqref{eq:splitting-of-vertex}):
\[
V^{k-1, I} = V^{l-1, I'} \oplus \Span\{\ba_{i; I}\ |\ i \in I \setminus I'\}.
\]

Given any bijection $\sigma:I \to \{1, \cdots, k\}$, since the simplex $\Delta^{k-1}_{0}$ chosen before Definition~\ref{defi:multisections} is also $(k-1)$-dimensional and centered at the origin, there is a linear isomorphism from $V^{k-1, I}$ to $E^{k-1}$ mapping $\ba_{i; I}$ to $v_{\sigma(i); k}$ for each $i\in I$. Thus, for all $p \in \mathring{M}_{I}$, it follows from Definition~\ref{defi:multisections}(2) that there exists a chart 
\[
\varphi: \widetilde{U}\subset \RR^{n + 2 -k} \times {V}^{k-1, I} \to U \subset M
\]
such that $p \in U$, that $(0, 0) \in \widetilde{U}$, and that
\begin{equation}\label{eq:modified-adapted-chart}
\left\{
\begin{array}{ll}
U \cap M_{i} = \varphi(\widetilde{U} \cap (\RR^{n + 2 - k} \times {V}^{k-1, I}_{i})),& \text{ if } i \in I,\\
\overline{U} \cap M_{i} = \emptyset, & \text{ if }i \not\in I.
\end{array}
\right.
\end{equation}
The fact that these charts respect the $n$-section in the sense of~\eqref{eq:modified-adapted-chart} allows us to patch them together along $\mathring{M}_{I}$. The result is the following, which we prove in Appendix~\ref{subsec:fixed-stratum}.

\begin{prop}\label{prop:fixed-stratum-patching}
Let $\Omega$ be any relatively open subset with compact closure in $\mathring{M}_{I}$, and $\cC$ a closed subset of $M$ not intersecting $\mathring{M}_{I}$. Then there exists a diffeomorphism $G$ from a neighborhood $\cW_1$ of $\Omega$ in $M$ onto $\Omega \times B^{k-1, I}_{\rho}$ for some $\rho > 0$, such that 
\[
\overline{\cW_{1}} \cap \cC = \emptyset = \overline{\cW_{1}} \cap M_{i}\quad\text{for all }i \not\in I,
\]
that $G(q) = (q, 0)$ for all $q \in \Omega$, and that
\[
G(\cW_1 \cap M_{i}) = \Omega \times (B_{\rho} \cap V^{k-1, I}_{i}) \quad \text{for all }i \in I.
\]
\end{prop}
\begin{proof}
We prove this in Appendix~\ref{subsec:fixed-stratum}, where the result is restated as Proposition~\ref{prop:global-multisection-chart}. The proof uses some of the notation introduced in Appendix~\ref{subsec:notation-for-patching}.\\
\end{proof}

We present the remainder of the construction in the form of a proof sketch, distilled from the content of Appendices~\ref{subsec:different-strata-I} through~\ref{subsec:different-strata-III}, for the following statement. 

\begin{prop}\label{prop:multisection-to-sweepout}
Let $n = 3$ or $4$, and assume that $M^{n+1}$ admits an $n$-section of some genus $g\geq 0$, with the central surface denoted $S_0$. Then there exist a bounded domain $X \subset V^{n-1}$ and a Lipschitz continuous map 
\[
\Psi: \overline{X}\times S_0 \to M
\]
such that $B^{n-1}_{r} \subset X$ for some $r > 0$, that $\overline{X}$ is bi-Lipschitz equivalent to $I^{n-1}$, and that $\Psi$ satisfies
\begin{enumerate}
\item $x \mapsto \Psi(x, \cdot)$ is a continuous map from $X$ into $C^{1}(S_0; M)$.
\item $\Psi|_{B^{n-1}_{r} \times S_0}$ is a smooth diffeomorphism onto an open set in $M$. Moreover,
\[
\Psi(B^{n-1}_{r} \times S_0) \cap \Psi((\overline{X} \setminus B^{n-1}_{r}) \times S_0) = \emptyset.
\]
\item For all $x \in \partial X$, we have $\cH^{1}(\Psi(\{x\} \times S_0)) < \infty$. Moreover, if $(x_i)$ is a sequence in $X$ converging to $x$, then
\[
\lim_{i \to \infty}A(\Psi(x_{i}, \cdot)) = 0.
\]
\end{enumerate}
\end{prop}
\begin{proof}[{\bf Sketch of proof}]
We only address the case $n=4$. Given $m \in \{3, 4, 5\}$, by the \textit{$m$-dimensional stratum} of the given quadrisection of $M$, we mean the union $\bigcup_{|I| = 6 - m}M_{I}$. Also, we refer to the desired map $\Psi$ loosely as an \emph{$S_0$-sweepout} of $M$. We shall build $\Psi$ in several stages. To start, we invoke Proposition~\ref{prop:fixed-stratum-patching}, with $I = \{1, 2, 3, 4\}$ and $\Omega = S_0$, to obtain a neighborhood $\cU_0$ of $S_0$ in $M$ and a diffeomorphism 
\[
f_0: \cU_0 \to B^{3}_{1} \times S_0
\]
satisfying $f_0(p) = (0, p)$ for all $p \in S_0$, and that
\[
f_0(\cU_0 \cap M_{i}) = (B_{1} \cap V^{3}_{i}) \times S_0 \quad\text{for all }i \in \{1, 2, 3, 4\}.
\]
Thus $f_0^{-1}$ can be regarded as providing an $S_0$-sweepout of a neighborhood of $S_0$ itself in $M$. The remainder of the proof consists of three steps, where we successively extend this sweepout to neighborhoods of the $3$ and $4$-dimensional strata, and finally to the $5$-dimensional stratum.

\vspace{0.5em}
\noindent\textbf{Step 1: Construction near the $3$-dimensional stratum.}
\vskip 2mm

This step corresponds largely to Appendix~\ref{subsec:different-strata-I}, and is divided further into two parts.
\vskip 2mm
\noindent\textit{(1) Parametrization of collar, and utilization of handlebody structure.} 
\vskip 2mm

Given $i \in \{1, 2, 3, 4\}$, let $\{i\}^{c} = \{1, 2, 3, 4\} \setminus \{i\}$. Note that $M_{\{i\}^{c}}$ is a smooth $3$-submanifold with boundary equal to $S_0$, and that a collar of the latter is parametrized by the map
\[
(t, p) \mapsto f_0^{-1}(t \frac{\ba_{i}}{|\ba_{i}|}, p),\quad\text{for }(t, p) \in [0, 7/8] \times S_0.
\]
Recalling also that $M_{\{i\}^{c}}$ is a $1$-handlebody, and using the building blocks described in Appendix~\ref{subsec:retraction}, particularly Lemma~\ref{lemm:Omega-S-basics} and Proposition~\ref{prop:bh-properties}, we may extend the above to a Lipschitz map 
\[
h_{\{i\}^{c}}:[0, 2] \times S_0 \to M_{\{i\}^{c}}
\]
such that (see below~\eqref{eq:1234-parametrize-collar})
\vskip 1mm
\begin{itemize}
\item $t \mapsto h_{\{i\}^{c}}(t, \cdot)$ is a continuous map from $[0, 2)$ into $C^1(S_0; M)$,
\vskip 1mm
\item $\cH^{1}(h_{\{i\}^{c}}(\{2\} \times S_0)) < \infty$, and $\lim_{t \to 2^{-}}A(h_{\{i\}^{c}}(t, \cdot)) = 0$,
\vskip 1mm
\item $h_{\{i\}^{c}}([0, 2] \times S_0) = M_{\{i\}^{c}}$, and $h_{\{i\}^{c}}(\{t\} \times S_0) \cap h_{\{i\}^{c}}(\{t'\} \times S_0) = \emptyset$ whenever $t \neq t'$.
\end{itemize}
\vskip 1mm
Again applying Proposition~\ref{prop:fixed-stratum-patching}, this time with $I = \{i\}^{c}$, and with $\Omega$ taken to be the set
\[
M_{\{i\}^{c}}^*: = M_{\{i\}^{c}} \setminus h_{\{i\}^{c}}([0, 1/8] \times S_0) = h_{\{i\}^{c}}((1/8, 2] \times S_0), 
\]
we obtain a neighborhood $\cU_{\{i\}^{c}}$ of $M_{\{i\}^{c}}^*$, along with a diffeomorphism
\[
f_{\{i\}^{c}}: \cU_{\{i\}^{c}} \to B_{1}^{2, \{i\}^{c}} \times M_{\{i\}^{c}}^*,
\]
such that $f_{\{i\}^{c}}(p) = (0, p)$ for all $p \in M_{\{i\}^{c}}^*$, and that 
\[
\cU_{\{i\}^{c}} \cap M_{j} = f_{\{i\}^{c}}^{-1}((B_1 \cap V^{2, {\{i\}^{c}}}_{j}) \times M_{\{i\}^{c}}^{*})\quad \text{for all }j \neq i.
\]
Furthermore, we can also arrange so that
\[
\overline{\cU_{\{i\}^{c}}} \cap M_{i} = \emptyset, \quad\text{and}\quad 
\overline{\cU_{\{i\}^{c}}} \cap \overline{\cU_{\{j\}^{c}}} = \emptyset\quad\text{whenever }i\neq j.
\]
\vskip 2mm
\noindent\textit{(2) Patching with sweepout near $S_0$ }
\vskip 2mm
To continue, again take $i \in \{1, 2, 3, 4\}$, and abbreviate $\{i\}^{c}$ as $I$. Consider the assignments
\begin{equation}\label{eq:step-1-assignment-1}
(x, t, p) \mapsto f_0^{-1}(x + t\frac{\ba_{i}}{|\ba_{i}|}, p),\quad\text{for }(x, t, p ) \in  B^{2, {I}}_{1/8} \times (0, 7/8) \times S_0,
\end{equation}
and 
\begin{equation}\label{eq:step-1-assignment-2}
(x, t, p) \mapsto f_{{I}}^{-1}(x, h_{{I}}(t, p)),\quad\text{for }(x, t, p) \in B^{2, {I}}_{1} \times (1/8, 2] \times S_0.
\end{equation}
We denote~\eqref{eq:step-1-assignment-2} by $f^{-1}_{{I}} \circ \phi_{{I}}$ (see~\eqref{eq:collar-spread-out}), and write~\eqref{eq:step-1-assignment-1} still as $f_0^{-1}$ by slight abuse of notation. The two maps agree with each other on $\{0\} \times (1/8, 7/8) \times S_0$ (see~\eqref{eq:agreement-on-central-leaf}). Moreover, both are diffeomorphisms when restricted to $ B^{2, {I}}_{1/8} \times (1/8, 7/8) \times S_0$. Thus, fixing an interval $A_0$ satisfying
\[
(1/3, 2/3) \Subset A_0 \Subset (1/8, 7/8), 
\]
we obtain, for sufficiently small $r_0  \in (0, 1/8)$, a well-defined injective local diffeomorphism
\[
\theta_{{I}}: B_{r_0}^{2, {I}} \times A_0 \times S_0 \longrightarrow B^{2, {I}}_{1/8} \times (1/8, 7/8) \times S_0,
\]
characterized by the relation
\[
f_{I}^{-1}\circ \phi_{I} \big( \theta_{{I}}(x, t, p) \big) = f_{0}^{-1} (x, t, p) \quad\text{(see below~\eqref{eq:transition-inclusion-2})}.
\]
It can be shown that $\theta_{{I}}$ restricts to the identity on $\{0\} \times A_0 \times S_0$, and that it respects the decomposition
\[
V^{2, {I}} = \bigcup_{j \in {I}}V^{2, {I}}_{j}
\]
in the $x$-variable (see Claim~\ref{claim:theta-properties}). Thanks primarily to these observations, upon using cut-off functions to interpolate between $\theta_{I}$ and the identity (see~\eqref{eq:F-I-definition}), and composing the result with $f_{I}^{-1} \circ \phi_I$, we obtain the map denoted $f_{I}^{-1}\circ \phi_{I} \circ F_{I}$ appearing on the left-hand side of~\eqref{eq:F-endpoints}, which respects the quadrisection structure (see Claim~\ref{claim:F-I-multisection} for the actual statement), reduces to~\eqref{eq:step-1-assignment-1} when $t \leq 3/8$ and to~\eqref{eq:step-1-assignment-2} when $t \geq 5/8$ (see~\eqref{eq:F-endpoints}), and restricts to a diffeomorphism on $B_{r_1}^{2, I} \times (1/3, 2/3) \times S_0$ for some $r_1 < r_0$ (Claim~\ref{claim:F-diffeo}). In view of the orthogonal decomposition~\eqref{eq:decomposition-wrt-lower-simplex}, which implies that 
\begin{equation}\label{eq:decomposition-for-step-1}
V^{3} = V^{2, \{i\}^{c}} \oplus \Span\{\ba_{i}\},
\end{equation}
upon patching together $(f_{I}^{-1}\circ \phi_{I} \circ F_{I})|_{B_{r_1}^{2,I} \times (1/3, 2/3) \times S_0}$ and $(f_{I}^{-1} \circ \phi_{I})|_{B_{r_1}^{2, I} \times (5/8, 2] \times S_0}$, and inverting the correspondence $(x, t) \mapsto x + t\frac{\ba_{i}}{|\ba_{i}|}$, we obtain a map defined on the product with $S_0$ of the tube
\begin{equation}\label{eq:cylinders-in-outline}
C_{r_1}^{I}((1/3, 2]): = \big\{x + t\frac{\ba_{i}}{|\ba_{i}|}\ \big|\ x \in B^{2, I}_{r_{1}},\ t \in (1/3, 2] \big\}.
\end{equation}

Fixing $\rho \in (1/3, 3/8)$, letting $i$ run through $\{1, 2, 3, 4\}$, and choosing $\sigma \in (0, \rho/3)$ subject to the requirements~\eqref{eq:sigma-threshold-1} and~\eqref{eq:sigma-threshold-2}, we may join each of the maps on $C^{\{i\}^{c}}_{\sigma}((1/3,2]) \times S_0$ obtained at the end of the previous paragraph with $f_0^{-1}|_{B^{3}_{\rho} \times S_0}$ to obtain the map 
\[
G:\cN \times S_0 \to M
\]
which is the subject of Proposition~\ref{prop:1234-with-123}. Here, given $R \in [3/8, 2]$, we define (see~\eqref{eq:cN-definition})
\[
\cN_{R} = B_{\rho}^{3} \cup \bigcup_{i = 1}^{4}C_{\sigma}^{\{i\}^{c}}((0, R])
,\quad \mathring{\cN}_{R} =  B_{\rho}^{3} \cup \bigcup_{i = 1}^{4}C_{\sigma}^{\{i\}^{c}}((0, R)),
\]
dropping $R$ from the notation when $R = 2$. The actual definition of $G$ is given in~\eqref{eq:G-in-core} and~\eqref{eq:G-on-tube}, while its important properties are listed in Proposition~\ref{prop:1234-with-123}, and verified in Claim~\ref{claim:G-properties}. We shall not repeat them here.

\vskip 2mm
\noindent\textbf{Step 2: Construction near the 4-dimensional stratum}
\vskip 2mm
This step corresponds to Appendix~\ref{subsec:different-strata-II}. We break the description into three parts.
\vskip 2mm
\noindent\textit{(1) Rounding out corners}
\vskip 2mm
Given $J \subset \{1,2, 3, 4\}$ with $|J| = 2$, let 
\[
\{i_0, i_1\}: = \{1, 2, 3, 4\}\setminus J.
\]
For $r \in (0, \frac{\rho}{3})$, we use Proposition~\ref{prop:1234-with-123}(d) to write $M_{J}$ as
\begin{equation}\label{eq:step-2-splitting}
M_{J} = (M_{J}\setminus  G(\overline{B_{r}} \times S_0)) \cup G((V^{3}_{J} \cap B_{3r}) \times S_0).
\end{equation}
The first set on the right-hand side is a smooth $4$-submanifold of $M \setminus G(\overline{B_{r}} \times S_0)$ with boundary given by $\partial M_{J}\setminus  G(\overline{B_{r}} \times S_0)$ (see the proof of Lemma~\ref{lemm:smoothed-strata}), while the second is a neighborhood in $M_{J}$ of $S_0$, the latter consisting of corner points of $\partial M_{J}$, so to speak. Notice that, in the present case $V^{3}_{J}$ has the form
\[
V^{3}_{J} = \{c_1\ba_{i_0} + c_{2}\ba_{i_1}\ |\ c_1, c_2 \geq 0\},
\]
whose relative boundary in $\Span\{\ba_{i_0},\ba_{i_1}\}$, denoted  $\partial V^{3}_{J}$, is the graph of a piecewise linear function. With $\alpha > 0$ to be chosen, by a standard procedure detailed in Appendix~\ref{subsec:smoothing}, which involves mollifying the graphing function and flowing by a suitable normal vector field (see Proposition~\ref{prop:smoothing-flow} and Remark~\ref{rmk:Psi-distance-bounds} for the end product), we obtain a universal constant $C$, a smooth domain $\Omega \subset V^{3}_{J}$, and a smooth vector field $\xi:\partial\Omega \to \Span\{\ba_{i_0},\ba_{i_1}\}$, where $\partial\Omega$ denotes boundary relative to $\Span\{\ba_{i_0},\ba_{i_1}\}$, such that, with $\rho_1: = C\alpha$, we have
\[
\Omega \setminus B^{3}_{\rho_1}  = V^{3}_{J} \setminus B^{3}_{\rho_1},\quad \partial\Omega \setminus B_{\rho_1}^3  = \partial V^{3}_{J} \setminus B_{\rho_1}^3,
\]
and that $\xi$ coincides on $\partial V^{3}_{J} \setminus B^3_{\rho_1}$ with the unit normal pointing into $V^{3}_{J}$. Moreover, the map
\[
\Psi: (s, z) \mapsto z + s\xi(z)
\]
is a diffeomorphism from $(-\alpha, \alpha) \times \partial\Omega$ onto an open set in $\Span\{\ba_{i_0} ,\ba_{i_1}\}$, and satisfies 
\[
\Psi(s, z) \in \Omega\quad\text{if and only if }s \geq 0.
\]
We then fix $\alpha$ small enough to have $\rho_1 < \frac{1}{8}\min\{\rho, \sigma\}$ (\eqref{eq:alpha-choice} and \eqref{eq:rho1-choice}), and define
\[
\widetilde{M}_{J} = \big( M_{J} \setminus G(\overline{B_{\rho_1}} \times S_0) \big) \cup G((\Omega \cap B_{3\rho_1}) \times S_0) \quad\text{(see~\eqref{eq:4-strata-smoothing})}.
\]
It is then shown in~\eqref{eq:difference-from-smoothing} and Lemma~\ref{lemm:smoothed-strata} that $\widetilde{M}_{J}$ is a compact, smooth $4$-submanifold of $M$ with boundary given by replacing $M_{J}$ and $\Omega$ in the right-hand side above with $\partial M_{J}$ and $\partial\Omega$, respectively. Moreover, by Remark~\ref{rmk:portions} we have
\[
\partial \widetilde{M}_{J} = G((\partial\Omega \cap \cN) \times S_0) \quad\text{(see~\eqref{eq:smoothing-contained-in-G-image})},
\]
so we may loosely think of $\partial\widetilde{M}_{J}$ as being swept out by copies of $S_0$.
\vskip 2mm
\noindent\textit{(2) Parametrization of collar}
\vskip 2mm
Our next task is to construct an analogue of $f_0^{-1}$, or, more precisely, of the map~\eqref{eq:step-1-assignment-1}. That is, we seek an injective local diffeomorphism
\[
g_{J}: B_{\alpha}^{1, J} \times (-\alpha, \alpha) \times \partial\widetilde{M}_{J} \to M
\]
whose restriction to $\{0\} \times [0, \alpha) \times \partial\widetilde{M}_{J}$ parametrizes a collar of $\partial \widetilde{M}_{J}$ in $\widetilde{M}_{J}$. For that purpose, we split $\partial \widetilde{M}_{J}$ into three pieces as
\begin{equation}\label{eq:partial-M-J-splitting-step-2}
\begin{split}
\partial\widetilde{M}_{J} =\ & G((\mathring{\cN}_{2/3} \cap \partial\Omega ) \times S_0) \\
& \cup h_{J \cup\{i_{0}\}}((5/8, 2] \times S_0) \big) \cup h_{J \cup\{i_{1}\}}((5/8, 2] \times S_0) \big) 
\end{split}
\end{equation}
(see~\eqref{eq:boundary-expression-for-g-J}). Then, in~\eqref{eq:g-J-definition-case-1} and~\eqref{eq:g-J-definition-case-2}, according to the above decomposition, we define $g_{J}$ by patching three maps together, and prove in Proposition~\ref{prop:g-J-properties} that it is indeed well-defined and has the properties described above. In particular, Proposition~\ref{prop:g-J-properties}(a) gives the following convenient formula:
\begin{equation}\label{eq:g-J-formula-step-2}
g_{J}(x, s, G(y, p)) = G(x + \Psi(s, y), p),
\end{equation}
for all $(x, s) \in B^{1, J}_{\alpha} \times (-\alpha, \alpha)$, and $(y, p) \in (\partial \Omega \cap \cN) \times S_0$. A key fact used in the proof of Proposition~\ref{prop:g-J-properties} is that, as a consequence of the properties of $\Psi$ and the orthogonal decomposition~\eqref{eq:decomposition-wrt-lower-simplex}, or more relevantly~\eqref{eq:decomposition-simplex-special-case}, the map
\begin{equation}\label{eq:Psi-hat-in-outline}
\widehat{\Psi}:(x, s, y) \mapsto x + \Psi(s, y)
\end{equation}
is a diffeomorphism from $B^{1, J}_{\alpha} \times (-\alpha, \alpha) \times \partial\Omega$ onto an open set in $V^{3}$. This subsequently plays the role of the correspondence $(x, t) \mapsto x + t\frac{\ba_{i}}{|\ba_{i}|}$ leading to the map~\eqref{eq:step-1-assignment-1}.

Fixing parameters $0 < \mu_0 < \mu_1 < \cdots < \mu_5 < 1$ whose choice we explain at the start of (3) below, and noting that 
\[
(s, q) \mapsto g_{J}(0, s, q) \quad\text{for }(s, q)\in [0, \mu_{4}\alpha] \times \partial\widetilde{M}_{J}
\]
parametrizes a collar of $\partial\widetilde{M}_{J}$ in $\widetilde{M}_{J}$, we obtain, using the $1$-handlebody structure of the latter and the results from Appendix~\ref{subsec:retraction}, a Lipschitz map 
\[
h_{J}:[0, \mu_5\alpha] \times \partial\widetilde{M}_{J} \to \widetilde{M}_{J}
\]
which extends the above map, and has properties analogous to those enjoyed by $h_{\{i\}^{c}}$ from step 1, part (1), except that in the second property, the statement regarding the mapping area should be replaced by the condition~\eqref{eq:h-J-shrink-area} on the action of $dh_{J}$ on $2$-vectors tangent to $\partial\widetilde{M}_{J}$. We then apply Proposition~\ref{prop:fixed-stratum-patching} with $I$ taken to be $J$, and $\Omega$ taken to be
\[
M_{J}^* : = \widetilde{M}_{J} \setminus h_{J}([0, \frac{\mu_0\alpha}{2}] \times \partial\widetilde{M}_{J}),
\]
to obtain a neighborhood $\cU_{J}$ of $M_{J}^{*}$ and a diffeomorphism 
\[
f_{J}: \cU_{J} \to B^{1, J}_{1} \times M_{J}^{*},
\]
satisfying properties analogous to the map $f_{\{i\}^{c}}$ from (1) in step 1.

\vskip 2mm
\noindent\textit{(3) Patching with sweepout near $3$-dimensional stratum}
\vskip 2mm

We first explain how the constants $\mu_0$ through $\mu_5$ mentioned above are found. Given $r > 0 $, we let $\cA_{r}$ denote the distance neighborhood
\[
\cA_{r} : = \{y \in V^{3}\ |\ |y - x| < r, \ \text{for some }x \in \cup_{i = 1}^{4}V^{3}_{\{i\}^{c}}\},
\]
which coincides with a union of tubes, namely 
\[
\cA_{r} = \cup_{i = 1}^{4}C^{\{i\}^{c}}_{r}([0, \infty)) \quad\text{(see~\eqref{eq:distance-neighborhood-1-strata})}.
\]
By the distance estimates in Remark~\ref{rmk:Psi-distance-bounds}, combined with~\eqref{eq:distance-attained}, we can find positive constants $\mu_0$, $\mu_1$, $\tau$, and $\delta$, independent of the index set $J$, such that $\mu_0 < \tau - 2\delta < \tau + 2\delta < \mu_1 < 1$, and that
\[
\begin{split}
\Omega \cap \cA_{\mu_0\alpha} \subset \ &\Psi([0, \mu_0\alpha) \times \partial\Omega) \subset \Omega \cap \cA_{(\tau - 2\delta)\alpha},\\
\Omega \cap \cA_{(\tau + 2\delta)\alpha} \subset \ &\Psi([0, \mu_1\alpha) \times \partial\Omega) \subset \Omega \cap \cA_{\alpha}.
\end{split}
\]
(See equations~\eqref{eq:mu1-mu0-inequalities} through~\eqref{eq:nested-smoothing-2}.) We then choose $\mu_2 < \cdots < \mu_{5}$ from $(\mu_{1}, 1)$ arbitrarily. Returning to the maps found in (2) above, and recalling that we are viewing $g_{J}$ as an analogue of~\eqref{eq:step-1-assignment-1}, we introduce also the following analogue of~\eqref{eq:step-1-assignment-2}:
\begin{equation}\label{eq:step-2-assignment-2}
(f_{J}^{-1}\circ\phi_{J})(x, s, q) := f_{J}^{-1}(x, h_{J}(s, q)),
\end{equation}
for $(x, s, q) \in B_{\alpha}^{1, J} \times (\frac{\mu_0\alpha}{2}, \mu_5\alpha] \times \partial\widetilde{M}_{J}$. In a way completely parallel to (2) in step 1, we obtain a sufficiently small $r_0 < \alpha$ so that a well-defined transition map
\[
\theta_{J}: B_{r_0}^{1, J} \times (\mu_0\alpha, \mu_{3}\alpha) \times \partial\widetilde{M}_{J} \to B^{1, J}_{\alpha} \times (\frac{\mu_0\alpha}{2}, \mu_{4}\alpha) \times \partial\widetilde{M}_{J}
\]
results from imposing the relation 
\[
f_{J}^{-1} \circ \phi_{J} \circ \theta_{J} = g_{J}.
\]
Choosing an interval $A_{1}(\alpha)$ satisfying
\[
(\mu_{1}\alpha, \mu_{2}\alpha) \Subset A_{1}(\alpha)\Subset (\mu_{0}\alpha, \mu_{3}\alpha),
\]
we can again interpolate between $\theta_{J}$ and the identity (see~\eqref{eq:F-J-for-4-strata}). Doing so leads to the map $f_{J}^{-1} \circ \phi_{J} \circ F_{J}$ on the left-hand side in~\eqref{eq:4-strata-F-endpoints}, which transitions from $g_{J}$ to the map~\eqref{eq:step-2-assignment-2}, and restricts to a diffeomorphism on $B^{1, J}_{r_2} \times A_{1}(\alpha) \times \partial\widetilde{M}_{J}$ for some $r_2 < r_0$. It is then possible to join $ (f_{J}^{-1}\circ\phi_{J} \circ F_{J}^{-1})|_{B^{1, J}_{r_2} \times A_{1}(\alpha) \times \partial\widetilde{M}_{J}}$ with $g_{J}|_{B^{1, J}_{r_{2}} \times (0, \mu_{1}\alpha) \times \partial\widetilde{M}_{J}}$ and $ (f_{J}^{-1}\circ \phi_{J})|_{B^{1, J} \times (\mu_{2}\alpha, \mu_{5}\alpha] \times \partial\widetilde{M}_{J}}$ to obtain the map
\[
H_{J}: B^{1, J}_{r_2} \times (0, \mu_{5}\alpha] \times \partial\widetilde{M}_{J} \to M \quad\text{(see~\eqref{eq:H-J-definition})},
\]
which is the subject of Proposition~\ref{prop:H-J-properties}. Notice that, by the relation~\eqref{eq:smoothing-contained-in-G-image} already mentioned above, it makes sense to consider
\begin{equation}\label{eq:H-J-G-in-outline}
(x, s, z, p) \mapsto H_{J}(x, s, G(z, p)), 
\end{equation}
for $(x, s, z) \in B^{1, J}_{r_2} \times (0, \mu_{5}\alpha] \times (\partial\Omega \cap \cN)$, and $p \in S_0$.

We are now at a point where we need to introduce domains analogous to~\eqref{eq:cylinders-in-outline}. Specifically, given $0 < r \leq \alpha$ and $E \subset [0, \alpha)$, and marking the $J$-dependence of $\Psi$, $\Omega$, and $\widehat{\Psi}$ by subscripts, we let
\[
\begin{split}
C^{J}_{r}(E) :=\ & \widehat{\Psi}_{J}(B^{1, J}_{r} \times E \times (\partial\Omega_{J} \cap \cN)) \\
=\ & \{x + \Psi_{J}(s, z)\ |\ x \in B^{1, J}_{r},\ s \in E,\ z \in \partial\Omega_{J} \cap \cN \}.
\end{split}
\]
Fixing $h > 0$ satisfying~\eqref{eq:N'-height-threshold}, our choice of $\mu_0$ and $\mu_1$ guarantees, among other things, that the sets $\{C^{J}_{h}((\mu_0\alpha, \mu_5\alpha])\}_{|J| = 2}$ are mutually disjoint (see~\eqref{eq:C-J-positive-distance}). Moreover, for each $J \subset\{1, 2, 3, 4\}$ with $|J| = 2$, thanks primarily to the formula~\eqref{eq:g-J-formula-step-2}, it turns out that composing~\eqref{eq:H-J-G-in-outline} with $\widehat{\Psi}^{-1} \times \id_{S_0}$ yields a map on $C^{J}_{h}((\mu_0\alpha, \mu_{5}\alpha]) \times S_0$ that agrees with $G_{(1)}|_{(\cN \cap \cA_{\tau\alpha}) \times S_0}$ on the overlap of their domains. Thus, letting
\[
\cN' := (\cN \cap \cA_{\tau\alpha}) \cup \bigcup_{|J| = 2}C^{J}_{h}((\mu_0\alpha, \mu_5\alpha]) \quad\text{(see~\eqref{eq:N'-definition})},
\]
we define, via the formulas~\eqref{eq:G-1-definition-case-1} and~\eqref{eq:G-1-definition-case-2}, the map $G_{(1)}: \cN' \times S_0 \to M$ whose relevant properties are established in Proposition~\ref{prop:G-1-properties}.

\vskip 2mm

\noindent\textbf{Step 3: Construction in the $5$-dimension stratum}

\vskip 2mm
This step corresponds to Appendix~\ref{subsec:different-strata-III}, and is more or less parallel to step 2, but technically more involved in several places. We again give our description in three parts.

\vskip 2mm
\noindent\textit{(1) Rounding out corners}
\vskip 2mm
Fix $j_0 \in \{1, 2, 3, 4\}$, and write $\{j_0\}^{c}$ as $\{i_0, i_1, i_2\}$. Similar to~\eqref{eq:step-2-splitting}, we use Proposition~\ref{prop:1234-with-123}(d) to write, for sufficiently small $r> 0$, 
\[
M_{j_0} = \big(M_{j_0} \setminus G((\cN \cap \overline{\cA_{r}}) \times S_0)\big) \cup G( (\cN \cap \cA_{3r} \cap V^{3}_{j_0} ) \times S_0 ).
\]
Since $\cup_{i= 1}^{4} M_{\{i\}^{c}} \subset G((\cN \cap \overline{\cA_{r}}) \times S_0)$, which can be deduced from Proposition~\ref{prop:1234-with-123}(c), we see that in the right-hand side above, the first set is a smooth domain in $M \setminus G((\cN \cap \overline{\cA_{r}}) \times S_0)$, and we want to modify the second set, as in (1) of Step 2. This time we note that
\[
V^{3}_{j_0} =\ \{ c_0\ba_{i_0} + c_1\ba_{i_1} + c_2\ba_{i_2}\ |\ \text{each }c_{\lambda} \geq 0\},
\]
so that $\partial V^{3}_{j_0}$ is the graph of a piecewise linear function over a copy of $\RR^2$. With parameters $\rho_2$ and $\alpha'$ chosen to satisfy~\eqref{eq:rho2-choice} and~\eqref{eq:alpha'-threshold}, we apply Proposition~\ref{prop:smoothing-flow} to get a smooth domain $\Omega \subset V^{3}_{j_0}$ along with a smooth vector field $\xi: \partial\Omega \to V^{3}$ with the property that
\[
\Omega \setminus \cA_{\rho_2} = V^{3}_{j_0} \setminus \cA_{\rho_2}, \quad \partial\Omega \setminus \cA_{\rho_2} = \partial V^{3}_{j_0} \setminus \cA_{\rho_2},
\]
that $\xi$ coincides on $\partial V^{3}_{j_0} \setminus \cA_{\rho_2}$ with the unit normal pointing into $V^{3}_{j_0}$, and that
\[
\Psi: (s, y) \mapsto y + s\xi(y)
\]
restricts to a diffeomorphism from $(-\alpha', \alpha') \times \partial\Omega$ onto an open set in $V^{3}$, with $\Psi(s, y) \in \Omega$ if and only if $s \geq 0$. We then let
\[
\widetilde{M}_{j_0} = \big(M_{j_0} \setminus G((\cN \cap \overline{\cA_{\rho_2}}) \times S_0)\big) \cup G( (\cN \cap \cA_{3\rho_2} \cap \Omega ) \times S_0 ).
\]

At this point, the construction becomes more involved than that in step 2, primarily since $G$ is not a diffeomorphism on all of $\cN \times S_0$, but is only guaranteed to be so in $\mathring{\cN}_{7/8} \times S_0$ (Proposition~\ref{prop:1234-with-123}(a)). The way forward is to further write
\begin{equation}\label{eq:step-3-decomp-1}
\begin{split}
G( (\cN \cap \cA_{3\rho_2} \cap \Omega) \times S_0 ) =\ & G( (\mathring{\cN}_{2/3} \cap \cA_{3\rho_2} \cap \Omega) \times S_0 )\\
&\cup G( ((\cN \setminus \cN_{5/8}) \cap \cA_{3\rho_2} \cap \Omega) \times S_0 ),
\end{split}
\end{equation}
and take into account the product structure provided by Proposition~\ref{prop:smoothing-product-structure}, which yields, for each $\lambda\in \{0, 1, 2\}$, a domain $\Omega_{\lambda}' \subset V^{2, \{i_{\lambda}\}^{c}}_{j_0}$ such that if we identify $V^{3}$ with $V^{2, \{i_{\lambda}\}^{c}} \times \RR$ by 
\[
(x, t) \mapsto x + t \frac{\ba_{i_\lambda}}{|\ba_{i_\lambda}|},
\]
then $\Omega \cap (V^{2, \{i_{\lambda}\}^{c}} \times [\rho_2, \infty)) = \Omega_{\lambda}' \times [\rho_{2}, \infty)$ (see~\eqref{eq:product-structure-specialized}). This in turn implies
\begin{equation}\label{eq:step-3-decomp-2}
G( ((\cN \setminus \cN_{5/8}) \cap \cA_{3\rho_2} \cap \Omega) \times S_0 ) = \cup_{\lambda \in \{0, 1, 2\}} f_{\{i_{\lambda}\}^{c}}^{-1}((B^{2, \{i_{\lambda}\}^{c}}_{3\rho_2} \cap \Omega_{\lambda}') \times M_{\{i_{\lambda}\}^{c}}^{**}),
\end{equation}
where $M_{\{i_{\lambda}\}^{c}}^{**}$ stands for $h_{\{i_{\lambda}\}^{c}}((5/8, 2] \times S_0)$ (see~\eqref{eq:product-structure-for-image-Omega}). Largely with the help of~\eqref{eq:step-3-decomp-1} and~\eqref{eq:step-3-decomp-2}, we prove in Proposition~\ref{prop:5-strata-smoothing}, and the remarks before it, that $\widetilde{M}_{j_0}$ as defined above is a compact, smooth domain in $M$ which is contained in $M_{j_0}$.

\vskip 2mm
\noindent\textit{(2) Parametrization of collar}
\vskip 2mm
To continue, take the constants $\mu_0$ through $\mu_{5}$, as well as $\tau$ and $\delta$, from the start of step 2, part (3). Then $\partial\widetilde{M}_{j_0}$ turns out to admit the following expression analogous to~\eqref{eq:partial-M-J-splitting-step-2}:
\[
\begin{split}
\partial\widetilde{M}_{j_0} = \ & 
G(( \mathring{\cN}_{2/3} \cap \cA_{(\tau - \delta)\alpha} \cap \partial\Omega) \times S_0)\nonumber\\
& \cup \big( \cup_{\lambda \in \{0, 1, 2\}} f_{\{i_\lambda\}^{c}}^{-1}( (B_{(\tau - \delta)\alpha} \times \partial\Omega_{\lambda}') \times M_{\{i_{\lambda}\}^{c}}^{**} )\big)\nonumber\\
& \cup \big(\cup_{\lambda \in \{0, 1, 2\}} h_{\{j_0, i_{\lambda}\}}\big((\mu_0 \alpha, \mu_4\alpha) \times \partial\widetilde{M}_{\{j_0, i_{\lambda}\}}\big) \big)\nonumber\\
& \cup \big( \cup_{\lambda \in \{0, 1,2\}} M_{\{j_0, i_{\lambda}\}}^{**} \big)\quad\text{(see~\eqref{eq:four-part-decomp})},
\end{split}
\]
where $M_{\{j_0, i_{\lambda}\}}^{**}$ stands for $h_{\{j_0, i_{\lambda}\}}((\mu_{2}\alpha, \mu_{5}\alpha] \times \partial \widetilde{M}_{\{j_0, i_{\lambda}\}})$, with $h_{\{j_0, i_{\lambda}\}}$ and $\widetilde{M}_{\{j_0, i_{\lambda}\}}$ being taken from step 2. Accordingly, in Definition~\ref{defi:g-j-definition}, we introduce the maps $g_{j_0}^{(1)}$ through $g_{j_0}^{(4)}$, one for the product with $(-\alpha', \alpha')$ of each part on the right-hand side above, and prove in Lemma~\ref{lemm:g-j-target}, Proposition~\ref{prop:g-j-patching}, and Proposition~\ref{prop:g-j-injective}, that they patch together to give an injective local diffeomorphism
\[
g_{j_0}:(-\alpha', \alpha') \times \partial\widetilde{M}_{j_0} \to M,
\]
with the property that $g_{j_0}(s, q) \in \widetilde{M}_{j_0}$ if and only if $s \geq 0$. Also, we mention that it now takes two formulas to cover the role played by~\eqref{eq:g-J-formula-step-2} in step 2. These are~\eqref{eq:g-j-expression-in-G} and~\eqref{eq:g-j-expression-in-H}, which we do not reproduce here.

Next we define, for $s > 0$, the distance neighborhood
\[
\cA^{(2)}_{s} = \{y \in V^{3}\ |\ |y - x| < s,\ \text{for some }x \in \cup_{|J| = 2}V^{3}_{J}\}.
\]
Then, similar to how $\mu_0$, $\mu_1$, $\tau$, and $\delta$ were chosen, we can find universal constants $\nu_{1}$, $\nu_{2}$, $\tau'$, and $\eta$, such that $0 < \nu_{1} < \tau' - 2\eta < \tau' + 2\eta < \nu_{2} < 1$, and that 
\[
\begin{split}
&\Omega_{j_0}  \cap \cA_{\nu_{1}\alpha'}^{(2)} \subset  \Psi_{j_0}([0, \nu_1\alpha') \times \partial\Omega ) \subset \Omega_{j_0} \cap \cA_{(\tau' - 2\eta)\alpha'}^{(2)},\\
&\Omega_{j_0} \cap \cA_{(\tau' + 2\eta)\alpha'}^{(2)} \subset \Psi_{j_0}([0, \nu_2\alpha') \times \partial\Omega_{j_0}) \subset \Omega_{j_0}  \cap \cA_{\alpha'}^{(2)},
\end{split}
\]
where we have added subscripts to indicate the $j_0$-dependence of $\Omega$ and $\Psi$ (see~\eqref{eq:2-strata-neighborhood-squeezed-1} through~\eqref{eq:smoothing-agree-2-strata-nbhd}). Choose in addition $\nu_{3}, \nu_{4} \in (\nu_{2}, 1)$ such that $\nu_{3}< \nu_{4}$, and note that
\[
(s, q) \mapsto g_{j_0}(s, q), \quad\text{for }(s, q) \in [0, \nu_{3}\alpha'] \times \partial\widetilde{M}_{j_0}
\]
parametrizes a collar of $\partial\widetilde{M}_{j_0}$ in $\widetilde{M}_{j_0}$. The latter being a $1$-handlebody, we again obtain a Lipschitz extension 
\[
h_{j_0}:[0, \nu_{4}\alpha'] \times \partial\widetilde{M}_{j_0} \to \widetilde{M}_{j_0}
\]
of the above map, satisfying the properties listed below~\eqref{eq:collar-5-strata}.
\vskip 2mm
\noindent\textit{(3) Patching with sweepout near $4$-dimensional stratum}
\vskip 2mm
Recall the domain $\cN'$ and the map $G_{(1)}:\cN' \times S_0 \to M$ from step 2, part (3). For each $j_0 \in \{1, 2, 3, 4\}$, with the help of the expression for $\partial\Omega_{j_0} \cap \cN'$ derived in Lemma~\ref{lemm:Psi-action-on-N'}, and the definition of $G_{(1)}$, we find in Proposition~\ref{prop:M-j-as-G-1-image} that $\partial\widetilde{M}_{j_0}$ can be expressed as
\[
\partial\widetilde{M}_{j_0} = G_{(1)} ((\partial\Omega_{j_0} \cap \cN') \times S_0).
\]
In particular, it makes sense to consider the map 
\[
(s, z, p) \mapsto h_{j_0}(s, G_{(1)}(z, p))\quad\text{for }(s, z, p)\in (\nu_1\alpha', \nu_{4}\alpha'] \times (\partial\Omega_{j_0} \cap \cN') \times S_0,
\]
which together with the inverse of $(s, z) \mapsto \Psi_{j_0}(s, z)$ yields a map defined on $\Psi_{j_0}((\nu_{1}\alpha', \nu_{4}\alpha'] \times (\partial\Omega_{j_0} \cap \cN')) \times S_0$. Letting
\[
\cN'' = (\cN' \cap \cA^{(2)}_{\tau'\alpha'}) \cup \big(\cup_{j_0 = 1}^{4} \Psi_{j_0}((\nu_{1}\alpha', \nu_{4}\alpha'] \times (\partial\Omega_{j_0} \cap \cN'))\big),
\]
we then join these previous maps with the restriction of $G_{(1)}$ to $\cN' \cap \cA^{(2)}_{\tau'\alpha'}$ to obtain a map 
\[
G_{(2)}:\cN'' \times S_0 \to M,
\]
in the following fashion (see~\eqref{eq:G-2-definition-case-1} and~\eqref{eq:G-2-definition-case-2}): if $y \in \cN' \cap \cA^{(2)}_{\tau'\alpha'}$, we let 
\[
G_{(2)}(y, p) = G_{(1)}(y, p),
\]
while if $y  = \Psi_{j_0}(s, z)$ for some $j_0\in \{1, 2, 3, 4\}$ and $(s, z) \in (\nu_1\alpha', \nu_{4}\alpha'] \times (\partial\Omega_{j_0} \cap \cN')$, we let 
\[
G_{(2)}(y, p) = h_{j_0}(s, G_{(1)}(z, p)).
\]
Proposition~\ref{prop:G-2-properties} establishes the properties of $G_{(2)}$, including the fact that it is well-defined, while in Lemma~\ref{lemm:cN''-compact} we show that $\cN''$ is compact, and characterize its interior, which coincides with the set denoted $\mathring{\cN}''$ in~\eqref{eq:N''-definition}. We complete the proof of Proposition~\ref{prop:multisection-to-sweepout} by taking $X = \mathring{\cN}''$, $\Psi = G_{(2)}$, and $r = \frac{\nu_{1}\alpha'}{2}$.\\
\end{proof} 

Combining all the above results, we can now prove Theorem~\ref{thm:existence-of-sweepout}, whose statement we reproduce below.

\begin{thm*}
Let $M$ be a closed, connected, oriented Riemannian manifold, isometrically embedded in some Euclidean space. Suppose for some $m \in \{2, 3\}$ that the homology group $H_{m+2}(M, \ZZ)$ is nontrivial. Then there exists a closed, oriented surface $S$ of genus $g > 1$, and an element $\bv_{0} \in \cP(m, S)$, such that 
\[
\cW([\bv_0]) > 0.
\]
\end{thm*}
\begin{proof}
Let $\alpha$ be a nontrivial class in $H_{m+2}(M,\ZZ)$. Since $m + 2 \leq 5$, by a famous result of Thom \cite[page 55, Theorem II.27]{Thom54} (see also the English translation \cite[page 172, Theorem II.27]{Novikov-Taimanov07}), the class $\alpha$ is realizable by a submanifold. That is (see~\cite[page 28]{Thom54} or~\cite[page 142]{Novikov-Taimanov07}), there exists a closed, oriented $(m+2)$-submanifold $N$ such that, with $h: N \to M$ denoting the inclusion map, there holds
\[
[h_{\#}\current{N}] = \alpha.
\]
Without loss of generality, we assume that $N$ is also connected. Applying Theorem \ref{thm:existence of multisection} and Proposition \ref{prop:multisection-to-sweepout} to $N$ yields a closed surface $S$ of genus $g > 1$, a domain $X \subset \RR^{m}$, and a Lipschitz map $\Psi: \overline{X} \times S \to N$, satisfying the requirements of Corollary \ref{coro:non-trivial-sweepout-in-P}. Consequently there exists $\bv_0 \in \cP(m, S)$ such that $\cW([\bv_0]) \geq \mu_{m, M} > 0$. This finishes the proof.
\end{proof}

\section{Facts on the Teichm\"uller space}\label{sec:Teichmuller}
For the proof of Theorem \ref{thm:min-max-existence} we employ several models for the Teichm\"uller space of closed surfaces with genus $g$. As background, we discuss in Section \ref{subsec:quasi-conformal} quasiconformal maps and their relationship to conformal classes of metrics, along with a number of standard \textit{a priori} estimates. In Section \ref{subsec:DE-extension} we review the Douady--Earle extension mentioned in Section \ref{subsec:statements}, which among other things produces quasiconformal maps of $\HH$ out of certain homeomorphisms of $\partial \HH$ in a way that respects the action of $\Aut{\HH}$. In Sections \ref{subsec:Fuchsian} and \ref{subsec:hyperbolic}, we very briefly review some basic facts about Fuchsian models and hyperbolic metrics. Finally, in Section \ref{subsec:Teichmuller} we recall the models of Teichm\"uller space that we need.

\subsection{Quasiconformal maps}\label{subsec:quasi-conformal}
There are a number of equivalent ways to define quasiconformal maps. For a discussion of these, see for instance~\cite[pages 78, 81 and 86]{Imayoshi-Taniguchi}, and~\cite[Section 4.5]{Hubbard}, and~\cite[Chapters I and IV]{Lehto-Virtanen1973}. We adopt here an analytical definition.
\begin{defi}\label{defi:quasiconformal-maps}
Let $\Omega$ be a domain in $\CC$ and suppose $K \geq 1$. A \emph{$K$-quasiconformal map} on $\Omega$ is an orientation-preserving homeomorphism $f$ from $\Omega$ onto an open set in $\CC$ such that the distributional derivatives $f_z, f_{\overline{z}}$ lie in $L^2_{\loc}(\Omega)$, and that
\begin{equation}\label{eq:qc-diff-ineq}
\big|f_{\overline{z}}\big| \leq \frac{K - 1}{K + 1}\big|f_{z}\big|,\quad \text{almost everywhere on }\Omega.
\end{equation}
By a \emph{quasiconformal map} we simply mean a map that is $K$-quasiconformal for some $K \geq 1$. 
\end{defi}

Classical conformal maps are $1$-quasiconformal. Conversely, any $1$-quasiconformal map, being distributionally holomorphic, is a conformal map in the classical sense by elliptic regularity. (See for instance~\cite[page 82, Corollary 2]{Imayoshi-Taniguchi}.) To motivate the next definition, we recall the following standard property.

\begin{lemm}\label{lemm:non-vanishing-Jacobian}
Let $f$ be a $K$-quasiconformal map on a domain $\Omega \subset \CC$. Then $f_{z} \neq 0$ almost everywhere.
\end{lemm}
\begin{proof}
This is~\cite[Theorem 5(ii)]{Ahlfors-Bers}. See also Appendix~\ref{appendix:quasiconformal}.
\end{proof}
\begin{defi}\label{defi:complex-dilatation}
Given a $K$-quasiconformal map $f$, by the previous lemma, letting
\[
\mu_f := f_{\overline{z}}/f_z
\]
defines a bounded measurable function on $\Omega$, called the \emph{complex dilatation} of $f$, which satisfies $\|\mu_f\|_{\infty} \leq \frac{K - 1}{K + 1} < 1$.
\end{defi}

Given a domain $\Omega \subset \CC$ and a complex-valued function $\mu \in L^{\infty}(\Omega)$ such that $\|\mu\|_{\infty; \Omega} < 1$, the problem of finding quasiconformal maps with complex dilatation $\mu$ leads to the famous Beltrami equation, which takes the form
\begin{equation}\label{eq:Beltrami-eq}
f_{\overline{z}} = \mu f_z,
\end{equation}
and is closely related to the construction of isothermal coordinates. Specifically, given a Riemannian metric $g$ on $\Omega$, recall that there is a unique way to put it into the form
\begin{equation}\label{eq:metric-complex}
g = \lambda |dz + \mu d\overline{z}|^2,
\end{equation}
where $\lambda$ is a positive function, while $\mu$ takes complex values and satisfies $|\mu| < 1$. More explicitly, if $g = E dx^2 + 2F dx dy + G dy^2$, then $\mu$ and $\lambda$ are given by (see for example~\cite[page 20]{Imayoshi-Taniguchi}):
\begin{equation}\label{eq:lambda-mu}
\lambda = \frac{1}{4}(E + G + 2\sqrt{EG - F^2}),\quad \mu = \frac{E - G + 2i F}{E + G + 2\sqrt{EG - F^2}}.
\end{equation}
In particular, $\mu$ depends only on the conformal class of $g$. Also, the eigenvalues of $g$ can be expressed in terms of $\lambda$ and $\mu$ as
\begin{equation}\label{eq:eigenvalue-lambda-mu}
\sigma_1 = \lambda(1 - |\mu|)^2,\quad \sigma_2 = \lambda(1 + |\mu|)^2.
\end{equation}
Now suppose $f: \Omega \to \CC$ is a quasiconformal map which solves~\eqref{eq:Beltrami-eq} with $\mu$ being the one from~\eqref{eq:metric-complex}. Then, at least on a formal level, we see that the pullback $f^*(|dz|^2)$ is a metric conformal to $g$ where $f_{z} \neq 0$. In other words, the inverse of $f$ can be viewed as providing an isothermal coordinate chart.

To continue, we recall a few additional properties of quasiconformal maps, before stating the fundamental existence and uniqueness theorems, due to Morrey~\cite{Morrey1938} and Ahlfors--Bers~\cite{Ahlfors-Bers}, for the Beltrami equation~\eqref{eq:Beltrami-eq}. 

\begin{lemm}[\cite{Ahlfors-Bers}, Theorem 5]
\label{lemm:qc-standard-properties}
Let $f$ be a quasiconformal map on a domain $\Omega$, and write $k: = \|\mu_f\|_{\infty; \Omega}$, so that $k \in [0, 1)$.
\vskip 1mm
\begin{enumerate}
\item[(a)] There exists $p = p(k)> 2$ such that $\nabla f \in L_{\loc}^{p}(\Omega)$.
\vskip 1mm
\item[(b)] Both $f$ and $f^{-1}$ take sets of (two-dimensional Lebesgue) measure zero to sets of measure zero.
\vskip 1mm
\item[(c)] $f^{-1}$ is quasiconformal on $f(\Omega)$, with complex dilatation given by
\begin{equation}\label{eq:inverse-Beltrami-coefficient}
\mu_{f^{-1}} = -\Big( \mu_f \cdot \frac{f_{z}}{\overline{f}_{\overline{z}}} \Big) \circ f^{-1}.
\end{equation}
In fact, the distributional derivatives of $f^{-1}$ satisfy
\begin{equation}\label{eq:f-inverse-derivatives}
(f^{-1})_{z} =\Big( \frac{1}{(1 - |\mu_{f}|^2)f_{z}}\Big) \circ f^{-1},\quad (f^{-1})_{\overline{z}} =\Big( \frac{-\mu_f}{(1 - |\mu_{f}|^2)\overline{f}_{\overline{z}}}\Big) \circ f^{-1}.
\end{equation}
\vskip 1mm
\item[(d)] Let $g$ be another quasiconformal map on $\Omega$. Then $g \circ f^{-1}:f(\Omega) \to g(\Omega)$ is quasiconformal. Also, 
\begin{equation}\label{eq:composition-Beltrami-coefficient}
\mu_{g \circ f^{-1}} = \Big(\frac{\mu_g - \mu_f}{1 - \overline{\mu_f}\mu_g} \cdot \frac{f_{z}}{\overline{f}_{\overline{z}}}\Big) \circ f^{-1}.
\end{equation}
\vskip 1mm
\item[(e)] If $\mu$ is $C^1$ on $\Omega$, then $f$ is an orientation-preserving $C^1$-diffeomorphism from $\Omega$ onto $f(\Omega)$.
\end{enumerate}
\end{lemm}
\begin{proof}
Each statement either is contained in~\cite[Theorem 5]{Ahlfors-Bers}, or follows easily from the proof thereof. We elaborate on some of the conclusions, as well as the choice of $p(k)$, in Appendix~\ref{appendix:quasiconformal}. See Remark~\ref{rmk:choice-of-p} for the latter point.
\end{proof}

\begin{rmk}\label{rmk:K-K'-quasiconformal}
From parts (c) and (d), we infer the well-known fact that the inverse of a $K$-quasiconformal map is still $K$-quasiconformal, while the composition of a $K$-quasiconformal map and a $K'$-quasiconformal map, whenever it makes sense, is $K'K$-quasiconformal.
\end{rmk}

\begin{thm}[\cite{Ahlfors-Bers}, Theorem 6; \cite{Imayoshi-Taniguchi}, Theorem 4.30]
\label{thm:qc-existence}
Given a complex-valued function $\mu \in L^{\infty}(\CC)$ with $k : = \|\mu\|_{\infty; \CC} < 1$, there exists a unique quasiconformal map $f^\mu$ from $\CC$ onto $\CC$ which has complex dilatation $\mu$ and satisfies
\begin{equation}\label{eq:qc-normalization}
f^{\mu}(0) = 0,\ \ f^{\mu}(1) = 1,\ \ \lim_{|z| \to \infty}|f^{\mu}(z)| = \infty.
\end{equation}
This unique solution has the property that $\nabla f^\mu \in L^p_{\loc}$, where $p = p(k) > 2$ is as in Lemma~\ref{lemm:qc-standard-properties}. Moreover, if $\mu$ is $C^1$ on some open set $\Omega \subset \CC$, then $f^{\mu}$ is an orientation-preserving $C^1$-diffeomorphism from $\Omega$ onto $f^{\mu}(\Omega)$.
\end{thm}
\begin{proof}
The existence and uniqueness statement is~\cite[Theorem 6]{Ahlfors-Bers}, where we note that the ``$\mu$-conformal homeomorphisms'' according to Ahlfors--Bers are quasiconformal maps in the sense of Definition~\ref{defi:quasiconformal-maps} (see Remark~\ref{rmk:mu-conformal-converse}). The integrability statement follows from Lemma~\ref{lemm:qc-standard-properties}(a). The final statement is Lemma~\ref{lemm:qc-standard-properties}(e).
\end{proof}

\begin{thm}[\cite{Ahlfors-Bers}, Lemma 18 and Theorem 9]
\label{thm:qc-dependence}
Given $k \in [0, 1)$, suppose $(\mu_n)$ is a sequence of complex-valued functions in $L^{\infty}(\CC)$ such that $\|\mu_n\|_{\infty; \CC} \leq k$ for all $n$, and that 
\[
\mu_n \to \mu \text{ almost everywhere on }\CC.
\]
Then, with $p = p(k)$ as in Theorem~\ref{thm:qc-existence}, we have for all $R > 0$ that
\[
\|\nabla f^{\mu_n} - \nabla f^\mu\|_{p; \bB_R} +  \|f^{\mu_n} - f^\mu\|_{\infty; \bB_R} \to 0 \text{ as } n \to \infty.
\]
\end{thm}
\begin{proof}
For use in later arguments, we reproduce a part of the proof in Appendix~\ref{subsec:beltrami-proofs}, taking for granted a crucial $C^0 \cap W^{1, p}$-estimate (Proposition \ref{prop:qc-estimates}) for solutions of~\eqref{eq:Beltrami-eq} produced by Theorem \ref{thm:qc-existence}.
\end{proof}

We also recall a version of Theorem~\ref{thm:qc-existence} that gives quasiconformal maps from the upper half-plane $\HH$ to itself. This is Theorem \ref{thm:qc-existence-H} below. The uniqueness statement therein again involves a certain three-point condition, which makes sense due to Proposition \ref{prop:qc-extension}.

\begin{prop}[\cite{Imayoshi-Taniguchi}, Proposition 4.31]
\label{prop:qc-extension}
Every quasiconformal map from $\bB$ onto itself extends to a homeomorphism of $\overline{\bB}$.
\end{prop}
\begin{proof}
Let $g$ be a quasiconformal map from $\bB$ onto itself. The proof of~\cite[Lemma 14]{Ahlfors-Bers} yields a quasiconformal map $\varphi$ from $\CC$ onto itself such that 
\[
\varphi(\partial \bB) = \partial\bB, \quad \varphi(\bB) = \bB,
\]
and that $\mu_{\varphi} = \mu_{g}$ on $\bB$. From Lemma~\ref{lemm:qc-standard-properties}(d), we infer that the composition $g \circ \varphi^{-1}|_{\bB}$ is a conformal map from $\bB$ to itself; that is, an element of $\Aut(\bB)$. In particular it extends to a homeomorphism of $\overline{\bB}$, and hence so does $g = (g \circ \varphi^{-1}) \circ \varphi$.
\end{proof}

Next consider the biholomorphic map $F:\HH \to \bB$ defined by
\begin{equation}\label{eq:H-to-D}
F(z) = \frac{z - i}{z + i}. 
\end{equation}
Given a quasiconformal map $w: \HH \to \HH$, by Lemma \ref{lemm:qc-standard-properties}, the composition 
\begin{equation}\label{eq:conjugation-by-F}
\theta_F(w):= F\circ w \circ F^{-1}: \bB \to \bB
\end{equation}
is still quasiconformal, and thus extends, by Proposition \ref{prop:qc-extension}, to a homeomorphism of $\overline{\bB}$, which we continue to write as $\theta_{F}(w)$. Observe that, for $\theta_F(w)|_{S^1}$ to fix the points $-1$, $-i$, and $1$, it is necessary and sufficient that
\begin{equation}\label{eq:qc-H-boundary-condition}
\lim_{z \to 0}w(z) = 0,\quad  \lim_{z \to 1}w(z) = 1,\quad \lim_{|z| \to \infty}|w(z)| = \infty.
\end{equation}
Following standard terminology, we say that \emph{$w$ fixes $0$, $1$, and $\infty$} whenever the above condition holds. If $w_1$ and $w_2$ are quasiconformal maps of $\HH$ that fix $0$, $1$, and $\infty$ in this sense, then so does $w_1 \circ w_2^{-1}$.

\begin{thm}[\cite{Hubbard}, pages 170 to 171; \cite{Imayoshi-Taniguchi}, Proposition 4.33]
\label{thm:qc-existence-H}
Suppose $\mu$ is a complex-valued function in $L^{\infty}(\HH)$ and that $\|\mu\|_{\infty; \HH} < 1$. Define \begin{equation}\label{eq:qc-reflection}
\widehat{\mu}(z) = \left\{
\begin{array}{ll}
\overline{\mu(\overline{z})}, & \text{ if }\im z  < 0,\\
\mu(z), & \text{ if }\im z > 0,
\end{array}
\right.
\end{equation}
and let $f^{\widehat{\mu}}$ be the map produced by Theorem~\ref{thm:qc-existence}. Then we have the following.
\vskip 1mm
\begin{enumerate}
\item[(a)] $f^{\widehat{\mu}}(\HH) = \HH$ and $f^{\widehat{\mu}}(\partial\HH) = \partial\HH$.
\vskip 1mm
\item[(b)] The restriction $w^{\mu}: = (f^{\widehat{\mu}})|_{\HH}$ is the unique quasiconformal map from $\HH$ to itself which has complex dilatation $\mu$ and fixes $0$, $1$, and $\infty$.
\vskip 1mm 
\item[(c)] If $\mu$ is $C^1$ on some open set $\Omega \subset \HH$, then $w^{\mu}$ is an orientation-preserving $C^1$-diffeomorphism from $\Omega$ to $w^{\mu}(\Omega)$.
\end{enumerate}
\end{thm}
\begin{proof}
With the help of the uniqueness statement in Theorem~\ref{thm:qc-existence}, we have 
\[
f^{\widehat{\mu}}(z) = \overline{f^{\widehat{\mu}}(\overline{z})} \quad\text{for all }z \in \CC,
\]
which together with the bijectivity of $f^{\widehat{\mu}}$ implies that $f^{\widehat{\mu}}(\CC \setminus \partial \HH) = \CC \setminus \partial \HH$. A connectedness argument then shows that $f^{\widehat{\mu}}(\HH)$ is either $\HH$ itself or the lower half-plane. To rule out the latter case, we mollify $\widehat{\mu}$ to obtain a sequence $(\lambda_n)$ of $C^1$-functions on $\CC$ such that $\|\lambda_n\|_{\infty} \leq \|\mu\|_{\infty}$ for all $n$, and that 
\[
\lambda_n \to \widehat{\mu} \quad\text{almost everywhere on }\CC.
\]
Upon replacing $\lambda_n(z)$ by $\frac{1}{2}(\lambda_n(z) + \overline{\lambda_n(\overline{z})})$ if needed, we may repeat the previous argument to see that $f^{\lambda_n}(\partial\HH) = \partial\HH$, and that $f^{\lambda_n}(\HH)$ can only be either $\HH$ or the lower half-plane. Now, since $f^{\lambda_n}$ has non-vanishing Jacobian by the last part of Theorem~\ref{thm:qc-existence}, we get from~\eqref{eq:qc-normalization} that it necessarily maps $\partial\HH$ onto itself in a strictly increasing fashion. The positivity of $Jf^{\lambda_n}$ then forces $f^{\lambda_n}(\HH) = \HH$. Passing to the limit using Theorem~\ref{thm:qc-dependence}, we find that, for instance, $f^{\widehat{\mu}}(i) \in \HH \cup \partial\HH$, and hence $f^{\widehat{\mu}}(\HH) = \HH$ by what we established at the start of the proof. 

For (b), since $f^{\widehat{\mu}}$ satisfies~\eqref{eq:qc-normalization}, it is clear that~\eqref{eq:qc-H-boundary-condition} holds for $w^{\mu}$. For uniqueness, we note that if $w:\HH \to \HH$ is a quasiconformal map satisfying the conditions in the statement, then $w^{\mu} \circ w^{-1}$ belongs to $\Aut(\HH)$ by Lemma~\ref{lemm:qc-standard-properties}(d), and moreover fixes $0$, $1$, and $\infty$, so consequently it must be the identity map.

Part (c) follows from the last assertion of Theorem~\ref{thm:qc-existence}.
\end{proof}
It is also important for us to have an analogue of Theorem \ref{thm:qc-dependence} that applies to the inverses of the maps produced by Theorem \ref{thm:qc-existence}. This is carried out in the work \cite{Zhou10} of the second named author. The relevant results are summarized below. We begin by recalling in Lemma \ref{lemm:conjugate-Beltrami} a differential equation for these inverses that is more convenient, for the purpose of deriving estimates, than the one coming from~\eqref{eq:inverse-Beltrami-coefficient}. The main continuous dependence result is Proposition \ref{prop:qc-H-dependence}.
\begin{lemm}[\cite{Zhou10}, Equation (63)]
\label{lemm:conjugate-Beltrami}
Given $\mu \in L^{\infty}(\CC)$ with $\|\mu\|_{\infty; \CC} < 1$, let $f:\Omega \to \CC$ be a quasiconformal map with complex dilatation $\mu$, and denote by $h: f(\Omega) \to \Omega$ its inverse. Then $h$ satisfies
\begin{equation}\label{eq:conjugate-Beltrami}
h_{\overline{z}} = -(\mu\circ h) \overline{h}_{\overline{z}}, \quad \text{almost everywhere on }f(\Omega).
\end{equation}
\end{lemm}
\begin{proof} 
By Lemma~\ref{lemm:qc-standard-properties}(c), we have
\[
h_{z}\circ f = \frac{\overline{f}_{\overline{z}}}{(1 - |\mu|^2)|f_z|^2},\ \ h_{\overline{z}} \circ f = -\frac{\mu f_z}{(1 - |\mu|^2)|f_z|^2},
\]
which implies that 
\[
h_{\overline{z}} \circ f = -\mu \cdot (\overline{h}_{\overline{z}}\circ f), \text{ almost everywhere on }\Omega.
\]
As quasiconformal maps take measure zero sets to measure zero sets, we obtain the desired conclusion upon composing both sides with $h$.
\end{proof}
\begin{prop}[\cite{Zhou10}, Lemma 6.3]
\label{prop:qc-H-dependence}
Suppose $(\mu_n)$ is a sequence of $C^1$-functions on $\HH$ that satisfies 
\[
\|\mu_n\|_{\infty; \HH} \leq k < 1\quad \text{for all }n,
\] 
and converges in $C^1_{\loc}(\HH)$ to some $\mu$. Let $w^{\mu_n}$ and $w^{\mu}$ be the quasiconformal maps produced by Theorem \ref{thm:qc-existence-H}, and denote their inverses by $v^{\mu_n}$ and $v^{\mu}$. Then, with $p = p(k)$ being the exponent in Theorem~\ref{thm:qc-existence}, the following conclusions hold.
\vskip 1mm
\begin{enumerate}
\item[(a)] We have
\begin{equation}\label{eq:qc-H-dependence}
\|v^{\mu_n} - v^{\mu}\|_{\infty; K} + \|\nabla v^{\mu_n} - \nabla v^{\mu} \|_{p; K} \to 0 \text{ as }n \to \infty,
\end{equation}
for all compact set $K \subset \HH$. 
\vskip 1mm
\item[(b)] In fact, $\nabla v^{\mu_n}$ converges to $\nabla v^{\mu}$ in $C^0_{\loc}(\HH)$. Moreover, we have
\begin{equation}\label{eq:uniform-W2p-on-compact}
\sup_{n}\|\nabla v^{\mu_n}\|_{1, p; K} < \infty,\quad \text{for all compact set } K \subset \HH.
\end{equation}
\vskip 1mm
\item[(c)] $w^{\mu_n}$ converges to $w^{\mu}$ in $C^1_{\loc}(\HH)$.
\end{enumerate}
\end{prop}
\begin{proof}
See Appendix~\ref{subsec:beltrami-proofs}. 
\end{proof}

By differentiating the Beltrami equation, we can upgrade the mode of convergence in Proposition~\ref{prop:qc-H-dependence} under stronger regularity assumptions. For example we have the following pair of results, whose proofs involve only standard techniques, and are gathered in Appendix~\ref{subsec:inverse-beltrami-proofs}. 
\begin{prop}\label{prop:qc-smooth-dependence}
Let $(w_n)$ and $w$ be smooth quasiconformal diffeomorphisms from $\HH$ onto itself that fix the points $0$, $1$, and $\infty$, and let $(\mu_n)$ and $\mu$ denote their respective complex dilatations. Assume that for some $k \in [0, 1)$ we have $\|\mu_n\|_{\infty; \HH} \leq k$ for all $n$, and that 
\[
\mu_n \to \mu \text{ smoothly locally on }\HH \text{ as }n \to \infty.
\]
Then $w_n \to w$ smoothly locally on $\HH$.
\end{prop}
\begin{proof}
See Appendix~\ref{subsec:inverse-beltrami-proofs}.
\end{proof}

\begin{prop}\label{prop:metric-inverse-convergence}
Under the assumptions of Proposition~\ref{prop:qc-smooth-dependence}, we have that $w_n^{-1}$ converges smoothly locally on $\HH$ to $w^{-1}$.
\end{prop}
\begin{proof}
See Appendix~\ref{subsec:inverse-beltrami-proofs}.
\end{proof}
Finally we explain how the notion of quasiconformality extends to mappings between closed Riemann surfaces, following \cite[Section 5.1.1]{Imayoshi-Taniguchi}. Suppose $\Sigma_0$ and $\Sigma$ are closed Riemann surfaces of genus $g$, and let $f$ be an orientation-preserving homeomorphism from $\Sigma_0$ to $\Sigma$. By a \emph{lift} of $f$ we mean a homeomorphism $\widetilde{f}: \HH \to \HH$ for which there are covering maps $p_0: \HH \to \Sigma_0$ and $p: \HH \to \Sigma$ such that 
\[
f \circ p_0 = p \circ \widetilde{f}.
\]
When we want to emphasize the choice of the covering map $p_0$, we refer to $\widetilde{f}$ as a lift of $f$ \emph{with respect to $p_0$}. 
\begin{equation}\label{eq:qc-lift}
\begin{tikzcd}
\mathbb{H} \arrow[d, "p_0"']  
\arrow[rr, "\widetilde{f}", dotted] &  & \mathbb{H} \arrow[d, "p"] \\
\Sigma_0 \arrow[rr, "f"']                                                                  &  & \Sigma                   
\end{tikzcd}
\end{equation}
Observe that the set of all the lifts of $f$ is given by $\{\alpha \circ \widetilde{f} \circ \beta\}_{\alpha, \beta \in \Aut(\HH)}$, while its lifts with respect to $p_0$ consists of $\{\alpha \circ \widetilde{f}\}_{\alpha \in \Aut(\HH)}$.
\begin{defi}\label{defi:qc-surfaces}
Let $f:\Sigma_0 \to \Sigma$ be an orientation-preserving homeomorphism as above.
\vskip 1mm
\begin{enumerate}
\item[(i)] $f$ is said to be \emph{K-quasiconformal} if it admits a lift that is $K$-quasiconformal on $\HH$ in the sense of Definition \ref{defi:quasiconformal-maps}. In this case, by Remark~\ref{rmk:K-K'-quasiconformal} and the observations made right after~\eqref{eq:qc-lift}, all the lifts of $f$ are $K$-quasiconformal.  
\vskip 1mm
\item[(ii)] $f$ is said to be \emph{quasiconformal} if it is $K$-quasiconformal for some $K \geq 1$. By Lemma~\ref{lemm:qc-standard-properties}, it is not hard to see that quasiconformality in the current sense is again preserved by forming compositions or taking inverses.
\end{enumerate}
\end{defi}
It is important to single out a particular lift as canonical. The next lemma reduces this task to Theorem~\ref{thm:qc-existence-H}.
\begin{lemm}\label{lemm:lift-criterion}
Let $f: \Sigma_0 \to \Sigma$ be a quasiconformal map in the above sense, and fix a covering map $p_0: \HH \to \Sigma_0$. Suppose also that $\widetilde{f}_0:\HH \to \HH$ is a lift of $f$ with respect to $p_0$, and write $\mu$ for its complex dilatation. Then, a given quasiconformal map $\widetilde{f}:\HH \to \HH$ is a lift of $f$ with respect to $p_0$ if and only if $\mu_{\widetilde{f}} = \mu$.
\end{lemm}
\begin{proof}
By Lemma~\ref{lemm:qc-standard-properties}(d) we have
\[
\mu_{\widetilde{f}} = \mu \Longleftrightarrow \widetilde{f}\circ \widetilde{f}_{0}^{-1} \in \Aut(\HH).
\]
The result then follows from the second observation made after~\eqref{eq:qc-lift}.
\end{proof}
\begin{defi}\label{defi:canonical-lift}
Given a quasiconformal map $f:\Sigma_0 \to \Sigma$, and a covering map $p_0:\HH \to \Sigma_0$, by Lemma~\ref{lemm:lift-criterion} and Theorem~\ref{thm:qc-existence-H}, there is a unique lift of $f$ with respect to $p_0$ that fixes $0$, $1$, and $\infty$, which we henceforth refer to as the \emph{canonical lift} of $f$ with respect to $p_0$. When the choice of $p_0$ is clear from the context, we omit it from the terminology.
\end{defi}

The next lemma records a standard fact that we need later. Suppose $S$ and $S'$ are two closed, oriented surfaces of genus $g$, and that $f: S \to S'$ is a smooth, orientation-preserving diffeomorphism. Equip $S$ and $S'$, respectively, with conformal classes of metrics $[\gamma_0]$ and $[\gamma]$, and let $J_0$ and $J$ denote the complex structures they induce.
\begin{lemm}[\cite{Hubbard}, Section 4.8]
\label{lemm:diffeomorphisms-are-qc}
In the above setting, $f:(S, [\gamma_0]) \to (S', [\gamma])$ is a quasiconformal map in the sense of Definition \ref{defi:qc-surfaces}. Moreover, in terms of the following bundle maps from $TS$ to $TS'$:
\begin{equation}\label{eq:1-form-decomposition}
d^{1, 0}f = \frac{1}{2}(df - J \circ df \circ J_0),\quad d^{0, 1}f = \frac{1}{2}(df + J \circ df \circ J_0),
\end{equation}
the complex dilataion $\mu$ of any lift of $f$ satisfies 
\begin{equation}\label{eq:diffeo-qc-bound}
\|\mu\|_{\infty; \HH} = \sup_{x \in S} \frac{|(d^{0, 1}f)_{x}|}{|(d^{1, 0}f)_{x}|} < 1,
\end{equation}
where the norms on the right-hand side can be computed with respect to any representatives of $[\gamma_0]$ and $[\gamma]$.
\end{lemm}
\begin{proof}
We include a proof for use later in Section~\ref{subsec:hyperbolic}. The idea is taken from~\cite[Section 4.8]{Hubbard}. Suppose $\widetilde{f}$ is a lift of $f$, and let $p_0$ and $p$ be as in diagram \eqref{eq:qc-lift}, with $\Sigma_0 = (S, [\gamma_0])$ and $\Sigma = (S', [\gamma])$. Then $\widetilde{f}$ preserves orientation, and thus
\begin{equation}\label{eq:tilde-f-orientation}
|\widetilde{f}_{z}|^2 - |\widetilde{f}_{\overline{z}}|^2 > 0 \quad \text{ on }\HH.
\end{equation}
In particular, setting 
\[
\mu: = \widetilde{f}_{\overline{z}}/ \widetilde{f}_{z}
\]
gives a well-defined smooth function on $\HH$ such that $|\mu| < 1$ everywhere. We next verify \eqref{eq:diffeo-qc-bound}, which will imply that $f$ is a quasiconformal map from $(S, [\gamma_0])$ to $(S', [\gamma])$ and thus conclude the proof.

Denote by $\be_1, \be_2$ the standard basis of $\RR^2$, and let $I:\RR^2 \to \RR^2$ be the standard complex structure, so that $I(\be_1) = \be_2$ and $I(\be_2) = -\be_1$. Viewing $\widetilde{f}$ as a mapping between subsets of $\RR^2$, we set
\begin{equation}\label{eq:splitting-for-qc-rmk}
d^{1, 0}\widetilde{f} = \frac{1}{2}(d\widetilde{f} - I \circ d\widetilde{f} \circ I),\quad d^{0, 1}\widetilde{f} = \frac{1}{2}(d\widetilde{f} + I \circ d\widetilde{f} \circ I).
\end{equation}
Also, given $v \in \RR^{2}$ and $a + bi \in \CC$, with $a, b$ real, we write 
\[
(a + i b) \cdot v : = a v + b I(v).
\]
Then the equation $\widetilde{f}_{\overline{z}} = \mu \widetilde{f}_{z}$ translates into
\begin{equation}\label{eq:Beltrami-in-R2}
(d^{0, 1}\widetilde{f})_{\zeta}(\be_{1}) = \mu(\zeta) \cdot (d^{1, 0}\widetilde{f})_{\zeta}(\be_{1})\quad\text{for all }\zeta \in \HH.
\end{equation}
Turning to the bundle maps defined by \eqref{eq:1-form-decomposition}, since $f$ preserves orientation, we have, similar to \eqref{eq:tilde-f-orientation}, that 
\begin{equation}\label{eq:f-orientation}
|d^{1, 0}f|^2 - |d^{0, 1}f|^2 > 0 \quad\text{on }S. 
\end{equation}
The compactness of $S$ then implies
\begin{equation}\label{eq:downstairs-sup-bound}
k: = \sup_{x\in S}\frac{|(d^{0, 1}f)_{x}|}{|(d^{1, 0}f)_{x}|} < 1.
\end{equation}
On the other hand, using
\[
d^{1, 0}f \circ J_0 = J\circ d^{1, 0}f,
\]
we infer from~\eqref{eq:f-orientation} that $d^{1, 0}f:TS \to TS'$ is a bundle isomorphism. Notice also that, by the above relation and its analogue for $d^{0, 1}f$, we have
\begin{equation}\label{eq:1,0-part-norm-formula}
|(d^{1, 0}f)_{x}|^2 |V|^2 = 2|(d^{1, 0}f)_{x}(V)|^2, \quad \text{ for all }x \in S,\ V \in T_{x}S,
\end{equation}
and a similar formula holds for $d^{0, 1}f$.

To relate the quantity~\eqref{eq:downstairs-sup-bound} to $\|\mu\|_{\infty; \HH}$, notice that the vertical maps in \eqref{eq:qc-lift} are holomorphic, which implies that
\begin{equation}\label{eq:splitting-preserved-vertical}
dp \circ d^{1, 0}\widetilde{f} = d^{1, 0}f \circ dp_0,\quad dp \circ d^{0, 1}\widetilde{f} = d^{0, 1}f \circ dp_0.
\end{equation}
Combining this with \eqref{eq:Beltrami-in-R2} leads to
\begin{equation}\label{eq:Beltrami-computation-1}
\begin{split}
(d^{0, 1}f \circ dp_0)(\be_{1}) =\ & (dp \circ d^{0, 1}\widetilde{f})(\be_1) = dp (\mu \cdot d^{1, 0}\widetilde{f}(\be_{1}))\\
=\ & (\re \mu) (d^{1, 0}f\circ dp_0)(\be_1) + (\im \mu) (J \circ d^{1, 0}f \circ dp_0)(\be_{1}).
\end{split}
\end{equation}
For later use, we observe from the first line of~\eqref{eq:Beltrami-computation-1} that
\[
(d^{0, 1}f \circ dp_0)(\be_{1})  =  (dp \circ d^{1, 0}\widetilde{f})(\mu \cdot \be_{1}) = (d^{1, 0}f \circ dp_0) (\mu \cdot \be_{1}).
\]
As a result, for all $\zeta \in \HH$, we have
\begin{equation}\label{eq:Beltrami-as-antilinear-map}
\mu(\zeta) \cdot \be_{1} = (dp_0)^{-1} \circ (d^{1, 0}f)^{-1} \circ d^{0, 1}f \circ dp_0(\be_{1}),
\end{equation}
where $dp_0$ stands for $(dp_0)_{\zeta}$, while $d^{1, 0}f$ and $d^{0, 1}f$ are evaluated at $p_0(\zeta)$. Returning to~\eqref{eq:Beltrami-computation-1}, we find after taking the norm and recalling~\eqref{eq:1,0-part-norm-formula} that
\[
|\mu(\zeta)|^2 = \frac{|d^{0, 1}f(dp_0(\be_1))|^2}{|d^{1, 0}f(dp_0(\be_1))|^2} = \frac{|(d^{0, 1}f)_{p_0(\zeta)}|^2}{|(d^{1, 0}f)_{p_0(\zeta)}|^2}, \quad \text{ for all }\zeta \in \HH,
\]
which combines with~\eqref{eq:downstairs-sup-bound} to give~\eqref{eq:diffeo-qc-bound}. This finishes the proof.
\end{proof}

\subsection{The Douady--Earle extension}\label{subsec:DE-extension}
Returning to quasiconformal maps in the sense of Definition \ref{defi:quasiconformal-maps}, a special class of such maps from $\bB$ to itself arise when certain homeomorphisms of $\partial \bB = S^1$ are extended into $\bB$ in a way that respects the action of $\Aut(\bB)$. This procedure is the celebrated Douady--Earle extension, which we recall in this section. Our main references are~\cite{Douady-Earle1986} and \cite[Section 5.1]{Hubbard}. Given a Borel measure $\mu$ on $S^1$ and a homeomorphism $f:S^1 \to S^1$, we denote by $f_{*}\mu$ the measure given by
\[
f_{*}\mu(E) = \mu(f^{-1}(E)), \text{ for all Borel set }E \subset S^1.
\]
The Douady--Earle extension rests on two ingredients, namely harmonic measures on $S^1 = \partial \bB$ and the notion of conformal barycenter. For $z \in \bB$, recall that the harmonic measure on $S^1$ with respect to $z$ is
\[
\eta_{z}  = \frac{1}{2\pi}\frac{1 - |z|^2}{|z -\zeta|^2} |d\zeta|.
\]
Thus, integrating a given $f \in C^0(S^1)$ with respect to $\eta_{z}$ yields the value at $z$ of the harmonic function on $\bB$ that restricts to $f$ on the boundary. Next, given a non-negative measure $\mu$ on $S^1$ with 
\[
\mu(S^1) = 1, \quad \mu(\{x\}) = 0 \text{ for all }x \in S^1,
\]
by~\cite[Proposition 5.1.3]{Hubbard}, the assignment $\xi_{\mu}: \bB \rightarrow \CC \simeq \RR^2$ defined by 
\[
\xi_{\mu}(z) = (1 - |z|^2)\int_{S^1}\frac{\zeta - z}{1 - \overline{z}\zeta} \ d\mu(\zeta),
\]
gives the unique vector field on $\bB$ with the property that $\xi_{\mu}(0) = \int_{S^1}\zeta\ d\mu(\zeta)$ and that
\[
(D\gamma)_{z}(\xi_{\mu}(z)) = \xi_{\gamma_* \mu}(\gamma(z)), \quad \text{ for all }\gamma \in \Aut(\bB).
\]
Moreover, it turns out that $\xi_{\mu}$ has a unique zero in $\bB$, which is denoted $B(\mu)$, and called the \emph{conformal barycenter} of the measure $\mu$. 
\begin{defi}
\label{defi:DE-extension}
Given a homeomorphism $f: S^1 \to S^1$, its Douady--Earle extension $\sE_{\bB}(f):\overline{\bB} \to \overline{\bB}$ is defined by 
\[
\sE_{\bB}(f)(z) = \left\{
\begin{array}{ll}
B(f_{*}\eta_{z}), & \text{ if }z \in \bB \\
f(z),     &  \text{ if }z \in S^1.
\end{array}
\right.
\]
\end{defi}
Thus, given $z \in \bB$, the image point $w = \sE_{\bB}(f)(z)$ is characterized as the unique solution in $\bB$ to the equation
\[
0 = \int_{S^1}\frac{f(\zeta) - w}{1 - \overline{w}f(\zeta)} \cdot \frac{1 - |z|^2}{|z - \zeta|^2} |d\zeta|.
\]
In particular, we see that
\begin{equation}\label{eq:DE-Aut-B}
\sE_{\bB}(\gamma|_{S^1}) = \gamma,\quad \text{for all }\gamma \in \Aut(\bB).
\end{equation}
The next proposition collects the basic properties of the Douady--Earle extension.
\begin{prop}\label{prop:DE-facts}
Given a homeomorphism $f: S^1 \to S^1$, the Douady--Earle extension has the following properties.
\vskip 1mm
\begin{enumerate}
\item[(a)] For all $\gamma_1, \gamma_2 \in \Aut(\bB)$, we have $\sE_{\bB}(\gamma_1 \circ f \circ \gamma_2) = \gamma_1 \circ \sE_{\bB}(f) \circ \gamma_2$.
\vskip 1mm
\item[(b)] $\sE_{\bB}(f):\overline{\bB} \to \overline{\bB}$ is a homeomorphism, and $\sE_{\bB}(f)|_{\bB}: \bB \to \bB$ is a diffeomorphism.
\vskip 1mm
\item[(c)] Suppose $f_i, f$ are homeomorphisms from $S^1$ to itself, and that $f_{i}$ converges uniformly to $f$ on $S^1$. Then $\sE_{\bB}(f_i) \to \sE_{\bB}(f)$ in $C^0(\overline{\bB}) \cap C^{\infty}_{\loc}(\bB)$.
\vskip 1mm
\item[(d)] Given $K \geq 1$, there exists $K^* \geq 1$ such that if $f$ admits a $K$-quasiconformal extension to $\bB$, then $\sE_{\bB}(f)$ is $K^*$-quasiconformal.
\end{enumerate}
\end{prop}
\begin{proof}
All four statements can be found in \cite{Douady-Earle1986}. 
Part (a) is stated at the bottom of page 27. Part (b) is Theorem 1 on page 28. Part (c) is Proposition 2 on page 30, and part (d) is Remark (1) on page 32.
\end{proof}

Next, recall the biholomorphic map $F:\HH \to \bB$ defined~\eqref{eq:H-to-D}, as well as the fact that, given a quasiconformal map $w: \HH \to \HH$, the composition $\theta_F(w)= F\circ w \circ F^{-1}$ is a quasiconformal map from $\bB$ to itself, and extends to a homeomorphism of $\overline{\bB}$. In particular, it makes sense to define $\sE_{\HH}(w): \HH \to \HH$ by
\begin{equation}\label{eq:DE-extension-H}
\sE_{\HH}(w): = F^{-1} \circ \sE_{\bB}( \theta_F(w)|_{S^1}) \circ F.
\end{equation}
The next two propositions summarize the properties of $\sE_{\HH}$ that are important for us.
\begin{prop}\label{prop:DE-H-properties}
Let $w:\HH \to \HH$ be a quasiconformal map. Then the following are true.
\vskip 1mm
\begin{enumerate}
\item[(a)] $\sE_{\HH}(w)$ is a smooth quasiconformal diffeomorphism from $\HH$ to itself.
\vskip 1mm
\item[(b)] Given $\gamma_1, \gamma_2 \in \Aut(\HH)$, there holds 
\[
\sE_{\HH}(\gamma_1 \circ w \circ \gamma_2) = \gamma_1 \circ \sE_{\HH}(w) \circ \gamma_2.
\]
\vskip 1mm
\item[(c)] Suppose $\theta_F(w)|_{S^1} = \theta_{F}(\gamma)|_{S^1}$ for some $\gamma\in \Aut(\HH)$. Then $\sE_{\HH}(w) = \gamma$.
\end{enumerate}
\end{prop}
\begin{proof}
Part (a) follows from Proposition~\ref{prop:DE-facts}(b)(d). Next, with $\gamma_1, \gamma_2$ as in part (b), notice that $\theta_F(\gamma_1), \theta_F(\gamma_2)$ lie in $\Aut(\bB)$, and that
\[
\theta_F(\gamma_1 \circ w \circ \gamma_2) = \theta_F(\gamma_1) \circ \theta_F(w) \circ \theta_F(\gamma_2).
\]
We get part (b) from these observations together with Proposition \ref{prop:DE-facts}(a). Finally, under the hypothesis of part (c), since $\theta_F(\gamma) \in \Aut(\bB)$, we have by~\eqref{eq:DE-extension-H} and~\eqref{eq:DE-Aut-B} that 
\[
\theta_F(\sE_{\HH}(w)) = \sE_{\bB}(\theta_{F}(w)|_{S^1}) = \theta_{F}(\gamma), 
\]
which immediately gives the desired equality.
\end{proof}
\begin{prop}\label{prop:DE-H-continuous}
Suppose $\mu_n: \HH \to \CC$ is a sequence of continuous functions converging uniformly locally on $\HH$ to some $\mu$, and that there exists some $k \in [0, 1)$ such that 
\[
\|\mu_n\|_{\infty; \HH} \leq k \text{ for all }n.
\]
Let $w^{\mu_n}, w^\mu$ denote the quasiconformal maps yielded by Theorem~\ref{thm:qc-existence-H}. Then, we have
\[
\sE_{\HH}(w^{\mu_n}) \to \sE_{\HH}(w^{\mu}) \text{ in }C^{\infty}_{\loc}(\HH).
\]
\end{prop}
\begin{proof}
Reflect $\mu_n$ and $\mu$ across $\partial \HH$ as indicated in~\eqref{eq:qc-reflection} to obtain bounded measurable functions $\widehat{\mu}_n$ and $\widehat{\mu}$ on $\CC$ with $L^{\infty}$-norms bounded by $k$, and denote by $f_n$ and $f$ the quasiconformal maps given by applying Theorem~\ref{thm:qc-existence}, so that
\[
w^{\mu_n} = f_n|_{\HH}, \quad w^\mu = f|_{\HH}.
\]
Since $\widehat{\mu}_n \to \widehat{\mu}$ almost everywhere on $\CC$, Theorem \ref{thm:qc-dependence} implies that $f_n \to f$ in $C^0_{\loc}(\partial \HH)$, and hence
\begin{equation}\label{eq:DE-continuous-boundary-1}
\theta_F(w^{\mu_n}) \to \theta_F(w^{\mu}) \quad \text{ in }C^{0}_{\loc}(S^1 \setminus \{1\}).
\end{equation}
To get convergence near $1$ as well, consider $\sigma \in \Aut(\HH)$ given by 
\[
\sigma(z) = \frac{1}{1 - z},
\]
and observe that it sends $(0, 1, \infty)$ to $(1, \infty, 0)$. In particular, the quasiconformal maps
\[
v_n = \sigma \circ w^{\mu_n}\circ \sigma^{-1}, \quad v = \sigma \circ w^{\mu} \circ \sigma^{-1},
\]
continue to fix $0$, $1$, and $\infty$. On the other hand, with the help of~\eqref{eq:composition-Beltrami-coefficient}, we see that the complex dilatations of these maps satisfy
\[
\mu_{v_n} = \big( \mu_n \cdot \frac{\sigma_{z}}{\overline{\sigma}_{\overline{z}}} \big) \circ \sigma^{-1} \longrightarrow \big( \mu \cdot \frac{\sigma_{z}}{\overline{\sigma}_{\overline{z}}} \big) \circ \sigma^{-1} = \mu_{v}, \quad\text{in }C^{0}_{\loc}(\HH).
\]
Thus, applying to the sequence $(v_n)$ the argument leading to~\eqref{eq:DE-continuous-boundary-1}, and noting that 
\[
\theta_{F}(v_n) = \theta_F(\sigma) \circ \theta_F(w^{\mu_n}) \circ \theta_F(\sigma)^{-1},
\]
and similarly for $\theta_{F}(v)$, we arrive at
\begin{equation}\label{eq:DE-continuous-boundary-2}
\theta_F(\sigma) \circ \theta_F(w^{\mu_n}) \circ \theta_F(\sigma)^{-1} \to \theta_F(\sigma) \circ \theta_{F}(w^{\mu})\circ \theta_F(\sigma)^{-1}\quad \text{ in }C^0_{\loc}(S^1 \setminus \{1\}).
\end{equation}
Since $1 = \theta_F(\sigma)^{-1}(-1)$, we see upon combining~\eqref{eq:DE-continuous-boundary-2} with~\eqref{eq:DE-continuous-boundary-1} that $\theta_F(w^{\mu_n}) \to \theta_F(w^{\mu})$ uniformly on all of $S^1$. It follows from Proposition~\ref{prop:DE-facts}(c) and the definition~\eqref{eq:DE-extension-H} that $\sE_{\HH}(w^{\mu_n}) \to \sE_{\HH}(w^{\mu})$ in $C^{\infty}_{\loc}(\HH)$, as claimed.
\end{proof}

\subsection{Fuchsian models and fundamental domains}\label{subsec:Fuchsian}
Let $\Sigma$ be a closed Riemann surface of genus at least $2$. Given a covering map $p: \HH \to \Sigma$, we refer to the group $\Gamma$ of deck transformations as a \emph{Fuchsian model} of $\Sigma$, which is among other things a discrete subgroup of $\Aut(\HH)$ with respect to uniform convergence on compact subsets of $\HH$. It is a standard result that if $\alpha_n$ is a sequence in $\Gamma$ that converges uniformly locally on $\HH$ to a limit that lies in $\Aut(\HH)$, then eventually the sequence coincides with the limit. (See the proof of~\cite[Lemma 2.16]{Imayoshi-Taniguchi}.) 

Another standard property of $\Gamma$ is that all its elements except the identity are hyperbolic (\cite[Lemma 2.22]{Imayoshi-Taniguchi}), where recall that an element of $\Aut(\HH)\setminus \{\id\}$ is classified as \emph{hyperbolic} if it fixes exactly two points on $\RR \cup \{\infty\}$. It can be shown, as remarked on \cite[page 47]{Imayoshi-Taniguchi}, that upon considering instead $\gamma \cdot \Gamma \cdot \gamma^{-1}$ for a suitable $\gamma \in \Aut(\HH)$, which amounts to switching from $p$ to the covering map $p \circ \gamma^{-1}$, we may assume additionally that each of $0$, $1$, and $\infty$ is fixed by some element of $\Gamma \setminus \{\id\}$. In the remainder of this paper, whenever Fuchsian models appear, we assume they satisfy this last condition.

For our purposes, yet another important fact is the existence of a fundamental domain for the action of $\Gamma$ on $\HH$. Specifically, denoting by $d_{\hyp}$ the distance function on $\HH$ induced by the hyperbolic metric $g_{\hyp}$, and fixing a reference point $w \in \HH$, the \emph{Dirichlet region} for $\Gamma$ with center $w$ is given by (see~\cite[Definition 9.4.1]{Beardon})
\[
\cD(w) = \bigcap_{\gamma \in \Gamma\setminus \{\id\}} \{z \in \HH\ |\ d_{\hyp}(z, w) < d_{\hyp}(z, \gamma(w))\}.
\]
It can be shown that $\cD(w)$ is a fundamental domain for $\Gamma$ which is locally finite. (See~\cite[Theorem 9.4.2]{Beardon}, as well as Definitions 9.3.1, 9.1.1, and 9.2.3 therein.) Since the quotient of $\HH$ by the action of $\Gamma$, being biholomorphic to the Riemann surface we started with, is compact, it follows that $\cD(w)$ is bounded with respect to $d_{\hyp}$ and thus has compact closure in $\HH$ (\cite[page 217, Exercise 3]{Beardon}). The next proposition summarizes what we need from the discussion in this paragraph. We denote by $\Gamma z$ the orbit of $z$ under the action of $\Gamma$.
\begin{prop}[\cite{Beardon}]
\label{prop:fundamental-domain}
Let $\Gamma$ be a Fuchsian model for a closed Riemann surface with genus at least $2$. There exists a compact set $K \subset \HH$ such that 
\[
K \cap \Gamma z \neq\emptyset,\quad \text{ for all }z \in \HH.
\]
\end{prop}
\subsection{Hyperbolic metrics}\label{subsec:hyperbolic}
Let $S$ be a closed oriented surface with genus $g > 1$, and fix a reference metric $\gamma_0$ from the collection $\Met_{-1}$, which, we recall, denotes the set of metrics on $S$ with constant curvature $-1$, equipped with the topology coming from smooth convergence of metrics. Recall also that $\Diff_0$ denotes the group of orientation-preserving diffeomorphisms of $S$ which are isotopic to the identity. These act on $\Met_{-1}$ by 
\[
\begin{split}
\Diff_0 \times \Met_{-1} & \longrightarrow \Met_{-1}\\
(f, \gamma)\ \ \ \ \ \  & \longmapsto \ f^* \gamma.
\end{split}
\]
As in \eqref{eq:hyperbolic-harmonic-id}, we set
\[
\Met^*_{-1} := \{\gamma \in \Met_{-1}\ |\  \id: (S, \gamma_0) \to (S, \gamma) \text{ is a harmonic map}\}.
\]
A well-known fact, recalled as Proposition \ref{prop:Met-*-representative} below, is that each orbit of the action of $\Diff_0$ meets $\Met_{-1}^*$ exactly once. Its proof is based on the following theorem due to Schoen and Yau \cite{Schoen-Yau1978}.
\begin{thm}[\cite{Schoen-Yau1978}, Corollary on page 272]
\label{thm:canonical-harmonic-mapping}
Let $S_1$ and $S_2$ be closed, oriented surfaces of genus at least $2$, equipped respectively with metrics $\gamma_1$ and $\gamma_2$, the latter assumed to have non-positive curvature. Let $\varphi: S_1 \to S_2$ be a mapping with degree $1$. Then there exists a unique harmonic map from $(S_1, \gamma_1)$ to $(S_2, \gamma_2)$ homotopic to $\varphi$. Moreover this harmonic map is an orientation-preserving diffeomorphism.
\end{thm}

Given $\gamma \in \Met_{-1}$, we write $f_{\gamma}$ for the harmonic diffeomorphism from $(S, \gamma_0)$ to $(S, \gamma)$ obtained by applying Theorem~\ref{thm:canonical-harmonic-mapping} with $\varphi$ being the identity map. Notice that, if $(\gamma_n)$ is a sequence in $\Met_{-1}$ converging smoothly to some $\gamma \in \Met_{-1}$, then $(f_{\gamma_n})$ converges smoothly to $f_{\gamma}$. This follows from combining the uniqueness part of Theorem \ref{thm:canonical-harmonic-mapping} with standard \textit{a priori} estimates for harmonic maps into closed, non-positively curved targets. (Pointwise gradient estimates can be obtained with the help of~\cite[equation (2.3) and Theorem 2.9]{Schoen-harmonic}. Higher-order estimates can then be derived in a routine manner.)

\begin{prop}[\cite{Earle-Eells1969}, page 36]
\label{prop:Met-*-representative}
In the above notation we have the following. 
\vskip 1mm
\begin{enumerate}
\item[(a)] Given an orbit $\cO$ of the action of $\Diff_0$, the metric $f_{\gamma}^*\gamma$ does not depend on the choice of representative $\gamma \in \cO$, and is the unique element of $\Met_{-1}^* \cap\ \cO$. 
\vskip 1mm
\item[(b)] The mapping $\cO \mapsto f_{\gamma}^*\gamma$, where $\gamma$ is any element of $\cO$, is a well-defined homeomorphism from $\Met_{-1}/\Diff_0$ onto $\Met_{-1}^*$. 
\end{enumerate}
\end{prop}
\begin{proof}
We outline the proof for the reader's convenience. Given $\varphi \in \Diff_0$ and $\gamma \in \Met_{-1}$, since $\varphi:(S, \varphi^*\gamma) \to (S, \gamma)$ is an isometry, and is homotopic to the identity, by the uniqueness statement in Theorem \ref{thm:canonical-harmonic-mapping}, there holds
\begin{equation}\label{eq:canonical-harmonic-pull-back}
\varphi \circ f_{\varphi^*\gamma} = f_{\gamma}.
\end{equation}
This has two implications. First, $f_{\gamma}^*\gamma$ is unchanged when $\gamma$ is replaced by $\varphi^*\gamma$, for any $\varphi \in \Diff_0$. Second, $f_{\varphi^*\gamma} = \id$ if and only if $\varphi = f_{\gamma}$. The two conclusions of (a) follow respectively from these observations.

Thanks to part (a), the map in (b) is a well-defined bijection. Also, the inverse mapping, which sends an element $\sigma \in \Met_{-1}^*$ to its orbit, is continuous by the definition of the quotient topology. To see that the mapping itself is continuous, suppose $\cO_n$ is a sequence of orbits converging to an orbit $\cO$. Up to taking a subsequence, we can find representatives $\gamma_n \in \cO_{n}$ and $\gamma \in \cO$ such that $(\gamma_n)$ converges smoothly to $\gamma$ as $n \to \infty$. The remarks after Theorem~\ref{thm:canonical-harmonic-mapping} then implies that 
\[
f_{\gamma_n}^*\gamma_n \to f_{\gamma}^*\gamma.
\]
The same argument shows that every subsequence of $\cO_n$ has a further subsequence along which the unique element in $\Met_{-1}^{*} \cap \cO_{n}$ converges to that in $\Met_{-1}^{*} \cap \cO$. This proves the required continuity.
\end{proof}

We end this section with another more or less standard result. It is to have an important role to play in Sections \ref{sec:conformal-reparametrization} and \ref{sec:harmonic-replacement}. 
\begin{prop}\label{prop:metric-qc-harmonic}
Fix a covering map $p_0: \HH \to (S, [\gamma_0])$. Suppose $(\sigma_n)$ is a sequence in $\Met_{-1}^*$ that converges smoothly to some $\sigma \in \Met_{-1}^*$. For each $n$, let 
\[
\widetilde{f}_n: \HH \to \HH
\]
be the canonical lift of $\id: (S, [\gamma_0]) \to (S, [\sigma_n])$ with respect to $p_0$, and define $\widetilde{f}$ similarly, with $\sigma$ in place of $\sigma_n$. Also, write $\mu_{n}$ and $\mu$ for the complex dilatations of $\widetilde{f}_n$ and $\widetilde{f}$. Then we have:
\vskip 1mm
\begin{enumerate}
\item[(a)] $(\mu_n)$ converges smoothly locally on $\HH$ to $\mu$. Moreover, there exists $k \in [0, 1)$ such that
\begin{equation}\label{eq:id-lift-Beltrami-bound}
\|\mu_n\|_{\infty; \HH} \leq k, \text{ for all }n.
\end{equation}
\vskip 1mm
\item[(b)] $(\widetilde{f}_n)$ and $(\widetilde{f}_n^{-1})$ converge smoothly locally on $\HH$ to $\widetilde{f}$ and $\widetilde{f}^{-1}$, respectively.
\end{enumerate}
\end{prop}
\begin{proof}
By Proposition~\ref{prop:qc-smooth-dependence} and Proposition~\ref{prop:metric-inverse-convergence}, we only need to prove part (a). Write $p_n$ for the covering map through which $\id \circ p_0$ is lifted to $\widetilde{f}_n$, so that we have the following diagram:
\begin{equation}\label{eq:id-harmonic-lift}
\begin{tikzcd}
\mathbb{H} \arrow[r, "\widetilde{f}_n"] \arrow[d, "p_0"'] & \mathbb{H} \arrow[d, "p_n"]  \\
(S, [\gamma_0]) \arrow[r, "\id"]                                 & (S, [\sigma_n])                       
\end{tikzcd}
\end{equation}
Also, let $I$ be the standard complex structure on $\RR^2$ as in the proof of Lemma \ref{lemm:diffeomorphisms-are-qc}. Following \eqref{eq:splitting-for-qc-rmk}, we define 
\[
d^{1, 0}\widetilde{f}_n = \frac{1}{2}(d\widetilde{f}_n - I \circ d\widetilde{f}_n \circ I),\quad \quad d^{0, 1}\widetilde{f}_n = \frac{1}{2}(d\widetilde{f}_n + I \circ d\widetilde{f}_n \circ I),
\]
so that
\begin{equation}\label{eq:Beltrami-in-R2-for-n}
d^{0, 1}\widetilde{f}_{n}(\be_{1}) = \mu_n \cdot d^{1, 0}\widetilde{f}_{n}(\be_{1}).
\end{equation}
On the other hand, let $J_n$, $J$, and $J_0$ be the complex structures on the Riemann surfaces $(S, [\sigma_n])$, $(S, [\sigma])$, and $(S, [\gamma_0])$, respectively. In analogy with~\eqref{eq:1-form-decomposition}, we define the following bundle maps:
\begin{equation}\label{eq:id-decomposition}
\begin{array}{rlrl}
F^{1, 0}_n : TS & \longrightarrow \ TS &\ \   F^{0, 1}_n : TS & \longrightarrow\  TS\\
v& \longmapsto \frac{1}{2}(v - J_nJ_0 v) &\ \  v & \longmapsto \frac{1}{2}(v + J_nJ_0 v),
\end{array}
\end{equation}
and also write $F^{1, 0}$ and $F^{0, 1}$ for the maps defined by the same formulas with $J_n$ replaced by $J$. Since both $J_n$ and $J_0$ are compatible with the orientation on $S$, we have 
\[
|F^{1, 0}_n|_{\gamma_0, \sigma_n} > |F^{0, 1}_n|_{\gamma_0, \sigma_n} \quad \text{everywhere on }S,
\]
where the subscripts indicate the metrics used to compute the norm. It follows, using also the relation $F^{1,0}_n \circ J_0 = J_n \circ F^{1, 0}_n$, that $F^{1, 0}_n$ is a bundle isomorphism. The exact same reasoning applies to $F^{1,0}$, and thus it makes sense to define
\[
\nu_n := (F^{1, 0}_n)^{-1} \circ F^{0, 1}_n, \quad \nu := (F^{1, 0})^{-1} \circ F^{0, 1}.
\]
Using~\eqref{eq:Beltrami-in-R2-for-n}, and also recalling \eqref{eq:splitting-preserved-vertical}, which in the present case gives
\[
dp_n \circ d^{1, 0}\widetilde{f}_n = F^{1, 0}_n \circ dp_0,\quad dp_n \circ d^{0, 1}\widetilde{f}_n = F^{0, 1}_n \circ dp_0,
\]
we obtain by the computations leading to~\eqref{eq:Beltrami-as-antilinear-map} that
\begin{equation}\label{eq:Beltrami-up-down}
\mu_n(z) \cdot \be_{1} = ((dp_0)_{z})^{-1} \circ (\nu_n)_{p_0(z)} \circ (dp_0)_{z}(\be_{1}),\quad \text{for all }z \in \HH.
\end{equation}
Now, since $(\sigma_n)$ converges smoothly to $\sigma$, we see, for example from~\eqref{eq:J-from-c}, that $(J_n)$ converges smoothly to $J$, and consequently $(F^{1, 0}_n)$ and $(F^{0, 1}_n)$ converge smoothly to $F^{1, 0}$ and $F^{0, 1}$. From this we deduce that $(\nu_n)$ converges smoothly to $\nu$, which together with \eqref{eq:Beltrami-up-down}, as well as the analogous relation between $\mu$ and $\nu$, implies that 
\[
\mu_n \to \mu\quad \text{in }C^{\infty}_{\loc}(\HH).
\]
This proves the first conclusion of part (a). For the second conclusion, we use \eqref{eq:diffeo-qc-bound} to get
\[
\| \mu_n \|_{\infty; \HH} = \sup_{x\in S}\frac{|(F_n^{0, 1})_{x}|_{\gamma_0, \sigma_n}}{|(F^{1, 0}_{n})_{x}|_{\gamma_0, \sigma_n}},\quad\quad \| \mu \|_{\infty; \HH} = \sup_{x\in S}\frac{|(F^{0, 1})_{x}|_{\gamma_0, \sigma}}{|(F^{1, 0})_{x}|_{\gamma_0, \sigma}}.
\]
Noting that $|F_n^{0, 1}|_{\gamma_0, \sigma_n}$ and $|F_n^{1, 0}|_{\gamma_0, \sigma_n}$ converge uniformly on $S$ to $|F^{0, 1}|_{\gamma_0, \sigma}$ and $|F^{1, 0}|_{\gamma_0, \sigma}$, respectively, with $|F^{1, 0}|_{\gamma_0, \sigma}$ having a positive lower bound, we deduce 
\[
\|\mu_n\|_{\infty; \HH} \to \|\mu\|_{\infty; \HH}.
\]
Since  $\|\mu\|_{\infty; \HH} < 1$, as is each $\|\mu_n\|_{\infty; \HH}$, we obtain $k \in [0, 1)$ so that \eqref{eq:id-lift-Beltrami-bound} holds. As noted at the beginning of the proof, having established part (a), we get part (b) from Proposition~\ref{prop:qc-smooth-dependence} and Proposition~\ref{prop:metric-inverse-convergence}.
\end{proof}

\subsection{Models of the Teichm\"uller space}\label{subsec:Teichmuller}
Fix a closed Riemann surface $\Sigma_0$ with genus $g > 1$ along with a covering map $p_0: \HH \to \Sigma_0$, and let $\Gamma_0$ denote its deck transformation group. As mentioned in Section \ref{subsec:Fuchsian}, without loss of generality we assume that each of $0$, $1$, and $\infty$ is fixed by some element of $\Gamma_0 \setminus \{\id\}$. Next, denote by $S$ the oriented surface underlying $\Sigma_0$, and by $\gamma_0$ the unique hyperbolic metric compatible with the conformal structure of $\Sigma_0$. In particular $p_0^*\gamma_0$ coincides with the hyperbolic metric $g_{\hyp}$.

The first model for the Teichm\"uller space we recall consists of equivalence classes of pairs $(\Sigma, f)$ where $\Sigma$ is a closed Riemann surface with genus $g$, and $f: \Sigma_0 \to \Sigma$ is a quasiconformal map according to Definition \ref{defi:qc-surfaces}. Two such pairs $(\Sigma, f), (\Sigma', f')$ are said to be equivalent if the composition
\[
f\circ (f')^{-1}: \Sigma' \to \Sigma
\]
is homotopic to a biholomorphic map. The set of equivalence classes is denoted $T(\Sigma_0)$. The standard topology on $T(\Sigma_0)$ is induced by the \emph{Teichm\"uller distance}. Its definition, which we do not repeat here (see for example~\cite[page 125]{Imayoshi-Taniguchi}), has the following immediate consequence: 
\vskip 2mm
\emph{Let $\tau_1 = [\Sigma_1,f_1]$ and $\tau_2 = [\Sigma_2, f_2]$ be two points in $T(\Sigma_0)$. Then $e^{d(\tau_1, \tau_2)} \leq K$ if and only if for all $K' > K$ there exists a $K'$-quasiconformal map from $\Sigma_1$ to $\Sigma_2$ homotopic to $f_2 \circ (f_1)^{-1}$.}
\vskip 2mm
The second Teichm\"uller space model we shall use is based on canonical lifts (with respect to $p_0$) of quasiconformal maps as described in Definition \ref{defi:canonical-lift}, or rather the homomorphisms from $\Gamma_0$ into $\Aut(\HH)$ that they induce. As preparation we state a pair of results that later serve to relate this second model back to $T(\Sigma_0)$.
\begin{lemm}[\cite{Imayoshi-Taniguchi}, Section 5.1]
\label{lemm:lift-Beltrami-coeff}
Let $\widetilde{f}: \HH \to \HH$ be a quasiconformal map and denote its complex dilatation by $\mu$. The following are equivalent.
\vskip 1mm
\begin{enumerate}
\item[(i)] There holds
\begin{equation}\label{eq:mu-equivariant}
\mu = (\mu\circ \gamma) \frac{\overline{\gamma}_{\overline{z}}}{\gamma_{z}},\quad \text{for all } \gamma \in \Gamma_0.
\end{equation}
\vskip 1mm
\item[(ii)] There holds 
\begin{equation}\label{eq:conjugation-aut}
\widetilde{f} \circ \gamma \circ \widetilde{f}^{-1} \in \Aut(\HH), \quad \text{for all } \gamma \in \Gamma_0.
\end{equation}
\vskip 1mm
\item[(iii)] There is a closed Riemann surface $\Sigma$ of genus $g$, and a quasiconformal map $f: \Sigma_0 \to \Sigma$, such that $\widetilde{f}$ is a lift of $f$ with respect to $p_0$.
\end{enumerate}
\end{lemm}
\begin{proof}
That (ii) is equivalent to (i) can be seen by applying~\eqref{eq:composition-Beltrami-coefficient} to $\widetilde{f} \circ (\widetilde{f} \circ \gamma^{-1})^{-1}$ and to $\widetilde{f} \circ \gamma^{-1}$. Next, assuming (ii), then the subgroup $\Gamma := \{\widetilde{f} \circ \gamma \circ \widetilde{f}^{-1}\}_{\gamma \in \Gamma_0}$ of $\Aut(\HH)$ inherits from $\Gamma_0$ the property of being discrete and containing no elements of finite order aside from $\id$. Consequently $\Sigma: = \HH/\Gamma$ is a Riemann surface in such a way that the projection $\pi:\HH \to \Sigma$ is holomorphic (\cite[Proposition 1.8.14]{Hubbard}). It can then be shown that $\pi \circ \widetilde{f}$ descends along $p_0$ to an orientation-preserving homeomorphism from $\Sigma_0$ to $\Sigma$ which, having $\widetilde{f}$ as a lift, is quasiconformal. In other words (iii) holds. Finally, assuming (iii), then there is a covering map $p:\HH \to \Sigma$ so that 
\[
p \circ \widetilde{f} = f \circ p_0.
\]
Given $\gamma \in \Gamma_0$, composing both sides above with $\gamma$ on the right shows that $\widetilde{f}\circ \gamma$ is again a lift of $f$ with respect to $p_0$, and hence must have the form
\[
\widetilde{f}\circ \gamma = \alpha \circ \widetilde{f},
\]
for some $\alpha \in \Aut(\HH)$. This proves (ii).
\end{proof}
\begin{lemm}[\cite{Imayoshi-Taniguchi}, Lemmas 5.1 and 5.2]
\label{lemm:marking-Fuchsian}
Given pairs $(\Sigma_1, f_1)$ and $(\Sigma_2, f_2)$, each representing an element of $T(\Sigma_0)$, let $\widetilde{f_1}$ and $\widetilde{f_2}$ denote the canonical lifts of $f_1$ and $f_2$ with respect to $p_0$, respectively. Then the following are equivalent:
\vskip 1mm
\begin{enumerate}
\item[(i)] $\widetilde{f}_1 \circ \gamma \circ \widetilde{f}_1^{-1} = \widetilde{f}_2 \circ \gamma \circ \widetilde{f}_2^{-1}$ for all $\gamma \in \Gamma_0$.
\vskip 1mm
\item[(ii)] $[\Sigma_1, f_1] = [\Sigma_2, f_2]$ in $T(\Sigma_0)$.
\vskip 1mm
\item[(iii)] $\theta_{F}(\widetilde{f}_1)|_{S^1} = \theta_{F}(\widetilde{f}_2)|_{S^1}$. (Here $\theta_F$ is given by~\eqref{eq:conjugation-by-F}, and the boundary values make sense by Proposition~\ref{prop:qc-extension}.)
\end{enumerate}
\end{lemm}

Following \cite[Section 5.1.2]{Imayoshi-Taniguchi}, whenever a quasiconformal map $\widetilde{f}: \HH \to \HH$ satisfies any one of the equivalent conditions in Lemma \ref{lemm:lift-Beltrami-coeff}, we associate to it the (injective) homomorphism $\theta_{\widetilde{f}}:\Gamma_0 \to \Aut(\HH)$ defined by
\begin{equation}\label{eq:theta-f-definition}
\begin{split}
\theta_{\widetilde{f}}: \gamma \ &\longmapsto \widetilde{f}\circ \gamma\circ \widetilde{f}^{-1}.
\end{split}
\end{equation}
With the above preparation, we define
\[
T(\Gamma_0)= \{\theta_{\widetilde{f}}\ |\ \widetilde{f}: \HH \to \HH  \text{ is quasiconformal and satisfies~\eqref{eq:qc-H-boundary-condition} and~\eqref{eq:conjugation-aut}}\}.
\]
We mention that by means of the so-called \emph{Fricke coordinates}, one can define an injective mapping into $\RR^{6g-6}$ from $T(\Gamma_0)$, thereby introducing a topology on the latter. Next, the spaces $T(\Gamma_0)$ and $T(\Sigma_0)$ are related by the map
\begin{equation}\label{eq:marking-to-Fuchsian}
\begin{split}
\cH: T(\Sigma_0) & \longrightarrow T(\Gamma_0)\\
[\Sigma, f] & \longmapsto \theta_{\widetilde{f}},
\end{split}
\end{equation}
where $\widetilde{f}$ denotes the canonical lift of $f$ with respect to $p_0$. This is a well-defined injection by Lemma \ref{lemm:marking-Fuchsian}, and is also surjective by the implication (ii) $\Rightarrow$ (iii) in Lemma \ref{lemm:lift-Beltrami-coeff}. That is, the map $\cH$ is a bijection. (See \cite[Proposition 5.3]{Imayoshi-Taniguchi}.) 

We pause to state two lemmas that we need in Sections \ref{subsec:admissible-subsets} and \ref{subsec:area-energy}, respectively. The proofs are basically exercises in using the definitions we just recalled. 
\begin{lemm}[See \cite{Imayoshi-Taniguchi}, page 125]
\label{lemm:qc-Fuchsian-action}
Let $(\tau_n)$ be a sequence converging to some $\tau$ in $T(\Sigma_0)$. Then there exists a sequence $\widetilde{\varphi}_n: \HH \to \HH$ of quasiconformal maps such that
\[
\widetilde{\varphi}_n \to id,\quad (\widetilde{\varphi}_n)^{-1} \to \id \quad \text{in }C^0_{\loc}(\HH),
\]
and that
\[
\cH(\tau_n)(\alpha) = \widetilde{\varphi}_n \circ \cH(\tau)(\alpha) \circ (\widetilde{\varphi}_{n})^{-1}, \quad \text{ for all }\alpha \in \Gamma_0.
\]
\end{lemm}
\begin{proof}
Define $K_n = (1 + \frac{1}{n})\cdot e^{d(\tau_n, \tau)}$, and choose representatives $(\Sigma_n, f_n)$ for $\tau_n$ and $(\Sigma, f)$ for $\tau$. Then, for each $n$ there exists a $K_n$-quasiconformal map $\varphi_n: \Sigma \to \Sigma_n$ homotopic to $f_n \circ f^{-1}$. Next, denote the canonical lifts of $f_n$, $f$ and $\varphi_n \circ f$ with respect to $p_0$ by $\widetilde{f}_n$, $\widetilde{f}$ and $\widetilde{g}_n$, respectively, so that in particular
\begin{equation}\label{eq:cH-represented}
\cH(\tau_n) = \theta_{\widetilde{f}_n}, \quad \cH(\tau) = \theta_{\widetilde{f}}.
\end{equation}
Recalling the definition of $T(\Sigma_0)$, we see that $(\Sigma_n, \varphi_n \circ f)$ also represents $\tau_n$. Combining this with the implication (ii) $\Rightarrow$ (i) of Lemma \ref{lemm:marking-Fuchsian}, upon letting $\widetilde{\varphi}_n := \widetilde{g_n} \circ \widetilde{f}^{-1}$, we have
\begin{equation}\label{eq:action-adjusted}
\theta_{\widetilde{f}_n}(\alpha) = \theta_{\widetilde{g}_n}(\alpha) = \widetilde{\varphi}_{n} \circ \theta_{\widetilde{f}}(\alpha) \circ (\widetilde{\varphi}_n)^{-1},\quad \text{ for all }\alpha \in \Gamma_0.
\end{equation}
In view of \eqref{eq:action-adjusted} and \eqref{eq:cH-represented}, it remains to prove that $\widetilde{\varphi}_n$ and its inverse both converge to $\id$ uniformly locally on $\HH$. To that end, notice that $\widetilde{\varphi}_n$ is indeed a lift of $\varphi_n$ as the notation suggests, and thus is $K_n$-quasiconformal in the sense of Definition \ref{defi:quasiconformal-maps}. Since $K_n \to 1$, it follows that if we let $\mu_n$ be the complex dilatation of $\widetilde{\varphi}_n$, then 
\[
\|\mu_n\|_{\infty; \HH} \to 0.
\]
From this, and the fact that $\widetilde{\varphi}_n$ fixes $0$, $1$, and $\infty$, we can invoke Theorem \ref{thm:qc-dependence} and Theorem~\ref{thm:qc-existence-H} to conclude that $\widetilde{\varphi}_{n} \to \id$ uniformly locally on $\HH$. Noting, by Lemma~\ref{lemm:qc-standard-properties}(c) and the observation made after~\eqref{eq:qc-H-boundary-condition}, that $(\widetilde{\varphi}_n)^{-1}$ is also $K_n$-quasiconformal, and that it, too, fixes $0$, $1$, and $\infty$, we may repeat this last argument to show that $(\widetilde{\varphi}_n)^{-1} \to \id$ uniformly locally on $\HH$, and we are done.
\end{proof}
\begin{lemm}\label{lemm:qc-to-Teichmuller}
Suppose $(F_n)_{n \in \NN}$ and $F$ are quasiconformal $C^1$-diffeomorphisms from $\HH$ to itself, each satisfying~\eqref{eq:qc-H-boundary-condition} and~\eqref{eq:conjugation-aut}. Assume further that $F_n \to F$ in $C^1_{\loc}(\HH)$. Then $\cH^{-1}(\theta_{F_n}) \to \cH^{-1}(\theta_{F})$ in $T(\Sigma_0)$.
\end{lemm}
\begin{proof}
Define $\tau_n = \cH^{-1}(\theta_{F_n})$ and $\tau = \cH^{-1}(\theta_{F})$. Then, by Lemmas~\ref{lemm:lift-Beltrami-coeff} and~\ref{lemm:marking-Fuchsian}, $\tau_n$ and $\tau$ admit representatives $(\Sigma_n, f_n)$ and $(\Sigma, f)$ so that $F_n$ and $F$ are, respectively, the canonical lifts of $f_n$ and $f$ with respect to $p_0$. It follows that $\theta_{F}(\Gamma_0)$ is a Fuchsian model for $\Sigma$. Also, $F_n \circ F^{-1}$ is a lift of $f_n \circ f^{-1}:\Sigma \to \Sigma_n$, and thus, defining 
\[
\mu_n: = \text{complex dilatation of }F_n \circ F^{-1},
\]
we have $\|\mu_n\|_{\infty;\HH} < 1$, and that
\[
e^{d(\tau_n, \tau)} \leq \frac{1 + \|\mu_n\|_{\infty; \HH}}{1 - \|\mu_n\|_{\infty; \HH}}.
\]
To estimate the norm of $\mu_{n}$, we observe the following. First, our assumption that $F_n$ converges to $F$ in $C^{1}_{\loc}(\HH)$ implies that $(F_n \circ F^{-1})_{z} \to 1$ and $(F_n \circ F^{-1})_{\overline{z}} \to 0$ uniformly locally, from which we deduce
\begin{equation}\label{eq:mu-n-C0-loc-converge}
\mu_n \to 0\quad \text{in }C^{0}_{\loc}(\HH).
\end{equation}
Secondly, because 
\[
(F_n \circ F^{-1}) \circ \theta_{F}(\alpha) \circ (F_n \circ F^{-1})^{-1} = \theta_{F_n}(\alpha) \in \Aut(\HH), \quad \text{for any }\alpha \in \Gamma_0,
\]
we get from the implication (ii) $\Rightarrow$ (i) of Lemma \ref{lemm:lift-Beltrami-coeff}, or rather its proof, that 
\begin{equation}\label{eq:abs-mu-n-invariant}
|\mu_n| \circ \theta_{F}(\alpha) = |\mu_n|,\quad \text{for all }\alpha \in \Gamma_0.
\end{equation} 
Now, as noted earlier, $\theta_{F}(\Gamma_0)$ is a Fuchsian model for the closed Riemann surface $\Sigma$. As such, Proposition \ref{prop:fundamental-domain} applies, which together with \eqref{eq:abs-mu-n-invariant} gives some compact set $K \subset \HH$ such that 
\[
\|\mu_n\|_{\infty; \HH} = \|\mu_n\|_{\infty; K}.
\]
Combining this with \eqref{eq:mu-n-C0-loc-converge} yields $\|\mu_n\|_{\infty; \HH} \to 0$, and hence $e^{d(\tau_n, \tau)} \to 1$, which gives the asserted convergence.
\end{proof}

Resuming our discussion of the Teichm\"uller space, the third model we need is the space $\Met_{-1}^{*}$ from Section~\ref{subsec:hyperbolic}, which is related to $T(\Sigma_0)$ by the map
\begin{equation}\label{eq:met-to-teich}
\begin{split}
\fc: \Met_{-1}^* \longrightarrow\ & T(\Sigma_0)\\
\sigma \longmapsto\ & [(S, [\sigma]), \id],
\end{split}
\end{equation}
where by $(S, [\sigma])$ we mean the Riemann surface consisting of $S$ equipped with the conformal class $[\sigma]$. Note that the pair $((S, [\sigma]), \id)$ does represent an element of $T(\Sigma_0)$, thanks to Lemma \ref{lemm:diffeomorphisms-are-qc}. Also, using Theorem \ref{thm:canonical-harmonic-mapping} and Proposition \ref{prop:Met-*-representative}(a), it can be shown that $\fc$ is a bijection.

It is a deep theorem that the spaces $T(\Sigma_0)$, $T(\Gamma_0)$ and $\Met_{-1}^*$ are all homeomorphic to $\RR^{6g-6}$. In particular, we will need the following aspect of this equivalence.
\begin{thm}[See for instance \cite{Daskalopoulos-Wentworth2007}, Theorem 2.11]
\label{thm:Teichmuller-equivalent}
The map $\fc:\Met_{-1}^{*} \to T(\Sigma_0)$ is a homeomorphism, where the topology on $\Met_{-1}^*$ is the one induced by smooth convergence of metrics, while that on $T(\Sigma_0)$ comes from the Teichm\"uller distance.
\end{thm}
\section{Analytical tools (I): Conformal reparametrization}\label{sec:conformal-reparametrization}
Let $M$ be a closed Riemannian manifold isometrically embedded into some Euclidean space $\RR^{q}$. Fix, as in Section \ref{subsec:Teichmuller}, a Riemann surface $\Sigma_0 = (S, [\gamma_0])$ with genus $g > 1$, where $S$ is the underlying closed oriented surface, and $\gamma_0$ is a hyperbolic metric. Choose also a covering map $p_0: \HH \to \Sigma_0$, and denote its deck transformation group by $\Gamma_0$. Again, without loss of generality we assume that each of $0$, $1$, and $\infty$ is fixed by some element of $\Gamma_0\setminus \{\id\}$. We mention once again that the pullback metric $p_0^*\gamma_0$ coincides with $g_{\hyp}$, the hyperbolic metric on $\HH$.

The purpose of this section is to construct the conformal reparametrization map $\Upsilon_{\eta}$ and the homotopy $\Xi_{\eta}$, mentioned in the outline in Section \ref{subsec:statements}. This involves passing to the universal cover $\HH$ and back. Section \ref{subsec:mappings} sets the stage, and the actual construction occupies Section \ref{subsec:area-energy}. The main result is Proposition \ref{prop:conformal-reparametrization}.

\subsection{Spaces of mappings}\label{subsec:mappings}
We work with two spaces whose elements are pairs consisting of a point in the Teichm\"uller space and a mapping into $M$. One of them, already defined in \eqref{eq:cM'-definition}, is the set 
\[
\cM' = \Met_{-1}^{*} \times (C^0 \cap W^{1, 2})(S; M),
\]
equipped with the product topology induced by smooth convergence of metrics and $C^0 \cap W^{1, 2}$ convergence of mappings. Thus, to say that a sequence $(\sigma_n, v_n)$ in $\cM'$ converges to $(\sigma, v) \in \cM'$ means that $\sigma_n \to \sigma$ smoothly on $S$, and that 
\[
\|v_n - v\|_{\infty; S} + \|dv_n - dv\|_{2; S} \to 0,
\]
where we use the reference metric $\gamma_0$ to compute the second term.

Next, given $\tau \in T(\Sigma_0)$, recall that the map $\cH$ in~\eqref{eq:marking-to-Fuchsian} assigns to it an injective homomorphism from $\Gamma_0$ into $\Aut(\HH)$, whose image we denote by $\Gamma_{\tau}$; that is,
\begin{equation}\label{eq:Gamma-tau-definition}
\Gamma_{\tau} : = \cH(\tau)(\Gamma_0).
\end{equation}
A map $u \in (C^0 \cap W^{1, 2})_{\loc}(\HH; M)$ is said to be \emph{$\Gamma_{\tau}$-invariant} if
\[
u \circ \gamma = u\quad \text{for all }\gamma \in \Gamma_{\tau}.
\]
We then consider the space
\begin{equation}\label{eq:cM-definition}
\cM = \{(\tau, u) \in T(\Sigma_0) \times  (C^0 \cap W^{1, 2})_{\loc}(\HH; M)\ |\ u \text{ is }\Gamma_{\tau} \text{-invariant }\}.
\end{equation}
The topology on $T(\Sigma_0) \times  (C^0 \cap W^{1, 2})_{\loc}(\HH; M)$ is given by the Teichm\"uller distance in the first factor, and by $C^0 \cap W^{1, 2}$-convergence on compact sets in the second factor. We equip $\cM$ with the subspace topology.

We next define a map from $\cM$ to $\cM'$ that will turn out to be a homeomorphism (Proposition \ref{prop:coming-down} below). Given $(\tau, u) \in \cM$, in terms of the bijection $\fc$ defined in \eqref{eq:met-to-teich}, we let 
\[
\sigma = \fc^{-1}(\tau) \in \Met^*_{-1},
\]
so that $\tau = [(S, [\sigma]), \id]$. Letting $\widetilde{f}$ denote the canonical lift of $\id:(S, [\gamma_0]) \to (S, [\sigma])$ with respect to $p_0$, we consider the following diagram, where $v$ is to be defined shortly:
\begin{equation}\label{eq:Phi-diagram}
\begin{tikzcd}
\mathbb{H} \arrow[r, "\widetilde{f}"] \arrow[d, "p_0"']    & \mathbb{H} \arrow[d, "p"] \arrow[r, "u"] & M \\
(S, [\gamma_0]) \arrow[r, "\id"] 
& {(S, [\sigma])}  \arrow[ru, "v"']      &  
\end{tikzcd}
\end{equation}
From~\eqref{eq:marking-to-Fuchsian} we have $\cH(\tau) = \theta_{\widetilde{f}}$, and thus $\theta_{\widetilde{f}}(\Gamma_0) = \Gamma_\tau$. As $u$ is $\Gamma_{\tau}$-invariant, it follows that $u \circ \widetilde{f}$ is $\Gamma_0$-invariant, and hence descends along $p_0$ to a map $v: S \to M$, which satisfies 
\begin{equation}\label{eq:Phi-u-v-relation}
v \circ p_0 = u \circ \widetilde{f},\quad \text{or equivalently}\quad v \circ p = u.
\end{equation}
Notice that since $\widetilde{f}$ is a smooth diffeomorphism, standard results on Sobolev spaces (for example \cite[Theorem 3.41]{Adams-Fournier}) shows that $u \circ \widetilde{f}$ again lies in $C^0 \cap W^{1, 2}$ locally on $\HH$, which implies by~\eqref{eq:Phi-u-v-relation} that $v \in C^0 \cap W^{1, 2}(S; M)$. We then define a map $\Phi: \cM \to \cM'$ by 
\begin{equation}\label{eq:Phi-definition}
\begin{split}
\Phi: (\tau, u) & \longmapsto  (\sigma, v).
\end{split}
\end{equation}
Conversely, given $(\sigma, v) \in  \cM'$, we define $\Psi(\sigma, v)$ to be the pair $(\tau, u) \in \cM$ where
\begin{equation}\label{eq:Psi-definition}
\tau := \fc(\sigma) = [(S, [\sigma]), \id],\quad\quad  u := v \circ p_0 \circ \widetilde{f}^{-1} = v \circ p,
\end{equation}
where, as before, $\widetilde{f}$ is the canonical lift of $\id:(S, [\gamma_0]) \to (S, [\sigma])$ with respect to $p_0$.
\begin{prop}\label{prop:coming-down}
The maps $\Phi: \cM \to  \cM'$ and $\Psi:  \cM'\to \cM$ are inverses of each other, and are both continuous.
\end{prop}
\begin{proof}
That $\Phi$ and $\Psi$ are inverses of each other can be checked by tracing their definitions. Next, to see that $\Phi$ is continuous, suppose $(\tau_n, u_n)$ is a sequence in $\cM$ converging to some $(\tau, u) \in \cM$ and let 
\[
(\sigma_n, v_n) = \Phi(\tau_n, u_n),\quad (\sigma, v) = \Phi(\tau, u).
\]
Since $\tau_n\to \tau$ in $T(\Sigma_0)$ by assumption, it follows from Theorem~\ref{thm:Teichmuller-equivalent} that $\sigma_n$ converges smoothly to $\sigma$. To show that $v_n$ converges to $v$ in $C^0 \cap W^{1, 2}$, we adopt the notation in diagram~\eqref{eq:Phi-diagram} and also let $\widetilde{f}_n$ denote the canonical lift of $\id:\Sigma_0 \to (S, [\sigma_n])$ with respect to $p_0$, as indicated in the following diagram:
\begin{equation}\label{eq:n-Phi-diagram}
\begin{tikzcd}
\mathbb{H} \arrow[r, "\widetilde{f}_n"] \arrow[d, "p_0"']            & \mathbb{H} \arrow[d, "p_n"] \arrow[r, "u_n"] & M \\
(S, [\gamma_0]) \arrow[r, "\id"]  & {(S, [\sigma_n])}          \arrow[ru, "v_n"']                       &  
\end{tikzcd}
\end{equation}
Since $\sigma_n$ converges smoothly to $\sigma$, Proposition~\ref{prop:metric-qc-harmonic} implies that $\widetilde{f}_n$ converges smoothly locally on $\HH$ to $\widetilde{f}$. In particular, given a compact subset $K$ of $\HH$, there exist another compact set $K' \subset \HH$ and a constant $A > 0$ such that for large enough $n$, we have for all $x \in K$ that 
\begin{equation}\label{eq:area-formula-conditions}
\widetilde{f}_n(x), \widetilde{f}(x) \in K' \quad \text{and} \quad
A^{-1} \leq \det(d\widetilde{f}_n(x)) \leq A.
\end{equation}
With the help of the first property in~\eqref{eq:area-formula-conditions}, we have $u_n \circ \widetilde{f}_n \to u \circ \widetilde{f}$ uniformly on $K$. To get $L^2$-convergence of the gradients, we note that 
\[
\begin{split}
\|(du_n)_{\widetilde{f}_n} d\widetilde{f}_n - (du)_{\widetilde{f}} d\widetilde{f}\|_{2; K} \leq\ & \underbrace{\|(du_n)_{\widetilde{f}_n} - (du)_{\widetilde{f}_n}  \|_{2; K}}_{(I)} \cdot  \|d\widetilde{f}_n\|_{\infty; K}\\
&+ \underbrace{\|(du)_{\widetilde{f}_n} \|_{2; K}}_{(II)} \cdot \|d\widetilde{f}_n - d\widetilde{f}\|_{\infty; K}\\
&+ \underbrace{\|(du)_{\widetilde{f}_n}  - (du)_{\widetilde{f}}\|_{2; K}}_{(III)} \cdot \|d\widetilde{f}\|_{\infty; K}
\end{split}
\]
Using the area formula and \eqref{eq:area-formula-conditions}, we see that $(I)$ tends to $0$, while $(II)$ is bounded uniformly in $n$. By approximating $du$ in $L^2(K')$ with continuous maps, and again using the area formula, we get that $(III)$ tends to $0$ as well. From these we infer that $d(u_n \circ \widetilde{f}_n)$ converges in $L^{2}_{\loc}$ to $d(u \circ \widetilde{f})$. Since $v_n \circ p_0 = u_n \circ \widetilde{f}_n$, from what we just proved we conclude that $v_n \to v$ in $(C^{0} \cap W^{1, 2})(S; M)$.

To prove that $\Psi$ is continuous, we let $(\sigma_n, v_n)$ be a sequence in $\cM'$ that converges to some $(\sigma, v) \in \cM'$, and define
\[
(\tau_n, u_n) = \Psi(\sigma_n, v_n),\quad (\tau, u) = \Psi(\sigma, v).
\]
Theorem~\ref{thm:Teichmuller-equivalent} implies that $\tau_n \to \tau$ in $T(\Sigma_0)$. Next, since $p_0$ is a local diffeomorphism, our assumption on $(v_n)$ implies that $v_n \circ p_0$ converges to $v \circ p_0$ in $W^{1, 2} \cap C^0$ locally on $\HH$. On the other hand, by Proposition~\ref{prop:metric-qc-harmonic}, the diffeomorphisms $\widetilde{f}_n^{-1}$ converge smoothly locally on $\HH$ to $\widetilde{f}^{-1}$. Recalling that $u_n = v_n \circ p_0 \circ \widetilde{f}_n^{-1}$ while $u = v \circ p_0 \circ \widetilde{f}^{-1}$, and arguing as in the previous paragraph, we get that $u_n \to u$ in $C^0 \cap W^{1, 2}$ locally on $\HH$.
\end{proof}
\begin{rmk}\label{rmk:convergence-p}
Suppose $(\sigma_n, v_n)$ is a sequence in $\cM'$ converging to some $(\sigma, v) \in \cM'$. In the notation of~\eqref{eq:Phi-diagram} and~\eqref{eq:n-Phi-diagram}, we have that 
\[
p_n = p_0 \circ \widetilde{f}_n^{-1}, \quad p = p_0 \circ \widetilde{f}^{-1}.
\]
Proposition~\ref{prop:metric-qc-harmonic} then shows that $(p_n)$ converges smoothly locally to $p$ as a sequence of smooth maps from $\HH$ into $S$.
\end{rmk}
\subsection{Area and energy}\label{subsec:area-energy}
Given $(\sigma, v) \in \cM'$, the mapping area $A(v)$ and the Dirichlet energy $E(\sigma, v)$ are defined as in Section~\ref{subsec:notation}. Here we are interested in continuously assigning to each $C^2$-map $v: S \to M$ a pair $(\sigma, \widetilde{v}) \in \cM'$ in such a way that $A(\widetilde{v}) = A(v)$, and that $E(\sigma, \widetilde{v})$ is close to $A(\widetilde{v})$. We begin with a preliminary lemma, before describing the construction. The objects $\gamma_0$, $p_0$ and $\Gamma_0$ in the statement are those fixed before the start of Section \ref{subsec:mappings}. 
\begin{lemm}\label{lemm:metric-to-qc}
Suppose $(\gamma_n)$ is a sequence of $C^{1}$-metrics on $S$ converging in the $C^1$-topology to a limiting metric $\gamma$, and express the pullback metrics $p_0^*\gamma_n$ and $p_0^*\gamma$ according to \eqref{eq:metric-complex} as
\begin{equation}\label{eq:pull-back-by-p0}
p_0^*\gamma_n = \lambda_n |dz + \mu_n d\overline{z}|^2,\ \ p_0^*\gamma = \lambda |dz + \mu d\overline{z}|^2.
\end{equation}
Let $w^{\mu_n}$ and $w^{\mu}$ be the quasiconformal maps produced by Theorem~\ref{thm:qc-existence-H}. Then we have
\vskip 1mm
\begin{enumerate}
\item[(a)] $w^{\mu} \circ \alpha \circ (w^{\mu})^{-1} \in \Aut(\HH)$ for all $\alpha \in \Gamma_0$. The same holds with $w^{\mu}$ replaced by $w^{\mu_n}$.
\vskip 1mm
\item[(b)] $\mu_n \to \mu$ in $C^1_{\loc}(\HH)$. Moreover, there exists some $k \in [0, 1)$ such that
\begin{equation}\label{eq:Beltrami-coeff-bound}
\|\mu_n\|_{\infty; \HH} \leq k,\quad \text{for all }n.
\end{equation}
\vskip 1mm
\item[(c)] $w^{\mu_n}$ and $(w^{\mu_n})^{-1}$ converge in $C^1_{\loc}(\HH)$ to $w^{\mu}$ and $(w^{\mu})^{-1}$, respectively.
\end{enumerate}
\end{lemm}
\begin{proof}
We first note that the $C^1$-convergence of $(\gamma_n)$ to a limiting metric implies the existence of some $\Lambda \geq 1$ such that  
\begin{equation}\label{eq:metric-surface-bound}
\Lambda^{-1} \gamma_0 \leq \gamma_n \leq \Lambda \gamma_0,\quad \text{for all }n.
\end{equation}
Also, from the formula~\eqref{eq:lambda-mu}, we see that $\mu$ and $(\mu_n)$ are $C^{1}$-functions on $\HH$, and hence $w^{\mu}$ and $(w^{\mu_n})$ are $C^1$-diffeomorphisms from $\HH$ onto itself.

To prove (a), take an arbitrary $\alpha \in \Gamma_0$ and use $p_0 \circ \alpha = p_0$ to get
\[
\alpha^*(\lambda |dz + \mu d\overline{z}|^2) = (p_0\circ \alpha)^*\gamma = p_0^*\gamma = \lambda |dz + \mu d\overline{z}|^2,
\]
from which we obtain after a straightforward computation that
\[
\mu = (\mu \circ \alpha) \frac{\overline{\alpha}_{\overline{z}}}{\alpha_{z}}.
\]
The assertion on $w^{\mu}$ then follows from the implication (i) $\Rightarrow$ (ii) of Lemma \ref{lemm:lift-Beltrami-coeff}. The proof for $w^{\mu_n}$ is exactly the same. 

For (b), since $p_0^* \gamma_n$ converges to $p_0^* \gamma$ in $C^1_{\loc}(\HH)$ by our assumption on $(\gamma_n)$, we get from~\eqref{eq:lambda-mu} that $(\mu_n)$ converges to $\mu$ in $C^1_{\loc}(\HH)$. Next, since $p_0^* \gamma_0$ coincides with $g_{\hyp}$, the latter having the form $\rho^2 |dz|^2$ for some positive smooth function $\rho$, we obtain from~\eqref{eq:metric-surface-bound} a uniform bound of the form 
\[
\Lambda^{-1}\rho^2|dz|^2 \leq p_0^*\gamma_n \leq \Lambda \rho^2|dz|^2.
\]
Combining this with \eqref{eq:eigenvalue-lambda-mu} gives
\[
\frac{1 + |\mu_n|}{1 - |\mu_n|} \leq \Lambda \text{ on }\HH,
\]
which immediately implies \eqref{eq:Beltrami-coeff-bound} with $k = \frac{\Lambda - 1}{\Lambda + 1}$. This proves part (b). Part (c) follows from part (b) and Proposition~\ref{prop:qc-H-dependence}.
\end{proof}
For each $\eta \in (0, 1)$, we construct two maps, denoted respectively by
\[
\Upsilon_{\eta}: C^2(S; M) \to \cM',\quad \quad \Xi_{\eta}: [0, 1] \times C^2(S; M)  \to C^1(S; M),
\]
through the following procedure. 
\vskip 1mm
\begin{enumerate}
\item[(1)] We first define $\Upsilon_{\eta}$. With $g_M$ denoting the Riemannian metric on $M$, given $v \in C^2(S; M)$, we consider the $C^1$-metric $\gamma := v^*g_M + \eta \gamma_0$ on $S$, and express the pullback metric $p_0^*\gamma$ according to~\eqref{eq:metric-complex} as 
\[
p_0^*\gamma = \lambda |dz + \mu d\overline{z}|^2.
\]
By the argument leading to~\eqref{eq:Beltrami-coeff-bound} we get $\|\mu\|_{\infty; \HH} < 1$. Theorem \ref{thm:qc-existence-H} then yields a quasiconformal $C^1$-diffeomorphism $w^{\mu}: \HH \to \HH$. By Lemma \ref{lemm:metric-to-qc}(a) and the fact that $w^{\mu}$ fixes $0$, $1$, and $\infty$, the homomorphism $\theta_{w^{\mu}}:\Gamma_0 \to\Aut(\HH)$ given by~\eqref{eq:theta-f-definition} belongs to $T(\Gamma_0)$. We then use the map $\cH$ in~\eqref{eq:marking-to-Fuchsian} to define
\begin{equation}\label{eq:Upsilon-tau-definition}
\tau := \cH^{-1}(\theta_{w^{\mu}}).
\end{equation}
To augment $\tau$ to a pair in $\cM$, note that the map
\begin{equation}\label{eq:Upsilon-u-definition}
u: = v \circ p_0 \circ (w^{\mu})^{-1}: \HH \to M
\end{equation}
is of class $C^0 \cap W^{1, 2}$ (in fact it is $C^1$), and satisfies $u \circ \gamma = u$ for all $\gamma$ in $\theta_{w^{\mu}}(\Gamma_0)$, the latter coinciding with $\Gamma_{\tau}$ by \eqref{eq:Upsilon-tau-definition}. Thus $(\tau, u) \in \cM$, and we set
\begin{equation}\label{eq:upsilon-definition}
\Upsilon_{\eta}(v) := \Phi(\tau, u),
\end{equation}
where $\Phi: \cM \to \cM'$ is the map in~\eqref{eq:Phi-definition}.
\vskip 1mm
\item[(2)] With $v$, $\tau$, and $u$ as above, we next define $\Xi_{\eta}(t, v)$ for $t \in [0, 1]$. To start, let 
\[
\sigma = \fc^{-1}(\tau)\in \Met_{-1}^{*},
\]
so that $\tau = [(S, [\sigma]), \id]$ by~\eqref{eq:met-to-teich}. Denote by $\widetilde{f}$ the canonical lift of $\id:(S, [\gamma_0]) \to (S, [\sigma])$ with respect to $p_0$, and write $p$ for the covering map such that 
\[
p \circ \widetilde{f} = \id \circ p_0.
\]
By \eqref{eq:Upsilon-tau-definition} and the definition of $\cH$, we then have 
\begin{equation}\label{eq:f-w-theta-coincide}
\theta_{\widetilde{f}} = \cH(\tau) = \theta_{w^{\mu}}.
\end{equation}
The desired path $\{\Xi_{\eta}(t, v)\}_{t \in [0, 1]}$ will be obtained from another path $\{\widetilde{H}_t\}_{t \in [0, 1]}$ of quasiconformal $C^1$-diffeomorphisms that connects $w^{\mu}$ to $\widetilde{f}$ and satisfies $\theta_{\widetilde{H}_{t}} = \cH(\tau)$ for all $t \in [0, 1]$. To construct the maps $\widetilde{H}_{t}$, we follow an argument of Earle and McMullen~\cite{Earle-McMullen1988} as it is presented in the proof of~\cite[Proposition 6.4.9]{Hubbard}. To start, let 
\[
\nu = \text{complex dilatation of $w^{\mu} \circ \widetilde{f}^{-1}$},
\]
and observe the following two properties. First, since the complex dilatation of $\widetilde{f}$ is smooth, while that of $w^{\mu}$ is $C^{1}$, we get from~\eqref{eq:composition-Beltrami-coefficient} that $\nu \in C^1_{\loc}(\HH)$. Secondly, by~\eqref{eq:f-w-theta-coincide}, there holds for all $\alpha \in \theta_{\widetilde{f}}(\Gamma_0)$ that
\[
(w^{\mu} \circ \widetilde{f}^{-1}) \circ \alpha \circ (w^{\mu} \circ \widetilde{f}^{-1})^{-1} = \alpha \in \Aut(\HH),
\]
and hence by the implication (ii) $\Rightarrow$ (i) in Lemma~\ref{lemm:lift-Beltrami-coeff}, we have
\begin{equation}\label{eq:nu-invariant}
\nu = (\nu \circ \alpha) \frac{\overline{\alpha}_{\overline{z}}}{\alpha_{z}},\quad \text{for all }\alpha \in \theta_{\widetilde{f}}(\Gamma_0).
\end{equation}
Given $t \in [0, 1]$, in the notation of Theorem \ref{thm:qc-existence-H}, we define
\begin{equation}\label{eq:tilde-g-definition}
\widetilde{g}_{t}: = w^{(1-t)\nu}:\HH \to \HH,
\end{equation}
which is a quasiconformal $C^1$-diffeomorphism since $\nu$ is $C^1$ on $\HH$. Moreover, by \eqref{eq:nu-invariant} and the implication (i) $\Rightarrow$ (ii) of Lemma~\ref{lemm:lift-Beltrami-coeff}, it satisfies
\begin{equation}\label{eq:g_t-Fuchsian}
\theta_{\widetilde{g}_t}(\alpha) \in \Aut(\HH),\quad \text{for all }\alpha \in \theta_{\widetilde{f}}(\Gamma_0).
\end{equation}
Observe also that $\widetilde{g}_{1}= \id$, while from the uniqueness part of Theorem~\ref{thm:qc-existence-H} we have $\widetilde{g}_0 = w^{\mu} \circ \widetilde{f}^{-1}$. Furthermore, by \eqref{eq:f-w-theta-coincide} and the implication (i) $\Rightarrow$ (iii) of Lemma~\ref{lemm:marking-Fuchsian}, we see that, in terms of the notation~\eqref{eq:conjugation-by-F},
\[
\theta_{F}(\widetilde{g}_0)|_{S^1} = \id_{S^1}.
\]

We next adjust each $\widetilde{g}_{t}$ so that the counterpart of $\theta_{\widetilde{g}_{t}}$ along the new path is independent of $t$. Specifically, in the notation of Section~\ref{subsec:DE-extension}, we set
\[
\widetilde{h}_{t}: = \sE_{\HH}((\widetilde{g}_{t})^{-1}),
\]
which by Proposition \ref{prop:DE-H-properties}(a) is a smooth quasiconformal diffeomorphism from $\HH$ to itself. Because $\theta_F((\widetilde{g}_{0})^{-1})$ and $\theta_F((\widetilde{g}_{1})^{-1})$ both restrict to the identity on $S^1$, by Proposition \ref{prop:DE-H-properties}(c) we have
\[
\widetilde{h}_{0} = \id =  \widetilde{h}_{1}.
\]
Also, by~\eqref{eq:g_t-Fuchsian} together with Proposition \ref{prop:DE-H-properties}(b), we get 
\begin{equation}\label{eq:h-t-Fuchsian}
(\widetilde{h}_{t})^{-1} \circ \alpha \circ \widetilde{h}_{t}  = \theta_{\widetilde{g}_{t}}(\alpha),\quad \text{for all }\alpha \in \theta_{\widetilde{f}}(\Gamma_0).
\end{equation}
Letting $\widetilde{H}_{t} = \widetilde{h}_{t} \circ \widetilde{g}_{t} \circ \widetilde{f}$ for $t \in [0, 1]$, we see that each $\widetilde{H}_{t}$ is a quasiconformal $C^1$-diffeomorphism from $\HH$ to itself, and that
\begin{equation}\label{eq:upstairs-endpoints}
\widetilde{H}_{0} = w^{\mu},\quad \widetilde{H}_{1} = \widetilde{f}.
\end{equation}
Moreover, by~\eqref{eq:h-t-Fuchsian}, there holds
\begin{equation}\label{eq:conjugation-along-path}
\widetilde{H}_{t} \circ \beta = \theta_{\widetilde{f}}(\beta) \circ \widetilde{H}_{t},\quad \text{for all }\beta \in \Gamma_0.
\end{equation}
Recalling that $\theta_{\widetilde{f}}(\Gamma_0) = \cH(\tau)(\Gamma_0)$ is the deck transformation group of $p$, we obtain a $C^{1}$-diffeomorphism $H_{t}: S \to S$ such that
\[
p \circ \widetilde{H}_{t} = H_t \circ p_0.
\]
In particular, by~\eqref{eq:upstairs-endpoints} we have
\[
H_1 = \id, \quad\quad (H_0)^{-1} \circ p = p_0 \circ (w^{\mu})^{-1}.
\]
We then define
\begin{equation}\label{eq:Xi-definition}
\Xi_{\eta}(t, v) := v \circ (H_0)^{-1} \circ H_{t}.
\end{equation}
Note that this indeed lies in $C^1(S; M)$. Also, with the help of~\eqref{eq:Upsilon-u-definition}, we see that
\begin{equation}\label{eq:Xi-alt-definition}
\Xi_{\eta}(t, v) \circ p_0 = u \circ \widetilde{H}_{t}.
\end{equation}
The diagram below is our attempt to summarize the relationship between some of the maps in the preceding construction. Recall that $\widetilde{H}_{0} = w^{\mu}$, while $\widetilde{H}_{1} = \widetilde{f}$ and $H_1 = \id$.
\vskip 1em
\begin{equation}\label{eq:upsilon-diagram-big}
\begin{tikzcd}
\HH \arrow[d, "p_0"] \arrow[rr, "\widetilde{H}_{t}"] && \HH \arrow[d, "p"] \arrow[rr, "(w^{\mu})^{-1}"] \arrow[rrr, "u", bend left = 35] && \HH \arrow[d, "p_0"] & M \\
S \arrow[rr, "H_{t}"] \arrow[rrrrru, "{\Xi_{\eta}(t, v)}", bend right = 65] && S \arrow[rr, "(H_0)^{-1}"]  && S \arrow[ru, "v"] & 
\end{tikzcd}
\end{equation}
\end{enumerate}
The next proposition establishes the properties of $\Upsilon_{\eta}$ and $\Xi_{\eta}$.
\begin{prop}\label{prop:conformal-reparametrization}
We have the following:
\vskip 1mm
\begin{enumerate}
\item[(a)] For all $\eta \in (0, 1)$, the map $\Upsilon_{\eta}: C^2(S; M) \to \cM'$ is continuous. 
\vskip 1mm
\item[(b)] For all $\eta \in (0, 1)$, the map $\Xi_{\eta}: [0, 1] \times C^2(S; M) \to C^1(S; M)$ is continuous, and satisfies
\begin{equation}\label{eq:Xi-area-constant}
A(\Xi_{\eta}(t, v)) = A(v),\quad \text{for all }(t, v) \in [0, 1] \times C^2(S; M).
\end{equation}
Moreover, given $v \in C^2(S; M)$, writing $\Upsilon_{\eta}(v)$ as $(\sigma, \widetilde{v})$, we have
\begin{equation}\label{eq:reparametrization-homotopy-endpoint}
\Xi_{\eta}(0, v) = v,\ \ \Xi_{\eta}(1, v) = \widetilde{v}.
\end{equation}
\vskip 1mm
\item[(c)] Given $\ep, \Lambda > 0$, there exists $\eta \in (0, 1)$ depending only on $\ep, \Lambda$, and $\Vol(S, \gamma_0)$ such that 
\begin{equation}\label{eq:reparametrization-area-energy}
A(v) \leq E(\Upsilon_{\eta}(v)) \leq A(v) + \ep,
\end{equation}
for all $v \in C^2(S; M)$ satisfying $E(\gamma_0, v) \leq \Lambda$. 
\end{enumerate}
\end{prop}
\begin{proof}
For part (a), to see that $\Upsilon_{\eta}$ is continuous, take a sequence $(v_n)$ converging in $C^2(S; M)$ to some $v$. With $\Psi: \cM' \to \cM$ being the inverse of $\Phi$ defined by~\eqref{eq:Psi-definition}, we let
\[
(\tau_n, u_n) = \Psi(\Upsilon_{\eta}(v_n)),\ \ (\tau, u) = \Psi(\Upsilon_{\eta}(v)).
\]
By Proposition~\ref{prop:coming-down}, to conclude $\Upsilon_{\eta}(v_n) \to \Upsilon_{\eta}(v)$ in $\cM'$, it is enough to prove that $(\tau_n, u_n) \to (\tau, u)$ in $\cM$. To that end, following the notation used in the preceding construction, we let
\[
\gamma_n = v_n^*g_M + \eta \gamma_0,\ \ \gamma = v^*g_M + \eta \gamma_0,
\]
define $\mu_n$ and $\mu$ by the relation~\eqref{eq:pull-back-by-p0}, and apply Theorem~\ref{thm:qc-existence-H} to get quasiconformal maps $w^{\mu_n}$ and $w^{\mu}$. Recalling~\eqref{eq:Upsilon-tau-definition} and~\eqref{eq:Upsilon-u-definition}, we see that $u_n$ and $\tau_n$, and likewise $u$ and $\tau$, can be expressed as follows:
\begin{equation}\label{eq:upsilon-continuity-u-expression}
u_n = v_n \circ p_0 \circ (w^{\mu_n})^{-1}, \quad u = v \circ p_0 \circ (w^{\mu})^{-1},
\end{equation}
\begin{equation}\label{eq:upsilon-continuity-tau-expression}
\tau_n = \cH^{-1}(\theta_{w^{\mu_n}}), \quad \tau = \cH^{-1}(\theta_{w^{\mu}}).
\end{equation}
Noting that $(\gamma_n)$ is a sequence of $C^{1}$-metrics converging in $C^1$ to $\gamma$, we have by Lemma~\ref{lemm:metric-to-qc} that
\begin{equation}\label{eq:reparam-w-convergence}
\mu_n \to \mu, \quad w^{\mu_n} \to w^{\mu },\quad (w^{\mu_n})^{-1} \to (w^{\mu})^{-1} \quad \text{in }C^{1}_{\loc}(\HH).
\end{equation}
From the last convergence and \eqref{eq:upsilon-continuity-u-expression} we deduce
\begin{equation}\label{eq:reparam-u-convergence}
u_n  \to   u \text{ in }C^{1}_{\loc}(\HH).
\end{equation}
On the other hand, by Lemma~\ref{lemm:metric-to-qc}(a) and the second convergence in~\eqref{eq:reparam-w-convergence}, we may invoke Lemma~\ref{lemm:qc-to-Teichmuller} to see from~\eqref{eq:upsilon-continuity-tau-expression} that $\tau_n \to \tau$ in $T(\Sigma_0)$, and we are done with part (a).

For part (b), to see that $\Xi_{\eta}$ is continuous, we suppose that $t_n \to t$ in $[0, 1]$ and that $v_n \to v$ in $C^2(S; M)$, and use the same set of notations as in the proof of (a). Then from the conclusions of that proof we have
\begin{equation}\label{eq:reparam-part-b-on-a}
\tau_n \to \tau \text{ in }T(\Sigma_0),\quad\quad u_n \to u \text{ in }C^{1}_{\loc}(\HH).
\end{equation}
Next, consider the elements $\sigma_n = \fc^{-1}(\tau_n)$ and $\sigma = \fc^{-1}(\tau)$ of $\Met_{-1}^{*}$. Recall that these are characterized by the relation
\[
\tau_n = [(S, [\sigma_n]), \id], \quad \tau = [(S, [\sigma]), \id].
\]
As before, denote by $\widetilde{f}_n$ and $\widetilde{f}$ the canonical lifts with respect to $p_0$ of $\id:(S, [\gamma_0]) \to (S, [\sigma_n])$ and $\id: (S, [\gamma_0]) \to (S, [\sigma])$. By Theorem \ref{thm:Teichmuller-equivalent} and the first statement in \eqref{eq:reparam-part-b-on-a}, we see that $\sigma_n \to \sigma$ smoothly on $S$, and thus, by Proposition~\ref{prop:metric-qc-harmonic}, we have
\begin{equation}\label{eq:reparam-f-convergence}
\mu_{\widetilde{f}_n} \to \mu_{\widetilde{f}}, \quad \widetilde{f}_n \to \widetilde{f}, \quad \widetilde{f}_n^{-1} \to \widetilde{f}^{-1} \quad \text{in $C^{\infty}_{\loc}(\HH)$}.
\end{equation}
Next, we define 
\[
\nu_n = \text{complex dilatation of } w^{\mu_n} \circ \widetilde{f}_n^{-1}.
\]
By Lemma~\ref{lemm:metric-to-qc} and Proposition~\ref{prop:metric-qc-harmonic}, respectively, we get the bounds~\eqref{eq:Beltrami-coeff-bound} and~\eqref{eq:id-lift-Beltrami-bound} on the complex dilatations of $w^{\mu_n}$ and $\widetilde{f}_n$, which together with~\eqref{eq:composition-Beltrami-coefficient} shows that
\begin{equation}\label{eq:nu-Beltrami-coefficient-away-from-1}
\sup_{n}\|\nu_n\|_{\infty; \HH} < 1.
\end{equation}
Furthermore, by~\eqref{eq:reparam-f-convergence} and the first convergence in~\eqref{eq:reparam-w-convergence}, and again using~\eqref{eq:composition-Beltrami-coefficient}, we have $\nu_n \to \nu$ in $C^{1}_{\loc}(\HH)$. Consequently, writing $\widetilde{g}_n$ and $\widetilde{g}_t$, respectively, for the quasiconformal $C^1$-diffeomorphisms $w^{(1 - t_n)\nu_n}$ and $w^{(1 - t)\nu}$ given by Theorem~\ref{thm:qc-existence-H}, we have by Proposition~\ref{prop:qc-H-dependence} that
\[
\widetilde{g}_n \to \widetilde{g}_{t},\quad (\widetilde{g}_{n})^{-1} \to (\widetilde{g}_{t})^{-1} \quad \text{in }C^{1}_{\loc}(\HH).
\]
In particular, the complex dilatations of $(\widetilde{g}_n)^{-1}$ converge uniformly locally on $\HH$ to that of $(\widetilde{g}_{t})^{-1}$, and have $L^{\infty}$-norms less than $1$ by a fixed amount independent of $n$ thanks to~\eqref{eq:nu-Beltrami-coefficient-away-from-1} and~\eqref{eq:composition-Beltrami-coefficient}. Since $(\widetilde{g}_n)^{-1}$ and $(\widetilde{g}_{t})^{-1}$ all fix $0$, $1$, and $\infty$, upon letting 
\[
\widetilde{h}_n = \sE_{\HH}((\widetilde{g}_n)^{-1}), \quad \widetilde{h}_{t} = \sE_{\HH}((\widetilde{g}_{t})^{-1})
\]
and using Proposition~\ref{prop:DE-H-continuous}, we get that $\widetilde{h}_n \to \widetilde{h}_{t}$ in $C^{\infty}_{\loc}(\HH)$. Putting everything together, then, we arrive at
\[
\widetilde{H}_n: = \widetilde{h}_n \circ \widetilde{g}_n \circ \widetilde{f}_n \to \widetilde{h}_{t} \circ \widetilde{g}_{t} \circ \widetilde{f} =: \widetilde{H}_{t}\quad \text{in }C^1_{\loc}(\HH).
\]
Combining this with \eqref{eq:Xi-alt-definition} and the analogous fact $\Xi_{\eta}(t_n, v_n) \circ p_0 = u_n \circ \widetilde{H}_n$, and also recalling \eqref{eq:reparam-part-b-on-a}, we conclude that
\[
\Xi_{\eta}(t_n, v_n) \circ p_0 \to \Xi_{\eta}(t, v) \circ p_0 \text{ in }C^1_{\loc}(\HH; M),
\]
which implies that $\Xi_{\eta}(t_n, v_n)$ converges to $\Xi_{\eta}(t, v)$ in $C^1(S; M)$. Moving on to the remaining conclusions of (b), the property~\eqref{eq:Xi-area-constant} follows from~\eqref{eq:Xi-definition} and the fact that $(H_0)^{-1} \circ H_t$ is a diffeomorphism from $S$ to itself. The first relation in~\eqref{eq:reparametrization-homotopy-endpoint} is also clear from~\eqref{eq:Xi-definition}. Next, upon letting $(\tau, u) = \Psi(\Upsilon_{\eta}(v))$ as before, we have by assumption that
\[
(\sigma, \widetilde{v}) = \Upsilon_{\eta}(v) = \Phi(\tau, u),
\]
so that, referring to the diagram \eqref{eq:Phi-diagram} illustrating the definition of $\Phi$, the map $\widetilde{v}$ must satisfy $\widetilde{v} \circ p_0 = u \circ \widetilde{f}$. Comparing this with \eqref{eq:Xi-alt-definition} (with $t = 1$), we get
\[
\widetilde{v} \circ p_0 = \Xi_{\eta}(1, v) \circ p_0,
\]
which gives the second relation in~\eqref{eq:reparametrization-homotopy-endpoint}.  We are done with part (b).

For part (c), take $v \in C^2(S; M)$ such that $E(\gamma_0, v) \leq \Lambda$, and define $\gamma = v^*g_M + \eta \gamma_0$ as before, with $\eta \in (0, 1)$ to be determined. By \eqref{eq:reparametrization-homotopy-endpoint} and \eqref{eq:Xi-definition} (with $t = 1$), we have 
\[
\Upsilon_{\eta}(v) = (\sigma, v \circ (H_0)^{-1}).
\]
Writing $F$ for $(H_0)^{-1}$, we note the standard fact that 
\begin{equation}\label{eq:A-E-lowerbound}
A(v) = A(v \circ F) \leq E(\sigma, v\circ F),
\end{equation}
which gives the first inequality in \eqref{eq:reparametrization-area-energy}. To obtain the second inequality, we recall that $p_0^*\gamma = \lambda |dz + \mu d\overline{z}|^2$, which implies that $(w^{\mu})^{-1}: (\HH, |dz|^2) \to (\HH, p_0^*\gamma)$ is a conformal map, and thus so is
\[
F:(S, \sigma) \to (S, \gamma).
\]
The conformal invariance of the Dirichlet energy then gives
\[
\begin{split}
E(\sigma, v\circ F) =\ & E(F^*\gamma, v\circ F) = E(\gamma, v) = \int_{S} \big(\frac{\tr_{\gamma}v^*g_M}{2}\big) \vol_{\gamma}.
\end{split}
\]
As $\gamma \geq v^*g_M \geq 0$ in the sense of symmetric $2$-tensors, we have $\tr_{\gamma}v^*g_M \leq \tr_{\gamma}\gamma = 2$, and hence 
\begin{equation}\label{eq:A-E-almost-done}
\begin{split}
E(\sigma, v \circ F) \leq\ & \int_{S} \vol_{\gamma} = \int_{S} \sqrt{\det\big[ g_M(dv(e_i), dv(e_j)) + \eta \delta_{ij} \big]} \vol_{\gamma_0},
\end{split}
\end{equation}
where $e_1, e_2$ is any oriented local $\gamma_0$-orthonormal frame on $S$. Expanding the determinant in the last integral yields
\[
\det\big[ g_M(dv(e_i), dv(e_j)) + \eta \delta_{ij} \big] = \det\big[ g_M(dv(e_i), dv(e_j))\big] + \eta |dv|_{\gamma_0}^2 + \eta^2.
\]
Substituting this into~\eqref{eq:A-E-almost-done} and using the elementary inequality $\sqrt{a + b} \leq \sqrt{a} + \sqrt{b}$, we get
\begin{equation}\label{eq:A-E-upperbound}
\begin{split}
E(\sigma, v\circ F) \leq\ & A(v) + \sqrt{\eta}\int_{S}\sqrt{ |dv|_{\gamma_0}^2 + \eta} \vol_{\gamma_0}\\
\leq\ & A(v) + \eta^{\frac{1}{2}} \cdot \big(\Vol(S, \gamma_0)\big)^{\frac{1}{2}} \cdot \Big( 2E(\gamma_0, v) + \eta \Vol(S, \gamma_0)\Big)^{\frac{1}{2}}\\
\leq\ & A(v) + \eta^{\frac{1}{2}}\cdot \big(\Vol(S, \gamma_0)\big)^{\frac{1}{2}} \cdot \Big( 2\Lambda + \Vol(S, \gamma_0) \Big)^{\frac{1}{2}},
\end{split}
\end{equation}
where to get the second line we used H\"older's inequality, and the last inequality is a consequence of our assumption that $E(\gamma_0, v) \leq \Lambda$, and the fact that $\eta \in (0, 1)$ to begin with. We finish the proof of~\eqref{eq:reparametrization-area-energy} upon taking $\eta$ to be sufficiently small depending on $\ep, \Lambda$, and $\Vol(S, \gamma_0)$, and recalling~\eqref{eq:A-E-lowerbound}.
\end{proof}

\section{Analytical tools (II): Harmonic replacement}\label{sec:harmonic-replacement}
The goal of this section is to define an energy decreasing process that operates on continuous families of pairs from the space $\cM'$ (recall \eqref{eq:cM'-definition}), which will be used in the proof of Theorem \ref{thm:min-max-existence} to ``tighten'' minimizing sequences in $[\bv_0]$ so as to extract min-max sequences with good compactness properties. We largely follow the steps in the foundational work of Colding and Minicozzi \cite{Colding-Minicozzi08b}, and also make use of subsequent refinements due to the second named author \cite{Zhou10,Zhou17b}. 

In Section \ref{subsec:local-considerations} we review harmonic replacement for small-energy maps defined on the unit disk $\bB$ in the plane. In Section \ref{subsec:admissible-subsets} we prove some more or less obvious facts about the action of Fuchsian groups on disjoint collections of disks in $\HH$. Based on these, in Section \ref{subsec:global-considerations} we extend the harmonic replacement operation first to pairs in the space $\cM$, and then to those in $\cM'$. In the same section we also discuss issues such as continuous dependence. Section \ref{subsec:iterated-replacement} is devoted to generalizing two crucial estimates from \cite{Colding-Minicozzi08b} on iterated harmonic replacements. These feed into the construction of the said energy decreasing process in Section \ref{subsec:energy-decreasing}, which proceeds by patching together local constructions, and thus iterated replacements occur on the overlap regions. The end results are Propositions \ref{prop:choice-of-disks} and \ref{prop:energy-decreasing}. 

\subsection{Local considerations}\label{subsec:local-considerations}
Recall that $\bB$ denotes the unit disk in $\RR^2$, and let $(M, g_M)$ be a closed Riemannian manifold isometrically embedded in some Euclidean space $\RR^q$. The following \textit{convexity estimate} is a deep result due to Colding and Minicozzi, based on earlier works by H\'elein and Wente, and is the main reason why harmonic replacements are well-defined and useful in min-max constructions. 
\begin{thm}[\cite{Colding-Minicozzi08b}, Theorem 3.1]
\label{thm:convexity}
There exists $\ep_0 > 0$ depending only on $M$ with the following property. For any weakly harmonic map $u \in W^{1, 2}(\bB; M)$ satisfying $\int_{\bB}|\nabla u|^2 \leq \ep_0$ and any $v \in W^{1, 2}(\bB; M)$ such that $v|_{\partial \bB} = u|_{\partial \bB}$ in the trace sense, there holds
\[
\frac{1}{4}\int_{\bB}|\nabla u - \nabla v|^2 \leq \frac{1}{2}\int_{\bB}|\nabla v|^2 - \frac{1}{2}\int_{\bB}|\nabla u|^2.
\]
\end{thm}
With $\ep_0$ as in Theorem~\ref{thm:convexity}, we consider
\begin{equation}\label{eq:small-energy-class}
\cA_{\ep_0} = \big\{u \in C^0(\overline{\bB}; M) \cap W^{1, 2}(\bB; M)\ |\ \int_{\bB}|\nabla u|^2 \leq \ep_0  \big\}.
\end{equation}
Theorem \ref{thm:replacement-existence-uniqueness} below is the existence and uniqueness result that underlies the definition, to be given shortly, of harmonic replacements. 
\begin{thm}[\cite{Colding-Minicozzi08b}, first half of Corollary 3.4]
\label{thm:replacement-existence-uniqueness}
Given $u \in \cA_{\ep_0}$, there exists a unique weakly harmonic map 
$h \in W^{1, 2}(\bB; M)$ such that $\int_{\bB}|\nabla h|^2 \leq \ep_0$, and that $h|_{\partial \bB} = u|_{\partial \bB}$ in the trace sense. The map $h$ minimizes the Dirichlet energy among $W^{1, 2}$-maps into $M$ that agree with $u$ on $\partial \bB$ in the trace sense. Moreover $h$ is continuous on $\overline{\bB}$ and smooth in $\bB$.
\end{thm}
The map $h$ given by Theorem~\ref{thm:replacement-existence-uniqueness} is the so-called \textit{harmonic replacement} of $u$, and we denote it by $\cR(u, \bB)$. The next theorem addresses the dependence of the harmonic replacement on the original map.
\begin{thm}[\cite{Colding-Minicozzi08b}, second half of Corollary 3.4]
\label{thm:replacement-continuity-1}
Suppose $(u_n)$ is a sequence in $\cA_{\ep_0}$ converging in $C^0(\overline{\bB}) \cap W^{1, 2}(\bB)$ to $u$. Then $\cR(u_n, \bB) \to \cR(u, \bB)$ in $C^0(\overline{\bB}) \cap W^{1, 2}(\bB)$. Moreover, given any $u, v \in \cA_{\ep_0}$ there holds the following estimate:
\begin{equation}\label{eq:replaced-energy}
\begin{split}
&\Big|\int_{\bB}|\nabla \cR(u, \bB)|^2 - \int_{\bB}|\nabla \cR(v, \bB)|^2\Big|\\
&\leq  C\|u - v\|_{\infty}\Big( \int_{\bB}|\nabla u|^2 + |\nabla v|^2\Big) + C\|\nabla u - \nabla v\|_2 \Big(\int_{\bB}|\nabla u|^2 + |\nabla v|^2 \Big)^{\frac{1}{2}}.
\end{split}
\end{equation}
where $C$ depends only on $M$.
\end{thm}
For the next result we slightly extend the definition of the harmonic replacement. Given an open set $\Omega \subset \RR^2$, suppose $u:\Omega \to M$ is a map of class $C^0 \cap W^{1, 2}$ locally, and that, for some disk $B$ with closure contained in $\Omega$, we have
\begin{equation}\label{eq:manifold-small-energy}
\int_{B} |\nabla u|^2 \leq \ep_0.
\end{equation}
Then $u$ can be harmonically replaced on $B$, and we denote the resulting map by $\cR(u, B)$; in other words,
\begin{equation}\label{eq:replacement-general}
\cR(u, B) := \left\{
\begin{array}{ll}
u, & \text{ on }\Omega \setminus B,\\
\cR(u \circ \varphi, \bB) \circ \varphi^{-1}, & \text{ on }B,
\end{array}
\right.
\end{equation}
where $\varphi$ is a conformal diffeomorphism from $\bB$ onto $B$ (any such $\varphi$ necessarily extends to a homeomorphism from $\overline{\bB}$ onto $\overline{B}$), the specific choice being irrelevant thanks to the uniqueness statement in Theorem~\ref{thm:replacement-existence-uniqueness}, and the fact that both the Dirichlet energy and the harmonicity condition are conformally invariant in dimension two. By the regularity statement in Theorem~\ref{thm:replacement-existence-uniqueness}, the new map $\cR(u, B)$ is still in $(C^0 \cap W^{1, 2})_{\loc}(\Omega; M)$. Also, since the domain of $u$, namely $\Omega$, is usually clear from the context, we omit it from the notation. Theorem~\ref{thm:replacement-continuity} below is an extension, due to the second named author, of Theorem \ref{thm:replacement-continuity-1}.
\begin{thm}[\cite{Zhou10}, Corollary 4.2]
\label{thm:replacement-continuity}
Suppose $(u_n)$ is a sequence in $\cA_{\ep_0}$ converging in $C^0(\overline{\bB}) \cap W^{1, 2}(\bB)$ to $u$. Suppose also that we are given closed disks $\overline{\bB_{R_{n}}(z_n)}$ and $\overline{\bB_{R}(z)}$ in $\bB$, with $R_n > 0$ for all $n$, such that $R_n \to R$ and $z_n \to z$ as $n \to \infty$. Then we have the following.
\vskip 1mm
\begin{enumerate}
\item[(a)] If $R > 0$, then
\[
\cR(u_n, \bB_{R_n}(z_n)) \to \cR(u, \bB_{R}(z)) \text{ in }C^0(\overline{\bB}) \cap W^{1, 2}(\bB).
\]
\vskip 1mm
\item[(b)] If $R = 0$, then
\[
\cR(u_n, \bB_{R_n}(z_n)) \to u \text{ in }C^0(\overline{\bB}) \cap W^{1, 2}(\bB).
\]
\end{enumerate}
\end{thm}
\begin{proof}
Some preparation is in order. We first extend $u_n$ and $u$ to all of $\RR^2$ by letting
\[
u_n(x) = u_n(|x|^{-2} x),\quad u(x) = u(|x|^{-2}x),\quad \text{for }x \in \RR^2 \setminus \bB.
\]
The extended maps, still denoted $u_n$ and $u$, lie in $(C^0 \cap W^{1, 2})_{\loc}(\RR^2; M)$, and moreover we have that $u_n$ converges to $u$ in $C^{0} \cap W^{1, 2}$ on compact subsets of $\RR^2$. Consequently, introducing the affine maps $\lambda_n(x) = z_n + R_n x$ and $\lambda(x) = z + Rx$, we have
\begin{equation}\label{eq:convergence-of-scaled-maps}
u_n \circ \lambda_n \to u \circ \lambda \quad \text{in }C^0 \cap W^{1, 2}, \text{ locally on }\RR^2.
\end{equation}
Indeed, given a compact set $K \subset \RR^2$, since $R_n \to R$ and $z_n \to z$, there exists $L > 0$ such that $\lambda_n(x), \lambda(x) \in \bB_{L}$ for all $n \in \NN$ and $x \in K$. From this we see that
\[
u_n \circ \lambda_n \to u \circ \lambda \quad \text{uniformly on }K.
\]
Next, given any $\ep > 0$, standard real analysis yields some $V \in C^0_{c}(\bB_{L + 1}; \RR^{2 \times q})$ such that $\|\nabla u - V\|_{2; \bB_{L + 1}} < \ep$. We then write
\[
\begin{split}
\nabla(u_n \circ \lambda_n)- \nabla(u \circ \lambda) =\ & R_n\cdot  (\nabla u_n \circ \lambda_n - \nabla u\circ \lambda_n) + R_n\cdot (\nabla u \circ \lambda_n - V \circ \lambda_n)\\
&+ (R_n - R) \cdot (V \circ \lambda_n) + R \cdot (V\circ \lambda_n - V \circ \lambda) \\
& + R \cdot (V \circ \lambda - \nabla u \circ \lambda).
\end{split}
\]
Using the area formula, the convergence of $\nabla u_n$ to $\nabla u$ in $L^2_{\loc}(\RR^2)$, and the continuity and boundedness of $V$, we infer from the above decomposition that 
\[
\limsup_{n \to \infty}\|\nabla(u_n \circ \lambda_n) - \nabla (u \circ \lambda)\|_{2; K} \leq 2\|\nabla u - V\|_{2; \bB_{L}} < 2\ep.
\]
Thus $\nabla (u_n \circ \lambda_n)\to \nabla (u \circ \lambda)$ in $L^2(K)$, and we have established \eqref{eq:convergence-of-scaled-maps}.

To prove part (a), we let $h_n = \cR(u_n, \bB_{R_n}(z_n))$ and $h = \cR(u, \bB_{R}(z))$, and notice by \eqref{eq:replacement-general} that
\begin{equation}\label{eq:scaled-maps}
h_n \circ \lambda_n = \left\{
\begin{array}{ll}
\cR(u_n \circ \lambda_{n}, \bB) &\text{ on }\bB,\\
u_n\circ \lambda_{n} &\text{ on }\RR^2 \setminus \bB,
\end{array}
\right.
\quad  h \circ \lambda = \left\{
\begin{array}{ll}
\cR(u \circ \lambda, \bB) &\text{ on }\bB,\\
u \circ \lambda&  \text{ on }\RR^2\setminus \bB.
\end{array}
\right.
\end{equation}
Both maps are of class $C^0 \cap W^{1, 2}$ on compact subsets of $\RR^2$. Moreover, by \eqref{eq:convergence-of-scaled-maps} and Theorem~\ref{thm:replacement-continuity-1}, we infer that
\begin{equation}\label{eq:scaled-replacement-convergence}
h_n \circ \lambda_n \to h \circ \lambda \quad \text{in }C^0 \cap W^{1, 2} \text{ locally on }\RR^2.
\end{equation}
Since $R > 0$, we can then think of $h_n$ and $h$ respectively as $(h_n \circ \lambda_n) \circ (\lambda_n)^{-1}$ and $(h \circ \lambda) \circ \lambda^{-1}$, and argue as in the proof of \eqref{eq:convergence-of-scaled-maps} to get that $h_n \to h$ in $C^0 \cap W^{1,2}$ on compact subsets of $\RR^2$. This proves (a).

For part (b), we again let $h_n = \cR(u_n, \bB_{R_n}(z_n))$, and note by \eqref{eq:convergence-of-scaled-maps} and Theorem~\ref{thm:replacement-continuity-1} that 
\begin{equation}\label{eq:scaled-replacement-b-prepare}
u_n \circ \lambda_n \to u(z), \quad h_n \circ \lambda_n \to u(z) \quad \text{in }(C^0 \cap W^{1, 2})_{\loc}(\RR^2),
\end{equation}
where $u(z)$ denotes the constant function taking that value. To deduce from this that $h_n \to u$ uniformly on $\overline{\bB}$, we argue by contradiction. Suppose that for some $\ep > 0$ there exists a subsequence of $h_n$, which we do not relabel, and a sequence of points $x_n \in \overline{\bB}$, such that 
\begin{equation}\label{eq:scaled-replacement-b-C0}
\big| h_n(x_n) - u(x_n)\big| \geq \ep, \quad \text{for all }n.
\end{equation}
If $|x_n - z_n| \leq R_n$ for all large enough $n$, then necessarily $x_n \to z$, and we get a contradiction to \eqref{eq:scaled-replacement-b-C0} upon noting that
\[
h_n(x_n) - u(x_n) = \big((h_n \circ\lambda_n)(\frac{x_n - z_n}{R_n}) - u(z)\big) - \big( u(x_n) - u(z) \big),
\]
and recalling that $h_n \circ \lambda_n$ converges to $u(z)$ on $\overline{\bB}$ uniformly. On the other hand, if along some further subsequence we have $|x_n - z_n| > R_n$ for all $n$, then since $h_n$ agrees with $u_n$ outside of $\bB_{R_n}(z_n)$, we have
\[
h_n(x_n) - u(x_n) = u_n(x_n) - u(x_n),
\]
which again leads to a contradiction with \eqref{eq:scaled-replacement-b-C0}, since $u_n \to u$ uniformly on $\overline{\bB}$. Thus we conclude that 
\[
h_n  \to u \text{ uniformly on }\overline{\bB}.
\]
It remains to establish that $\nabla h_n \to \nabla u$ in $L^2(\bB)$. Again using the fact that $h_n = u_n$ outside of $\bB_{R_n}(z_n)$, together with the conformal invariance of the energy, we find that
\[
\begin{split}
\int_{\bB}|\nabla u - \nabla h_n|^2 \leq\ & 2\int_{\bB} |\nabla u - \nabla u_n|^2 + 2\int_{\bB_{R_n}(z_n)}|\nabla u_n - \nabla h_n|^2\\
=\ &  2\int_{\bB} |\nabla u - \nabla u_n|^2  + 2 \int_{\bB} |\nabla (u_n \circ \lambda_n) - \nabla(h_n \circ \lambda_n)|^2.
\end{split}
\]
By \eqref{eq:scaled-replacement-b-prepare} and our assumption on $u_n$, both integrals on the second line tend to $0$, and we are done.
\end{proof}
\subsection{Admissible subsets of $\HH$}\label{subsec:admissible-subsets}
We adopt the setting of Section~\ref{subsec:mappings}. Thus, let the Riemann surface $\Sigma_0 = (S, [\gamma_0])$, the covering map $p_0: \HH \to \Sigma_0$, and the group $\Gamma_0$ of deck transformations be as in the start of Section~\ref{sec:conformal-reparametrization}. Also, with $\cH$ given by~\eqref{eq:marking-to-Fuchsian}, we recall from~\eqref{eq:Gamma-tau-definition} that
\[
\Gamma_\tau = \cH(\tau)(\Gamma_0),\quad \text{for }\tau \in T(\Sigma_0).
\]
\begin{defi}\label{defi:admissible-set}
Given $\tau \in T(\Sigma_0)$, a compact set $K \subset \HH$ is said to be \textit{$\tau$-admissible} if 
\[
\gamma(K) \cap K = \emptyset,\quad \text{for all }\gamma \in \Gamma_{\tau} \setminus \{\id\}.
\]
\end{defi}
Before making the next definition, we recall that whether a subset $B \subset \HH$ is a disk does not depend on the choice between $g_{\euc}$ and $g_{\hyp}$, whereas its radius and center do. Below, when rescaling a disk concentrically, we always do so using $g_{\hyp}$. In other words, if $\zeta$ and $s$ denote the (hyperbolic) center and radius of a closed disk $B$ in $\HH$, and if $\exp^{\HH}$ is the exponential map of $(\HH, g_{\hyp})$, then 
\begin{equation}\label{eq:hyp-scaling-exp}
B = \exp^{\HH}_{\zeta}\big(\overline{\bB_{s}(0)}\big),\quad \lambda B = \exp^{\HH}_{\zeta}\big(\overline{\bB_{\lambda s}(0)}\big).
\end{equation}
Also, given any $\gamma \in \Aut(\HH)$, we have
\begin{equation}\label{eq:scaling-commute-with-mobius}
\gamma(\lambda B) = \lambda \gamma(B).
\end{equation}
\begin{defi}\label{defi:admissible-disks-upstairs}
A finite collection $\sD$ of mutually disjoint closed disks in $\HH$ is said to be $\tau$-admissible if their union is so in the sense of Definition \ref{defi:admissible-set}, in which case we define
\begin{equation}\label{eq:orbited-disks-definition}
\Gamma_\tau(\sD)  = \{\gamma(B)\ |\ \gamma \in \Gamma_{\tau}, B \in \sD\},
\end{equation}
which is itself a collection of mutually disjoint closed disks. Given in addition some $\lambda > 0$, we define
\begin{equation}\label{eq:scaled-disks-definition}
\lambda \sD = \{\lambda B\ |\ B \in \sD\}.
\end{equation}
Note that if $\sD$ is $\tau$-admissible for some $\tau \in T(\Sigma_0)$, then so is $\lambda \sD$ as long as $\lambda \leq 1$.
\end{defi}

Lemmas \ref{lemm:admissible-set-limit}, \ref{lemm:collection-properties} and \ref{lemm:disk-relation} below contain a number of important, but perhaps obvious, properties of $\tau$-admissible sets. For use in their proofs, we need the following.
\begin{lemm}\label{lemm:proper-action}
Let $(\tau_n)$ be a sequence converging to some $\tau$ in $T(\Sigma_0)$. Suppose there exist a compact set $K_0 \subset \HH$ and elements $\gamma_n$ of $\Gamma_{\tau_n}$ such that
\begin{equation}\label{eq:pinned-by-compact}
\gamma_n(K_0) \cap K_0 \neq \emptyset,\quad \text{for all }n.
\end{equation}
Then, passing to a subsequence if necessary, there exists $\alpha \in \Gamma_0$ such that $\gamma_n = \cH(\tau_n)(\alpha)$ for all large enough $n$, and that $\gamma_n \to \cH(\tau)(\alpha)$ uniformly locally on $\HH$.
\end{lemm}
\begin{proof}
Take a biholomorphic map from $\HH$ to $\bB$, say the map $F$ in \eqref{eq:H-to-D}, and define $\psi_n = F\circ \gamma_n \circ F^{-1}$. By the assumption \eqref{eq:pinned-by-compact}, we have
\[
\psi_n(F(K_0)) \cap F(K_0) \neq \emptyset, \quad \text{for all }n.
\]
Since $F(K_0)$ is a compact subset of $\bB$, standard facts about $\Aut(\bB)$ yields a subsequence of $(\psi_n)$, which we do not relabel, that converges in $C^0_{\loc}(\bB) $ to some $\psi \in \Aut(\bB)$, and hence $\gamma_n$ converges in $C^0_{\loc}(\HH)$ to $F^{-1}\circ \psi \circ F$, which lies in $\Aut(\HH)$. Using the definition of $\Gamma_{\tau_n}$ to write 
\begin{equation}\label{eq:gamma-n-expression-1}
\gamma_n = \cH(\tau_n)(\alpha_n),
\end{equation}
for some $\alpha_n \in \Gamma_0$, and letting $(\widetilde{\varphi}_n)$ be the sequence of quasiconformal maps given by Lemma \ref{lemm:qc-Fuchsian-action}, we have 
\begin{equation}\label{eq:gamma-n-expression-2}
\cH(\tau)(\alpha_n) = (\widetilde{\varphi}_n)^{-1} \circ \gamma_n \circ \widetilde{\varphi}_n,
\end{equation}
which converges uniformly locally on $\HH$ to a limit that lies in $\Aut(\HH)$. Since $\cH(\tau)(\Gamma_0) = \Gamma_{\tau}$ is a discrete subgroup of $\Aut(\HH)$, we infer (see the proof of \cite[Lemma 2.16]{Imayoshi-Taniguchi}) that there exists $N \in \NN$ such that for all $n \geq N$ we have $\cH(\tau)(\alpha_n) = \cH(\tau)(\alpha_N)$, and hence $\alpha_n = \alpha_N$ since $\cH(\tau)$ is injective. Letting $\alpha = \alpha_N$, we get the first conclusion from  \eqref{eq:gamma-n-expression-1}, and the second conclusion follows from \eqref{eq:gamma-n-expression-2} combined with the fact that $\widetilde{\varphi}_n$ and $(\widetilde{\varphi}_n)^{-1}$ both converge to $\id$ locally uniformly on $\HH$.
\end{proof}

The following simple result guarantees, roughly speaking, that $\tau$-admissibility persists under perturbation of both the point $\tau$ and the compact set $K$. Note that, given in addition a sequence of compact sets $(K_n)$ in $\HH$, whether it converges to $K$ in the Hausdorff distance does not depend on the choice between $g_{\hyp}$ and $g_{\euc}$ in defining the distance metric on $\HH$.
\begin{lemm}\label{lemm:admissible-set-limit}
Suppose $(\tau_n)$ is a sequence converging to $\tau$ in $T(\Sigma_0)$, and $(K_n)$ is a sequence of compact subsets of $\HH$ converging to some compact set $K \subset \HH$ in the Hausdorff distance. Assume in addition that $K$ is $\tau$-admissible. Then $K_n$ is $\tau_n$-admissible for large enough $n$.
\end{lemm}
\begin{proof}
Suppose by contradiction that, up to taking a subsequence, we can find for each $n$ some $\gamma_n \in \Gamma_{\tau_n} \setminus \{\id\}$ and $z_n \in K_n$ such that $\gamma_n(z_n) \in K_n$. Since $K_n \to K$ in Hausdorff distance, passing to a subsequence again, there exist $z, w \in K$ such that
\begin{equation}\label{eq:pinned-pointwise}
z_n \to z,\quad \gamma_n(z_n) \to w.
\end{equation}
On the other hand, from the convergence of $K_n$ to $K$, we also get a compact set $K' \subset \HH$ that contains $K_n$ for all $n$. Thus \eqref{eq:pinned-by-compact} holds with $K_0 = K'$, and we may apply Lemma \ref{lemm:proper-action} to get some $\alpha \in \Gamma_0$ such that, taking a further subsequence if needed, we have
\[
\gamma_n \to \cH(\tau)(\alpha)=:\gamma, \quad \text{in }C^0_{\loc}(\HH).
\]
Combining this with \eqref{eq:pinned-pointwise} and the fact that $z_n \in K'$ for all $n$, we deduce that $\gamma(z) = w \in K$, and thus $\gamma(K) \cap K \neq \emptyset$. Since $\gamma \in \Gamma_{\tau}$ and $K$ is assumed to be $\tau$-admissible, we deduce that $\gamma =\id$, and hence $\alpha = \id$. This however implies, by Lemma \ref{lemm:proper-action} again, that eventually $\gamma_n = \id$, a contradiction, and we are done.
\end{proof}
Lemmas~\ref{lemm:collection-properties} and~\ref{lemm:disk-relation} below establish the basic properties of collections of the form $\Gamma_{\tau}(\sD)$, which are locations where harmonic replacement occurs, as we explain in the next section.
\begin{lemm}\label{lemm:collection-properties}
Suppose $\tau \in T(\Sigma_0)$ and that $\sD$ is a finite, disjoint, $\tau$-admissible collection of closed disks. Then we have the following.
\vskip 1mm
\begin{enumerate}
\item[(a)] The assignment below defines a bijection from $\Gamma_{\tau} \times \sD$ to $\Gamma_{\tau}(\sD)$:
\begin{equation}\label{eq:acted-collection}
(\gamma, B) \longmapsto \gamma(B).
\end{equation}
In fact, if $(\gamma, B), (\gamma', B') \in \Gamma_{\tau} \times \sD$ are such that $\gamma(B) \cap \gamma'(B') \neq \emptyset$, then $\gamma = \gamma'$ and $B = B'$.
\vskip 1mm
\item[(b)] Any compact subset $K \subset \HH$ intersects at most finitely many disks from $\Gamma_{\tau}(\sD)$.
\end{enumerate}
\end{lemm}
\begin{proof}
For part (a), we prove the last conclusion first. Suppose for some $(\gamma, B), (\gamma', B') \in \Gamma_{\tau} \times \sD$ we have $\gamma(B) \cap \gamma'(B') \neq \emptyset$. Then $\gamma(\cup_{D \in \sD} D)$ intersects $\gamma'(\cup_{D \in \sD} D)$, so the $\tau$-admissibility of $\cup_{D \in \sD} D$ implies that $\gamma = \gamma'$, and consequently $B \cap B' \neq \emptyset$. Since $\sD$ is a disjoint collection to begin with, this forces $B = B'$. Therefore $(\gamma, B) = (\gamma', B')$ as asserted. What we just proved shows that the map \eqref{eq:acted-collection} is injective. That it is surjective follows straight from the definition of $\Gamma_{\tau}(\sD)$.

For part (b), assume by contradiction that there exist a sequence of mutually distinct disks $B_n$ in $\Gamma_{\tau}(\sD)$ such that
\[
K \cap B_n \neq\emptyset\quad \text{for all }n.
\]
Since $\sD$ is finite, up to taking a subsequence we may assume that each $B_n$ is of the form $\gamma_n(B)$ for some $\gamma_n \in \Gamma_{\tau}$ and some $B \in \sD$. By Lemma \ref{lemm:proper-action} applied to the constant sequence $\tau$ and with $K_0 = K \cup B$, we get that eventually $\gamma_n$ all coincide with each other, which contradicts the assumption that the disks $B_n = \gamma_n(B)$ are mutually distinct.
\end{proof}

\begin{lemm}\label{lemm:disk-relation}
Given $\tau \in T(\Sigma_0)$, let $\sD$ be a finite, disjoint, $\tau$-admissible collection of disks, and let $(\tau_n)$ be a sequence in $T(\Sigma_0)$ converging to $\tau$. Then we have the following.
\vskip 1mm
\begin{enumerate}
\item[(a)] For all $(\gamma, B) \in \Gamma_{\tau} \times \sD$, there exists $N \in \NN$ such that for each $n \geq N$, exactly one disk in $\Gamma_{\tau_n}(\frac{1}{2}\sD)$ intersects $\gamma(B)$, and moreover this disk is contained in $\gamma(\frac{3}{4}B)$.
\vskip 1mm
\item[(b)] Given a compact set $K \subset \HH$, let $\cF = \{(\gamma, B) \in \Gamma_{\tau} \times \sD\ |\ \gamma(B) \cap K \neq \emptyset \}$. Then there exists $N \in \NN$ such that for all $n \geq N$, any member of $\Gamma_{\tau_n}(\frac{1}{2}\sD)$ that intersects $K$ is contained in $\gamma(\frac{3}{4}B)$ for some $(\gamma, B) \in \cF$.
\end{enumerate}
\end{lemm}
\begin{proof}
For part (a), given $(\gamma, B) \in \Gamma_{\tau} \times \sD$, by the definition of $\Gamma_{\tau}$ there exists a unique $\alpha \in\Gamma_0$ such that $\gamma = \cH(\tau)(\alpha)$. Thus, letting $\gamma_n = \cH(\tau_n)(\alpha)$, we see from Lemma \ref{lemm:qc-Fuchsian-action} that
\[
\gamma_n  \to \gamma,\quad \text{in }C^0_{\loc}(\HH),
\]
which yields some $N_1 \in \NN$ so that
\begin{equation}\label{eq:disk-relation-claim-1}
\gamma_n(\frac{1}{2}B) \subset \gamma(\frac{3}{4}B),\quad \text{for all }n \geq N_1.
\end{equation}
As $\gamma_n \in \Gamma_{\tau_n}$, this gives a member of $\Gamma_{\tau_n}(\frac{1}{2}\sD)$ which is contained in $\gamma(\frac{3}{4}B)$. 

It remains to prove that, increasing $N_1$ if necessary, any disk in $\Gamma_{\tau_n}(\frac{1}{2}\sD)$ that intersects $\gamma(B)$ must coincide with $\gamma_n(\frac{1}{2}B)$, as long as $n \geq N_1$. Assume by contradiction that along a subsequence of $(\tau_n)$ that we do not relabel, we can find $B_n \in \Gamma_{\tau_n}(\frac{1}{2}\sD)$ so that 
\begin{equation}\label{eq:disk-relation-contradiction}
B_n \neq \gamma_n(\frac{1}{2}B),\ \ B_n \cap \gamma(B) \neq\emptyset.
\end{equation}
Since $\sD$ is finite, again passing to a subsequence if needed, we may assume that each $B_n$ has the form $\delta_n(\frac{1}{2}B')$ for some $\delta_n \in \Gamma_{\tau_n}$ and some $B' \in \sD$, the latter not depending on $n$. Applying Lemma \ref{lemm:proper-action} with $K_0 = B' \cup \gamma(B)$, we obtain a further subsequence of $(\delta_n)$ that converges uniformly locally on $\HH$ to some $\delta \in \Gamma_{\tau}$, which satisfies by \eqref{eq:disk-relation-contradiction} that
\begin{equation}\label{eq:disk-relation-limit}
\delta(\frac{1}{2}B') \cap \gamma(B) \neq \emptyset.
\end{equation}
Moreover, still from Lemma \ref{lemm:proper-action}, we know that $\delta = \cH(\tau)(\beta)$ for some $\beta \in \Gamma_0$, and that eventually $\delta_n = \cH(\tau_n)(\beta)$. On the other hand, by the second conclusion of Lemma~\ref{lemm:collection-properties}(a), the condition \eqref{eq:disk-relation-limit} forces $\delta = \gamma$ and $B = B'$. Recalling that $\gamma = \cH(\tau)(\alpha)$, we conclude that $\alpha = \beta$, and hence 
\[
\delta_n = \cH(\tau_n)(\alpha) = \gamma_n,
\]
for all sufficiently large $n$. It follows that eventually $B_n = \delta_n(\frac{1}{2}B') = \gamma_n(\frac{1}{2}B)$, a contradiction to the first condition in~\eqref{eq:disk-relation-contradiction}. This proves part (a).

For part (b), assuming that the statement does not hold and using the finiteness of $\sD$, we get, after passing to a subsequence of $(\tau_n)$, some $B' \in \sD$ and $\delta_n \in \Gamma_{\tau_n}$ such that
\[
\delta_n(\frac{1}{2}B') \cap K \neq\emptyset,\quad \text{and that}\quad \delta_{n}(\frac{1}{2}B') \not\subset \gamma(\frac{3}{4}B) \text{ for all }(\gamma, B) \in \cF.
\]
Applying Lemma \ref{lemm:proper-action} as in the proof of (a), up to taking a further subsequence, we have that $(\delta_n)$ converges in $C^0_{\loc}(\HH)$ to some $\delta \in \Gamma_{\tau}$ which satisfies $\delta(B') \cap K \neq \emptyset$. In particular, $(\delta, B')\in \cF$, and we arrive at a contradiction upon noting that since $\delta_n \to \delta$ on compact subsets of $\HH$, eventually we have $\delta_n(\frac{1}{2}B') \subset \delta(\frac{3}{4}B')$.
\end{proof}
\subsection{Harmonic replacement on collections of disks}\label{subsec:global-considerations}
Given $(\tau, u)$ in the space $\cM$ defined by \eqref{eq:cM-definition} in Section~\ref{subsec:mappings}, and suppose $\sD$ is a finite, disjoint, $\tau$-admissible collection of closed disks in $\HH$ satisfying
\begin{equation}\label{eq:replacement-location} 
\int_{\bigcup_{B \in \sD}B}|\nabla u|^2 < \ep_0,
\end{equation}
where $\ep_0$ is the constant in Theorem \ref{thm:convexity}. By the conformal invariance of the energy, as well as the definition of $\cM$, for all $\gamma \in\Gamma_{\tau}$ we have
\[
\int_{\bigcup_{B \in \sD}\gamma(B)} |\nabla u|^2 = \int_{\bigcup_{B \in \sD}B} |\nabla u|^2 < \ep_0.
\]
It follows that $u$ can be harmonically replaced on each disk in the collection $\Gamma_{\tau}(\sD)$ defined in \eqref{eq:orbited-disks-definition}, and we denote the resulting map by $\cR(\tau, u, \sD)$. In other words, we set
\begin{equation}\label{eq:u-replacement}
\cR(\tau, u, \sD) =
\left\{
\begin{array}{ll}
u, & \text{ on }\HH \setminus \bigcup_{B \in \Gamma_{\tau}(\sD)}B,\\
\cR(u, B), & \text{ on each }B \in \Gamma_{\tau}(\sD).
\end{array}
\right.
\end{equation}
\begin{lemm}\label{lemm:replacement-invariant}
Given $(\tau, u) \in \cM$ and a finite, disjoint, $\tau$-admissible collection $\sD$ of disks so that \eqref{eq:replacement-location} holds, define $\cR(\tau, u, \sD)$ as above. Then $(\tau, \cR(\tau, u, \sD)) \in \cM$.
\end{lemm}
\begin{proof}
By the regularity statement in Theorem \ref{thm:replacement-existence-uniqueness}, the map $\cR(\tau, u, \sD)$ is $C^0 \cap W^{1, 2}$ locally on $\HH$. It remains to check that
\[
\cR(\tau, u, \sD)\circ \gamma = \cR(\tau, u, \sD) \text{ for all }\gamma \in \Gamma_{\tau}.
\]
To see that, take any $\gamma \in \Gamma_{\tau}$. For all $x$ which does not lie in $\cup_{B \in \Gamma_{\tau}(\sD)}B$, neither does $\gamma(x)$ since $\cup_{B \in \Gamma_{\tau}(\sD)}B$ is invariant under $\Gamma_{\tau}$. Recalling that $u \circ \gamma = u$, we get
\[
\cR(\tau, u, \sD)(\gamma(x)) = u(\gamma(x)) = u(x) = \cR(\tau, u, \sD)(x).
\]
On the other hand, given $B \in \Gamma_{\tau}(\sD)$, we take $\lambda$ to be any conformal diffeomorphism from $\bB$ onto $\Int(B)$ and recall the definition~\eqref{eq:replacement-general} to get
\[
\cR(\tau, u, \sD) = \cR(u \circ \gamma\circ \lambda, \bB) \circ (\gamma \circ \lambda)^{-1}\quad \text{on }\gamma(B).
\]
Thus, for all $x \in B$, evaluating the above at $\gamma(x)$ and using again that $u \circ \gamma = u$, we infer that
\[
\cR(\tau, u, \sD)(\gamma(x)) = \cR(u \circ \lambda, \bB)(\lambda^{-1}(x)) = \cR(u, B)(x) = \cR(\tau, u, \sD)(x).
\]
This completes the proof.
\end{proof}
Our next goal is to establish the continuous dependence of $(\tau, \cR(\tau, u, \sD))$ on $(\tau, u)$ as the latter varies in $\cM$. 
\begin{prop}\label{prop:replacement-H-continuity}
Suppose $(\tau_n, u_n)$ is a sequence in $\cM$ converging to some $(\tau, u) \in \cM$ and that $\sD$ is a finite, disjoint, $\tau$-admissible collection of closed disks in $\HH$ such that~\eqref{eq:replacement-location} holds. Assume also that $(r_n)$ is a sequence in $(0, \frac{1}{2}]$ converging to some $r$. 
\vskip 1mm
\begin{enumerate}
\item[(a)] If $r > 0$, then $(\tau_n, \cR(\tau_n, u_n, r_n\sD))$ converges to $(\tau, \cR(\tau, u, r\sD))$ in $\cM$.
\vskip 1mm
\item[(b)] If $r = 0$, then $(\tau_n, \cR(\tau_n, u_n, r_n\sD))$ converges to $(\tau, u)$ in $\cM$.
\end{enumerate}
\end{prop}
\begin{proof}
By the definition of the topology on $\cM$ we have $\tau_n \to \tau$ in $T(\Sigma_0)$. Thus in both parts it suffices to establish the $C^0 \cap W^{1, 2}$ convergence of $\cR(\tau_n, u_n, r_n\sD)$ on compact subsets of $\HH$.

For part (a), by Lemma~\ref{lemm:admissible-set-limit}, and the $W^{1, 2}$-convergence of $u_n$ to $u$ on the compact set $\cup_{B \in \sD}B$, we have for all sufficiently large $n$ that $\sD$ is $\tau_n$-admissible, and that \eqref{eq:replacement-location} holds with $u_n$ in place of $u$. Thus, eventually it makes sense to consider $\cR(\tau_n, u_n, r_{n}\sD)$. To save space, we define
\[
h_n = \cR(\tau_n, u_n, r_{n}\sD),\ \ h = \cR(\tau, u, r\sD).
\]
We first prove that $h_n \to h$ in $C^0 \cap W^{1, 2}$ on each disk in $\Gamma_{\tau}(\sD)$. Given such a disk $\gamma(B)$, where $\gamma\in \Gamma_{\tau}$ and $B \in \sD$, since $\gamma(rB)$ is the unique element of $\Gamma_{\tau}(r\sD)$ that intersects $\gamma(B)$, we have in terms of the notation \eqref{eq:replacement-general} that
\begin{equation}\label{eq:h-gamma-B}
h(x) = \cR(u, \gamma(rB))(x),\quad \text{for all }x \in \gamma(B). 
\end{equation}
On the other hand, Lemma~\ref{lemm:disk-relation} and its proof yields elements $\gamma_n \in \Gamma_{\tau_n}$ converging in $C^0_{\loc}(\HH)$ to $\gamma$, such that eventually
\vskip 1mm
\begin{enumerate}
\item[(i)] the disk $\gamma_n(\frac{1}{2}B)$ is contained in $\gamma(\frac{3}{4}B)$, and
\vskip 1mm
\item[(ii)] no other member of the collection $\Gamma_{\tau_n}(\frac{1}{2}\sD)$ intersects $\gamma(B)$. 
\end{enumerate}
Consequently, using in addition~\eqref{eq:scaling-commute-with-mobius}, we see that $\gamma_n(r_n B)$ is the only disk in $\Gamma_{\tau_n}(r_n \sD)$ that intersects $\gamma(B)$, so that, again in the notation \eqref{eq:replacement-general}, we have
\begin{equation}\label{eq:hn-gamma-B}
h_n(x) =\cR(u_n, \gamma_n(r_n B))(x), \quad \text{for all }x \in \gamma(B).
\end{equation}
Next, fix any conformal diffeomorphism $\lambda: \bB \to \Int(B)$. With the help of (i) above, there exist $R, R_n > 0$ and $z, z_n \in \bB$ such that
\begin{equation}\label{eq:parametrized-disks}
\overline{\bB_{R}(z)} = \lambda^{-1}(rB),\quad  \overline{\bB_{R_n}(z_n)} = (\lambda^{-1}\circ\gamma^{-1} \circ \gamma_n)(r_n B),
\end{equation}
and that both closed disks are subsets of $\bB$. From \eqref{eq:h-gamma-B} and \eqref{eq:hn-gamma-B}, we have for all $x \in \bB$ that
\[
(h\circ\gamma\circ \lambda)(x) = \cR(u \circ \gamma \circ \lambda, \bB_{R}(z))(x), \quad (h_n\circ\gamma\circ \lambda)(x) = \cR(u_n \circ \gamma \circ \lambda, \bB_{R_n}(z_n))(x).
\]
Writing $B=\exp^{\HH}_{\zeta}(\overline{\bB_\rho(0)})$ as in \eqref{eq:hyp-scaling-exp}, and noting that 
\[
(\gamma^{-1} \circ \gamma_n \circ \exp_{\zeta}^{\HH})(r_n \ \cdot\ ) \longrightarrow \exp_{\zeta}^{\HH}(r\ \cdot\ ),\quad \text{uniformly on }\overline{\bB_{\rho}(0)},
\]
we have $\overline{\bB_{R_n}(z_n)} \to \overline{\bB_{R}(z)}$ with respect to the Hausdorff distance. In particular 
\[
R_n \to R, \quad z_n \to z \quad \text{as }n \to \infty.
\]
Note also that $u_n \circ \gamma \circ \lambda \to u \circ \gamma \circ \lambda$ in $C^0(\overline{\bB}) \cap W^{1, 2}(\bB)$, and that eventually $u_n \circ \gamma \circ \lambda$ belongs to the class $\cA_{\ep_0}$ (see~\eqref{eq:small-energy-class}), because
\[
\lim_{n \to \infty}\int_{\bB}|\nabla (u_n \circ \gamma \circ \lambda)|^2 = \int_{\bB}|\nabla (u \circ \gamma \circ \lambda)|^2 =\int_{B}|\nabla (u \circ \gamma)|^2 = \int_{B} |\nabla u|^2 < \ep_0,
\]
where the last inequality and the equality that precedes it follow, respectively, from~\eqref{eq:replacement-location} and from the $\Gamma_{\tau}$-invariance of $u$. We can now apply Theorem~\ref{thm:replacement-continuity}(a), from which we get
\begin{equation}\label{eq:convergence-on-disk}
h_n \to h\quad \text{in }C^0 \cap W^{1, 2} \text{ on }\gamma(B).
\end{equation}
Now take any compact subset $K \subset \HH$ and define
\[
\cF : = \{(\gamma, B) \in \Gamma_{\tau} \times \sD\ |\ \gamma(B) \cap K \neq \emptyset \}.
\]
Then $h = u$ on $K \setminus \cup_{(\gamma, B) \in \cF}\,\gamma(\frac{3}{4}B)$, whereas Lemma~\ref{lemm:disk-relation}(b) and our assumption $r_n \leq \frac{1}{2}$ implies that $h_n = u_n$ on $K \setminus \cup_{(\gamma, B) \in \cF}\,\gamma(\frac{3}{4}B)$ for all sufficiently large $n$, and thus
\[
h_n \to u = h,\quad \text{in }C^0 \cap W^{1, 2}\big(K \setminus \cup_{(\gamma, B) \in \cF}\, \gamma(\frac{3}{4}B)\big).
\] 
Having already shown that $h_n \to h$ in $C^0 \cap W^{1, 2}$ on each disk in $\Gamma_{\tau}(\sD)$, and seeing that $\cF$ is finite by Lemma~\ref{lemm:collection-properties}, we conclude that
\[
h_n \to h \text{ in }C^0 \cap W^{1, 2}(K).
\]

For part (b), we essentially repeat the above argument. Given $\gamma \in \Gamma_{\tau}$ and $B \in \sD$, we replace the map $h$ in~\eqref{eq:h-gamma-B} by $u$ itself. In \eqref{eq:parametrized-disks} we instead set $R = 0$ and let $z$ be $\lambda^{-1}(\zeta)$, where $\zeta$ is the (hyperbolic) center of $B$, and then observe that 
\[
(\gamma^{-1} \circ \gamma_n \circ \exp_{\zeta}^{\HH})(r_n \ \cdot\ ) \longrightarrow \text{ constant map }\zeta,\quad \text{uniformly on }\overline{\bB_{\rho}(0)},
\]
so that $R_n \to 0$, while $z_n \to z$. Invoking Theorem~\ref{thm:replacement-continuity}(b) now leads to 
\[
h_n  \to u \quad \text{in }C^0 \cap W^{1, 2} \text{ on }\gamma(B).
\]
The remainder of the proof of (a) goes through without further change.
\end{proof}
We next discuss harmonic replacement for pairs in the space $\cM'$ defined by \eqref{eq:cM'-definition}. To set the stage, take $\sigma \in \Met^*_{-1}$, and let $\tau$ be the point in $T(\Sigma_0)$ represented by $((S, [\sigma]), \id)$. Then we have the following diagram, which has already appeared in Section~\ref{subsec:mappings}:
\begin{equation}\label{eq:coming-down-diagram-again}
\begin{tikzcd}
\mathbb{H} \arrow[r, "\widetilde{f}"] \arrow[d, "p_0"']    & \mathbb{H} \arrow[d, "p"]  \\
\Sigma_0 = (S, [\gamma_0]) \arrow[r, "\id"] & {(S, [\sigma])}                          
\end{tikzcd}
\end{equation}
Here $\widetilde{f}$ is the canonical lift of $\id: \Sigma_0 \to (S, [\sigma])$ with respect to $p_0$. Also, in the notation \eqref{eq:Gamma-tau-definition}, the deck transformation group of $p: \HH \to (S, [\sigma])$ coincides with $\Gamma_{\tau}$. Given a closed geodesic disk $B$ in $(S, \sigma)$, by a \emph{lift of $B$ with respect to $p$} we mean a closed disk $\widetilde{B} \subset \HH$ such that $p(\widetilde{B}) = B$ and that $p$ restricts to an isometry with respect to $g_{\hyp}$ on a neighborhood of $\widetilde{B}$. Since $p: (\HH, g_{\hyp}) \to (S, \sigma)$ is a Riemannian covering map, the set of lifts of $B$ is non-empty, and form a disjoint collection permuted by the action of $\Gamma_{\tau}$.

\begin{rmk}\label{rmk:upstairs-and-downstairs}
The following facts, which are straightforward to verify, link the current discussion to the earlier one concerning $\cM$. 
\vskip 1mm
\begin{enumerate}
\item[(1)]  If $\fB$ is a finite, disjoint collection of closed geodesic disks in $(S, \sigma)$, and to each $B \in \fB$ we assign a choice of lift $\widetilde{B}$ with respect to $p$, then $\{\widetilde{B}\ |\ B \in \fB\}$ is a disjoint collection of closed disks which is $\tau$-admissible in the sense of Definition \ref{defi:admissible-disks-upstairs}.
\vskip 1mm
\item[(2)] Starting with a finite, disjoint, $\tau$-admissible collection $\sD$ of closed disks in $\HH$, then Lemma \ref{lemm:admissible-set-limit} yields some $\delta > 0$ such that $p$ is injective on $\bigcup_{\widetilde{B} \in \sD}(1+\delta)\widetilde{B}$. Consequently $\{p(\widetilde{B})\ |\ \widetilde{B} \in \sD\}$ is a disjoint collection of closed geodesic disks in $(S, \sigma)$, and each $\widetilde{B}$ is a lift of its projection $p(\widetilde{B})$.
\end{enumerate}
\end{rmk}
One final piece of notation (for use in Corollary \ref{coro:replacement-downstairs-continuity} below, among other places): given a finite, disjoint collection $\fB$ of geodesic disks in $(S, \sigma)$ as well as some $\lambda \in (0, 1]$, we set
\[
\lambda \fB = \{\lambda B\ |\ B \in \fB \},
\]
where of course $\lambda B$ means the geodesic disk in $(S, \sigma)$ that is concentric with $B$ and has $\lambda$ times its radius. Notice that, with respect to $p$, if $\widetilde{B}$ is a lift of $B$, then $\lambda \widetilde{B}$ is a lift of $\lambda B$.

\begin{defi}\label{defi:v-replacement}
Let $\ep_0$ be the energy threshold from Theorem \ref{thm:convexity}. 
Given $(\sigma, v) \in \cM'$ and a finite, disjoint collection $\fB$ of closed geodesic disks in $(S, \sigma)$ such that
\begin{equation}\label{eq:replacement-location-downstairs}
\sum_{B \in \fB}\int_{B} |\nabla v|_{\sigma}^2 \vol_{\sigma} < \ep_0,
\end{equation} 
we let $(\tau ,u) = \Psi(\sigma, v)$, so that $u = v \circ p$ in the notation of \eqref{eq:Phi-diagram}, and choose for each $B \in \fB$ a lift $\widetilde{B}$ with respect to $p$. Then Remark~\ref{rmk:upstairs-and-downstairs} shows that the collection
\[
\sD := \{\widetilde{B}\ |\ B \in \fB\}
\]
of disks is $\tau$-admissible, and by~\eqref{eq:replacement-location-downstairs} we see that $u$ satisfies \eqref{eq:replacement-location} with this choice of $\sD$. As a result $\cR(\tau, u, \sD)$ is defined, which together with $\tau$ form an element of $\cM$ by Lemma~\ref{lemm:replacement-invariant}. We then define $\cR(\sigma, v, \fB)$ by the following relation:
\begin{equation}\label{eq:harmonic-replacement-downstairs}
(\sigma, \cR(\sigma, v, \fB)) = \Phi(\tau, \cR(\tau, u, \sD)).
\end{equation}
\end{defi}
We make some remarks about this definition before proceeding.
\begin{rmk}\label{rmk:equivalent-lift}
Since any two lifts of the same geodesic disk in $(S, \sigma)$ are related by an action of $\Gamma_{\tau}$, it follows that in the above definition, the collection $\Gamma_{\tau}(\sD)$ (see~\eqref{eq:orbited-disks-definition}), and hence the map $\cR(\tau, u, \sD)$, does not depend on the choice of lifts, even though $\sD$ does. Consequently $\cR(\sigma, v, \fB)$ is well-defined.
\end{rmk}
\begin{rmk}\label{rmk:v-replacement}
In view of \eqref{eq:harmonic-replacement-downstairs} and the second relation in \eqref{eq:Phi-u-v-relation}, we have $\cR(\sigma, v, \fB) \circ p = \cR(\tau, u, \sD)$. Consequently, 
\begin{equation}\label{eq:replacement-up-down-relation}
\cR(\sigma, v, \fB)(x) = \left\{
\begin{array}{ll}
v(x), &  \text{ if }x \not \in \cup_{B \in \fB}B, \\
\cR(u, \widetilde{B}) \circ (p|_{\widetilde{B}})^{-1}(x), & \text{ if }x \in B \text{ for some }B \in \fB,
\end{array}
\right.
\end{equation}
where in the second case, $\widetilde{B}$ can be replaced by any other lift of $B$. As $p: (\widetilde{B}, g_{\euc}) \to (B, \sigma)$ is a conformal diffeomorphism for each $B \in \fB$, we obtain from Theorems \ref{thm:convexity} and \ref{thm:replacement-existence-uniqueness} the following.
\vskip 1mm
\begin{enumerate}
\item[(1)] $\cR(\sigma, v, \fB)\big|_{B}$ is the unique harmonic map from $(B, \sigma)$ to $M$ which has $E(\sigma,\cdot)$ energy at most $\ep_0$ and agrees with $v$ on $\partial B$. Also it has the least $E(\sigma, \cdot)$ energy among $W^{1, 2}$-maps from $B$ into $M$ that agrees with $v$ on $\partial B$.
\vskip 1mm
\item[(2)] By the convexity estimate of Theorem \ref{thm:convexity} we have
\[
\begin{split}
\frac{1}{4}\int_{S}|\nabla v - \nabla \cR(\sigma, v, \fB)|^2_{\sigma}\vol_{\sigma} =\ & \frac{1}{4}\sum_{B\in \fB}\int_{\widetilde{B}} |\nabla u - \nabla \cR(u, \widetilde{B})|^2\\
\leq\ &\frac{1}{2}\sum_{B \in \fB}\int_{\widetilde{B}} |\nabla u|^2 - |\nabla \cR(u, \widetilde{B})|^2\\
= \ & E(\sigma, v) - E(\sigma, \cR(\sigma, v, \fB)).
\end{split}
\]
\end{enumerate}
\end{rmk}

We end this section with the counterpart of Proposition~\ref{prop:replacement-H-continuity} in the context of $\cM'$. Below we suppose that $(\sigma_n, v_n)$ is a sequence in $\cM'$ converging to some $(\sigma, v) \in \cM'$, and adopt the notation in the diagrams~\eqref{eq:Phi-diagram} and~\eqref{eq:n-Phi-diagram}. In particular we write $(\tau_n, u_n) : = \Psi(\sigma_n, v_n)$ and $(\tau, u) : = \Psi(\sigma, v)$.
\begin{coro}\label{coro:replacement-downstairs-continuity}
In the above setting, let $\fB$ be a finite, disjoint collection of geodesic disks in $(S, \sigma)$ such that~\eqref{eq:replacement-location-downstairs} holds, and suppose that for each $B \in \fB$, a lift $\widetilde{B}$ with respect to $p$ is chosen. Also, let $(r_n)$ be a sequence in $(0, \frac{1}{2}]$ converging to some $r$. Then, for all large enough $n$, 
\[
\fB_n := \{p_n(\widetilde{B})\ |\ B \in \fB\}
\]
is a finite, disjoint collection of closed geodesic disks in $(S, \sigma_n)$. Moreover, $\cR(\sigma_n, v_n, r_n\fB_n)$ is defined, and we have:
\vskip 1mm
\begin{enumerate}
\item[(a)] If $r > 0$, then $(\sigma_n, \cR(\sigma_n, v_n, r_n\fB_n)) \to (\sigma, \cR(\sigma, v, r\fB))$ in $\cM'$.
\vskip 1mm
\item[(b)] If $r  = 0$, then $(\sigma_n, \cR(\sigma_n, v_n, r_n\fB_n)) \to (\sigma, v)$ in $\cM'$.
\end{enumerate}
\end{coro}
\begin{proof}
By Remark~\ref{rmk:upstairs-and-downstairs}, we see that the collection $\sD := \{\widetilde{B}\ |\ B \in \fB\}$ of disks in $\HH$ is $\tau$-admissible. Moreover, the smallness condition~\eqref{eq:replacement-location} holds with respect to $u = v \circ p$. Since $\tau_n \to \tau$ in $T(\Sigma_0)$ thanks to Proposition~\ref{prop:coming-down}, we see from Lemma~\ref{lemm:admissible-set-limit} and Remark~\ref{rmk:upstairs-and-downstairs} that eventually $\sD$ is $\tau_n$-admissible, and that indeed $\fB_n$ is a disjoint collection of geodesic disks in $(S, \sigma_n)$. On the other hand, Proposition \ref{prop:coming-down} also gives $u_n \to u$ in $C^0\cap W^{1,2}$ locally on $\HH$, and hence, by \eqref{eq:replacement-location} and the fact that $u_n = v_n \circ p_n$, we have for large enough $n$ that
\[
\sum_{\widetilde{B} \in \sD} \int_{p_n(\widetilde{B})}|\nabla v_n|_{\sigma_n}^2 \vol_{\sigma_n} = \sum_{\widetilde{B} \in \sD} \int_{\widetilde{B}} |\nabla u_n|^2 < \ep_0.
\]
From this energy bound and the fact that $r_n \fB_n = \{p_n(r_n \widetilde{B})\ |\ \widetilde{B} \in \sD\}$, we see that eventually $\cR(\sigma_n, v_n, r_n\fB_{n})$ is defined, and that
\[
(\sigma_n, \cR(\sigma_n, v_n, r_n \fB_n)) = \Phi(\tau_n, \cR(\tau_n, u_n, r_n \sD)).
\]
In the case $r > 0$, we similarly have
\[
(\sigma, \cR(\sigma, v, r \fB)) = \Phi(\tau, \cR(\tau, u, r \sD)).
\]
Both assertions are now direct consequences of Proposition~\ref{prop:coming-down} and Proposition~\ref{prop:replacement-H-continuity}. 
\end{proof}
\subsection{Iterated harmonic replacements}\label{subsec:iterated-replacement}
In this section we recall the generalizations obtained by the second named author in~\cite{Zhou17b} of the estimates due to Colding and Minicozzi~\cite{Colding-Minicozzi08b} regarding iterated harmonic replacements. These are stated as Proposition \ref{prop:iterated-replacement-estimates}. The parts of the argument whose adaptation to our setting require some effort are isolated as Proposition \ref{prop:iterated-replacement}, the proof of which takes up most of the length of this section.

In addition to $(\HH, g_{\hyp})$ we shall have occasion to work with the Poincar\'e disk $(\bB, g_{-1}: = \frac{4}{(1 - |z|^2)^2}g_{\euc})$. Standard facts imply that there exists $\rho_0 \in (0, \frac{1}{4})$ such that 
\begin{equation}\label{eq:radius-comparable}
\bB_{\frac{r}{3}}(0) \subset B_{g_{-1}}(0, r) \subset \bB_{\frac{r}{2}}(0), \text{ for all }r \in [0, \rho_0],
\end{equation}
where $B_{g_{-1}}(0, r)$ denotes the open geodesic disk with respect to $g_{-1}$ having the indicated center and radius. Also, letting $F:\HH \to \bB$ be the bihomorphic map defined by \eqref{eq:H-to-D} and writing $G$ for $F^{-1}$, it is well-known that 
\begin{equation}\label{eq:B-to-H-isometry}
G:(\bB, g_{-1}) \to (\HH, g_{\hyp})\quad \text{is an isometry.}
\end{equation}

\begin{prop}\label{prop:iterated-replacement}
Suppose $(\sigma, v_1), (\sigma, v_2) \in \cM'$ and that $\fB$ is a finite, disjoint collection of closed geodesic disks in $(S, \sigma)$, each with radius at most $\rho_0$, and that 
\begin{equation}\label{eq:iterated-replacement-small-energy}
\sum_{B\in \fB} \int_{B} |\nabla v_1|_{\sigma}^2 + |\nabla v_2|^2_{\sigma} \vol_{\sigma} < \ep_0.
\end{equation}
Assume further that $v_1$ and $v_2$ agree somewhere on $\partial(\mu B)$ for all $B \in \fB$ and $\mu \in [\frac{1}{2}, 1]$. Then we have 
\begin{equation}\label{eq:iterated-replacement}
\begin{split}
&E(\sigma, v_2) - E(\sigma, \cR(\sigma, v_2, \frac{1}{2}\fB))\\
&\leq E(\sigma, v_1) - E(\sigma, \cR(\sigma, v_1, \fB)) + \frac{1}{\kappa} \Big(\sum_{B \in \fB} \int_{B} |\nabla v_1 - \nabla v_2|_{\sigma}^2 \vol_{\sigma}\Big)^{\frac{1}{2}},
\end{split}
\end{equation}
where the constant $\kappa$ depends only on $M$. 
\end{prop}
Before giving the proof of Proposition \ref{prop:iterated-replacement} we recall the following construction from~\cite{Colding-Minicozzi08b}.
\begin{lemm}[\cite{Colding-Minicozzi08b}, Lemma 3.11]
\label{lemm:interpolation}
There exists $\tau > 0$ depending on $M$ with the following property. Suppose $f, g$ are maps in $C^0 \cap W^{1, 2}(\partial \bB; M)$ that agree somewhere on $\partial \bB$ and satisfy
\[
\|f' - g'\|_{2; \partial \bB} \leq \tau.
\]
Then there exist some $\rho \in (0, \frac{1}{2}]$ and an interpolating map $w \in C^0 \cap W^{1, 2}(\bB \setminus \bB_{1 - \rho}; M)$ such that 
\[
w((1-\rho)x) = f(x) , \quad w(x) = g(x), \quad \text{for all }x \in \partial \bB,
\]
and that
\[
\int_{\bB \setminus \bB_{1 - \rho}}|\nabla w|^2 \leq C\big(\|f'\|_{2; \partial \bB} + \|g'\|_{2; \partial \bB}\big)\cdot\|f' - g'\|_{2; \partial \bB},
\]
where $C$ is a numerical constant.
\end{lemm}
\begin{proof}[Proof of Proposition~\ref{prop:iterated-replacement}]
Some preliminary remarks are in order. Define $(\tau, u_i) = \Psi(\sigma, v_i)$ for $i = 1, 2$, so that in particular $\tau = [(S, [\sigma]), \id]$. As before, we write $\widetilde{f}$ for the canonical lift of $\id:(S, [\gamma_0]) \to (S, [\sigma])$ with respect to $p_0$, and let $p$ be the covering map such that $p \circ \widetilde{f} = p_0$. We then have the following diagram, which is basically a repeat of \eqref{eq:Phi-diagram}:
\begin{equation}\label{eq:replacement-lift}
\begin{tikzcd}
\mathbb{H} \arrow[r, "\widetilde{f}"] \arrow[d, "p_0"']    & \mathbb{H} \arrow[d, "p"] \arrow[r, "u_i"] & M \\
\Sigma_0 = (S, [\gamma_0]) \arrow[r, "\id"] 
& {(S, [\sigma])}     \arrow[ru, "v_i"']  &
\end{tikzcd}
\end{equation}
Fix, for each $B \in \fB$, a lift $\widetilde{B}$ with respect to $p$. Noting that $\partial (\mu B) = p(\partial(\mu \widetilde{B}))$, we see that
\begin{equation}\label{eq:u1-u2-somewhere}
u_1 = u_2 \text{ somewhere on }\partial(\mu\widetilde{B}), \text{ for all }B \in \fB, \ \mu \in [\frac{1}{2}, 1].
\end{equation}
Next, from the fact that $p:(\widetilde{B}, g_{\euc}) \to (B, \sigma)$ is a conformal diffeomorphism, we have
\begin{equation}\label{eq:W12-difference-upstairs}
\int_{B} |\nabla v_1 - \nabla v_2|_{\sigma}^2 \vol_{\sigma} = \int_{\widetilde{B}} |\nabla u_1 - \nabla u_2|^2.
\end{equation}
Taking into account also~\eqref{eq:replacement-up-down-relation}, we get
\begin{equation}\label{eq:energy-drop-upstairs}
E(\sigma, v_1) - E(\sigma, \cR(\sigma, v_1, \fB)) = \frac{1}{2}\sum_{B \in \fB} \big[\int_{\widetilde{B}}|\nabla u_1|^2 - \int_{\widetilde{B}}|\nabla \cR(u_1, \widetilde{B})|^2\big],
\end{equation}
and a similar expression holds for $E(\sigma, v_2) - E(\sigma, \cR(\sigma, v_2, \frac{1}{2}\fB))$, which in particular implies 
\[
E(\sigma, v_2) - E(\sigma, \cR(\sigma, v_2, \frac{1}{2}\fB)) \leq \frac{1}{2}\sum_{B \in \fB}\int_{B}|\nabla v_2|_{\sigma}^2 \vol_{\sigma} < \frac{\ep_0}{2},
\]
and thus it suffices to prove~\eqref{eq:iterated-replacement} under the extra assumption 
\begin{equation}\label{eq:extra-smallness}
\sum_{B \in \fB} \int_{B} |\nabla v_1 - \nabla v_2|_{\sigma}^2 \vol_{\sigma} < \frac{\tau^2}{A^2},
\end{equation}
where $\tau$ is the threshold given by Lemma~\ref{lemm:interpolation}, and $A$ is a universal constant to be determined.

To that end, fix any $B \in \fB$, and choose $\gamma \in \Aut(\HH)$ such that the center of $\gamma^{-1}(\widetilde{B})$ with respect to $g_{\hyp}$ is at $i \in \HH$, in which case there exists some some $R \leq \rho_0$ such that 
\[
(\gamma \circ G)^{-1}(\widetilde{B}) = \overline{B_{g_{-1}}(0, R)}.
\]
Next, define 
\[
\widehat{u}_{i} = u_{i} \circ \gamma \circ G,\quad \text{for }i = 1, 2.
\]
Noting that $\gamma \circ G: (\bB, g_{\euc}) \to (\HH, g_{\euc})$ is a conformal diffeomorphism, and also using \eqref{eq:radius-comparable} and the co-area formula, we get some $r \in [\frac{R}{4}, \frac{R}{3}]$ such that $\widehat{u}_1$ and $\widehat{u}_2$ are of class $C^0 \cap W^{1, 2}$ when restricted to $\partial \bB_{r}$, and that 
\begin{equation}\label{eq:coarea-bound-1}
r\int_{\partial \bB_{r}} |\nabla \widehat{u}_1 - \nabla \widehat{u}_2|^2 \leq 10\int_{\bB_{\frac{R}{3}}} |\nabla \widehat{u}_1 - \nabla \widehat{u}_2|^2 \leq 10\int_{\widetilde{B}}|\nabla u_1 - \nabla u_2|^2,
\end{equation}
\begin{equation}\label{eq:coarea-bound-2}
r\int_{\partial \bB_{r}} |\nabla \widehat{u}_1|^2 + |\nabla\widehat{u}_2|^2 \leq 10\int_{\bB_{\frac{R}{3}}} |\nabla \widehat{u}_1|^2 + |\nabla \widehat{u}_2|^2 \leq 10\int_{\widetilde{B}}|\nabla u_1|^2 + |\nabla u_2|^2.
\end{equation}
We next want to apply Lemma~\ref{lemm:interpolation} to $\widehat{u}_1$ and $\widehat{u}_2$. For that, observe that since $r \in [\frac{R}{4}, \frac{R}{3}]$, we have by~\eqref{eq:radius-comparable} and~\eqref{eq:B-to-H-isometry} that, for some $\mu \in [\frac{1}{2}, 1]$, 
\begin{equation}\label{eq:good-disk-mu}
\overline{\bB_{r}} = \overline{B_{g_{-1}}(0, \mu R)} = (\gamma \circ G)^{-1}(\mu \widetilde{B}).
\end{equation}
This together with \eqref{eq:u1-u2-somewhere} shows that 
\[
\widehat{u}_1 = \widehat{u}_2,\quad \text {somewhere on }\partial \bB_{r}.
\]
Recalling also \eqref{eq:coarea-bound-1}, \eqref{eq:W12-difference-upstairs}, and \eqref{eq:extra-smallness}, we see that upon requiring, say, $A > 10$, we can apply Lemma~\ref{lemm:interpolation} (suitably scaled) to $\widehat{u}_1$ and $\widehat{u}_2$, yielding some $\rho \in (0, \frac{1}{2}]$ and $w \in C^0 \cap W^{1, 2}(\bB_r \setminus \bB_{(1 - \rho)r}; M)$ such that 
\[
w((1 - \rho)x) = \widehat{u}_2(x),\quad w(x) = \widehat{u}_1(x), \quad \text{for all }x \in \partial \bB_{r},
\]
and that
\begin{equation}\label{eq:interpolation-in-proof}
\int_{\bB_{r} \setminus \bB_{(1 - \rho)r}} |\nabla w|^2 \leq C\Big( \int_{\widetilde{B}} |\nabla u_1|^2 + |\nabla u_2|^2 \Big)^{\frac{1}{2}} \cdot \Big( \int_{\widetilde{B}} |\nabla u_1 - \nabla u_2|^2 \Big)^{\frac{1}{2}},
\end{equation}
where in getting \eqref{eq:interpolation-in-proof} we also used \eqref{eq:coarea-bound-1} and \eqref{eq:coarea-bound-2}. To continue, we let
\[
\widehat{v}(x) = \left\{
\begin{array}{ll}
\widehat{u}_1(x), & \text{ if }x \in B_{g_{-1}}(0, R) \setminus \bB_{r},\\
w(x), & \text{ if }x \in \bB_{r} \setminus \bB_{(1 - \rho)r},\\
\cR(u_2, \mu \widetilde{B})\circ \gamma\circ G\big(\frac{x}{1 -\rho}\big), & \text{ if }x \in \bB_{(1 - \rho)r}.
\end{array}
\right.
\]
With the help of \eqref{eq:good-disk-mu} and the boundary behavior of $w$, it is not hard to see that $\widehat{v}$ belongs to $C^0 \cap W^{1, 2}$ on $B_{g_{-1}}(0, R)$, so that 
\[
\widehat{v} \circ (\gamma \circ G)^{-1} \in C^0 \cap W^{1,2}(\widetilde{B}; M)
\]
Since also $\widehat{v} \circ (\gamma \circ G)^{-1}  = u_1$ on $\partial \widetilde{B}$, it follows from the energy minimizing property of $\cR(u_1, \widetilde{B})$ that
\[
\begin{split}
\int_{\widetilde{B}} |\nabla \cR(u_1, \widetilde{B})|^2 \leq\ & \int_{\widetilde{B}}|\nabla (\widehat{v} \circ G^{-1} \circ \gamma^{-1})|^2 =  \int_{B_{g_{-1}}(0, R)} |\nabla \widehat{v}|^2\\
=\ & \int_{B_{g_{-1}}(0, R)\setminus \bB_{r}} |\nabla \widehat{u}_1|^2 + \int_{\bB_{r}} |\nabla (\cR(u_2, \mu\widetilde{B}) \circ \gamma \circ G)|^2 + \int_{\bB_{r}\setminus \bB_{(1-\rho)r}}|\nabla w|^2\\
\leq\ & \int_{\widetilde{B}} |\nabla u_1|^2 - \int_{\mu\widetilde{B}}|\nabla u_1|^2 + \int_{\mu\widetilde{B}} |\nabla \cR(u_2, \mu \widetilde{B})|^2 \\
&+ C\Big( \int_{\widetilde{B}} |\nabla u_1|^2 + |\nabla u_2|^2 \Big)^{\frac{1}{2}} \cdot \Big( \int_{\widetilde{B}} |\nabla u_1 - \nabla u_2|^2 \Big)^{\frac{1}{2}},
\end{split}
\]
where for the last inequality we used the estimate \eqref{eq:interpolation-in-proof}. By a similar but much simpler comparison argument, this time using the energy minimizing property of $\cR(u_2, \mu \widetilde{B})$, we get
\[
\begin{split}
\int_{\mu \widetilde{B}} |\nabla \cR(u_2, \mu \widetilde{B})|^2 \leq\ & 
\int_{\mu\widetilde{B}} |\nabla u_2|^2 - \int_{\frac{1}{2}\widetilde{B}} |\nabla u_2|^2 + \int_{\frac{1}{2}\widetilde{B}}|\nabla \cR(u_2, \frac{1}{2}\widetilde{B})|^2.
\end{split}
\]
Note also the following consequence of H\"older's inequality:
\[
\int_{\mu\widetilde{B}}\big| |\nabla u_1|^2 - |\nabla u_2|^2 \big| \leq 2\Big( \int_{\widetilde{B}} |\nabla u_1|^2 + |\nabla u_2|^2 \Big)^{\frac{1}{2}} \cdot \Big( \int_{\widetilde{B}} |\nabla u_1 - \nabla u_2|^2 \Big)^{\frac{1}{2}}.
\]
Combining the previous three estimates and rearranging leads to
\[
\begin{split}
&\int_{\frac{1}{2}\widetilde{B}} |\nabla u_2|^2 - \int_{\frac{1}{2}\widetilde{B}} |\nabla \cR(u_2, \frac{1}{2}\widetilde{B})|^2\\
\leq\ & \int_{\widetilde{B}}|\nabla u_1|^2 - \int_{\widetilde{B}}|\nabla \cR(u_1, \widetilde{B})|^2 + C\Big( \int_{\widetilde{B}} |\nabla u_1|^2 + |\nabla u_2|^2 \Big)^{\frac{1}{2}} \cdot \Big( \int_{\widetilde{B}} |\nabla u_1 - \nabla u_2|^2 \Big)^{\frac{1}{2}}.
\end{split}
\]
Summing over $B \in \fB$ and using~\eqref{eq:energy-drop-upstairs} and its analogue for $E(\sigma, v_2) - E(\sigma, \cR(\sigma, v_2, \frac{1}{2}\fB))$, while also applying the Cauchy--Schwarz inequality to the contribution from the product term in the above estimate, and recalling~\eqref{eq:iterated-replacement-small-energy} and~\eqref{eq:W12-difference-upstairs}, we obtain~\eqref{eq:iterated-replacement} as asserted.
\end{proof}

Another preparatory result we need is a simple energy estimate which leads to a condition (see Remark~\ref{rmk:many-replacements}) under which harmonic replacement can be performed successively over multiple collections.
\begin{lemm}\label{lemm:energy-after-replacement}
Let $(\sigma, v) \in \cM'$ and suppose $\fB$ is a finite, disjoint collection of closed geodesic disks in $(S, \sigma)$ such that each has radius at most $\rho_0$ and that 
\[
\sum_{B \in \fB} \int_{B} |\nabla v|_{\sigma}^2 \vol_{\sigma} < \ep_0.
\]
Then, letting $w = \cR(\sigma, v, \fB)$, we have for all Borel set $A \subset S$ that 
\[
\int_{A}|\nabla w|_{\sigma}^2 \vol_{\sigma} \leq \int_{A\setminus (\cup_{B \in \fB}B)} |\nabla v|_{\sigma}^2 \vol_{\sigma} + \sum_{B \in \fB} \int_{B} |\nabla v|_{\sigma}^2 \vol_{\sigma}.
\]
\end{lemm}
\begin{proof}
Since $w = v$ on the complement of $\cup_{B\in \fB}B$, and since $w$ minimizes $E(\sigma, \cdot)$ on each $B \in \fB$ among maps agreeing with it on $\partial B$, we have
\[
\begin{split}
\int_{A} |\nabla w|_{\sigma}^2 \vol_{\sigma} =\  & \int_{A \cap (\cup_{B \in \fB}B)} |\nabla w|_{\sigma}^2 \vol_{\sigma} + \int_{A \setminus (\cup_{B \in \fB}B)} |\nabla v|_{\sigma}^2 \vol_{\sigma}\\
\leq\ & \sum_{B \in \fB}\int_{B} |\nabla w|_{\sigma}^2 \vol_{\sigma} + \int_{A \setminus (\cup_{B \in \fB}B)} |\nabla v|_{\sigma}^2 \vol_{\sigma}\\
\leq\ & \sum_{B \in \fB}\int_{B} |\nabla v|_{\sigma}^2 \vol_{\sigma} + \int_{A \setminus (\cup_{B \in \fB}B)} |\nabla v|_{\sigma}^2 \vol_{\sigma}.
\end{split}
\]
This gives the asserted estimate.
\end{proof}
\begin{rmk}\label{rmk:many-replacements}
Suppose $(\sigma, v) \in \cM$ and that for some $L \geq 2$ we have finite, disjoint collections $\fB_1, \cdots, \fB_{L}$ of geodesic disks in $(S, \sigma)$, each with radius at most $\rho_0$, such that 
\[
\sum_{B \in \fB_i} \int_{B} |\nabla v|_{\sigma}^2 \vol_{\sigma} < \ep < \frac{\ep_0}{3^{L-1}} \text{ for }i = 1, \cdots, L.
\]
Then for any distinct $i, j \in \{1, \cdots, L\}$, we have by Lemma~\ref{lemm:energy-after-replacement} that
\[
\sum_{B \in \fB_{j}} \int_{B}|\nabla \cR(\sigma, v, \fB_i)|_{\sigma}^2 \vol_{\sigma} < 2\ep < \frac{\ep_0}{3^{L-2}},
\]
and therefore harmonic replacement of $\cR(\sigma, v, \fB_{i})$ on $\fB_{j}$ is permitted, and we denote the resulting map by $\cR(\sigma, v, \fB_{i}, \fB_{j})$. More generally, for all $k = 2, \cdots, L$ and distinct $i_1, \cdots, i_k \in \{1, \cdots, L\}$, we may perform successive harmonic replacements on $\fB_{i_1}, \cdots, \fB_{i_k}$ to obtain what we subsequently denote by $\cR(\sigma, v, \fB_{i_1}, \cdots, \fB_{i_k})$.
\end{rmk}

We now state the estimates due to Colding--Minicozzi for iterated replacements, as generalized by the second named author~\cite{Zhou10,Zhou17b}. (See also~\cite{Laurain-Petrides2019}.)
\begin{prop}
\label{prop:iterated-replacement-estimates}
Suppose $(\sigma, v) \in \cM'$ and that $\fB_1, \fB_2$ are two finite, disjoint collections of closed geodesic disks in $(S, \sigma)$, each with radius at most $\rho_0$, such that 
\begin{equation}\label{eq:iterated-replacement-smallness}
\sum_{B \in \fB_i} \int_{B} |\nabla v|_{\sigma}^2 \vol_{\sigma} < \frac{\ep_0}{3},\quad \text{for }i = 1, 2.
\end{equation}
Then, letting $v_1 = \cR(\sigma, v, \fB_{1})$, we have
\begin{equation}\label{eq:iterated-replacement-estimate-1}
\begin{split}
&E(\sigma, v_1) - E(\sigma, \cR(\sigma, v_1, \fB_2)) \\
& \geq E(\sigma, v) - E(\sigma, \cR(\sigma, v, \frac{1}{2}\fB_{2})) - \frac{1}{\kappa}\big[ E(\sigma, v) - E(\sigma, v_1) \big]^{\frac{1}{2}},
\end{split}
\end{equation}
\begin{equation}\label{eq:iterated-replacement-estimate-2}
\begin{split}
&E(\sigma, v_1) - E(\sigma, \cR(\sigma, v_1, \frac{1}{2}\fB_2))\\
& \leq E(\sigma, v) - E(\sigma, \cR(\sigma, v, \fB_2)) +  \frac{1}{\kappa}\big[ E(\sigma, v) - E(\sigma, v_1) \big]^{\frac{1}{2}},
\end{split}
\end{equation}
where again the constant $\kappa$ depends only on $M$.
\end{prop}
\begin{proof}
Having established Proposition \ref{prop:iterated-replacement}, we can deduce this proposition by arguing as in \cite{Colding-Minicozzi08b}. See Appendix~\ref{appendix:estimates-iterated-replacement}.
\end{proof}
The estimates in Proposition \ref{prop:iterated-replacement-estimates} may be iterated further to yield the following corollary.
\begin{coro}[See also \cite{Laurain-Petrides2019}, proof of Theorem 4.1]
\label{coro:iterated-iterated-replacement}
Suppose $(\sigma, v) \in \cM'$, and that for some $L \geq 2$ we have finite, disjoint collections $\fB_{1}, \cdots, \fB_{L}$ of closed geodesic disks in $(S, \sigma)$ with radius at most $\rho_0$, such that 
\begin{equation}\label{eq:iterated-iterated-smallness}
\sum_{B\in \fB_j} \int_{B} |\nabla v|_{\sigma}^2 \vol_{\sigma} < \frac{\ep_0}{3^{L-1}}, \text{ for each }j = 1, \cdots, L.
\end{equation}
Then, in the notation of Remark~\ref{rmk:many-replacements}, we have
\begin{equation}\label{eq:many-iteration-1}
\begin{split}
&E(\sigma, v) - E(\sigma, \cR(\sigma, v, 2^{-L + 1}\fB_{L}))\\
\leq\ & E(\sigma , \cR(\sigma, v, \fB_1, \cdots, \fB_{L-1})) - E(\sigma, \cR(\sigma, v, \fB_1, \cdots, \fB_L)) \\
& + \frac{L - 1}{\kappa}\big[ E(\sigma, v) - E(\sigma, \cR(\sigma, v, \fB_1, \cdots, \fB_{L-1})) \big]^{\frac{1}{2}}\\
\leq\ & \big( 1 + \frac{L - 1}{\kappa} \big)\big[ E(\sigma, v) - E(\sigma, \cR(\sigma, v, \fB_1, \cdots, \fB_{L})) \big]^{\frac{1}{2}},
\end{split}
\end{equation}
and that
\begin{equation}\label{eq:many-iteration-2}
\begin{split}
&E(\sigma, \cR(\sigma, v, \fB_1, \cdots, \fB_{L-1})) - 
E(\sigma, \cR(\sigma, v, \fB_1, \cdots, \fB_{L - 1}, 2^{-L + 1}\fB_{L}))\\
\leq\ & E(\sigma, v) - E(\sigma, \cR(\sigma, v, \fB_{L})) + \frac{L - 1}{\kappa}\big[ E(\sigma, v) - E(\sigma, \cR(\sigma, v, \fB_1, \cdots, \fB_{L-1})) \big]^{\frac{1}{2}}
\end{split}
\end{equation}
\end{coro}
\begin{proof}
For the sake of completeness we include a proof in Appendix~\ref{appendix:estimates-iterated-replacement}.
\end{proof}
\subsection{Energy decreasing process for families of mappings}\label{subsec:energy-decreasing}
The main results of this section are Propositions \ref{prop:choice-of-disks} and \ref{prop:energy-decreasing}, in which harmonic replacement is applied to continuous families of pairs in $\cM'$. We first make a preliminary definition and state a technical result (Proposition \ref{prop:e-semicontinuous}) whose proof we present in Appendix \ref{appendix:semi-continuity}. Let $\ep_0 > 0$ be the constant in Theorem~\ref{thm:convexity}. Given a pair $(\sigma, v) \in \cM'$, an energy threshold $\ep \in (0, \frac{\ep_0}{3})$, and a scaling factor $\lambda \in (0, 1)$, we define the \emph{maximal energy drop} with respect to these parameters to be
\begin{equation}\label{eq:maximal-drop}
e(\sigma, v, \ep, \lambda) = \sup\big\{E(\sigma, v) - E(\sigma, \cR(\sigma, v, \lambda\fB))\big\},
\end{equation}
where the supremum is taken over all finite disjoint collections $\fB$ of closed geodesic disks in $(S, \sigma)$ such that each has radius at most $\rho_0$, and that
\[
\sum_{B \in \fB} \int_{B} |\nabla v|_{\sigma}^2 \vol_{\sigma} < \ep.
\]
The following semicontinuity property of the maximal energy drop is key to the entire construction in this section.
\begin{prop}[\cite{Colding-Minicozzi08b}, Lemma 3.20; \cite{Zhou17b}, Lemma 4.11]
\label{prop:e-semicontinuous}
Suppose $(\sigma, v) \in \cM'$ and that $v$ is not a harmonic map with respect to $\sigma$. Suppose also that $(\sigma_n, v_n)$ is a sequence in $\cM'$ converging to $(\sigma ,v)$. Then we have
\[
\limsup_{n \to \infty} e(\sigma_n, v_n, \frac{\ep}{3}, \frac{\lambda}{4}) \leq e(\sigma, v, \ep, \lambda).
\]
\end{prop}
\begin{proof}
The proof involves the same ideas as in the work of Colding and Minicozzi~\cite{Colding-Minicozzi08b}, but the execution in our context is rather technical. We give the details in Appendix~\ref{appendix:semi-continuity}.
\end{proof}
To continue, fix some $d \in \NN$. A given subset of $\cM'$ of the form $\{(\sigma_t, v_t)\}_{t \in X}$, where $X \subset \RR^d$, is said to be a \emph{continuous family} provided there is a continuous map $(\bsig, \bv): X \to \cM'$ such that $(\sigma_{t}, v_{t}) = (\bsig(t), \bv(t))$ for all $t \in X$. We will mostly take $X$ to be either $\Int(I^d)$ or $I \times \Int(I^d)$, where recall that $I = [0, 1]$, while for $\delta \in (0, \frac{1}{4})$ we let $I_{\delta} = [\delta, 1 - \delta]$. Whenever a continuous family $\{(\sigma_t, v_t)\}$ is given, we often write
\[
(\tau_t, u_t) := \Psi(\sigma_t, v_t),
\]
which varies continuously in $\cM$ with $t$ by Proposition \ref{prop:coming-down}, and let $f_{t}$ be the canonical lift of $\id: \Sigma_0 \to (S, [\sigma_{t}])$ with respect to $p_0$.
\begin{equation}\label{eq:diagram-family}
\begin{tikzcd}
\mathbb{H} \arrow[r, "\widetilde{f}_t"] \arrow[d, "p_0"']    & \mathbb{H} \arrow[d, "p_t"] \arrow[r, "u_t"] & M \\
\Sigma_0 = (S, [\gamma_0]) \arrow[r, "\id"]  & {(S, [\sigma_t])} \arrow[ru, "v_t"']  &  
\end{tikzcd}
\end{equation}

We also need a version of the Besicovitch covering theorem. Specifically, the following statement can be established by following the proof of \cite[Section 1.5.2, Theorem 2]{Evans-Gariepy}. We use $U(t, r)$ to denote the open ball in $\RR^d$ with center $t$ and radius $r$. Given $U = U(t, r)$, we write $\lambda U$ for $U(t, \lambda r)$.
\begin{lemm}[\cite{Evans-Gariepy}, Section 1.5.2]
\label{lemm:besicovitch}
There exists a constant $N = N(d)$ with the following property. Given a bounded set $X \subset \RR^d$, suppose to each $t \in X$ there is associated a radius $r_t > 0$, such that 
\[
\sup_{t \in X}r_{t} < \infty.
\]
Then there exists an at most countable subset $\{t_{j}\}_{j = 1}^{J} \subset X$, where $J$ could be $\infty$, such that
\vskip 1mm
\begin{enumerate}
\item[(a)] $X \subset \bigcup_{j = 1}^{J} U(t_j, r_{t_j})$.
\vskip 1mm
\item[(b)] The collection of closures, namely $\big\{ \overline{U(t_j, r_{t_j})} \big\}_{j = 1}^{J}$, can be partitioned into $N$ subcollections of mutually disjoint balls.
\end{enumerate}
\end{lemm}
\begin{proof}
Below, by Claim \#1, Claim \#2, and so on, we mean the claims appearing in the proof of~\cite[Section 1.5.2, Theorem 2]{Evans-Gariepy}. Following the argument behind Claims \#1 through \#4, we obtain an at most countable subset $\{t_{j}\}_{j = 1}^{J}$ of $X$ such that, writing $r_{j}$ for $r_{t_j}$, we have
\begin{enumerate}
\item[(i)] $t_j \not\in U(t_i, r_i)$ and $r_j \leq \frac{4}{3}r_i$, whenever $j > i$.
\vskip 1mm
\item[(ii)] $\overline{U(t_j, \frac{r_j}{3})} \cap \overline{U(t_i, \frac{r_i}{3})} = \emptyset$, whenever $j \neq i$.
\vskip 1mm
\item[(iii)] $X \subset \bigcup_{j = 1}^{J} U(t_j, r_j)$.
\end{enumerate}
To prove that $\big\{ \overline{U(t_j, r_{j})} \big\}_{j = 1}^{J}$ can be partitioned in the fashion described in (b), for which we continue to follow~\cite{Evans-Gariepy}, requires only the first two properties listed above, and in fact (ii) follows from (i) since the latter implies 
\[
|t_i - t_j| \geq r_i > \frac{r_i}{3} + \frac{r_j}{3}, \quad \text{whenever }j > i,
\]
as explained in Claim \#2. At any rate, we next fix $k > 1$ and define 
\[
\begin{split}
I =\ & \{j \in \{1, \cdots, k-1\}\ |\ \overline{U(t_j, r_j)} \cap \overline{U(t_k, r_k)} \neq\emptyset\}, \\
K =\ & \{j \in I\ |\ r_j \leq 3r_k\}.
\end{split}
\]
The proof of Claim \#5 shows that $\overline{U(t_j, \frac{r_j}{3})} \subset \overline{U(t_k, 5r_k)}$ for all $j \in K$, and that combining this with (i) and (ii) yields $20^{d}$ as an upper bound for the cardinality of $K$.

On the other hand, given $i, j \in I \setminus K$ with $i \neq j$, and assuming without loss of generality that $t_k = 0$ and $|t_i| \leq |t_j|$, in place of the inequalities above Claim \#6a, we now have
\[
3r_k < r_i \leq |t_i| \leq r_i + r_k,\quad 3r_k < r_j \leq |t_j| \leq r_j + r_k.
\]
Letting $\theta \in [0, \pi]$ denote the angle between the vectors $t_i, t_j \in \RR^d$ and feeding the above inequalities into the proof of Claims \#6a through \#6c, we see that
\[
\cos\theta \leq 
\left\{
\begin{array}{ll}
1/2, & \text{ if }|t_i - t_j| \geq |t_j|,\\
5/6, & \text{ if }|t_j| > |t_i - t_j| \geq r_j,\\
61/64, & \text{ if }r_j > |t_i - t_j|.
\end{array}
\right.
\]
The proof of Claim \#7 then goes through without change to give a dimensional bound on the cardinality of $I \setminus K$. Steps 14 and 15 of the proof of \cite[Section 1.5.2, Theorem 2]{Evans-Gariepy} now gives us the desired partition in part (b).
\end{proof}

\begin{prop}[\cite{Colding-Minicozzi08b}, Lemma 3.24;~\cite{Zhou17b}, Lemma 4.12]
\label{prop:choice-of-disks}
Let $\ep \in (0, \frac{\ep_0}{3})$ and $W > 0$ be given, where $\ep_0$ is the threshold from Theorem~\ref{thm:convexity}. Suppose $\{(\sigma_t, v_t)\}_{t \in \Int(I^d)}$ is a continuous family in $\cM'$ satisfying for some $\delta \in (0, \frac{1}{20})$ that
\vskip 1mm
\begin{enumerate}
\item[(r1)] $\emptyset \neq \{t \in \Int(I^d)\ |\ A(v_t) \geq \frac{W}{2}\} \subset I_{5\delta}^d$.
\vskip 1mm
\item[(r2)] $A(v_t) \geq \frac{W}{2}$ implies that the map $v_t: (S, \sigma_t) \to M$ is not harmonic.
\end{enumerate}
Then there exist finitely many points $t_1, \cdots, t_m \in I^d_{5\delta}$ and, associated with each $t_j$, the following objects:
\begin{itemize}
\item an open disk $U_{j} \subset \RR^d$ centered at $t_j$ with $2U_{j} \subset \Inte(I_{4\delta}^d)$,
\vskip 1mm
\item a continuous function $r_j:I^d \to [0, 1]$ supported in $2 U_{j}$,
\vskip 1mm
\item a finite disjoint collection $\fB_{j} = \{B_{j, \alpha}\}_{\alpha \in A_j}$ of closed geodesic disks in $(S, \sigma_{t_j})$ with radii at most $\rho_0$, along with a choice of lift $\widetilde{B}_{j, \alpha}$ for each $B_{j, \alpha}$ with respect to $p_{t_j}$,
\end{itemize}
such that the following hold.
\vskip 1mm
\begin{enumerate}
\item[(a)] For all $t \in \Int(I^d)$, the number of $j$'s for which $r_j(t) > 0$ is at most $N$, the constant from Lemma \ref{lemm:besicovitch}.
\vskip 1mm
\item[(b)] If $t \in 2U_{j}$, then $\cup_{\alpha \in A_j}\widetilde{B}_{j, \alpha}$ is $\tau_t$-admissible in the sense of Section \ref{subsec:admissible-subsets}, and we have
\begin{equation}\label{eq:3N-energy-bound}
\sum_{\alpha \in A_j}\int_{p_t(\widetilde{B}_{j, \alpha})} |\nabla v_t|_{\sigma_t}^2 \vol_{\sigma_t} < \frac{\ep}{3^{N}}.
\end{equation}
Moreover, letting $\fB_{j, t} = \{p_t(\widetilde{B}_{j, \alpha})\}_{\alpha \in A_j}$, we have
\begin{equation}\label{eq:energy-drop-lower-bound}
E(\sigma_t, v_t) - E(\sigma_t, \cR(\sigma_t, v_t, 2^{-N}\fB_{j, t})) \geq \frac{1}{8}e(\sigma_t, v_t, \frac{\ep}{3^{N + 2}}, 2^{-N - 2}).
\end{equation}
\vskip 1mm
\item[(c)] If $t \in \Int(I^d)$ is such that $A(v_t) \geq \frac{W}{2}$, then there exists $j = j(t)$ such that $r_j(t) = 1$.
\end{enumerate}
\end{prop}
\begin{proof}
By assumption (r1), letting 
\[
C := \{t \in \Int(I^d)\ |\ A(v_{t}) \geq \frac{W}{2}\}
\]
yields a non-empty, compact subset of $\Inte(I_{4\delta}^d)$, while assumption (r2) guarantees that for each $t \in C$ we have
\begin{equation}\label{eq:choice-of-disk-not-harmonic}
e(\sigma_t, v_t, \frac{\ep}{3^{N + 1}}, 2^{-N }) > 0,
\end{equation}
and consequently there exists a finite disjoint collection $\fB_t = \{B_{t, \alpha}\}_{\alpha \in A_t}$ of disks in $(S, \sigma_t)$ such that each has radius at most $\rho_0$, and that
\begin{equation}\label{eq:t-collection-small-energy}
\sum_{\alpha \in A_t} \int_{B_{t, \alpha}} |\nabla v_t|_{\sigma_t}^2 \vol_{\sigma_t} < \frac{\ep}{3^{N+1}},
\end{equation}
\begin{equation}\label{eq:t-collection-almost-max-drop}
E(\sigma_t, v_t) - E(\sigma_t, \cR(\sigma_t, v_t, 2^{-N}\fB_{t})) \geq \frac{1}{2}e(\sigma_t, v_t, \frac{\ep}{3^{N + 1}}, 2^{-N}).
\end{equation}
Now fix a lift $\widetilde{B}_{t, \alpha}$ for each $B_{t, \alpha}$ with respect to $p_{t}$. In particular, the compact set $\cup_{\alpha \in A_{t}} \widetilde{B}_{t, \alpha}$ is $\tau_{t}$-admissible by Remark \ref{rmk:upstairs-and-downstairs}. Since the continuity of $s \mapsto \tau_s$ allows us to invoke Lemma~\ref{lemm:admissible-set-limit}, and since $C \subset \Int(I_{4\delta}^{d})$, we can find an open disk $U_t \subset \RR^d$ centered at $t$ such that 
\[
2U_{t} \subset \Inte(I_{4\delta}^d),
\]
and that, for all $s \in 2U_t$,
\begin{equation}\label{eq:close-admissible}
\cup_{\alpha \in A_t}\widetilde{B}_{t, \alpha} \text{ is $\tau_s$-admissible},
\end{equation}
in which case, by Remark \ref{rmk:upstairs-and-downstairs} again, we get a finite, disjoint collection of geodesic disks in $(S, \sigma_s)$, each having radius at most $\rho_0$, upon letting
\[
\fB_{t, s} = \{p_s(\widetilde{B}_{t, \alpha})\}_{\alpha \in A_t}.
\]
Furthermore, expressing~\eqref{eq:t-collection-small-energy} in terms of $u_{t}$, and using the continuity of $s \mapsto u_s$ in $(C^0 \cap W^{1, 2})_{\loc}(\HH)$, we get after shrinking $U_t$ if necessary that
\begin{equation}\label{eq:close-energy}
\sum_{\alpha} \int_{p_s(\widetilde{B}_{t, \alpha})} |\nabla v_s|_{\sigma_s}^2 \vol_{\sigma_s} = \sum_{\alpha}\int_{\widetilde{B}_{t, \alpha}}|\nabla u_s|^2 < \frac{\ep}{3^N},\quad \text{for all }s \in 2U_{t}.
\end{equation}
Using again the continuity of $s \mapsto u_s$, this time combined with Theorem~\ref{thm:replacement-continuity-1}, we have for all $\alpha \in A_t$ that
\[
\lim_{s \to t}\int_{\widetilde{B}_{t, \alpha}} |\nabla u_s|^2 - |\nabla \cR(u_s, 2^{-N}\widetilde{B}_{t, \alpha})|^2 = \int_{\widetilde{B}_{t, \alpha}} |\nabla u_t|^2 - |\nabla \cR(u_t, 2^{-N}\widetilde{B}_{t, \alpha})|^2.
\]
Combining this with~\eqref{eq:t-collection-almost-max-drop}, and also using Proposition~\ref{prop:e-semicontinuous} and~\eqref{eq:choice-of-disk-not-harmonic}, upon shrinking $U_t$ further, we can also arrange that for all $s \in 2U_{t}$ there holds
\begin{equation}\label{eq:close-drop}
\begin{split}
E(\sigma_s, v_s) - E(\sigma_s, \cR(\sigma_s, v_s, 2^{-N}\fB_{t, s})) \geq\ & \frac{1}{4}e(\sigma_t, v_t, \frac{\ep}{3^{N + 1}}, 2^{-N}) \\
\geq\ & \frac{1}{8}e(\sigma_s, v_s, \frac{\ep}{3^{N + 2}}, 2^{-N - 2}).
\end{split}
\end{equation}
We now apply Lemma \ref{lemm:besicovitch} to the collection $\{U_{t}\}_{t \in C}$, and use the compactness of $C$ to extract a finite subcovering. The result is that there exist points $t_1, \cdots, t_m \in C$ so that 
\vskip 1mm
\begin{enumerate}
\item[(i)] $C \subset \bigcup_{j =1}^m U_{t_j}$,
\vskip 1mm
\item[(ii)] $\{\overline{U_{t_1}}, \cdots, \overline{U_{t_m}}\}$ can be partitioned into $N$ collections of mutually disjoint disks.
\end{enumerate}
For each $j \in \{1, \cdots, m\}$, we let
\[
\fB_j := \fB_{t_j}, \quad U_j := U_{t_j},
\]
and choose a continuous function $r_j:I^d \to [0, 1]$ satisfying
\[
r_j(t) = 1 \text{ if }t \in U_{j},\quad \supp(r_j) \subset 2U_{j}.
\]
As $\{\overline{U_{1}}, \cdots, \overline{U_{m}}\}$ is a finite collection of compact sets in $\RR^d$, we can also arrange that
\[
\supp(r_i) \cap \supp(r_j) = \emptyset\quad \text{ whenever }\overline{U_{i}} \cap \overline{U_{j}} = \emptyset.
\]
Point (ii) above then implies that $\{\supp(r_1), \cdots, \supp(r_m)\}$ can likewise be partitioned into $N$ subcollections, each mutually disjoint. In particular, for each $t \in \Int(I^d)$, the number of $j$'s for which $r_j(t) > 0$ is at most $N$, which verifies conclusion (a). Next, when $t \in 2U_{j}$, the admissibility assertion in conclusion (b) follows from \eqref{eq:close-admissible}, while the estimates \eqref{eq:3N-energy-bound} and \eqref{eq:energy-drop-lower-bound} follow from \eqref{eq:close-energy} and \eqref{eq:close-drop}, respectively. Finally, if $t \in C$, then property (i) above yields some $j$ such that $t \in U_{j}$, in which case $r_j(t) = 1$. This proves (c).\\
\end{proof}
\begin{prop}\label{prop:energy-decreasing}
There exists a continuous, increasing function $\Theta:[0, \infty) \to [0, \infty)$, depending only on $M$ and $d$, which satisfies $\Theta(0) = 0$ and has the following property. Let $\ep \in (0, \frac{\ep_0}{3})$, $W> 0$, $\{(\sigma_t, v_t)\}_{t \in \Int(I^d)}$, and $\delta \in (0, \frac{1}{20})$ be as in the hypotheses of Proposition~\ref{prop:choice-of-disks}, so that both (r1) and (r2) hold. Then there exists a continuous family $\{(\sigma_t, \widehat{v}_{s, t})\}_{(s, t) \in [0, 1] \times \Int(I^d)}$ in $\cM'$ such that
\begin{equation}\label{eq:energy-decreasing-bc}
\widehat{v}_{s,t} = v_t, \text{ if $s =0$ or $t \not\in \Int(I_{4\delta}^d$}).
\end{equation}
In addition, writing $\widehat{v}_t$ for $\widehat{v}_{1, t}$, we have:
\vskip 1mm
\begin{enumerate}
\item[(a)] $E(\sigma_t, \widehat{v}_t) \leq E(\sigma_t, v_t)$  for all $t\in \Int(I^d)$.
\vskip 1mm
\item[(b)] If $A(v_t) \geq \frac{W}{2}$ and if $\fB = \{B_i\}$ is a finite, disjoint collection of closed geodesic disks in $(S, \sigma_t)$, each with radius at most $\rho_0$, such that 
\[
\sum_{i}\int_{B_i} |\nabla \widehat{v}_t|_{\sigma_t}^2 \vol_{\sigma_t} < \frac{\ep}{3^{N + 2}},
\]
then
\[
\begin{split}
\int_{S} |\nabla \widehat{v}_t - \nabla \cR(\sigma_t, \widehat{v}_t, 2^{-2N - 2}\fB)|_{\sigma_t}^2 \vol_{\sigma_t} \leq\ & 4\big[ E(\sigma_t, \widehat{v}_t) - E(\sigma_t, \cR(\sigma_t, \widehat{v}_t, 2^{-2N - 2}\fB))\big]\\ 
\leq\ & \Theta(E(\sigma_t, v_t) - E(\sigma_t, \widehat{v}_t)).
\end{split}
\]
\end{enumerate}
\end{prop}
\begin{proof}
Let $t_j \in I_{5\delta}^{d}$, $U_{j} \subset \RR^d$, $r_j:I^d \to [0, 1]$, and $\{(B_{j, \alpha}, \widetilde{B}_{j, \alpha})\}_{\alpha \in A_j}$ (for $j = 1, \cdots, m$) be as given by Proposition~\ref{prop:choice-of-disks}. If $j$ and $t$ are such that $t \in 2U_j$, we define, as before,
\[
\fB_{j, t} = \{p_t(\widetilde{B}_{j, \alpha})\}_{\alpha \in A_j},
\]
which makes sense by Proposition~\ref{prop:choice-of-disks}(b), while if $t \not\in 2U_j$ we set $\fB_{j, t} = \emptyset$. Recalling that $\supp(r_j) \subset 2U_j$, and that for each $t$, at most $N$ of the $r_j(t)$'s are non-zero, we see from Remark~\ref{rmk:many-replacements} and the bound \eqref{eq:3N-energy-bound} in Proposition~\ref{prop:choice-of-disks}(b) that the following is well-defined for all $k \in \{1, \cdots, m\}$ and $(s, t) \in [0, 1] \times \Int(I^d)$:
\[
\widehat{v}^{k}_{s, t} = \cR(\sigma_t, v_t, \frac{s}{2}r_1(t)\fB_{1, t}, \cdots, \frac{s}{2}r_k(t)\fB_{k, t}).
\]
We then define $\widehat{v}_{s, t}: = \widehat{v}^m_{s, t}$, and observe that~\eqref{eq:energy-decreasing-bc} holds since $2U_1, \cdots, 2U_{m}$ are all contained in $\Int(I_{4\delta}^d)$. Also, since harmonic replacement does not increase the Dirichlet energy, we get conclusion (a) as well.

Next we argue by induction that the map $(s, t) \mapsto (\sigma_t, \widehat{v}^{k}_{s, t})$ from $[0, 1] \times \Int(I^d)$ to $\cM'$ is continuous for all $k \in \{1, \cdots, m\}$. For $k = 1$, by considering separately the cases $t \not\in \supp(r_1)$ and $t \in 2U_1$, and using Corollary~\ref{coro:replacement-downstairs-continuity} in the latter case, we get that 
\[
(s, t) \mapsto (\sigma_t, \cR(\sigma_t, v_t, \frac{s}{2} \cdot r_1(t)\fB_{1, t})),
\]
is continuous as a map into $\cM'$. Suppose by induction that $(s, t) \mapsto (\sigma_t, \widehat{v}^{k}_{s, t})$ is continuous for some $k \in \{1, \cdots, m-1\}$, and take $(s_0, t_0) \in [0, 1] \times \Int(I^d)$. In the case $t_0 \not\in 2U_{k + 1}$, since $r_{k+1}$ is supported in $2U_{k + 1}$, we can find a neighborhood $U$ of $t_0$ such that
\[
\widehat{v}^{k + 1}_{s, t} = \widehat{v}^{k}_{s, t},\quad \text{for all }s \in [0, 1],\ t \in U.
\]
Thus in this case $(s, t) \mapsto (\sigma_t, \widehat{v}^{k + 1}_{s, t})$ is continuous at $(s_0, t_0)$ by the induction hypothesis. On the other hand, in the case $t_0 \in 2U_{k + 1}$, we have that $\cup_{\alpha} \widetilde{B}_{k + 1, \alpha}$ is $\tau_{t_0}$-admissible by Proposition \ref{prop:choice-of-disks}(b), and we distinguish the following two subcases:
\vskip 1mm
\begin{enumerate}
\item[(i)] If only at most $N-1$ of $r_1(t_0), \cdots, r_k(t_0)$ are positive, then by the energy bound \eqref{eq:3N-energy-bound} and Remark~\ref{rmk:many-replacements}, we must have
\[
\sum_{\alpha \in A_{k + 1}} \int_{p_{t_0}(\widetilde{B}_{k + 1, \alpha})} |\nabla \widehat{v}^{k}_{s_0, t_0}|_{\sigma_{t_0}}^2 \vol_{\sigma_{t_0}} < \frac{2^{N-1}\ep}{3^{N}} < \frac{\ep}{3},
\]
and hence, in view of the relation
\[
\widehat{v}^{k + 1}_{s, t} = \cR(\sigma_t, \widehat{v}^{k}_{s, t}, \frac{s}{2} r_{k + 1}(t)\fB_{k + 1, t}),
\]
as well as how $\fB_{k + 1, t}$ is defined in terms of $\{\widetilde{B}_{k + 1, \alpha}\}_{\alpha \in A_{k +1}}$ when $t \in 2U_{k + 1}$, we may use Corollary~\ref{coro:replacement-downstairs-continuity} to get the continuity of $(s, t) \mapsto (\sigma_t, \widehat{v}^{k + 1}_{s, t})$ at $(s_0, t_0)$. 
\vskip 1mm
\item[(ii)] If $N$ of $r_1(t_0), \cdots, r_{k}(t_0)$ are positive, then by continuity the same remains true for $t$ close to $t_0$, in which case Proposition~\ref{prop:choice-of-disks}(a) forces $r_{k + 1}(t)$ to vanish. Thus we find, as in the case where $t_0 \not\in 2U_{k + 1}$, that $\widehat{v}^{k + 1}_{s, t} = \widehat{v}^{k}_{s, t}$ for all $t$ near $t_0$ and $s \in [0, 1]$.
\end{enumerate}
This completes the inductive argument. Recalling that $\widehat{v}_{s, t} = \widehat{v}^m_{s, t}$ by definition, we see that $(s, t) \mapsto (\sigma_t, \widehat{v}_{s,t})$ is continuous from $[0, 1] \times \Int(I^d)$ to $\cM'$. 

It remains to establish (b). Given $t$ such that $A(v_t) \geq \frac{W}{2}$, by Proposition~\ref{prop:choice-of-disks}, we may assume that the set of $j$'s for which $r_j(t) >0$ is given by 
\[
1 \leq i_1 < \cdots < i_M \leq m,
\]
for some $M \leq N$, with the index $j(t)$ in Proposition~\ref{prop:choice-of-disks}(c) occurring as $i_L$ for some $L \in \{1, \cdots, M\}$. In particular, we have $t \in 2U_{i_{j}}$ for all $j \in \{1, \cdots, M\}$. Thanks to the upper bound \eqref{eq:3N-energy-bound}, we may apply the estimate~\eqref{eq:many-iteration-1} from Corollary \ref{coro:iterated-iterated-replacement}, where we take $L$ to be as present, $\fB_{j}$ to be $\frac{r_{i_j}(t)}{2}\fB_{i_j, t}$ for $j = 1, \cdots, L-1$, and $\fB_{L}$ to be $\frac{1}{2}\fB_{i_L, t}$. Combining the resulting inequality with the lower bound \eqref{eq:energy-drop-lower-bound}, we deduce that
\begin{equation}\label{eq:energy-decreasing-total-drop}
\begin{split}
&\frac{1}{8}e(\sigma_t, v_t, \frac{\ep}{3^{N + 2}}, 2^{-N-2}) \\
\leq\ & E(\sigma_t, v_t) - E(\sigma_t, \cR(\sigma_t, v_t, 2^{-N}\fB_{i_{L}, t}))\\
\leq\ & E(\sigma_t, v_t) - E(\sigma_t, \cR(\sigma_t, v_t, 2^{-L}\fB_{i_{L}, t}))\\
\leq\ & \big( 1 + \frac{L-1}{\kappa} \big)\big[ E(\sigma_t, v_t) - E(\sigma_t, \cR(\sigma_t, v_t, \frac{r_{i_1}(t)}{2}\fB_{i_1, t}, \cdots, \frac{1}{2}\fB_{i_L, t})) \big]^{\frac{1}{2}}\\
\leq\ & \big( 1 + \frac{N-1}{\kappa} \big)\big[ E(\sigma_{t}, v_{t}) - E(\sigma_t, \widehat{v}_t)\big]^{\frac{1}{2}}.
\end{split}
\end{equation}
Now, given $\fB$ as in the statement of part (b), applying instead the estimate \eqref{eq:many-iteration-2} from Corollary \ref{coro:iterated-iterated-replacement}, with $L = M+1$, $\fB_{M+1} = 2^{-N-2}\fB$, and $\fB_{j} = \frac{r_{i_j}(t)}{2}\fB_{i_j, t}$ for $j = 1, \cdots, M$, we obtain
\[
\begin{split}
&E(\sigma_t, \widehat{v}_t) - E(\sigma_t, \cR(\sigma_t, \widehat{v}_t, 2^{-2N - 2}\fB))\\
\leq\ & E(\sigma_t, \widehat{v}_t) - E(\sigma_t, \cR(\sigma_t, \widehat{v}_t, 2^{-M - N - 2}\fB))\\
\leq\ & E(\sigma_t, v_t) - E(\sigma_t, v_t, 2^{-N - 2}\fB) + \frac{M}{\kappa}\big[E(\sigma_t, v_t) - E(\sigma_t,  \widehat{v}_t) \big]^{\frac{1}{2}}\\
\leq\ & e(\sigma_t, v_t, \frac{\ep}{3^{N + 2}}, 2^{-N-2}) + \frac{N}{\kappa}\big[ E(\sigma_t, v_t) - E(\sigma_t,  \widehat{v}_t) \big]^{\frac{1}{2}},
\end{split}
\]
where the last line follows from our assumptions on the collection $\fB$ and the definition of $e(\sigma_t, v_t, \frac{\ep}{3^{N _+ 2}}, 2^{-N-2})$. Combining this with~\eqref{eq:energy-decreasing-total-drop} gives
\begin{equation}\label{eq:drop-after-decreasing}
E(\sigma_t, \widehat{v}_t) - E(\sigma_t, \cR(\sigma_t, \widehat{v}_t, 2^{-2N - 2}\fB)) \leq 9\big( 1 + \frac{N}{\kappa} \big)\big[ E(\sigma_t, v_t) - E(\sigma_t,  \widehat{v}_t) \big]^{\frac{1}{2}}.
\end{equation}
To finish, recall from Remark \ref{rmk:v-replacement} that, as a consequence of Theorem~\ref{thm:convexity}, we have
\begin{equation}\label{eq:drop-bounds-distance}
\begin{split}
\frac{1}{4}\int_{S}|\nabla \widehat{v}_t - \nabla \cR(\sigma_t, \widehat{v}_t, 2^{-2N - 2}\fB)|_{\sigma_t}^2 \vol_{\sigma_t} \leq E(\sigma_t, \widehat{v}_t) - E(\sigma_t, \cR(\sigma_t, \widehat{v}_t, 2^{-2N - 2}\fB)).
\end{split}
\end{equation}
Combining~\eqref{eq:drop-bounds-distance} with~\eqref{eq:drop-after-decreasing} gives the asserted estimate in part (b), with 
\[
\Theta(s) = 36\big( 1 + \frac{N}{\kappa} \big)\sqrt{s}.
\]
\end{proof}
\section{Proof of min-max theorem}\label{sec:proof-of-main}
In this section we prove Theorem~\ref{thm:min-max-existence}. First we recall the setting: $S$ is a closed, oriented surface with genus $g > 1$, equipped with a reference metric $\gamma_0$; $M$ is a closed oriented Riemannian $n$-manifold, isometrically embedded in some Euclidean space $\RR^q$; $\cV$ is a tubular neighborhood of $M$ in $\RR^d$; $\Pi: \cV \to M$ stands for the nearest-point projection. We choose $d_0 > 0$ so that 
\[
\{y \in \RR^q\ |\ \dist(y, M) \leq 4d_0\} \subset \cV.
\]
Also, fixing $m \in \NN$, we let $\cP = \cP(m, S)$ be as in Definition~\ref{defi:class-cP}, and recall that the equivalence class of a given $\bv_{0}\in \cP$ with respect to the relation ``$\sim$'' defined in Section~\ref{subsec:statements} is denoted $[\bv_{0}]$. 

Theorem~\ref{thm:min-max-existence} is a direct consequence of the following result. Recall that harmonic maps from the $2$-sphere are necessarily weakly conformal. (See for example~\cite[Lemma 4.25]{Colding-Minicozzi11}.)
\begin{thm}\label{thm:main-existence}
Suppose $\bv_{0} \in \cP$ is such that $\cW([\bv_{0}]) > 0$. There exist a closed, possibly disconnected Riemann surface $(S_0, c_0)$ with genus at most $g$, a weakly conformal harmonic map $u_0: (S_0, c_0) \to M$, and a finite collection of harmonic $2$-spheres $\{u_{j}\}_{j \in J}$ in $M$, such that 
\[
\cW([\bv_{0}]) = E(h_0, u_0) + \sum_{j \in J}E(g_{S^2}, u_j).
\]
Here $g_{S^2}$ denotes the round metric, while $h_0$ is any representative of the conformal class $c_0$.
\end{thm}
The proof occupies the remainder of Section~\ref{sec:proof-of-main}. First, in Section~\ref{subsec:good-min-max}, starting with a class $[\bv_{0}]$ in $\cP$ with $\cW([\bv_{0}]) > 0$, by the tools developed in Sections~\ref{sec:conformal-reparametrization} and~\ref{sec:harmonic-replacement}, we obtain a min-max sequence $(v_n)$ of maps $S \to M$ with area attaining the width $\cW([\bv_{0}])$, along with hyperbolic metrics $(\sigma_n)$ on $S$, so that $v_n$ is becoming conformal with respect to $\sigma_n$ in a certain sense, and enjoys good compactness properties on disks with small energy. The implications of these two properties for the extraction of bubbles and the analysis of neck regions are derived in Section~\ref{subsec:compactness-properties}. These pave the way for Section~\ref{subsec:convergence}, where we first recall certain aspects of the compactification of the space of hyperbolic structures, from which the limiting Riemann surface $(S_0, c_0)$ arises, before completing the proof of Theorem~\ref{thm:main-existence} by standard bubble tree convergence arguments.

\subsection{Producing a good min-max sequence}\label{subsec:good-min-max}
Recall that $I = [0, 1]$ while $I_{\delta} = [\delta, 1 - \delta]$. In particular, 
\[
I_{\delta'} \subset I_{\delta} \subset I, \quad \text{if }0 <\delta < \delta'.
\]
The main result of this section is Proposition~\ref{prop:good-min-max-seq}, which rests primarily on Proposition~\ref{prop:conformal-reparametrization} and Proposition~\ref{prop:energy-decreasing}. A key assumption in the latter, namely the absence of non-constant harmonic slices, is arranged for using Proposition~\ref{prop:area-min-seq}, whose proof in turn depends on Lemma~\ref{lemm:cut-and-paste}. On the other hand, Lemma~\ref{lemm:smooth-sweepout} contains the mollification process that makes the sweepouts sufficiently smooth for the various operations in this section. Provided $\cW([\bv_{0}]) > 0$, the condition~\eqref{eq:no-area-at-boundary} built into the definition of $\cP$ means that we need only perform these operations away from $\partial I^m$, thereby keeping the sweepouts in the class $[\bv_{0}]$ throughout. 
\begin{lemm}[\cite{Colding-Minicozzi08b}, Lemma D.1]
\label{lemm:smooth-sweepout}
Given $\bv \in \cP$, together with constants $\delta\in (0, \frac{1}{4})$ and $\mu \in (0, d_0)$, there exists $\bw \in [\bv]$ satisfying
\vskip 1mm
\begin{enumerate}
\item[(a)] $\bw(t) \in C^2(S; M)$ for all $t \in I_{\delta}^m$, and $\bw|_{I_{\delta}^m}$ is continuous as a map into $C^2(S; M)$.
\vskip 1mm
\item[(b)] $\bw(t) = \bv(t)$ for all $t \in I^m \setminus I^m_{\delta/2}$.
\vskip 1mm
\item[(c)] $\|\bw(t) - \bv(t)\|_{C^0 \cap W^{1, 2}} < \mu$ for all $t \in I_{\delta/2}^m$.
\end{enumerate}
\end{lemm}
\begin{proof}
See Appendix~\ref{appendix:operations}.
\end{proof}
For the next lemma, we define, for $K > 0$, 
\[
\cX_{K} = \{u \in C^2(\bB_2; M)\ |\ \|Du\|_{\infty; \bB_2} \leq K\}.
\]
\begin{lemm}\label{lemm:cut-and-paste}
Fix $\tau \in (0, \frac{1}{10})$, say $\tau = \frac{1}{12}$. Given $K > 0$ and $\mu \in (0, 1)$, there exist a radius $\overline{\rho} \in (0, 1]$ depending only on $K, \mu, d_0$, and $\|d\Pi\|_{\infty; \cV}$, along with a continuous map 
\[
T:(0, \overline{\rho}] \times \cX_{K} \to C^2(\bB_2; M)
\]
with the following properties:
\vskip 1mm
\begin{enumerate}
\item[(a)] For all $(\rho, u) \in (0, \overline{\rho}] \times \cX_{K}$, we have that $T(\rho, u)(x) = u(0)$ on a neighborhood of $0$, while $T(\rho, u)(x) = u(x)$ if $|x| \geq (1 + \tau)\overline{\rho}$. 
\vskip 1mm
\item[(b)] $A(T(\rho, u); \bB_2) < A(u; \bB_2) + \mu$, for all $(\rho, u) \in (0, \overline{\rho}] \times \cX_{K}$.
\vskip 1mm
\item[(c)] Setting $T(0, u) = u$ extends $T$ to a continuous map 
\[
[0, \overline{\rho}] \times  \cX_{K} \to (C^0 \cap W^{1, 2})(\bB_{2}; M).
\]
\end{enumerate}
\end{lemm}
\begin{proof}
See Appendix~\ref{appendix:operations}.
\end{proof}
\begin{prop}\label{prop:area-min-seq}
Suppose $\bv_0 \in \cP$ is such that $W : = \cW([\bv_0]) > 0$. Given $\ep > 0$, there exist $\bw \in [\bv_0]$ and some $\delta \in (0, \frac{1}{20})$ so that the following hold:
\vskip 1mm
\begin{enumerate}
\item[(a)] $\bw(t) \in C^2(S; M)$ for all $t \in I_{2\delta}^m$, and $\bw|_{I_{2\delta}^m}$ is continuous as a map into $C^2(S; M)$. 
\vskip 1mm
\item[(b)] For all $t \in I_{2\delta}^m$, the map $\bw(t)$ is constant on some open subset of $S$.
\vskip 1mm
\item[(c)] $\sup_{t \in \Int(I^m)}A(\bw(t)) < W + \frac{\ep}{2}$.
\vskip 1mm
\item[(d)] $A(\bw(t)) \geq \frac{3W}{8}$ implies $t \in I_{5\delta}^m$.
\end{enumerate}
\end{prop}
\begin{proof}
By the definition~\eqref{eq:width-defi} of $W$, there exists $\bv \in [\bv_0]$ such that 
\begin{equation}\label{eq:v-area-sup}
\sup_{t \in \Int(I^m)}A(\bv(t)) < W + \frac{\ep}{8}.
\end{equation}
Since $\bv$ agrees with $\bv_0$ near $\partial I^m$, and since $W > 0$ by assumption, upon recalling the boundary condition \eqref{eq:no-area-at-boundary} in the definition of $\cP$, we see it is possible to find $\delta \in (0, \frac{1}{20})$ that satisfies
\begin{equation}\label{eq:large-area-set-of-v}
\bv(t) = \bv_0(t) \quad \text{and} \quad A(\bv_0(t)) < \frac{W}{8},\quad \text{ for all }t \in \Int(I^m) \setminus I_{5\delta}^m.
\end{equation}
Since $\bv(I_{\delta/2}^m)$ is a compact subset of $W^{1, 2}(S; M)$, with respect to the reference metric $\gamma_0$ on $S$ we have
\[
\sup_{t \in I_{\delta/2}^m} E(\gamma_0, \bv(t)) < \infty.
\]
Thanks to this and~\cite[page 2565, estimate (A.4)]{Colding-Minicozzi08b}, by invoking Lemma~\ref{lemm:smooth-sweepout} with a small enough $\mu$, we get some $\bv_1 \in [\bv_0]$ such that $\bv_1|_{I_{\delta}^m}$ is a continuous map into $C^2(S; M)$, that $\bv_{1}(t) = \bv(t)$ for all $t \in I^m \setminus I_{\delta/2}^{m}$, and that 
\[
|A(\bv(t)) - A(\bv_1(t))| < \frac{1}{8}\min\{\ep, W\}, \text{ for all } t \in \Int(I^m).
\]
This last estimate together with~\eqref{eq:v-area-sup} and~\eqref{eq:large-area-set-of-v} gives
\begin{equation}\label{eq:v_1-area-sup}
\sup_{t \in \Int(I^m)}A(\bv_1(t)) < W + \frac{\ep}{4},
\end{equation}
\begin{equation}\label{eq:large-area-set-of-v1}
\{t\in \Int(I^m)\ |\ A(\bv_1(t)) \geq \frac{W}{4}\} \subset I_{5\delta}^m.
\end{equation}
Next take any $p_0 \in S$ and a chart $\varphi:\bB_3 \to S$ with $\varphi(0) = p_0$. For $t \in I_{\delta}^m$ we let 
\[
\widetilde{\bv}_1(t) = \bv_1(t) \circ \varphi.
\]
Note that $\{\widetilde{\bv}_1(t)\}_{t \in I_{\delta}^m}$ is a continuous family in $C^2(\bB_2; M)$, and in particular
\begin{equation}\label{eq:uniform-gradient}
K : = \sup_{t \in I_{\delta}^m,\ x \in \bB_2} \big|D[\widetilde{\bv}_1(t)](x)\big| < \infty.
\end{equation}
With this choice of $K$, and with $\mu$ taken to be $\frac{1}{8}\min\{\ep, W\}$, we apply Lemma~\ref{lemm:cut-and-paste} to obtain $\overline{\rho} \in (0, 1]$ and $T:[0, \overline{\rho}] \times  \cX_{K} \to (C^0 \cap W^{1, 2})(\bB_{2}; M)$, and set
\[
\rho(t) = \overline{\rho} \cdot \frac{\dist(t, I^m \setminus I_{3\delta/2}^m)}{\dist(t, I^m \setminus I_{3\delta/2}^m) +  \dist(t, I_{2\delta}^m)}.
\]
For $(s, t) \in [0, 1] \times I^m$, we define a map $\bh(s, t): S \to M$ as follows. If $t \in I_{\delta}^m$, we let
\[
\bh(s, t)(x) =\left\{
\begin{array}{ll}
T( s \cdot \rho(t), \widetilde{\bv}_1(t))\circ \varphi^{-1}(x), & \text{ if }x \in \varphi(\bB_2),\\
\bv_1(t)(x), & \text{ if }x \in S \setminus \varphi(\bB_2),
\end{array}
\right.
\]
while if $t \in I^m \setminus I_{\delta}^m$, we set $\bh(s, t) = \bv_1(t)$. Notice that
\begin{equation}\label{eq:h-bc}
\bh(s, t)= \bv_1(t), \text{ if either $t \not\in \Inte(I_{3\delta/2}^m)$ or $s = 0$},
\end{equation}
while when $t \in \Inte(I_{3\delta/2}^m)$ and $s > 0$, we see with the help of the second statement in Lemma~\ref{lemm:cut-and-paste}(a) that each $\bh(s, t)$ lies in $C^2(S; M)$. Furthermore, by the continuity of $T|_{(0, \overline{\rho}] \times \cX_{K}}$ into $C^2(\bB_2; M)$, the map
\[
\bh: (0, 1] \times \Int(I_{\frac{3\delta}{2}}^m) \to C^2(S; M)
\]
is continuous. On the other hand, using Lemma~\ref{lemm:cut-and-paste}(c) as well as~\eqref{eq:h-bc}, we see that so is
\[
\bh:[0, 1] \times \Int(I^m) \to (C^0 \cap W^{1, 2})(S; M).
\]
This last continuity property and~\eqref{eq:h-bc} shows that 
\[
\bw : = \bh(1, \cdot) \in [\bv_1] =  [\bv_0].
\]
Next we observe some of the properties of $\bw$. When $t \not \in \Inte(I_{3\delta/2}^m)$, we have $\bw(t)= \bv_1(t)$ by~\eqref{eq:h-bc}. When $t \in \Inte(I_{3\delta/2}^m)$, we have:
\vskip 1mm
\begin{enumerate}
\item[(i)] $\bw(t)\in C^2(S; M)$. Moreover, $\bw: \Inte(I_{3\delta/2}^m) \to C^2(S; M)$ is continuous.
\vskip 1mm
\item[(ii)] By Lemma~\ref{lemm:cut-and-paste}(a), each $\bw(t)$ is constant near $p_0 = \varphi(0)$, and agrees with $\bv_1(t)$ on $S\setminus \varphi(\bB_{(1 + \tau)\overline{\rho}})$.
\vskip 1mm
\item[(iii)] By Lemma~\ref{lemm:cut-and-paste}(b), and the diffeomorphism invariance of the mapping area,
\begin{equation}\label{eq:v2-area}
\begin{split}
A(\bw(t)) =\ & A(\bw(t); S \setminus \varphi(\bB_2)) + A(T(\rho(t), \widetilde{\bv}_1(t)); \bB_{2})\\
\leq\ & A(\bv_1(t); S \setminus \varphi(\bB_2)) + A(\widetilde{\bv}_1(t); \bB_{2}) + \frac{1}{8}\min\{W, \ep\}\\
=\ & A(\bv_1(t)) + \frac{1}{8}\min\{W, \ep\}.
\end{split}
\end{equation}
\end{enumerate}
From property (iii),~\eqref{eq:v_1-area-sup},~\eqref{eq:large-area-set-of-v1}, and the fact that $\bw = \bv_1$ outside of $\Inte(I_{3\delta/2}^m)$, we get 
\[
\sup_{t \in \Int(I^m)}A(\bw(t)) < W + \frac{\ep}{2},
\]
\[
\{t\in \Int(I^m)\ |\ A(\bw(t)) \geq \frac{3W}{8}\} \subset I_{5\delta}^m.
\]
To finish, observe that we have just verified (c) and (d), while conclusions (a) and (b) follow, respectively, from (i) and (ii) above.
\end{proof}
We now come to the main result of this section. For the statement we recall that $\ep_0$ is the energy threshold from Theorem \ref{thm:convexity}, $\rho_0 \in (0, \frac{1}{4})$ is the radius fixed at the start of Section~\ref{subsec:iterated-replacement}, and $N$ is given by Lemma \ref{lemm:besicovitch}.
\begin{prop}\label{prop:good-min-max-seq}
Again suppose $\bv_0 \in \cP$ is such that $W : = \cW([\bv_0]) > 0$. Given $\ep_1 \in (0, \frac{\ep_0}{3})$, there exist a sequence $(\sigma_n, v_n)$ in $\cM'$, and a sequence of positive numbers $\delta_n$ tending to zero, so that:
\vskip 1mm
\begin{enumerate}
\item[(a)] $\lim_{n \to \infty}A(v_n) = W$.
\vskip 1mm
\item[(b)] $\lim_{n \to \infty}\big[E(\sigma_n, v_n) - A(v_n) \big] = 0$.
\vskip 1mm
\item[(c)] For all $n$ and any finite, disjoint collection of geodesic disks $\fB = \{B_i\}$ in $(S, \sigma_n)$, each with radius at most $\rho_0$, satisfying
\begin{equation}\label{eq:good-min-max-seq-balls}
\sum_{i}\int_{B_i} |\nabla v_n|_{\sigma_n}^2 \vol_{\sigma_n} < \frac{\ep_1}{3^{N + 2}},
\end{equation}
we have 
\begin{equation}\label{eq:good-min-max-seq-convexity}
\begin{split}
&\int_{S} |\nabla v_{n} - \nabla \cR(\sigma_n, v_n, 2^{-2N-2}\fB)|_{\sigma_{n}}^2 \vol_{\sigma_{n}} \\
&\leq 4\big[ E(\sigma_n, v_n) - E(\sigma_n, \cR(\sigma_n, v_n, 2^{-2N-2}\fB)) \big] \leq \delta_n.
\end{split}
\end{equation}
\end{enumerate}
\end{prop}
\begin{proof}
Given any $\ep < \frac{W}{8}$,  we let $\bw \in [\bv_0]$ and $\delta \in (0, \frac{1}{20})$ be as given by Proposition~\ref{prop:area-min-seq}. Since $\bw(I_{\delta}^m)$ is a compact subset of $W^{1,2}(S; M)$, we have
\[
\Lambda: = \sup_{t \in I_\delta^m} E(\gamma_0, \bw(t)) < \infty.
\]
With $\ep$ and $\Lambda$ as above, Proposition~\ref{prop:conformal-reparametrization} yields some constant $\eta \in (0, 1)$ together with continuous maps $\Upsilon_{\eta}: C^2(S; M) \to \cM'$ and $\Xi_{\eta}: [0, 1] \times C^2(S; M) \to C^1(S; M)$ so that, upon letting
\begin{equation}\label{eq:energy-decreasing-sigma-defi}
(\bsig(t), \widetilde{\bw}(t)) = \Upsilon_{\eta}(\bw(t))\quad \text{for }t \in I_{2\delta}^m,
\end{equation}
we have that 
\begin{equation}\label{eq:w-hat-area-energy}
A(\bw(t)) \leq E(\bsig(t), \widetilde{\bw}(t)) < A(\bw(t)) + \frac{\ep}{4},\quad \text{for all }t \in I_{2\delta}^m,
\end{equation}
while $\Xi_{\eta}$ satisfies 
\begin{equation}\label{eq:Xi-and-hat}
\Xi_{\eta}(0, \bw(t)) = \bw(t), \quad \Xi_{\eta}(1, \bw(t)) = \widetilde{\bw}(t),\quad \text{for }t \in I_{2\delta}^m.
\end{equation}
Moreover, referring back to the definition~\eqref{eq:Xi-definition}, each $\Xi_{\eta}(s, \bw(t))$ consists of a diffeomorphism of $S$ to itself followed by $\bw(t)$. This has two consequences: First, as already noted in Proposition~\ref{prop:conformal-reparametrization},
\begin{equation}\label{eq:w-hat-area}
A(\Xi_{\eta}(s, \bw(t)))  =  A(\bw(t)),\quad\text{for all }(s, t) \in [0, 1] \times I_{2\delta}^m.
\end{equation}
Secondly, for all $t \in I_{2\delta}^m$, the map $\widetilde{\bw}(t)$ still is constant on some open subset of $S$ thanks to Proposition~\ref{prop:area-min-seq}(b). The unique continuation property of harmonic maps (see for instance \cite[Theorem 1]{Sampson1978}) then implies that, for all $t \in I_{2\delta}^m$,
\begin{equation}\label{eq:w-hat-no-harmonic}
\widetilde{\bw}(t):(S, \bsig(t)) \to M \text{ is not harmonic unless it is constant on $S$.}
\end{equation}
Next, take a continuous function $f:I^m \to [0, 1]$ such that
\[
f = 1 \text{ on }I_{4\delta}^m, \quad f = 0 \text{ on }I^m \setminus I_{3\delta}^m,
\]
and set 
\begin{equation}\label{eq:energy-decreasing-homotopy-defi}
\bh(s, t) = \left\{
\begin{array}{ll}
\Xi_{\eta}(f(t)s, \bw(t)), & \text{ if }(s, t) \in [0, 1] \times I_{2\delta}^m,\\
\bw(t), & \text{ if } (s, t) \in [0, 1] \times (I^m \setminus I_{2\delta}^m).
\end{array}
\right.
\end{equation}
Then $\bh$ maps $[0, 1] \times \Int(I^m)$ continuously into $(C^0 \cap W^{1, 2})(S; M)$. Since $\bh(s, t) = \bw(t)$ whenever $t \not\in I_{3\delta}^m$ or $s = 0$, we infer that
\begin{equation}\label{eq:w1-in-class}
\bw_1: = \bh(1, \cdot) \in [\bw] = [\bv_0].
\end{equation}
For all $t \in I^{m}_{4\delta}$, by~\eqref{eq:w-hat-area-energy}, the second property in~\eqref{eq:Xi-and-hat}, and our choice of $f$, we have
\begin{equation}\label{eq:w1-area-energy}
A(\bw(t)) \leq E(\bsig(t), \bw_1(t)) < A(\bw(t)) + \frac{\ep}{4}, \text{ for all }t \in I_{4\delta}^m.
\end{equation}
Also, with the help of~\eqref{eq:w-hat-area}, we have
\begin{equation}\label{eq:w1-area-agree}
A(\bw_1(t)) = A(\bw(t)), \text{ for all }t \in \Int(I^m).
\end{equation}
Finally, we extend the family of metrics $\{\bsig(t)\}_{t \in I_{2\delta}^m}$ to $I^m$ by using the radial retraction from $I^m \setminus I_{2\delta}^m$ to $\partial I_{2\delta}^m$. Then $\{(\bsig(t), \bw_1(t))\}_{t \in \Int(I^m)}$ is a continuous family in $\cM'$ in the sense defined prior to Lemma~\ref{lemm:besicovitch}.

We next verify that the family $\{(\bsig(t), \bw_1(t))\}_{t \in \Int(I^m)}$, together with $W = \cW([\bv_0])$ and the constant $\delta > 0$ from the beginning of the proof, satisfy the conditions (r1) and (r2) appearing in the hypotheses of Propositions~\ref{prop:choice-of-disks} and~\ref{prop:energy-decreasing}. Indeed, recall from~\eqref{eq:w1-in-class} that $\bw_1 \in [\bv_0]$, so that, by the definition of the width, there holds
\[
\sup_{t \in \Int(I^m)}A(\bw_1(t)) \geq W > 0.
\]
This, along with~\eqref{eq:w1-area-agree} and Proposition~\ref{prop:area-min-seq}(d), shows that (r1) holds. Next, given $t \in \Int(I^m)$ such that $A(\bw_1(t)) \geq \frac{W}{2} > 0$, by (r1), which we just checked, we must have $t \in I_{5\delta}^m$, in which case $\bw_1(t) = \widetilde{\bw}(t)$ by our choice of $f$ and~\eqref{eq:Xi-and-hat}. The previous area lower bound then prevents $\widetilde{\bw}(t)$ from being constant, so we get from~\eqref{eq:w-hat-no-harmonic} that
\[
\bw_1(t):(S, \bsig(t)) \to M \text{ is not a harmonic map.}
\]
This verifies (r2), and Proposition~\ref{prop:energy-decreasing} is applicable to $\{(\bsig(t), \bw_1(t))\}_{t \in \Int(I^m)}$, with ``$\ep$'' taken to be the given $\ep_1$ here. We denote the resulting new family by $\{(\bsig(t), \widehat{\bw}_1(t))\}$, and extend $\widehat{\bw}_1$ to $I^m$ by letting
\[
\widehat{\bw}_1(t) = \bw_1(t), \text{ for all }t \in \partial I^m.
\]
By~\eqref{eq:energy-decreasing-bc} and the continuity statement that precedes it, we see that $\widehat{\bw}_1 \in [\bw_1] =  [\bv_0]$, and that in fact
\[
\widehat{\bw}_1(t) = \bw_1(t), \text{ for all }t \in I^m \setminus \Int(I_{4\delta}^m).
\]
To prepare for the construction of the desired sequences, choose any $t_* \in \Int(I^m)$ such that 
\[
A(\widehat{\bw}_1(t_*)) > \sup_{t \in \Int(I^m)}A(\widehat{\bw}_1(t)) - \frac{\ep}{4} \geq W - \frac{\ep}{4} \geq \frac{W}{2},
\]
where the second inequality uses $\widehat{\bw}_1 \in [\bv_0]$, and the last inequality holds because $\ep < \frac{W}{8}$. If $t_* \not\in I_{4\delta}^m$, then $\widehat{\bw}_1(t_*) = \bw_1(t_*)$, which together with~\eqref{eq:w1-area-agree} gives
\[
\frac{W}{2} < A(\widehat{\bw}_1(t_*)) = A(\bw_1(t_*)) = A(\bw(t_*)).
\]
This forces $t_* \in I_{5\delta}^m$ by Proposition~\ref{prop:area-min-seq}(d), a contradiction. Thus $t_*$ must lie in $I_{4\delta}^m$. We can then estimate
\[
\begin{split}
W - \frac{\ep}{4} < A(\widehat{\bw}_1(t_*)) \leq\ & E(\bsig(t_*), \widehat{\bw}_1(t_*))\\
\leq\ & E(\bsig(t_*), \bw_1(t_*)) < A(\bw(t_*)) + \frac{\ep}{4} < W + \frac{3\ep}{4}.
\end{split}
\]
Here, the passage to the second line uses Proposition~\ref{prop:energy-decreasing}(a), and the last two inequalities uses, respectively,~\eqref{eq:w1-area-energy} and Proposition~\ref{prop:area-min-seq}(c). In particular, using again~\eqref{eq:w1-area-agree}, we have
\begin{equation}\label{eq:area-lower-bound-for-decreasing}
A(\bw_1(t_*)) = A(\bw(t_*)) > W - \frac{\ep}{2} > \frac{W}{2},
\end{equation}
that is, the area lower bound required in Proposition~\ref{prop:energy-decreasing}(b) holds at $t_*$.

Next we apply the above argument with $\ep = \frac{1}{n}$ for all sufficiently large $n$, and denote the various resulting objects in terms of $n$ as follows:
\[
\begin{split}
(\sigma_{n}, v_{n}): =\ & (\bsig(t_*), \widehat{\bw}_1(t_*)), \\
\delta_{n} :=\ & \Theta(E(\bsig(t_*), \bw_1(t_*)) - E(\bsig(t_*), \widehat{\bw}_1(t_*))),
\end{split}
\]
where $\Theta$ is the function given by Proposition~\ref{prop:energy-decreasing}. Then from the above string of inequalities we gather that
\[
|A(v_{n}) - W| < \frac{1}{n}, \quad E(\sigma_n, v_n) - A(v_n) < \frac{1}{n},
\]
so that conclusions (a) and (b) hold. Recalling also that $\Theta$ is increasing, continuous, and that $\Theta(0) = 0$, we have
\[
\delta_n \leq \Theta(\frac{1}{n}) \to 0 \text{ as }n \to \infty.
\]
Finally, given a finite, disjoint collection of geodesic disks $\fB$ satisfying the hypothesis of part (c), by~\eqref{eq:area-lower-bound-for-decreasing} and Proposition~\ref{prop:energy-decreasing}(b) we see that~\eqref{eq:good-min-max-seq-convexity} holds. The proof is complete.
\end{proof}
\subsection{Compactness properties of good min-max sequences}\label{subsec:compactness-properties}
We state without proof the following two standard results from the regularity theory of harmonic maps in dimension two.
\begin{thm}[$\ep$-regularity]
\label{thm:ep-regularity}
Given $\Lambda > 0$, there exists a positive constant $\ep_{\reg}$, depending only on $M$ and $\Lambda$, with the following property. Given a $C^1$-Riemannian metric $g$ on $\bB$ satisfying
\[
\Lambda^{-1}g_{\euc} \leq g \leq \Lambda g_{\euc},\quad |D g| \leq \Lambda, \text{ on }\bB,
\]
if $u \in W^{1, 2}(\bB, M)$ minimizes energy with respect to $g$ on compact subsets of $\bB$, and satisfies 
\[
\int_{\bB}|Du|^2< \ep_{\reg},
\]
then $u \in C^{1, \alpha}_{\loc}(\bB)$ for all $\alpha \in (0, 1)$, and 
\begin{equation}\label{eq:ep-reg-estimate}
\|Du\|_{\infty; \bB_{\frac{1}{2}}}^2 + [Du]_{\alpha; \bB_{\frac{1}{2}}}^2 \leq C\int_{\bB}|Du|^2,
\end{equation}
where $C$ depends only on $M, \Lambda$ and $\alpha$. If in addition $g$ is smooth on $\bB$, then so is $u$.
\end{thm}
\begin{thm}[energy gap]
\label{thm:harmonic-S2-gap}
Let $g_{S^2}$ denote the standard metric on $S^2$ with constant curvature $1$. There exists $\ep_{\text{gap}} >0$ depending only on $M$ such that if $u:(S^2, g_{S^2}) \to M$ is any non-constant harmonic map then
\[
E(g_{S^2}, u;S^2) \geq \ep_{\text{gap}}.
\]
\end{thm} 
Below we let $\ep_{\reg}$ be the threshold from Theorem~\ref{thm:ep-regularity} with $\Lambda = 2$, and recall that $\ep_0$ is given by Theorem \ref{thm:convexity}. Assuming that $\bv_0 \in \cP$ satisfies 
\[
W: = \cW([\bv_0]) > 0,
\]
we apply Proposition \ref{prop:good-min-max-seq} with $\ep_1$ taken to satisfy
\begin{equation}\label{eq:ep1-threshold}
0 < \ep_1 < \frac{1}{5}\min\{\ep_0, \ep_{\reg}, \ep_{\text{gap}}\},
\end{equation}
to obtain sequences $(\sigma_n, v_n) \in \cM'$ and $\delta_n > 0$ with the properties stated in the conclusion there. The proof of Theorem \ref{thm:main-existence}, carried out in the next section, amounts to showing that a subsequence of $(\sigma_n, v_n)$ bubble-tree converges. This section serves as preparation. In Lemma~\ref{lemm:ep-reg-almost-harmonic} and Proposition~\ref{prop:domain-convergence} we establish, based primarily on Proposition~\ref{prop:good-min-max-seq}(c), a ``small energy $\Rightarrow$ compactness'' result that we use to locate concentration points and extract bubbles. Then, Propositions~\ref{prop:neck} and~\ref{prop:collar-bubble} constitute our main tools for working with neck regions between bubbles. Certain assumptions common to both Propositions~\ref{prop:neck} and~\ref{prop:collar-bubble} are singled out as Definition~\ref{defi:well-prepared}.
\begin{lemm}[\cite{Zhou17b}, Lemma 5.8]
\label{lemm:ep-reg-almost-harmonic}
Let $\Omega$ be a two-dimensional open manifold, exhausted by an increasing sequence $(\Omega_{n})$ of open subsets. Suppose for each $n$ there exist $\lambda_n \in (0, 1]$ and a diffeomorphism  $\varphi_n$ from $\Omega_n$ onto an open subset of $S$ such that the metrics $g_n: = \lambda_{n}^{-2}\varphi_n^*\sigma_n$ converge smoothly locally on $\Omega$ to some limiting metric $g$. Letting $\rho_0$ be the radius fixed at the start of Section~\ref{subsec:iterated-replacement}, if $\overline{B_{g}(x, t)} \subset (\Omega, g)$ is a geodesic disk satisfying
\begin{equation}\label{eq:eventually-small-radius}
\limsup_{n \to \infty}\lambda_n t \leq \rho_0,
\end{equation}
\begin{equation}\label{eq:eventually-no-concentration}
\liminf_{n \to \infty}\int_{B_{g}(x, t)} |\nabla (v_n \circ \varphi_n)|_{g_n}^2 \vol_{g_n} < \frac{\ep_1}{3^{N + 2}},
\end{equation}
then a subsequence of $(v_n \circ \varphi_n)$ converges strongly in $W^{1, 2}$ on $B_{g}(x, \frac{t}{2^{2N + 4}})$ to a smooth harmonic map with respect to $g$.
\end{lemm}
\begin{proof}
We abbreviate $v_n \circ \varphi_n$ as $u_n$. Given $B_{g}(x, t)$ as in the statement, eventually $B_{g_n}(x, \frac{2}{3}t)$ is a geodesic disk with respect to $g_n$, and moreover
\begin{equation}\label{eq:radius-comparable-2}
B_{g}(x, \frac{r}{2}) \subset B_{g_n}(x, \frac{2}{3}r) \subset B_{g}(x, r), \text{ for all }r \leq t
\end{equation}
Since $\varphi_n: (\Omega_n, \lambda_n^2 g_n) \to (\varphi_n(\Omega_n), \sigma_n)$ is an isometry, we infer that $B_n : = B_{\sigma_n}(\varphi_n(x), \frac{2}{3}\lambda_n t)$ is a geodesic disk in $(S, \sigma_n)$. By~\eqref{eq:eventually-small-radius}, eventually we have $\frac{2}{3}\lambda_n t < \rho_0$. Moreover, by~\eqref{eq:eventually-no-concentration}, up to taking a subsequence we also have
\[
\int_{B_n} |\nabla v_n|_{\sigma_n}^2 \vol_{\sigma_n} \leq \int_{B_{g}(x, t)} |\nabla u_n|_{g_n}^2 \vol_{g_n} < \frac{\ep_1}{3^{N + 2}},\quad \text{for all large enough }n.
\]
Thus, letting $k_n$ be the harmonic replacement of $v_n$ on $2^{-2N-2}B_n$ with respect to $\sigma_n$, we deduce from Proposition~\ref{prop:good-min-max-seq}(c) that 
\[
\int_{2^{-2N-2}B_n} |\nabla v_n - \nabla k_n|_{\sigma_n}^2 \vol_{\sigma_n} \leq \delta_n.
\]
Next define 
\[
h_n = k_n \circ \varphi_n,
\]
which is an energy-minimizing harmonic map on $B_{g_n}(x, \frac{t}{3 \cdot 2^{2N + 1}})$ with respect to $g_n$, and satisfies eventually that 
\begin{align}
&\int_{B_{g_n}(x, \frac{t}{3 \cdot 2^{2N + 1}})} |\nabla h_n|^2_{g_n} \vol_{g_n}  \leq \int_{B_{g_n}(x, \frac{t}{3 \cdot 2^{2N + 1}})} |\nabla u_n|^2_{g_n} \vol_{g_n}  < \frac{\ep_1}{3^{N + 2}}, \label{eq:hn-small-energy}\\
&\int_{B_{g_n}(x, \frac{t}{3 \cdot 2^{2N + 1}})} |\nabla u_n - \nabla h_n|_{g_n}^2 \vol_{g_n} \leq \delta_n.\label{eq:W12-close}
\end{align}
Since $u_n - h_n = 0$ on $\partial B_{g_n}(x, \frac{t}{3 \cdot 2^{2N + 1}})$, using~\eqref{eq:radius-comparable-2} and the domain monotonicity of the lowest Dirichlet eigenvalue (of $-\Delta_g$), there exists a constant $C$ independent of $n$ such that
\[
\int_{B_{g_n}(x, \frac{t}{3 \cdot 2^{2N + 1}})} |u_n - h_n|^2 \vol_{g} \leq C\int_{B_{g_n}(x, \frac{t}{3 \cdot 2^{2N + 1}})} |\nabla u_n - \nabla h_n|_{g}^2 \vol_{g}.
\]
Estimating the right-hand side with~\eqref{eq:W12-close} and the convergence of $g_n$ to $g$, we get for all sufficiently large $n$ that
\begin{equation}\label{eq:L2-convergence}
\int_{B_{g_n}(x, \frac{t}{3 \cdot 2^{2N + 1}})} |u_n - h_n|^2 \vol_{g} \leq 4C\int_{B_{g_n}(x, \frac{t}{3 \cdot 2^{2N + 1}})} |\nabla u_n - \nabla h_n|_{g_n}^2 \vol_{g_n} \leq 4C\delta_n,
\end{equation}
Next, again by~\eqref{eq:radius-comparable-2} and the local smooth convergence of $g_n$ to $g$, we infer from \eqref{eq:hn-small-energy} and \eqref{eq:W12-close}, respectively, that the following two strings of estimates hold for all large enough $n$:
\begin{align}
&\int_{B_{g}(x, \frac{t}{2^{2N + 3}})} |\nabla h_n|_{g}^2 \vol_{g} \leq 4\int_{B_{g_n}(x, \frac{t}{3 \cdot 2^{2N + 1}})} |\nabla u_n|^2_{g_n} \vol_{g_n}  < \frac{\ep_1}{3^{N}} < \ep_{\reg}, \label{eq:hn-small-energy-g}\\
&\int_{B_{g}(x, \frac{t}{2^{2N + 3}})} |\nabla u_n - \nabla h_n|_{g}^2 \vol_{g} \leq 4\int_{B_{g_n}(x, \frac{t}{3 \cdot 2^{2N + 1}})} |\nabla u_n - \nabla h_n|_{g_n}^2 \vol_{g_n} \leq 4\delta_n. \label{eq:W12-close-g}
\end{align}
Now, if $\vartheta:\bB \to U \subset B_{g}(x, \frac{t}{2^{2N+3}})$ is an isothermal chart with respect to $g$ centered at some point in $B_{g}(x, \frac{t}{2^{2N + 4}})$, then eventually Theorem~\ref{thm:ep-regularity} is applicable to $h_n \circ \vartheta$, resulting in $C^{0, \alpha}$-bounds for $D(h_n\circ \vartheta)$ on $\bB_{\frac{1}{2}}$ that are independent of $n$. Combining this observation with a standard covering argument yields a subsequence of $(h_n)$ that converges in $C^{1}$ on $B_{g}(x, \frac{t}{2^{2N + 4}})$ to a smooth harmonic map $h$ with respect to $g$. By~\eqref{eq:W12-close-g} and~\eqref{eq:L2-convergence}, we deduce that
\[
u_n \to h\quad \text{strongly in }W^{1, 2} \text{ on }B_{g}(x, \frac{t}{2^{2N + 4}}),
\]
which concludes the proof.
\end{proof}
\begin{prop}\label{prop:domain-convergence}
Under the same hypotheses as in Lemma~\ref{lemm:ep-reg-almost-harmonic}, there exist a finite set $\cS \subset \Omega$, along with a smooth, weakly conformal harmonic map $u: (\Omega, g) \to M$ with finite energy, such that, up to taking a subsequence, we have
\[
u_n : = v_n\circ \varphi_n \to u \text{ in }W^{1, 2}_{\loc}(\Omega\setminus \cS; M),
\]
and that, as Radon measures on $\Omega$, 
\[
\frac{1}{2}|\nabla u_n|_{g_n}^2 \vol_{g_n} \to \frac{1}{2}|\nabla u|_{g}^2 \vol_g + \sum_{x \in \cS}m_x \delta_x \quad\text{in weak$^*$ sense},
\]
where $\delta_x$ denotes the delta measure supported at $x$, and each $m_x$ is at least $\frac{\ep_1}{2 \cdot 3^{N + 2}}$.
\end{prop}
\begin{proof}
By Proposition~\ref{prop:good-min-max-seq}(a)(b) and the conformal invariance of the Dirichlet energy, we have
\[
E(g_n, u_n; \Omega_n) = E(\sigma_n, v_n; \varphi_n(\Omega_n)) \leq W + o(1).
\]
The desired conclusion, except for the weakly conformal property, now follows by standard arguments from this uniform energy bound and Lemma~\ref{lemm:ep-reg-almost-harmonic}, as well as the removable singularity theorem for finite-energy harmonic maps~\cite{Sacks-Uhlenbeck81}. To see that $u$ is weakly conformal, take any compact set $K \subset \Omega \setminus \cS$ and note that as soon as $n$ is large enough so that $K \subset \Omega_n$, we have 
\[
\begin{split}
E(g_n, u_n; K) - A(u_n; K) =\ & E(\sigma_n, v_n; \varphi_n(K)) - A(v_n; \varphi_n(K))\\
\leq\ & E(\sigma_n, v_n) - A(v_n).
\end{split}
\]
Letting $n \to \infty$ and using Proposition~\ref{prop:good-min-max-seq}(b) gives
\[
E(g, u; K) \leq A(u; K), \quad \text{for all compact }K \subset \Omega \setminus \cS.
\]
From this we infer that $u$ is weakly conformal on $\Omega \setminus \cS$, and thus on $\Omega$ by smoothness.\\
\end{proof}
Next we introduce some additional notation. Given $a < b$, we define $\bC_{a, b} = [a, b] \times S^1$. By $g_{\pro}$, we mean the product metric on $\RR \times S^1$; that is
\[
g_{\pro} : =  dt^2 + d\theta^2.
\]
For $s \in \RR$, we write $\boldsymbol{\tau}_{s}$ for the translation map
\begin{equation}\label{eq:tau-translation}
\begin{split}
\boldsymbol{\tau}_{s}: \RR \times S^1 \longrightarrow\ & \RR \times S^1\\
(t, \theta) \longmapsto\ & (t + s, \theta),
\end{split}
\end{equation}
which is closely related to the dilation maps
\[
\begin{split}
\boldsymbol{\delta}_{x, r}: \RR^2 \longrightarrow\ & \RR^2\\
y \longmapsto \ & x + ry.
\end{split}
\]
Specifically, denoting by $\mathbf{c}:\RR \times S^1 \to \RR^2 \setminus \{0\}$ the map taking $(t, \theta)$ to the point whose polar coordinates are $(e^{t}, \theta)$, we have
\begin{equation}\label{eq:c-d-t-relation}
\mathbf{c} \circ \boldsymbol{\tau}_{s} = \boldsymbol{\delta}_{0, e^{s}} \circ \mathbf{c}.
\end{equation}

Following standard practice, in identifying further bubble regions, as well as the neck regions between them, we conformally transform the relevant domains into cylinders of the kind described in Definition~\ref{defi:well-prepared} below. Two types of conformal transformations are employed. One of them is explained in Lemma~\ref{lemm:neck-assumptions}. The other, used specifically in collar regions where the metrics $\sigma_n$ are degenerating, appears in Section~\ref{subsec:convergence} (see especially \eqref{eq:collar-isometry}).
\begin{defi}\label{defi:well-prepared}
Given a positive sequence $(T_n)$ with $\lim_{n \to \infty}T_n = \infty$ and, for each $n$, a diffeomorphism $\psi_n$ from $\bC_{-T_n, T_n}$ into $S$, a positive number $\lambda_n$, and a smooth function $f_n:[-T_n, T_n] \to (0, \infty)$, we say that the sequence $\{(\bC_{-T_n, T_n}, \psi_n, \lambda_n, f_n)\}$ is \emph{well-prepared} if the following hold:
\vskip 1mm
\begin{enumerate}
\item[(p1)] $\limsup_{n \to \infty}\lambda_nT_n < \infty$, while
\[
\sup_{|t| \leq T_n}\frac{\lambda_n}{f_n(t)} \to 0 \quad \text{as }n \to \infty.
\]
\vskip 1mm
\item[(p2)] Let $g_n = \lambda_n^{-2}\psi_n^*\sigma_n$. There exists a positive sequence $(\alpha_n)$ converging to $0$ such that, for all sufficiently large $n$,
\[
(1 + \alpha_n)^{-1}  \cdot g_{\pro} \leq f_n^2 \cdot g_n \leq (1 + \alpha_n) \cdot g_{\pro}, \quad \text{ on }\bC_{-T_n, T_n}.
\]
\vskip 1mm
\item[(p3)] With $g_n$ as in (p2), for all $T > 0$ and any sequence $t_n \in [-T_n + T, T_n - T]$, we have
\[
(f_n(t_n))^2 \cdot \boldsymbol{\tau}_{t_n}^* g_n \to g\quad \text{smoothly on }\bC_{-T, T} \text{ as }n \to \infty,
\]
where $g$ is some Riemannian metric on $\bC_{-T, T}$ which is conformal to $g_{\pro}$.
\end{enumerate}
\end{defi}
\begin{lemm}\label{lemm:neck-assumptions}
Given $R > 0$, suppose for each $n$ there exist $a_n \in (0, 1]$ and a diffeomorphism $\varphi_n$ from $\bB_R \subset \RR^2$ into $S$ such that $h_n:= a_n^{-2}\varphi_n^*\sigma_n$ converge smoothly locally on $\bB_R$ to a limiting metric $h$ satisfying $h(0) = g_{\euc}$. Given a sequence $(y_n)$ in $\bB_R$ converging to the origin, and sequences of radii $(\rho_n), (s_n)$ satisfying
\begin{equation}\label{eq:neck-assumptions-scales}
\lim_{n \to \infty}\rho_n = \lim_{n \to \infty}\frac{s_n}{\rho_n} = 0,
\end{equation}
define $T_n = \log\sqrt{\frac{\rho_n}{s_n}}$, $\mu_n = \sqrt{\rho_n s_n}$, and let
\[
\lambda_n = \mu_n a_n, \quad \psi_n = \varphi_n \circ \boldsymbol{\delta}_{y_n, \mu_n} \circ \mathbf{c}, \quad f_n(t) = e^{-t}.
\]
Then $\{(\bC_{-T_n, T_n}, \psi_n, \lambda_n, f_n)\}$ is well-prepared.
\end{lemm}
\begin{proof}
By~\eqref{eq:neck-assumptions-scales} we have $\lim_{n \to\infty}T_n = \infty$. In particular, eventually 
\[
\lambda_n T_n + \sup_{|t| \leq T_n}\frac{\lambda_n}{f_n(t)} \leq 2\lambda_n e^{T_n} \leq 2\mu_n e^{T_n} = 2\rho_n.
\]
Since $\rho_n \to 0$, we see that (p1) in Definition \ref{defi:well-prepared} holds. Next, by the definition of $\psi_n$ and $h_n$ we have
\begin{equation}\label{eq:neck-assumptions-gn-hn}
g_n: = \lambda_n^{-2}\psi_n^*\sigma_n = \mathbf{c}^*\big( \mu_n^{-2}\boldsymbol{\delta}_{y_n,\mu_n}^*h_n \big).
\end{equation}
For all $x \in \RR^2$ such that  $e^{-T_n} \leq |x| \leq e^{T_n}$, using the assumption $h(0) = g_{\euc}$, we find that
\begin{equation}\label{eq:dilation-convergence}
\begin{split}
\big|(\mu_n^{-2}\boldsymbol{\delta}_{y_n, \mu_n}^{*}h_{n})_{ij}(x) - \delta_{ij}\big| =\ & \big|(h_n)_{ij}(y_n + \mu_n x) - h_{ij}(0)\big| \\
\leq\ & \sup_{y \in \bB_{\rho_n}(y_n)}|(h_n)_{ij}(y) - h_{ij}(0)|,
\end{split}
\end{equation}
Since $\rho_n \to 0$ and $y_n \to 0$, and since $h_n \to h$ uniformly on compact subsets of $\bB_R$, we infer that 
\[
\mu_n^{-2}\boldsymbol{\delta}_{y_n, \mu_n}^{*}h_{n} \to g_{\euc}\quad \text{uniformly on }\bB_{e^{T_n}} \setminus \bB_{e^{-T_n}}.
\]
In particular there exists a positive sequence $\alpha_n \to 0$ such that eventually
\[
(1 + \alpha_n)^{-1}\cdot g_{\euc} \leq \mu_n^{-2}\boldsymbol{\delta}_{y_n, \mu_n}^{*}h_{n} \leq (1 + \alpha_n) \cdot g_{\euc}, \quad \text{on }\bB_{e^{T_n}} \setminus \bB_{e^{-T_n}}.
\]
Pulling back by $\mathbf{c}$ and recalling that $\mathbf{c}^*g_{\euc} = e^{2t}g_{\pro} = (f_n(t))^{-2}g_{\pro}$, we get
\[
(1 + \alpha_n)^{-1}\cdot g_{\pro} \leq (f_n(t))^2\cdot g_n \leq (1 + \alpha_n) \cdot g_{\pro},\quad \text{on }\bC_{-T_n, T_n}.
\]
which verifies (p2). Finally, given $T > 0$ and a sequence $t_n\in [-T_n + T, T_n - T]$, by~\eqref{eq:c-d-t-relation} and~\eqref{eq:neck-assumptions-gn-hn} we have
\[
(f_n(t_n))^2 \cdot \boldsymbol{\tau}_{t_n}^* g_n = \mathbf{c}^*\big((\mu_ne^{t_n})^{-2} \cdot  \boldsymbol{\delta}_{y_n, \mu_ne^{t_n}}^*h_n\big).
\]
By a consideration similar to~\eqref{eq:dilation-convergence}, and noting also that $\mu_n e^{t_n + T} \leq \mu_n e^{T_n} = \rho_n \to 0$, it is not hard to deduce from the above relation that $(f_n(t_n))^2 \cdot \boldsymbol{\tau}_{t_n}^* g_n$ converges smoothly on $\bC_{-T, T}$ to $\mathbf{c}^*g_{\euc} = e^{2t}g_{\pro}$. This establishes (p3). 
\end{proof}

For the next result, based on ideas from \cite[Appendix B]{Colding-Minicozzi08b} and \cite[Section 5]{Zhou10}, we need to fix two additional parameters. By~\cite[Proposition B.1 and Remark B.3]{Colding-Minicozzi08b}, there exist constants $\ep_{\text{cyl}} \in (0, 1)$ and $L \geq 1$ depending only on $M$ such that if $v$ is either a harmonic map from $(\bC_{-3L, 3L}, g_{\pro})$ into $M$ satisfying $E(g_{\pro}, v; \bC_{-3L, 3L}) \leq \ep_{\text{cyl}}$, or a harmonic function from $(\bC_{-3L, 3L}, g_{\pro})$ into $\RR^q$, then we necessarily have
\begin{equation}\label{eq:CM-estimate}
\int_{\bC_{-L, L}} |v_{\theta}|^2 \leq \frac{1}{31}\int_{\bC_{-2L, 2L}}|\nabla v|^2.
\end{equation}
In what follows we fix $\ep_2$ so that
\begin{equation}\label{eq:ep2-choice}
0 < \ep_2 < \min\{\frac{\ep_1}{2\cdot 3^{N + 2}},  \ep_{\text{cyl}} \}.
\end{equation}
\begin{prop}[\cite{Colding-Minicozzi08b}, Appendix B; \cite{Zhou10}, Section 5]
\label{prop:neck}
Suppose $\{(\bC_{-T_n, T_n}, \psi_n, \lambda_n, f_n)\}$ is well-prepared according to Definition~\ref{defi:well-prepared}. Define
\[
g_n = \lambda_n^{-2} \psi_n^*\sigma_n, \quad u_n = v_n \circ \psi_n,
\]
and assume that for all $T > 0$, we have 
\begin{equation}\label{eq:no-energy-on-finite-part}
\sup_{-T_n + T \leq t \leq T_n - T} E(g_{n}, u_n; \bC_{t-T, t+T}) \to 0 \text{ as }n \to \infty.
\end{equation}
Then in fact
\begin{equation}\label{eq:no-energy-on-neck}
\lim_{n \to \infty}E(g_{n}, u_n; \bC_{-T_n, T_n}) = 0.
\end{equation}
\end{prop}
\begin{proof}
As is the case with Proposition~\ref{prop:e-semicontinuous}, the ideas behind the proof come mostly from the work of Colding--Minicozzi~\cite{Colding-Minicozzi08b}, but carrying them out requires some care to details. We refer the reader to Appendix~\ref{appendix:long-cylinders}.
\end{proof}
\begin{prop}\label{prop:collar-bubble}
Again suppose $\{(\bC_{-T_n, T_n}, \psi_n, \lambda_n, f_n)\}$ is well-prepared according to Definition~\ref{defi:well-prepared} and let 
\[
g_n = \lambda_n^{-2} \psi_n^*\sigma_n, \quad u_n = v_n \circ \psi_n.
\]
Assume that 
\begin{equation}\label{eq:iv-weakened}
\lim_{n \to \infty}E(g_{n}, u_n; \bC_{-T_n, -T_n + T} \cup \bC_{T_n - T, T_n })  = 0, \textit{ for all }T > 0,
\end{equation}
and that
\begin{equation}\label{eq:there-is-energy}
\limsup_{n \to \infty}\sup_{t \in [-T_n + 1, T_n - 1]}E(g_{n}, u_n; \bC_{t-1, t+1}) > 0.
\end{equation}
Then there exists $Q_1 \in \NN$ so that, up to taking a subsequence, we can find, for each $n$, concentric closed intervals $I_{n, i}' \subset I_{n, i}$ in $[-T_n, T_n]$ for $i = 1, \cdots, Q_1$ such that $I_{n, 1}, \cdots, I_{n, Q_1}$ are mutually disjoint, and that the following hold.
\vskip 1mm
\begin{enumerate}
\item[(a)] For each $i$, 
\[
\lim_{n \to \infty}\diam I'_{n, i} = \lim_{n \to \infty} \dist(I'_{n, i}, \partial I_{n, i}) = \lim_{n \to \infty}\dist(I_{n, i}, \{-T_n, T_n\}) = \infty.
\]
\vskip 1mm
\item[(b)] Whenever $i \neq j$, 
\[
\lim_{n \to \infty} \dist(I_{n, i}, I_{n, j}) = \infty.
\]
\vskip 1mm
\item[(c)] For each $i$, we have
\[
\lim_{n \to \infty}E(g_{n}, u_n; (I_{n, i}\setminus I_{n, i}') \times S^1) = 0.
\]
\vskip 1mm
\item[(d)] Writing $t_{n, i}$ for the midpoint of $I_{n, i}$, there exists a Radon measure $\mu_{i}$ on $\RR \times S^1$, with positive and finite total mass, such that 
\[
\boldsymbol{\tau}_{t_{n, i}}^*\big(\frac{1}{2}|\nabla u_{n}|_{g_{n}}^2 \vol_{g_{n}}\big) \to \mu_{i}\quad\text{in weak$^*$ sense on }\RR \times S^1,
\]
and that
\[
\lim_{n \to \infty}E(g_{n}, u_n; I'_{n, i} \times S^1) = \mu_i(\RR \times S^1).
\]
\vskip 1mm
\item[(e)] $\sup \big\{ E(g_{n}, u_n; I \times S^1) \ |\  I \subset [-T_n, T_n]\setminus (\cup_{i=1}^{Q_1}I_{n, i}), \diam I = 2 \big\} \to 0$ as $n \to \infty$.
\end{enumerate} 
\end{prop}
\begin{proof}
The idea of proof is entirely standard. One way of executing the argument is included in Appendix~\ref{appendix:long-cylinders} for the convenience of the reader.
\end{proof}
\begin{rmk}\label{rmk:no-leftover}
Under the assumptions of Proposition~\ref{prop:collar-bubble}, we may apply Proposition~\ref{prop:neck} to each connected component of $[-T_n, T_n] \setminus \big(\cup_{i = 1}^{Q_1}I_{n, i} \big)$. Specifically, denote the closures of these components by 
\[
[s_{n,i} - L_{n,i}, s_{n, i} + L_{n, i}],\quad i = 0, \cdots, Q_1,
\]
and fix some $i \in \{0, \cdots, Q_1\}$. By Proposition~\ref{prop:collar-bubble}(a)(b), we have $\lim_{n \to \infty}L_{n, i} = \infty$, and it is not hard to check that the sequence
\[
\{(\bC_{-L_{n, i}, L_{n, i}},\ \psi_n \circ \boldsymbol{\tau}_{s_{n, i}},\ \lambda_n,\ f_n(s_{n, i} + \cdot))\}_{n \in \NN}
\]
is well-prepared. Moreover, given $T > 0$, iterating Proposition~\ref{prop:collar-bubble}(e) shows that the assumption~\eqref{eq:no-energy-on-finite-part} holds with $T_n$ replaced by $L_{n, i}$, and with $g_n, u_n$ replaced by their pullbacks with respect to $\boldsymbol{\tau}_{s_{n, i}}$. Thus we can invoke Proposition~\ref{prop:neck} to see that
\begin{equation}\label{eq:no-leftover}
\lim_{n \to \infty}E(g_n, u_n; \bC_{s_{n,i} - L_{n,i}, s_{n, i} + L_{n, i}}) = 0, \quad \text{for }i = 0, 1, \cdots, Q_1,
\end{equation}
which combines with Proposition~\ref{prop:collar-bubble}(c)(d) to give
\begin{equation}\label{eq:transition-energy}
\lim_{n \to \infty}E(g_n, u_n; \bC_{-T_n, T_n}) = \sum_{i = 1}^{Q_1}\mu_i(\RR \times S^1).
\end{equation}
\end{rmk}
\begin{rmk}\label{rmk:form-of-mu-i}
Continuing to work under the hypotheses of Proposition~\ref{prop:collar-bubble}, we observe the following about the limiting measures $\mu_{i}$ in conclusion (d). Fixing $i \in \{1, \cdots, Q_1\}$ and denoting by $t_{n, i}$ the midpoint of $I_{n, i}$, we let
\[
\begin{split}
\widehat{u}_n :=\ & u_n \circ \boldsymbol{\tau}_{t_{n, i}} = v_n \circ (\psi_n \circ \boldsymbol{\tau}_{t_{n, i}}),\\
\widehat{g}_n :=\ & (f_n(t_{n, i}))^2 \cdot \boldsymbol{\tau}_{t_{n, i}}^* g_n = \big( \frac{\lambda_n}{f_n(t_{n, i})} \big)^{-2} \cdot (\psi_n \circ \boldsymbol{\tau}_{t_{n, i}})^*\sigma_n.
\end{split}
\]
These are defined on domains that exhaust $\RR \times S^1$, since $\diam I_{n, i} \to \infty$. Also, by (p1) and (p3) of Definition~\ref{defi:well-prepared}, we have that $\lim_{n \to \infty}\frac{\lambda_n}{f_n(t_{n, i})} = 0$, and that $\widehat{g}_n$ converges in $C^{\infty}_{\loc}(\RR \times S^1)$ to a limiting metric conformal to $g_{\pro}$. Noting that 
\[
\boldsymbol{\tau}_{t_{n, i}}^*\big( \frac{1}{2}|\nabla u_n|_{g_n}^2 \vol_{g_n} \big) = \frac{1}{2}|\nabla \widehat{u}_n|_{\widehat{g}_n}^2 \vol_{\widehat{g}_n},
\]
we obtain from Proposition~\ref{prop:domain-convergence} a finite-energy, weakly conformal harmonic map $w_i: (\RR \times S^1, g_{\pro}) \to M$, and a finite set $\cA_i \subset \RR \times S^1$, such that
\begin{equation}\label{eq:form-of-mu-i}
\mu_{i} =\frac{1}{2} |\nabla w_i|_{g_{\pro}}^2 \vol_{g_{\pro}} + \sum_{x \in \cA_{i}} m_{i, x}\delta_{x},
\end{equation}
where $m_{i, x} \geq \frac{\ep_1}{2 \cdot 3^{N + 2}} > \ep_2$ for each $x \in \cA_{i}$. Moreover, by the removable singularity theorem~\cite{Sacks-Uhlenbeck81} and Theorem~\ref{thm:harmonic-S2-gap}, we have $E(g_{\pro}, w_i) \geq \ep_{\text{gap}} > \ep_2$, provided $w_i$ is not a constant.
\end{rmk}


\subsection{Bubble tree convergence}\label{subsec:convergence}
We recall certain elements of the Deligne--Mumford compactification, or rather its differential geometric description, following the presentation in~\cite[Appendix A]{Rupflin-Topping-Zhu2013}. A more comprehensive reference is~\cite{Hummel}. We do not reproduce the full statement of the compactness theorem, but only mention what for us are the relevant consequences. With $(\sigma_n, v_n)$ as in the previous section, up to taking a subsequence, exactly one of the following two alternatives occur. 
\vskip 2mm
\noindent\textbf{Case 1:} There exist orientation-preserving diffeomorphisms $\varphi_n: S \to S$ such that $\varphi_n^* \sigma_n$ converges smoothly on $S$ to a limiting metric $\sigma$. 
\vskip 2mm
\noindent\textbf{Case 2:} There exist:
\vskip 1mm
\begin{enumerate}
\item[(i)] a positive integer $Q \in \{1, 2, \cdots, 3 \cdot \text{genus}(S) - 3\}$;
\vskip 1mm
\item[(ii)] a closed, possibly disconnected Riemann surface $(\widehat{S}, [\hat{h}])$, along with a collection $\cC = \{p_{i}^+, p_{i}^{-}\}_{i = 1}^{Q}$ of $2Q$ points in $\widehat{S}$, and a complete hyperbolic metric $h$ on $S': = \widehat{S} \setminus \cC$ which is conformal to the restriction of $\hat{h}$;
\vskip 1mm
\item[(iii)] for each $n$, a collection $\{\gamma_{n, i}\}_{i = 1}^{Q}$ of mutually disjoint, homotopically non-trivial, simple closed geodesics in $(S, \sigma_n)$, along with an orientation-preserving diffeomorphism 
\[
F_n: S\setminus \cup_{i = 1}^{Q}\gamma_{n, i}\to S'
\]
which extends to a continuous map from $S$ into the space obtained from $\widehat{S}$ by identifying each $p_{i}^+$ with $p_i^{-}$, in such a way that $F_n(\gamma_{n, i}) = \{p_i^+, p_i^-\}$ for $i =1, \cdots, Q$.
\end{enumerate}
\vskip 1mm
These objects have the following further properties:
\vskip 1mm
\begin{enumerate}
\item[(c1)] Let $l_{n, i}$ denote the length of $\gamma_{n, i}$ with respect to $\sigma_n$. Then $\lim_{n \to \infty}l_{n, i} = 0$ for all $i$.
\vskip 1mm
\item[(c2)] Let $G_n = F_n^{-1}$. Then $G_n^*\sigma_n \to h$ smoothly locally on $S'$.
\vskip 1mm
\item[(c3)] (See~\cite[equation (A.5)]{Rupflin-Topping-Zhu2013}) For all $d \in (0, \sinh^{-1}(1))$, define 
\[
\Sigma^{d} = \{x \in S'\ |\ \inj(x, h) \geq d\}, \quad\quad \Sigma_n^{d} = \{x \in S\ |\ \inj(x, \sigma_n) \geq d\}.
\]
Then, given $0 < d_1 < d < d_2 < \sinh^{-1}(1)$, we have $\Sigma^{d_2} \subset G_n^{-1}(\Sigma_n^{d}) \subset \Sigma^{d_1}$, for all sufficiently large $n$.
\end{enumerate}
\begin{rmk}
By (c1), and the standard characterization of the $\sinh^{-1}(1)$-thin part of hyperbolic surfaces (see~\cite[Lemmas A.4 and A.5 and Proposition A.6]{Rupflin-Topping-Zhu2013} for a summary), the $\sigma_n$-distance between $\Sigma_n^{d}$ and $\cup_{i = 1}^{Q}\gamma_{n, i}$ tends to infinity as $n \to \infty$ while $d \in (0, \sinh^{-1}(1))$ is kept fixed. It follows that eventually $G_n^{-1}(\Sigma_n^{d})$ coincides with $\{x \in S'\ |\ \inj(x, G_n^*\sigma_n) \geq d\}$, and hence (c3) does follow from~\cite[equation (A.5)]{Rupflin-Topping-Zhu2013}.
\end{rmk}
We now restate Theorem~\ref{thm:main-existence} in a more precise form.
\begin{thm*}
Suppose $\bv_0 \in \cP$ is such that $W: = \cW([\bv_0]) > 0$, and let $(\sigma_n, v_n)$ be the sequence given by Proposition~\ref{prop:good-min-max-seq}, with $\ep_1$ satisfying~\eqref{eq:ep1-threshold}. Take $(S_0, c_0)$ to be $(S, [\sigma])$ in Case 1, and $(\widehat{S}, [\hat{h}])$ in Case 2. Up to taking a subsequence, there exist a weakly conformal harmonic map $u_0:(S_0, c_0) \to M$ and finitely many harmonic $2$-spheres $\{w_j\}_{j\in J}$ in $M$ such that 
\begin{equation}\label{eq:energy-identity}
(W = \ )\ \lim_{n \to \infty}E(\sigma_n, v_n) = E(h_0, u_0; S_0) + \sum_{j\in J}E(g_{S^2}, w_j),
\end{equation}
where $h_0$ is any representative of the conformal class $c_0$.
\end{thm*}
We only explain the proof in Case 2. The other case is essentially contained in it. We do not give a full account of how to construct the bubble tree or locate neck regions, as these are completely standard. Instead we focus the role played by what might be called the ``almost conformal'' (Proposition~\ref{prop:good-min-max-seq}(b)) and ``almost harmonic'' (Proposition~\ref{prop:good-min-max-seq}(c)) properties of the good min-max sequence $(\sigma_n, v_n)$. 
\vskip 1em
\noindent\textit{Step 1: Base map.} By Proposition~\ref{prop:domain-convergence} with $\varphi_n = G_n$, $\Omega_n = \Omega = S'$ and $\lambda_n = 1$, we get a finite set $\cB \subset S'$, and a smooth, weakly conformal harmonic map $u_0: (S', h) \to M$ with finite energy, such that, up to taking a subsequence, we have
\begin{equation}\label{eq:level-0-map}
u_n: = v_n \circ G_n \to u_0 \text{ in }W^{1, 2}_{\loc}(S' \setminus \cB),
\end{equation}
and that, with $h_n := G_n^* \sigma_n$, 
\begin{equation}\label{eq:level-0-measure}
\frac{1}{2}|\nabla u_n|_{h_n}^2 \vol_{h_n} \to \frac{1}{2}|\nabla u_0|^2_{h} \vol_{h} + \sum_{x \in \cB} m_{x} \delta_{x},
\end{equation}
as Radon measures on $S'$, where 
\begin{equation}\label{eq:mx-prelim-lower-bound}
m_{x} \geq \frac{\ep_1}{2\cdot 3^{N + 2}} > \ep_2, \quad\text{for all }x \in \cB,
\end{equation}
the second inequality being a consequence of~\eqref{eq:ep2-choice}. By point (ii) above, we can regard $u_0$ as a weakly conformal harmonic map on $(\widehat{S}\setminus \cC, \widehat{h})$ with finite energy. The removable singularity theorem then implies that $u_0$ extends smoothly to $\widehat{S}$, and further we have
\begin{equation}\label{eq:base-map-energy}
E(\widehat{h}, u_0; \widehat{S}) = E(h, u_0; S').
\end{equation}
Next, given any open set $U$ with compact closure in $S'$, with the help of Proposition~\ref{prop:good-min-max-seq}(a)(b) we have
\[
E(h_n, u_n; U) = E(\sigma_n, v_n; G_n(U)) \leq W + o(1),
\]
From this and~\eqref{eq:level-0-measure}, we deduce that
\begin{equation}\label{eq:loss-at-collars}
\tau: = W - \big(  E(h, u_0; S') + \sum_{x \in \cB} m_{x} \big) \geq 0.
\end{equation}
It remains to analyze what $m_{x}$ and $\tau$ are composed of. These are done in Steps 2 to 3 and Steps 4 to 5, respectively. Throughout the argument we make frequent use of the maps introduced before Definition~\ref{defi:well-prepared}.
\vskip 1em
\noindent\textit{Step 2: Top bubbles and transition regions.} 
If $\cB = \{x_1, \cdots, x_N\}$ is non-empty, we choose $\rho_1 > 0$ so that
\[
\rho_1 < \frac{1}{8}\inj(x_{i}, h), \text{ for all }i\ ;\quad \overline{B_{h}(x_i, 4\rho_1)} \cap \overline{B_{h}(x_j, 4\rho_1)} = \emptyset, \text{ whenever }i \neq j.
\]
Fixing some $i \in \{1, \cdots, N\}$ and an orthonormal basis for $T_{x_i}S'$, we use the exponential map $\exp_{x_i}^{h}$ to identify $B_h(x_i, 4\rho_1)$ with $\bB_{4\rho_1} \subset \RR^2$, and drop all explicit mention of $\exp_{x_i}^{h}$ in what follows. In particular we denote $(\exp_{x_i}^{h})^*h_n$ and $(\exp_{x_i}^{h})^*h$ still by $h_n$ and $h$, and similarly for $u_n$ and $u_0$. The metrics $h_n$, understood this way, then satisfy
\begin{equation}\label{eq:hn-converge-for-top-bubble}
h_n \to h\text{ in }C^{\infty}(\bB_{4\rho_1}),\quad h_{ij}(0) = \delta_{ij}.
\end{equation}
By the standard procedure for extracting top bubbles (see for example~\cite{Brezis-Coron1985} or~\cite{Wang-houston1996}), and in view of the strict inequality~\eqref{eq:mx-prelim-lower-bound}, we can find a subsequence of $(u_n, h_n)$, along with $y_n \to 0$ and $r_n \to 0$, such that
\begin{equation}\label{eq:top-bubble-region}
\ep_2 = E(h_n, u_n; \bB_{r_n}(y_n)) \geq E(h_n, u_n; \bB_{r_n}(y)), \quad \text{for all }y \in \bB_{2\rho_1}.
\end{equation}
Note that eventually $\bB_{\frac{3\rho_1}{2}}(y_n) \subset \bB_{2\rho_1}$, and that
\begin{equation}\label{eq:blown-up-to-g-euc}
r_n^{-2}\boldsymbol{\delta}_{y_n, r_n}^*h_n \to g_{\euc}\quad\text{in }C^{\infty}_{\loc}(\RR^2).
\end{equation}
Combining these observations with~\eqref{eq:top-bubble-region}, and recalling that $\ep_2 < \frac{\ep_1}{3^{N + 2}}$ by~\eqref{eq:ep2-choice}, we may apply Lemma~\ref{lemm:ep-reg-almost-harmonic} with the following choices
\[
\Omega =\RR^2, \quad  \Omega_n =\bB_{\frac{3\rho_1}{2r_n}} ,\quad \lambda_n = r_n, \quad \varphi_n = G_n \circ \boldsymbol{\delta}_{y_n, r_n}, \quad t = 1,
\]
to obtain a harmonic map $w_{i}: (\RR^2, g_{\euc}) \to M$ so that, along a further subsequence,
\begin{equation}\label{eq:level-0-top-bubble}
u_n \circ \boldsymbol{\delta}_{y_n, r_n} \to w_i \text{ strongly in }W^{1, 2}_{\loc}(\RR^2).
\end{equation}
With the help of this convergence and the one in~\eqref{eq:level-0-measure}, we have
\begin{equation}\label{eq:wi-energy-upper}
E(g_{\euc}, w_{i}) \leq m_{x_i}.
\end{equation}
Thus, since~\eqref{eq:blown-up-to-g-euc},~\eqref{eq:level-0-top-bubble}, together with the equality in~\eqref{eq:top-bubble-region} prevents $w_{i}$ from being constant, we deduce from the removable singularity theorem and Theorem~\ref{thm:harmonic-S2-gap} that
\begin{equation}\label{eq:wi-energy-lower}
E(g_{\euc}, w_i) \geq \ep_{\text{gap}}.
\end{equation}
By another standard argument based on~\eqref{eq:level-0-measure},~\eqref{eq:level-0-top-bubble}, as well as the fact that $|y_n| + r_n\to 0$, we get a still further subsequence of $(u_n, h_n)$, along with radii $(\rho_n)$ and $(s_n)$, so that the following hold:
\begin{equation}\label{eq:level-0-transition-scales}
\lim_{n \to \infty}\frac{r_n}{s_n} = \lim_{n \to \infty}\frac{s_n}{\rho_n} = \lim_{n \to \infty}\rho_n = 0,
\end{equation}
\begin{equation}\label{eq:level-0-transition}
m_{x_i} - E(g_{\euc}, w_i) = \lim_{n \to \infty} E(h_n, u_n; \bB_{\rho_n}(y_n) \setminus \bB_{s_n}(y_n)),
\end{equation}
\begin{equation}\label{eq:level-0-transition-tips}
\lim_{n \to \infty}E(h_n, u_n; \bB_{\rho_n}(y_n) \setminus \bB_{L^{-1}\rho_n}(y_n)) = \lim_{n \to \infty}E(h_n, u_n; \bB_{L s_n}(y_n) \setminus \bB_{s_n}(y_n)) = 0, 
\end{equation}
where the last property holds for any $L > 1$.
\vskip 1em
\noindent\textit{Step 3: Neck regions and further bubbles.} 
To further break down the right-hand side of~\eqref{eq:level-0-transition}, we let
\[
\lambda_n = \sqrt{\rho_n s_n}, \quad T_n = \log\sqrt{\frac{\rho_n}{s_n}},
\]
By~\eqref{eq:level-0-transition-scales},~\eqref{eq:hn-converge-for-top-bubble}, and Lemma~\ref{lemm:neck-assumptions} (applied with $a_n = 1$ and $\varphi_n = G_n$), upon letting 
$\psi_n = G_n \circ \boldsymbol{\delta}_{y_n, \lambda_n}\circ \mathbf{c}$ and $f_n(t) = e^{-t}$, we see that $\{(\bC_{-T_n, T_n}, \psi_n, \lambda_n, f_n)\}_{n \in \NN}$ is well-prepared in the sense of Definition~\ref{defi:well-prepared}. Next, define
\[
\begin{split}
\widetilde{u}_n :=\ & v_n \circ \psi_n = u_n \circ \boldsymbol{\delta}_{y_n, \lambda_n} \circ \mathbf{c},\\
\widetilde{h}_n :=\ & \lambda_n^{-2}\psi_n^*\sigma_n = \lambda_n^{-2}(\boldsymbol{\delta}_{y_n, \lambda_n}\circ \mathbf{c})^{*}h_n.
\end{split}
\]
Then, given $T > 0$, we have by~\eqref{eq:level-0-transition-tips} that 
\begin{equation}\label{eq:level-0-cylinder-left}
E(\widetilde{h}_n, \widetilde{u}_{n}; \bC_{-T_n, -T_n + T}) = E(u_n, h_n; \bB_{s_n e^{T}}(y_n) \setminus \bB_{s_n}(y_n)) \to 0,
\end{equation}
and, similarly,
\begin{equation}\label{eq:level-0-cylinder-right}
\lim_{n \to \infty}E(\widetilde{h}_n, \widetilde{u}_{n}; \bC_{T_n-T, T_n}) = 0.
\end{equation}
These together show that the condition~\eqref{eq:iv-weakened} in Proposition~\ref{prop:collar-bubble} holds. Now if~\eqref{eq:there-is-energy} fails, then assumption~\eqref{eq:no-energy-on-finite-part} of Proposition~\ref{prop:neck} holds, and we get
\[
\begin{split}
0 = \lim_{n \to \infty}E(\widetilde{h}_n, \widetilde{u}_n; \bC_{-T_n, T_n }) =\ & \lim_{n \to \infty} E(h_n, u_n; \bB_{\rho_n}(y_n) \setminus \bB_{s_n}(y_n))\\
=\ & m_{x_i} - E(g_{\euc}, w_i).
\end{split}
\]
On the other hand if~\eqref{eq:there-is-energy} holds, then by Proposition~\ref{prop:collar-bubble} and Remarks~\ref{rmk:no-leftover} and~\ref{rmk:form-of-mu-i}, we get $Q_i \in \NN$, along with Radon measures $\mu_{i1}, \cdots, \mu_{iQ_{i}}$ on $\RR \times S^1$, such that 
\begin{equation}\label{eq:m_xi-decomposition-1}
m_{x_i} - E(g_{\euc}, w_i) = \sum_{j = 1}^{Q_i}\mu_{ij}(\RR \times S^1).
\end{equation}
Moreover, for each $j \in \{1, \cdots, Q_i\}$, there exist a finite-energy harmonic map $w_{ij}: (\RR \times S^1, g_{\pro}) \to M$, and a finite set $\cB_{ij} \subset \RR \times S^1$, such that $\mu_{ij}$ is given by
\begin{equation}\label{eq:m_xi-decomposition-2}
\mu_{ij} = \frac{1}{2}|\nabla w_{ij}|_{g_{\pro}}^2 \vol_{g_{\pro}} + \sum_{x \in \cB_{ij}}m_{ij, x}\delta_{x}.
\end{equation}
In particular, with the help of~\eqref{eq:loss-at-collars} and~\eqref{eq:wi-energy-lower} we have for all $x \in \cB_{ij}$ that
\[
m_{ij, x} \leq \mu_{ij}(\RR \times S^1)  \leq m_{x_i} - E(g_{\euc}, w_i) \leq  W - \ep_{\text{gap}}.
\] 

Whenever $\cB_{ij}$ is non-empty, recalling that it arises along with $w_{ij}$ from Proposition~\ref{prop:domain-convergence}, and noting how the hypotheses of the latter are verified in Remark~\ref{rmk:form-of-mu-i}, we find ourselves in a position to repeat the argument in Step 2 and Step 3. Doing so results in a decomposition for each $m_{ij,x}$ analogous to~\eqref{eq:m_xi-decomposition-1} and~\eqref{eq:m_xi-decomposition-2}, where in the counterpart of~\eqref{eq:m_xi-decomposition-2}, the amount of energy concentrating at each further bubbling point is now bounded by $W - 2\ep_{\text{gap}}$. Thus, after at most $\lfloor\frac{W}{\ep_{\text{gap}}} \rfloor + 1$ iterations, we end up with finitely many harmonic maps from the round $S^2$ into $M$ whose energy sum up to $m_{x_i}$. Recalling~\eqref{eq:loss-at-collars}, we have thus shown that
\begin{equation}\label{eq:end-of-step3}
W = \tau + E(h, u_0; S') + \text{(energy of finitely many harmonic $S^2$ into $M$)}.
\end{equation}
\vskip 1em
\noindent\textit{Step 4: Finding collar regions.} To prove the energy identity, it remains to show that $\tau$ either vanishes or is again the sum of energies of finitely many harmonic maps from $(S^2, g_{S^2})$ into $M$. To that end, we denote by $\mu$ the measure on the right-hand side of~\eqref{eq:level-0-measure}, so that
\[
\mu = \frac{1}{2}|\nabla u_0|^2_{h} \vol_{h} + \sum_{x \in \cB} m_{x} \delta_{x}.
\]
In the notation of (c3) above, we observe that each $\Sigma^{d}$ is compact, and exhausts $S'$ as $d \to 0$. Hence by the definition of $\tau$ in~\eqref{eq:loss-at-collars}, we can choose a sequence $d_n \to 0$ such that 
\begin{equation}\label{eq:exhaust-W-tau}
W - \tau \geq \mu(\Sigma^{2d_n}) \geq W - \tau - \frac{1}{n}. 
\end{equation}
Still using (c3), together with and the convergence of measures $\frac{1}{2}|\nabla u_n|_{h_n}^2 \vol_{h_n} \to \mu$, as well as (c1), we obtain a subsequence $(h_{n_k}, u_{n_k})$ such that
\vskip 1mm
\begin{enumerate}    
\item[(s1)] $\Sigma^{2d_k} \subset G_{n_k}^{-1}(\Sigma_{n_k}^{d_k}) \subset \Sigma^{\frac{d_k}{2}}$,
\vskip 1mm
\item[(s2)] $\Sigma^{\frac{2d_k}{k}} \subset G_{n_k}^{-1}(\Sigma_{n_k}^{\frac{d_k}{k}}) \subset \Sigma^{\frac{d_k}{2k}}$,
\vskip 1mm
\item[(s3)] $|E(h_{n_k}, u_{n_k}; \Sigma^{2d_k}) - \mu(\Sigma^{2d_k})| \leq \frac{1}{k}$.
\vskip 1mm
\item[(s4)] $|E(h_{n_k}, u_{n_k}; \Sigma^{\frac{d_k}{2k}} \setminus \Sigma^{2d_k}) - \mu(\Sigma^{\frac{d_k}{2k}} \setminus \Sigma^{2d_k})| \leq \frac{1}{k}$.
\vskip 1mm
\item[(s5)] $\frac{l_{n_k, i}}{d_k} \leq \frac{1}{k^2}$, for all $i = 1, \cdots, Q$.
\end{enumerate}
\vskip 1mm
By Proposition~\ref{prop:good-min-max-seq}(a)(b) we may also arrange that 
\begin{equation}\label{eq:energy-convergence}
\big|E(\sigma_{n_k}, v_{n_k}) - W\big|\leq \frac{1}{k}.
\end{equation}
Then, denoting $(h_{n_k}, u_{n_k})$ still by $(h_n, u_n)$, we have the following lower bound
\begin{equation}\label{eq:thick-lower}
\begin{split}
E(h_n, u_n; G_n^{-1}(\Sigma_{n}^{d_n})) \geq\ & E(h_n, u_n; \Sigma^{2d_n})\\
\geq\ & \mu(\Sigma^{2d_n}) - \frac{1}{n} \geq W - \tau - \frac{2}{n},
\end{split}
\end{equation}
together with an upper bound
\begin{equation}\label{eq:thick-upper}
\begin{split}
E(h_n, u_n; G_n^{-1}(\Sigma_{n}^{d_n}))\leq \ &E(h_n, u_n; \Sigma^{\frac{d_n}{2n}})\\
\leq\ & \mu(\Sigma^{\frac{d_n}{2n}}) + \frac{2}{n} \leq W - \tau + \frac{2}{n}.
\end{split}
\end{equation}
Moreover, with the help of~\eqref{eq:exhaust-W-tau}, we see that
\begin{equation}\label{eq:small-transition}
\begin{split}
E(\sigma_n, v_n; \Sigma_n^{\frac{d_n}{n}}\setminus \Sigma_n^{d_n}) = E(h_n, u_n; G_{n}^{-1}(\Sigma_n^{\frac{d_n}{n}}\setminus \Sigma_n^{d_n}) )\leq\ & E(h_n, u_n; \Sigma^{\frac{d_n}{2n}}\setminus \Sigma^{2d_n}) \\
\leq\ & \mu(\Sigma^{\frac{d_n}{2n}}\setminus \Sigma^{2d_n}) + \frac{1}{n} \leq \frac{2}{n}.
\end{split}
\end{equation}
From~\eqref{eq:thick-lower},~\eqref{eq:thick-upper} and~\eqref{eq:energy-convergence}, we get
\begin{equation}\label{eq:collar-captured}
E(\sigma_n, v_n; S\setminus \Sigma_n^{d_n}) = E(\sigma_n, v_n) - E(h_n, u_n, G_n^{-1}(\Sigma_n^{d_n})) = \tau + O(\frac{1}{n}).
\end{equation}
Standard results on closed hyperbolic surfaces (see for instance~\cite[Lemmas A.4, A.5 and Proposition A.6]{Rupflin-Topping-Zhu2013}) yields the following description of the connected components of $S \setminus \Sigma_n^{d_n}$. Specifically, given $0 < l < 2d < 2\sinh^{-1}(1)$ we define
\[
T(l, d) = \frac{\pi^2}{l} - \frac{2\pi}{l}\sin^{-1}\Big( \frac{\sinh(\frac{l}{2})}{\sinh d} \Big),
\]
and then set
\[
T_{n, i} = T(l_{n, i}, d_n),\ \ T_{n, i}' = T(l_{n, i}, \frac{d_n}{n}).
\]
Then eventually $S \setminus \Int\Sigma_n^{d_n}$ consists of connected components $K_{n, 1}, \cdots, K_{n, Q}$, where for each $i \in \{1, \cdots, Q\}$ there is a diffeomorphism
\begin{equation}\label{eq:collar-isometry}
\psi_{n, i}: \bC_{-T_{n, i}, T_{n, i}} \to K_{n, i}
\end{equation}
such that $\psi_{n, i}(\{0\} \times S^1) = \gamma_{n, i}$, and that upon letting
\begin{equation}\label{eq:lambda-f-definition}
\lambda_{n, i} = \frac{l_{n, i}}{2\pi},\quad    f_{n, i}(t) = \cos\big( \frac{l_{n, i}t}{2\pi} \big),
\end{equation}
we have
\begin{equation}\label{eq:collar-pullback-metric}
g_{n, i}: = \lambda_{n, i}^{-2}\psi_{n, i}^*\sigma_n = (f_{n, i})^{-2} \cdot (dt^2 + d\theta^2).
\end{equation}
Noting also that  $S \setminus \Int \Sigma_n^{\frac{d_n}{n}} = \bigcup_{i=1}^{Q} \psi_{n, i}(\bC_{-T_{n, i}', T_{n, i}'})$, and defining 
\[
w_{n, i} = v_{n, i} \circ \psi_{n, i},
\]
we find that~\eqref{eq:collar-captured} and~\eqref{eq:small-transition} translate respectively into
\begin{equation}\label{eq:cylinder-captured}
\lim_{n \to \infty}\sum_{i = 1}^{Q}E(g_{n, i}, w_{n, i}; \bC_{-T_{n,i}, T_{n, i}}) = \tau,
\end{equation}
\begin{equation}\label{eq:small-transition-cylinder}
\sum_{i = 1}^{Q}E(g_{n, i}, w_{n, i}; \bC_{-T_{n, i}, -T_{n, i}'} \cup \bC_{T_{n, i}', T_{n, i}}) \leq \frac{2}{n}.
\end{equation}
Fixing $i \in \{1, \cdots, Q\}$ and dropping it from the subscripts, we next verify that the sequence $\{\bC_{-T_n, T_n}, \psi_n, \lambda_n, f_n\}$ is well-prepared according to Definition~\ref{defi:well-prepared}. Since $d_n$ and $l_n$ both converge to $0$ as $n \to \infty$, and since so does $\frac{l_n}{d_n/n}$ by (s5) above, it is elementary to see that
\begin{equation}\label{eq:Tn-to-infinity}
\lim_{n \to \infty}T_n' = \infty,\quad  \lim_{n \to \infty}(T_n - T_n') = \infty.
\end{equation}
To check the first part of property (p1), we note that $\lambda_n T_n \leq \frac{\pi}{2}$, while for the second part of (p1) we observe
\[
f_n(t) \geq \frac{\sinh(\frac{l_n}{2})}{\sinh d_n},\quad \text{for all }t \in [-T_n, T_n],
\]
and hence
\[
\begin{split}
\sup_{|t|\leq T_n}\frac{\lambda_n}{f_n(t)} \leq\ & \frac{l_n}{2\pi} \cdot \frac{\sinh d_n}{\sinh(\frac{l_n}{2})} = \frac{1}{\pi} \cdot \frac{l_n/2}{\sinh(l_n/2)}\cdot \sinh d_n \to 0, \text{ as }n \to \infty.
\end{split}
\]
Next, property (p2) in Definition~\ref{defi:well-prepared} follows straight from~\eqref{eq:collar-pullback-metric}, and in fact we can take $\alpha_n = 0$ in this case. To check (p3), take an arbitrary $T > 0$ and a sequence $(t_n)$ such that $|t_n|\leq T_n - T$. It suffices to prove that
\begin{equation}\label{eq:p3-reduced}
\frac{f_n(t_n + t)}{f_n(t_n)} \to 1 \text{ smoothly on }  [-T, T].
\end{equation}
To that end, notice that
\begin{equation}\label{eq:sum-expanded}
\frac{f_n(t_n + t)}{f_n(t_n)} =  \cos\big(\frac{l_{n}t}{2\pi}\big) - \tan\big( \frac{l_n t_n}{2\pi} \big) \cdot \sin\big( \frac{l_n t}{2\pi} \big).
\end{equation}
Since $l_n$ tends to $0$ and $T$ is fixed, we have
\begin{equation}\label{eq:cos-to-1}
\cos\big( \frac{l_n t}{2\pi} \big) \to 1\quad \text{smoothly on }[-T, T].
\end{equation}
For the second term on the right-hand side of~\eqref{eq:sum-expanded}, we write
\begin{equation}\label{eq:p3-rewrite}
\tan\big( \frac{l_n t_n}{2\pi} \big) \cdot \sin\big( \frac{l_n t}{2\pi} \big) = \sinh(\frac{l_n}{2})  \tan\big( \frac{l_n t_n}{2\pi}  \big) \cdot \frac{\sin\big(\frac{l_nt}{2\pi} \big)}{\sinh(\frac{l_n}{2})}.
\end{equation}
For all $k \in \NN \cup \{0\}$, the function $t \mapsto \frac{\sin\big(\frac{l_nt}{2\pi} \big)}{\sinh(\frac{l_n}{2})}$ has $C^k$-norm on $[-T, T]$ bounded independently of $n$. On the other hand, since $|t_n| \leq T_n$, we get
\[
|\sinh(\frac{l_n}{2})  \tan\big( \frac{l_n t_n}{2\pi} \big)| \leq \sinh(\frac{l_n}{2}) \cot\big( \sin^{-1}\big( \frac{\sinh(\frac{l_n}{2})}{\sinh d_n} \big)\big) \leq \sinh d_n \to 0, \text{ as }n \to \infty,
\]
which together with the previous observation and~\eqref{eq:p3-rewrite} shows that 
\[
\tan\big( \frac{l_n t_n}{2\pi} \big) \cdot \sin\big( \frac{l_n t}{2\pi} \big) \to 0\quad\text{smoothly on }[-T, T].
\]
Recalling~\eqref{eq:cos-to-1}, we deduce that~\eqref{eq:p3-reduced} holds, and thus so does (p3). 

Having shown that $\{\bC_{-T_n, T_n}, \psi_n, \lambda_n, f_n\}$ is well-prepared according to Definition \ref{defi:well-prepared}, we note further that the second statement in~\eqref{eq:Tn-to-infinity} together with~\eqref{eq:small-transition-cylinder} fulfills the condition~\eqref{eq:iv-weakened} in Proposition~\ref{prop:collar-bubble}; that is, for all $T > 0$, there holds
\begin{equation}\label{eq:no-energy-at-tips}
\lim_{n \to \infty}E(g_{n}, w_n; \bC_{-T_n, -T_n + T}\cup \bC_{T_n - T, T_n}) = 0.
\end{equation}
\vskip 1em
\noindent\textit{Step 5: No energy loss at the collars.} We will be brief about this last step, since it essentially addresses the same dichotomy as the one we had in Step 3 after stating~\eqref{eq:level-0-cylinder-left} and~\eqref{eq:level-0-cylinder-right}. Given $i \in \{1, \cdots, Q\}$, if it happens that
\[
\limsup_{n \to \infty} \sup_{t \in [-T_{n, i} + 1, T_{n, i} - 1]} E(g_{n, i}, w_{n, i}; \bC_{t-1, t+1}) = 0,
\]
then we deduce from Proposition~\ref{prop:neck} that
\[
\lim_{n \to \infty}E(g_{n, i}, w_{n,i}; \bC_{-T_{n, i}, T_{n, i}}) = 0.
\]
On the other hand, if 
\[
\limsup_{n \to \infty} \sup_{t \in [-T_{n, i} + 1, T_{n, i} - 1]} E(g_{n, i}, w_{n, i}; \bC_{t-1, t+1}) > 0,
\]
we instead invoke Proposition~\ref{prop:collar-bubble} along with the observations~\eqref{eq:transition-energy} and~\eqref{eq:form-of-mu-i} to put ourselves in a situation similar to the end of Step 3. After a finite number of iterations of the arguments in Steps 2 and 3, we can express $\lim_{n \to \infty}E(g_{n, i}, w_{n,i}; \bC_{-T_{n, i}, T_{n, i}})$ as the total energy of finitely many harmonic maps from the round $S^2$ into $M$. Recalling~\eqref{eq:cylinder-captured}, we get what we set out to prove about $\tau$ at the start of Step 4, which combines with~\eqref{eq:end-of-step3} to give~\eqref{eq:energy-identity}, and concludes the proof of Theorem~\ref{thm:main-existence}. As a consequence we get Theorem~\ref{thm:min-max-existence}.

\appendix
\section{Patching of charts adapted to multisection}\label{appendix:patching-adapted}
This appendix supplements Section~\ref{subsec:multisections}. Appendices~\ref{subsec:smoothing} and~\ref{subsec:retraction} contain auxiliary technical results used heavily in subsequent constructions. Appendix~\ref{subsec:fixed-stratum} concerns primarily the proof of Proposition~\ref{prop:fixed-stratum-patching}. A summary of the content of Appendices~\ref{subsec:different-strata-I} through~\ref{subsec:different-strata-III} can be found in the proof sketch of Proposition~\ref{prop:multisection-to-sweepout} given in Section~\ref{subsec:multisections}.
\subsection{Notation}\label{subsec:notation-for-patching}
We repeat and expand the list of notation introduced in Section~\ref{subsec:multisections} after Definition~\ref{defi:multisections}. Give an integer $n \geq 3$, and a non-empty subset $I \subset \{1, \cdots, n\}$, with $\be_1, \cdots, \be_n$ being the standard basis on $\RR^{n}$, we define
\[
\bc_{I} := \frac{1}{|I|}\sum_{i \in I}\be_{i},\quad \text{and}\quad \ba_{i; I} := \be_{i} - \bc_{I}\quad (\text{for }i \in I),
\]
omitting $I$ from $\bc_I$ and $\ba_{i; I}$ when $I = \{1, \cdots, n\}$. It is an elementary fact that
\begin{equation}\label{eq:projected-vertices-relation}
\sum_{i \in I}\lambda_i \ba_{i; I} = 0 \quad\text{if and only if}\quad \lambda_i = \lambda_{j}\text{ for any }i, j \in I.
\end{equation}
Then, with $k := |I|$, we let $\Delta^{k-1, I}$ be the convex hull of $\{\ba_{i; I}\ |\ i \in I\}$, which is a $(k-1)$-dimensional simplex, and write $V^{k-1, I}$ for the $(k-1)$-plane in $\RR^{n}$ containing $\Delta^{k-1, I}$, which satisfies 
\[
V^{k-1, I} = \Span\{\ba_{i; I}\ |\ i \in I\}.
\] 
The origin in $V^{k-1, I}$ is denoted $\0^{k-1}$ when we want to emphasize the dimension of the space. Open balls in $V^{n-1}$ centered at $\0^{n-1}$ are denoted $B^{n-1}_{r}$, or sometimes simply $B_{r}$, while intersections of the form $V^{k-1, I} \cap B_{r}$ are often abbreviated as $B^{k- 1, I}_{r}$. Given $i_0 \in I$, we let
\[
V^{k-1, I}_{i_0} : = \{ \sum_{i \in I \setminus \{i_0\}} c_i \ba_{i; I}\ |\ \text{each }c_i \geq 0 \},
\]
and notice that, given $c_i \in \RR$ for each $i\in I$, we have as a result of~\eqref{eq:projected-vertices-relation} that
\begin{equation}\label{eq:sector-characterization}
\sum_{i \in I}c_i\ba_{i; I} \in V^{k-1, I}_{i_0} \quad\text{if and only if}\quad c_{i_0} = \min_{i \in I}c_i.
\end{equation}
In particular, there holds
\begin{equation}\label{eq:sectors-cover}
V^{k-1, I} = \cup_{i\in I}V^{k-1, I}_{i}.
\end{equation}
For any non-empty $I' \subset I$, we define $V^{k-1, I}_{I'} := \cap_{i \in I'} V^{k-1, I}_{i}$ and note by~\eqref{eq:sector-characterization} that it can be expressed as
\begin{equation}\label{eq:sectors}
V^{k-1, I}_{I'} = \big\{\sum_{i \in I \setminus I'}c_i \ba_{i; I}\ |\ \text{each }c_i \geq 0\big\}.
\end{equation}
Similarly, letting $\mathring{V}^{k-1, I}_{I'} : = V^{k-1, I}_{I'} \setminus \big( \cup_{i \in I\setminus I'} V^{k-1, I}_{i} \big)$, we have 
\[
\mathring{V}^{k-1, I}_{I'} = \left\{
\begin{array}{ll}
\big\{\sum_{i \in I \setminus I'}c_i \ba_{i; I}\ |\ \text{each }c_i > 0\big\},& \text{ if }I' \neq I,\\
\\
V^{k-1, I}_{I} = \{\0^{k-1}\}, & \text{ if }I' = I.
\end{array}
\right.
\]
Note also the following analogue of~\eqref{eq:sectors-cover}:
\begin{equation}\label{eq:interior-sectors-components}
V^{k-1, I} \setminus  \big( \cup_{I' \subset I,\ |I'| = 2}V^{k-1, I}_{I'} \big) = \cup_{i \in I}\mathring{V}^{k-1, I}_{i}.
\end{equation}

Next, suppose $\emptyset \neq I' \subsetneq I$, with $|I| = k$ and $|I'| = l$. We mention a few more elementary facts for later use. First, from the definitions of $\bc_{I'}$ and $\bc_{I}$ we have 
\begin{equation}\label{eq:center-distance}
|\bc_{I'} - \bc_{I}| = \sqrt{|I'|^{-1} - |I|^{-1}} =: \lambda_{|I'|, |I|},
\end{equation}
along with the following alternative expressions for $\bc_{I'} - \bc_{I}$:
\begin{equation}\label{eq:change-center}
\bc_{I'} - \bc_{I} = \frac{1}{|I'|}\sum_{i \in I'}\ba_{i; I} = -\frac{1}{|I'|}\sum_{i \in I\setminus I'}\ba_{i; I}.
\end{equation}
Second, there holds the orthogonal decomposition
\begin{equation}\label{eq:decomposition-wrt-lower-simplex}
V^{k-1, I} =   V^{l-1, I'} \oplus \Span\{\ba_{i; I}\ |\ i \in I \setminus I'\},
\end{equation}
and, for each $i \in I'$, the splitting of $\ba_{i; I}$ with respect to this decomposition is given by
\begin{equation}\label{eq:splitting-of-vertex}
\ba_{i; I} = \ba_{i; I'} + (\bc_{I'} - \bc_{I}).
\end{equation}
Finally, given in addition some $\emptyset \neq J \subsetneq I'$ together with $x \in V^{l-1, I'}$ and $z \in V^{k-1, I}_{I'}$, we have 
\begin{equation}\label{eq:projection-of-sectors}
x \in V^{l-1, I'}_{J}\ \text{(resp., $\mathring{V}^{l-1, I'}_{J}$)}\quad \Longleftrightarrow\quad x + z \in V^{k-1, I}_{J} \ \text{(resp., $\mathring{V}^{k-1, I}_{J}$)},
\end{equation}
which can be seen by writing 
\[
x= \sum_{i \in I'} c_i\ba_{i; I'}, \quad\text{with }\min_{i \in I'}c_i = 0,
\]
and then using~\eqref{eq:splitting-of-vertex} and~\eqref{eq:change-center} to get
\[
x = \sum_{i\in I'} c_i \ba_{i; I} + \frac{1}{|I'|}\big( \sum_{j \in I'}c_j \big) \cdot \sum_{i \in I \setminus I'} \ba_{i; I}.
\]
A similar argument shows that if $z \in \mathring{V}^{k-1, I}_{I'}$ and $x \in V^{l-1, I'}$, then 
\begin{equation}\label{eq:projection-of-sectors-2}
x + z \in V^{k-1, I}_{i} \Longrightarrow i \in I'. 
\end{equation}

\subsection{A standard approximation result}\label{subsec:smoothing}
Take integers $n \geq 3$ and $2 \leq k \leq n-1$, along with some $I \subset \{1, \cdots, n\}$ with $|I| = k$. Also, let $J := \{1, \cdots, n\} \setminus I$. The constructions in this section mostly take place within the subspace
\[
W: = \Span\{\ba_{i}\ |\ i\in I\},
\]
and center around the cone $V^{n-1}_{J}$. Recall from~\eqref{eq:sectors} that the latter consists of linear combinations of $\{\ba_{i}\}_{i \in I}$ involving only non-negative coefficients, while its relative boundary in $W$, denoted $\partial V_{J}^{n-1}$, satisfies
\begin{equation}\label{eq:V-J-boundary-definition}
\partial V^{n-1}_{J} = \cup_{i \in I}V^{n-1}_{J \cup \{i\}}, 
\end{equation}
and is smooth away from
\begin{equation}\label{eq:V-J-corner-definition}
\partial_{2} V^{n-1}_J := \cup_{I' \subset I,\ |I'| = 2}V^{n-1}_{J \cup I'}.
\end{equation}

Our goal in this section is to produce a smooth domain $\Omega \subset V^{n-1}_{J}$ which agrees with $V^{n-1}_{J}$ outside of a prescribed distance neighborhood of $\partial_{2}V^{n-1}_{J}$, and whose end in the direction of each $\ba_{i}$, so to speak, exhibits a product structure with respect to the orthogonal decomposition $W = \Span\{\ba_{i}\} \oplus \big( \Span\{\ba_{i}\} \big)^{\perp}$. In addition, the construction also yields a collar neighborhood of $\partial\Omega$ in $\Omega$ where the distance from each leaf to $\partial V^{n-1}_{J}$ satisfies certain upper and lower bounds. The main results of this section are Proposition~\ref{prop:smoothing-flow}, Proposition~\ref{prop:smoothing-product-structure}, and Remark~\ref{rmk:Psi-distance-bounds}.

We start by making a few observations and fixing some additional notation. For each $i \in I$, the unit normal to $V^{n-1}_{J \cup \{i\}} \setminus \partial_{2}V^{n-1}_{J} = \mathring{V}^{n-1}_{J \cup \{i\}}$ that points into $V^{n-1}_{J}$ is given by $\frac{\ba_{i; J \cup \{i\}}}{|\ba_{i; J \cup \{i\}}|}$. Indeed, by~\eqref{eq:splitting-of-vertex} and~\eqref{eq:change-center} we have
\begin{equation}\label{eq:inward-pointing}
\ba_{i; J \cup \{i\}} = \ba_{i} + \frac{1}{n-k + 1}\sum_{j \in I \setminus \{i\}}\ba_{j},
\end{equation}
which implies by a direct computation that
\begin{equation}\label{eq:normal-to-faces}
\ba_{i; J \cup \{i\}} \cdot \ba_{j} = \frac{n-k}{n-k + 1}\,\delta_{ij},\quad\text{for all }j \in I.
\end{equation}
In particular $\ba_{i; J \cup \{i\}}$ is normal to $\mathring{V}^{n-1}_{J \cup \{i\}}$ since the tangent space to the latter is everywhere equal to $\Span\{\ba_{j}\}_{j \in I \setminus \{i\}}$. That $\ba_{i; J \cup \{i\}}$ points into $V^{n-1}_{J}$ follows from the positivity of the coefficients on the right-hand side of~\eqref{eq:inward-pointing}. 

Next we introduce a splitting of $W$ in order to express $\partial V^{n-1}_{J}$ as a graph. Define
\[
\bu := \frac{\bc_I - \bc}{|\bc_I - \bc|} = \frac{1}{\lambda\cdot k}\sum_{i \in I}\ba_{i},
\]
where $\lambda = \lambda_{k, n} = \sqrt{k^{-1} - n^{-1}}$, and the second equality follows from~\eqref{eq:center-distance} and the first expression in~\eqref{eq:change-center}. Recalling also~\eqref{eq:splitting-of-vertex} and the second expression in~\eqref{eq:change-center}, we have the orthogonal decomposition 
\begin{equation}\label{eq:W-orthogonal-decomp}
W = V^{k-1, I} \oplus \Span\{\bu\}.
\end{equation}
Accordingly, given $x \in V^{k-1, I}$ and $y \in \RR$, we write $(x, y)$ for the point $x + y\bu \in W$, so that in particular $\ba_{i} = (\ba_{i; I}, \lambda)$. Also, letting $\bb_{i; I} := -(\lambda k) \ba_{i; I}$, we have
\begin{equation}\label{eq:splitting-vertex-consequence}
\ba_{i; J \cup \{i\}} = \frac{n-k}{(n-k+1)\lambda k}(-\bb_{i; I}, 1).
\end{equation}
Combining this with~\eqref{eq:normal-to-faces}  shows that, given any $(x, y) = \sum_{i \in I}c_i\ba_{i} \in W$, the coefficients satisfy
\begin{equation}\label{eq:coefficient-sign}
c_i = \frac{1}{\lambda k}(y - \bb_{i; I} \cdot x).
\end{equation}
Thus, defining $f: V^{k-1, I} \to [0,\infty)$ by 
\begin{equation}\label{eq:boundary-sector-graph-definition}
f(x) = \max\limits_{i \in I}\bb_{i; I} \cdot x,
\end{equation}
which takes only non-negative values since $\sum_{i \in I}\bb_{i; I}\cdot x = 0$, we get
\begin{equation}\label{eq:graph-sectors}
V^{n-1}_{J}\ \text{(resp., $\mathring{V}^{n-1}_{J}$)} = \{(x, y) \in W\ |\  y \geq f(x)\ \text{(resp., $y > f(x)$)}\},
\end{equation}
and 
\begin{equation}\label{eq:graph-sectors-boundary}
\partial V^{n-1}_{J} =  \{(x, y) \in W \ |\  y = f(x) \}.
\end{equation}

To continue, we define $\theta = \theta_{n, k} \in (0, \frac{\pi}{2})$ to be the angle satisfying
\begin{equation}\label{eq:slope-of-sector}
\sin\theta = \frac{1}{\sqrt{|\bb_{i; I}|^2 + 1}} =  \sqrt{\frac{n}{(n-k + 1)k}}.
\end{equation}
Note that, for all $x, x' \in V^{k-1, I}$,
\begin{equation}\label{eq:f-Lipschitz}
|f(x) - f(x')| \leq (\cot\theta)|x - x'|.
\end{equation}
Also, given $(x, y) \in V_{J}^{n-1}$, we have
\begin{equation}\label{eq:distance-to-translation}
\dist((x, y), \partial V^{n-1}_J) = (y - f(x))\sin\theta.
\end{equation}
The next lemma gathers some additional properties of $f$.
\begin{lemm}\label{lemm:piecewise-affine}
Given $x \in V^{k-1, I}$ and $i_0\in I$, the following are equivalent.
\vskip 1mm
\begin{enumerate}
\item[(i)] $x \in V^{k-1, I}_{i_0}$. 
\vskip 1mm
\item[(ii)] $f(x) =\bb_{i_0; I} \cdot x$.
\vskip 1mm
\item[(iii)] $(x, f(x)) \in V^{n-1}_{J\cup \{i_0\}}$.
\end{enumerate}
\end{lemm}
\begin{proof}
Assume (i). Then $x = \sum_{i \in I \setminus \{i_0\}}c_i \ba_{i; I}$ for some non-negative constants $\{c_i\}_{i \in I \setminus \{i_0\}}$, so that, with $y := \lambda \cdot \sum_{i \in I \setminus \{i_0\}}c_i$, we have $(x, y)= \sum_{i \in I \setminus  \{i_0\}}c_i\ba_{i}$. Since each $c_i \geq 0$, it follows from~\eqref{eq:coefficient-sign} that
\[
\bb_{i_0; I} \cdot x = y \geq \bb_{i; I} \cdot x \quad \text{for all }i \in I,
\]
so (ii) holds. Next, assume (ii) and express $(x, f(x)) \in W$ as $\sum_{i \in I}c_i\ba_{i}$. Then we have by~\eqref{eq:coefficient-sign} that 
\[
c_{i_0} = 0 \leq c_{i} \quad \text{for all }i \in I,
\]
which gives (iii). That (iii) implies (i) is a direct consequence of the relation $\ba_{i} = (\ba_{i; I}, \lambda)$. 
\end{proof}

For the next results, we choose an arbitrary orthonormal basis of $V^{k-1, I}$ and identify it with $\RR^{k-1}$. Let $\zeta$ be a radially symmetric smooth function which is positive in the unit ball $B^{k-1} \subset \RR^{k-1}$, vanishes outside, and satisfies $\int_{B^{k-1}}\zeta = 1$. For $\ep > 0$, define
\begin{equation}\label{eq:f-ep-definition}
f_{\ep}(x) = \int_{\RR^{k-1}}\zeta_{\ep}(y)f(x-y)dy,
\end{equation}
where as usual $\zeta_\ep$ denotes $\frac{1}{\ep^{k-1}}\zeta(\frac{\cdot}{\ep})$. Then each $f_{\ep}$ is a smooth convex function. Also, since $f$ is positive homogeneous of degree $1$, a change of variables shows that
\begin{equation}\label{eq:f-ep-scaling-property}
f_{\ep}(x) = \ep \cdot f_{1}(\frac{x}{\ep}).
\end{equation}
Next, using the fact that $\zeta(y) = \zeta(-y)$, we have
\[
f_{\ep}(x) =  \int_{B^{k-1}}\zeta(y) \cdot \frac{1}{2}[f(x - \ep y) + f(x + \ep y)] dy,
\]
which together with the convexity of $f$ and the Lipschitz bound~\eqref{eq:f-Lipschitz} implies that
\begin{equation}\label{eq:f-ep-f-bound}
f(x) \leq f_{\ep}(x) \leq f(x) + (\cot\theta)\ep,
\end{equation}
and that
\begin{equation}\label{eq:f-ep-monotone}
f_{\ep}(x) \leq f_{\ep'}(x) \leq f_{\ep}(x) + (\cot\theta)(\ep' -\ep ), \quad\text{whenever }0 < \ep \leq \ep'.
\end{equation}
\begin{lemm}\label{lemm:mollification-of-convex}
The functions $f_{\ep}$ have the following further properties.
\vskip 1mm
\begin{enumerate}
\item[(a)] Suppose $\overline{B^{k-1}_{\ep}(x)} \subset \mathring{V}^{k-1, I}_{i_0}$ for some $i_0 \in I$. Then we have
\[
f_{\ep}(x) = f(x) = \bb_{i_0; I} \cdot x, \quad \text{and}\quad \nabla f_{\ep}(x) = \nabla f(x) =  \bb_{i_0; I}.
\]
\vskip 1mm
\item[(b)] For all $x \in \RR^{k-1}$, the vector $(-\nabla f_{\ep}(x), 1)$ lies in $\mathring{V}^{n-1}_{J}$, and we have the bounds
\begin{equation}\label{eq:f-ep-estimates}
|\nabla f_{\ep}(x)| \leq \cot\theta,\quad |\nabla^2 f_{\ep}(x)| \leq \frac{C}{\ep},
\end{equation}
where $C$ depends only on $n, k$, and the choice of $\zeta$. Also, $\nabla f_{\ep}(\0^{k-1}) = \0^{k-1}$.
\vskip 1mm
\item[(c)] Suppose for some $i_0 \in I$ that $\overline{B_{\ep}^{k-1}(x)} \cap V_{i_0}^{k-1, I} = \emptyset$.
Then we have 
\[
(-\nabla f_{\ep}(x), 1) \cdot \ba_{i_0} = 0.
\]
\vskip 1mm
\item[(d)] For all $x$ as in part (c) and $t \geq 0$, we have
\[
f_{\ep}(x + t\ba_{i_0; I}) = f_{\ep}(x) + \lambda t,\quad \nabla f_{\ep}(x + t\ba_{i_0; I}) = \nabla f_{\ep}(x).
\]
\end{enumerate}
\end{lemm}
\begin{proof}
For part (a), let $x_0 \in \RR^{k-1}$ and $i_0 \in I$ be such that $\overline{B^{k-1}_{\ep}(x_0)} \subset \mathring{V}^{k-1, I}_{i_0}$. Then Lemma~\ref{lemm:piecewise-affine} implies
\[
f(x) = \bb_{i_0; I} \cdot x, \quad \text{for all }x \in B^{k-1}_{\ep}(x_0).
\]
Combining this with~\eqref{eq:f-ep-definition} and the evenness of $\zeta$, we obtain the first conclusion of part (a). For the second conclusion, we instead use
\begin{equation}\label{eq:grad-f-ep}
\nabla f_{\ep}(x_0) = \int_{B^{k-1}_{\ep}(x_0)}\zeta_{\ep}(x_0 - y)\nabla f(y) dy,
\end{equation}
which is valid because $f$ is Lipschitz, and thus of class $W^{1, \infty}$.

For part (b), Lemma~\ref{lemm:piecewise-affine} together with~\eqref{eq:grad-f-ep} gives 
\begin{equation}\label{eq:grad-f-ep-convex-hull}
\nabla f_{\ep}(x) = \sum_{i \in I} \Big(\int_{B^{k-1}_{\ep}(x) \cap V^{k-1, I}_{i}} \zeta_{\ep}(x - y) dy\Big) \bb_{i; I},
\end{equation}
which implies that $(-\nabla f_{\ep}(x), 1)$ lies in the convex hull of $\{(-\bb_{i; I}, 1)\}_{i \in I}$. That the latter is contained in $\mathring{V}^{n-1}_{J}$ can be seen at once from~\eqref{eq:inward-pointing} and~\eqref{eq:splitting-vertex-consequence}. The estimates in~\eqref{eq:f-ep-estimates} are both straightforward consequences of~\eqref{eq:grad-f-ep} and the Lipschitz bound~\eqref{eq:f-Lipschitz}. For the last conclusion of part (b), we apply~\eqref{eq:grad-f-ep-convex-hull} at $x = \0^{k-1}$, and use the radial symmetry of $\zeta$ and the fact that $\sum_{i \in I}\bb_{i; I} = 0$. 

Next, under the hypothesis of part (c), we see that the right-hand side of~\eqref{eq:grad-f-ep-convex-hull} reduces to a convex combination of $\{\bb_{i; I}\}_{i \in I \setminus \{i_0\}}$, and the conclusion follows since $\ba_{i_0}$ is orthogonal to $(-\bb_{i; I}, 1)$ for all $i \in I \setminus \{i_0\}$ by~\eqref{eq:normal-to-faces} and~\eqref{eq:splitting-vertex-consequence}.

For part (d), given $x' \in B_{\ep}^{k-1}(x)$, our assumption on $x$, along with~\eqref{eq:sectors-cover}, implies that $x' \in V^{k-1, I}_{i}$ for some $i \neq i_0$. By~\eqref{eq:sectors}, we then have $x' + t\ba_{i_0; I} \in V^{k-1, I}_{i}$ for all $t \geq 0$, so that, by Lemma~\ref{lemm:piecewise-affine}, 
\[
\begin{split}
f(x' + t\ba_{i_0; I}) =\ & \bb_{i; I} \cdot (x' + t\ba_{i_0; I})\\
=\ & f(x') - (\lambda k)(\ba_{i; I} \cdot \ba_{i_0; I}) t = f(x') +\lambda t.
\end{split}
\]
This being true for all $x' \in B_{\ep}^{k-1}(x)$, we deduce from~\eqref{eq:f-ep-definition} the first asserted equality. Since $\{\ba_{i; I}\}_{i \in I \setminus \{i_0\}}$ is a basis for $V^{k-1, I}$, we see that $V^{k-1, I}_{i_0}$ is closed, and hence the condition on $x$ in (c) defines an open subset of $V^{k-1, I}$. We get the second of the asserted equalities in (d) upon differentiating the first with respect to $x$.
\end{proof}
To continue, we let
\[
\cD_{\ep} = \{(x, y)\ |\ y \geq f_{\ep}(x)\},\quad \partial\cD_{\ep} = \{(x, f_{\ep}(x))\ |\ x \in \RR^{k-1}\},
\]
and define a vector field $\eta_{\ep}:\partial\cD_{\ep} \to W$ by
\begin{equation}\label{eq:eta-vector-definition}
\eta_{\ep}(x, f_{\ep}(x)) = (\nabla f_{\ep}(x), -1),\quad\text{for all }x \in \RR^{k-1}.
\end{equation}
Then, fixing also some $\alpha > 0$, we define $\Phi_{\ep, \alpha}:\RR \times \partial\cD_{\ep} \to W$ by
\[
\Phi_{\ep, \alpha}(s, z) := z + \frac{\alpha}{\sin\theta} \bu + (s\sin\theta)\eta_{\ep}(z), 
\]
In other words, for all $x \in \RR^{k-1}$, we have
\begin{equation}\label{eq:Phi-for-smoothing-parametrized}
\Phi_{\ep, \alpha}(s, x, f_{\ep}(x)) = \big(x + (s\sin\theta)\nabla f_{\ep}(x),\ f_{\ep}(x) + \frac{\alpha}{\sin\theta} - s\sin\theta \big).
\end{equation}
Note by~\eqref{eq:f-ep-scaling-property} that $\partial \cD_{\ep} = \ep \cdot \partial\cD_{1}$, and that $\eta_{\ep}(z) = \eta_{1}(\frac{z}{\ep})$ for all $z \in \partial \cD_{\ep}$. As a result, we have
\[
\Phi_{\ep, \alpha}(s, z) = \ep \cdot \Phi_{1, \frac{\alpha}{\ep}}(\frac{s}{\ep}, \frac{z}{\ep}) \quad\text{for }(s, z) \in \RR \times \partial\cD_{\ep}.
\]
Thus, working first with $\ep = 1$ and using the convexity of $f_1$ along with the estimates in Lemma~\ref{lemm:mollification-of-convex}(b), and then invoking the above scaling property, we obtain $\sigma_0(n, k) \in (0, 1)$ such that $\Phi_{\ep, \alpha}$ restricts to a diffeomorphism from $(-\sigma_0\ep, \infty) \times \partial \cD_{\ep}$ onto an open set in $W$. Again by~\eqref{eq:f-ep-estimates}, upon decreasing $\sigma_0$ if necessary, we also have
\begin{equation}\label{eq:mollified-hessian-bound}
|\nabla^2 f_{\ep}(x)| \leq \frac{1}{2\sigma_0 \ep} \quad\text{for all }x \in \RR^{k-1}.
\end{equation}
We then consider the following embedded smooth hypersurfaces in $W$:
\[
\Gamma_{\ep, \alpha}(s) := \Phi_{\ep, \alpha}(\{s\} \times \partial\cD_{\ep}), \quad s \in (-\sigma_0\ep, \infty).
\]
In what follows, we often write $\Phi_{\ep, \alpha}$ and $\eta_{\ep}$ simply as $\Phi$ and $\eta$, respectively.
\begin{lemm}\label{lemm:Phi-leaves-are-graphical}
Given $\ep, \alpha > 0$, there exists a smooth function 
\[
f_{\ep,\alpha}:(-\sigma_0\ep, \infty) \times \RR^{k-1} \to \RR
\]
with the following properties:
\vskip 1mm
\begin{enumerate}
\item[(a)] For all $s \in (-\sigma_0\ep, \infty)$, the set $\Gamma_{\ep, \alpha}(s)$ is the graph of $f_{\ep, \alpha}(s, \cdot)$.
\vskip 1mm
\item[(b)] The assignment
\[
(s, x) \to (x, f_{\ep, \alpha}(s, x))
\]
defines a diffeomorphism from $(-\sigma_0\ep, \infty) \times \RR^{k-1}$ onto $\Phi\big( (-\sigma_0\ep, \infty) \times \partial\cD_{\ep} \big)$.
\vskip 1mm
\item[(c)] $f_{\ep, \alpha}(s, x)$ is strictly decreasing in $s$ for fixed $x$.
\end{enumerate}
\end{lemm}
\begin{proof}
Using the bound~\eqref{eq:mollified-hessian-bound} when $s \in (-\sigma_0\ep, 0)$, and the convexity of $f_{\ep}$ when $s \geq 0$, we see that, in the sense of symmetric matrices, there holds
\begin{equation}\label{eq:projection-derivative}
I_{k-1} + (s\sin\theta )\nabla^2f_{\ep}(x) \geq \frac{I_{k-1}}{2} \quad\text{for all }(s, x) \in (-\sigma_0\ep, \infty) \times \RR^{k-1},
\end{equation}
which implies that the map
\begin{equation}\label{eq:graphical-reparametrization}
\begin{array}{rrcl}
F: &  (-\sigma_0\ep, \infty) \times \RR^{k-1} &  \longrightarrow &  (-\sigma_0 \ep, \infty) \times \RR^{k-1} \\
&  (s, x) &  \longmapsto & (s,\ x + (s\sin\theta)\nabla f_{\ep}(x)),
\end{array}
\end{equation}
is a local diffeomorphism, so in particular its image is open in $(-\sigma_0 \ep, \infty) \times \RR^{k-1}$. Using the boundedness of $\nabla f_{\ep}$ on $\RR^{k-1}$, it is straightforward to see that the image of $F$ is also relatively closed, and thus must be all of $(-\sigma_0 \ep, \infty) \times \RR^{k-1}$. To see that $F$ is injective, suppose $(s, x), (s', x') \in (-\sigma_0\ep, \infty) \times \RR^{k-1}$ are such that $F(s, x) = F(s', x')$. Immediately we get $s = s'$, so that 
\[
x + (s\sin\theta)\nabla f_{\ep}(x) = x' + (s\sin\theta)\nabla f_{\ep}(x').
\]
Rearranging and taking the inner product with $x - x'$ yields
\begin{equation}\label{eq:projection-injective}
|x  - x'|^2 = -(s\sin\theta)(x - x') \cdot (\nabla f_{\ep}(x) - \nabla f_{\ep}(x')).
\end{equation}
When $s \in (-\sigma_0 \ep, 0)$, we use~\eqref{eq:mollified-hessian-bound} to infer from the above that
\[
|x - x'|^2 \leq \frac{1}{2}|x - x'|^2,
\]
which implies $x = x'$. When $s \geq 0$, we use instead the convexity of $t \mapsto f_{\ep}(tx + (1-t)x')$ to see that the right-hand side of~\eqref{eq:projection-injective} is non-positive, which again forces $x = x'$. To sum up, we have proved that the map $F$ defined in~\eqref{eq:graphical-reparametrization} is a bijective local diffeomorphism, and hence a global diffeormorphism. 

To continue, note that the inverse of $F$ must have the form
\begin{equation}\label{eq:X-definition-by-F}
F^{-1}(s, x) = (s, X(s, x)),
\end{equation}
for some smooth function $X: (-\sigma_0\ep, \infty) \times \RR^{k-1} \to \RR^{k-1}$ which satisfies by~\eqref{eq:graphical-reparametrization} that
\begin{equation}\label{eq:X-characterization}
X(s, x) + (s\sin\theta) \nabla f_{\ep}(X(s, x)) = x.
\end{equation}
We then define $f_{\ep, \alpha}:(-\sigma_0, \infty) \times \RR^{k-1} \to \RR$ by
\begin{equation}\label{eq:graphing-function-smooth}
f_{\ep, \alpha}(s, x) = f_{\ep}(X(s, x)) + \frac{\alpha}{\sin\theta} - s\sin\theta,
\end{equation}
which clearly yields a smooth function. Moreover, by~\eqref{eq:X-characterization} and~\eqref{eq:Phi-for-smoothing-parametrized}, there holds for all $(s, x)\in (-\sigma_0\ep, \infty) \times \RR^{k-1}$ that
\begin{equation}\label{eq:Phi-X-relation}
\Phi(s, X(s, x), f_{\ep}(X(s, x))) = (x, f_{\ep, \alpha}(s, x)), 
\end{equation}
from which we get parts (a) and (b) upon recalling that $(s, x) \mapsto (s, X(s, x)) = F^{-1}(s, x)$ is a diffeomorphism from $(-\sigma_0\ep, \infty) \times \RR^{k-1}$ onto itself. 

For part (c), we differentiate~\eqref{eq:X-characterization} with respect to $s$ and use~\eqref{eq:projection-derivative} to get
\[
(\sin\theta) \nabla f_{\ep}(X) \cdot \pa{X}{s} \leq -\frac{1}{2}\cdot \big| \pa{X}{s} \big|^2.
\]
Combining this with the result of differentiating~\eqref{eq:graphing-function-smooth} with respect to $s$, we obtain
\begin{equation}\label{eq:graphing-order}
\pa{f_{\ep, \alpha}}{s} = \nabla f_{\ep}(X) \cdot \pa{X}{s} - \sin\theta  < 0.
\end{equation}
The proof is complete.
\end{proof}

\begin{lemm}\label{lemm:foliation-from-flow}
Let $\theta = \theta_{n, k} \in (0, \frac{\pi}{2})$ be the angle defined by~\eqref{eq:slope-of-sector}. Given $\alpha > 0$, there exist constants $\tau_0(n, k) \in (\cos^2\theta, \frac{1}{1 + \sin^2\theta})$ and $C_0(n, k) > 2$ such that, provided $\ep$ satisfies 
\begin{equation}\label{eq:foliation-ep-threshold}
0 < \ep < \alpha \cdot \min\big\{1,\,\frac{\tau_0 - \cos^2\theta}{\cos\theta} \big\},
\end{equation}
the maps $f_{\ep, \alpha}$ have the following additional properties.
\vskip 1mm
\begin{enumerate}
\item[(a)] For all $0 \leq s \leq \alpha$ and $x \in \RR^{k-1}$,
\begin{equation}\label{eq:graphing-function-relations}
f(x) + \frac{\alpha - s}{\sin\theta} \leq f_{\ep, \alpha}(s, x) \leq f(x) + \frac{\tau_{0}\alpha}{\sin\theta} + (\alpha - s)\sin\theta.
\end{equation}
\vskip 1mm
\item[(b)] Suppose $\overline{B_{\frac{\alpha \tau_0}{\cos\theta}}(x_0)} \subset \mathring{V}^{k-1, I}_{i_0}$ for some $i_0 \in I$.
Then for any $s \in [0, \alpha]$, and any $x \in \RR^{k-1}$ sufficiently close to $x_0$, we have 
\[
f_{\ep, \alpha}(s, x) = f(x) + \frac{\alpha - s}{\sin\theta} = \bb_{i_0; I} \cdot x + \frac{\alpha - s}{\sin\theta}.
\]
\vskip 1mm
\item[(c)] Denoting by $\cK_{\alpha}$ the $\frac{C_0\alpha}{\sin\theta}$-neighborhood of $\partial_2 V^{n-1}_{J}$ in $W$, we have
\begin{equation}\label{eq:smoothing-agree}
\{(x, y) \in W\ |\ y \geq f_{\ep, \alpha}(\alpha, x)\} \setminus \cK_{\alpha} = V^{n-1}_{J}\setminus \cK_{\alpha},
\end{equation}
and similarly
\begin{equation}\label{eq:smoothing-boundary-agree}
\{(x, y) \in W\ |\ y = f_{\ep, \alpha}(\alpha, x)\} \setminus \cK_{\alpha} = \partial V^{n-1}_{J}\setminus \cK_{\alpha}.
\end{equation}
\end{enumerate}
\end{lemm}
\begin{proof}
Let $X:(-\sigma_0\ep, \infty) \times \RR^{k-1} \to \RR^{k-1}$ be as in the previous proof. For part (a), observe first that, by the relation~\eqref{eq:X-characterization} and Lemma~\ref{lemm:mollification-of-convex}(b), there holds
\begin{equation}\label{eq:x-X-distance}
|x - X(s, x)| \leq s\cos\theta,
\end{equation}
and hence
\[
|f_{\ep}(x) - f_{\ep}(X(s, x))| \leq \cot\theta \cdot (s\cos\theta) = \frac{s}{\sin\theta} - s\sin\theta.
\]
Combining this with~\eqref{eq:graphing-function-smooth} yields 
\[
f_{\ep, \alpha}(s, x) \geq f_{\ep}(x) - |f_{\ep}(x) - f_{\ep}(X(s, x))| + \frac{\alpha}{\sin\theta } - s\sin\theta \geq f_{\ep}(x)+ \frac{\alpha - s}{\sin\theta}.
\]
Since $f_{\ep} \geq f$, this gives the first inequality in part (a). For the second inequality, by the convexity of $f_{\ep}$ and the relation~\eqref{eq:X-characterization}, we have
\[
\begin{split}
f_{\ep}(x) \geq\ & f_{\ep}(X(s, x)) + (x - X(s, x)) \cdot \nabla f_{\ep}(X(s, X))\\
=\ & f_{\ep}(X(s, x)) + (s\sin\theta)|\nabla f_{\ep}(X(s, x))|^2,
\end{split}
\]
which together with~\eqref{eq:graphing-function-smooth} and~\eqref{eq:f-ep-f-bound} gives
\[
\begin{split}
f_{\ep, \alpha}(s, x) \leq\ & f_{\ep}(x) + \frac{\alpha}{\sin\theta} - (s\sin\theta)(1 + |\nabla f_{\ep}(X(s, x))|^2)\\
\leq\ & f(x) + (\cot\theta)\ep + \frac{\alpha}{\sin\theta} - s\sin\theta\\
\leq\ & f(x) + (\cot\theta)\ep + \frac{\alpha \cos^2\theta}{\sin\theta}  + (\alpha - s)\sin\theta.
\end{split}
\]
Upon letting $\tau_0 = \frac{1}{2}\big( \cos^2\theta + \frac{1}{1 + \sin^2\theta} \big)$ and imposing on $\ep$ the bound~\eqref{eq:foliation-ep-threshold}, we obtain the second inequality in~\eqref{eq:graphing-function-relations}.

For part (b), we define 
\[
\delta = \frac{\alpha\tau_0}{\cos\theta} - (\ep + \alpha\cos\theta).
\]
Then by~\eqref{eq:x-X-distance}, whenever $s \in [0, \alpha]$ and $|x - x_0| < \delta$, we have
\[
\overline{B_{\ep}(X(s, x))} \subset \overline{B_{\ep + s\cos\theta}(x)} \subset \mathring{V}^{k-1, I}_{i_0}.
\]
For such $s$ and $x$, upon applying Lemma~\ref{lemm:mollification-of-convex}(a) to both $x$ and $X(s, x)$, and using~\eqref{eq:X-characterization}, we obtain
\[
f_{\ep}(X(s, x)) = f(x) - (s\sin\theta)|\bb_{i_0; I}|^2 = \bb_{i_0;I} \cdot x - (s\sin\theta)|\bb_{i_0; I}|^2.
\]
Combining these with~\eqref{eq:graphing-function-smooth} and recalling~\eqref{eq:slope-of-sector} yields both of the asserted equalities. 

Moving on to part (c), we choose 
\[
C_0 = 2 + \frac{\tau_0}{\cos\theta}.
\]
By~\eqref{eq:graphing-function-relations} applied with $s = \alpha$, we have
\begin{equation}\label{eq:two-sided-bound-on-f-ep-a}
f(x) \leq f_{\ep, \alpha}(\alpha, x) \leq f(x) + \frac{\tau_0 \alpha}{\sin\theta},
\end{equation}
which along with~\eqref{eq:graph-sectors} gives the inclusion ``$\subset$'' in~\eqref{eq:smoothing-agree}. For the reverse inclusion, take $(x, y)$ such that 
\[
y \geq f(x), \quad (x, y) \not\in \cK_{\alpha}.
\]
In the case where
\begin{equation}\label{eq:away-from-corners}
\overline{B_{\frac{\alpha \tau_0}{\cos\theta}}^{k-1}(x)} \cap \big( \cup_{I' \subset I,\ |I'| = 2}V^{k-1, I}_{I'} \big) = \emptyset,
\end{equation}
since $\overline{B_{\frac{\alpha \tau_0}{\cos\theta}}^{k-1}(x)}$ is connected, we obtain from~\eqref{eq:interior-sectors-components} some $i_0$ such that the hypothesis of part (b) holds. Consequently $f_{\ep, \alpha}(\alpha, x) = f(x)$, so that $(x, y)$ lies in the left-hand side of~\eqref{eq:smoothing-agree} as desired. On the other hand, if~\eqref{eq:away-from-corners} fails, then there exists $I' \subset I$ with $|I'| = 2$, along with some $x' \in V^{k-1, I}_{I'}$, such that
\[
|x - x'| \leq \frac{\alpha \tau_0}{\cos\theta}.
\]
By Lemma~\ref{lemm:piecewise-affine}, we see that $(x', f(x')) \in \partial_2 V^{n-1}_{J}$, while~\eqref{eq:f-Lipschitz} implies that
\[
|f(x) - f(x')| \leq \frac{\alpha\tau_0}{\cos\theta}\cdot \cot\theta.
\]
Consequently, we have
\[
d((x, f(x)), \partial_2 V^{n-1}_{J}) \leq \frac{\alpha}{\sin\theta} \cdot \frac{\tau_0}{\cos\theta}.
\]
Recalling our choice of $C_0$ above, we deduce that
\[
(x, t) \in \cK_{\alpha} \quad\text{for all }t \in [f(x), f(x) + \frac{\alpha}{\sin\theta}).
\]
In view of~\eqref{eq:two-sided-bound-on-f-ep-a}, as well as our assumption about $(x, y)$, we conclude that $y > f_{\ep, \alpha}(\alpha, x)$, so again $(x, y)$ lies in the left-hand side of~\eqref{eq:smoothing-agree}. 

It remains to prove~\eqref{eq:smoothing-boundary-agree}. To that end, we note from the above argument that, given $x \in \RR^{k-1}$, if either $(x, f(x))$ or $(x, f_{\ep, \alpha}(\alpha, x))$ is assumed to lie outside of $\cK_{\alpha}$, then $x$ must satisfy~\eqref{eq:away-from-corners}, in which case $f(x) = f_{\ep, \alpha}(\alpha, x)$. This observation suffices to establish~\eqref{eq:smoothing-boundary-agree}. The proof is complete.
\end{proof}
Combining the previous two lemmas, we get the following result.
\begin{prop}\label{prop:smoothing-flow}
Given $\alpha > 0$, there exist a smooth domain $\Omega \subset V^{n-1}_{J}$ whose boundary relative to $W$ we denote by $\partial\Omega$, a smooth vector field $\xi:\partial\Omega \to \mathring{V}^{n-1}_{J} \cap B_{1}$, and some $\alpha_1 > \alpha$ such that, with $\tau_0$, $C_0$, and $\cK_{\alpha}$ having the same meaning as in Lemma~\ref{lemm:foliation-from-flow}, the following hold.
\vskip 1mm
\begin{enumerate}
\item[(a)] $\Omega \setminus \cK_{\alpha} = V^{n-1}_{J}\setminus \cK_{\alpha}$, and $\partial\Omega \setminus \cK_{\alpha} = \partial V^{n-1}_{J}\setminus \cK_{\alpha}$. 
\vskip 1mm
\item[(b)] Given $i_0 \in I$, we have for all $z \in V^{n-1}_{J \cup\{i_0\}}\setminus \cK_{\alpha}$ that
\begin{equation}\label{eq:xi-coincide-with-inward-normal}
\xi(z) = \frac{\ba_{i_0; J \cup \{i_0\}}}{|\ba_{i_0; J \cup \{i_0\}}|} = \text{inward unit normal to $\partial\Omega$}.
\end{equation}
On the other hand, for all $z \in \partial\Omega$ which lies outside of the $\frac{C_0\alpha}{\sin\theta}$-neighborhood of $V^{n-1}_{J \cup \{i_0\}}$, we have $\xi(z) \cdot \ba_{i_0} = 0$.
\vskip 1mm
\item[(c)] The map $\Psi$ defined by 
\begin{equation}\label{eq:flow-Psi-defi}
\Psi(s, z) =  z + s\xi(z)
\end{equation}
is a diffeomorphism from $(-\infty, \alpha_1) \times \partial \Omega$ onto an open set in $W$. 
\end{enumerate}
\vskip 2mm
In fact, there exists a smooth function $\bphi:(-\infty, \alpha_{1})  \times \RR^{k-1} \to \RR$ related to the above objects in the following way.
\vskip 2mm
\begin{enumerate}
\item[(d)] The map $(s, x) \mapsto (x, \boldsymbol{\varphi}(s, x))$ is a diffeomorphism from $(-\infty, \alpha_{1}) \times \RR^{k-1}$ onto $\Psi((-\infty, \alpha_{1}) \times \partial\Omega)$. Also, $s \mapsto \bphi(s, x)$ is strictly increasing in $s$ for fixed $x$.
\vskip 1mm
\item[(e)] We have
\[
\begin{split}
\Omega = \ &\{(x, y) \in W\ |\ y \geq \bphi(0, x)\};\\
\Psi(\{s\} \times \partial\Omega) =\ & \{(x, \bphi(s, x))\ |\ x \in \RR^{k-1}\}, \quad \text{for all }s \in (-\infty, \alpha_1).
\end{split}
\]
\vskip 1mm
\item[(f)] For all $s \in [0, \alpha]$ and $x \in \RR^{k-1}$, we have
\begin{equation}\label{eq:graphing-function-relations-1}
f(x) + \frac{s}{\sin\theta} \leq \boldsymbol{\varphi}(s, x) \leq f(x) + \frac{\tau_0\alpha}{\sin\theta} + s\sin\theta.
\end{equation}
Moreover, provided $x \in \RR^{k-1}$ satisfies~\eqref{eq:away-from-corners}, we have 
\begin{equation}\label{eq:graphing-function-relations-2}
\boldsymbol{\varphi}(s, x) = f(x) + \frac{s}{\sin\theta}\quad\text{for all }s \in [0, \alpha].
\end{equation}
\end{enumerate}
\end{prop}
\begin{rmk}\label{rmk:Psi-distance-bounds}
The following further consequences of Proposition~\ref{prop:smoothing-flow} will be useful later.
\vskip 1mm
\begin{enumerate}
\item From conclusions (d) and (e), we deduce that $\Psi(s, z) \in \Omega$ if and only if $s\geq 0$. 
\vskip 1mm
\item By~\eqref{eq:graphing-function-relations-1} and~\eqref{eq:distance-to-translation}, there holds for all $s \in [0, \alpha]$ and $z \in \partial\Omega$ that
\begin{equation}\label{eq:Psi-distance-estimates}
s \leq \dist(\Psi(s, z), \partial V^{n-1}_{J}) \leq \tau_0\alpha + s\sin^2\theta.
\end{equation}
Letting $B_{s}(\partial V^{n-1}_{J})$ denote the $s$-distance neighborhood of $\partial V^{n-1}_{J}$ in $V^{n-1}$, then by the first inequality in~\eqref{eq:graphing-function-relations-1} together with the intermediate value theorem, we have
\begin{equation}\label{eq:distance-neighborhood-contained}
\Omega \cap B_{s}(\partial V^{n-1}_{J}) \subset \Psi([0, s) \times \partial\Omega),
\end{equation}
and likewise $\Omega \cap \overline{B_{s}(\partial V^{n-1}_{J})} \subset \Psi([0, s] \times \partial\Omega)$.
\vskip 1mm
\item Still using~\eqref{eq:graphing-function-relations-1} and~\eqref{eq:distance-to-translation}, we have for all $s \in [\tau_0\alpha, \alpha]$, that
\begin{equation}\label{eq:smoothing-agree-outside-distance-nbhd}
\begin{split}
\Omega \setminus B_{s}(\partial V^{n-1}_{J}) =\ & \{(x, y) \in W\ |\ y \geq \max\{\bphi(0, x), f(x) + \frac{s}{\sin\theta}\}\}\\
=\ & \{(x, y) \in W\ |\ y \geq f(x) + \frac{s}{\sin\theta}\} = V^{n-1}_{J} \setminus B_{s}(\partial V^{n-1}_{J}).
\end{split}
\end{equation}
\end{enumerate}
\end{rmk}
\begin{proof}[Proof of Proposition~\ref{prop:smoothing-flow}]
Fix any $\ep$ satisfying~\eqref{eq:foliation-ep-threshold}, and let $X$ and $f_{\ep, \alpha}$ be as in Lemma~\ref{lemm:Phi-leaves-are-graphical} and its proof. (See~\eqref{eq:graphical-reparametrization},~\eqref{eq:X-definition-by-F}, and~\eqref{eq:graphing-function-smooth}.) With $\alpha_{1}$ taken to be $\alpha + \sigma_0 \ep$, we define
\begin{equation}\label{eq:bphi-definition}
\boldsymbol{\varphi}(s, x) = f_{\ep, \alpha}(\alpha - s, x) \quad\text{for }(s, x) \in (-\infty, \alpha_{1}) \times \RR^{k-1}.
\end{equation}
Then part (f) holds thanks to parts (a) and (b) of Lemma~\ref{lemm:foliation-from-flow}, while by Lemma~\ref{lemm:Phi-leaves-are-graphical}, we get that $s \mapsto \bphi(s, x)$ is strictly increasing, and that $(s, x) \mapsto (x, \boldsymbol{\varphi}(s, x))$ is a diffeomorphism from $(-\infty, \alpha_{1}) \times \RR^{k-1}$ onto its image. We determine this image towards the end of the proof, but apart from that we are done with (d). To continue, we define
\[
\Omega = \{(x, y) \in W\ |\ y \geq \boldsymbol{\varphi}(0, x)\},
\]
so that $\partial\Omega$ is the graph of $\bphi(0, \cdot)$, and that by~\eqref{eq:graph-sectors} and the first inequality in~\eqref{eq:graphing-function-relations-1} we have 
\[
\Omega \subset V^{n-1}_J,
\]
as required. The assertion of part (a) is then a restatement of Lemma~\ref{lemm:foliation-from-flow}(c). Noting that $\Phi_{\ep, \alpha}(\alpha, \cdot)$ is a diffeomorphism from $\partial\cD_{\ep}$ onto $\partial\Omega$, we define $\xi: \partial\Omega \to W$ by the following relation:
\begin{equation}\label{eq:xi-definition}
\xi(\Phi_{\ep, \alpha}(\alpha, z)) = -(\sin\theta)\eta_{\ep}(z), \quad\text{for }z \in \partial\cD_{\ep}.
\end{equation}
In terms of the function $X:(-\sigma_0\ep, \infty) \times \RR^{k-1} \to \RR^{k-1}$ from the proof of Lemma~\ref{lemm:Phi-leaves-are-graphical}, particularly the relation~\eqref{eq:Phi-X-relation}, we can express the above definition as
\begin{equation}\label{eq:xi-X-relation}
\xi(x, \boldsymbol{\varphi}(0, x)) = (\sin\theta) (-\nabla f_{\ep}(X(\alpha, x)), 1), \quad\text{for } x \in \RR^{k-1}.
\end{equation}
Lemma~\ref{lemm:mollification-of-convex}(b) then implies that $\xi$ takes values in $\mathring{V}^{n-1}_{J} \cap B_{1}$, as required. 

To verify first statement in part (b), take any $z \in V^{n-1}_{J \cup \{i_0\}} \setminus \cK_{\alpha}$. By~\eqref{eq:graph-sectors-boundary}, Lemma~\ref{lemm:piecewise-affine}, and the final paragraph of the proof of Lemma~\ref{lemm:foliation-from-flow}, we can express it as
\[
z = (x, f(x)) = (x, \bphi(0, x)),
\]
for some $x \in V^{k-1, I}_{i_0}$ that also satisfies~\eqref{eq:away-from-corners}. In particular, we see from~\eqref{eq:interior-sectors-components} that $\overline{B_{\frac{\alpha \tau_0}{\cos\theta}}(x)} \subset \mathring{V}^{k-1, I}_{i_0}$. As in the proof of Lemma~\ref{lemm:foliation-from-flow}(b), we then deduce using~\eqref{eq:foliation-ep-threshold} and~\eqref{eq:x-X-distance} that
\[
\overline{B_{\ep}(X(\alpha, x))} \subset \mathring{V}^{k-1, I}_{i_0}.
\]
It follows from Lemma~\ref{lemm:mollification-of-convex}(a) that $\nabla f_{\ep}(X(\alpha, x)) = \bb_{i_0; I}$, which together with~\eqref{eq:xi-X-relation} shows that 
\[
\xi(z) = (\sin\theta)(-\bb_{i_0; I}, 1).
\]
Upon recalling~\eqref{eq:slope-of-sector} and~\eqref{eq:splitting-vertex-consequence}, we get the first equality in~\eqref{eq:xi-coincide-with-inward-normal}. Using again the inclusion $\overline{B_{\frac{\alpha \tau_0}{\cos\theta}}(x)} \subset \mathring{V}^{k-1, I}_{i_0}$, along with Lemma~\ref{lemm:foliation-from-flow}(b), we see that 
\[
\bphi(0, x') = \bb_{i_0; I} \cdot x',
\]
for all $x'$ near $x$, and hence the right-hand side in the above expression for $\xi(z)$ coincides with the inward unit normal to $\partial\Omega$ at $z$. This proves the first assertion in (b). To continue, take instead some $z \in \partial\Omega \setminus B_{\frac{C_0\alpha}{\sin\theta}}\big( V^{n-1}_{J \cup \{i_0\}} \big)$, and write it as 
\[
z = (x, \bphi(0,x)) \quad\text{for some } x \in \RR^{k-1}.
\]
Suppose by contradiction that $\overline{B_{\frac{\tau_0\alpha}{\cos\theta}}(x)} \cap V^{k-1, I}_{i_0} \neq\emptyset$. Then by Lemma~\ref{lemm:piecewise-affine} together with~\eqref{eq:f-Lipschitz} we have 
\[
\dist((x, f(x)), V^{n-1}_{J\cup\{i_0\}}) \leq \frac{\alpha}{\sin\theta}\cdot \frac{\tau_0}{\cos\theta},
\]
which along with~\eqref{eq:graphing-function-relations-1} and the choice of $C_0$ leads to a contradiction to our assumption about $z$. Thus we must have 
\[
\overline{B_{\frac{\tau_0\alpha}{\cos\theta}}(x)} \cap V^{k-1, I}_{i_0} = \emptyset,
\]
in which case~\eqref{eq:x-X-distance} and~\eqref{eq:foliation-ep-threshold} together gives $\overline{B_{\ep}(X(\alpha, x))} \cap V^{k-1, I}_{i_0} = \emptyset$, and we get from Lemma~\ref{lemm:mollification-of-convex}(c) and~\eqref{eq:xi-X-relation} that
\[
\xi(z) \cdot \ba_{i_0} = 0,
\]
as asserted, and we are done with (b). 

Turning to $\Psi$ as defined by~\eqref{eq:flow-Psi-defi}, after a straightforward computation using~\eqref{eq:xi-definition} and the definition of $\Phi_{\ep, \alpha}$, we have
\[
\Psi(\alpha - s, \Phi_{\ep, \alpha}(\alpha, z)) = \Phi_{\ep, \alpha}(s, z),
\]
which together with our choice of $\alpha_{1}$ proves (c). Recalling~\eqref{eq:bphi-definition} and Lemma~\ref{lemm:Phi-leaves-are-graphical}, we see from the above that $\Psi(\{s\} \times \partial\Omega)$ coincides with the graph of $\boldsymbol{\varphi}(s, \cdot)$. This completes parts (d) and (e). Having already verified part (f) right after defining $\bphi$, we are done. 
\end{proof}

We next observe that, for each $i_0 \in I$, certain portions of $V^{n-1}_{J}$ and $\Omega$ split into products with respect to the decomposition of $W$ into $\Span\{\ba_{i_0}\}$ and its orthogonal complement. Fixing $i_0 \in I$ and writing $\{i_0\}^{c}$ for $\{1, \cdots, n\}\setminus \{i_0\}$, we define
\[
W' = \Span\big\{\ba_{i; \{i_0\}^{c}}\ |\ i \in I \setminus \{i_0\}\big\}.
\]
Note that $V_{J}^{n-2, \{i_0\}^{c}} \subset W'$. Then, in analogy with~\eqref{eq:V-J-boundary-definition} and~\eqref{eq:V-J-corner-definition}, we let
\begin{equation}\label{eq:rel-boundary-notation}
\begin{split}
\partial V^{n-2, \{i_0\}^{c}}_{J} =\ & \cup_{i \in I \setminus \{i_0\}} V^{n-2, \{i_0\}^{c}}_{J \cup \{i\}},\\
\partial_2 V^{n-2, \{i_0\}^{c}}_{J} =\ & \cup_{I' \subset I \setminus \{i_0\},\ |I'| = 2} V^{n-2, \{i_0\}^{c}}_{J \cup I'}.
\end{split}
\end{equation}
The first set in~\eqref{eq:rel-boundary-notation} coincides with the relative boundary of $V^{n-2, \{i_0\}^{c}}_{J}$ in $W'$.
\begin{lemm}\label{lemm:sectors-product-structure}
In the above notation, and recalling also that $W = \Span\{\ba_{i}\ |\ i \in I\}$, we have
\vskip 1mm
\begin{enumerate}
\item[(a)] $W = W' \oplus \Span\{\ba_{i_0}\}$;
\vskip 1mm
\item[(b)] Given $y \in W'$ and $t  > 0$, we have $\dist(y + t\ba_{i_0}, V^{n-1}_{J \cup \{i_0\}}) \geq t\sqrt{\frac{n-1}{n}}$. 
\vskip 1mm
\item[(c)] With $y$ and $t$ as in part (b), we have
\[
y + t\ba_{i_0} \in V^{n-1}_{J}\quad \text{if and only if}\quad y \in V^{n-2, \{i_0\}^{c}}_{J}.
\]
The same equivalence holds with $V^{n-1}_{J}$ and $V^{n-2, \{i_0\}^{c}}_{J}$ replaced, respectively, by $\partial V^{n-1}_{J}$ and $\partial V^{n-2, \{i_0\}^{c}}_{J}$.
\vskip 1mm
\item[(d)] Given $r > 0$, $y \in W'$, and $t \geq r \sqrt{\frac{n}{n-1}}$, we have
\[
y + t\ba_{i_0} \in V^{n-1}_{J}\cap B_{r}(\partial V^{n-1}_{J}) \quad\text{if and only if}\quad y \in V^{n-2, \{i_0\}^{c}}_{J}\cap B_{r}(\partial V^{n-2, \{i_0\}^{c}}_{J}).
\]
\end{enumerate}
\end{lemm}
\begin{proof}
For part (a), by~\eqref{eq:splitting-of-vertex} and~\eqref{eq:change-center}, we have whenever $i \neq i_0$ that
\begin{equation}\label{eq:part-a-identity}
\ba_{i; \{i_0\}^{c}} = \ba_{i} + \frac{1}{n-1}\ba_{i_0}.
\end{equation}
Applying this with $i \in I \setminus \{i_0\}$ shows that $W = W' + \Span\{\ba_{i_0}\}$. The orthogonality of this decomposition can be seen by dotting both sides above with $\ba_{i_0}$. 

For part (b), take any point $z$ in $V^{n-1}_{J \cup \{i_0\}}$, which, by the definition of the latter and the relation~\eqref{eq:part-a-identity}, can be written as
\[
z = \sum_{i \in I \setminus \{i_0\}}c_i \ba_{i} =\Big( \sum_{i \in I \setminus \{i_0\}}c_i \ba_{i; \{i_0\}^{c}}\Big) - \big(\frac{1}{n-1}\sum_{i \in I \setminus \{i_0\}} c_i\big) \ba_{i_0},
\]
where each $c_i$ is non-negative. From the orthogonal decomposition in part (a), we infer that 
\[
|y + t\ba_{i_0} - z| \geq \big( t + \frac{1}{n-1}\sum_{i \in I \setminus \{i_0\}} c_i \big)|\ba_{i_0}| \geq t|\ba_{i_0}|.
\]
This gives part (b). 

For part (c), we use the definition of $W'$ to write $y = \sum_{i \in I \setminus \{i_0\}} \lambda_{i}\ba_{i; \{i_0\}^{c}}$, and use~\eqref{eq:part-a-identity} to get
\[
y + t\ba_{i_0} = \Big(\sum_{i \in I \setminus \{i_0\}} \lambda_{i}\ba_{i}\Big) + \big(t + \frac{1}{n-1}\sum_{i\in I \setminus \{i_0\}}\lambda_i \big) \ba_{i_0},
\]
from which we get the asserted equivalence, since $t > 0$ by assumption. For the boundary version, we need only note in addition that $y + t\ba_{i_0}$ cannot lie in $V^{n-1}_{J \cup \{i_0\}}$, since for that to happen we would need $\lambda_{i} \geq 0$ for each $i \in I\setminus \{i_0\}$ and $t + \frac{1}{n-1}\sum_{i\in I \setminus \{i_0\}}\lambda_i = 0$, which contradicts $t > 0$.

For part (d), suppose first that $y \in V^{n-2, \{i_0\}^{c}}_{J}$ and that 
\[
|y - z| < r \quad\text{for some }z \in \partial V^{n-2, \{i_0\}^{c}}_{J}.
\]
Then by part (c) we have $z + t\ba_{i_0} \in \partial V^{n-1}_{J}$, and consequently $y + t\ba_{i_0} \in V^{n-1}_{J} \cap B_{r}(\partial V^{n-1}_{J})$. This proves the implication ``$\Leftarrow$''. Conversely, suppose that 
\[
y + t\ba_{i_0} \in V^{n-1}_{J} \cap B_{r}(\partial V^{n-1}_{J}),
\]
so that $y \in V^{n-2, \{i_0\}^{c}}_{J}$ to start with. If in addition $y + t\ba_{i_0}$ lies in $\partial V^{n-1}_{J}$, then $y \in \partial V^{n-2, \{i_0\}^{c}}_{J}$ by part (c), and we are done. If $y + t\ba_{i_0}$ lies not in $\partial V^{n-1}_{J}$, then its distance to the latter must be realized somewhere on $\partial V^{n-1}_{J} \setminus \partial_{2}V^{n-1}_{J}$. This implies the existence of some $i \in I$ and $z \in V^{n-1}_{J \cup \{i\}}$ such that 
\begin{equation}\label{eq:closest-point-conditions}
|y + t\ba_{i_0} - z| < r,\quad y + t\ba_{i_0} - z \perp \Span\{\ba_{j}\ |\ j \in I \setminus \{i\}\}.
\end{equation}
The first condition, along with part (b) and our assumption that $t\sqrt{\frac{n-1}{n}} \geq r$, implies that $i \neq i_0$, in which case we infer from the second condition and part (a) that $y + t\ba_{i_0} - z \in W'$, so that $z - t\ba_{i_0} \in W'$, because $y \in W'$ to start with. Noting that
\[
(z - t\ba_{i_0}) + t\ba_{i_0} = z \in \partial V^{n-1}_{J},
\]
we get from part (c) that $z - t\ba_{i_0} \in \partial V^{n-2, \{i_0\}^{c}}_{J}$. The first condition in~\eqref{eq:closest-point-conditions} then gives $y \in B_{r}(\partial V^{n-2, \{i_0\}^{c}}_{J})$. The proof is complete.
\end{proof}

We end this section with a number of properties of the domain $\Omega$ that are analogous to those of $V^{n-1}_{J}$ established in Lemma~\ref{lemm:sectors-product-structure}. Still we write $x + y\bu$ as $(x, y)$ for $x \in V^{k-1, I}$ and $y \in \RR$, and identify $V^{k-1,I}$ with $\RR^{k-1}$.
\begin{prop}\label{prop:smoothing-product-structure}
Given $\alpha > 0$ and $i_0 \in I$, there exist a smooth domain $\Omega'_{i_0} \subset V^{n-2, \{i_0\}^{c}}_{J}$ whose boundary relative to $W'$ we denote by $\partial\Omega'_{i_0}$, and a smooth vector field $\xi': \partial\Omega'_{i_0} \to W'$, such that the objects produced by Proposition~\ref{prop:smoothing-flow} satisfy in addition the following.
\vskip 1mm
\begin{enumerate}
\item[(a)] Given $y \in W'$ and $t \geq 2\sqrt{\frac{n}{n-1}}\cdot\frac{C_0\alpha}{\sin\theta} =: T_0$, where $C_0$ is the constant in Lemma~\ref{lemm:foliation-from-flow}, we have 
\[
y + t\ba_{i_0} \in \Omega \quad\text{if and only if} \quad y \in \Omega'_{i_0}.
\]
The same equivalence holds with $\Omega$ and $\Omega'_{i_0}$ replaced by $\partial\Omega$ and $\partial\Omega'_{i_0}$, respectively.
\vskip 1mm
\item[(b)] For all $y \in \partial\Omega'_{i_0}$ and $t \geq T_0$, we have $\xi(y + t\ba_{i_0}) = \xi'(y)$. Moreover, the map 
\begin{equation}\label{eq:Psi-on-slice}
(s, y) \mapsto y + s\xi'(y)
\end{equation}
is a diffeomorphism from $(-\alpha_1, \alpha_1) \times \partial\Omega_{i_0}'$ onto an open set in $W'$, and the image point lies in $\Omega_{i_0}'$ if and only if $s \geq 0$.
\vskip 1mm
\item[(c)] Denoting by $\cK'$ the $\frac{C_0\alpha}{\sin\theta}$-neighborhood of $\partial_2 V^{n-2, \{i_0\}^{c}}_{J}$ in $W'$, we have 
\[
\Omega'_{i_0} \setminus \cK' = V^{n-2, \{i_0\}^{c}}_{J} \setminus \cK', \quad\quad \partial\Omega'_{i_0} \setminus \cK' = \partial V^{n-2, \{i_0\}^{c}}_{J} \setminus \cK'.
\]
Moreover, outside of $\cK'$ the vector field $\xi'$ coincides with the inward pointing unit normal to $\partial\Omega'_{i_0}$.
\end{enumerate}
\end{prop}
\begin{proof}
We fix the same $\ep$ as in the proof of Proposition~\ref{prop:smoothing-flow}, and use freely the notation from that proof. To begin, given $y \in W'$, we define 
\[
A_{y} = \{t \geq \frac{T_0}{2}\ |\ y + t\ba_{i_0} \in \Omega\},
\]
which is closed in $[\frac{T_0}{2}, \infty)$. Observe also that, by Lemma~\ref{lemm:sectors-product-structure}(b), we have
\begin{equation}\label{eq:product-structure-distance}
\dist(y + t\ba_{i_0}, V^{n-1}_{J\cup\{i_0\}}) \geq \frac{C_0\alpha}{\sin\theta}\quad\text{for all }t \geq \frac{T_0}{2}.
\end{equation}
To prove that $A_{y}$ is relatively open, suppose $t_* \in A_{y}$, so that
\[
z := y + t_*\ba_{i_0} \in \Omega.
\]
If $z$ lies in $\Omega \setminus \partial\Omega$, which is open relative to $W$, then certainly $y + t\ba_{i_0} \in \Omega$ for all $t$ sufficiently close to $t_*$. If instead $z \in \partial\Omega$, then, upon writing $z = (x, \bphi(0, x))$, we have from~\eqref{eq:product-structure-distance} and the proof of Proposition~\ref{prop:smoothing-flow}(b) that 
\[
\overline{B_{\frac{\tau_0\alpha}{\cos\theta}}(x)} \cap V^{k-1, I}_{i_0} = \emptyset,
\]
where $\tau_0$ is given by Lemma~\ref{lemm:foliation-from-flow}. Again letting $\delta = \frac{\alpha\tau_0}{\cos\theta} - (\ep + \alpha\cos\theta)$, we deduce upon recalling~\eqref{eq:foliation-ep-threshold} and~\eqref{eq:x-X-distance} that
\begin{equation}\label{eq:X-ep+delta-ball-inclusion}
\overline{B_{\ep+\delta}(X(\alpha, x))} \subset \RR^{k-1}\setminus V^{k-1, I}_{i_0}.
\end{equation}
Since $|\ba_{i_0; I}| = \sqrt{\frac{k-1}{k}} < 1$, we can then apply Lemma~\ref{lemm:mollification-of-convex}(d) to $X(\alpha, x) - \delta\ba_{i_0; I}$ to obtain, for all $s \in [-\delta, \infty)$,
\[
\begin{split}
&\nabla f_{\ep}(X(\alpha, x) + s\ba_{i_0; I}) = \nabla f_{\ep}(X(\alpha, x)) \quad\text{and}\\
&f_{\ep}(X(\alpha, x) + s\ba_{i_0; I}) = f_{\ep}(X(\alpha, x)) + s\lambda.
\end{split}
\]
The first relation together with~\eqref{eq:X-characterization} shows that the point $x' := X(\alpha, x) + s\ba_{i_0; I}$ satisfies
\[
x' + (\alpha\sin\theta) \nabla f_{\ep}(x') = x + s\ba_{i_0; I},
\]
which forces $X(\alpha, x + s\ba_{i_0; I}) = X(\alpha, x) + s\ba_{i_0; I}$. Combining this with the second relation above, and recalling~\eqref{eq:graphing-function-smooth} and the definition of $\bphi$, we conclude that
\[
\bphi(0, x + s\ba_{i_0; I}) = \bphi(0, x) + \lambda s,\quad\text{for all }s \geq -\delta,
\]
and consequently
\begin{equation}\label{eq:A-y-open}
\begin{split}
z + s\ba_{i_0} = \ &  (x + s\ba_{i_0; I}, \bphi(0, x) + s\lambda)\\
=\ &  (x + s\ba_{i_0; I}, \bphi(0, x + s\ba_{i_0; I})) \in \partial \Omega.
\end{split}
\end{equation}
This being true for all $s \geq -\delta$, we see in particular that $t_*$ lies in the relative interior of $A_y$. Having shown that $A_{y}$ is both open and closed in $[\frac{T_0}{2}, \infty)$, we conclude that $A_{y}$ is either empty or all of $[\frac{T_0}{2}, \infty)$. From this dichotomy, we see that
\begin{equation}\label{eq:Omega'-expression}
\begin{split}
\Omega_{i_0}' :=\ & \{y \in W'\ |\ A_{y} = [\frac{T_0}{2}, \infty)\}\\ =\ & \{y \in W'\ |\ y + t\ba_{i_0} \in \Omega \quad\text{for some } t \geq \frac{T_0}{2}\}.
\end{split}
\end{equation}
In particular $\Omega_{i_0}'$ is closed in $W'$, and is contained in $V^{n-2, \{i_0\}^{c}}_{J}$ thanks to the inclusion $\Omega \subset V^{n-1}_{J}$ and Lemma~\ref{lemm:sectors-product-structure}(c). From~\eqref{eq:Omega'-expression}, we deduce also that, for each $t > \frac{T_0}{2}$, 
\begin{equation}\label{eq:boundary-Omega'-relation}
\Omega_{i_0}' \setminus \partial\Omega_{i_0}' = \{y \in W'\ |\ y + t\ba_{i_0} \in \Omega \setminus \partial\Omega\}.
\end{equation}
This together with~\eqref{eq:Omega'-expression} implies, again for each $t > \frac{T_0}{2}$, that
\begin{equation}\label{eq:boundary-Omega'-as-transverse-intersection}
\partial\Omega_{i_0}'+ t\ba_{i_0} = \partial\Omega \cap \big( W'  + t \ba_{i_0} \big).
\end{equation}
In particular, given $y \in \partial\Omega_{i_0}'$ and $t > \frac{T_0}{2}$, the tangent space $T_{y + t\ba_{i_0}}\partial\Omega$ contains $\ba_{i_0}$. The latter being orthogonal to $W'$ by Lemma~\ref{lemm:sectors-product-structure}(a), we see that the right-hand side above is a transverse intersection. From this, and counting dimensions, we conclude that $\partial\Omega_{i_0}'$ is an embedded smooth hypersurface in $W'$. Equation~\eqref{eq:boundary-Omega'-as-transverse-intersection} then implies that 
\begin{equation}\label{eq:product-structure}
\begin{array}{rcl}
\partial\Omega_{i_0}' \times (\frac{T_0}{2}, \infty) & \longrightarrow & \partial\Omega \cap \big(  W' + (\frac{T_{0}}{2}, \infty) \cdot \ba_{i_0} \big)\\
(y, t) & \longmapsto & y + t\ba_{i_0},
\end{array}
\end{equation}
is a diffeomorphism, with the inverse mapping given by restricting the orthogonal projections onto $W'$ and $\Span\{\ba_{i_0}\}$.

Towards defining $\xi'$, take any $y \in \partial\Omega_{i_0}'$, and observe by~\eqref{eq:boundary-Omega'-as-transverse-intersection} and~\eqref{eq:product-structure-distance} that
\[
z: = y + T_0\ba_{i_0} \in \partial\Omega\setminus B_{\frac{C_0\alpha}{\sin\theta}}(V^{n-1}_{J \cup \{i_0\}}).
\]
Thus, expressing $z$ as $(x, \bphi(0, x))$ for some $x \in \RR^{k-1}$, we again have the inclusion~\eqref{eq:X-ep+delta-ball-inclusion}. Tracing the argument leading up to~\eqref{eq:A-y-open}, and recalling~\eqref{eq:xi-X-relation}, we get for all $s \geq 0$ that
\begin{equation}\label{eq:calculation-leading-to-xi'}
\begin{split}
\xi(z + s\ba_{i_0}) = \ &\xi(x + s\ba_{i_0; I}, \bphi(0, x + s\ba_{i_0; I}))\\
=\ & (\sin\theta)(-\nabla f_{\ep}(X(\alpha, x + s\ba_{i_0; I})), 1)\\
=\ &  (\sin\theta)(-\nabla f_{\ep}(X(\alpha, x)), 1) = \xi(z).
\end{split}
\end{equation}
This suggests defining $\xi'$ to be
\[
\xi'(y) = \xi(y + T_0\ba_{i_0})\quad \text{for }y \in \partial\Omega'_{i_0}.
\]
Recalling that $y + T_0\ba_{i_0} \in \partial\Omega\setminus B_{\frac{C_0\alpha}{\sin\theta}}(V^{n-1}_{J \cup \{i_0\}})$, we get from Proposition~\ref{prop:smoothing-flow}(b) that $\xi'(y) \perp \ba_{i_0}$, so that, by Lemma~\ref{lemm:sectors-product-structure}(a), $\xi'$ takes values in $W'$, as required.

We next verify (a), (b), and (c). For part (a), the asserted equivalences follow immediately from~\eqref{eq:Omega'-expression} and~\eqref{eq:boundary-Omega'-as-transverse-intersection}. For part (b), the equality involving $\xi(y + t\ba_{i_0})$ is a consequence of the calculation~\eqref{eq:calculation-leading-to-xi'}. In particular, for all $s \in (-\alpha_1, \alpha_1)$, $y \in \partial\Omega_{i_0}'$, and $t \geq T_0$, we have
\begin{equation}\label{eq:Psi-slice-expression}
\Psi(s, y + t\ba_{i_0}) = s\xi'(y) + y + t\ba_{i_0}.
\end{equation}
Turning to the map~\eqref{eq:Psi-on-slice}, we claim that, given $s \in (-\alpha_1, \alpha_1)$, $z \in \partial\Omega$, and $t_{1} > T_0 + \alpha_1\sqrt{\frac{n}{n-1}}$, we have the following equivalence:
\begin{equation}\label{eq:Psi-slice-as-image}
\Psi(s, z) \in W' + t_{1}\ba_{i_0} \quad\text{if and only if}\quad z \in W' + t_{1}\ba_{i_0}. 
\end{equation}
Indeed, suppose first that $z \in W' + t_{1}\ba_{i_0}$. Then since $z \in \partial\Omega$ to start with, we may use~\eqref{eq:boundary-Omega'-as-transverse-intersection} to write 
\[
z  =  y + t_1\ba_{i_0},\quad\text{for some }y \in \partial\Omega_{i_0}'.
\]
Since $t_1 > T_0$, we get from~\eqref{eq:Psi-slice-expression} that $\Psi(s, z) \in W' + t_{1}\ba_{i_0}$, which proves the implication ``$\Leftarrow$'' in~\eqref{eq:Psi-slice-as-image}. Conversely, starting with $\Psi(s, z) \in W' + t_1\ba_{i_0}$, by the triangle inequality, and the bounds $|\xi(z)| \leq 1$ and $|s| < \alpha_{1}$, we have 
\[
\begin{split}
t_2 : = z \cdot \frac{\ba_{i_0}}{|\ba_{i_0}|^2} \geq\ & \Psi(s, z) \cdot \frac{\ba_{i_0}}{|\ba_{i_0}|^2} -  \frac{|z - \Psi(s, z)|}{|\ba_{i_0}|}\\
>\ & t_1 - \alpha_1\sqrt{\frac{n}{n-1}} > T_0.
\end{split}
\]
Thus, as in the proof of the implication``$\Leftarrow$'', we get some $y \in \partial\Omega_{i_0}'$ such that 
\[
z = y + t_2\ba_{i_0},
\]
in which case~\eqref{eq:Psi-slice-expression} gives $\Psi(s, z) \in W' + t_2\ba_{i_0}$. Comparing this with our assumption on $\Psi(s, z)$ forces $t_2 = t_1$, and consequently $z \in W' + t_1\ba_{i_0}$. This finishes the proof of~\eqref{eq:Psi-slice-as-image}. Fixing a choice of $t_{1} > T_0 + \alpha_1\sqrt{\frac{n}{n-1}}$, we deduce from~\eqref{eq:Psi-slice-as-image} and~\eqref{eq:boundary-Omega'-as-transverse-intersection} that
\[
\begin{split}
\Psi\big( (-\alpha_1, \alpha_1) \times (\partial\Omega_{i_0}' + t_1\ba_{i_0}) \big) = \Psi((-\alpha_1, \alpha_1) \times \partial\Omega) \cap (W' + t_1\ba_{i_0}).
\end{split}
\]
Since $\Psi$ is a diffeomorphism on $(-\alpha_1, \alpha_1) \times \partial\Omega$, with $\Psi((-\alpha_1, \alpha_1) \times \partial\Omega)$ being open in $W$, and since, thanks to~\eqref{eq:Psi-slice-expression}, the map~\eqref{eq:Psi-on-slice} is none other than
\[
(s, y) \mapsto \Psi(s, y + t_1\ba_{i_0}) - t_1\ba_{i_0}, 
\]
we conclude that it sends $(-\alpha_1,\alpha_1) \times \partial\Omega_{i_0}'$ diffeomorphically onto $\big(\Psi((-\alpha_1, \alpha_1) \times \partial\Omega) - t_1\ba_{i_0} \big)\cap W'$. The last statement of part (b) follows from part (a), the formula~\eqref{eq:Psi-slice-expression}, and Remark~\ref{rmk:Psi-distance-bounds}.

For part (c), we first show that
\begin{equation}\label{eq:product-exterior}
y + T_0\ba_{i_0} \in W\setminus \cK_{\alpha} \quad\text{provided }y \in W' \setminus \cK',
\end{equation}
where recall that $\cK_{\alpha}$ is the $\frac{C_0\alpha}{\sin\theta}$-neighborhood of $\partial_{2}V^{n-1}_{J}$ in $W$. To see this, take $y \in W' \setminus \cK'$, along with a subset $I' \subset I$ of size $2$, and some point $z \in V^{n-1}_{J \cup I'}$. In the case $i_0 \in I'$, we have $V^{n-1}_{J \cup I'} \subset V^{n-1}_{J \cup \{i_0\}}$, and it follows from~\eqref{eq:product-structure-distance} that 
\[
|y + T_0\ba_{i_0} - z| \geq \frac{C_0\alpha}{\sin\theta}.
\]
On the other hand, if $i_0 \not\in I'$, we may express $z$ as
\[
z = c_{i_0}\ba_{i_0} + \sum_{i \in I \setminus (I' \cup \{i_0\})}c_i\ba_{i}  \in \Span\{\ba_{i_0}\} + \Big(\sum_{i \in I \setminus (I' \cup \{i_0\})}c_i\ba_{i; \{i_0\}^{c}}\Big),
\]
where each $c_i$ is non-negative, and the second step uses~\eqref{eq:part-a-identity}. Since the term in parentheses belongs to $V^{n-2, \{i_0\}^{c}}_{J \cup I'}$, we deduce by Lemma~\ref{lemm:sectors-product-structure}(a) and the assumption $y \in W' \setminus \cK'$ that
\[
|y + T_0\ba_{i_0} - z| \geq |y - \big(\sum_{i \in I \setminus (I' \cup \{i_0\})}c_i\ba_{i; \{i_0\}^{c}}\big)| \geq \frac{C_0\alpha}{\sin\theta}.
\]
Thus~\eqref{eq:product-exterior} holds indeed, which together with part (a), Proposition~\ref{prop:smoothing-flow}(a), and Lemma~\ref{lemm:sectors-product-structure}(c) shows that, under the hypothesis $y \in W' \setminus \cK'$, we have $y \in \Omega_{i_0}'$ if and only if $y \in V^{n-2, \{i_0\}^{c}}_{J}$. This gives the first of the two equalities of sets in part (c), and the argument for the second one is the same. Finally, given $y \in \partial\Omega_{i_0}' \setminus \cK'$, by~\eqref{eq:product-exterior} and part (a) we have 
\[
y + T_0\ba_{i_0} \in \partial\Omega \setminus \cK_{\alpha}.
\]
The definition of $\xi'$ and Proposition~\ref{prop:smoothing-flow}(b) then implies that $\xi'(y)$ coincides with the inward unit normal to $\partial\Omega$ at $y + T_0\ba_{i_0}$. Using the product structure given by~\eqref{eq:product-structure}, we deduce that $\xi'(y)$ is orthogonal to $T_y\partial\Omega_{i_0}'$ in $W'$. That it points into $\Omega_{i_0}'$ at $y$ follows from the last statement of part (b). The proof is now complete.
\end{proof}
\subsection{Sweeping out $1$-handlebodies}\label{subsec:retraction}
Take any integer $m$ with $m \geq 2$. It is a standard fact that, given an orientable $1$-handlebody $H$ of dimension $m + 1$, there exists a subset $\Gamma \subset H$ of finite $1$-dimensional Hausdorff measure, together with a Lipschitz map 
\[
[0, 1] \times \partial H \to H
\]
which 
\vskip 1mm
\begin{itemize}
\item acts as projection onto the second factor on $\{0\} \times \partial H$, and takes $\{1\} \times \partial H$ to $\Gamma$, 
\vskip 1mm
\item restricts to a diffeomorphism from $[0, 1) \times \partial H$ to $H \setminus \Gamma$,
\vskip 1mm
\item shrinks the length of $2$-vectors tangent to $\partial H$ uniformly to zero as $t \to 1^{-}$.
\end{itemize}
The purpose of this section is to indicate how one might obtain such a Lipschitz map, by carrying out the construction on what we consider to be the essential building blocks of $H$. At the end of the section we elaborate on the last point above regarding the length of $2$-vectors, which in fact is a statement about the derivative of the map. (See Remark~\ref{rmk:retraction-end-remark}.)

To be more specific about the scope of what we do here, note that since $m + 1 \geq 3$, upon writing $\ell$ for the genus of $H$, which determines its diffeomorphism type~\cite[Corollary VI.11.4]{Kosinski1993}, we can assume without loss of generality that $H$ is the boundary connected sum of $\ell$ copies of a distance neighborhood of $S^{1} \times \{\0^{m-1}\}$ in $\RR^{m + 1} = \RR^{2} \times \RR^{m-1}$. Thus, roughly speaking, issues of smoothness aside, $H$ consists of two types of building blocks, namely distance neighborhoods of $(\text{arc in }\RR^{2}) \times \{\0^{m-1}\}$, and those of $(\text{figure ``$Y$''}) \times \{\0^{m-1}\}$. The remainder of this section is devoted to constructing a smooth version, denoted $\Omega$, of the latter type of distance neighborhoods, along with a Lipschitz map 
\[
\bh: [0, 1] \times \partial \Omega \to \Omega
\]
having properties corresponding to those enumerated at the beginning. The main results are Lemma~\ref{lemm:Omega-S-basics} and Proposition~\ref{prop:bh-properties}. 

Let $\theta: = \frac{\pi}{3}$ and consider
\[
\bv_{1} = (\sin\theta, \cos\theta), \quad \bv_{2} = (-\sin\theta, \cos\theta), \quad \bv_{3} = (0, -1),
\]
along with
\[
Y = \cup_{k = 1}^{3}[0, \infty) \cdot \bv_{k},\quad\quad Y^{-} =  \cup_{k = 1}^{3}[0, \infty) \cdot (-\bv_{k}).
\]
Also, let $\bw_{\lambda}$ denote the result of rotating $\bv_{\lambda}$ counterclockwise by $\frac{\pi}{2}$. Given $\lambda \in \{1, 2, 3\}$, and employing indices modulo $3$, we define
\[
\begin{split}
\cV_{\lambda} :=\ & [0, \infty) \cdot \bv_{\lambda} + [0, \infty) \cdot \bv_{\lambda + 1} = \{a\bv_{\lambda} + b\bv_{\lambda + 1}\ |\ a, b \geq 0\},
\end{split}
\]
and write $\mathring{\cV}_{\lambda}$ for the interior of $\cV_{\lambda}$. We shall construct the domain $\Omega$ by putting together four pieces (see~\eqref{eq:Omega-for-retraction}), contained respectively in $\mathring{\cV}_{\lambda} \times \RR^{m-1}$ ($\lambda = 1, 2, 3$) and $B_{\rho}(Y) \times \RR^{m-1}$, where $\rho > 0$ is to be determined. Accordingly, the map $\bh: [0, 1] \times \partial\Omega  \to \Omega$ is also defined in a piecewise manner (see~\eqref{eq:bh-definition}). In fact, among the regions $\cV_{\lambda} \times \RR^{m-1}$, we only consider the case $\lambda = 1$, since the construction for $\lambda = 2$ or $3$ turn out to be the same up to suitably rotating the $\RR^{2}$ factor in $\RR^{m + 1} = \RR^{2} \times \RR^{m-1}$. 

Thus we focus for the moment on $\cV_{1}$, written as $\cV$ until we begin to assemble the pieces in~\eqref{eq:Omega-for-retraction} and~\eqref{eq:S-for-retraction}. Our first task is to construct a more convenient pair of coordinates on $\cV$, which is a two-stage process. The important results along the way are Lemma~\ref{lemm:foliation-by-translated-mollification} and Lemma~\ref{lemm:rotated-graph}, and the end product is the map $\widetilde{\bF}$ in~\eqref{eq:tilde-bF-definition}. To start, let
\begin{equation}\label{eq:f-for-retraction}
f(x) = (\cot \theta)|x|,
\end{equation}
so that we have
\[
\cV = \{(x, y) \in \RR^2\ |\ x \in \RR,\ y \geq f(x)\}.
\]
Next, taking $\zeta$ to be a smooth even function that is positive on $(-1, 1)$, vanishes elsewhere, and integrates to $1$ on $\RR$, we again define $f_{\ep}$ by~\eqref{eq:f-ep-definition}, which yields a smooth, convex function which is even, and satisfies
\begin{equation}\label{eq:f-ep-agree}
f_{\ep}(x) = f(x) \quad\text{provided }|x| \geq \ep.
\end{equation}
Also, the scaling property~\eqref{eq:f-ep-scaling-property}, the bounds~\eqref{eq:f-ep-f-bound} and~\eqref{eq:f-ep-monotone}, and the derivative estimates~\eqref{eq:f-ep-estimates} all continue to hold. 
\begin{lemm}\label{lemm:foliation-by-translated-mollification}
Define $\bF: \RR \times [0, \infty) \to \cV$ by 
\begin{equation}\label{eq:bF-definition}
\bF(x, \ep) = \left\{ 
\begin{array}{ll}
(x, \ep + f_{\ep}(x)) , & \text{ if }\ep > 0,\\
(x, f(x)), & \text{ if }\ep = 0.
\end{array}
\right.
\end{equation}
\vskip 1mm
\begin{enumerate}
\item[(a)] $\bF$ is a bi-Lipschitz map from $\RR \times [0,\infty)$ onto $\cV$. The Lipschitz constants of $\bF$ and $\bF^{-1}$ are both bounded from above in terms of $\theta$.
\vskip 1mm
\item[(b)] $\bF$ restricts to a diffeomorphism from $\RR \times (0, \infty)$ onto $\mathring{\cV}$. 
\vskip 1mm
\item[(c)] For all $(x, \ep) \in \RR \times [0, \infty)$ and $t \geq 0$, we have
\begin{equation}\label{eq:bF-scaling-property}
\bF(tx, t\ep) = t \bF(x, \ep).
\end{equation}
\end{enumerate}
\end{lemm}
\begin{proof}
That $\bF$ takes values in $\cV$ and sends $\RR \times (0, \infty)$ into $\mathring{\cV}$ can be seen from~\eqref{eq:f-ep-f-bound}. By~\eqref{eq:f-ep-monotone} and the first inequality in~\eqref{eq:f-ep-estimates}, we have for all $x, x' \in \RR$ and $\ep, \ep' > 0$ that
\[
|\bF(x, \ep) - \bF(x', \ep')| \leq (1 + \cot\theta)(|x - x'| + |\ep - \ep'|),
\]
showing that $\bF$ is Lipschitz on $\RR \times (0, \infty)$. Since $\lim_{\ep \to 0^+}\bF(x, \ep) = \bF(x, 0)$, the above estimate extends to $\RR \times [0, \infty)$. Next, with the help of the strict monotonicity of $\ep \mapsto f_{\ep}(x) + \ep$ implied by~\eqref{eq:f-ep-monotone}, as well as the intermediate value theorem, we see that $\bF$ is a bijection from $\RR \times (0,\infty)$ onto $\mathring{\cV}$, and also maps $\RR \times [0,\infty)$ bijectively onto $\cV$.

To see that $\bF^{-1}$ is Lipschitz, take $p_1 := \bF(x_1, \ep_1)$ and $p_2 := \bF(x_2, \ep_2)$, with $\ep_1, \ep_2 > 0$. The definition of $\bF$ immediately gives
\[
|x_1 - x_2| \leq |p_1 - p_2|.
\]
Assuming without loss of generality that $\ep_{2} \geq \ep_{1}$, and using~\eqref{eq:f-ep-monotone} and~\eqref{eq:f-ep-estimates}, we also have
\begin{equation}\label{eq:inverse-bF-Lipschitz}
\begin{split}
|p_2 - p_1| \geq\ & f_{\ep_{2}}(x_2) + \ep_{2} - f_{\ep_{1}}(x_1) - \ep_1 \\
\geq\ & f_{\ep_{1}}(x_2) + \ep_{2} - f_{\ep_{1}}(x_1) - \ep_1 \geq -(\cot\theta)|x_1 - x_2| + \ep_{2} - \ep_{1},
\end{split}
\end{equation}
Combining this with the previous estimate on $|x_1 - x_2|$ yields
\[
|x_{1} - x_{2}| + |\ep_{2} - \ep_{1}| \leq (2 + \cot\theta)\,|\bF(x_1, \ep_1) - \bF(x_2, \ep_2)|.
\]
As before, this estimate continues to hold on $\RR \times [0,\infty)$, from which we infer that $\bF^{-1}$ is Lipschitz on $\bF(\RR \times [0, \infty)) = \cV$. This proves part (a).

For part (b), that $\bF$ is smooth on $\RR \times (0,\infty)$ can be seen from the relation~\eqref{eq:f-ep-scaling-property}. Moreover, the monotonicity property~\eqref{eq:f-ep-monotone} implies that $\pa{f_{\ep}}{\ep} + 1 \geq 1$, so that $d\bF$ is invertible everywhere on $\RR \times (0, \infty)$. Having already seen above that $\bF|_{\RR \times (0, \infty)}$ is a bijection onto $\mathring{\cV}$, we are done with part (b).

For part (c), when $t = 0$, both sides of the stated equality is $(0, 0)$. Next, when $\ep = 0$, the result follows straight from~\eqref{eq:f-for-retraction}. Finally, when $\ep > $ and $t > 0$, we note from~\eqref{eq:f-ep-scaling-property} that
\[
t\ep + f_{t\ep}(tx) = t\ep + t\ep f_{1}(\frac{x}{\ep}) = t\ep + t f_{\ep}(x),
\]
which completes the proof.
\end{proof}
Next, take another smooth cut-off function $\varphi:\RR \to [0, 1]$ such that
\begin{equation}\label{eq:central-cutoff}
\varphi(x) = \left\{
\begin{array}{ll}
0, & \text{ if }x \leq 1/4\\
1, & \text{ if }x \geq 3/4
\end{array}
\right.;\quad \quad 0 \leq \varphi' \leq 3.
\end{equation}
Taking $\ell : =  4\cot\theta$, we define, for $x \in \RR$ and $\ep > 0$,
\[
g_{\ep}(x) = f_{\ep}(\varphi(\frac{|x|}{\ep})x) + \ep \ell.
\]
By our choice of $\varphi$, as well as~\eqref{eq:f-ep-agree}, we have
\begin{equation}\label{eq:g-particular}
\begin{split}
g_{\ep}(x) =\left\{
\begin{array}{ll}
f_{\ep}(0) + \ep \ell = g_{\ep}(0), & \text{if }|x| \leq \ep/4, \\
f(x) + \ep\ell, & \text{if }|x| \geq \ep. 
\end{array}
\right.
\end{split}
\end{equation}
Also~\eqref{eq:f-ep-scaling-property} implies that
\[
g_{\ep}(x) = \ep g_{1}(\frac{x}{\ep}).
\]
Recalling from the convexity and evenness of $f_{\ep}$ that $f'_{\ep}(x)\geq f'_{\ep}(0) = 0$ for all $x \geq 0$, and using~\eqref{eq:f-ep-estimates} to bound $f'_{\ep}$ from above, we have for all $x \geq 0$ and $\ep > 0$ that
\begin{equation}\label{eq:g-ep-bounds}
0 \leq g'_{\ep}(x)\leq 4\cot\theta,\quad \cot\theta \leq \pa{g_{\ep}}{\ep}(x) \leq 5\cot\theta,
\end{equation}
where for the second estimate we used also~\eqref{eq:f-ep-monotone} and our choice of $\ell$. In particular, for all $x > 0$ and $\ep > 0$, we have 
\begin{equation}\label{eq:g-ep-f}
g_{\ep}(x) >  \lim_{\ep \to 0^+}g_{\ep}(x) = (\cot\theta)x.
\end{equation}

For the next result, we choose $\mu: = \frac{\csc\theta}{f_{1}(0) + \ell} \in (0, \frac{\sec\theta}{4})$, so that $g_{\mu}(0) = \csc\theta$. The scaling property noted below~\eqref{eq:g-particular} implies that
\begin{equation}\label{eq:g-ep-scaling}
g_{\ep\mu}(x) = \ep g_{\mu}(\frac{x}{\ep}),\quad\text{for all }x \in \RR,\ \ep > 0.
\end{equation}
\begin{lemm}\label{lemm:rotated-graph}
For each $\ep > 0 $, there exists a smooth function $\psi_{\ep}:(\ep \cot\theta, \infty) \to (0, \infty)$ which is related to $g_{\mu\ep}$ by
\begin{equation}\label{eq:psi-g-graph-relation}
\{e^{i\theta}\cdot (x + i g_{\mu\ep}(x))\ |\ x  > 0\} = \{- \psi_{\ep}(t) + it\ |\ t > \ep\cot\theta \}.
\end{equation}
Moreover, define, for $(x, y) \in \mathring{\cV}$,
\[
\bP(x, y) = \left\{
\begin{array}{ll}
(\psi_{x}(y), y), & \text{ if }x > 0,\\
(0, y), & \text{ if }x = 0,\\
(-\psi_{|x|}(y),y), & \text{ if } x < 0.
\end{array}
\right.
\]
Then $\bP$ is a diffeomorphism from $\mathring{\cV}$ onto itself, and the Lipschitz constants of both $\bP$ and $\bP^{-1}$ on $\mathring{\cV}$ are bounded from above in terms of $\theta$.
\end{lemm}
\begin{proof}
Define, for $x \geq 0$ and $\ep  > 0$,
\begin{equation}\label{eq:bt-definition}
\bt(x, \ep) = g_{\mu\ep}(x)\cos\theta + x\sin\theta.
\end{equation}
By the first inequality in~\eqref{eq:g-ep-bounds} we have 
\begin{equation}\label{eq:dt-dx-lower-bound}
\pa{\bt}{x} \geq \sin\theta.
\end{equation}
Consequently, $\bt(x, \ep)$ increases strictly from $\ep\cot\theta$ to $\infty$ as $x$ increases from $0$ to $\infty$, and the map $(x, \ep) \mapsto (\bt(x, \ep), \ep)$ takes $(0, \infty) \times (0, \infty)$ diffeomorphically onto $\{(t, \ep)\ |\ \ep > 0, \ t > \ep\cot\theta\}$. Expressing its inverse as $(t, \ep) \mapsto (\bx(t, \ep), \ep)$, and letting
\begin{equation}\label{eq:psi-ep-definition}
\psi_{\ep}(t) := g_{\mu\ep}(\bx(t, \ep))\sin\theta - \bx(t, \ep) \cos\theta,
\end{equation}
we get
\begin{equation}\label{eq:psi-g-relation}
it - \psi_{\ep}(t) = e^{i\theta}\cdot (\bx(t, \ep) + i g_{\mu\ep}(\bx(t, \ep))),
\end{equation}
from which~\eqref{eq:psi-g-graph-relation} immediately follows. Also, from~\eqref{eq:g-ep-f} and~\eqref{eq:psi-ep-definition}, we see that $\psi_{\ep}$ is indeed positive.

Towards establishing the asserted properties of $\bP$, we make several observations. First, by the relation $\bt(\bx(t, \ep), \ep) = t$ and~\eqref{eq:psi-ep-definition}, we have that
\begin{equation}\label{eq:bP-target}
t - (\cot\theta)\psi_{\ep}(t) = \bx(t, \ep)\csc\theta > 0.
\end{equation}
Combining $\bt(\bx(t, \ep), \ep) = t$ instead with the scaling property~\eqref{eq:g-ep-scaling}, we get $\bx(t, \ep) = \ep  \bx(\frac{t}{\ep}, 1)$, and hence 
\begin{equation}\label{eq:psi-scaling-property}
\psi_{\ep}(t) = \ep \psi_{1}(\frac{t}{\ep}).
\end{equation}
Next, given $(t_1, \ep_1)$ and $(t_2, \ep_2)$ satisfying 
\[
t_i > \ep_i \cot\theta > 0, \quad\text{for }i = 1, 2,
\]
we let $x_i = \bx(t_i, \ep_i)$, and assume without loss of generality that $x_2 \geq x_1$. Then, from~\eqref{eq:dt-dx-lower-bound} and the estimate for $\pa{g_{\ep}}{\ep}$ in~\eqref{eq:g-ep-bounds}, we have
\[
\begin{split}
t_2 - t_1 = \bt(x_2, \ep_2) - \bt(x_1, \ep_1) \geq\ & (\sin\theta)(x_2 - x_1) + \bt(x_1, \ep_2) - \bt(x_1, \ep_1)\\
=\ & (\sin\theta)(x_2 - x_1) + g_{\mu\ep_2}(x_1) - g_{\mu\ep_1}(x_1)\\
\geq\ & (\sin\theta)(x_2 - x_1)  - C_{\theta}|\ep_2 - \ep_1|.
\end{split}
\]
Rearranging to isolate $x_2 - x_1$ on one side, we obtain the following Lipschitz estimate for $\bx$ on $\{(t, \ep)\ |\ t > \ep\cot\theta > 0\}$:
\begin{equation}\label{eq:bx-Lipschitz}
|\bx(t_1, \ep_1) - \bx(t_2, \ep_2)| \leq C_{\theta} (|t_2 - t_1| + |\ep_2 - \ep_1|).
\end{equation}
Combining this with~\eqref{eq:psi-ep-definition}, and using both estimates in~\eqref{eq:g-ep-bounds}, we infer that
\begin{equation}\label{eq:psi-ep-Lipschitz}
|\psi_{\ep_1}(t_1) - \psi_{\ep_2}(t_2)| \leq C_{\theta} ( |t_1 - t_2| + |\ep_1 - \ep_2|),
\end{equation}
whenever $t_i > \ep_i \cot\theta > 0$. To continue, we differentiate the relation $\bt(\bx(t, \ep), \ep) = t$ with respect to $\ep$, and combine the result with~\eqref{eq:g-ep-bounds} to obtain upper and lower bounds on $\pa{\bx}{\ep}$. Upon also differentiating~\eqref{eq:psi-ep-definition} with respect to $\ep$, we get after some straightforward calculation that 
\begin{equation}\label{eq:d-psi-d-ep}
c_1 \leq \pa{\psi_{\ep}}{\ep} = -\pa{\bx}{\ep} \cdot \sec\theta \leq c_2,
\end{equation}
where $c_1, c_2$ depend only on $\theta$. To describe the boundary behavior, so to speak, of $\psi_{\ep}$, we use the strict monotonicity of $\bx(\cdot, 1)$ to get $\gamma_{1} < \gamma_{2}$ in $(\cot\theta, \infty)$ such that 
\[
\bx(\gamma_{1}, 1) = \frac{\mu}{4},\quad \bx(\gamma_{2}, 1) = \mu.
\]
Taking $x = \bx(\gamma_{1}, 1)$ and $\ep = 1$ in~\eqref{eq:bt-definition}, and recalling~\eqref{eq:g-particular}, we see that $\gamma_{1}$ is given more explicitly by
\begin{equation}\label{eq:gamma1-explicit}
\gamma_{1} = \cot\theta + \frac{\mu \sin\theta}{4},
\end{equation}
With $\gamma_{1}, \gamma_{2}$ thus chosen, we again use~\eqref{eq:g-particular} and~\eqref{eq:bt-definition}, as well as the definition of $\psi_{\ep}$, to find that
\begin{equation}\label{eq:psi-boundary-values}
\psi_{\ep}(t) = \left\{
\begin{array}{ll}
\ep \csc^2\theta - t\cot\theta, &\text{ if } \ep \cot\theta < t \leq \gamma_{1}\ep, \\
\ep\mu \ell\sin\theta, &  \text{ if }t \geq \gamma_{2}\ep.
\end{array}
\right.
\end{equation}

Now, by~\eqref{eq:bP-target} and the positivity of $\psi_{\ep}$, the map $\bP$ indeed takes values in $\mathring{\cV}$. Next, using~\eqref{eq:g-ep-bounds} and~\eqref{eq:g-ep-f}, and arguing as in the proof of Lemma~\ref{lemm:foliation-by-translated-mollification}, we see that the assignment
\[
(x,\ep) \mapsto (x, g_{\mu\ep}(x))
\]
defines a diffeomorphism from $\{x >0, \ep > 0\}$ onto $\mathring{\cV} \cap \{x > 0\}$, which together with~\eqref{eq:psi-g-relation} shows that 
\[
(t, \ep) \mapsto (\psi_{\ep}(t), t)
\]
is a diffeomorphism from $\{(t, \ep)\ |\ \ep > 0,\ t > \ep\cot\theta\}$ onto $\mathring{\cV} \cap \{x > 0\}$. From this we deduce that $\bP$ maps $\mathring{\cV} \cap \{x > 0\}$ and $\mathring{\cV} \cap \{x < 0\}$ diffeomorphically onto themselves. Taking also the second line of~\eqref{eq:psi-boundary-values} into account, which implies
\begin{equation}\label{eq:P-in-smaller-cone}
\bP(x, y) = ((\mu\ell \sin\theta) \cdot x, y) \quad\text{whenever }y \geq \gamma_{2}|x|,\ y \neq 0,
\end{equation}
we conclude that $\bP$ is a diffeomorphism from $\mathring{\cV}$ onto itself. Finally, the asserted Lipschitz continuity of $\bP$ follows from~\eqref{eq:psi-ep-Lipschitz},~\eqref{eq:P-in-smaller-cone}, and the convexity of $\mathring{\cV}$. As for the inverse, the expression~\eqref{eq:P-in-smaller-cone} implies that
\[
|\bP^{-1}(z_1) - \bP^{-1}(z_2)| \leq C_{\theta}|z_{1} - z_{2}|
\]
for all $z_{1}, z_{2} \in \{(x, y) \ |\ y> 0,\ y \geq \frac{\gamma_{2}|x|}{\mu\ell\sin\theta}\}$.
Using~\eqref{eq:d-psi-d-ep} and~\eqref{eq:psi-ep-Lipschitz}, and following how~\eqref{eq:inverse-bF-Lipschitz} was proved, we obtain a Lipschitz estimate of the above form when $z_{1}, z_{2} \in \mathring{\cV} \cap \{x > 0\}$, and when $z_{1}, z_{2} \in \mathring{\cV} \cap \{x < 0\}$. These last observations are enough to give the claimed Lipschitz continuity of $\bP^{-1}$. The proof is complete.
\end{proof}
\begin{rmk}\label{rmk:bP-extension}
By the first line in~\eqref{eq:psi-boundary-values}, we have whenever $x > 0$ and $\frac{y}{x} \in (\cot\theta, \gamma_{1}]$ that 
\begin{equation}\label{eq:bP-near-boundary-V}
\bP(x, y) = (x\csc^2\theta - y\cot\theta, y).
\end{equation}
With the help of~\eqref{eq:gamma1-explicit} and the inequality $\mu < \frac{\sec\theta}{4}$, it makes sense to define 
\begin{equation}\label{eq:b-gamma-defi}
b(\gamma): = \frac{\gamma}{\csc^2\theta - \gamma\cot\theta}\quad\text{for }\gamma \in (\cot\theta, \gamma_{1}].
\end{equation}
A direct computation using~\eqref{eq:bP-near-boundary-V} then shows that, whenever $x > 0$ and $\frac{y}{x} \in (\cot\theta, b(\gamma_1)]$, there holds
\begin{equation}\label{eq:inverse-bP-near-boundary-V}
\bP^{-1}(x, y) = ((x + y\cot\theta)\sin^2\theta, y).
\end{equation}
Consequently, for all $\gamma \in (\cot\theta, \gamma_1]$, we have
\begin{equation}\label{eq:bP-image-of-cones}
\bP\big(\{x > 0,\ \cot\theta < \frac{y}{x} \leq  \gamma\}\big) = \{x > 0,\ \cot\theta < \frac{y}{x} \leq b(\gamma)\}.
\end{equation}
Also, from~\eqref{eq:bP-near-boundary-V},~\eqref{eq:inverse-bP-near-boundary-V}, and their analogues for the case $x < 0$, as well as the Lipschitz bounds on $\bP$ and $\bP^{-1}$ proved above, it follows that extending both to be the identity on $\partial \cV$ yields a pair of Lipschitz maps from $\cV$ to itself that remain inverses of each other. Further, the extended maps, which we still denote by $\bP$ and $\bP^{-1}$, both restrict to be smooth on $\cV \setminus \{\0^2\}$.
\end{rmk}

\begin{rmk}\label{rmk:bP-scaling}
The map $\bP$, as extended in the previous remark, enjoys the following scaling property: for all $t \geq 0$ and $(x, y) \in \cV$, there holds
\begin{equation}\label{eq:bP-scaling-property}
\bP(tx, ty) = t\bP(x, y).
\end{equation}
Indeed, when $(x, y)\in \partial \cV$, this follows from~\eqref{eq:f-for-retraction} and the fact that $\bP|_{\partial\cV}$ is the identity. Next, when $(x, y) \in \mathring{\cV}$ and $x = 0$, we get the result straight from the definition of $\bP$. Finally, in the case $(x, y) \in \mathring{\cV}$ and $x \neq 0$, if $t = 0$ then~\eqref{eq:bP-scaling-property} holds trivially since both sides reduce to $(0, 0)$, while if $t > 0$, then we have by~\eqref{eq:psi-scaling-property} that
\[
\frac{tx}{|tx|} \psi_{|tx|}(ty) = \frac{x}{|x|} \cdot t|x|\psi_{1}(\frac{y}{|x|}) = t\cdot  \frac{x}{|x|} \psi_{|x|}(y),
\]
which implies~\eqref{eq:bP-scaling-property}.
\end{rmk}

To continue, define
\begin{equation}\label{eq:tilde-bF-definition}
\widetilde{\bF}(x, \ep) = \bP(\bF(x, \ep)), \quad (x, \ep) \in \RR \times [0,\infty).
\end{equation}
By Lemmas~\ref{lemm:foliation-by-translated-mollification} and~\ref{lemm:rotated-graph}, as well as Remark~\ref{rmk:bP-extension}, we see that $\widetilde{\bF}$ is a bi-Lipschitz map from $\RR \times [0, \infty)$ onto $\cV$, and restricts to a diffeomorphism from $\RR \times (0, \infty)$ onto $\mathring{\cV}$. Also, there exists $L_{\theta} > 1$ such that for any distinct points $(x_1, \ep_1), (x_2, \ep_2) \in \RR \times [0, \infty)$, there holds
\begin{equation}\label{eq:tilde-bF-biLip-bounds}
L_{\theta}^{-1} \leq \frac{|\widetilde{\bF}(x_1, \ep_1) - \widetilde{\bF}(x_2, \ep_{2})|}{|x_1 - x_2| + |\ep_1 - \ep_2|} \leq L_{\theta}.
\end{equation}
Recall that $f(1) = \cot\theta < \gamma_1$, and take $\delta_0 \in (0, 1)$ such that 
\begin{equation}\label{eq:delta-0-choice}
\gamma_0: = f_{\delta_0}(1) + \delta_0 < \gamma_1.
\end{equation}
With $\delta_0$ thus chosen, we define, for $\lambda = 1, 2, 3$,
\[
\begin{split}
K_{\lambda} :=\ & \{q \in \RR^2\ |\ |q - (q\cdot \bv_{\lambda})\bv_{\lambda}| \leq \delta_0 (q \cdot \bv_{\lambda})  \}\\
=\ & \{a\bv_{\lambda} + b\bw_{\lambda}\ |\ a \geq 0,\ |b| \leq \delta_0 a\},
\end{split}
\]
where $\bv_{\lambda}$ and $\bw_{\lambda}$ are defined at the start of the section. Note that since $\delta_0 < 1 < \tan\theta$, the sets $\{K_{\lambda}\}_{\lambda \in \{1, 2, 3\}}$ have pairwise disjoint interiors. 

Lemma~\ref{lemm:tilde-bF-boundary-behavior}, Remark~\ref{rmk:tilde-bF-image-cone}, and Remark~\ref{rmk:tilde-bF-image-K2} below contain additional properties of $\tilde{\bF}$ that we later use heavily.
\begin{lemm}\label{lemm:tilde-bF-boundary-behavior}
The map $\widetilde{\bF}$ has the following properties.
\vskip 1mm
\begin{enumerate}
\item[(a)] With $L_{\theta}$ as in~\eqref{eq:tilde-bF-biLip-bounds}, we have for all $r > 0$ that
\begin{equation}\label{eq:tilde-bF-image-of-dist-nbhd}
\cV \cap B_{L_{\theta}^{-1}r}(Y) \subset \widetilde{\bF}(\RR \times [0, r)) \subset \cV \cap B_{L_{\theta}r}(Y),
\end{equation}
and that
\begin{equation}\label{eq:tilde-bF-image-of-central-nbhd}
\cV \cap B_{L_{\theta}^{-1}r}(Y^{-})\subset \widetilde{\bF}((-r, r) \times [0, \infty)).
\end{equation}
\vskip 1mm
\item[(b)] Whenever $0 \leq \ep \leq \delta_0 x$, there holds
\begin{equation}\label{eq:tilde-bF-expression}
\begin{split}
\widetilde{\bF}(x, \ep) =\ & (x - \ep \cot\theta, x\cot\theta + \ep) = (x\csc\theta) \bv_{1} + (\ep\csc\theta)\bw_{1}.
\end{split}
\end{equation}
Consequently, 
\begin{equation}\label{eq:tilde-bF-image-cones}
\widetilde{\bF}(\{x \geq 0,\ 0 \leq \ep \leq \delta_0 x\}) = \cV \cap K_{1}.
\end{equation}
\end{enumerate}
\end{lemm}
\begin{proof}
First, it is elementary that
\[
\cV\cap B_{L_{\theta}^{-1}r}(Y) = \cV \cap B_{L_{\theta}^{-1}r}(\partial\cV),
\]
whereas
\[
\cV \cap B_{L_{\theta}^{-1}r}(Y^{-}) = \{(x, y) \in \cV\ |\ |x| < L_{\theta}^{-1}r\}.
\]
Now take any $(x, y) \in \cV \cap B_{L_{\theta}^{-1}r}(\partial\cV)$ and let $(X, \ep) = \widetilde{\bF}^{-1}(x, y)$. By definition there exists some $x_0 \in \RR$ such that 
\[
|(x, y) - (x_0, f(x_0))| < L_{\theta}^{-1}r. 
\]
Using the first inequality in~\eqref{eq:tilde-bF-biLip-bounds} then leads to
\[
\begin{split}
\ep \leq\ & L_{\theta}|\widetilde{\bF}(X, \ep) - \widetilde{\bF}(x_0, 0)| = L_{\theta}|(x, y) - (x_0, f(x_0))| < r,
\end{split}
\]
and the first inclusion in~\eqref{eq:tilde-bF-image-of-dist-nbhd} follows. For the second inclusion, we simply note that 
\[
|\widetilde{\bF}(x, \ep) - (x, f(x))| = |\widetilde{\bF}(x, \ep) - \widetilde{\bF}(x, 0)| \leq L_{\theta}\ep,
\]
where the last step uses the second inequality in~\eqref{eq:tilde-bF-biLip-bounds}. Since $\widetilde{\bF}(x,\ep) \in \cV$ to start with, we are done with~\eqref{eq:tilde-bF-image-of-dist-nbhd}. 

Next we prove~\eqref{eq:tilde-bF-image-of-central-nbhd}. Suppose $(x, y) \in \cV$ is such that $|x| < L_{\theta}^{-1}r$, and again write $(X, \ep) = \widetilde{\bF}^{-1}(x, y)$, the goal this time being to show that $|X| <r $. If $x = 0$, then $X = 0$ and we are done. If $x \neq 0$, in which case $y > 0$, then upon letting $\delta = \frac{y}{f_1(0) + 1}$, so that $y = f_{\delta}(0) + \delta$, we have
\[
L_{\theta}^{-1}r > |(x, y) - (0, y)| = |\widetilde{\bF}(X, \ep) - \widetilde{\bF}(0, \delta)| \geq L_{\theta}^{-1}|X|,
\]
where the last step uses~\eqref{eq:tilde-bF-biLip-bounds}. Thus we have proved that~\eqref{eq:tilde-bF-image-of-central-nbhd} holds. 

For part (b), we first claim that, with $\gamma_0$, $\delta_0$ as above, there holds
\begin{equation}\label{eq:bF-images-of-cones}
\bF(\{x \geq 0,\ 0 \leq \ep \leq \delta_0 x\}) = \{x \geq 0,\ (\cot\theta)x \leq y \leq \gamma_0x\}.
\end{equation}
Indeed, suppose $0 \leq \ep \leq \delta_0 x$ and let $(x, y) = \bF(x, \ep)$. In the case $\ep = 0$, we simply have $y = f(x) =  (\cot\theta)x$, while if $\ep > 0$, in which case $x > 0$ also, we use the monotonicity of $\ep \mapsto f_{\ep}(x) + \ep$ along with the scaling property~\eqref{eq:f-ep-scaling-property} to see that
\[
\begin{split}
(\cot\theta)x < y = f_{\ep}(x) + \ep \leq\ & f_{\delta_0x}(x) + \delta_0 x\\
=\ & x\cdot(f_{\delta_0}(1) + \delta_0) = \gamma_0 x.
\end{split}
\]
Conversely, suppose $0 \leq (\cot\theta)x \leq y\leq \gamma_0 x$, and let $(x, \ep) = \bF^{-1}(x, y)$. If $y = (\cot\theta)x$ then necessarily $\ep = 0$. On the other hand, if $y > (\cot\theta)x$, then both $x > 0$ and $\ep > 0$, so by the equality $\gamma_0 x = f_{\delta_0 x}(x) + \delta_0 x$ computed above, and the monotonicity noted before it, we must have $\ep \leq \delta_0 x$. Thus we have proved~\eqref{eq:bF-images-of-cones}, and in fact the argument shows that 
\[
\bF(\{0 < \ep \leq \delta_0 x\}) = \{(\cot\theta)x < y \leq \gamma_0x\}.
\]
We proceed to verify~\eqref{eq:tilde-bF-expression}. Given $(x, \ep)$ such that $0 < \ep \leq \delta_0 x$, since $\gamma_0 < \gamma_1$, we see from the above equality and~\eqref{eq:bP-near-boundary-V} that
\[
\widetilde{\bF}(x, \ep) = (x\csc^2\theta - (f_{\ep}(x) + \ep)\cot\theta, f_{\ep}(x) + \ep).
\]
Recalling that $\delta_0 < 1$, we obtain~\eqref{eq:tilde-bF-expression} in the case $\ep > 0$ upon substituting~\eqref{eq:f-ep-agree} into the above. The case $\ep = 0$ is obvious from~\eqref{eq:f-for-retraction} and the fact that $\bP|_{\partial\cV}$ is the identity. Finally, we prove~\eqref{eq:tilde-bF-image-cones}. The inclusion ``$\subset$'' follows straight from the second equality in~\eqref{eq:tilde-bF-expression}. For the reverse inclusion, take $q \in \cV \cap K_{1}$ and write 
\[
q = a\bv_{1} + b\bw_{1},
\]
for some $a, b \in \RR$. Since $q \in K_{1}$, we have $|b| \leq \delta_0a$, while $q \in \cV$ implies $b \geq 0$. Thus, letting $x = a\sin\theta$ and $\ep = b\sin\theta$, we get $0 \leq \ep \leq \delta_0 x$, and hence
\[
q = x\cdot (1, \cot\theta) + \ep \cdot (-\cot\theta, 1) = \widetilde{\bF}(x, \ep).
\]
This gives the inclusion ``$\supset$'' in~\eqref{eq:tilde-bF-image-cones}. 
\end{proof}
\begin{rmk}\label{rmk:tilde-bF-image-cone}
Considering what part (b) and~\eqref{eq:tilde-bF-image-of-central-nbhd} together imply for points in $\cV \cap K_1 \cap B_{(4L_{\theta})^{-1}}(Y^{-})$, and extending that to $K_1 \cap B_{(4L_{\theta})^{-1}}(Y^{-})$ by reflection across $\Span\{\bv_{1}\}$, we deduce that any $q \in K_1 \cap B_{(4L_{\theta})^{-1}}(Y^{-})$ has the form
\[
q = a\bv_{1} + b\bw_{1},
\]
for some $a, b$ satisfying that $0 \leq a\sin\theta < \frac{1}{4}$, and that $|b| \leq \delta_0 a$.
\end{rmk}
\begin{rmk}\label{rmk:tilde-bF-image-K2}
From the definition of $\bP$ and the fact that each $f_{\ep}$, as well as $f$ itself, is an even function, it follows that for all $(X, \ep) \in \RR \times [0,\infty)$, if we let $(x, y): = \widetilde{\bF}(X, \ep)$, then 
\[
\widetilde{\bF}(-X, \ep) = (-x, y).
\]
Combining this with Lemma~\ref{lemm:tilde-bF-boundary-behavior}(b) shows that, whenever $0 \leq \ep \leq -\delta_0 X$, there holds
\[
\widetilde{\bF}(X, \ep) = (X + \ep\cot\theta, -X\cot\theta + \ep) = -(X\csc\theta)\bv_2 - (\ep\csc\theta)\bw_{2},
\]
and, since $\cV \cap K_2$ happens to be the image of $\cV \cap K_1$ under $(x, y) \mapsto (-x, y)$, that
\[
\widetilde{\bF}(\{X \leq 0,\ 0 \leq \ep \leq -\delta_0 X\}) = \cV \cap K_2.
\]
\end{rmk}

We next construct what ends up modeling the cross section of $\Omega$ when intersected with $m$-planes of the form $(a\bv_{\lambda} + \Span\{\bw_{\lambda}\}) \times \RR^{m-1}$ for sufficiently large $a > 0$ (see Lemma~\ref{lemm:Omega-S-basics}(b) below). Writing
\[
|w|_{\infty}: = \max\{|w^{1}|, \cdots, |w^{m}|\},\quad\text{for }w \in \RR^{m},
\]
and denoting by $\eta$ the mollifier introduced before~\eqref{eq:f-ep-definition} in the case $k - 1= m$, we define
\[
N(\cdot) := \eta * |\cdot|_{\infty}\quad\text{on }\RR^{m}.
\]
Similar to~\eqref{eq:f-ep-f-bound}, by the convexity of $|\cdot|_{\infty}$ and the fact that $|\cdot|_{\infty} \leq |\cdot|$, we have
\begin{equation}\label{eq:reg-norm-bound}
|w|_{\infty} \leq N(w) \leq |w|_{\infty} + 1.
\end{equation}
Also, $N$ is a convex function, and is even with respect to each coordinate. The latter means that
\begin{equation}\label{eq:N-even-in-each}
N(w^1, \cdots, -w^{i}, \cdots, w^{m}) = N(w^1, \cdots, w^{i}, \cdots, w^{m}),
\end{equation}
for each $i \in \{1, \cdots, m\}$. Next, given $i, j \in \{1, \cdots, m\}$, we have
\[
\partial_{i} (|\cdot|_{\infty}) = \pm \delta_{ij} \quad\text{a.e. on }\{v \in \RR^{m}\ |\ |v|_{\infty} = \pm v_{j}\}.
\]
Substituting this into~\eqref{eq:grad-f-ep} yields
\begin{equation}\label{eq:grad-N-formula}
\partial_{i}N(w) = \int_{B_{1}(w) \cap \{|v|_{\infty} = v^{i}\}} \eta(w - v) dv - \int_{B_{1}(w) \cap \{|v|_{\infty} = -v^{i}\}} \eta(w - v) dv.
\end{equation}
\begin{lemm}\label{lemm:N-gradient-radial}
For all $w \in \RR^m$ such that $|w|_{\infty} > 3$, there holds
\begin{equation}\label{eq:N-gradient-radial}
w \cdot \nabla N(w) > 0.
\end{equation}
\end{lemm}
\begin{proof}
If $i \in \{1, \cdots, m\}$ is such that $|w^i| < 1$, then by the triangle inequality, we have for all $v \in B_{1}(w)$ that $|v^{i}| < 2 < |v|_{\infty}$, so that 
\[
B_{1}(w) \cap \{|v|_{\infty} = \pm v^{i}\} = \emptyset.
\]
On the other hand, if $|w^{i}| \geq 1$, then $v^i$ has the same sign as $w^i$ for all $v \in B_{1}(w)$, so that 
\[
B_{1}(w) \cap \{|v|_{\infty} = -\frac{w^i}{|w^i|}\cdot v^i\} = \emptyset.
\]
Combining these observations with~\eqref{eq:grad-N-formula} and summing over $i$ gives $w \cdot \nabla N(w) \geq 0$. 

To get strict inequality, choose $i$ such that $|w|_{\infty} = |w^i|$, and assume without loss of generality that $w^i > 0$. Then, as already mentioned, we have
\[
B_{1}(w) \cap \{|v|_{\infty} = -v^i\} = \emptyset.
\]
Noting also that
\[
B_{\frac{\delta}{2}}(w + \delta e_{i}) \subset B_{1}(w) \cap \{|v|_{\infty} = v^i\}, 
\]
for all sufficiently small $\delta >0$, and that $\eta$ is positive on $B_1$, we deduce from~\eqref{eq:grad-N-formula} that $w^i \pa{N}{w^i}(w) > 0$. Consequently, the inequality we got at the end of the previous paragraph must be strict. The proof is complete.
\end{proof}
Combining~\eqref{eq:reg-norm-bound}, the convexity of $N$, and Lemma~\ref{lemm:N-gradient-radial}, we see that
\begin{equation}\label{eq:cC-for-retraction}
\cC := \{w \in \RR^{m}\ |\ N(w) \leq 5\}
\end{equation}
is a non-empty, compact, convex set with smooth boundary, the latter given by $\{N(w) = 5\}$. More specifically,~\eqref{eq:reg-norm-bound} implies
\begin{equation}\label{eq:cC-inclusions}
\cC \subset \{|w|_{\infty} \leq 5\},\quad \partial\cC \subset \{4 \leq |w|_{\infty} \leq 5\}.
\end{equation}

At this point we digress to mention the following standard fact, which we use a number of times in the remainder of this section.
\begin{lemm}\label{lemm:boundary-IVT}
Let $X$ and $W$ be smooth $d$-dimensional manifolds with boundary, and suppose $F: X \to W$ is a smooth, injective immersion satisfying $F(\partial X) \subset \partial W$. Then $F$ maps $X$ diffeomorphically onto an open set in $W$.
\end{lemm}
\begin{proof}
It suffices to prove that $F$ is a local diffeomorphism. Given $x_0 \in X$, if it is an interior point, then since $dF_{x_0}$ is invertible, the inverse function theorem guarantees that $F(x_0) \in W \setminus\partial W$ (see~\cite[page 96, Problem 4-2]{Lee2013}), and that $F$ is a local diffeomorphism at $x_0$. Next, if $x_0$ is a boundary point, in which case $F(x_0) \in \partial W$ by assumption, upon introducing coordinate charts, we obtain a neighborhood $U$ of $0$ in $\RR^{d}$, along with a smooth map $f$ on $U \cap \{x^{d} \geq 0\}$ with values in $\{x^{d} \geq 0\}$, such that $f(0) = 0$, that $(df)_{0}$ is invertible, and that 
\begin{equation}\label{eq:f-preserve-half-space}
f(U \cap \{x^{d} = 0\}) \subset \{x^{d} = 0\},\quad f(U \cap \{x^{d} > 0\}) \subset \{x^{d} > 0\},
\end{equation}
where the second condition follows from our initial observation that $F$ preserves interior points. Our task is then reduced to finding a pair of neighborhoods of the origin in $\{x^{d} \geq 0\}$ so that $f$ maps one diffeomorphically onto the other. To that end, we first use the smoothness of $f$ (in the sense of~\cite[page 34]{Lee2013}) and the inverse function theorem to get, up to shrinking $U$, a diffeomorphism $\tilde{f}$ from $U$ onto a neighborhood $V \subset \RR^{d}$ of the origin such that 
\begin{equation}\label{eq:tilde-f-f-agree}
\tilde{f} = f\quad \text{on }U \cap \{x^{d} \geq 0\}.
\end{equation}
On the other hand, the first inclusion in~\eqref{eq:f-preserve-half-space}, along with the invertibility of $(df)_0$, implies that the inverse function theorem is also applicable, at the origin, to $f|_{U\cap \{x^{d} = 0\}}$ viewed as a map into $\{x^{d} = 0\}$. Hence we get open sets $U'$ and $V'$ in $\RR^{d}$, both containing the origin, such that 
\[
f: U' \cap \{x^{d} = 0\} \to V' \cap \{x^{d} = 0\},\quad\text{diffeomorphically}.
\]
Assuming, without loss of generality, that $U' \subset U$, and choosing an open ball $B \subset \RR^{d}$ centered at $0$ such that $B \subset U' \cap \tilde{f}^{-1}(V')$, we claim that
\begin{equation}\label{eq:tilde-f-preserve-half-space}
\tilde{f}(B \cap \{x^{d} = 0\}) = \tilde{f}(B) \cap \{x^{d} = 0\}.
\end{equation}
The inclusion ``$\subset$'' is a straightforward consequence of~\eqref{eq:tilde-f-f-agree} and the first fact in~\eqref{eq:f-preserve-half-space}. For the inclusion ``$\supset$'', take $y \in \tilde{f}(B) \cap \{x^{d} = 0\}$ and define 
\[
x_1 = \tilde{f}^{-1}(y), \quad\quad x_2 = (f|_{U' \cap \{x^{d} = 0\}})^{-1}(y),
\]
where the definition for $x_2$ makes sense since $\tilde{f}(B) \subset V'$. Since both $x_1, x_2$ belong to $U$, with $x_2$ lying additionally in $\{x^{d} = 0\}$, we have
\[
\tilde{f}(x_2) = f(x_2) = y = \tilde{f}(x_1).
\]
The injectivity of $\tilde{f}$ on $U$ then forces $x_1$ to coincide with $x_2$, so that $\tilde{f}^{-1}(y) \in B \cap \{x^{d} = 0\}$. This finishes the proof of~\eqref{eq:tilde-f-preserve-half-space}, and as a result we get
\[
\{x^{d} = 0\} \cap \tilde{f}(B \cap \{x^{d} < 0\}) = \tilde{f}(B \cap \{x^{d} = 0\}) \cap \tilde{f}(B \cap \{x^{d} < 0\})  = \emptyset.
\]
The connectedness of $\tilde{f}(B \cap \{x^{d} < 0\})$ then implies that it is contained either in $\{x^{d} < 0\}$ or in $\{x^{d} > 0\}$. As the latter alternative is ruled out by~\eqref{eq:f-preserve-half-space}, or rather what those inclusions imply for $\pa{f^{d}}{x^{d}}(0)$, upon using~\eqref{eq:f-preserve-half-space} once more and recalling~\eqref{eq:tilde-f-f-agree}, we arrive at
\[
f(B \cap \{x^{d} \geq 0\}) = \tilde{f}(B \cap \{x^{d} \geq 0\}) = \tilde{f}(B) \cap \{x^{d} \geq 0\}.
\]
We conclude the proof upon noting that the right-most set is open in $\{x^{d} \geq 0\}$, and that $f:B \cap \{x^{d} \geq 0\} \to \tilde{f}(B) \cap \{x^{d} \geq 0\}$ is smoothly inverted by $(\tilde{f}^{-1})|_{\tilde{f}(B) \cap \{x^{d} \geq 0\}}$.\\
\end{proof}
Returning from the above digression, we next establish some additional properties of $\cC$.
\begin{lemm}\label{lemm:cC-properties}
The set $\cC$ has the following additional properties.
\vskip 1mm
\begin{enumerate}
\item[(a)] $\cC': = \{w' \in \RR^{m-1}\ |\ N(0, w') \leq 5\}$ is a non-empty convex set with smooth boundary in $\RR^{m-1}$, and 
\[
\cC \cap \{|w^{1}| < 1\} = (-1, 1) \times \cC', \quad 
\partial \cC \cap \{|w^{1}| < 1\} = (-1, 1) \times \partial \cC'.
\]
\vskip 1mm
\item[(b)] The assignment $(r, w) \mapsto rw$ restricts to a Lipschitz map from $[0, 1] \times \partial \cC$ onto $\cC$, and a diffeomorphism from $(0, 1] \times \partial \cC$ onto $\cC \setminus \{\0^{m}\}$.
\end{enumerate}
\end{lemm}
\begin{proof}
From~\eqref{eq:reg-norm-bound} we see that $\cC'$ is non-empty. That $\cC'$ is compact and convex is also clear, since $\cC$ has these properties and $\{0\} \times \cC' = \cC \cap (\{0\} \times \RR^{m-1})$. That $\partial\cC'$ is smooth and coincides with $\{w' \in \RR^{m-1}\ |\ N(0, w') = 5\}$ follows from applying Lemma~\ref{lemm:N-gradient-radial} to see that 
\begin{equation}\label{eq:slice-regular}
w' \cdot \nabla (N(0, \cdot)) > 0,\quad \text{whenever }N(0, w') = 5.
\end{equation}
To continue, given a point $w \in \RR^{m}$, we write $w'$ to mean the element of $\RR^{m-1}$ obtained by dropping $w^{1}$. Suppose $w \in \RR^{m}$ is such that $|w^{1}| < 1$ and $|(0, w')|_{\infty} > 3$. Then, whenever $v \in B_{1}(w)$, we have by the triangle inequality that $|v^{1}| < |(0, v')|_{\infty}$, and hence $|v|_{\infty} = |(0, v')|_{\infty}$. Substituting this into the definition of $N$ yields
\[
\begin{split}
N(w) =\ & \int_{\RR^{m}} \eta(w^{1} - v^{1}, w' - v')|(0, v')|_{\infty}\, dv \\
=\ & \int_{\RR^{m-1}} \big(\int_{\RR}\eta(t, w' - v') dt\big) |(0, v')|_{\infty}\, dv'.
\end{split}
\]
From this we get
\begin{equation}\label{eq:N-independent-of-w1}
N(w) = N(0, w')\quad\text{whenever }|w^1| < 1,\ |(0, w')|_{\infty} > 3,
\end{equation}
which together with~\eqref{eq:reg-norm-bound} gives
\[
\{|w^{1}| < 1\} \setminus \cC = (-1, 1) \times (\RR^{m-1}\setminus \cC'),
\]
proving the first asserted equality. The second likewise follows from~\eqref{eq:N-independent-of-w1} and~\eqref{eq:reg-norm-bound}, and we are done with (a).

For part (b), since $\cC$ is convex, and contains $\0^{m}$ in its interior by~\eqref{eq:reg-norm-bound}, we see that $rw \in \cC$ for all $(r, w) \in [0, 1] \times \partial \cC$, and that $rw \neq \0^{m}$ in the case $r > 0$. Given $w \in \cC\setminus  \{\0^{m}\}$, again by~\eqref{eq:reg-norm-bound}, as well as the intermediate value theorem, there exists $r \in (0, 1]$ such that $N(r^{-1}w) = 5$, so that $w \in r \cdot \partial\cC$. On the other hand, let $(r_1, w_1), (r_2, w_2) \in (0, 1] \times \partial\cC$ be such that 
\[
r_1w_1 = r_2 w_2,
\]
and assume without loss of generality that $r_2 \geq r_1$. Then since 
\[
N(w_2) = 5 = N(w_1) = N(\frac{r_2}{r_1}w_2),
\]
and since $t \mapsto N(tw_2)$ is strictly increasing on $[1, \infty)$ thanks to Lemma~\ref{lemm:N-gradient-radial}, we must have $r_2 = r_1$, in which case $w_1 = w_2$. Thus we have shown that $(r,w) \mapsto rw$ takes $(0, 1] \times\partial\cC$ bijectively to $\cC\setminus \{\0^m\}$, and consequently also takes $[0, 1] \times \partial \cC$ onto $\cC$. Lipschitz continuity follows with the help of~\eqref{eq:cC-inclusions}, while Lemma~\ref{lemm:N-gradient-radial} gives invertibility of the derivative on $(0, 1] \times \partial \cC$. Since the map in question clearly takes $\{1\} \times \partial\cC$ bijectively to $\partial\cC$, an application of Lemma~\ref{lemm:boundary-IVT} completes the proof.
\end{proof}
\begin{rmk}\label{rmk:cC'-retraction}
Combining the conclusion of part (b) with the relations 
\[
\{0\} \times \cC' = \cC \cap (\{0\} \times \RR^{m-1}),\quad \{0\} \times \partial\cC' = \partial\cC \cap (\{0\} \times \RR^{m-1}),
\]
and using also~\eqref{eq:slice-regular}, we see that the assignment
\[
(r, w') \mapsto rw'
\]
is a Lipschitz map from $[0,1] \times \partial \cC'$ onto $\cC'$ which takes $(0, 1] \times \partial\cC'$ diffeomorphically onto $\cC'\setminus \{\0^{m-1}\}$.
\end{rmk}

For the next result, we take the same cut-off function $\varphi$ as in~\eqref{eq:central-cutoff}, and define 
\begin{equation}\label{eq:xi-for-retraction}
\xi(x, t) = t \varphi(|x|)x + (1-t)x,\quad (x, t) \in \RR \times [0, 1].
\end{equation}
In particular, notice that for all $t \in [0, 1]$, we have
\[
\xi(x, t) = 
\left\{
\begin{array}{ll}
x,& \text{ if }|x| \geq 1,\\
(1-t)x,& \text{ if }|x| \leq 1/4.
\end{array}
\right.
\]
\begin{lemm}\label{lemm:x-e-z-reparametrization}
The correspondence 
\begin{equation}\label{eq:x-e-z-reparametrization}
(x, t, w) \mapsto (\xi(x, t), (1-t)w)
\end{equation}
defines a Lipschitz map from $\RR \times [0, 1] \times \partial\cC$ onto $\RR \times \cC$, and moreover restricts to a diffeomorphism from $\RR \times [0, 1) \times \partial \cC$ onto $\RR \times (\cC \setminus \{\0^{m}\})$.
\end{lemm}
\begin{proof}
By direct computation we have
\[
\pa{\xi}{x} =t\cdot [ \varphi'(|x|)|x| + \varphi(|x|) ] + (1-t),\quad \pa{\xi}{t} = [\varphi(|x|)-1]x,
\]
and consequently, for $(x, t) \in \RR \times [0,1]$, 
\begin{equation}\label{eq:xi-partial-derivatives}
1-t \leq \pa{\xi}{x} \leq 4,\quad \quad |\pa{\xi}{t}| \leq 1.
\end{equation}
Combining these with~\eqref{eq:cC-inclusions}, we infer that the map in question is Lipschitz on $\RR \times [0, 1] \times \partial \cC$. Moreover, by the intermediate value theorem, each $\RR \times \{t\} \times \partial\cC$ is mapped onto $\RR \times (1-t)\cdot\partial \cC$, which along with Lemma~\ref{lemm:cC-properties}(b) implies that~\eqref{eq:x-e-z-reparametrization} takes $\RR \times [0, 1] \times \partial \cC$ onto $\RR \times \cC$, and takes $\RR \times [0, 1) \times \partial \cC$ onto $\RR \times (\cC\setminus \{\0^{m}\})$. 

Next, note from~\eqref{eq:xi-partial-derivatives} that $\pa{\xi}{x} \geq 1-t > 0$ on $\RR \times [0, 1)$. Thus, again using Lemma~\ref{lemm:cC-properties}(b), we get that~\eqref{eq:x-e-z-reparametrization} is an injective on $\RR \times [0, 1) \times \partial \cC$, and that its derivative is invertible everywhere on this latter set. As in Lemma~\ref{lemm:cC-properties}, we conclude the proof upon noting that~\eqref{eq:x-e-z-reparametrization} is a bijection from $\RR \times \{0\} \times \partial \cC$ to $\RR \times \partial\cC$, and using Lemma~\ref{lemm:boundary-IVT}. 
\end{proof}
Expressing $\widetilde{\bF}^{-1}$ in components as
\begin{equation}\label{eq:tilde-bF-inverse-components}
\widetilde{\bF}^{-1}(x, y) = (X(x, y), \tilde{\ep}(x, y)),
\end{equation}
we define
\[
\begin{split}
U =\ & \{(x, y, z) \in \cV \times \RR^{m-1}\ |\ (\tilde{\ep}(x, y), z) \in \cC\},\\
\Sigma =\ & \{(x, y, z) \in \cV \times \RR^{m-1}\ |\ (\tilde{\ep}(x, y), z) \in \partial \cC\},
\end{split}
\]
and write 
\[
\mathring{U} = U \cap (\mathring{\cV} \times \RR^{m-1}), \quad \mathring{\Sigma} = \Sigma \cap (\mathring{\cV} \times \RR^{m-1}).
\]
\begin{lemm}\label{lemm:U-Sigma-basics}
The above sets have the following basic properties.
\vskip 1mm
\begin{enumerate}
\item[(a)] The assignment
\begin{equation}\label{eq:tilde-bF-as-parametrization}
(x, \ep, z)  \mapsto (\widetilde{\bF}(x, \ep), z)
\end{equation}
is a bi-Lipschitz map from $\RR \times \{(\ep, z)\ |\ \ep \geq 0,\ (\ep, z) \in \partial \cC\}$ onto $\Sigma$, and a diffeomorphism from $\RR \times \{(\ep, z)\ |\ \ep > 0,\ (\ep, z) \in \partial \cC\}$ onto $\mathring{\Sigma}$. 
\vskip 1mm
\item[(b)] Part (a) continues to hold with both occurrences of $\partial\cC$ replaced by $\cC$, and with $\Sigma$ and $\mathring{\Sigma}$ replaced by $U$ and $\mathring{U}$, respectively. 
\vskip 1mm
\item[(c)] For all $r \in (0, 1]$, we have 
\[
\Sigma \cap (B_{L_{\theta}^{-1}r}(Y) \times \RR^{m-1})= (B_{L_{\theta}^{-1}r}(Y) \cap \cV) \times \partial \cC',
\]
and similarly 
\[
U \cap (B_{L_{\theta}^{-1}r}(Y) \times \RR^{m-1}) = (B_{L_{\theta}^{-1}r}(Y) \cap \cV) \times \cC'.
\]
\end{enumerate}
\end{lemm}
\begin{proof}
For the statement about $\Sigma$ in (a), we note that~\eqref{eq:tilde-bF-as-parametrization} is a bi-Lipschitz map from $\RR \times [0, \infty) \times \RR^{m-1}$ onto $\cV \times \RR^{m-1}$ by the discussion after~\eqref{eq:tilde-bF-definition}, and that $\Sigma$ is indeed the image of $\RR \times \{(\ep, z)\ |\ \ep \geq 0,\ (\ep, z) \in \partial \cC\}$. The statement about $\mathring{\Sigma}$ is proved in a similar way, except we use the fact that~\eqref{eq:tilde-bF-as-parametrization} takes $\RR \times (0, \infty) \times \RR^{m-1}$ diffeomorphically onto $\mathring{\cV}\times \RR^{m-1}$. This proves part (a). Similar observations give part (b), and we omit the details.

For (c), we first use part (a) together with Lemma~\ref{lemm:cC-properties}(a) to get
\[
\Sigma \cap (\widetilde{\bF}(\RR \times [0, 1)) \times \RR^{m-1} ) = \widetilde{\bF}(\RR \times [0, 1)) \times \partial \cC'.
\]
Recalling the first inclusion in~\eqref{eq:tilde-bF-image-of-dist-nbhd} and the assumption $r \leq 1$, we deduce that
\[
\Sigma\cap ((\cV \cap B_{L_{\theta}^{-1}r}(Y)) \times \RR^{m-1}) = (\cV \cap B_{L_{\theta}^{-1}r}(Y)) \times \partial \cC',
\]
which implies the first asserted equality since $\Sigma \subset \cV \times \RR^{m-1}$ to start with. The proof of the second one is similar, and again we omit the details.
\end{proof}

By Lemma~\ref{lemm:x-e-z-reparametrization} and Lemma~\ref{lemm:U-Sigma-basics}(a)(b), we obtain a well-defined map 
\[
\bg: [0, 1] \times \Sigma \to U
\]
by requiring that 
\begin{equation}\label{eq:bg-definition}
\bg(t, \widetilde{\bF}(x, \ep), z) = (\widetilde{\bF}(\xi(x, t), (1-t)\ep), (1-t) z),
\end{equation}
for all $t \in [0, 1]$, $x \in \RR$, and $(\ep, z) \in \partial \cC$ such that $\ep \geq 0$. Still by the previous two lemmas, we have $\bg([0, 1] \times \Sigma) = U$, and that 
\begin{equation}\label{eq:bg-lambda-diffeo}
\mathring{U} = \bg([0, 1) \times \mathring{\Sigma}) \quad\text{diffeomorphically.}
\end{equation}
Moreover, using $\widetilde{\bF}(0, 0) = (0, 0)$ along with the estimate~\eqref{eq:tilde-bF-biLip-bounds}, and observing from~\eqref{eq:xi-for-retraction} that $|\xi(x, t)| \leq |x|$, we get for all $R > 0$ that
\begin{equation}\label{eq:bg-scale-factor}
\begin{split}
\bg([0, 1] \times (\Sigma \cap (B_R \times \RR^{m-1}))) \subset\ & U \cap (B_{R\cdot L_\theta^2} \times \RR^{m-1})\\
\subset\ & (\cV \cap B_{R\cdot L_\theta^2}) \times \{|z|_{\infty} \leq 5\},
\end{split}
\end{equation}
where $B_{r} \subset \RR^{2}$ denotes an open disk centered at $\0^{2}$, and for the second inclusion we used~\eqref{eq:cC-inclusions} and the definition of $U$. To continue, note also that $\bg$ is a Lipschitz map, as can be seen by using~\eqref{eq:tilde-bF-biLip-bounds} and the Lipschitz continuity of the map~\eqref{eq:x-e-z-reparametrization}, followed by another application of~\eqref{eq:tilde-bF-biLip-bounds}. The resulting estimate is
\begin{equation}\label{eq:bg-Lipschitz}
|\bg(t_1, p_1) - \bg(t_2, p_2)| \leq C_{\theta}(|t_1 - t_2| + |p_1 -p_2|),
\end{equation}
for all $(t_1, p_1), (t_2, p_2) \in [0, 1] \times \Sigma$. Next, tracing the definitions, we have
\begin{equation}\label{eq:bg-for-t=1}
\bg(\{1\} \times \Sigma)  = \partial\cV \times \{\0^{m-1}\}.
\end{equation}
On the other hand, given $q \in \partial \cV$ and $z \in \partial\cC'$, upon writing $q = (x, f(x)) = \widetilde{\bF}(x, 0)$ for some unique $x \in \RR$, so that $|x| = |q|\sin\theta$, and also expressing~\eqref{eq:xi-for-retraction} as 
\[
\xi(x, t) = [t\varphi(|x|) + (1 - t)]x,
\]
so as to use the homogeneity property of $f$, we have
\begin{equation}\label{eq:bg-on-boundary}
\begin{split}
\bg(t, q, z)  =\ & (\xi(x, t), f(\xi(x,t)), (1-t)z)\\
=\ & ((1-t)q + t\varphi(|q|\sin\theta)q, (1-t) z).
\end{split}
\end{equation}
From this we infer that
\begin{equation}\label{eq:bg-image-boundary}
\bg(\{t\} \times \partial \cV \times \partial\cC') = \partial\cV \times (1-t)\cdot \partial\cC' \quad\text{for all }t \in [0, 1],
\end{equation}
and, using in addition Remark~\ref{rmk:cC'-retraction} and the lower bound for $\pa{\xi}{x}$ in~\eqref{eq:xi-partial-derivatives}, that
\begin{equation}\label{eq:bg-boundary-injective}
\bg|_{[0, 1) \times \partial\cV \times \partial\cC'}\quad \text{is injective}.
\end{equation}
\begin{prop}\label{prop:bg-behavior}
Let $\gamma_0$ and $\delta_0$ be the constants chosen before Lemma~\ref{lemm:tilde-bF-boundary-behavior}. Also, for $\lambda \in \{1, 2, 3\}$, let 
\[
\pi_{\lambda}:\RR^{2} \to \RR^{2}
\]
denote orthogonal projection onto $\Span\{\bv_{\lambda}\}$.
\vskip 1mm
\begin{enumerate}
\item[(a)] For all $t \in [0, 1]$ and $(q, z) \in \Sigma \cap (K_1 \times \RR^{m-1})$, where $K_{1}$ is defined above Lemma~\ref{lemm:tilde-bF-boundary-behavior}, we have 
\begin{equation}\label{eq:bg-near-boundary}
\bg(t, q, z) = ((1-t)q + t\varphi(|\pi_{1}(q)|\sin\theta) \pi_{1}(q),\ (1-t)z).
\end{equation}
\vskip 1mm
\item[(b)] For all $t \in [0, 1]$ and $(q, z) \in \Sigma \cap (B_{(4L_{\theta})^{-1}}(Y^{-}) \times \RR^{m-1})$, we have
\begin{equation}\label{eq:bg-near-center}
\bg(t, q, z) = ((1-t)q, (1-t)z).
\end{equation}
\vskip 1mm
\item[(c)] There exists $C_{\theta} > 0$ such that given $t \in [0, 1)$, upon writing $\bg_{t}$ for $\bg(t, \cdot)$, we have, for all $p \in \mathring{\Sigma}$ and $2$-plane $E \subset T_{p}\mathring{\Sigma}$, that
\[
|(d\bg_{t})_{p}(u_{1}) \wedge (d\bg_{t})_{p}(u_{2})| \leq C_{\theta}(1-t),
\]
where $u_1, u_2$ is any orthonormal basis of $E$, the norm on the left-hand side being independent of this choice.
\end{enumerate}
\end{prop}
\begin{proof}
For convenience, in this proof we write
\[
\mu: = \csc\theta.
\]
For part (a), since $q \in \cV \cap K_{1}$, upon letting $(X, \ep) = \widetilde{\bF}^{-1}(q)$, we obtain from Lemma~\ref{lemm:tilde-bF-boundary-behavior}(b) that $0 \leq \ep \leq \delta_0 X$, so that~\eqref{eq:tilde-bF-expression} gives
\[
(q, z) = (X\mu\bv_1 + \ep\mu\bw_1, z).
\]
Recalling from~\eqref{eq:xi-for-retraction} that, in the present situation,
\[
\xi(X, t) = t\varphi(|X|)X + (1-t)X \geq \frac{(1-t)\ep}{\delta_0},
\]
and again using~\eqref{eq:tilde-bF-expression}, we have
\[
\begin{split}
\bg(t, q, z) =\ & (\widetilde{\bF}(\xi(X, t), (1-t)\ep), (1-t)z)\\
=\ & (\xi(X,t)\mu\bv_{1} + (1-t)\ep\mu\bw_1 , (1-t)z)\\
=\ & ((1-t)q + t\varphi(|X|)X\mu\bv_{1}, (1-t)z).
\end{split}
\]
We arrive at the expression~\eqref{eq:bg-near-boundary} upon noting that $X\mu\bv_{1} = \pi_{1}(q)$.

For part (b), by assumption we have $q \in \cV \cap B_{(4L_{\theta})^{-1}}(Y^{-})$. Thus, again writing $(X, \ep) = \widetilde{\bF}^{-1}(q)$, we infer from~\eqref{eq:tilde-bF-image-of-central-nbhd} that $|X| < \frac{1}{4}$, in which case $\xi(X, t) = (1-t)X$. This together with the scaling properties~\eqref{eq:bF-scaling-property} and~\eqref{eq:bP-scaling-property} implies that
\[
\widetilde{\bF}(\xi(X, t),(1-t)\ep) = \widetilde{\bF}((1-t)X, (1-t)\ep) = \bP((1-t)\bF(X, \ep)) = (1-t) \widetilde{\bF}(X, \ep).
\]
Substituting this into~\eqref{eq:bg-definition} gives~\eqref{eq:bg-near-center}.

Coming to part (c), for convenience we denote the map~\eqref{eq:tilde-bF-as-parametrization} by $\widehat{\bF}$. Given $p \in \mathring{\Sigma}$, upon expressing $p$ as $\widehat{\bF}(x, \ep, z)$, where $(\ep, z) \in \partial \cC$ and $\ep > 0$, we see from Lemma~\ref{lemm:U-Sigma-basics}(a) and~\eqref{eq:bg-definition} that 
\[
(d\widehat{\bF})_{(x, \ep, z)}(\{0\} \times T_{(\ep, z)}\partial\cC) \subset T_{p} \mathring{\Sigma},
\]
and that, for all $v \in T_{(\ep, z)}\partial\cC$ and $t \in [0, 1)$,
\begin{equation}\label{eq:bg-differential}
(d\bg_{t})_{p}\big( (d\widehat{\bF})_{(x, \ep, z)}(0, v)  \big)  = (1-t) \cdot (d\widehat{\bF})_{(\xi(x, t), (1-t)\ep, (1-t)z)}(0, v).
\end{equation}
Now, by counting dimensions, any given $2$-plane $E \subset T_{p}\mathring{\Sigma}$ must intersect $(d\widehat{\bF})_{(x, \ep, z)}(\{0\} \times T_{(\ep, z)}\partial \cC)$ nontrivially, and hence possess an orthonormal basis $u_1, u_2$ where 
\[
u_{1} = \frac{(d\widehat{\bF})_{(x, \ep, z)}(0,v)}{|(d\widehat{\bF})_{(x, \ep, z)}(0,v)|} \quad\text{for some non-zero }v \in T_{(\ep, z)}\partial \cC.
\]
With $t \in [0, 1)$, we infer from~\eqref{eq:bg-differential} and~\eqref{eq:tilde-bF-biLip-bounds} that
\[
|(d\bg_{t})_{p}(u_1)| = (1-t)\cdot \frac{|(d\widehat{\bF})_{(\xi(x, t), (1-t)\ep,(1-t)z)}(0,v)|}{|(d\widehat{\bF})_{(x, \ep, z)}(0,v)|} \leq C_{\theta}(1-t),
\]
while from~\eqref{eq:bg-Lipschitz} we get $|(d\bg_{t})_{p}(u_2)| \leq C_{\theta}$. This gives the asserted estimate.
\end{proof}
\begin{rmk}\label{rmk:bg-behavior-K2}
Following the proof of part (a), but using instead Remark~\ref{rmk:tilde-bF-image-K2},  we see that
\[
\bg(t, q, z) = ((1-t)q + t\varphi(|\pi_{2}(q)|\sin\theta)\pi_{2}(q),\ (1-t)z),
\]
for all $t \in [0, 1]$ and $(q, z) \in \Sigma \cap (K_{2} \cap \RR^{m-1})$.
\end{rmk}

\begin{prop}\label{prop:U-end-product}
There exists a universal constant $a_0 > 1$ such that for all $a > a_0$, the following hold.
\vskip 1mm
\begin{enumerate}
\item[(a)] We have 
\[
U \cap ( (a\bv_1 + \Span\{\bw_1\}) \times \RR^{m-1} ) = \{(a\bv_{1}+ (\ep\csc\theta)\bw_{1}, z)\ |\ (\ep, z) \in \cC,\ \ep \geq 0\}.
\]
Moreover this continues to hold with $U$ replaced by $\Sigma$, and with $\cC$ replaced by $\partial\cC$.
\vskip 1mm
\item[(b)] Given $(\ep, z) \in \partial\cC$ with $\ep\geq 0$, we have
\[
\bg(t, a\bv_1 + (\ep\csc\theta)\bw_{1}, z ) = (a\bv_{1} + (1-t)(\ep\csc\theta)\bw_{1}, (1-t)z),
\]
for all $t \in [0, 1]$.
\end{enumerate}
\end{prop}
\begin{proof}
Letting $L_{\theta} > 1$ and $\delta_0 \in (0, 1)$ be, respectively, the constants occurring in~\eqref{eq:tilde-bF-biLip-bounds} and~\eqref{eq:delta-0-choice}, we choose
\[
a_0 = \frac{10 L_{\theta}}{\delta_0 \sin\theta}.
\]
Observe by~\eqref{eq:cC-inclusions} that, for all $a > a_0$ and $(\ep, z) \in \cC$ with $\ep \geq 0$, we have $0 \leq \ep < \delta_0 \cdot a\sin\theta$, and hence~\eqref{eq:tilde-bF-expression} applies to give
\[
a\bv_{1} + (\ep\csc\theta)\bw_{1} = \widetilde{\bF}(a\sin\theta, \ep),
\]
which together with Lemma~\ref{lemm:U-Sigma-basics}(b) implies the inclusion ``$\supset$'' in part (a). For the reverse inclusion, take any $p \in U \cap ( (a\bv_1 + \Span\{\bw_1\}) \times \RR^{m-1} )$. By Lemma~\ref{lemm:U-Sigma-basics}(b) we can write it as
\begin{equation}\label{eq:U-end-p-expression}
p = (\widetilde{\bF}(X, \ep), z) =: (q, z),
\end{equation}
where $X \in \RR$, $\ep \geq 0$, and $(\ep, z) \in \cC$. Since $q \in a\bv_{1} + \Span\{\bw_1\}$, there is some $b \in \RR$ such that
\begin{equation}\label{eq:U-end-q-expression}
q = a\bv_{1} + b \bw_{1}.
\end{equation}
Using~\eqref{eq:cC-inclusions} and the second inclusion in~\eqref{eq:tilde-bF-image-of-dist-nbhd}, we infer from~\eqref{eq:U-end-p-expression} that 
\[
q \in \cV \cap B_{6L_{\theta}}(\partial\cV).
\]
Thus, in the expression~\eqref{eq:U-end-q-expression}, we must have $b \geq 0$, and that 
\begin{equation}\label{eq:U-end-consequence-of-dist}
a\cos\theta + b\sin\theta \leq  |a\sin\theta - b\cos\theta|\cdot \cot\theta + \frac{6L_{\theta}}{\sin\theta}.
\end{equation}
Suppose by contradiction that $a\sin\theta - b\cos\theta < 0$. Then the above becomes
\[
a\cos\theta + b\sin\theta \leq \frac{b\cos^2\theta}{\sin\theta} - a\cos\theta + \frac{6L_{\theta}}{\sin\theta}.
\]
Recalling that $b \geq 0$ and $\theta = \frac{\pi}{3}$, we deduce that $a \leq \frac{6L_{\theta}}{\sin\theta}$, contradicting our choice of $a_0$. Consequently $a\sin\theta - b\cos\theta \geq 0$, and we obtain from~\eqref{eq:U-end-consequence-of-dist} that
\[
b \leq 6L_{\theta} < \delta_0 a,
\]
which allows us to use~\eqref{eq:tilde-bF-expression} in~\eqref{eq:U-end-q-expression} to get $q = \widetilde{\bF}(a\sin\theta, b\sin\theta)$. Comparing this with~\eqref{eq:U-end-p-expression} gives $b = \ep \csc\theta$, so that 
\[
p = (q, z) = (a\bv_{1} + (\ep\csc\theta)\bw_{1}, z),
\]
where $(\ep, z) \in \cC$ and $\ep \geq 0$. This gives the inclusion ``$\subset$'' in part (a). The stated variant involving $\Sigma$ and $\partial\cC$ can be established by a similar argument, and we omit the details.

Given $(\ep, z) \in \partial\cC$ with $\ep \geq 0$, the second conclusion of part (a) implies that $(a\bv_{1} + (\ep\csc\theta)\bw_{1}, z) \in \Sigma$, so the statement of (b) makes sense. Next, since $\ep < \delta_0\cdot a\sin\theta$, as noted in the proof of (a), we get from~\eqref{eq:tilde-bF-expression} that
\[
(a\bv_{1} + (\ep\csc\theta) \bw_{1}, z) = (\widetilde{\bF}(a\sin\theta, \ep), z).
\]
Noting also that $a\sin\theta > 1$, and thus $\xi(a\sin\theta, t) = a\sin\theta > \frac{(1-t)\ep}{\delta_0}$, we have
\[
\begin{split}
\bg(t, a\bv_{1} + (\ep\csc\theta)\bw_{1}, z) = \ &(\widetilde{\bF}(a\sin\theta, (1-t)\ep), (1-t)z)\\
=\ & (a\bv_{1} + (1-t)(\ep\csc\theta)\bw_{1}, (1-t)z),
\end{split}
\]
and the proof is complete.
\end{proof}
\begin{rmk}\label{rmk:U-end-E2}
Using Remark~\ref{rmk:tilde-bF-image-K2}, and following the above proof, we have for all $a > a_0$ that
\[
U \cap ( (a\bv_{2} + \Span\{\bw_{2}\}) \times \RR^{m-1}) = \{(a\bv_2 - (\ep\csc\theta)\bw_2, z)\ |\ \ep \geq 0,\ (\ep, z) \in \cC\},
\]
and a similar equality holds with $U$ and $\cC$ replaced by $\Sigma$ and $\partial\cC$, respectively. Moreover, given $(\ep, z) \in \partial\cC$ with $\ep \geq 0$, there holds
\[
\bg(t, a\bv_{2} - (\ep\csc\theta)\bw_{2}, z) = (a\bv_{2} - (1-t)(\ep\csc\theta)\bw_{2}, (1-t)z),
\]
for all $t \in [0, 1]$.
\end{rmk}
To continue, take any $r \in (0, 1)$ which satisfies $r < \frac{\delta_0}{8\sin\theta}$, and define
\[
\rho_{1} = L_{\theta}^{-1}r, \quad \rho_{2} = (8L_{\theta})^{-1}.
\]
The choice of $r$ guarantees that
\begin{equation}\label{eq:outside-central-nbhd-in-cone}
B_{\rho_{1}}(Y) \setminus B_{\rho_{2}}(Y^{-}) \subset \mathring{K}_{1} \cup \mathring{K}_{2} \cup \mathring{K}_{3},
\end{equation}
where $\mathring{K}_{\lambda}$ stands for the interior of $K_{\lambda}$. Moreover, given $\lambda \in \{1, 2, 3\}$, we have
\begin{equation}\label{eq:outside-central-nbhd-expression}
B_{\rho_{1}}(Y) \cap ( \mathring{K}_{\lambda} \setminus B_{\rho_{2}}(Y^{-})) \subset \{a \bv_{\lambda} + b \bw_{\lambda}\ |\ a > \frac{\rho_1}{\delta_0},\ |b| < \rho_1\}.
\end{equation}
We then repeat the constructions up to this point with $\cV = \cV_{1}$ replaced by $\cV_{2}$ and $\cV_{3}$, and mark by subscripts the dependence of the resulting objects on the choice of $\lambda \in \{1, 2, 3\}$. Thus, for instance, we write 
\[
\bg_{1}:[0, 1] \times \Sigma_{1} \to U_{1}
\]
for the map defined by~\eqref{eq:bg-definition}. With this notation, we let
\begin{equation}\label{eq:Omega-for-retraction}
\Omega := (B_{\rho_{1}}(Y) \times \cC') \cup \mathring{U}_{1} \cup \mathring{U}_{2} \cup \mathring{U}_{3}, 
\end{equation}
and similarly define
\begin{equation}\label{eq:S-for-retraction}
S := (B_{\rho_{1}}(Y) \times \partial\cC') \cup \mathring{\Sigma}_{1} \cup \mathring{\Sigma}_{2} \cup \mathring{\Sigma}_{3}.
\end{equation}

Before proving, in Lemma~\ref{lemm:Omega-S-basics} below, that $\Omega$ is a smooth domain in $\RR^{m + 1}$ with boundary equal to $S$, we make another digression, this time to state a standard fact about patching together submanifolds with boundary.
\begin{lemm}\label{lemm:submanifold-construction}
Let $X$ be a subset of a smooth manifold $M^{n}$. Suppose there exist $d\in \NN$ and an open covering $\{\cO_{i}\}_{i \in \Lambda}$ of $M$, with $\Lambda$ an arbitrary index set, such that for all $i \in \Lambda$, the intersection $X_{i}: = X \cap \cO_{i}$ is an embedded $d$-dimensional submanifold of $\cO_{i}$ with boundary. Then $X$ is an embedded $d$-submanifold of $M$ with boundary. Moreover, denoting by $\mathring{X}_{i}$ and $\partial X_{i}$, respectively, the interior and boundary of $X_{i}$ as a submanifold of $\cO_{i}$, and defining $\mathring{X}$ and $\partial X$ similarly with respect to $M$, we have 
\[
\mathring{X} \cap \cO_{i} = \mathring{X}_{i},\quad \partial X \cap \cO_{i} = \partial X_{i}, \quad\text{for all }i\in \Lambda.
\]
\end{lemm}
\begin{proof}
We give a proof since the result is also used in subsequent sections. Given any $p \in X$, there is some $i \in \Lambda$ such that $p \in X_{i}$, in which case the first part of~\cite[Theorem 5.51]{Lee2013} gives a chart 
\[
\varphi:\widetilde{V} \subset \RR^{n} \to V \subset \cO_{i}
\]
such that $p \in V$, and that $\varphi^{-1}(V \cap X_{i})$ is either a $d$-slice or a $d$-dimensional half-slice of $\widetilde{V}$ in the terminology of~\cite[pages 101 and 122]{Lee2013}. Noting that $V \cap X_{i} = V \cap X$, and that $\varphi$ is also a chart for $M$ since $\cO_{i}$ is open, we may invoke~\cite[Theorem 5.51]{Lee2013} again to conclude that $X$ is a smooth $d$-submanifold of $M$ with boundary. Moreover, the previous argument actually implies
\[
\mathring{X}_{i} \subset \mathring{X} \cap \cO_{i}, \quad \partial X_{i} \subset \partial X \cap \cO_{i},\quad\text{for all }i \in \Lambda.
\]
Since $\mathring{X}_{i} \sqcup \partial X_{i} = X_{i} = (\mathring{X} \cap \cO_{i}) \sqcup (\partial X \cap \cO_{i})$, with both ends being disjoint unions, neither of the inclusions above can be proper. This finishes the proof.\\
\end{proof}
Returning to the main line of argument, we next establish the smoothness of $\Omega$ as a domain in $\RR^{m + 1}$, among other important properties.
\begin{lemm}\label{lemm:Omega-S-basics}
$\Omega$ is a smooth domain in $\RR^{m + 1}$ with boundary given by $S$. Moreover, the following hold.
\vskip 1mm
\begin{enumerate}
\item[(a)] We have
\[
\Omega = (Y \times \cC') \sqcup \mathring{U}_{1} \sqcup \mathring{U}_{2} \sqcup \mathring{U}_{3} = U_1 \cup U_2 \cup U_3,
\]
and that
\[
S = (Y \times \partial\cC') \sqcup \mathring{\Sigma}_{1} \sqcup \mathring{\Sigma}_{2} \sqcup \mathring{\Sigma}_{3} = \Sigma_1 \cup \Sigma_2 \cup \Sigma_3.
\]
\vskip 1mm
\item[(b)] Let $a_0$ be the threshold in Proposition~\ref{prop:U-end-product}. Then, for all $a > a_0$ and $\lambda \in \{1, 2, 3\}$, we have
\[
\Omega\cap ( (a\bv_{\lambda} + \Span\{\bw_{\lambda}\}) \times \RR^{m-1}) = \{(a\bv_{\lambda} + b\bw_{\lambda}, z)\ |\ (b\sin\theta, z) \in \cC\},
\]
and that
\[
S \cap ( (a\bv_{\lambda} + \Span\{\bw_{\lambda}\}) \times \RR^{m-1}) = \{(a\bv_{\lambda} + b\bw_{\lambda}, z)\ |\ (b\sin\theta, z) \in \partial\cC\}.
\]
\end{enumerate}
\end{lemm}
\begin{proof}
From Lemma~\ref{lemm:U-Sigma-basics}(c), we deduce that, for $\lambda = 1, 2, 3$,
\begin{equation}\label{eq:U-lambda-intersect-center}
\mathring{U}_{\lambda} \cap (B_{\rho_{1}}(Y) \times \RR^{m-1}) = (B_{\rho_{1}}(Y) \cap \mathring{\cV}_{\lambda}) \times \cC',
\end{equation}
which implies that
\begin{equation}\label{eq:Omega-smooth-center}
\Omega \cap (B_{\rho_{1}}(Y) \times \RR^{m-1})  = B_{\rho_{1}}(Y) \times \cC',
\end{equation}
and that
\begin{equation}\label{eq:Omega-smooth-branches}
\Omega \cap (\mathring{\cV}_{\lambda} \times \RR^{m-1}) = \mathring{U}_{\lambda}.
\end{equation}
Notice that the right-hand side in~\eqref{eq:Omega-smooth-center} and~\eqref{eq:Omega-smooth-branches} are smooth domains in $B_{\rho_1}(Y) \times \RR^{m-1}$ and $\mathring{\cV}_{\lambda} \times \RR^{m-1}$, respectively, with (relative) boundary given by $B_{\rho_1}(Y) \times \partial\cC'$ and $\mathring{\Sigma}_{\lambda}$, the latter characterization being a consequence of Lemma~\ref{lemm:U-Sigma-basics}. Since $B_{\rho_{1}}(Y) \times \RR^{m-1}$ along with $\mathring{\cV}_{\lambda} \times \RR^{m-1}$ for $\lambda = 1, 2, 3$ form an open covering of $\RR^{m+1}$, we get from Lemma~\ref{lemm:submanifold-construction} that $\Omega$ is a smooth domain in $\RR^{m + 1}$, with its boundary satisfying
\[
\partial \Omega \cap (B_{\rho_{1}}(Y) \times \RR^{m-1})  = B_{\rho_{1}}(Y) \times \partial\cC',
\]
and, for $\lambda = 1, 2, 3$, 
\[
\partial\Omega \cap (\mathring{\cV}_{\lambda} \times \RR^{m-1}) = \mathring{\Sigma}_{\lambda}.
\]
In particular, we get $\partial\Omega = S$. This proves the first assertion of the lemma.

For part (a), the first stated expression for $\Omega$ follows from~\eqref{eq:U-lambda-intersect-center}, and the fact that $B_{\rho_1}(Y) \setminus Y \subset \cup_{\lambda = 1}^{3}\mathring{\cV}_{\lambda}$. For the second, we note in addition that 
\[
U_{\lambda} \setminus \mathring{U}_{\lambda} = U_{\lambda} \cap (\partial \cV_{\lambda} \times \RR^{m-1})  = \partial \cV_{\lambda} \times \cC',
\]
where the last step uses Lemma~\ref{lemm:U-Sigma-basics}(c). Using instead that $\Sigma_{\lambda}\setminus \mathring{\Sigma}_{\lambda} = \partial\cV_{\lambda} \times \partial\cC'$, and that
\[
\mathring{\Sigma}_{\lambda} \cap (B_{\rho_{1}}(Y) \times \RR^{m-1}) = (B_{\rho_{1}}(Y) \cap \mathring{\cV}_{\lambda}) \times \partial\cC',
\]
which again come from Lemma~\ref{lemm:U-Sigma-basics}(c), we get the stated expressions for $S$.

For part (b), to avoid repetition, we only address the first equality, and only in the case $\lambda = 2$. To start, note that $(a\bv_{2} + \Span\{\bw_2\}) \cap \cV_{3} = \emptyset$, and hence by part (a) we have
\begin{equation}\label{eq:Omega-end}
\Omega \cap ( (a\bv_{2} + \Span\{\bw_{2}\})  \times \RR^{m-1}) = (U_1 \cup U_2) \cap ( (a\bv_{2} + \Span\{\bw_{2}\}) \times \RR^{m-1}).
\end{equation}
Next, Proposition~\ref{prop:U-end-product}(a) implies that 
\[
U_2 \cap ( (a\bv_{2} + \Span\{\bw_2\})  \times \RR^{m-1}) = \{ (a\bv_2 + b\bw_2, z)\ |\ b \geq 0,\ (b\sin\theta, z)\in \cC \}.
\]
On the other hand, by Remark~\ref{rmk:U-end-E2}, along with the symmetry property of $\cC$ coming from~\eqref{eq:N-even-in-each}, we get
\[
\begin{split}
U_1 \cap ( ( a\bv_{2} + \Span\{\bw_{2}\} ) \times \RR^{m-1} ) =\ & \{(a\bv_{2} +b \bw_{2}, z)\ |\ b \leq 0,\ (-b\sin\theta, z) \in \cC\} \\
=\ & \{ (a\bv_2 + b\bw_2, z)\ |\ b\leq 0,\ (b\sin\theta, z)\in \cC \}.
\end{split}
\]
Combining the previous two observations with~\eqref{eq:Omega-end}, we arrive at the desired conclusion.
\end{proof}

We next construct a map from $[0, 1] \times B_{\rho_1}(Y) \times \partial\cC'$ into $B_{\rho_1}(Y) \times \cC'$ which we will patch with $\bg_{1}$, $\bg_{2}$, and $\bg_{3}$ to obtain the desired $\bh$ mentioned at the start of the section. As an ingredient, we first define 
\[
\chi:[0, 1] \times B_{\rho_{1}}(Y) \to B_{\rho_{1}}(Y) 
\]
by letting
\begin{equation}\label{eq:bg-c-definition}
\chi(t, q) = \left\{
\begin{array}{ll}
(1-t)q, & \text{if } q \in B_{\rho_{1}}(Y) \cap B_{2\rho_{2}}(Y^{-}),\\
(1-t)q + t\varphi(|\pi_{\lambda}(q)|\sin\theta)\pi_{\lambda}(q), & \text{if }q \in B_{\rho_{1}}(Y) \cap  \mathring{K}_{\lambda}\setminus \overline{B_{\rho_{2}}(Y^{-})}\\
& \text{for some }\lambda \in \{1, 2, 3\},
\end{array}
\right.
\end{equation}
where recall that $\pi_{\lambda}$ denotes orthogonal projection onto $\Span\{\bv_{\lambda}\}$. In view of~\eqref{eq:outside-central-nbhd-in-cone}, any $q \in B_{\rho_1}(Y)$ falls into at least one of the cases considered above. 
\begin{prop}\label{prop:bg-c-properties}
The map $\chi$ is well-defined, and indeed takes values in $B_{\rho_1}(Y)$. Moreover, letting 
\[
\bg_{c}(t, q, z) := (\chi(t, q), (1-t)z), \quad\text{for }(t, q, z) \in [0, 1] \times B_{\rho_1}(Y) \times \partial\cC'
\]
yields a map into $B_{\rho_1}(Y) \times \cC'$ that is Lipschitz on compact subsets of $[0, 1] \times B_{\rho_1}(Y) \times \partial\cC'$, and has the following additional properties.
\vskip 1mm
\begin{enumerate}
\item[(a)] For each $\lambda \in \{1, 2, 3\}$, there holds
\begin{equation}\label{eq:bg-c-bg-lambda-agree}
\bg_{c}(t, q, z) = \bg_{\lambda}(t,q, z),
\end{equation}
whenever $t \in [0, 1]$, and $(q, z) \in (B_{\rho_1}(Y) \cap \cV_{\lambda}) \times \partial\cC' = (B_{\rho_1}(Y) \times \partial\cC') \cap \Sigma_{\lambda}$.
\vskip 1mm
\item[(b)] There exists a universal constant $C$ such that given $t \in [0, 1)$ and $(q, z) \in B_{\rho_1}(Y) \times \partial \cC'$, along with a $2$-plane $E \subset T_{q}\RR^2 \times T_{z}\partial\cC'$, upon abbreviating $\bg_{c}(t, \cdot)$ as $\bg_{c, t}$, there holds
\begin{equation}\label{eq:bg-c-area-estimate}
\big|(d\bg_{c, t})_{(q,z)}(u_1) \wedge (d\bg_{c, t})_{(q, z)}(u_2)\big| \leq C(1-t),
\end{equation}
where $u_1, u_2$ is any orthonormal basis of $E$, the choice being irrelevant for the norm on the left-hand side.
\end{enumerate}
\end{prop}
\begin{proof}
We first prove that $\chi$ is well-defined. Given $\lambda \in \{1, 2, 3\}$ and some $q \in B_{\rho_{1}}(Y) \cap \mathring{K}_{\lambda} \cap ( B_{2\rho_{2}}(Y^{-}) \setminus \overline{B_{\rho_{2}}(Y^{-})})$, by Remark~\ref{rmk:tilde-bF-image-cone} we have $|\pi_{\lambda}(q)|\sin\theta < \frac{1}{4}$, which together with~\eqref{eq:central-cutoff} implies that 
\[
(1-t)q + t\varphi(|\pi_{\lambda}(q)|\sin\theta)\pi_{\lambda}(q)= (1-t)q.
\]
Since, by what we observed right before Lemma~\ref{lemm:tilde-bF-boundary-behavior}, the sets $B_{\rho_1}(Y) \cap \mathring{K}_{\lambda} \setminus \overline{B_{\rho_{2}}(Y^{-})}$ for $\lambda = 1, 2, 3$ are pairwise disjoint, we have shown that the definition~\eqref{eq:bg-c-definition} is consistent on all possible overlap regions. Next, we claim that 
\begin{equation}\label{eq:chi-image-ball}
\chi([0, 1] \times (B_{\rho_1}(Y) \cap B_{R})) \subset B_{\rho_{1}}(Y) \cap B_{R},
\end{equation}
for all $R > 0$. Indeed, since $0 \leq \varphi \leq 1$, we see from~\eqref{eq:bg-c-definition} that 
\begin{equation}\label{eq:chi-preserive-BR}
|\chi(t, q)| \leq |q|, \quad\text{for all } t \in [0, 1],\ q \in B_{\rho_1}(Y).
\end{equation}
As for the distance from $\chi(t, q)$ to $Y$, in the case $q \in B_{\rho_1}(Y) \cap B_{2\rho_2}(Y^{-})$, since $Y$ is preserved by scaling, we have
\[
\dist(\chi(t, q), Y) < (1-t)\rho_1 \leq \rho_1,
\]
whereas in the case $q \in B_{\rho_{1}}(Y) \cap  \mathring{K}_{\lambda}\setminus \overline{B_{\rho_{2}}(Y^{-})}$ for some $\lambda \in \{1, 2, 3\}$, we write 
\begin{equation}\label{eq:chi-expression-decomposed}
\chi(t, q) = (1-t)(q - \pi_{\lambda}(q)) + [(1-t)+ t\varphi(|\pi_{\lambda}(q)|\sin\theta)] \pi_{\lambda}(q),
\end{equation}
and notice from~\eqref{eq:outside-central-nbhd-expression} that $\pi_{\lambda}(q) \in Y$ while $|q - \pi_{\lambda}(q)| < \rho_1$. Consequently we get
\[
\dist(\chi(t, q), Y) \leq (1-t) |q - \pi_{\lambda}(q)| < \rho_1.
\]
Combining this with the distance estimate in the previous case, as well as~\eqref{eq:chi-preserive-BR}, we get~\eqref{eq:chi-image-ball} as claimed. In particular, $\chi$ does indeed send $[0, 1] \times B_{\rho_1}(Y)$ into $B_{\rho_1}(Y)$. 

We also need a Lipschitz estimate on $\chi$. Using the first line of~\eqref{eq:bg-c-definition} and the boundedness of $B_{\rho_1}(Y) \cap B_{2\rho_2}(Y^{-})$, we control the Lipschitz semi-norm of $\chi$ on $[0, 1] \times (B_{\rho_1}(Y) \cap B_{2\rho_2}(Y^{-}))$ by some universal constant. To see that a similar estimate holds for $\chi$ restricted to each $[0, 1] \times ((B_{\rho_{1}}(Y) \setminus \overline{B_{\rho_{2}}(Y^{-})}) \cap \mathring{K}_{\lambda})$, we recall the definition~\eqref{eq:xi-for-retraction}, and use the fact that $\pi_{\lambda}(q) = (q\cdot \bv_{\lambda})\bv_{\lambda}$, with $q \cdot \bv_{\lambda}$ being positive by~\eqref{eq:outside-central-nbhd-expression}, to rewrite~\eqref{eq:chi-expression-decomposed} as
\begin{equation}\label{eq:chi-expression-decomposed-2}
\chi(t, q) = (1-t)(q - \pi_{\lambda}(q)) + \xi(q\cdot \bv_{\lambda}\sin\theta, t)\cdot \csc\theta\cdot\bv_{\lambda},
\end{equation}
and use again the inequality $|q - \pi_{\lambda}(q)| < \rho_1$, along with the bounds~\eqref{eq:xi-partial-derivatives} on the derivatives of $\xi$. Recalling~\eqref{eq:chi-image-ball}, and noting that any two points in $[0, 1] \times B_{\rho_1}(Y)$ which are sufficiently close together must fall into the same case among those considered in~\eqref{eq:bg-c-definition}, we obtain, for each $R > 0$, 
\[
|\chi(t_1, q_1) - \chi(t_2, q_2)| \leq C_{R}(|t_1 - t_2| + |q_1 - q_2|),
\]
whenever $(t_1, q_2), (t_2, q_2) \in [0, 1] \times (B_{\rho_1}(Y) \cap B_{R})$.

We now turn to the map $\bg_{c}$. From~\eqref{eq:chi-image-ball} and Remark~\ref{rmk:cC'-retraction}, we see that $\bg_{c}$ takes values in $B_{\rho_1}(Y) \times \cC'$ as asserted. Moreover, combining Remark~\ref{rmk:cC'-retraction} with the estimate at the end of the previous paragraph shows that $\bg_{c}$ is Lipschitz on $[0, 1] \times (B_{\rho_1}(Y) \cap B_R) \times \partial\cC'$ for each $R > 0$. For part (a), the equality of sets at the end of the statement follows from Lemma~\ref{lemm:U-Sigma-basics}(c). To prove the relation~\eqref{eq:bg-c-bg-lambda-agree}, take $t \in [0, 1]$, $z \in \partial\cC'$, and $q \in B_{\rho_1}(Y) \cap \cV_{\lambda}$. In the case $q \in B_{2\rho_2}(Y^{-})$, we have $(q, z) \in \Sigma_{\lambda} \cap (B_{2\rho_2}(Y^{-}) \times \RR^{m-1})$, and thus Proposition~\ref{prop:bg-behavior}(b) gives
\[
\bg_{\lambda}(t, q, z) = ((1-t)q, (1-t)z) = \bg_{c}(t, q, z).
\]
On the other hand, if $q \in \mathring{K}_{\mu} \setminus \overline{B_{\rho_{2}}(Y^{-})} $ for some $\mu \in \{1, 2, 3\}$, then necessarily $\mu = \lambda$ or $\lambda + 1$, and we get~\eqref{eq:bg-c-bg-lambda-agree} from Proposition~\ref{prop:bg-behavior}(a) and Remark~\ref{rmk:bg-behavior-K2}.

For part (b), the estimate~\eqref{eq:bg-c-area-estimate} obviously holds when $q \in B_{\rho_1}(Y) \cap B_{2\rho_2}(Y^{-})$, since in this case $(d\bg_{c, t})_{(q, z)}$ coincides with multiplication by $(1-t)$ on all of $T_{q}\RR^{2} \times T_{z}\partial\cC'$. On the other hand, when $q \in B_{\rho_1}(Y) \cap \mathring{K}_{\lambda} \setminus \overline{B_{\rho_2}(Y^{-})}$ for some $\lambda \in \{1, 2, 3\}$, we express it as 
\[
q  = a_0 \bv_{\lambda} + b_0\bw_{\lambda},
\]
so that $|b_0| < \rho_1 < \delta_0 a_0$ by~\eqref{eq:outside-central-nbhd-expression}. Then for all $(a, b)$ sufficiently close to $(a_0, b_0)$, it follows from~\eqref{eq:chi-expression-decomposed-2} that
\begin{equation}\label{eq:bg-c-defi-alt}
\chi(t, a\bv_{\lambda} + b\bw_{\lambda}) =  \xi(a\sin\theta, t)\cdot \csc\theta\cdot \bv_{\lambda} + (1-t)b\bw_{\lambda},
\end{equation}
which implies that
\[
(d\bg_{c, t})_{(q, z)}(\bv_{\lambda}, 0) = \pa{\xi}{x}(|\pi_{\lambda}(q)|\sin\theta, t)\cdot (\bv_{\lambda}, 0),
\]
and that $(d\bg_{c, t})_{(q, z)}$ restricts to multiplication by $(1-t)$ on $\Span\{\bw_{\lambda}\} \times T_{z}\partial\cC'$. For dimension reasons, the latter must have nontrivial intersection with the given $2$-plane $E$. Thus, recalling also the bound~\eqref{eq:xi-partial-derivatives} on $\pa{\xi}{x}$, we conclude again that~\eqref{eq:bg-c-area-estimate} holds.
\end{proof}

\begin{prop}\label{prop:bh-properties}
Define $\bh: [0, 1] \times S \to \Omega$ by letting
\begin{equation}\label{eq:bh-definition}
\bh(t, p) = \left\{
\begin{array}{ll}
\bg_{c}(t, p), & \text{ if }p \in B_{\rho_{1}}(Y) \times \partial \cC',\\
\bg_{\lambda}(t, p), & \text{ if }p \in \mathring{\Sigma}_{\lambda}.
\end{array}
\right.
\end{equation}
Then $\bh$ is well-defined, maps $[0, 1] \times S$ onto $\Omega$, and is Lipschitz on compact subsets. Moreover, the following hold.
\vskip 1mm
\begin{enumerate}
\item[(a)] $\bh(0, \cdot) = \id_{S}$, and $\bh(\{1\} \times S) = Y \times \{\0^{m-1}\}$.
\vskip 1mm
\item[(b)] $\bh$ maps $[0, 1) \times S$ diffeomorphically onto $\Omega \setminus (Y \times \{\0^{m-1}\})$.
\vskip 1mm
\item[(c)] For all $(t, p) \in [0, 1) \times S$ and $2$-plane $E \subset T_{p}S$, letting $\bh_{t} = \bh(t, \cdot)$ we have 
\begin{equation}\label{eq:bh-area-estimate}
|(d\bh_{t})_{p}(u_1) \wedge (d\bh_{t})_{p}(u_2)| \leq C(1-t),
\end{equation}
where $C$ is a universal constant, and as before $u_1, u_2$ is any orthonormal basis for $E$.
\vskip 1mm
\item[(d)] Let $a_0$ be the threshold from Proposition~\ref{prop:U-end-product}. Then for all $a > a_0$, $\lambda \in \{1, 2, 3\}$, and $p' \in \Span\{\bw_{\lambda}\} \times \RR^{m-1}$ such that $p: = (a\bv_{\lambda}, \0^{m-1}) + p' \in S$, there holds
\[
\bh(t, p) =  (a\bv_{\lambda}, \0^{m-1}) + (1-t)p'.
\]
\end{enumerate}
\end{prop}
\begin{proof}
Since $\Sigma_{\lambda} \cap \mathring{\Sigma}_{\mu} = \emptyset$ whenever $\lambda \neq \mu$, we get from Proposition~\ref{prop:bg-c-properties}(a) that the definition~\eqref{eq:bh-definition} is consistent, and that in fact
\begin{equation}\label{eq:bh-bg-lambda-agree}
\bh(t, p) = \bg_{\lambda}(t, p) \quad\text{for all }(t, p) \in [0, 1] \times \Sigma_{\lambda}.
\end{equation}
Recalling that $U_{\lambda} = \bg_{\lambda}([0, 1] \times \Sigma_{\lambda})$ (see below~\eqref{eq:bg-definition}), we infer from~\eqref{eq:bh-bg-lambda-agree} and the expressions in Lemma~\ref{lemm:Omega-S-basics} that $\Omega = \bh([0, 1] \times S)$.

Towards proving the assertion about the Lipschitz continuity of $\bh$, note that if $p = (q, z) \in S \setminus (B_{\rho_1}(Y) \times \partial\cC')$, then there is a unique $\lambda \in \{1, 2, 3\}$ such that 
\[
p \in \mathring{\Sigma}_{\lambda} \setminus  (B_{\rho_1}(Y) \times \partial\cC') = \mathring{\Sigma}_{\lambda} \setminus  (B_{\rho_1}(Y) \times \RR^{m-1}),
\]
where the equality can be deduced from Lemma~\ref{lemm:U-Sigma-basics}(c). In particular $q \in \mathring{\cV}_{\lambda} \setminus B_{\rho_1}(Y)$, and consequently we have
\[
p' \in \mathring{\Sigma}_{\lambda}, \quad\text{for all }p' \in S \text{ satisfying }|p' - p| < \rho_1.
\]
Indeed, if $p' = (q', z') \in S \setminus \mathring{\Sigma}_{\lambda}$, then $q' \in \cV_{\lambda+1} \cup \cV_{\lambda+2}$ by Lemma~\ref{lemm:Omega-S-basics}(a), and thus 
\[
|p - p'| \geq |q - q'| \geq \rho_1.
\]

Now, thanks respectively to Proposition~\ref{prop:bg-c-properties} and~\eqref{eq:bg-Lipschitz}, we see that $\bh$ restricts to be Lipschitz on $[0, 1] \times (B_{\rho_1}(Y) \cap B_{R}) \times \partial \cC'$ for all $R > 0$, and on each $[0, 1] \times \Sigma_{\lambda}$. Combining this with the discussion at the end of the previous paragraph, and using also the $L^{\infty}$-bounds coming from~\eqref{eq:bg-scale-factor},~\eqref{eq:chi-image-ball}, and~\eqref{eq:cC-inclusions}, we conclude that $\bh$ is a Lipschitz map on $[0, 1] \times (S \cap ( B_R \times \RR^{m-1}) )$ for each $R > 0$. 

For part (a), inspecting the definitions~\eqref{eq:bg-definition},~\eqref{eq:xi-for-retraction}, and~\eqref{eq:bg-c-definition}, we see that $\bh(0, p) = p$ for all $p \in S$. The second assertion follows from~\eqref{eq:bh-bg-lambda-agree} and~\eqref{eq:bg-for-t=1}.

For part (b), we observe first that, as observed below~\eqref{eq:Omega-smooth-branches} in the proof of Lemma~\ref{lemm:Omega-S-basics}, we have 
\[
B_{\rho_{1}}(Y) \times \partial\cC' = S\cap (B_{\rho_{1}}(Y) \times \RR^{m-1}) , \quad \mathring{\Sigma}_{\lambda} = S \cap (\mathring{\cV}_{\lambda} \times \RR^{m-1}),
\]
so each of the regions appearing in~\eqref{eq:bh-definition} is open relative to $S$. To continue, by~\eqref{eq:bh-bg-lambda-agree} and~\eqref{eq:bg-lambda-diffeo}, we see that $d\bh$ is invertible everywhere on $[0, 1) \times \mathring{\Sigma}_{\lambda}$. That the same is true on $[0, 1) \times B_{\rho_1}(Y) \times \partial\cC'$ can be seen with the help of~\eqref{eq:bg-c-definition}, the expression~\eqref{eq:bg-c-defi-alt}, and Remark~\ref{rmk:cC'-retraction}. Furthermore, using again~\eqref{eq:bg-lambda-diffeo}, along with~\eqref{eq:bg-boundary-injective} and~\eqref{eq:bg-image-boundary}, we see that $\bh$ is injective when restricted to each of 
\[
[0, 1) \times Y \times \partial\cC', \quad [0, 1) \times \mathring{\Sigma}_{1}, \quad [0, 1) \times \mathring{\Sigma}_{2}, \quad [0, 1) \times \mathring{\Sigma}_{3}
\]
and that, under $\bh$, these have images equal respectively to
\[
Y \times (\cC' \setminus \{\0^{m-1}\}),\quad \mathring{U}_{1},\quad \mathring{U}_{2},\quad \mathring{U}_{3}.
\]
It follows that $\bh$ is a bijection from $[0, 1) \times S$ onto $\Omega \setminus (Y \times \{\0^{m-1}\})$. Having already verified the invertibility of $d\bh$ everywhere on $[0, 1) \times S$, we are done with (b) upon recalling part (a) and using Lemma~\ref{lemm:boundary-IVT}.

For part (c), we emphasize once again that each region involved in~\eqref{eq:bh-definition} is open in $S$. Now, in the case $p = (q, z) \in B_{\rho_1}(Y) \times \partial\cC'$, we have $T_{p}S = T_{q}\RR^2 \times T_{z}\partial\cC'$ and that $(d\bh_{t})_{p} = (d\bg_{c, t})_{(q, z)}$, and the result follows from Proposition~\ref{prop:bg-c-properties}(b). In the case $p \in \mathring{\Sigma}_{\lambda}$ for some $\lambda \in \{1, 2, 3\}$, we have instead that $T_{p}S = T_{p}\mathring{\Sigma}_{\lambda}$ and $(d\bh_{t})_{p} = (d\bg_{\lambda, t})_{p}$, and we get~\eqref{eq:bh-area-estimate} from Proposition~\ref{prop:bg-behavior}(c).

Finally, let $p$ be as in the statement of (d), and assume without loss of generality that $\lambda = 2$. By Lemma~\ref{lemm:Omega-S-basics}(b), we get some $b \in \RR$ and $z \in \RR^{m-1}$ such that $(b\sin\theta, z) \in \partial\cC$, and that
\[
p = (a\bv_2 + b\bw_{2}, z).
\]
If $b \geq 0$, then with the help of Lemma~\ref{lemm:Omega-S-basics}(a), we have
\[
p \in S \cap (\cV_{2} \times \RR^{m-1}) = \Sigma_{2},
\]
in which case~\eqref{eq:bh-bg-lambda-agree} and Proposition~\ref{prop:U-end-product}(b) together implies
\[
\bh(t, p) = \bg_{2}(t, p) = (a\bv_{2} + (1-t)b\bw_{2}, (1-t)z),
\]
proving the asserted formula. If $b < 0$, then we get $p \in \Sigma_1$, so that $\bh(t, p) = \bg_{1}(t, p)$, and the conclusion follows instead from Remark~\ref{rmk:U-end-E2}. The proof is complete.
\end{proof}

\begin{rmk}\label{rmk:retraction-end-remark}
One may paraphrase Proposition~\ref{prop:bh-properties}(c) by saying that $\bh$ shrinks the length of $2$-vectors tangent to $S$ uniformly to $0$ as $t \to 1^{-}$, and this is how the statement involving $2$-vectors made at the beginning of the section is to be interpreted. Also, Lemma~\ref{lemm:Omega-S-basics} and Proposition~\ref{prop:bh-properties}(d) guarantees that the pair $(\Omega, \bh)$ can be smoothly patched together with other similar building blocks to yield the map $[0, 1] \times \partial H \to H$ mentioned at the start. The listed properties of the latter would then emerge from parts (a), (b), and (c) of Proposition~\ref{prop:bh-properties}, respectively.
\end{rmk}
\subsection{Patching within fixed stratum}\label{subsec:fixed-stratum}
In this section we prove Proposition~\ref{prop:fixed-stratum-patching}. Suppose $M^{n + 1}$, where $n \geq 3$, is a closed, connected, oriented Riemannian manifold isometrically embedded in some Euclidean space $\RR^{Q}$, and that it admits an $n$-section in the sense of Section~\ref{subsec:multisections}. Take a non-empty $I \subset \{1, \cdots, n\}$ and denote its length by $k$. Recall from Definition~\ref{defi:multisections} that, by definition,
\begin{equation}\label{eq:boundary-sector}
\mathring{M}_{I} = M_I \setminus \big( \cup_{i \not\in I}M_{i} \big),\quad \partial M_{I} = M_{I} \setminus \mathring{M_{I}},
\end{equation}
and that the former consists exactly of points $p$ for which $\{1 \leq i \leq n\ |\ p \in M_{i}\} = I$. Let $\Delta^{k-1, I}, V^{k-1, I}$ be the objects defined at the start of Appendix~\ref{appendix:patching-adapted}. As noted in Section~\ref{subsec:multisections}, it follows from Definition~\ref{defi:multisections}(2) that, for all $p \in \mathring{M}_{I}$, there exists a chart
\[
\varphi: \widetilde{U}\subset \RR^{n + 2 -k} \times {V}^{k-1, I} \to U \subset M
\]
such that $(\0^{n+2 - k}, \0^{k-1}) \in \widetilde{U}$ and $p \in U$, and that~\eqref{eq:modified-adapted-chart} holds, namely
\[
\left\{
\begin{array}{ll}
U \cap M_{i} = \varphi(\widetilde{U} \cap (\RR^{n + 2 - k} \times {V}^{k-1, I}_{i})),& \text{ if } i \in I,\\
\overline{U} \cap M_{i} = \emptyset, & \text{ if }i \not\in I.
\end{array}
\right.
\]
In this appendix we call such a $\varphi$ an adapted chart at $p$. Thus, as mentioned in Definition~\ref{defi:multisections}, from the property~\eqref{eq:modified-adapted-chart} and standard results about smooth manifolds (for example~\cite[Theorem 5.51]{Lee2013}), it can be checked that the set $\mathring{M}_{I} \cup \big( \cup_{i \not\in I} \mathring{M}_{I \cup \{i\}} \big)$ is an embedded $(n+2 - k)$-submanifold with boundary, and $\mathring{M}_{I}$ is the set of interior points. 

The next lemma is an immediate consequence of~\eqref{eq:modified-adapted-chart}.
\begin{lemm}\label{lemm:adapted-chart-differential}
Given adapted charts $\varphi_{a}: \widetilde{U}_{a} \to U_{a}$ at some $p \in \mathring{M}_{I}$, where $a = 1, 2$, let $\widetilde{V}_{a} = \varphi_{a}^{-1}(U_1 \cap U_2)$, and define $\tau : \widetilde{V}_{1} \to \widetilde{V}_{2}$ by 
\[
\tau = \varphi_{2}^{-1} \circ \varphi_{1}.
\]
Then for all $(z, {\bf 0}^{k-1}) \in \widetilde{V}_{1} \cap (\RR^{n + 2-k} \times \{\0^{k-1}\})$ we have the following.
\vskip 1mm
\begin{enumerate}
\item[(a)] $(d\tau)_{(z, {\0^{k-1}})}$ preserves $\RR^{n + 2 - k} \times \{\0^{k-1}\}$.
\vskip 1mm
\item[(b)] For each $i \in I$, there is some $\lambda > 0$ such that
\[
(d\tau)_{(z, {\bf 0}^{k-1})}({\bf 0}^{n + 2 - k}, \ba_{i; I})  - (\0^{n+2 - k}, \lambda \ba_{i; I}) \in \RR^{n + 2 - k} \times \{\0^{k-1}\}.
\]
\end{enumerate}
\end{lemm}
\begin{proof}
From~\eqref{eq:modified-adapted-chart} we have
\[
\tau(\widetilde{V}_{1}  \cap (\RR^{n + 2 - k} \times \{\0^{k-1}\})) \subset \RR^{n + 2 - k} \times \{\0^{k-1}\},
\]
which implies (a). For part (b), notice that, again by~\eqref{eq:modified-adapted-chart}, for all sufficiently small positive $t$, we have
\[
\tau(z, t\ba_{i; I}) \in \RR^{n + 2 - k} \times {V}^{k-1, I}_{I \setminus \{i\}} .
\]
Since $\tau(z, \0^{k-1}) \in \RR^{n + 2 - k} \times \{\0^{k-1}\}$, differentiating at $t =0$ shows that the desired conclusion holds for some $\lambda \geq 0$. That $\lambda$ cannot vanish follows from part (a) and the invertibility of $(d\tau)_{(z, \0^{k})}$. This proves part (b).
\end{proof}

\begin{coro}[\cite{BACGM2023}, page 4]
\label{coro:trivial-quotient-bundle}
The normal bundle $T^{\perp}\mathring{M}_{I}$ is trivial. More precisely, for any $i_0 \in I$ there exists a global frame $\{\xi_{i}\ |\ i\in I \setminus \{i_0\}\}$ of $T^{\perp}\mathring{M}_{I}$ which admits the following characterization: for all $p \in \mathring{M}_{I}$  and $i \in I \setminus \{i_0\}$, there holds 
\begin{equation}\label{eq:section-by-adapted-chart}
\xi_{i}(p) = \frac{(d\varphi)_{\varphi^{-1}(p)}({\bf 0}^{n + 2 - k}, \ba_{i;I})^{\perp}}{\| (d\varphi)_{\varphi^{-1}(p)}({\bf 0}^{n+ 2 - k}, \ba_{i;I})^{\perp} \|},
\end{equation}
where $\varphi:\widetilde{U} \to U$ is any choice of adapted chart whose image contains $p$.
\end{coro}
\begin{proof}
Given $p_0 \in \mathring{M}_{I}$, we first show that the right-hand side of~\eqref{eq:section-by-adapted-chart} is independent of the choice of adapted chart whose image contains $p_0$. Indeed, if $\psi: \widetilde{V} \to V$ is another such chart, then thanks to Lemma~\ref{lemm:adapted-chart-differential}(b), there exists $\lambda > 0$ such that 
\[
(d\varphi)_{\varphi^{-1}(p_0)}({\bf 0}^{n + 2 - k}, \ba_{i; I})  - \lambda \cdot (d\psi)_{\psi^{-1}(p_0)}({\bf 0}^{n + 2 - k}, \ba_{i;I}) \in T_{p_0}\mathring{M}_{I}.
\]
Taking the normal component on both sides shows that replacing $\varphi$ by $\psi$ makes no difference in~\eqref{eq:section-by-adapted-chart}, so in particular each $\xi_i$ is a well-defined smooth section of $T^{\perp}\mathring{M}_{I}$. It remains to show that $\{\xi_{i}(p_0)\}_{i \in I \setminus \{i_0\}}$ is a basis for $T_{p_0}^{\perp}\mathring{M}_{I}$ for all $p_0 \in \mathring{M}_{I}$. By the definition~\eqref{eq:section-by-adapted-chart} and standard linear algebra, this is equivalent to showing that the cosets
\[
(d\varphi)_{\varphi^{-1}(p_0)}({\bf 0}^{n + 2 - k}, \ba_{i;I}) + T_{p_0}\mathring{M}_{I}\quad \text{ for }i \in I \setminus \{i_0\},
\]
form a basis for $T_{p_0}M / T_{p_0}\mathring{M}_{I}$, which can be verified without much difficulty using the fact that $\{\ba_{i; I}\}_{i \in I \setminus \{i_0\}}$ is a basis for $V^{k-1, I}$, and that $(d\varphi)_{\varphi^{-1}(p_0)}$ is an isomorphism from $\RR^{n + 2 - k} \times V^{k-1, I}$ to $T_{p_0}M$ taking $\RR^{n + 2 - k} \times \{\0^{k-1}\}$ to $T_{p_0}\mathring{M}_{I}$.
\end{proof}

We can now prove Proposition~\ref{prop:fixed-stratum-patching}, restated as the following result.
\begin{prop}
\label{prop:global-multisection-chart}
Let $\Omega$ be any relatively open subset with compact closure in $\mathring{M}_{I}$, and $\cC$ a closed subset of $M$ not intersecting $\mathring{M}_{I}$. Then there exists a diffeomorphism $G$ from a neighborhood $\cW_1$ of $\Omega$ in $M$ onto a neighborhood $\cW_2$ of $\Omega \times \{\0^{k-1}\}$ in $\Omega \times V^{k-1, I}$, such that
\vskip 1mm
\begin{enumerate}
\item[(a)] $\overline{\cW_{1}} \cap \big[ (\cup_{i \not\in I}M_{i} ) \cup \cC \big] = \emptyset$.
\vskip 1mm
\item[(b)] $\cW_{2}$ has the form $\Omega \times B^{k-1, I}_{\rho}$ for some $\rho > 0$.
\vskip 1mm
\item[(c)] $G(q) = (q, \0^{k-1})$ for all $q \in \Omega$.
\vskip 1mm
\item[(d)] Given $q \in \cW_{1}$ and $i \in I$, we have that
\[
q \in M_{i} \quad\text{if and only if}\quad G(q) \in \Omega \times V^{k-1, I}_{i}.
\]
\end{enumerate}
\end{prop}
\begin{rmk}\label{rmk:global-multisection-chart}
Since $G$ is a bijection from $\cW_{1}$ onto $\cW_{2}$, part (d) above is equivalent to the statement that
\[
G(\cW_1 \cap M_{i}) = \cW_{2} \cap (\Omega \times V^{k-1, I}_{i}) \quad \text{for all }i \in I.
\]
\end{rmk}
\begin{proof}[Proof of Proposition~\ref{prop:global-multisection-chart}]
Fix any $i_0 \in I$. By Corollary~\ref{coro:trivial-quotient-bundle}, the sections $\{\xi_i\}_{i \in I \setminus \{i_0\}}$ defined by~\eqref{eq:section-by-adapted-chart} form a global frame for $T^{\perp}\mathring{M}_{I}$. With $\exp^{\perp}:T^{\perp}\mathring{M}_{I} \to M$ denoting the normal exponential map of $\mathring{M}_{I}$, we define $F : \mathring{M}_{I} \times V^{k-1, I} \rightarrow M$ by
\[
F(p, \sum_{i \in I\setminus \{i_0\}}v_i\ba_{i; I}) = \exp^{\perp}(p, \sum_{i\in I\setminus \{i_0\}} v_i\xi_{i}(p) ).
\]
Next, fix open subsets $\Omega_1$, $\Omega_{2}$, and $\Omega_{3}$ relative to $\mathring{M}_{I}$, all with compact closures in $\mathring{M}_{I}$, such that 
\[
\overline{\Omega} \subset \Omega_1,\quad \overline{\Omega_1} \subset \Omega_2,\quad\overline{\Omega_2} \subset \Omega_3.
\]
By the compactness of $\overline{\Omega_3}$, there exists $\eta > 0$ such that $F$ maps $\Omega_{3} \times B^{k-1, I}_{\eta}$ diffeomorphically onto a neighborhood $\cU_0$ of $\Omega_{3}$ in $M$. By~\eqref{eq:boundary-sector} and our assumption about $\cC$, decreasing $\eta$ if necessary, we may also arrange for $F(\overline{\Omega_{3}} \times \overline{B^{k-1, I}_{\eta}})$ to be disjoint from $\cC$ and from $M_{i}$ for each $i \not\in I$, so that 
\begin{equation}\label{eq:adjacent-sectors}
\overline{\cU_0} \cap (\cup_{i \not\in I}M_{i}) = \emptyset, \quad\quad \overline{\cU_0}\cap \cC = \emptyset.
\end{equation}
For the remainder of the proof, we write $F$ for the restriction $F|_{\Omega_3 \times B_{\eta}^{k-1, I}}$.

Since $\overline{\Omega_{2}}$ is compact, there exist finitely many adapted charts $\varphi_{a}: \widetilde{U}_{a} \to U_{a}$, for $a = 1, \cdots, L$, such that each $U_{a}$ has the form $F(Z_{a} \times B^{k-1, I}_{\delta_{a}})$ for some open set $Z_{a} \subset \Omega_{3}$ and radius $\delta_{a} \in (0, \eta)$, and that 
\[
\overline{\Omega_{2}} \subset \bigcup_{a = 1}^{L}Z_{a}.
\]
For all $a \in \{1, \cdots, L\}$, note that $\varphi_{a}^{-1} \circ F$ can be expressed as 
\begin{equation}\label{eq:za-wa-definitions}
\varphi_{a}^{-1}\big( F(p, v) \big) = (z_{a}(p, v), w_{a}(p, v)) \quad\text{for }(p, v) \in Z_a \times  B^{k-1, I}_{\delta_{a}},
\end{equation}
where $z_{a}$ and $w_{a}$ take values in $\RR^{n + 2 - k}$ and $V^{k-1, I}$, respectively. We then let $\{\zeta_a\}$ be a partition of unity on $\overline{\Omega_{2}}$ subordinate to $\{Z_a\}_{a = 1}^{L}$, put
\[
\eta_0:= \min_{1 \leq a \leq L}\delta_{a},
\]
and define $G_{1}: \Omega_{2} \times B^{k-1, I}_{\eta_0} \to \Omega_{2} \times V^{k-1, I}$ by letting
\begin{equation}\label{eq:G-defi}
G_{1}(p, v) = (p, \sum_{a = 1}^{L}\zeta_{a}(p) \cdot w_{a}(p, v)).
\end{equation}
This is a smooth map, with $\zeta_{a}(p) \cdot w_{a}(p, v)$ understood to mean $\0^{k-1}$ when $p \not\in Z_{a}$. Also, by~\eqref{eq:modified-adapted-chart} we have $w_{a}(p, \0^{k-1}) = \0^{k-1}$ for all $p \in Z_a$, and hence
\begin{equation}\label{eq:G1-id-on-Omega}
G_1(p, \0^{k-1}) = (p, \0^{k-1}) \quad \text{for all }p \in \Omega_{2}.
\end{equation}
To study the derivative of $G_1$ at such points, take $p \in \Omega_{2}$ and $a \in \{1, \cdots, L\}$ such that $p \in Z_{a}$. By~\eqref{eq:za-wa-definitions} we have
\[
\exp^{\perp}\big( p, \sum_{i \in I \setminus \{i_0\}} v_i \xi_i(p) \big) = \varphi_{a}(z_a(p, v), w_a(p, v)), \quad\text{where } v = \sum_{i \in I \setminus \{i_0\}} v_i \ba_{i; I}.
\]
Differentiating with respect to $v_i$ at $v = \0^{k-1}$, and using~\eqref{eq:section-by-adapted-chart}, we find that
\[
\begin{split}
(d\varphi_a)_{\varphi_a^{-1}(p)}({\bf 0}^{n+2-k}, \lambda_i\cdot \ba_{i;I})^{\perp}
=\ & (d\varphi_{a})_{\varphi_a^{-1}(p)}(\pa{z_a}{v_i}, \pa{w_a}{v_i}),
\end{split}
\]
where $\lambda_i = \big(\| (d\varphi_a)_{\varphi_a^{-1}(p)}({\bf 0}^{n+2-k}, \ba_{i;I})^{\perp} \|\big)^{-1}$, and $\pa{w_a}{v_i}$ stands for $(dw_{a})_{(p, \0^{k-1})}(0, \ba_{i; I})$. Since $(d\varphi_{a})_{\varphi_{a}^{-1}(p)}(\RR^{n+2 - k} \times \{\0^{k-1}\}) = T_{p}\mathring{M}_{I}$, we infer from the above that
\[
(dw_{a})_{(p, \0^{k-1})}(0, \ba_{i; I}) = \lambda_{i} \ba_{i; I}.
\]
Combining this with the definition of $G_1$, we get for each $i \in I \setminus \{i_0\}$ some $c_i > 0$ so that
\[
(dG_1)_{(p, \0^{k-1})}(\0^{n+2-k}, \ba_{i; I}) = c_i\cdot ({\bf 0}^{n+ 2-k}, \ba_{i; I}).
\]
This together with~\eqref{eq:G1-id-on-Omega} implies that $(dG_{1})_{(p, \0^{k-1})}$ is invertible for any $p \in \Omega_{2}$. Recalling that $\overline{\Omega_{1}}$ is a compact subset of $\Omega_{2}$, and using~\eqref{eq:G1-id-on-Omega} again, we obtain $\eta_1 < \eta_0$ such that $G_{1}$ maps $\Omega_{1} \times B^{k-1, I}_{\eta_1}$ diffeomorphically onto some neighborhood of $\Omega_{1} \times \{\0^{k-1}\}$ in $\Omega_{2} \times V^{k-1, I}$. Similar to what we did with $F$, below we write $G_1$ for the restriction $G_1|_{\Omega_{1} \times B_{\eta_1}^{k-1, I}}$.

Towards defining $\cW_{1}$, $\cW_{2}$, and $G$, we observe that given $J \subset I$ and $(p, v) \in (\Omega_{1} \times B_{\eta_1}^{k-1, I}) \cap F^{-1}(\mathring{M}_{J})$, as well as $a \in \{1, \cdots, L\}$ such that $p \in Z_{a}$, by the inclusion $U_a \subset \cU_0$ and~\eqref{eq:adjacent-sectors} we have
\[
F(p, v) \in U_{a} \cap \mathring{M}_{J} = U_{a} \cap \big( \cap_{i \in J} M_{i} \big) \setminus \big( \cup_{i \in I \setminus J}M_{i} \big),
\]
and hence~\eqref{eq:za-wa-definitions} and~\eqref{eq:modified-adapted-chart} implies
\[
 w_{a}(p, v) \in \big( \cap_{i \in J}{V}^{k-1, I}_{i} \big) \setminus \big( \cup_{i \in I \setminus J}{V}^{k-1, I}_{i}\big) =  \mathring{{V}}^{k-1, I}_{J}.
\]
Using the convexity of $\mathring{{V}}^{k-1, I}_{J}$, we arrive at
\begin{equation}\label{eq:G1-respects-multisection}
G_{1}(p, v) \in \Omega_{2} \times \mathring{{V}}^{k-1, I}_{J} \quad \text{for all }(p, v) \in  (\Omega_{1} \times B_{\eta_1}^{k-1, I}) \cap F^{-1}(\mathring{M}_{J}).
\end{equation}
Recalling from the previous paragraph that $G_1(\Omega_1 \times B_{\eta_1}^{k-1,I})$ is a neighborhood of $\Omega_1 \times \{\0^{k-1}\}$, we may use the compactness of $\overline{\Omega} \subset \Omega_1$ to choose $\rho > 0$ such that 
\[
\cW_{2} : = \Omega \times B_{\rho}^{k-1, I} \subset G_{1}(\Omega_{1} \times B^{k-1, I}_{\eta_{1}}).
\]
Notice that $\cW_{2}$ is mapped diffeomorphically by $F\circ G_{1}^{-1}$ onto a neighborhood of $\Omega$ contained in $\cU_0$, which we define to be $\cW_{1}$; that is, we set
\[
\cW_{1}: = (F\circ G_1^{-1})(\cW_{2}).
\]
Part (a) then follows from noting that $\cW_{1} \subset \cU_0$, and using~\eqref{eq:adjacent-sectors}, while part (b) holds by construction. Letting 
\[
G: = G_{1} \circ F^{-1},
\]
we see that it is a diffeomorphism from $\cW_{1}$ onto $\cW_{2}$, and that part (c) holds by~\eqref{eq:G1-id-on-Omega}. For part (d), we first prove that 
\begin{equation}\label{eq:global-multisection-chart-J}
G(\cW_1 \cap \mathring{M}_{J}) \subset  \cW_{2} \cap (\Omega \times \mathring{V}^{k-1, I}_{J}) \quad \text{for all }J \subset I.
\end{equation}
To that end, take $J \subset I$ and note that for all $q \in \cW_{1} \cap \mathring{M}_{J}$ we have
\[
(p, v) : = F^{-1}(q) \in  (\Omega_{1} \times B_{\eta_1}^{k-1, I}) \cap F^{-1}(\mathring{M}_{J}).
\]
Since $G(q)$ lies in $\cW_{2}$ and also coincides with $G_{1}(p, v)$, we infer from~\eqref{eq:G1-respects-multisection} that~\eqref{eq:global-multisection-chart-J} holds. Now take $q \in \cW_{1}$ and define 
\[
I' = \{1 \leq i \leq n\ |\ q \in M_{i}\}. 
\]
Then $q \in \cW_{1} \cap \mathring{M}_{I'}$, and by part (a) we have $I' \subset I$. We then use~\eqref{eq:global-multisection-chart-J} to get
\[
G(q) \in \Omega \times \mathring{V}^{k-1, I}_{I'}. 
\]
For any $i \in I$, if $q$ lies in $M_{i}$, then $i \in I'$, in which case $G(q)$ lies in $\Omega \times V^{k-1, I}_{i}$. On the other hand, if $q \not\in M_{i}$, then we have $i \in I \setminus I'$, and thus $G(q) \not\in \Omega \times V^{k-1, I}_{i}$. In other words, we have established the assertion of part (d), and the proof of Proposition~\ref{prop:global-multisection-chart} is complete.
\end{proof}
\begin{rmk}\label{rmk:extension}
Expressing $G:\cW_{1} \to \Omega \times B_{\rho}^{k-1, I}$ in components as $G(q) = (A(q), X(q))$, we see from Proposition~\ref{prop:global-multisection-chart}(c) that $A(q) = q$ for all $q \in \Omega$. Recalling that $M$ is embedded in $\RR^{Q}$, and denoting by $\Pi$ the nearest-point projection from a tubular neighborhood onto $M$, then the map from $\Pi^{-1}(\cW_{1}) \times B_{\rho}^{k-1, I}$ to $\cW_{1}$ given by
\[
(q, v) \mapsto  G^{-1}((A \circ \Pi)(q), v),
\]
is smooth, and restricts to $G^{-1}$ on $\Omega \times B_{\rho}^{k-1, I}$. Noting that the domain $\Pi^{-1}(\cW_{1}) \times B_{\rho}^{k-1, I}$ of the above map is open in $\RR^{Q} \times V^{k-1, I}$, while the target $\cW_{1}$ is embedded in $\RR^{Q}$, we infer that for any $\sigma < \rho$ and compact subset $K \subset \Omega$, the map $G^{-1}\big|_{K \times \overline{B_{\sigma}^{k-1, I}}}$ is Lipschitz in the usual sense defined with Euclidean distances.
\end{rmk}
\subsection{Patching across different strata (I)}\label{subsec:different-strata-I}
We continue to work in the setting of Appendix~\ref{subsec:fixed-stratum}, and write $S_0$ for the central surface $M_{\{1, \cdots, n\}}$. 
Below are the main building blocks for the construction in this section.
\vskip 1mm
\begin{enumerate}
\item[(i)] By Proposition~\ref{prop:global-multisection-chart}, switching the order of factors in the target for convenience, we obtain a neighborhood $\cU_0$ of $S_0$ in $M$, and a diffeomorphism $f_0 : \cU_0 \to B^{n-1}_{1} \times S_0$, such that 
\begin{equation}\label{eq:f1234-center}
f_{0}^{-1} (\0^{n-1}, x) = x,\quad \text{for all }x \in S_0,
\end{equation}
\begin{equation}\label{eq:f1234-compatible}
f_{0}^{-1}((B_{1} \cap V^{n-1}_{i}) \times S_0) = M_i \cap \cU_{0},\quad \text{for any }i\in \{1, 2, \cdots, n\},
\end{equation}
and we let $\cW_{0} := f^{-1}_{0}(B^{n-1}_{1/8} \times S_0)$.
\vskip 1mm
\item[(ii)] For each $I \subset \{1, \cdots, n\}$ with $|I| = n-1$, we let 
\[
M_{I}^{*} : = M_{I} \setminus \overline{\cW_{0}}.
\]
Since $\mathring{M}_I = M_I \setminus S_0$ and $S_0 \subset \cW_0$, we have
\[
M_{I}^{*}= \mathring{M}_{I} \setminus \overline{\cW_{0}} \subset  \mathring{M}_I \setminus \cW_0 = M_{I} \setminus \cW_0,
\]
which shows that $M_{I}^*$ is open in $\mathring{M}_I$, and is contained in a compact subset of the latter. In particular, Proposition~\ref{prop:global-multisection-chart} is applicable, and we obtain a neighborhood $\cU_{I}$ of $M_I^{*}$ in $M$, and a diffeomorphism $f_{I}: \cU_{I} \to B^{n-2, I}_{1} \times M^{*}_{I}$, satisfying that
\begin{equation}\label{eq:f123-center}
f_I^{-1}(\0^{n-2}, q) = q,\quad \text{for all }q \in M_{I}^{*},
\end{equation}
\begin{equation}\label{eq:f123-compatible}
f_{I}^{-1}((B_{1} \cap V^{n-2, I}_{i}) \times M_I^{*}) = M_{i} \cap \cU_{I}, \quad\text{for all }i \in I.
\end{equation}
By Proposition~\ref{prop:global-multisection-chart}(a), we can also assume that
\begin{equation}\label{eq:U-disjoint}
\overline{\cU_{I}} \cap M_{i} =\emptyset\quad\text{if }i \not\in I,\quad\quad \overline{\cU_{I}} \cap \overline{\cU_{J}} = \emptyset \quad \text{if }I \neq J.
\end{equation}
\vskip 1mm
\item[(iii)] Fix $i_0 \in \{1, \cdots, n\}$ and let $I = \{i_0\}^{c}$, where we recall that $\{i_0\}^{c}$ stands for $\{1, \cdots, n\} \setminus \{i_{0}\}$. By~\eqref{eq:decomposition-wrt-lower-simplex}, there holds the orthogonal decomposition $V^{n-1} = V^{n-2, I} \oplus \Span\{\ba_{i_0}\}$. Accordingly, we introduce the notation
\begin{equation}\label{eq:split-coordinates}
(x, t): = x + t\frac{\ba_{i_0}}{|\ba_{i_0}|},\quad\text{for } t\in \RR \text{ and }x \in V^{n-2, I}.
\end{equation}
Then~\eqref{eq:f1234-compatible} implies that, for all $r \in (0, 1]$,
\begin{equation}\label{eq:1234-parametrize-collar}
f_{0}^{-1}(\{\0^{n-2}\} \times [0, r) \times S_0) = M_{I} \cap f_{0}^{-1}(B^{n-1}_{r} \times S_0),
\end{equation}
and these are all collar neighborhoods of $S_0$ in the $3$-dimensional $1$-handlebody $M_{I}$. Choosing the one corresponding to $r = 7/8$, then as noted in Section~\ref{subsec:retraction}, we may use the sets $\Omega$ and $S$ from Lemma~\ref{lemm:Omega-S-basics} and the map $\bh$ from Proposition~\ref{prop:bh-properties} as essential building blocks, taking $m = 2$ in both results, to construct a map
\[
[0, 1] \times f_{0}^{-1}(\{\0^{n-2}\} \times \{7/8\} \times S_0) \longrightarrow M_{I} \setminus f_{0}^{-1}(\{\0^{n-2}\} \times [0, 7/8) \times S_0),
\]
which is Lipschitz continuous, among other things. Changing the domain to $[7/8, 2] \times S_0$ via the reparametrization 
\[
(t, p) \mapsto  (\frac{8t - 7}{9}, f_{0}^{-1}(\0^{n-2}, 7/8, p)),
\]
and joining the result with $f_{0}^{-1}(\0^{n-2},\cdot, \cdot)|_{[0, 7/8] \times S_0}$, we obtain a Lipschitz map 
\[
h = h_{I}:[0, 2] \times S_0 \to M_{I}
\]
which has the following properties:
\vskip 1mm
\begin{enumerate}
\item[(1)] $t \mapsto h(t, \cdot)$ is continuous from $[0, 2)$ into $C^{1}(S_0; M)$.
\vskip 1mm
\item[(2)] $\Gamma = \Gamma_{I}:= h(\{2\} \times S_0)$ has finite $\cH^{1}$-measure. Also, $\Area(h(t, \cdot)) \to 0$ as $t \to 2^{-}$.
\vskip 1mm
\item[(3)] $M_I = h([0, 2] \times S_0)$. Also, 
\begin{equation}\label{eq:h-disjoint-images}
h(\{t\} \times S_0) \cap h(\{t'\} \times S_0) = \emptyset,\quad\text{whenever }t \neq t'.
\end{equation}
\vskip 1mm
\item[(4)] $h(t, \cdot) = f_{0}^{-1}(\0^{n-2}, t, \cdot)$ for all $t \in [0, 7/8]$.
\end{enumerate}
\vskip 1mm

By properties (3) and (4), as well as~\eqref{eq:1234-parametrize-collar}, we have
\begin{equation}\label{eq:M_I*-expression}
h((1/8, 2] \times S_0) = M_I \setminus h([0, 1/8] \times S_0 ) = M_I^*.
\end{equation}
We then define
\[
\phi = \phi_{I}:B^{n-2, I}_{1} \times (1/8, 2]\times S_0 \to B^{n-2, I}_{1} \times M_{I}^{*}
\]
by
\begin{equation}\label{eq:collar-spread-out}
\phi(x, t, p) = (x, h(t, p)).
\end{equation}
By property (4) above, $\phi$ restricts to a diffeomorphism from $B^{n-2, I}_{1} \times (1/8, 7/8) \times S_0$ onto an open subset of $B^{n-2, I}_{1} \times M_{I}^*$, and that, given $(t, p) \in (1/8, 7/8) \times S_0$, there holds
\begin{equation}\label{eq:agreement-on-central-leaf}
(f_{I}^{-1}\circ \phi)(\0^{n-2}, t, p) = h(t, p) =  f_0^{-1}(\0^{n-2}, t, p).
\end{equation}
\end{enumerate}

\vskip 2mm
Given $I = \{i_0\}^{c}$ as above, along with a subset $A$ of $\RR$ and some $\sigma > 0$, we define 
\begin{equation}\label{eq:cylinder-definition}
\begin{split}
C_{\sigma}^{I}(A): =\ & B^{n-2, I}_{\sigma} + A\cdot \frac{\ba_{i_0}}{|\ba_{i_0}|} = \{x + t\frac{\ba_{i_0}}{|\ba_{i_0}|}\ |\ t \in A, \ x \in B^{n-2, I}_{\sigma}\}.
\end{split}
\end{equation}
When the notation~\eqref{eq:split-coordinates} is adopted, we write $C_{\sigma}^{I}(A)$ as $B^{n-2, I}_{\sigma} \times A$. Given positive numbers $\rho$, $\sigma$, and $R$, with $\sigma< \rho$, we let
\begin{equation}\label{eq:cN-definition}
\cN_{\rho, \sigma, R} := B_{\rho}^{n-1} \cup \bigcup_{|I| = n-1}C_{\sigma}^{I}((0, R]),\quad
\mathring{\cN}_{\rho, \sigma, R} : =  B_{\rho}^{n-1} \cup \bigcup_{|I| = n-1}C_{\sigma}^{I}((0, R)),
\end{equation}
dropping the subscript $R$ in the case $R = 2$. 
\begin{prop}\label{prop:1234-with-123}
There exist $\rho > \sigma > 0$ and a Lipschitz map $G:\cN_{\rho, \sigma} \times S_0 \to M$ such that $y \mapsto G(y, \cdot)$ is continuous from $\mathring{\cN}_{\rho, \sigma}$ into $C^{1}(S_0; M)$, and that the following hold.
\vskip 1mm
\begin{enumerate}
\item[(a)] $G|_{\mathring{\cN}_{\rho, \sigma, 7/8} \times S_0}$ is a diffeomorphism onto an open set in $M$. Also, given disjoint subsets $A, B \subset \cN_{\rho, \sigma}$ with $A \subset \cN_{\rho, \sigma, 2/3}$, we have 
\begin{equation}\label{eq:G-disjoint-images}
G(A \times S_0) \cap G(B \times S_0) = \emptyset.
\end{equation}
\vskip 1mm
\item[(b)] For all $y_0 \in \cN_{\rho, \sigma}\setminus \mathring{\cN}_{\rho, \sigma}$, we have $\cH^{1}(G(\{y_0\} \times S_0)) < \infty$, and that
\begin{equation}\label{eq:G-area-boundary}
\Area(G(y, \cdot)) \to 0\quad \text{as }y \to y_0 \text{ from within } \mathring{\cN}_{\rho, \sigma}.
\end{equation}
\vskip 1mm
\item[(c)]  For all $i_0 \in \{1, \cdots, n\}$, letting $I= \{i_0\}^{c}$, we have 
\begin{equation}\label{eq:handle-sweepout-on-core}
G(t \frac{\ba_{i_0}}{|\ba_{i_0}|}, p) = h_{I}(t, p),\quad \text{for }(t, p) \in [0, 2] \times S_0.
\end{equation}
In particular, $G(\0^{n-1}, p) = p$ for all $p \in S_0$.
\vskip 1mm
\item[(d)] Given $(y, p) \in \cN_{\rho, \sigma} \times S_0$ and $j \in \{1, \cdots, n\}$, we have 
\begin{equation}\label{eq:G-multisection-iff}
y \in V^{n-1}_{j} \quad \text{if and only if}\quad G(y, p) \in M_{j}.
\end{equation}
\end{enumerate}
\end{prop}
\begin{proof}
Define the following intervals
\begin{equation}\label{eq:intervals-choices}
A_2 = (3/8, 5/8),\quad A_1 = (1/3, 2/3),\quad A_0 = (7/24, 17/24), \quad A = (1/4, 3/4).
\end{equation}
Given $i_0 \in \{1, \cdots, n\}$, letting $I := \{i_0\}^{c}$ as above, we have by~\eqref{eq:agreement-on-central-leaf} that
\[
(f_{I}^{-1}\circ \phi)(\{\0^{n-2}\} \times \overline{A} \times S_0)\subset \cU_{I} \cap f_{0}^{-1}\big((B_{1}^{n-1} \setminus \overline{B^{n-1}_{1/8}}) \times S_0\big),
\]
so there exists some $\lambda_0 = \lambda_{0, I} \in (0, 1/8)$ such that 
\begin{equation}\label{eq:transition-inclusion-1}
\begin{split}
(f_{I}^{-1} \circ \phi)(\overline{B^{n-2, I}_{\lambda_0}} \times \overline{A} \times S_0) \subset\ & \cU_{I} \cap f_{0}^{-1}\big((B_{1}^{n-1} \setminus \overline{B^{n-1}_{1/8}}) \times S_0\big)\\
\subset\ & \cU_{I} \cap (\cU_0 \setminus \overline{\cW_0}).
\end{split}
\end{equation}
Using~\eqref{eq:agreement-on-central-leaf} again, we get $r_0 = r_{0, I} < \lambda_{0, I}$ such that
\begin{equation}\label{eq:transition-inclusion-2}
f_{0}^{-1}(\overline{B^{n-2, I}_{r_0}} \times \overline{A_0} \times S_0) \subset (f_{I}^{-1} \circ \phi)(B^{n-2, I}_{\lambda_0} \times A \times S_0).
\end{equation}
In particular, it makes sense to define the transition map
\[
\theta = \theta_{I} = (\overline{x}, \overline{t}, \overline{p}): B^{n-2, I}_{r_0} \times A_0 \times S_0 \longrightarrow B^{n-2, I}_{\lambda_0} \times A \times S_0
\]
by
\[
\theta : = (\phi|_{B^{n-2, I}_{1} \times (1/8, 7/8) \times S_0})^{-1} \circ f_I \circ f_0^{-1},
\]
which maps $B^{n-2, I}_{r_0} \times A_0 \times S_0$ diffeomorphically onto its image.
\begin{claim}\label{claim:theta-properties}
The map $\theta$ has the following properties.
\vskip 1mm
\begin{enumerate}
\item[(a)] For all $(t, p) \in A_0 \times S_0$, there holds $\theta(\0^{n-2}, t, p) = (\0^{n-2}, t, p)$.
\vskip 1mm
\item[(b)] Given $J \subset I$ and some point $(x, t, p)$ in the domain of $\theta$ satisfying $x \in \mathring{V}^{n-2, I}_{J}$, we have $\overline{x}(x, t, p) \in \mathring{V}^{n-2, I}_{J}$.
\vskip 1mm
\item[(c)] For any $i \in I$ and $(t, p) \in A_0 \times S_0$, there exists $c_{i} = c_{i}(t, p) > 0$ such that 
\begin{equation}\label{eq:x-bar-derivative}
(d\overline{x})_{(\0^{n-2}, t, p)}(\ba_{i; I}, 0, 0) = c_{i}\cdot \ba_{i; I}.
\end{equation}
\end{enumerate}
\end{claim}
\begin{proof}[Proof of Claim]
Part (a) is an immediate consequence of~\eqref{eq:agreement-on-central-leaf}. For part (b), when $J = I$, in which case $\mathring{V}^{n-2, I}_{J} = \{\0^{n-2}\}$, we get the desired conclusion from part (a), while when $J \subsetneq I$, we may apply~\eqref{eq:projection-of-sectors} to see that $(x, t) \in \mathring{V}^{n-1}_{J}$, so that by~\eqref{eq:transition-inclusion-2},~\eqref{eq:f1234-compatible}, and~\eqref{eq:U-disjoint}, we have
\[
f_{0}^{-1}(x, t, p) \in \cU_{I} \cap \mathring{M}_{J} = \cU_{I} \cap M_{J} \setminus \big( \cup_{i \in I \setminus J}M_{i} \big).
\]
Recalling from the definitions of $\theta$ and $\phi_{I}$ that 
\[
f_0^{-1}(x, t, p) = f_{I}^{-1}(\overline{x}, h_{I}(\overline{t}, \overline{p})),
\]
and then using~\eqref{eq:f123-compatible}, we deduce that
\[
\overline{x}(x, t, p) \in V^{n-2,I}_{J} \setminus \big( \cup_{i \in I \setminus J} V^{n-2, I}_{i}\big) = \mathring{V}^{n-2, I}_{J},
\]
as asserted. For part (c), we apply (b) with $J = I \setminus \{i\}$ and observe by (a) that $\overline{x}(\0^{n-2}, t, p) = \0^{n-2}$. These together imply that
\[
(d\overline{x})_{(\0^{n-2}, t, p)}(\ba_{i; I}, 0, 0) = c_i\cdot \ba_{i;I},\quad\text{for some }c_i \geq 0.
\]
In view of (a) and the invertibility of $(d\theta)_{(\0^{n-2}, t, p)}$, the constant $c_i$ cannot be zero. This proves part (c), and we are done.
\end{proof}
Next we choose a cut-off function $\zeta \in C^{\infty}(\RR; [0, 1])$ satisfying
\begin{equation}\label{eq:cut-off-towards-3-strata}
\zeta(t) = 1 \text{ if }t \leq 3/8,\quad \zeta(t) = 0 \text{ if }t \geq 5/8,
\end{equation}
and define 
\[
F = F_{I} = (\widetilde{x}, \widetilde{t}, \widetilde{p}): B^{n-2, I}_{r_0} \times A_0 \times S_0 \to B^{n-2, I}_{\lambda_0} \times A \times S_0
\]
by
\begin{equation}\label{eq:F-I-definition}
\begin{split}
\widetilde{x}(x, t, p) =\ & (1 - \zeta(t))\cdot x + \zeta(t)\cdot \overline{x}(x,t ,p),\\
\widetilde{t}(x, t, p) =\ & (1 - \zeta(t))\cdot t + \zeta(t)\cdot \overline{t}(x,t ,p),\\
\widetilde{p}(x, t, p) =\ & \overline{p}(\zeta(t)x, t, p).
\end{split}
\end{equation}
Since $B_{r_0}^{n-2, I}$ and $A_0$ are contained respectively in $B^{n-2, I}_{\lambda_0}$ and $A$, which are both convex sets, we see that $F$ does have the stated codomain. By Claim~\ref{claim:theta-properties}(a) we have 
\begin{equation}\label{eq:F-id-on-center}
F(\0^{n-2}, t, p) = (\0^{n-2}, t, p)\quad \text{for all } (t, p) \in A_0 \times S_0,
\end{equation}
which combines with~\eqref{eq:agreement-on-central-leaf} to give
\begin{equation}\label{eq:images-coincide-on-center}
(f_{I}^{-1} \circ \phi \circ F)(\0^{n-2}, t, p) = (f_{I}^{-1} \circ \phi)(\0^{n-2}, t, p) = f_0^{-1}(\0^{n-2}, t, p).
\end{equation}
Next, by our choice of $\zeta$ we have
\begin{equation}\label{eq:F-endpoints}
(f_{I}^{-1} \circ \phi \circ F)(x, t, p) = \left\{
\begin{array}{ll}
f_{0}^{-1}(x, t, p),& \text{ if }t \leq 3/8,\\
(f_{I}^{-1}\circ \phi)(x, t, p), & \text{ if }t \geq 5/8.
\end{array}
\right.
\end{equation}
\begin{claim}\label{claim:F-I-multisection}
Given $(x, t, p) \in B_{r_0}^{n-2, I} \times A_0 \times S_0$ and $j \in I$, we have 
\[
x \in V^{n-2, I}_{j} \quad \text{if and only if}\quad(f_{I}^{-1}\circ \phi \circ F)(x, t, p) \in M_{j}.
\]
\end{claim}
\begin{proof}
We first prove that 
\begin{equation}\label{eq:interpolating-map-multisection}
(f_{I}^{-1} \circ \phi \circ F)(x ,t, p) \in  \mathring{M}_J,
\end{equation}
whenever $J \subset I$ and $(x, t, p) \in (B^{n-2, I}_{r_0} \cap \mathring{V}^{n-2, I}_{J}) \times A_0 \times S_0$. Indeed, in this case we have by Claim~\ref{claim:theta-properties}(b) and the convexity of $\mathring{V}^{n-2, I}_{J}$ that
\[
\widetilde{x}(x, t, p)\in B^{n-2, I}_{\lambda_0} \cap \mathring{V}^{n-2, I}_{J},
\]
which implies by~\eqref{eq:f123-compatible} that
\[
(f_{I}^{-1} \circ \phi \circ F)(x ,t, p) \in \cU_I \cap M_{J} \setminus \big( \cup_{i \in I\setminus J}M_{i} \big).
\]
Combining this with~\eqref{eq:U-disjoint} yields~\eqref{eq:interpolating-map-multisection}. To derive the claim from this, take $(x, t, p) \in B^{n-2, I}_{r_0} \times A_0 \times S_0$ and define
\[
I' = \{i \in I\ |\ x \in V^{n-2, I}_{i}\}.
\]
Then we have $x \in \mathring{V}^{n-2, I}_{I'}$, and~\eqref{eq:interpolating-map-multisection} implies $(f_I^{-1}\circ \phi \circ F)(x, t, p) \in \mathring{M}_{I'}$. Consequently, given $j \in I$, both ends in the claimed equivalence reduce to $j \in I'$, and we are done.
\end{proof}

\begin{claim}\label{claim:F-diffeo}
There exists $r_1  = r_{1, I} < r_0$ such that $F$ maps $B^{n-2, I}_{r_1} \times A_1 \times S_0$ diffeomorphically onto a neighborhood of $\{\0^{n-2}\} \times A_1 \times S_0$ inside $B^{n-2, I}_{\lambda_0} \times A \times S_0$. 
\end{claim}
\begin{proof}[Proof of Claim]
By Claim~\ref{claim:theta-properties}(c), we have for all $(t, p) \in A_0 \times S_0$ that 
\[
(d\widetilde{x})_{(\0^{n-2}, t, p)}(\ba_{i; I}, 0, 0) = [(1 - \zeta(t)) + c_i \zeta(t)] \cdot \ba_{i; I},\quad\text{for }i \in I.
\]
Combining this with~\eqref{eq:F-id-on-center} shows that $(dF)_{(\0^{n-2}, t, p)}$ is invertible for all $(t, p) \in A_0 \times S_0$, and thus there exists $r_1 < r_0$ such that $F$ is a local diffeomorphism at each point in $\overline{B^{n-2, I}_{r_1}} \times \overline{A_1} \times S_0$. It remains to prove that, upon decreasing $r_1$ if necessary, we can also arrange $F$ to be injective on $\overline{B^{n-2, I}_{r_1}} \times \overline{A_1} \times S_0$. Indeed, suppose there is no such $r_1$, then there exist sequences $x_i, x_i' \to \0^{n-2}$ and $(t_i, p_i), (t_i', p_i')$ in $\overline{A_1} \times S_0$ such that 
\begin{equation}\label{eq:F-injective-contradiction}
(x_i, t_i, p_i) \neq (x_i', t_i', p_i'), \quad\text{but } F(x_i, t_i, p_i) = F(x_i', t_i', p_i').
\end{equation}
Passing to subsequences if needed, we get $(t, p), (t', p') \in \overline{A_1} \times S_0$ such that 
\[
(t_i, p_i) \to (t, p),\quad \text{and }(t_i', p_i') \to (t', p').
\]
Taking the limit in the second relation in~\eqref{eq:F-injective-contradiction} and using~\eqref{eq:F-id-on-center}, we find that $(t, p) = (t', p')$, and thus~\eqref{eq:F-injective-contradiction} implies that $F$ is not injective on any neighborhood of $(\0^{n-2}, t, p)$, contradicting the fact that $F$ is a local diffeomorphism at that point. The claim is proved.
\end{proof}
\begin{claim}\label{claim:image-of-F}
Let $r_1$ be as in Claim~\ref{claim:F-diffeo} and fix any $\rho$ in $(\frac{1}{3}, \frac{3}{8})$, that is, between the left endpoint of $A_1$ and that of $A_2$. There exists $r_2 = r_{2, I} < \min\{r_1, \frac{\rho}{3}\}$ such that the following hold.
\vskip 1mm
\begin{enumerate}
\item[(a)] The three sets below are pairwise disjoint:
\begin{equation}\label{eq:image-of-F-sets}
\begin{split}
&(f_{I}^{-1}\circ \phi)(\overline{B_{r_2}^{n-2, I}} \times (5/8, 2] \times S_0), \quad (f_{I}^{-1} \circ \phi \circ F)\big(\overline{B^{n-2, I}_{r_2}} \times \overline{A_2} \times S_0\big),\\
&f_{0}^{-1}\big((\overline{B_{\rho}^{n-1}} \times S_0) \cup (\overline{B_{r_2}^{n-2,I}} \times (1/3, 3/8)\times S_0)\big).
\end{split}
\end{equation}
\vskip 1mm
\item[(b)] $(f_{I}^{-1} \circ\phi)(\overline{B_{r_2}^{n-2, I}} \times [7/8, 2] \times S_0) \cap (f_{I}^{-1} \circ\phi)(\overline{B_{r_2}^{n-2, I}} \times [5/8, 2/3] \times S_0) = \emptyset$.
\vskip 1mm
\item[(c)] $\overline{B_{r_2, I}^{n-2}} \times [0, 1/3] \subset B^{n-1}_{\rho}$.
\end{enumerate}
\end{claim}
\begin{proof}[Proof of Claim]
Recall from~\eqref{eq:images-coincide-on-center} that, for all $(t, p) \in A_0 \times S_0$,
\[
(f_{I}^{-1} \circ \phi \circ F)(\0^{n-2}, t, p) = (f_{I}^{-1} \circ \phi)(\0^{n-2}, t, p) = f_0^{-1}(\0^{n-2}, t, p).
\]
By the equality between the left and right-most terms, and the fact that $f_0^{-1}$ is injective, as well as our choice of $\rho$ and $A_2$, we get
\begin{equation}\label{eq:disjoint-image-with-1234}
(f_{I}^{-1}\circ \phi\circ F)(\{\0^{n-2}\} \times \overline{A_2} \times S_0) \cap f_{0}^{-1}(\overline{B_{\rho}^{n-1}} \times S_0) = \emptyset.
\end{equation}
Combining~\eqref{eq:images-coincide-on-center} with~\eqref{eq:h-disjoint-images} instead, we have 
\begin{equation}\label{eq:disjoint-image-with-123-final}
(f_{I}^{-1}\circ \phi\circ F)(\{\0^{n-2}\} \times \overline{A_2} \times S_0) \cap (f_{I}^{-1} \circ \phi)(\{\0^{n-2}\} \times [2/3, 2] \times S_0) = \emptyset.
\end{equation}
Next, observe that
\[
\begin{split}
&(f_{I}^{-1}\circ \phi)(\{\0^{n-2}\} \times [5/8, 2] \times S_0)\\
=\ & (f_{I}^{-1}\circ \phi)(\{\0^{n-2}\} \times [5/8, 7/8) \times S_0) \cup (f_{I}^{-1}\circ \phi)(\{\0^{n-2}\} \times [7/8, 2] \times S_0) \\
=:\ & (I) \cup (II).
\end{split}
\]
Letting
\[
E := f_{0}^{-1}\big( (\overline{B_{\rho}^{n-1}} \times S_0) \cup (\{\0^{n-2}\} \times [1/3, 3/8]\times S_0)\big),
\]
we see by~\eqref{eq:agreement-on-central-leaf} and the injectivity of $f_0^{-1}$ that 
\begin{equation}\label{eq:set-I}
(I) \cap E = f_0^{-1}(\{\0^{n-2}\} \times [5/8, 7/8) \times S_0) \cap E = \emptyset.
\end{equation}
As for $(II)$, with the help of~\eqref{eq:1234-parametrize-collar} and properties (3) and (4) below it, we get
\[
(II) = M_{I} \setminus h([0, 7/8) \times S_0) = M_I \setminus f_{0}^{-1}(B^{n-1}_{7/8} \times S_0).
\]
Since $E \subset f_{0}^{-1}(B^{n-1}_{1/2} \times S_0)$, we infer that $(II) \cap E = \emptyset$, which together with~\eqref{eq:set-I} gives
\begin{equation}\label{eq:disjoint-I-II}
(f_{I}^{-1}\circ \phi)(\{\0^{n-2}\} \times [5/8, 2] \times S_0) \cap E = \emptyset.
\end{equation}
From~\eqref{eq:disjoint-image-with-1234},~\eqref{eq:disjoint-image-with-123-final}, and~\eqref{eq:disjoint-I-II}, along with a compactness argument, we deduce the existence of some $r_2 < \min\{r_1, \frac{\rho}{3}\}$ so that the following hold:
\vskip 1mm
\begin{enumerate}
\item[(D1)] $(f_{I}^{-1}\circ \phi\circ F)(\overline{B_{r_2}^{n-2, I}} \times \overline{A_2} \times S_0) \cap f_{0}^{-1}(\overline{B_{\rho}^{n-1}} \times S_0) = \emptyset$,
\vskip 1mm
\item[(D2)] $(f_{I}^{-1}\circ \phi\circ F)(\overline{B_{r_2}^{n-2, I}} \times \overline{A_2} \times S_0) \cap (f_{I}^{-1} \circ \phi)(\overline{B^{n-2, I}_{r_2}}\times [2/3, 2] \times S_0) = \emptyset$,
\vskip 1mm
\item[(D3)] $(f_{I}^{-1}\circ \phi)(\overline{B_{r_2}^{n-2, I}} \times [5/8, 2] \times S_0) \cap f_{0}^{-1}\big( (\overline{B_{\rho}^{n-1}} \times S_0) \cup (\overline{B_{r_2}^{n-2, I}} \times [1/3, 3/8]\times S_0)\big) = \emptyset$.
\end{enumerate}
\vskip 1mm

To prove conclusion (a), recall from the proof of Claim~\ref{claim:F-diffeo} that $f_I^{-1} \circ \phi \circ F$ is injective on $\overline{B_{r_2}^{n-2, I}}\times \overline{A_1} \times S_0$. Combining this with (D1), we get that in fact $(f_{I}^{-1}\circ \phi\circ F)(\overline{B_{r_2}^{n-2, I}} \times \overline{A_2} \times S_0)$ is disjoint from 
\[
\begin{split}
&f_{0}^{-1} (\overline{B_{\rho}^{n-1}} \times S_0) \cup (f_{I}^{-1} \circ \phi \circ F)\big(\overline{B_{r_2}^{n-2, I}} \times (1/3, 3/8)\times S_0\big) \\
&= f_{0}^{-1}\big( (\overline{B_{\rho}^{n-1}} \times S_0) \cup (\overline{B_{r_2}^{n-2, I}} \times (1/3, 3/8)\times S_0)\big),
\end{split}
\]
where for the second line we used~\eqref{eq:F-endpoints}. By the same argument, but using (D2) instead of (D1), we see that $(f_{I}^{-1}\circ \phi\circ F)(\overline{B_{r_2}^{n-2, I}} \times \overline{A_2} \times S_0)$ is also disjoint from
\[
\begin{split}
&(f_{I}^{-1} \circ \phi \circ F)\big(\overline{B^{n-2, I}_{r_2}}\times (5/8, 2/3) \times S_0\big)  \cup (f_{I}^{-1} \circ \phi)\big(\overline{B^{n-2, I}_{r_2}}\times [2/3, 2] \times S_0\big) \\
&=(f_{I}^{-1} \circ \phi)\big(\overline{B^{n-2, I}_{r_2}}\times (5/8, 2] \times S_0\big).
\end{split}
\]
Recalling also (D3), we see that the sets in~\eqref{eq:image-of-F-sets} are mutually disjoint. This proves (a). Next, by~\eqref{eq:h-disjoint-images} we have
\[
(f_{I}^{-1}\circ \phi)(\{\0^{n-2}\} \times [7/8, 2] \times S_0) \cap (f_{I}^{-1}\circ \phi)(\{\0^{n-2}\} \times [5/8, 2/3] \times S_0) = \emptyset,
\]
from which we get the assertion of (b) upon decreasing $r_2$ if needed. Finally, since $\rho > 1/3$, we get (c) by decreasing $r_2$ further. The claim is proved.\\
\end{proof}
Repeating the proof up to this point for each $I \subset\{1, \cdots, n\}$ having length $n-1$, with the intervals $A, A_0, A_1, A_2$ and the radius $\rho \in (\frac{1}{3}, \frac{3}{8})$ kept fixed, we obtain parameters $\lambda_{0, I}$ and $r_{0, I} > r_{1,I} > r_{2, I} > 0$ as above, and we fix any $\sigma > 0$ satisfying
\begin{equation}\label{eq:sigma-threshold-1}
\sigma \leq \min\{r_{2, I}\ |\ I \subset \{1, \cdots, n\},\ |I| = n-1\}.
\end{equation}
In particular, we have by Claim~\ref{claim:image-of-F}(c) that $C^{I}_{\sigma}([0, 1/3]) \subset B^{n-1}_{\rho}$ for each $|I| = n-1$. Recalling~\eqref{eq:cN-definition}, we obtain for $R \in [3/8, 2]$ that
\begin{equation}\label{eq:cN-definition-small-rho}
\cN_{\rho, \sigma, R} = B^{n-1}_{\rho} \cup \bigcup_{|I| = n-1}C_{\sigma}^{I}((1/3, R]),\quad \mathring{\cN}_{\rho, \sigma, R} = B^{n-1}_{\rho} \cup \bigcup_{|I| = n-1}C_{\sigma}^{I}((1/3, R)).
\end{equation}
Decreasing $\sigma$ if necessary, we also get
\begin{equation}\label{eq:sigma-threshold-2}
C^{I}_{2\sigma}((1/3, \infty)) \cap C^{J}_{2\sigma}([0, \infty)) = \emptyset \quad\text{whenever }I \neq J.
\end{equation}
In particular, for all $R \in [3/8, 2]$ and $|I| = n-1$, we have
\begin{equation}\label{eq:tube-intersection}
C^{I}_{\sigma}([0, \infty))\cap \cN_{\rho, \sigma, R} = C^{I}_{\sigma}([0, R]).
\end{equation}
We then define a map
\[
G: \cN_{\rho, \sigma}\times S_0 \to M
\]
in a piecewise manner as follows. First we let 
\begin{equation}\label{eq:G-in-core}
G(y, p) = f_0^{-1}(y, p),\quad\text{if }(y, p) \in B_{\rho}^{n-1} \times S_0.
\end{equation}
Then, given $I =\{i_0\}^{c}$ for some $i_0 \in \{1, \cdots, n\}$, identifying $C_{\sigma}^{I}((1/3, 2])$ with $B^{n-2, I}_{\sigma} \times (1/3, 2]$, we set
\begin{equation}\label{eq:G-on-tube}
G = \left\{
\begin{array}{ll}
f_{I}^{-1} \circ \phi_{I}\circ F_{I},  & \text{ on } B_{\sigma}^{n-2, I} \times (1/3, 2/3) \times S_0,\\
\\
f_{I}^{-1} \circ \phi_{I},  & \text{ on }B_{\sigma}^{n-2, I} \times (5/8, 2] \times S_0.
\end{array}
\right.
\end{equation}
With the help of~\eqref{eq:F-endpoints} and~\eqref{eq:sigma-threshold-2}, we see that $G$ is well-defined, and that in fact
\begin{equation}\label{eq:G-f0}
G = f_0^{-1} \quad \text{on } \mathring{\cN}_{\rho, \sigma, 3/8} \times S_0.
\end{equation}
Also, in view of the codomain of $F_I$ (see above~\eqref{eq:F-I-definition}) and the inclusion~\eqref{eq:transition-inclusion-1}, we have
\begin{equation}\label{eq:G-tube-middle-image}
G(C^{I}_{\sigma}((1/3, 2/3)) \times S_0) \subset \cU_I \cap \cU_0.
\end{equation}
We next verify that $G$ has all the desired properties.
\begin{claim}\label{claim:G-properties}
$G$ is a Lipschitz map, and the assignment $y \mapsto G(y, \cdot)$ is continuous from $\mathring{\cN}_{\rho, \sigma}$ to $C^{1}(S_0; M)$. Moreover, it has the properties asserted in the conclusions of Proposition~\ref{prop:1234-with-123}. 
\end{claim}
\begin{proof}[Proof of Claim]
Noting by~\eqref{eq:M_I*-expression} that $h_{I}([5/8, 2] \times S_0)$ is a compact subset of $M_{I}^{*}$, we infer from Remark~\ref{rmk:extension} and the second case in~\eqref{eq:G-on-tube} that $G$ is Lipschitz on each $C_{\sigma}^{\{i\}^{c}}((5/8, 2]) \times S_0$. On the other hand, extending $f_I^{-1}\circ \phi_I \circ F_I$ using the nearest-point projection onto $S_0$ defined on some tubular neighborhood of it in $\RR^{Q}$, we see, in a way similar to the end of Remark~\ref{rmk:extension}, that $f_I^{-1}\circ \phi_I \circ F_I$ is a Lipschitz map on $\overline{B_{\sigma}^{n-2, I}} \times [1/3, 2/3] \times S_0$, and hence so is $G$ on each $C_{\sigma}^{\{i\}^{c}}((1/3, 2/3)) \times S_0$. Since $G$ maps into a compact target, we deduce from the observations made thus far that, for each $i \in \{1, \cdots, n\}$,
\begin{equation}\label{eq:G-is-Lipschitz-on-tube}
G \quad\text{is Lipschitz on }C_{\sigma}^{\{i\}^{c}}((1/3, 2]) \times S_0,
\end{equation}
by considering separately pairs of points that are close together and those that are not. Next, using Remark~\ref{rmk:extension} again, this time combined with~\eqref{eq:G-f0}, we see that $G$ is Lipschitz when restricted to $\mathring{\cN}_{\rho,\sigma, 3/8} \times S_0$. From this along with~\eqref{eq:G-is-Lipschitz-on-tube}, it follows that $G$ is a Lipschitz map on its entire domain. We omit the details of this last inference, and only mention that, by the disjointness property~\eqref{eq:sigma-threshold-2}, the expression~\eqref{eq:cN-definition-small-rho}, and the inequality $\rho < 3/8$, we have 
\[
\cN_{\rho, \sigma} \setminus \mathring{\cN}_{\rho,\sigma, 3/8}  = \bigcup_{|I| = n-1}C^{I}_{\sigma}([3/8, 2]).
\]
Consequently, for any pair points in $\cN_{\rho, \sigma} \times S_0$ that are sufficiently close to each other, if one of them lies outside of $\mathring{\cN}_{\rho, \sigma, 3/8}$, then there is some $i \in \{1, \cdots, n\}$ so that both points lie in $C_{\sigma}^{\{i\}^{c}}((1/3, 2]) \times S_0$.

Next, from the definition it is clear that $G$ restricts to a diffeomorphism on each of $B_{\rho}^{n-1} \times S_0$, $B_{\sigma}^{n-2, I} \times (1/3, 2/3) \times S_0$, and $B_{\sigma}^{n-2, I} \times (5/8, 7/8) \times S_0$. In particular, we see that
\begin{equation}\label{eq:G-2/3-local-diffeo}
G|_{\mathring{\cN}_{\rho, \sigma, 7/8} \times S_0}\quad\text{is a local diffeomorphism},
\end{equation}
and that the assignment $y \mapsto G(y, \cdot)$ is continuous as a map into $C^{1}(S_0; M)$ at every point $y_0 \in \mathring{\cN}_{\rho, \sigma, 7/8}$. To extend this continuity property to all of $\mathring{\cN}_{\rho, \sigma}$, suppose $(y_k)$ is a sequence in $\mathring{\cN}_{\rho, \sigma}$ converging to some 
\[
y_0 = x_0 + t_0 \frac{\ba_{i}}{|\ba_{i}|} \in C^{\{i\}^{c}}_{\sigma}((5/8, 2)),
\]
where $i \in \{1, \cdots, n\}$, $x_0 \in B^{n-2, \{i\}^{c}}_{\sigma}$ and $t_0 \in (5/8, 2)$. Without loss of generality we can assume that $y_{k} \in C^{\{i\}^{c}}_{\sigma}((5/8, 2))$ for all $k$, so that writing 
\[
y_k = x_k + t_k\frac{\ba_{i}}{|\ba_{i}|}
\]
yields a sequence in $B^{n-2, \{i\}^{c}}_{\sigma} \times (5/8, 2)$ converging to $(x_0, t_0)$. In view of property (1) below~\eqref{eq:1234-parametrize-collar}, as well as the smooth extension of $f_I^{-1}$ provided by Remark~\ref{rmk:extension}, and the compactness of $h_{I}([5/8, 2] \times S_0)$, we conclude that
\[
G(y_{k}, \cdot) = f_{I}^{-1}(x_k, h_{I}(t_k, \cdot)) \longrightarrow f_{I}^{-1}(x_0, h_{I}(t_0, \cdot)) = G(y_0, \cdot),
\]
in $C^{1}(S_0; M)$ as $i \to \infty$. Thus we have shown that $y \mapsto G(y, \cdot)$ is continuous on each $C^{\{i\}^{c}}_{\sigma}((5/8, 2))$, and hence on all of $\mathring{\cN}_{\rho, \sigma}$ by what we observed right after~\eqref{eq:G-2/3-local-diffeo}.

Towards proving Proposition~\ref{prop:1234-with-123}(a), note that, for each $i \in \{1, \cdots, n\}$, since $f_{\{i\}^{c}}^{-1}$ takes values in $\cU_{\{i\}^{c}}$, we have from the definition of $G$ that 
\begin{equation}\label{eq:G-tube-image-UI}
G(C_{\sigma}^{\{i\}^{c}}((1/3, 2]) \times S_0) \subset \cU_{\{i\}^{c}}.
\end{equation}
The sets in $\{\cU_{\{i\}^{c}}\ |\ i = 1, \cdots, n\}$ being pairwise disjoint by~\eqref{eq:U-disjoint}, we infer from Claim~\ref{claim:image-of-F}(a) that the members of the collection
\begin{equation}\label{eq:G-image-limb}
\Big\{ C_{\sigma}^{\{i\}^{c}}([3/8, 5/8]) \times S_0  ,\quad  C_{\sigma}^{\{i\}^{c}}((5/8, 2]) \times S_0 \Big\}_{i= 1, \cdots, n}
\end{equation}
have pairwise disjoint images under $G$. Next, recalling~\eqref{eq:G-f0}, and using Claim~\ref{claim:image-of-F}(a),~\eqref{eq:G-tube-image-UI}, and~\eqref{eq:U-disjoint}, we see that each set in the collection~\eqref{eq:G-image-limb} has image disjoint from that of $\mathring{\cN}_{\rho, \sigma, 3/8} \times S_0$. Since this set along with those in the collection~\eqref{eq:G-image-limb} cover $\mathring{\cN}_{\rho, \sigma, 7/8} \times S_0$, and since each set in this covering has an intersection with $\mathring{\cN}_{\rho, \sigma, 7/8} \times S_0$ on which $G$ is injective, we conclude that actually $G$ is injective on all of $\mathring{\cN}_{\rho, \sigma, 7/8} \times S_0$. Recalling~\eqref{eq:G-2/3-local-diffeo} yields the first conclusion of (a). By the above discussion, along with Claim~\ref{claim:image-of-F}(b), we get
\begin{equation}\label{eq:prelim-for-G-disjoint-images}
G(C_{\sigma}^{\{i\}^{c}}([7/8, 2]) \times S_0 ) \cap G(\cN_{\rho,\sigma, 2/3} \times S_0) = \emptyset,\quad\text{for }i = 1, \cdots, n.
\end{equation}
Given disjoint sets $A, B \subset \cN_{\rho, \sigma}$ with $A \subset \cN_{\rho, \sigma, 2/3}$, upon writing 
\[
G(B \times S_0) = G((B \cap \mathring{\cN}_{\rho, \sigma, 7/8}) \times S_0) \cup G((B \setminus \mathring{\cN}_{\rho, \sigma, 7/8}) \times S_0), 
\]
we see from~\eqref{eq:prelim-for-G-disjoint-images} and the injectivity of $G|_{\mathring{\cN}_{\rho, \sigma, 7/8} \times S_0}$ that~\eqref{eq:G-disjoint-images} holds. This proves the second conclusion of (a).

Next, given $y_0 \in \cN_{\rho, \sigma}\setminus \mathring{\cN}_{\rho, \sigma}$, there exists $i \in \{1, \cdots, n\}$ and $x_0 \in B^{n-2, \{i\}^{c}}_{\sigma}$ such that $y_0 = (x_0, 2)$ in terms of~\eqref{eq:split-coordinates}, and thus
\[
G(\{y_0\} \times S_0) = f_{\{i\}^{c}}^{-1}(\{x\} \times  h_{\{i\}^{c}}(\{2\} \times S_0)),
\]
which has finite $\cH^{1}$-measure by property (2) listed below~\eqref{eq:1234-parametrize-collar}, and Remark~\ref{rmk:extension}. This is the first assertion of part (b). For the second assertion, note that any sequence $(y_i)$ in $\mathring{\cN}_{\rho, \sigma}$ tending to $y_0$ must eventually have the form 
\[
y_k = (x_k, t_k) \in B^{n-2, \{i\}^{c}}_{\sigma} \times (5/8, 2),
\]
where $x_k \to x_0$ and $t_k \to 2^{-}$. Letting $L$ denote the $C^1$-norm of $f_{\{i\}^{c}}^{-1}$ on $\overline{B_{\sigma}^{n-2, \{i\}^{c}}} \times h_{\{i\}^{c}}([5/8, 2] \times S_0)$, we get after a straightforward computation that
\[
\begin{split}
\Area(G(y_k, \cdot)) =\ & \Area\big(f_{\{i\}^{c}}^{-1}(x_k, h_{\{i\}^{c}}(t_k, \cdot))\big)\\
\leq\ & L^2 \cdot \Area\big(h_{\{i\}^{c}}(t_k, \cdot)\big) \to 0,\quad\text{as }k \to \infty,
\end{split}
\]
where the last convergence again uses the property (2) quoted above. This finishes the proof of part (b).

Moving on to (c), by~\eqref{eq:G-f0} and the last of the properties of $h$ appearing after~\eqref{eq:1234-parametrize-collar}, we get~\eqref{eq:handle-sweepout-on-core} for $(t, p) \in [0, \frac{3}{8}) \times S_0$. On the other hand, by~\eqref{eq:G-on-tube} and~\eqref{eq:F-id-on-center}, we get~\eqref{eq:handle-sweepout-on-core} for $(t, p) \in (\frac{1}{3}, 2] \times S_0$. Having established~\eqref{eq:handle-sweepout-on-core}, we then get the second conclusion of part (c) upon recalling~\eqref{eq:f1234-center}.

For part (d), take $(y, p) \in \cN_{\rho, \sigma} \times S_0$ along with $j \in \{1, \cdots, n\}$. In the case where $y \in \mathring{\cN}_{\rho, \sigma, 3/8}$, it follows from~\eqref{eq:G-f0} and~\eqref{eq:f1234-compatible} that $G(y, p) \in M_{j}$ if and only if $y \in V^{n-1}_{j}$. Next, in the case
\[
y \in C_{\sigma}^{I}([3/8, 2])\quad \text{for some }I =\{i_0\}^{c},
\]
we identify $C_{\sigma}^{I}([3/8, 2])$ with $B_{\sigma}^{n-2, I} \times [3/8, 2]$ and write $y = (x, t)$ as in~\eqref{eq:split-coordinates}, so that 
\[
G(y, p) = 
\left\{ 
\begin{array}{ll}
(f_{I}^{-1} \circ \phi_{I} \circ F_{I})(x, t, p),& \text{ if } t \in [3/8, 5/8],\\
(f_I^{-1}\circ \phi_I)(x, t, p), & \text{ if }t \in (5/8, 2]. 
\end{array}
\right.
\]
Now if $y \in V^{n-1}_{j}$, then since $t > 0$, we have by~\eqref{eq:projection-of-sectors-2} that $j \in I$, in which case~\eqref{eq:projection-of-sectors} gives $x \in V^{n-2, I}_{j}$, and it follows from Claim~\ref{claim:F-I-multisection}, or respectively~\eqref{eq:f123-compatible}, depending on whether $t \leq 5/8$ or $t > 5/8$, that $G(y, p) \in M_{j}$. Conversely, assuming $G(y, p) \in M_{j}$, then~\eqref{eq:G-tube-image-UI} and~\eqref{eq:U-disjoint} force $j \in I$. Again using Claim~\ref{claim:F-I-multisection} and~\eqref{eq:f123-compatible}, we get $x \in V^{n-2, I}_{j}$, and hence $y \in V^{n-1}_{j}$ by~\eqref{eq:projection-of-sectors}. This finishes the proof.
\end{proof}

In view of this last claim, the proof of Proposition~\ref{prop:1234-with-123} is complete.\\
\end{proof}
\subsection{Patching across different strata (II)}\label{subsec:different-strata-II}
We henceforth abbreviate $\cN_{\rho, \sigma, R}$ as $\cN_{R}$, dropping the subscripts $\rho$ and $\sigma$ since they have been fixed. When $R =2$, we drop it as well. The same applies to $\mathring{\cN}_{\rho, \sigma, R}$. To move towards the next step of the construction, pick any distinct $i_0, i_1 \in \{1, \dots, n\}$ and let
\[
I = \{i_0, i_1\},\quad J = \{1, \cdots, n\}\setminus \{i_0, i_1\},\quad  W = \Span\{\ba_{i_0}, \ba_{i_1}\}
\]
in the considerations of Section~\ref{subsec:smoothing}. Then the sets defined in~\eqref{eq:V-J-boundary-definition} and~\eqref{eq:V-J-corner-definition} reduce to
\begin{equation}\label{eq:strata-special-case}
\partial V^{n-1}_{J} = V^{n-1}_{J \cup\{i_0\} } \cup V^{n-1}_{J \cup \{i_1\}},\quad \partial_2 V^{n-1}_{J} = \{\0^{n-1}\}.
\end{equation}
With $C_0(n, |I|) > 2$ and $\theta(n, |I|) \in (0, \frac{\pi}{2})$ being given respectively by Lemma~\ref{lemm:foliation-from-flow} and~\eqref{eq:slope-of-sector}, we choose $\alpha > 0$ so that
\begin{equation}\label{eq:alpha-choice}
\alpha < \frac{\min\{\rho, \sigma\}}{8(C_0 + 1)}\sin\theta,
\end{equation}
and set
\begin{equation}\label{eq:rho1-choice}
\rho_1 : = \frac{C_0\alpha}{\sin\theta}\ \in (2\alpha, \frac{\min\{\rho, \sigma\}}{8}).
\end{equation}
With $\alpha$ thus chosen, we let $\Omega$, $\xi$, $\Psi$, and $\boldsymbol{\varphi}$ be the objects given by Proposition~\ref{prop:smoothing-flow}, where for now we do not emphasize their $J$-dependence in the notation. To reiterate some of their properties, the set $\Omega \subset V^{n-1}_{J}$ and its relative boundary in $W$ are described by 
\begin{align}
\Omega =\ & \{x + t\bu\ |\ x \in V^{1, I},\ t \geq \boldsymbol{\varphi}(0, x)\},\nonumber\\
\partial\Omega =\ & \{x + \boldsymbol{\varphi}(0, x)\bu\ |\ x \in V^{1, I}\},\label{eq:Omega-definition}
\end{align}
where $\bu$, defined at the start of Section~\ref{subsec:smoothing}, is in the present setting the unit vector in the direction of $\ba_{i_0} + \ba_{i_1}$. By our choice~\eqref{eq:rho1-choice} of $\rho_1$, we have
\begin{equation}\label{eq:smoothing-outside-B-rho}
\partial\Omega \setminus B^{n-1}_{\rho_1} = \partial V^{n-1}_{J} \setminus B^{n-1}_{\rho_1}, \quad \Omega \setminus B_{\rho_1}^{n-1} = V^{n-1}_{J} \setminus B_{\rho_1}^{n-1}.
\end{equation}
Next, the vector field $\xi:\partial\Omega \to W$ satisfies for all $y \in \partial\Omega$ that
\begin{equation}\label{eq:gauss-image}
\xi(y) \in \mathring{V}^{n-1}_{J}\quad \text{and}\quad |\xi(y)| \leq 1.
\end{equation}
In addition, for $\lambda \in \{0, 1\}$, on the component $V^{n-1}_{J \cup \{i_\lambda\}} \setminus B_{\rho_1}^{n-1}$ of $\partial V^{n-1}_{J} \setminus B^{n-1}_{\rho_1}$, we have
\begin{equation}\label{eq:normal-on-flat-part}
\xi = \frac{\ba_{i_\lambda; J \cup \{i_\lambda\}}}{|\ba_{i_\lambda; J \cup \{i_\lambda\}}|} = \frac{(-\bb_{i_{\lambda}; I}, 1)}{\sqrt{1 + |\bb_{ i_{\lambda}; I}|^2}},
\end{equation}
where the second equality follows from~\eqref{eq:splitting-vertex-consequence}. The map $\Psi$ is a diffeomorphism from $(-\infty, \alpha] \times \partial \Omega$ into $W$, and is given in terms of $\xi$ by the formula
\begin{equation}\label{eq:3-strata-smoothing-Psi}
\Psi(s, y) = y + s\xi(y).
\end{equation}
We note also the following orthogonal decompositions, again for $\lambda \in \{0, 1\}$. The first two come from~\eqref{eq:decomposition-wrt-lower-simplex}, while the last one uses in addition~\eqref{eq:inward-pointing}:
\begin{equation}\label{eq:decomposition-simplex-special-case}
\begin{split}
V^{n-1} =\ & V^{n-2, J\cup \{i_{\lambda}\}} \oplus \Span\{\ba_{i_{1-\lambda}}\}\\
=\ &  V^{n-3, J} \oplus \Span\{\ba_{i_{\lambda}; J \cup \{i_{\lambda}\}}\} \oplus \Span\{\ba_{i_{1-\lambda}}\}\\
=\ & V^{n-3, J} \oplus W.
\end{split}
\end{equation}

Below, we construct what may be considered the result of taking $M_{J}$ and rounding out the corners on its boundary. Specifically, define
\begin{equation}\label{eq:4-strata-smoothing}
\widetilde{M}_{J} = \big( M_{J} \setminus G(\overline{B_{\rho_1}} \times S_0) \big) \cup G((\Omega \cap B_{3\rho_1}) \times S_0),
\end{equation}
where by $B_{r}$ we mean $B^{n-1}_{r}$. Note that $\widetilde{M}_{J} \subset M_{J}$, which can be seen from the inclusion $\Omega \subset V^{n-1}_{J}$ and Proposition~\ref{prop:1234-with-123}(d). Moreover, $\widetilde{M}_{J}$ is compact. To see this, we use Proposition~\ref{prop:1234-with-123}(d), the injectivity of $G$ on $B_{3\rho_1} \times S_0$, and~\eqref{eq:smoothing-outside-B-rho} to get that
\begin{equation}\label{eq:difference-from-smoothing}
\begin{split}
M_J \setminus \widetilde{M}_{J} =\ & \big( M_{J} \cap G(\overline{B_{\rho_1}} \times S_0) \big)\setminus G( (\Omega \cap B_{3\rho_1}) \times S_0)\\
=\ & G(((V^{n-1}_{J}\setminus \Omega) \cap \overline{B_{\rho_1}}) \times S_0)\\
=\ & G(((V^{n-1}_{J}\setminus \Omega) \cap B_{\rho_1}) \times S_0) = M_{J} \cap G((B_{\rho_1} \setminus \Omega) \times S_0).
\end{split}
\end{equation}
Since $\Omega$ is closed relative to $W$, and hence relative to $V^{n-1}$ as well, the above proves $M_{J} \setminus \widetilde{M}_{J}$ to be open in $M_{J}$. The latter being compact, we conclude that so is $\widetilde{M}_{J}$.

\begin{lemm}\label{lemm:smoothed-strata}
$\widetilde{M}_{J}$ is an embedded, smooth, $4$-dimensional submanifold with boundary. Moreover, we have:
\vskip 1mm
\begin{enumerate}
\item[(a)] $\widetilde{M}_{J} \setminus G(\overline{B_{\rho_1}} \times S_0) = M_{J} \setminus G(\overline{B_{\rho_1}} \times S_0)$. 
\vskip 1mm
\item[(b)] $\widetilde{M}_{J} \cap G(B_{3\rho_1} \times S_0) = G((\Omega \cap B_{3\rho_1}) \times S_0)$.
\vskip 1mm
\item[(c)] Both (a) and (b) continue to hold with $\widetilde{M}_{J}$, $M_{J}$ and $\Omega$ replaced by $\partial\widetilde{M}_{J}$, $\partial M_{J}$ and $\partial\Omega$.
\end{enumerate}
\end{lemm}
\begin{proof}
Again by the injectivity of $G|_{B_{3\rho_1} \times S_0}$, the relation~\eqref{eq:smoothing-outside-B-rho}, and Proposition~\ref{prop:1234-with-123}(d), we have 
\begin{equation}\label{eq:smoothing-pieces-argument}
\begin{split}
G((\Omega \cap B_{3\rho_1}) \times S_0) \setminus G(\overline{B_{\rho_1}} \times S_0)=\ & G( (\Omega\cap (B_{3\rho_1}\setminus\overline{B_{\rho_1}})) \times S_0) \\
=\ & G( (V^{n-1}_{J}\cap (B_{3\rho_1}\setminus\overline{B_{\rho_1}})) \times S_0) \\
=\ & M_{J} \cap G((B_{3\rho_1} \setminus \overline{B_{\rho_1}}) \times S_0)\\
=\ & M_{J} \cap \big[ G(B_{3\rho_1}\times S_0) \setminus G(\overline{B_{\rho_1}} \times S_0) \big],
\end{split}
\end{equation}
from which we deduce both (a) and (b). Now consider the open covering of $M$ formed by
\[
\cO_{1} := M \setminus G(\overline{B_{\rho_1}} \times S_0), \quad \cO_{2} := G(B_{3\rho_1} \times S_0).
\]
Since $G|_{B_{3\rho_{1}} \times S_0}$ is a diffeomorphism, we get from part (b) that $\widetilde{M}_{J} \cap \cO_{2}$ is a smooth $4$-submanifold of $\cO_{2}$ with boundary given by $G((\partial\Omega \cap B_{3\rho_1}) \times S_0)$. On the other hand, since $S_0 \subset G(\overline{B_{\rho_1}} \times S_0)$ thanks to Proposition~\ref{prop:1234-with-123}(c), we have by part (a) that 
\[
\widetilde{M}_{J} \cap \cO_{1} = \big(\mathring{M}_{J} \cup \mathring{M}_{J \cup \{i_0\}} \cup \mathring{M}_{J \cup \{i_1\}}\big) \setminus G(\overline{B_{\rho_1}} \times S_0),
\]
which is a smooth $4$-submanifold of $\cO_{1}$ with boundary given by $\partial M_{J} \setminus G(\overline{B_{\rho_1}} \times S_0)$. Since $\cO_{1}$ and $\cO_{2}$ form an open covering of $M$, we conclude from Lemma~\ref{lemm:submanifold-construction} that $\widetilde{M}_{J}$ is a smooth $4$-submanifold of $M$, and that its set of boundary points indeed satisfies the analogues of (a) and (b), as asserted in (c). The proof is complete.
\end{proof}
\begin{prop}\label{prop:smoothed-strata-portions}
$\widetilde{M}_{J}$ has the following additional properties.
\begin{enumerate}
\item[(a)] For all $R \in [3/8, 2]$, we have 
\[
\partial\widetilde{M}_{J} \setminus G(\cN_{R} \times S_0) =  \cup_{\lambda \in \{0, 1\}}h_{J \cup \{i_{\lambda}\}}((R, 2] \times S_0) = G((\partial\Omega \cap \cN\setminus \cN_{R}) \times S_0),
\]
and the same holds with $\cN_{R}$ and $(R, 2]$ replaced by $\mathring{\cN}_{R}$ and $[R, 2]$, respectively. 
\vskip 1mm
\item[(b)] Let $E$ be any subset of $\cN$. Then 
\[
\partial\widetilde{M}_{J} \cap G(E \times S_0) = G((\partial\Omega \cap E) \times S_0).
\]
Moreover, this continues to hold with $\partial\widetilde{M}_{J}$ and $\partial\Omega$ replaced by $\widetilde{M}_{J}$ and $\Omega$, respectively.
\end{enumerate}
\end{prop}
\begin{proof}
For part (a) we only prove the first assertion, since the version for $\mathring{\cN}_R$ uses the same argument. To start, note that since $\overline{B_{\rho_1}} \subset \cN_{R}$, Lemma~\ref{lemm:smoothed-strata} gives
\begin{align*}
\partial \widetilde{M}_{J}\setminus G(\cN_{R} \times S_0) =\ & \partial M_{J}\setminus G(\cN_{R} \times S_0)\\
=\ & \cup_{\lambda \in \{0, 1\}} M_{J \cup \{i_{\lambda}\}}\setminus G(\cN_{R} \times S_0).
\end{align*}
By parts (d) and (c) of Proposition~\ref{prop:1234-with-123}, as well as~\eqref{eq:tube-intersection}, we have for $\lambda \in \{0, 1\}$ that
\[
M_{J \cup \{i_{\lambda}\}} \cap G(\cN_{R} \times S_0) = h_{J \cup \{i_{\lambda}\}}([0, R] \times S_0).
\]
Combining the previous two observations with $M_{J \cup \{i_{\lambda}\}} = h_{J \cup \{i_{\lambda}\}}([0, 2] \times S_0)$, and recalling~\eqref{eq:h-disjoint-images}, we get the first equality in part (a). By Proposition~\ref{prop:1234-with-123}(c) and~\eqref{eq:tube-intersection} again, we find that
\[
\cup_{\lambda \in \{0, 1\}}h_{J \cup \{i_{\lambda}\}}((R, 2] \times S_0) = G( ( \partial V^{n-1}_{J} \cap \cN \setminus \cN_{R} ) \times S_0),
\]
which along with~\eqref{eq:smoothing-outside-B-rho} gives the second equality in part (a).

For part (b), since 
\[
G((E \setminus \overline{B_{\rho_1}}) \times S_0) \cap G(\overline{B_{\rho_1}} \times S_0) = \emptyset
\]
by Proposition~\ref{prop:1234-with-123}(a), we deduce from Lemma~\ref{lemm:smoothed-strata}, Proposition~\ref{prop:1234-with-123}(d), and~\eqref{eq:smoothing-outside-B-rho} that
\[
\begin{split}
\partial\widetilde{M}_{J} \cap G((E \setminus \overline{B_{\rho_1}}) \times S_0) =\ & \partial M_J \cap G((E \setminus \overline{B_{\rho_1}}) \times S_0)\\
= \ & G((\partial \Omega \cap (E\setminus \overline{B_{\rho_1}}) ) \times S_0).
\end{split}
\]
Using Lemma~\ref{lemm:smoothed-strata} again, this time to see that 
\[
\partial\widetilde{M}_{J} \cap G((E \cap B_{3\rho_1}) \times S_0) = G((\partial \Omega \cap  E \cap B_{3\rho_1})\times S_0),
\]
we get the first conclusion of part (b). The version stated in the second conclusion follows from the same argument, and we omit its proof.
\end{proof}
\begin{rmk}\label{rmk:portions}
We observe the following.
\vskip 1mm
\begin{enumerate}
\item[(1)] The intersection $\partial\Omega \cap \cN_{R}$ is compact for all $R \in [3/8, 2]$, since we have by~\eqref{eq:smoothing-outside-B-rho} that
\begin{equation}\label{eq:N-R-intersect-Omega}
\begin{split}
\partial\Omega \cap \cN_{R} =\ & (\partial\Omega \cap \overline{B_{2\rho_{1}}}) \cup (\partial V^{n-1}_{J} \cap (\cN_{R} \setminus B_{\rho_{1}}))\\
=\ & (\partial\Omega \cap \overline{B_{2\rho_{1}}}) \cup \big( \cup_{\lambda \in \{0, 1\}}  [\rho_{1}, R] \cdot \frac{\ba_{i_{1-\lambda}}}{|\ba_{i_{1-\lambda}}|} \big).
\end{split}
\end{equation}
\vskip 1mm
\item[(2)] Parts (a) and (b) of Proposition~\ref{prop:smoothed-strata-portions} together imply that $\partial\widetilde{M}_{J}$ can be expressed as
\begin{equation}\label{eq:smoothing-contained-in-G-image}
\partial\widetilde{M}_{J} = G((\partial\Omega \cap \cN) \times S_0).
\end{equation}
Thus, given $R \in [3/8, 2/3]$, using also Proposition~\ref{prop:1234-with-123}(a), we find that
\[
G(( \partial\Omega \cap \cN \setminus \cN_{R} ) \times S_0) = \partial\widetilde{M}_{J} \setminus G((\partial\Omega \cap \cN_{R}) \times S_0),
\]
which together with item (1) above shows that the left-hand side is open relative to $\partial\widetilde{M}_{J}$. Noting in addition that, by~\eqref{eq:M_I*-expression}, the sets $h_{J\cup \{i_{\lambda}\}}((R, 2] \times S_0)$ for $\lambda = 0, 1$ have disjoint compact closures, we conclude that each is relatively open in $\partial\widetilde{M}_{J}$.
\vskip 1mm
\item[(3)] For all $R \in [3/8, 7/8]$, upon taking $E = \mathring{\cN}_{R}$ in Proposition~\ref{prop:smoothed-strata-portions}(b) and recalling Proposition~\ref{prop:1234-with-123}(a), we get that $G|_{(\partial\Omega \cap \mathring{\cN}_{R}) \times S_0}$ is a diffeomorphism onto $\partial\widetilde{M}_{J} \cap G(\mathring{\cN}_{R} \times S_0)$. 
\end{enumerate}
\end{rmk}

As another preliminary result, we make a simple observation about the map $\Psi$ from~\eqref{eq:3-strata-smoothing-Psi}.
\begin{lemm}\label{lemm:preservation-by-flow}
Given $x \in B^{n-3, J}_{\alpha}$, $s \in (-\alpha, \alpha)$, and $y \in \partial\Omega$, we have for any $R \in [3/8, 2]$ that
\begin{equation}\label{eq:Psi-image-intersect-N}
y \in \cN_{R} \Longleftrightarrow x + \Psi(s, y) \in \cN_{R}.
\end{equation}
A similar equivalence holds with $\cN_{R}$ replaced by $\mathring{\cN}_{R}$ on both sides.
\end{lemm}
\begin{proof}
By~\eqref{eq:rho1-choice} we have $B_{8\rho_1}\subset \cN_{R}$, and that $\alpha < \frac{\rho_1}{2}$. Taking also~\eqref{eq:gauss-image} into account, we deduce that, with $(x, s, y)$ as given, if either $y\not\in \cN_{R}$ or $x + \Psi(s, y)\not\in\cN_{R}$, then necessarily $y \in \partial\Omega \setminus B_{7\rho_1}$, so that by~\eqref{eq:smoothing-outside-B-rho} and~\eqref{eq:normal-on-flat-part} we obtain some  $\lambda \in \{0, 1\}$ and $t \geq 7\rho_1$ such that 
\[
y = t\frac{\ba_{i_{1-\lambda}}}{|\ba_{i_{1-\lambda}}|},\quad\quad x + \Psi(s, y) = x + s\frac{\ba_{i_{\lambda}; J \cup \{i_{\lambda}\}}}{|\ba_{i_{\lambda}; J \cup \{i_{\lambda}\}}|} + t\frac{\ba_{i_{1-\lambda}}}{|\ba_{i_{1-\lambda}}|}.
\]
In view of~\eqref{eq:decomposition-simplex-special-case}, and using the notation~\eqref{eq:cylinder-definition}, we see that $y$ and $x + \Psi(s, y)$ both lie in $C^{J \cup \{i_{\lambda}\}}_{2\alpha}(\{t\})$. By~\eqref{eq:tube-intersection}, since at least one of these two points is assumed to lie outside of $\cN_{R}$, we must have $t > R$, in which case, again by~\eqref{eq:tube-intersection}, both points are outside $\cN_{R}$. This proves the claimed equivalence. The argument is the same when $\cN_{R}$ is replaced by $\mathring{\cN}_{R}$ in~\eqref{eq:Psi-image-intersect-N}, and we omit the details.
\end{proof}

Next, we use Proposition~\ref{prop:smoothed-strata-portions} to write
\begin{equation}\label{eq:boundary-expression-for-g-J}
\begin{split}
\partial\widetilde{M}_{J} = \ & (\partial\widetilde{M}_{J} \cap G(\mathring{\cN}_{2/3} \times S_0)) \cup (\partial\widetilde{M}_{J} \setminus G(\cN_{5/8} \times S_0))\\
=\ & G((\mathring{\cN}_{2/3} \cap \partial\Omega ) \times S_0) \cup \big( \cup_{\lambda \in \{0, 1\}} h_{J \cup\{i_{\lambda}\}}((5/8, 2] \times S_0) \big).
\end{split}
\end{equation}
Notice that $G((\mathring{\cN}_{2/3} \cap \partial\Omega ) \times S_0)$ is open relative to $\partial\widetilde{M}_{J}$, and that, by Remark~\ref{rmk:portions}, so is each $h_{J \cup\{i_{\lambda}\}}((5/8, 2] \times S_0) \big)$. With $\alpha$ chosen as in~\eqref{eq:alpha-choice}, we define 
\[
g_J: B^{n-3, J}_{\alpha} \times (-\alpha, \alpha) \times \partial \widetilde{M}_{J} \to M
\]
by considering the following cases:
\vskip 1mm
\begin{enumerate}
\item[(1)] If $q \in G((\partial\Omega \cap \mathring{\cN}_{2/3}) \times S_0)$, we let $(y, p) = \big(G|_{\mathring{\cN}_{7/8} \times S_0}\big)^{-1}(q)$, and define
\begin{equation}\label{eq:g-J-definition-case-1}
g_{J}(x, s, q) = G(x + \Psi(s, y), p).
\end{equation}
\vskip 1mm
\item[(2)] If $q \in h_{J \cup \{i_{\lambda}\}}((5/8, 2] \times S_0)$ for some $\lambda \in \{0, 1\}$, necessarily unique by~\eqref{eq:M_I*-expression} and~\eqref{eq:U-disjoint}, we set
\begin{equation}\label{eq:g-J-definition-case-2}
g_{J}(x, s, q) = f_{J\cup\{i_\lambda\}}^{-1}(x+ s\frac{\ba_{i_{\lambda}; J \cup \{i_\lambda\}}}{| \ba_{i_{\lambda}; J \cup \{i_\lambda\}}|}, q).
\end{equation}
\end{enumerate}
By the first two decompositions in~\eqref{eq:decomposition-simplex-special-case} and, respectively, Lemma~\ref{lemm:preservation-by-flow}, we see that the maps on the right-hand side of~\eqref{eq:g-J-definition-case-2} and~\eqref{eq:g-J-definition-case-1} are indeed evaluated at points in their domains. That the two definitions agree on the overlap region will be addressed shortly. For brevity, in what follows we let
\begin{equation}\label{eq:v-lambda-definition}
\bv_{\lambda}: = \frac{\ba_{i_{\lambda}; J \cup \{i_\lambda\}}}{| \ba_{i_{\lambda}; J \cup \{i_\lambda\}}|}.
\end{equation}
\begin{prop}\label{prop:g-J-properties}
$g_{J}$ is a well-defined smooth map. Moreover the following hold.
\begin{enumerate}
\item[(a)] For all $(x, s) \in B^{n-3, J}_{\alpha} \times (-\alpha, \alpha)$ and $(y, p) \in (\partial\Omega \cap \cN) \times S_0$, we have 
\[
g_J(x, s, G(y, p)) = G(x + \Psi(s, y), p).
\]
\vskip 1mm
\item[(b)] $g_J(\0^{n-3}, 0, q) = q$ for all $q \in \partial\widetilde{M}_{J}$. Given in addition $(x, s) \in B^{n-3, J}_{\alpha} \times (-\alpha, \alpha)$, we have
\begin{equation}\label{eq:into-smoothing}
g_{J}(x, s, q) \in \widetilde{M}_{J} \quad\text{if and only if}\quad x  = \0^{n-3} \text{ and }s \geq 0.
\end{equation}
\vskip 1mm
\item[(c)] $g_J$ takes $B^{n-3, J}_{\alpha} \times (-\alpha, \alpha) \times \partial \widetilde{M}_{J}$ diffeomorphically onto a neighborhood of $\partial \widetilde{M}_{J}$ in $M$.
\vskip 1mm
\item[(d)] Given $(x, s, q) \in B^{n-3, J}_{\alpha} \times [0, \alpha) \times \partial \widetilde{M}_{J}$ and $j \in J$, we have
\begin{equation}\label{eq:g-J-multisection}
g_{J}(x, s, q) \in M_{j} \quad\text{ if and only if }\quad x \in V^{n-3, J}_{j}.
\end{equation}
\end{enumerate}
\end{prop}
\begin{proof}
We begin by observing that, given $(x, s) \in B^{n-3, J}_{\alpha} \times (-\alpha, \alpha)$, $(t, p) \in (5/8, 2] \times S_0$, and $\lambda \in \{0, 1\}$, as in the proof of Lemma~\ref{lemm:preservation-by-flow}, we have by~\eqref{eq:normal-on-flat-part} that
\[
\Psi(s, t\frac{\ba_{i_{1-\lambda}}}{|\ba_{i_{1-\lambda}}|}) = s\bv_{\lambda} + t\frac{\ba_{i_{1-\lambda}}}{|\ba_{i_{1-\lambda}}|},
\]
and hence, by the definition of $G$, more precisely the second case in~\eqref{eq:G-on-tube}, there holds
\begin{equation}\label{eq:g-J-well-defined}
f_{J \cup \{i_{\lambda}\}}^{-1}(x + s\bv_{\lambda}, h_{J \cup \{i_{\lambda}\}}(t, p)) = G(x + \Psi(s, t\frac{\ba_{i_{1-\lambda}}}{|\ba_{i_{1-\lambda}}|}), p).
\end{equation}
To prove that $g_{J}$ is well-defined, we need only check that~\eqref{eq:g-J-definition-case-1} agrees with~\eqref{eq:g-J-definition-case-2} when $q$ lies in $G((\partial\Omega \cap \mathring{\cN}_{2/3}) \times S_0) \cap  h_{J \cup \{i_{\lambda}\}}((5/8, 2] \times S_0)$ for some $\lambda \in \{0, 1\}$. Noting from Proposition~\ref{prop:1234-with-123}(c) and~\eqref{eq:tube-intersection} that
\[
h_{J \cup \{i_{\lambda}\}}((5/8, 2] \times S_0) = G( ( V^{n-1}_{J \cup \{i_{\lambda}\}} \cap \cN \setminus \cN_{5/8}) \times S_0),
\]
we have by Proposition~\ref{prop:1234-with-123}(a) and~\eqref{eq:smoothing-outside-B-rho} that
\[
G((\partial\Omega \cap \mathring{\cN}_{2/3}) \times S_0) \cap  h_{J \cup \{i_{\lambda}\}}((5/8, 2] \times S_0) = G((V^{n-1}_{J \cup \{i_{\lambda}\}} \cap \mathring{\cN}_{2/3} \setminus \cN_{5/8}) \times S_0).
\]
As a result, any $q$ belonging to the left-hand side above can be expressed as
\[
q = G(t\frac{\ba_{i_{1-\lambda}}}{|\ba_{i_{1-\lambda}}|}, p) = h_{J \cup \{i_{\lambda}\}}(t, p),
\]
for some $(t, p) \in (5/8, 2/3) \times S_0$. The value assigned to $g_{J}(x, s, q)$ according to~\eqref{eq:g-J-definition-case-2} then coincides with the left-hand side of~\eqref{eq:g-J-well-defined}, while since
\[
\big(G|_{\mathring{\cN}_{7/8} \times S_0}\big)^{-1}(q) = (t\frac{\ba_{i_{1-\lambda}}}{|\ba_{i_{1-\lambda}}|}, p),
\]
we see that the right-hand side of~\eqref{eq:g-J-well-defined} is exactly $g_{J}(x, s, q)$ defined according to~\eqref{eq:g-J-definition-case-1}. Thus we conclude that $g_J$ is well-defined.

For part (a), let $(x, s, y, p)$ be as in the statement. In the case where $y \in \partial\Omega \cap \mathring{\cN}_{2/3}$, the result follows directly from~\eqref{eq:g-J-definition-case-1}. In the case $y \in \partial\Omega \cap (\cN \setminus \cN_{5/8})$, by~\eqref{eq:smoothing-outside-B-rho} we can write
\[
y = t\frac{\ba_{i_{1-\lambda}}}{|\ba_{i_{1-\lambda}}|},\quad\text{for some }t \in (5/8, 2]\text{ and }\lambda \in \{0, 1\},
\]
in which case $G(y, p) = h_{J \cup \{i_{\lambda}\}}(t, p)$, and we get the desired equality from~\eqref{eq:g-J-definition-case-2} and~\eqref{eq:g-J-well-defined}.

For the first conclusion in part (b), we merely mention that~\eqref{eq:f123-center} is used when $q$ falls into case (2) in the definition of $g_J$. For the second conclusion, we take $x \in B^{n-3, J}_{\alpha}$, $s \in (-\alpha,\alpha)$ and $q \in \partial\widetilde{M}_{J}$, and consider separately the following two cases:
\vskip 1mm
\begin{enumerate}
\item[(i)] If $q = G(y, p)$ for some $(y, p) \in (\partial\Omega \cap \mathring{\cN}_{2/3}) \times S_0$, we have by Lemma~\ref{lemm:preservation-by-flow} that 
\begin{equation}\label{eq:preservation-by-flow-applied-1}
x + \Psi(s, y) \in \mathring{\cN}_{2/3}.
\end{equation}
Using also part (a), we see that in this case $g_{J}(x, s, q) \in \widetilde{M}_{J}$ if and only if
\[
G(x + \Psi(s, y), p) \in \widetilde{M}_{J} \cap G(\mathring{\cN}_{2/3} \times S_0) = G(( \Omega \cap \mathring{\cN}_{2/3})\times S_0),
\]
where the equality follows from Proposition~\ref{prop:smoothed-strata-portions}(b). Since $G$ is injective on $\mathring{\cN}_{7/8} \times S_0$, and since we have~\eqref{eq:preservation-by-flow-applied-1} in the background, the above in turn is equivalent to 
\[
x + \Psi(s, y) \in \Omega.
\]
This certainly holds when $x = \0^{n-3}$ and $s \geq 0$, thanks to the first part of Remark~\ref{rmk:Psi-distance-bounds}. Conversely, if $x + \Psi(s, y) \in \Omega$, then the inclusion $\Omega \subset W$ and the third decomposition in~\eqref{eq:decomposition-simplex-special-case} forces $x = \0^{n-3}$, in which case Remark~\ref{rmk:Psi-distance-bounds} gives $s \geq 0$.
\vskip 1mm
\item[(ii)] If $q = h_{J \cup \{i_{\lambda}\}}(t, p)$ for some $(t, p) \in (5/8, 2] \times S_0$ and $\lambda \in \{0, 1\}$, then Proposition~\ref{prop:1234-with-123}(c) gives $q =  G(t\frac{\ba_{i_{1-\lambda}}}{|\ba_{i_{1-\lambda}}|}, p)$. Since, by~\eqref{eq:smoothing-outside-B-rho},
\[
t\frac{\ba_{i_{1-\lambda}}}{|\ba_{i_{1-\lambda}}|} \in \partial V^{n-1}_{J} \cap (\cN \setminus \cN_{5/8}) = \partial \Omega \cap (\cN \setminus \cN_{5/8}),
\]
we deduce from part (a), Lemma~\ref{lemm:preservation-by-flow}, and Proposition~\ref{prop:1234-with-123}(a) that
\[
g_J(x, s, q) \in G((\cN \setminus \cN_{5/8}) \times S_0) = G(\cN \times S_0) \setminus G(\cN_{5/8} \times S_0).
\]
Thus, using also Lemma~\ref{lemm:smoothed-strata}, we see that in this case $g_{J}(x,s, q)$ lies in $\widetilde{M}_{J}$ if and only if it lies in $M_J$. Since $g_J(x, s, q)$ in this case can be computed using~\eqref{eq:g-J-definition-case-2}, this last statement further reduces to
\[
f_{J \cup \{i_{\lambda}\}}^{-1}(x + s\bv_{\lambda}, q) \in M_J.
\]
By~\eqref{eq:f123-compatible}, this occurs exactly when 
\[
x + s\bv_{\lambda} \in V^{n-2, J \cup \{i_{\lambda}\}}_{J} = [0, \infty) \cdot \bv_{\lambda},
\]
where the equality uses~\eqref{eq:sectors} and~\eqref{eq:v-lambda-definition}. Upon comparing the first two decompositions in~\eqref{eq:decomposition-simplex-special-case} we see that the above is equivalent to $x = \0^{n-3}$ and $s \geq 0$.
\end{enumerate}

For part (c), we define 
\[
\begin{split}
E_1 :=\ & B_{\alpha}^{n-3, J} \times (-\alpha, \alpha) \times G( (\partial\Omega \cap \mathring{\cN}_{2/3}) \times S_0),\\
E_{2, \lambda} :=\ & B_{\alpha}^{n-3, J} \times (-\alpha, \alpha)  \times h_{J \cup \{i_{\lambda}\}}((5/8, 2] \times S_0),\\
E_2 : =\ & E_{2, 0} \cup E_{2, 1}\\
=\ & B_{\alpha}^{n-3, J} \times (-\alpha, \alpha)  \times G((\partial\Omega \cap \cN \setminus \cN_{5/8}) \times S_0),
\end{split}
\]
where the second equality for $E_2$ is a consequence of Proposition~\ref{prop:smoothed-strata-portions}. Also, by what we noticed below~\eqref{eq:boundary-expression-for-g-J}, each of $E_{1}$, $E_{2, 0}$, and $E_{2, 1}$ is open relative to $B_{\alpha}^{n-3, J} \times (-\alpha, \alpha) \times \partial \widetilde{M}_{J}$. Next, with the help of the last decomposition in~\eqref{eq:decomposition-simplex-special-case} and Lemma~\ref{lemm:preservation-by-flow}, we see that the map 
\[
(x, s, y) \mapsto x + \Psi(s,y)
\]
restricts to a diffeomorphism on $B_{\alpha}^{n-3, J} \times (-\alpha, \alpha) \times (\partial\Omega \cap \mathring{\cN}_{2/3})$ onto an open subset of $\mathring{\cN}_{2/3}$. Combining this with Remark~\ref{rmk:portions}(3) and Proposition~\ref{prop:1234-with-123}(a), we deduce that 
\[
g_J|_{E_1} \quad\text{is a diffeomorphism.}
\]
On the other hand, note that $(x, s) \mapsto x+ s\bv_{\lambda}$ defines a diffeomorphism from $B^{n-3, J}_{\alpha} \times (-\alpha, \alpha)$ into $B_{2\alpha}^{n-2, J \cup \{i_{\lambda}\}}$, and that $h_{J \cup \{i_{\lambda}\}}((5/8, 2] \times S_0)$ is open relative to $M_{J \cup \{i_{\lambda}\}}^{*}$ by~\eqref{eq:h-disjoint-images} and~\eqref{eq:M_I*-expression}. Since $f_{J \cup\{i_{\lambda}\}}^{-1}$ is a diffeomorphism, we deduce that $g_J$ restricts to a diffeomorphism on each $E_{2, \lambda}$. Moreover, in view of~\eqref{eq:U-disjoint}, we have
\[
g_J(E_{2, 0}) \cap g_J(E_{2, 1}) \subset \cU_{J \cup \{i_0\}} \cap \cU_{J \cup \{i_1\}} = \emptyset,
\]
and hence 
\[
g_J|_{E_2}\quad\text{is a diffeomorphism}.
\]
To see that $g_{J}$ is a diffeomorphism on $E_1 \cup E_2$, which equals the whole domain of $g_J$ thanks to~\eqref{eq:smoothing-contained-in-G-image}, it remains to show that $E_1 \setminus E_2$ and $E_2 \setminus E_1$ have disjoint images under $g_J$. To that end, we use Proposition~\ref{prop:1234-with-123}(a) to obtain
\[
\begin{split}
E_1 \setminus E_2 = \ &  B_{\alpha}^{n-3, J} \times (-\alpha, \alpha) \times G( (\partial\Omega \cap \cN_{5/8}) \times S_0),\\
E_2 \setminus E_1 = \ & B_{\alpha}^{n-3, J} \times (-\alpha, \alpha) \times G( (\partial\Omega \cap \cN \setminus \mathring{\cN}_{2/3}) \times S_0).
\end{split}
\]
By part (a) and Lemma~\ref{lemm:preservation-by-flow}, followed by Proposition~\ref{prop:1234-with-123}(a), we conclude
\[
g_J(E_1 \setminus E_2) \cap g_J(E_2 \setminus E_1) \subset G(\cN_{5/8} \times S_0) \cap G((\cN \setminus \mathring{\cN}_{2/3}) \times S_0) = \emptyset,
\]
and we are done with part (c).

For (d), take $(x, s, q) \in B^{n-3, J}_{\alpha} \times [0, \alpha) \times \partial \widetilde{M}_{J}$ along with some $j \in J$. When $x = \0^{n-3}$, by part (b) and the inclusion $\widetilde{M}_{J} \subset M_{J}$, we see that both ends in~\eqref{eq:g-J-multisection} are true statements, so the asserted equivalence holds. Next, when $x \neq \0^{n-3}$, which can only occur if $n > 3$, we treat separately cases (1) and (2) in the definition of $g_J$. In the case where $q = G(y, p)$ for some $(y, p) \in (\partial\Omega \cap \mathring{\cN}_{2/3}) \times S_0$, by~\eqref{eq:gauss-image} and the inclusion $\partial\Omega \subset V^{n-1}_{J}$, we have
\[
\Psi(s, y) \in V^{n-1}_{J}.
\]
Consequently, from~\eqref{eq:projection-of-sectors}, we infer that $x \in V^{n-3, J}_{j}$ if and only if $x + \Psi(s, y) \in V^{n-1}_{j}$. By Proposition~\ref{prop:1234-with-123}(d) and~\eqref{eq:g-J-definition-case-1}, this last inclusion holds exactly when $g_J(x, s, q) \in M_j$. On the other hand, if $q \in h_{J \cup\{i_\lambda\}}((5/8, 2] \times S_0)$ for some $\lambda \in\{0, 1\}$, then since $s\bv_{\lambda} \in V^{n-2, J \cup \{i_\lambda\}}_{J}$, we have by~\eqref{eq:projection-of-sectors} that $x \in V^{n-3, J}_{j}$ if and only if $x + s\bv_{\lambda} \in V^{n-2, J\cup\{i_{\lambda}\}}_{j}$, which in turn is equivalent to $g_J(x, s, q) \in M_j$ by~\eqref{eq:f123-compatible} and~\eqref{eq:g-J-definition-case-2}. This finishes the proof of part (d).
\end{proof}

Recalling that $\sigma$ is fixed to satisfy both~\eqref{eq:sigma-threshold-1} and~\eqref{eq:sigma-threshold-2}, we define, for $r \in (0, \sigma)$, the distance neighborhood
\begin{equation}\label{eq:distance-neighborhood-1-strata}
\cA_{r} : = B_{r}(\cup_{i=1}^{n}V^{n-1}_{\{i\}^{c}}) = \cup_{i= 1}^{n} C^{\{i\}^{c}}_{r}([0, \infty)),
\end{equation}
where the ``$\subset$'' part of the second equality uses the fact that, since $\sum_{i=1}^{n}\ba_{i} = \0^{n-1}$, the distance to $\cup_{i = 1}^{n}V^{n-1}_{\{i\}^{c}}$ from a point outside of it cannot be realized at the origin. Given $R \in [3/8, 2]$, upon recalling~\eqref{eq:tube-intersection}, we deduce from~\eqref{eq:distance-neighborhood-1-strata} that 
\begin{equation}\label{eq:N-A-intersection}
\cN_{R} \cap \cA_{r} = \cup_{i = 1}^{n}C^{\{i\}^{c}}_{r}([0, R]), \quad\quad (\cN \setminus \cN_{R}) \cap \cA_{r} = \cup_{i = 1}^{n}C^{\{i\}^{c}}_{r}((R, 2]).
\end{equation}
Similarly, replacing $\cA_{r}$ by its closure, and using~\eqref{eq:tube-intersection} along with relations analogous to~\eqref{eq:distance-neighborhood-1-strata}, we find that
\begin{equation}\label{eq:N-A-closure-intersection}
\begin{split}
\cN_{R} \cap \overline{\cA_{r}} =\ & \cup_{i = 1}^{n}\big( \overline{B^{n-2, \{i\}^{c}}_{r}} + [0, \infty) \cdot \frac{\ba_{i}}{|\ba_{i}|} \big) \cap \cN_{R}\\
=\ & \cup_{i=1}^{n} \big( \overline{B^{n-2, \{i\}^{c}}_{r}} + [0, R] \cdot \frac{\ba_{i}}{|\ba_{i}|} \big).
\end{split}
\end{equation}
In particular, $\cN_{R} \cap \overline{\cA_{r}}$ is a compact set.
\begin{rmk}\label{rmk:distance-attained}
Given $J \subset \{1, \cdots, n\}$ with $1 \leq |J| \leq n-2$, together with some $i \in J$, we have
\begin{equation}\label{eq:cone-half-space-disjoint}
\big( V^{n-2, \{i\}^{c}} + [0, \infty) \cdot \frac{\ba_{i}}{|\ba_{i}|} \big) \cap V^{n-1}_{J} = \{\0^{n-1}\},
\end{equation}
which can be seen by expressing points in $V^{n-1}_{J}$ in the form~\eqref{eq:sectors}, and recalling that $\ba_{i} \perp V^{n-2, \{i\}^{c}}$ (see~\eqref{eq:decomposition-wrt-lower-simplex}) and that $\ba_{k} \cdot \ba_{l} = \delta_{kl} - \frac{1}{n}$. From~\eqref{eq:cone-half-space-disjoint}, we get the first of the two equalities below:
\begin{equation}\label{eq:distance-attained}
V^{n-1}_{J} \cap \cA_{r} = \cup_{i \not\in J} V^{n-1}_{J} \cap C^{\{i\}^{c}}_{r}([0, \infty)) = V^{n-1}_J \cap B_{r}\big( \cup_{i \not\in J}V^{n-1}_{\{i\}^{c}} \big),
\end{equation}
while second equality holds for the same reason as its counterpart in~\eqref{eq:distance-neighborhood-1-strata}, except we use instead the fact that $y \cdot \sum_{i \not\in J}\ba_{i} > 0$ for all $y \in V^{n-1}_{J} \setminus \{\0^{n-1}\}$. 
\end{rmk}
Below, we specialize back to the case $|J| = n-2$. For such $J$, the right-most set in~\eqref{eq:distance-attained} reduces to $V^{n-1}_{J} \cap B_{r}(\partial V^{n-1}_{J})$.
\begin{prop}\label{prop:distance-neighborhood-in-smoothing}
Suppose $\alpha$ satisfies~\eqref{eq:alpha-choice}, and let $\tau_0 = \tau_0(n, 2)$ be the constant from Lemma~\ref{lemm:foliation-from-flow}. Given $J \subset \{1, \cdots, n\}$ with $|J| = n-2$ and writing $J^{c} = \{i_{0}, i_{1}\}$, we have the following regarding the neighborhoods $\cA_{r}$.
\vskip 1mm
\begin{enumerate}
\item[(a)] For all $r \in (0, \alpha]$, 
\[
\begin{split}
&\widetilde{M}_{J} \cap G\big( (\cN \cap \cA_{r}) \times S_0 \big)\\
&\subset g_{J}(\{\0^{n-3}\} \times [0, r) \times \partial\widetilde{M}_{J}) \subset \widetilde{M}_{J} \cap G\big( (\cN \cap \cA_{\tau_0\alpha + r\sin^2\theta}) \times S_0 \big).
\end{split}
\]
\vskip 1mm
\item[(b)] For all $r \in [\tau_0\alpha, \alpha]$, we have 
\[
\widetilde{M}_{J} \setminus G((\cN \cap \cA_{r}) \times S_0) = M_{J} \setminus G((\cN \cap \cA_{r}) \times S_0).
\]
\end{enumerate}
\end{prop}
\begin{proof}
Using Proposition~\ref{prop:smoothed-strata-portions}, followed by the inclusion $\Omega \subset V^{n-1}_{J}$ and~\eqref{eq:distance-attained}, we have
\[
\widetilde{M}_{J} \cap G((\cN \cap \cA_{r})\times S_0)
= G(( \cN \cap \Omega \cap \cA_{r}) \times S_0) = G(( \cN \cap \Omega \cap B_{r}(\partial V^{n-1}_J) ) \times S_0).
\]
Combining this with the following consequence of~\eqref{eq:distance-neighborhood-contained} and Lemma~\ref{lemm:preservation-by-flow}:
\[
\cN \cap \Omega \cap B_{r}(\partial V^{n-1}_{J}) \subset \cN \cap \Psi([0, r) \times \partial\Omega) = \Psi([0, r) \times (\cN \cap \partial\Omega)),
\]
we get
\[
\widetilde{M}_{J} \cap G((\cN \cap \cA_{r})\times S_0) \subset G(\Psi([0, r)\times (\cN \cap \partial\Omega)) \times S_0).
\]
Noting from Proposition~\ref{prop:g-J-properties}(a) and~\eqref{eq:smoothing-contained-in-G-image} that
\[
G(\Psi([0, r)\times (\cN \cap \partial\Omega)) \times S_0) = g_{J}(\{\0^{n-3}\} \times [0, r) \times \partial\widetilde{M}_{J}),
\]
we arrive at the first inclusion in part (a). By the same argument, but using~\eqref{eq:Psi-distance-estimates} instead of~\eqref{eq:distance-neighborhood-contained} to see that
\[
\Psi([0, r) \times (\cN \cap \partial\Omega)) \subset \cN \cap \Omega \cap B_{\tau_0\alpha + r\sin^2\theta}(\partial V^{n-1}_{J}),
\]
we get the second inclusion in (a). 

For part (b), in view of the definition~\eqref{eq:4-strata-smoothing} of $\widetilde{M}_{J}$, and the corresponding expression for $M_J$, it suffices to prove that 
\begin{equation}\label{eq:3-strata-agree-reduced}
\begin{split}
&G((\Omega \cap B_{3\rho_1}) \times S_0) \setminus G((\cN \cap \cA_{r}) \times S_0)\\
&= G((V^{n-1}_{J} \cap B_{3\rho_1}) \times S_0) \setminus G((\cN \cap \cA_{r}) \times S_0).
\end{split}
\end{equation}
To that end, note that, by Proposition~\ref{prop:1234-with-123}(a) and the inclusion $B_{3\rho_1} \subset \cN_{2/3}$, we have
\[
G((\Omega \cap B_{3\rho_1}) \times S_0) \setminus G((\cN \cap \cA_{r}) \times S_0) = G(( B_{3\rho_1} \cap (\Omega \setminus \cA_{r}) ) \times S_0),
\]
and similarly
\[
G((V^{n-1}_{J} \cap B_{3\rho_1}) \times S_0) \setminus G((\cN \cap \cA_{r}) \times S_0) = G(( B_{3\rho_1} \cap (V^{n-1}_{J} \setminus \cA_{r}) ) \times S_0).
\]
Since $r \in [\tau_0\alpha, \alpha]$, we obtain~\eqref{eq:3-strata-agree-reduced} upon recalling from~\eqref{eq:smoothing-agree-outside-distance-nbhd} and~\eqref{eq:distance-attained} that
\[
\Omega \setminus \cA_{r} = V^{n-1}_{J} \setminus \cA_{r}.
\]
The proof is complete.
\end{proof}

Largely for use in the next section, we mention a few additional facts involving $G$ and the distance neighborhoods $\cA_{r}$.
\begin{lemm}\label{lemm:1-strata-nbhd-facts}
Given $0 < s < r \leq 4\alpha$ and $R \in [5/8, 2/3]$, the following hold.
\vskip 1mm
\begin{enumerate}
\item[(a)] We have 
\[
G((\cA_{r} \cap \cN \setminus \cN_{R}) \times S_0) = \cup_{i = 1}^{n} f^{-1}_{\{i\}^{c}}(B^{n-2, \{i\}^{c}}_{r} \times h_{\{i\}^{c}}((R, 2] \times S_0)).
\]
Moreover, this continues to hold with both $\cA_{r}$ and $B^{n-2, \{i\}^{c}}_{r}$ replaced by their closures.
\vskip 1mm
\item[(b)] We have that
\[
G( ( (\cA_{r}\setminus \cA_{s})\cap (\cN \setminus \cN_{R}) ) \times S_0 ) \cap 
G( ( \cA_{s} \cap \cN \setminus \cN_{R} ) \times S_0 ) = \emptyset.
\]
and that, as a result,
\[
G( ( (\cA_{r}\setminus \cA_{s}) \cap \cN ) \times S_0 ) \cap 
G( (\cA_{s}  \cap \cN ) \times S_0 ) = \emptyset.
\]
Each equality remains true if we replace either $\cA_{r}$ or $\cA_{s}$, or both, by the closure.
\vskip 1mm
\item[(c)] $G((\cA_r \cap \cN) \times S_0)$ is a neighborhood of $\cup_{i=1}^{n}M_{\{i\}^{c}}$ in $M$.
\end{enumerate}
\end{lemm}
\begin{proof}
Since $R \geq 5/8$, we get the first conclusion of part (a) upon using the second relation in~\eqref{eq:N-A-intersection} and recalling~\eqref{eq:G-on-tube}. For the stated variant involving closures, we simply observe that the reasoning leading to~\eqref{eq:N-A-closure-intersection} yields
\begin{equation}\label{eq:N-A-closure-intersection-2}
(\cN \setminus \cN_{R}) \cap \overline{\cA_{r}} = \cup_{i=1}^{n} \big( \overline{B^{n-2, \{i\}^{c}}_{r}} + (R, 2] \cdot \frac{\ba_{i}}{|\ba_{i}|} \big).
\end{equation}

For part (b), by combining~\eqref{eq:N-A-intersection} with~\eqref{eq:sigma-threshold-2}, and using again~\eqref{eq:G-on-tube}, we have
\[
G( ((\cA_{r}\setminus \cA_{s}) \cap (\cN \setminus \cN_{R}) ) \times S_0 ) = \cup_{i = 1}^{n}\, f_{\{i\}^{c}}^{-1}((B_r^{n-2, \{i\}^{c}} \setminus B_{s}^{n-2, \{i\}^{c}})\times h_{\{i\}^{c}}((R, 2] \times S_0)).
\]
Since the maps $f_{\{i\}^{c}}^{-1}$ have mutually disjoint images by~\eqref{eq:U-disjoint}, and since each of them is injective, we get the first equality in (b) from part (a) and the above. The second conclusion of (b) then follows from splitting $\cN$ as $(\cN \setminus \cN_{R}) \cup \cN_{R}$, and recalling Proposition~\ref{prop:1234-with-123}(a) and the injectivity of $G|_{\mathring{\cN}_{7/8} \times S_0}$. The version involving closures can be proved by essentially the same argument, using in addition~\eqref{eq:N-A-closure-intersection-2} in the process. We omit the details.

For part (c), it follows from Proposition~\ref{prop:1234-with-123}(c) that 
\[
\cup_{i = 1}^{n}M_{\{i\}^{c}} \subset G((\cA_{r} \cap \cN) \times S_0).
\]
To see that the right-hand side is an open set, we note from~\eqref{eq:h-disjoint-images} and~\eqref{eq:M_I*-expression} that each $h_{\{i\}^{c}}((5/8, 2] \times S_0)$ is open relative to $M_{\{i\}^{c}}^{*}$, which together with part (a) shows that $G((\cA_{r} \cap \cN \setminus \cN_{5/8}) \times S_0)$ is open in $M$, and we are done since $G( (\cA_{r} \cap \mathring{\cN}_{2/3}) \times S_0)$ is also an open subset of $M$.
\end{proof}
The map $g_J$ studied in Proposition~\ref{prop:g-J-properties} above is to serve as the analogue of $f_0^{-1}$ in point (i) at the start of Appendix~\ref{subsec:different-strata-I}, and we shall also need versions of the maps designated $f_I$, $h_I$ and $\phi_I$ there. Recalling that $\tau_0 \in (\cos^2\theta, \frac{1}{ 1 + \sin^2\theta})$, we choose constants $\mu_1, \mu_0$ so that
\begin{equation}\label{eq:mu1-mu0-choices}
\tau_0 < \mu_1 < \frac{1- \tau_0}{\sin^2\theta},\quad 0 < \mu_0 < \frac{\mu_1 - \tau_0}{2\sin^2\theta}.
\end{equation}
It follows that $\mu_1, \mu_0 < 1$, and thus
\begin{equation}\label{eq:mu1-mu0-inequalities}
\max\{\mu_0, \tau_0\} < \tau_0 + \mu_0\sin^2\theta < \tau_0 + 2\mu_0\sin^2\theta < \mu_1 < \tau_0 + \mu_1\sin^2\theta < 1.
\end{equation}
In addition, we fix $\tau, \delta > 0$ such that
\begin{equation}\label{eq:transition-parameter}
\begin{split}
&\tau_0 + \mu_{0}\sin^2\theta  < \tau-2\delta <  \tau + 2\delta < \tau_0 + 2\mu_{0}\sin^2\theta\\
&\tau_0 + \sin^2\theta  < 2 - \delta.
\end{split}
\end{equation}
As result, regarding the map~\eqref{eq:3-strata-smoothing-Psi}, we get from~\eqref{eq:Psi-distance-estimates} that
\begin{equation}\label{eq:nested-smoothing-1}
\begin{split}
&\Psi([0, \mu_0\alpha] \times \partial\Omega) \subset \Omega \cap B_{(\tau-2\delta)\alpha}( \partial V^{n-1}_{J}),\\
&\Psi([0, \mu_1\alpha] \times \partial\Omega) \subset \Omega \cap B_{\alpha}( \partial V^{n-1}_{J}),\\
&\Psi([0, \alpha) \times \partial\Omega) \subset \Omega \cap B_{(2-\delta)\alpha}( \partial V^{n-1}_{J}),
\end{split}
\end{equation}
while from~\eqref{eq:distance-neighborhood-contained} we have
\begin{equation}\label{eq:nested-smoothing-2}
\begin{split}
& \Omega \cap B_{\mu_0\alpha}( \partial V^{n-1}_{J}) \subset \Psi([0, \mu_0\alpha) \times \partial\Omega),\\
&\Omega \cap \overline{B_{\mu_0\alpha}( \partial V^{n-1}_{J})} \subset \Psi([0, \mu_0\alpha] \times \partial\Omega),\\
& \Omega \cap \overline{B_{(\tau+2\delta) \alpha}( \partial V^{n-1}_{J})} \subset \Psi([0, \mu_1\alpha) \times \partial\Omega).
\end{split}
\end{equation}

Next, choose $\mu_2 < \mu_3< \mu_4 < \mu_5$ from $(\mu_1, 1)$, along with open intervals $A(\alpha)$ and $A_1(\alpha)$ satisfying
\begin{equation}\label{eq:intervals-with-alpha}
\underbrace{(\mu_{1}\alpha, \mu_{2}\alpha)}_{=: A_{2}(\alpha)} \Subset A_{1}(\alpha) \Subset  \underbrace{(\mu_0\alpha, \mu_3\alpha)}_{=: A_{0}(\alpha)} \Subset A(\alpha) \Subset (\frac{\mu_0\alpha}{2}, \mu_4\alpha),
\end{equation}
and consider the sets
\begin{equation}\label{eq:tubes-around-smoothing}
\begin{split}
\cW_J :=\ & g_J(B_{\frac{\alpha}{2}}^{n-3, J} \times (-\frac{\mu_0\alpha}{2}, \frac{\mu_0\alpha}{2}) \times \partial \widetilde{M}_{J}),\\
\cV_J :=\ & g_J(B_{\frac{\alpha}{2}}^{n-3, J} \times (-\mu_4\alpha, \mu_4\alpha) \times \partial \widetilde{M}_{J}).
\end{split}
\end{equation}
By Proposition~\ref{prop:g-J-properties}(b) we have
\begin{equation}\label{eq:V-J-expression}
\begin{split}
\cV_{J} \cap \widetilde{M}_{J} =\ & g_{J}(\{\0^{n-3}\} \times [0, \mu_4\alpha) \times \partial\widetilde{M}_{J}),\\
\cW_{J} \cap \widetilde{M}_{J} =\ &  g_J(\{\0^{n-3}\} \times [0, \frac{\mu_0\alpha}{2})  \times \partial \widetilde{M}_{J}).
\end{split}
\end{equation}
To continue, we define
\[
M_{J}^{*} := \widetilde{M}_{J} \setminus \overline{\cW_{J}}.
\]
\begin{lemm}\label{lemm:M-J-*-properties}
$M_{J}^{*}$ is a relatively open subset with compact closure in $\mathring{M}_{J}$. Moreover, we have
\[
\overline{M_{J}^{*}} \cap G((\cN \cap \overline{\cA_{\frac{\mu_0\alpha}{4}}}) \times S_0) = \emptyset,
\]
and that
\[
\overline{M_{J}^{*}} \cap M_{i} = \emptyset,\quad\text{for all }i \not\in J.
\]
\end{lemm}
\begin{proof}
Since $\widetilde{M}_{J}$ is compact, as observed before Lemma~\ref{lemm:smoothed-strata}, while $\cW_{J}$ is open, we see that $\widetilde{M}_{J} \setminus\cW_{J}$ is compact. On the other hand, from~\eqref{eq:V-J-expression} and Proposition~\ref{prop:distance-neighborhood-in-smoothing}(a), we obtain
\begin{equation}\label{eq:deleted-distance}
\widetilde{M}_{J} \setminus \cW_{J} \subset \widetilde{M}_{J} \setminus G((\cN \cap \cA_{\frac{\mu_0 \alpha}{2}}) \times S_0),
\end{equation}
which along with the inclusion $\widetilde{M}_{J} \subset M_J$ and Lemma~\ref{lemm:1-strata-nbhd-facts}(c) gives
\begin{equation}\label{eq:intersection-not-in-J}
(\widetilde{M}_{J} \setminus \cW_{J}) \cap M_{i} = \emptyset \quad \text{for all }i \not\in J,
\end{equation}
and consequently
\[
\widetilde{M}_{J}\setminus \cW_{J} \subset \mathring{M}_{J}.
\]
Since $M_{J}^* \subset \widetilde{M}_{J} \setminus \cW_{J}$, and we have just seen the latter is a compact subset of $\mathring{M}_J$, we conclude that $M_{J}^{*}$ has compact closure in $\mathring{M}_J$. Moreover, from~\eqref{eq:deleted-distance} and respectively~\eqref{eq:intersection-not-in-J}, we infer the two disjointness assertions of the lemma.

It remains to prove that $M_J^*$ is open in $\mathring{M}_{J}$. To that end, recall that $\mathring{M}_{J} = M_{J} \setminus \partial M_J$, since in the present case $J \neq \{1, \cdots, n\}$. Thus Lemma~\ref{lemm:smoothed-strata} implies
\begin{equation}\label{eq:interior-relatively-open}
\begin{split}
\widetilde{M}_{J}\setminus \partial \widetilde{M}_{J} =\ & \big[(\widetilde{M}_{J}\setminus \partial \widetilde{M}_{J}) \setminus G(\overline{B_{\rho_1}} \times S_0)  \big] \cup \big[(\widetilde{M}_{J}\setminus \partial \widetilde{M}_{J}) \cap G(B_{3\rho_1} \times S_0)  \big] \\
=\ & \big[\mathring{M}_{J} \setminus G(\overline{B_{\rho_1}} \times S_0)\big] \cup G(((\Omega\setminus \partial\Omega)\cap B_{3\rho_1}) \times S_0).
\end{split}
\end{equation}
Next, in the notation of Proposition~\ref{prop:smoothing-flow}, we have by~\eqref{eq:graphing-function-relations-1} and~\eqref{eq:graph-sectors} that
\[
\Omega \setminus \partial\Omega = \{(x, y) \in W\ |\ y > \bphi(0, x)\} \subset \mathring{V}^{n-1}_{J},
\]
which shows in particular that $\Omega\setminus \partial\Omega$ can be expressed as $\mathring{V}^{n-1}_{J} \cap \cO$ for some open subset $\cO$ of $V^{n-1}$. Proposition~\ref{prop:1234-with-123}(d) then gives
\[
G(((\Omega\setminus \partial\Omega)\cap B_{3\rho_1}) \times S_0) = \mathring{M}_{J} \cap G((\cO \cap B_{3\rho_1}) \times S_0),
\]
and we conclude from~\eqref{eq:interior-relatively-open} that $\widetilde{M}_{J}\setminus \partial \widetilde{M}_{J}$ is open relative to $\mathring{M}_{J}$. Combining this with 
\[
M_J^{*} = (\widetilde{M}_{J}\setminus \partial\widetilde{M}_{J}) \setminus \overline{\cW_J},
\]
which holds since $\partial\widetilde{M}_{J} \subset \cW_{J}$ by Proposition~\ref{prop:g-J-properties}(b), we get that $M_J^{*}$ is open in $\mathring{M}_{J}$ as asserted. The proof is complete.
\end{proof}
By Lemma~\ref{lemm:M-J-*-properties} and Proposition~\ref{prop:global-multisection-chart}, there exist a neighborhood $\cU_J$ of $M_{J}^{*}$ and a diffeomorphism
\begin{equation}\label{eq:chart-near-4-strata}
f_J: \cU_J \to B_{1}^{n-3 , J} \times M_{J}^{*}
\end{equation}
satisfying the analogues of~\eqref{eq:f123-center} and~\eqref{eq:f123-compatible}. Furthermore, the sets $\cU_J$ can be chosen to satisfy
\begin{equation}\label{eq:U-disjoint-4-strata}
\overline{\cU_J} \cap M_{i} = \emptyset \quad \text{for all }i \not\in J,
\end{equation}
\begin{equation}\label{eq:U-disjoint-distance}
\overline{\cU_{J}} \cap G((\cN \cap \overline{\cA_{\frac{\mu_0\alpha}{4}}}) \times S_0) = \emptyset,
\end{equation}
and that 
\begin{equation}\label{eq:U-mutually-disjoint}
\overline{\cU_{J}} \cap \overline{\cU_{J'}} = \emptyset,\quad\text{whenever }J \neq J',
\end{equation}
where for~\eqref{eq:U-disjoint-distance} we used also the fact that $G((\cN \cap \overline{\cA_{\frac{\mu_0\alpha}{4}}}) \times S_0)$ is compact by~\eqref{eq:N-A-closure-intersection}.

Next consider the map from $[0, \mu_{4}\alpha) \times \partial \widetilde{M}_{J}$ to $\widetilde{M}_{J}$ given by
\[
(s, q) \mapsto g_J(\0^{n-3}, s, q).
\]
Proposition~\ref{prop:g-J-properties} shows that this indeed takes values in $\widetilde{M}_{J}$, and is an injective immersion. Recalling also that $g_{J}(\0^{n-3}, 0, q) = q$ for all $q \in \partial\widetilde{M}_{J}$, we deduce, from Lemma~\ref{lemm:boundary-IVT} for instance, that the above parametrizes a collar neighborhood of $\partial \widetilde{M}_{J}$ in $\widetilde{M}_{J}$. Since the latter is a ($4$-dimensional) $1$-handlebody, by the procedure sketched below~\eqref{eq:1234-parametrize-collar}, taking instead $m = 3$ in Lemma~\ref{lemm:Omega-S-basics} and Proposition~\ref{prop:bh-properties}, we get a Lipschitz map
\begin{equation}\label{eq:collar-pushout-4-strata}
h_{J}:[0, \mu_5\alpha] \times \partial \widetilde{M}_{J} \to \widetilde{M}_{J},
\end{equation}
such that 
\vskip 1mm
\begin{enumerate}
\item $t \mapsto h_{J}(t, \cdot)$ is continuous from $[0, \mu_{5}\alpha)$ into $C^{1}(\partial\widetilde{M}_{J}; M)$.
\vskip 2mm
\item  $\Gamma_J: = h_J(\{\mu_5\alpha\} \times\partial \widetilde{M}_{J})$ has finite $\cH^{1}$-measure. Moreover, writing $h_{J, t}$ for $h_{J}(t, \cdot):\partial\widetilde{M}_{J} \to M$, and setting 
\[
\|\Lambda^2 dh_{J, t}\|: = \sup\big\{ |(dh_{J, t})_{q}(u_1 \wedge u_2)|\ \big|\ q \in \partial\widetilde{M}_{J},\ u_1, u_2 \text{ orthonormal in }T_{q}\partial\widetilde{M}_{J}  \big\},
\]
then we have 
\begin{equation}\label{eq:h-J-shrink-area}
\|\Lambda^2 dh_{J, t}\| \longrightarrow 0\quad\text{as }t \to (\mu_5\alpha)^{-}.
\end{equation}
\vskip 2mm
\item $h_{J}([0, \mu_5\alpha] \times \partial\widetilde{M}_{J}) = \widetilde{M}_{J}$. Also,
\begin{equation}\label{eq:h-J-disjoint-images}
h_{J}(\{t\} \times \partial \widetilde{M}_{J}) \cap h_{J}(\{t'\} \times \partial \widetilde{M}_{J}) = \emptyset, \quad\text{whenever }t \neq t'.
\end{equation}
\vskip 2mm
\item $h_{J}(s, q) = g_{J}(\0^{n-3}, s, q)$ for all $(s, q) \in [0, \mu_4\alpha] \times \partial\widetilde{M}_{J}$. In particular, $h_J|_{(0, \mu_{4}\alpha) \times \partial\widetilde{M}_{J}}$ is a diffeomorphism onto an open set in $\widetilde{M}_{J} \setminus \partial\widetilde{M}_{J}$.
\end{enumerate}
\vskip 1mm
Combining properties (3) and (4) with~\eqref{eq:V-J-expression} yields
\begin{equation}\label{eq:h-J-image-outside-collar}
\begin{split}
h_J((\frac{\mu_0\alpha}{2}, \mu_{5}\alpha] \times \partial\widetilde{M}_{J}) =\ & \widetilde{M}_{J} \setminus h_J([0, \frac{\mu_0\alpha}{2}] \times \partial\widetilde{M}_{J}) = M^*_{J}.
\end{split}
\end{equation}
Imitating~\eqref{eq:collar-spread-out}, we define 
\[
\phi_J:B^{n-3, J}_{1} \times (\frac{\mu_0\alpha}{2}, \mu_5\alpha] \times \partial \widetilde{M}_{J} \to B^{n-3, J}_{1} \times M_J^{*}
\]
by
\[
\phi_J(x, s, q) = (x, h_J(s, q)).
\]
Then $\phi_J$ restricts to a diffeomorphism on $B^{n-3, J}_{1} \times (\frac{\mu_0\alpha}{2}, \mu_4\alpha) \times \partial\widetilde{M}_{J}$. Moreover, we have for all $(s, q) \in (\frac{\mu_{0}\alpha}{2}, \mu_4\alpha) \times \partial \widetilde{M}_{J}$ that
\begin{equation}\label{eq:4-strata-transition-center}
(f_J^{-1} \circ \phi_J)(\0^{n-3}, s, q) = h_J(s, q) = g_J(\0^{n-3}, s, q),
\end{equation}
and that these all lie in $\cU_J \cap \cV_J \setminus \overline{\cW_J}$. Following the argument leading up to the inclusions~\eqref{eq:transition-inclusion-1} and~\eqref{eq:transition-inclusion-2}, we obtain $\lambda_0 = \lambda_{0, J}$ and $r_0 = r_{0, J}$ in $(0, \frac{\alpha}{2})$ such that $r_0 < \lambda_0$ and that
\begin{equation}\label{eq:fJ-phiJ-inclusion}
(f_J^{-1} \circ \phi_J)(B^{n-3, J}_{\lambda_0} \times A(\alpha) \times \partial\widetilde{M}_{J}) \Subset \cU_J \cap \cV_J \setminus \overline{\cW_J},
\end{equation}
\begin{equation}\label{eq:gJ-inclusion}
g_J(B^{n-3, J}_{r_0} \times A_0(\alpha) \times \partial \widetilde{M}_{J}) \Subset (f_J^{-1}\circ \phi_J)(B^{n-3, J}_{\lambda_0} \times A(\alpha) \times \partial\widetilde{M}_{J}).
\end{equation}
We then define 
\[
\theta = \theta_J: B^{n-3, J}_{r_0} \times A_0(\alpha) \times \partial \widetilde{M}_{J} \to B^{n-3, J}_{\lambda_0} \times A(\alpha) \times \partial\widetilde{M}_{J}
\]
by 
\[
\theta := (\phi_J\big|_{B^{n-3, J}_{1} \times  (\frac{\mu_0\alpha}{2}, \mu_4\alpha) \times \partial\widetilde{M}_{J}})^{-1} \circ f_J \circ g_J.
 \]
As before we write $\theta$ as $(\overline{x}, \overline{s}, \overline{q})$ according to the product structure on its target.
\begin{lemm}\label{lemm:4-strata-theta-props}
The map $\theta$ defined above has the following properties.
\vskip 1mm
\begin{enumerate}
\item[(a)] $\theta(\0^{n-3}, s, q) = (\0^{n-3}, s, q)$, for all $(s, q) \in A_0(\alpha) \times \partial \widetilde{M}_{J}$.
\vskip 1mm
\item[(b)] Given a non-empty $J_1 \subset J$ and a point $(x, s, q)$ in the domain of $\theta$ satisfying $x \in \mathring{V}^{n-3, J}_{J_1}$, we have $\overline{x}(x,s, q)\in \mathring{V}^{n-3, J}_{J_1}$.
\end{enumerate}
\end{lemm}
\begin{proof}
Conclusion (a) follows at once from~\eqref{eq:4-strata-transition-center}. For part (b), since the domain of $\theta$ involves only positive values of $s$, under the assumption that $x \in \mathring{V}^{n-3, J}_{J_1}$, we have, by Proposition~\ref{prop:g-J-properties}(d) and the inclusions~\eqref{eq:fJ-phiJ-inclusion} and~\eqref{eq:gJ-inclusion} that
\[
g_J(x, s, q) \in  \cU_{J} \cap  M_{J_1} \setminus \big( \cup_{i \in J \setminus J_{1}}M_{i} \big). 
\]
The properties of $f_J$ and the definition of $\phi_J$ then gives 
\[
\overline{x}(x, s, q) \in V^{n-3, J}_{J_{1}} \setminus \big( \cup_{i \in J \setminus J_{1}}V^{n-3, J}_{i} \big) = \mathring{V}^{n-3, J}_{J_{1}},
\]
and we are done with part (b).
\end{proof}
To continue, take a cut-off function $\zeta:\RR \to [0, 1]$ satisfying
\[
\zeta(t) = 1\text{ if }t \leq \mu_1\alpha,\quad \zeta(t) = 0 \text{ if }t \geq \mu_{2}\alpha,
\]
and define 
\[
F_{J} = (\widetilde{x}, \widetilde{s}, \widetilde{q}):  B^{n-3, J}_{r_0} \times A_0(\alpha) \times \partial \widetilde{M}_{J} \to B^{n-3, J}_{\lambda_0} \times A(\alpha) \times \partial\widetilde{M}_{J}
\]
by the following formulas, which are analogous to~\eqref{eq:F-I-definition}:
\begin{equation}\label{eq:F-J-for-4-strata}
\begin{split}
\widetilde{x}(x, s, q) =\ & (1 - \zeta(s))\cdot x + \zeta(s)\cdot \overline{x}(x, s, q),\\
\widetilde{s}(x, s, q) =\ & (1 - \zeta(s))\cdot s + \zeta(s)\cdot \overline{s}(x, s, q),\\
\widetilde{q}(x, s, q) =\ & \overline{q}(\zeta(s)x, s, q).
\end{split}
\end{equation}
Notice that 
\begin{equation}\label{eq:FJ-id-on-center}
F_J(\0^{n-3}, s, q) = (\0^{n-3}, s, q) \quad\text{for all }(s, q) \in A_0(\alpha) \times \partial\widetilde{M}_{J},
\end{equation}
and that
\begin{equation}\label{eq:4-strata-F-endpoints}
(f_J^{-1} \circ \phi_J \circ F_J)(x ,s, q) = \left\{
\begin{array}{ll}
(f_J^{-1} \circ \phi_J)(x, s, q), & \text{ if }s \geq \mu_2\alpha,\\
g_J(x, s, q), & \text{ if }s \leq \mu_1\alpha.
\end{array}
\right.
\end{equation}
Following the proof of Claim~\ref{claim:F-I-multisection}, using Lemma~\ref{lemm:4-strata-theta-props}(b) and~\eqref{eq:U-disjoint-4-strata} in the process, it can be shown that for any $(x, s, q)\in B^{n-3, J}_{r_0} \times A_0(\alpha) \times \partial \widetilde{M}_{J}$ and $j \in J$, we have
\begin{equation}\label{eq:F-J-multisection}
(f_J^{-1}\circ\phi_J \circ F_J)(x, s, q) \in M_j \quad \text{if and only if}\quad x \in V^{n-3, J}_{j}.
\end{equation}
The proofs of Claim~\ref{claim:theta-properties}(c) and Claim~\ref{claim:F-diffeo} carry over as well, and we obtain some $r_1 = r_{1, J} < r_0$ such that 
\[
F_J|_{B^{n-3, J}_{r_1} \times A_1(\alpha) \times \partial\widetilde{M}_{J}}
\]
is a diffeomorphism onto a neighborhood of $\{\0^{n-3}\} \times A_1(\alpha) \times \partial\widetilde{M}_{J}$. In analogy with Claim~\ref{claim:image-of-F}, we also have the following result.
\begin{lemm}\label{lemm:images-4-strata}
There exists $r_2 = r_{2, J} < r_1$ such that the following hold.
\vskip 1mm
\begin{enumerate}
\item[(a)] These three sets are pairwise disjoint:
\[
(f_J^{-1}\circ \phi_J \circ F_J)\big( \overline{B^{n-3, J}_{r_2}} \times [\mu_1\alpha, \mu_2\alpha] \times\partial \widetilde{M}_{J} \big),\quad (f_J^{-1} \circ \phi_J)\big( \overline{B^{n-3, J}_{r_2}} \times (\mu_2\alpha, \mu_5\alpha] \times \partial\widetilde{M}_{J} \big),
\]
\[
g_J\big( \overline{B^{n-3, J}_{r_2}} \times [0, \mu_1\alpha) \times \partial\widetilde{M}_{J} \big).
\]
\vskip 1mm
\item[(b)] With $\tau$ as in~\eqref{eq:transition-parameter}, the first two sets in part (b) are disjoint from $G\big( (\cN \cap \cA_{\tau\alpha}) \times S_0 \big)$.
\vskip 1mm
\item[(c)] The set $g_J\big( \overline{B^{n-3, J}_{r_2}} \times [\mu_0\alpha, \mu_1\alpha] \times \partial\widetilde{M}_{J} \big)$ is disjoint from $G\big( (\cN \cap \cA_{\frac{\mu_0\alpha}{2}}) \times S_0 \big)$.
\end{enumerate}
\end{lemm}
\begin{proof}
Since $\mu_1\alpha$ lies in the interior of $A_1(\alpha)$ and $\tau < \mu_1$ by~\eqref{eq:mu1-mu0-inequalities} and~\eqref{eq:transition-parameter}, there is $\beta > 0$ such that 
\[
\beta \in A_1(\alpha) \cap (\tau\alpha, \mu_1\alpha).
\]
Next, by~\eqref{eq:FJ-id-on-center} and~\eqref{eq:4-strata-transition-center}, we have
\begin{equation}\label{eq:agree-on-center-for-images}
f_J^{-1}\circ \phi_J \circ F_J = f_{J}^{-1} \circ \phi_J = g_J\quad\text{on }\{\0^{n-3}\} \times A_0(\alpha) \times \partial \widetilde{M}_{J}.
\end{equation}
The injectivity of $g_{J}$ then gives
\begin{equation}\label{eq:images-for-4-strata-1}
(f_J^{-1} \circ \phi_J \circ F_J)\big( \{\0^{n-3}\} \times \overline{A_2(\alpha)} \times \partial\widetilde{M}_{J} \big) \cap g_J\big( \{\0^{n-3}\} \times [0, \beta] \times \partial \widetilde{M}_{J} \big) = \emptyset,
\end{equation}
while by~\eqref{eq:h-J-disjoint-images}, we have
\begin{equation}\label{eq:images-for-4-strata-2}
\begin{split}
&(f_J^{-1} \circ \phi_J \circ F_J)\big( \{\0^{n-3}\} \times \overline{A_2(\alpha)} \times \partial\widetilde{M}_{J} \big) \\&\cap (f_J^{-1}\circ \phi_J)\big( \{\0^{n-3}\} \times ((\mu_2\alpha, \mu_5\alpha]\setminus A_1(\alpha)) \times \partial\widetilde{M}_{J} \big) = \emptyset,
\end{split}
\end{equation}
Recalling that
\[
g_J\big( \{\0^{n-3}\} \times [0, \mu_1\alpha] \times \partial\widetilde{M}_{J}\big) = h_J\big([0, \mu_1\alpha] \times \partial\widetilde{M}_{J}\big),
\]
and using again~\eqref{eq:h-J-disjoint-images}, we have
\begin{equation}\label{eq:images-for-4-strata-3}
(f_J^{-1} \circ \phi_J)\big(\{\0^{n-3}\} \times [\mu_2\alpha, \mu_5\alpha] \times \partial\widetilde{M}_{J} \big)\cap g_J\big( \{\0^{n-3}\} \times [0, \mu_1\alpha] \times \partial\widetilde{M}_{J}\big) = \emptyset.
\end{equation}
In view of~\eqref{eq:images-for-4-strata-1} through~\eqref{eq:images-for-4-strata-3}, we obtain by a standard compactness argument some $r_2 = r_{2,J} < r_{1, J}$ such that the following pairs of sets are disjoint:
\vskip 1mm
\begin{enumerate}
\item[(i)] $(f_J^{-1} \circ \phi_J \circ F_J)\big( \overline{B^{n-3, J}_{r_2}} \times \overline{A_2(\alpha)} \times \partial\widetilde{M}_{J} \big)$ and $g_J\big( \overline{B^{n-3, J}_{r_2}} \times [0,  \beta] \times \partial \widetilde{M}_{J} \big)$.
\vskip 1mm
\item[(ii)] $(f_J^{-1} \circ \phi_J \circ F_J)\big( \overline{B^{n-3, J}_{r_2}} \times \overline{A_2(\alpha)} \times \partial\widetilde{M}_{J} \big)$ and $(f_J^{-1}\circ \phi_J)\big(\overline{B^{n-3, J}_{r_2}} \times ((\mu_2\alpha, \mu_5\alpha]\setminus A_1(\alpha)) \times \partial\widetilde{M}_{J} \big)$.
\vskip 1mm
\item[(iii)] $(f_J^{-1}\circ \phi_J)\big(\overline{B^{n-3, J}_{r_2}} \times [\mu_2\alpha, \mu_5\alpha] \times \partial\widetilde{M}_{J} \big)$ and $ g_J\big( \overline{B^{n-3, J}_{r_2}} \times [0, \mu_1\alpha] \times \partial\widetilde{M}_{J} \big)$.
\end{enumerate}
\vskip 1mm

For part (a), denoting the three sets in the statement by $E_1$, $E_2$, and $E_3$ for the moment, we see by (iii) that $E_2 \cap E_3 = \emptyset$. Using (i), along with the inclusion $(\beta, \mu_1\alpha) \subset A_1(\alpha)$, the transition property~\eqref{eq:4-strata-F-endpoints}, and the injectivity of $f_J^{-1} \circ \phi_J \circ F_J$ on $\overline{B^{n-3, J}_{r_1}} \times A_1(\alpha) \times \partial\widetilde{M}_{J}$, we get $E_1\cap E_3 = \emptyset$. By the same reasoning, but with (i) replaced by (ii), and $(\beta, \mu_1\alpha)$ replaced by $(\mu_2\alpha, \mu_5\alpha] \cap A_{1}(\alpha)$, we get $E_1 \cap E_2 = \emptyset$.

For parts (b) and (c), note that, by Proposition~\ref{prop:distance-neighborhood-in-smoothing}(a) and our choice of $\beta$, we have 
\[
\widetilde{M}_{J} \cap G((\cN \cap \overline{A_{\tau\alpha}}) \times S_0) \subset g_{J}(\{\0^{n-3}\} \times [0, \beta) \times \partial\widetilde{M}_{J}).
\]
which combines with~\eqref{eq:images-for-4-strata-1} and~\eqref{eq:images-for-4-strata-3}, respectively, to give
\[
(f_J^{-1} \circ \phi_J \circ F_J)\big( \{\0^{n-3}\} \times \overline{A_2(\alpha)} \times \partial\widetilde{M}_{J} \big) \cap G((\cN \cap \overline{A_{\tau\alpha}}) \times S_0) = \emptyset,
\]
and
\[
(f_J^{-1} \circ \phi_J)\big(\{\0^{n-3}\} \times [\mu_2\alpha, \mu_5\alpha] \times \partial\widetilde{M}_{J} \big) \cap G((\cN \cap \overline{A_{\tau\alpha}}) \times S_0)  = \emptyset.
\]
Here we also used the fact that the maps in~\eqref{eq:agree-on-center-for-images} send $\{\0^{n-3}\} \times A_0(\alpha) \times \partial \widetilde{M}_{J}$ into $\widetilde{M}_{J}$. Again by Proposition~\ref{prop:distance-neighborhood-in-smoothing}(a), we also have
\[
\begin{split}
&g_{J}\big( \{\0^{n-3}\} \times [\mu_0\alpha, \mu_1\alpha] \times \partial\widetilde{M}_{J} \big) \cap G\big( (\cN \cap \overline{\cA_{\frac{\mu_0\alpha}{2}}}) \times S_0 \big)\\
&\subset g_{J}\big( \{\0^{n-3}\} \times [\mu_0\alpha, \mu_1\alpha] \times \partial\widetilde{M}_{J} \big) \cap g_{J}\big( \{\0^{n-3}\} \times [0, \mu_0\alpha) \times \partial\widetilde{M}_{J} \big) = \emptyset.
\end{split}
\]
From these observations and another compactness argument, we get (b) and (c) upon decreasing $r_2$ if necessary. The proof is complete.
\end{proof}

Next we define $H_{J}: B^{n-3, J}_{r_2} \times (0, \mu_5\alpha] \times \partial \widetilde{M}_{J} \to M$ by
\begin{equation}\label{eq:H-J-definition}
H_J(x, s, q) = \left\{
\begin{array}{ll}
g_J(x, s, q), & \text{ if }s \in  (0, \mu_1\alpha), \\
(f_J^{-1} \circ \phi_J \circ F_J)(x, s, q), & \text{ if }s \in  A_{1}(\alpha),\\
(f_{J}^{-1}\circ \phi_{J})(x, s, q), & \text{ if }s \in  (\mu_2\alpha, \mu_5\alpha].
\end{array}
\right.
\end{equation}
By~\eqref{eq:4-strata-F-endpoints}, this is is well-defined. Arguing as at the start of the proof of Claim~\ref{claim:G-properties}, we see that $H_{J}$ is Lipschitz. Recalling that $f_{J}^{-1}$ maps into $\cU_{J}$, which together with~\eqref{eq:gJ-inclusion} and~\eqref{eq:intervals-with-alpha} implies that $g_{J}(B^{n-3, J}_{r_2} \times (\mu_0\alpha, \mu_1\alpha) \times \partial\widetilde{M}_{J})$ is contained in $\cU_J$, we see that
\begin{equation}\label{eq:H-J-image-property-1}
H_{J}(B^{n-3, J}_{r_2} \times (\mu_0\alpha, \mu_5\alpha] \times \partial\widetilde{M}_{J}) \subset \cU_{J}.
\end{equation}
On the other hand, by Lemma~\ref{lemm:images-4-strata}, and recalling from~\eqref{eq:mu1-mu0-inequalities} and~\eqref{eq:transition-parameter} that $\mu_0 < \tau$, we have 
\begin{equation}\label{eq:H-J-image-property-2}
H_{J}(B^{n-3, J}_{r_2} \times (\mu_0\alpha, \mu_5\alpha] \times \partial\widetilde{M}_{J}) \cap G\big( (\cN \cap \cA_{\frac{\mu_0\alpha}{2}}) \times S_0 \big) = \emptyset.
\end{equation}
\begin{prop}\label{prop:H-J-properties}
The assignment $(x, s) \mapsto H_{J}(x, s, \cdot)$ defines a continuous map from $B_{r_2}^{n-3, J} \times (0, \mu_{5}\alpha)$ to $C^{1}(\partial\widetilde{M}_{J}; M)$. Moreover, the following additional properties hold.
\vskip 1mm
\begin{enumerate}
\item[(a)] $H_{J}$ restricts to a diffeomorphism on $B^{n-3, J}_{r_2} \times (0, \mu_{4}\alpha) \times \partial\widetilde{M}_{J}$. Also, given disjoint subsets $A$ and $B$ of $(0, \mu_5\alpha]$, with $A \subset (0, \mu_{2}\alpha]$, we have 
\begin{equation}\label{eq:H-J-disjoint-images}
H_{J}(B_{r_2}^{n-3, J} \times A \times  \partial\widetilde{M}_{J}) \cap H_{J}(B_{r_2}^{n-3, J} \times B \times \partial\widetilde{M}_{J}) = \emptyset.
\end{equation}
\vskip 1mm
\item[(b)] We have $H_J(\0^{n-3}, s, q) = h_J(s, q)$ for all $(s, q) \in (0, \mu_5\alpha] \times \partial\widetilde{M}_{J}$. Moreover, for each $x \in B^{n-3, J}_{r_2}$, we have 
\[
\cH^{1}\big(H_J(\{x\} \times \{\mu_5\alpha\} \times \partial\widetilde{M}_{J})\big) < \infty.
\]
\vskip 1mm
\item[(c)] Write $\{i_0, i_1\}$ for $\{1, \cdots, n\}\setminus J$. Given $x \in B_{r_2}^{n-3, J}$, $s \in [\mu_0\alpha, \mu_5\alpha]$, and $\lambda \in\{0, 1\}$, we have 
\[
\cH^{1}\big(H_J(\{x\} \times \{s\} \times h_{J \cup \{i_{\lambda}\}}(\{2\} \times S_0))\big) < \infty.
\]
\vskip 1mm
\item[(d)] Given $(x, s, q)$ in the domain of $H_{J}$ and $i \in J$, we have
\[
H_J(x, s, q) \in M_{i} \quad \text{if and only if}\quad x \in V^{n-3, J}_{i}.
\]
\end{enumerate}
\end{prop}
\begin{proof}
Note that $H_{J}$ restricts to a diffeomorphism on each of the first two regions in~\eqref{eq:H-J-definition}, as well as the part of the third region where $s \in (\mu_{2}\alpha, \mu_{4}\alpha)$. This together with Lemma~\ref{lemm:images-4-strata}(a) shows that $H_{J}$ is an injective local diffeomorphism on $B^{n-3, J}_{r_2} \times (0, \mu_{4}\alpha) \times \partial\widetilde{M}_{J}$, which is the first assertion in part (a). Combining this with property (1) listed below~\eqref{eq:collar-pushout-4-strata}, along with the smoothness of $f_J^{-1}$ and the compactness of $h_{J}([\mu_2\alpha, \mu_5\alpha] \times \partial\widetilde{M}_{J})$, we infer, as in the proof of Claim~\ref{claim:G-properties}, that $(x, s) \mapsto H_J(x, s, \cdot)$ is continuous on $B_{r_2}^{n-3, J} \times (0, \mu_{5}\alpha)$ as a map into $C^{1}(\partial\widetilde{M}_{J}; M)$. For the second statement in (a), again by Lemma~\ref{lemm:images-4-strata}(a), we have
\[
H_{J}(B_{r_2}^{n-3, J} \times [\mu_{4}\alpha, \mu_{5}\alpha] \times \partial \widetilde{M}_{J}) \cap 
H_{J}(B_{r_2}^{n-3, J} \times (0, \mu_{2}\alpha] \times \partial \widetilde{M}_{J}) = \emptyset.
\]
Since $H_{J}$ is injective on $B_{r_2}^{n-3, J} \times (0, \mu_{4}\alpha) \times \partial \widetilde{M}_{J}$, we get~\eqref{eq:H-J-disjoint-images} as desired upon splitting $B$ into $B \cap (0, \mu_{4}\alpha)$ and $B \cap [\mu_{4}\alpha, \mu_{5}\alpha]$.

For the first conclusion of part (b), by property (4) listed below~\eqref{eq:collar-pushout-4-strata}, we get the asserted equality when $s \in (0, \mu_1\alpha)$. In the case $s \in [\mu_1\alpha, \mu_5\alpha]$, we use instead the definition of $\phi_J$ and the properties of $f_J$, recalling in addition~\eqref{eq:FJ-id-on-center} when $s \in [\mu_1\alpha, \mu_2\alpha]$. For the second conclusion of (b), by the definition of $H_J$ we have
\[
H_{J}(\{x\} \times \{\mu_5\alpha\} \times \partial\widetilde{M}_{J}) = f_{J}^{-1}(\{x\} \times \Gamma_{J}).
\]
Since $\Gamma_{J}$ has finite $\cH^{1}$-measure, we deduce with the help of Remark~\ref{rmk:extension} that so does $H_J(\{x\} \times \{\mu_5\alpha\} \times \partial\widetilde{M}_{J})$. 

For part (c), by Proposition~\ref{prop:smoothed-strata-portions} we have $h_{J \cup \{i_\lambda\}}(\{2\} \times S_0) \subset \partial\widetilde{M}_{J}$, so the set in question makes sense. It has finite $\cH^{1}$-measure because $h_{J \cup \{i_\lambda\}}(\{2\} \times S_0)$ does (see below~\eqref{eq:1234-parametrize-collar}), and because $H_{J}$ is Lipschitz. 

For part (d), we take any $(x, s, q)$ in the domain of $H_{J}$ and $i \in J$, and note that, in the three cases considered in~\eqref{eq:H-J-definition}, the asserted equivalence follows, respectively, from Proposition~\ref{prop:g-J-properties}(d),~\eqref{eq:F-J-multisection}, and the properties of $f_J$.
\end{proof}

We are almost in a position to define the analogue of the map $G$ in Proposition~\ref{prop:1234-with-123}. Writing $\Psi_J$ and $\Omega_J$ for $\Psi$ and $\Omega$ to emphasize their dependence on the choice of $J$, we consider the assignment
\begin{equation}\label{eq:Psi-J-augmented}
\widehat{\Psi}_J : (x, s, z) \mapsto x + \Psi_{J}(s, z),
\end{equation}
which defines a diffeomorphism from $V^{n-3, J} \times (-\infty, \alpha) \times \partial\Omega_{J}$ onto an open set in $V^{n-1}$. Then, given $0 < r \leq \alpha$, $E \subset [0, \alpha)$, and $R \in [3/8, 2]$, we define
\begin{equation}\label{eq:C-J-E-definition}
\begin{split}
C^{J}_{r, R}(E)  :=\ & \widehat{\Psi}_{J}( B^{n-3, J}_{r} \times E \times (\partial\Omega_J \cap \cN_{R})), \\
\mathring{C}^{J}_{r, R}(E)  :=\ & \widehat{\Psi}_{J}(B^{n-3, J}_{r} \times E \times (\partial\Omega_J \cap \mathring{\cN}_{R})).\\
\end{split}
\end{equation}
In the case $R = 2$, we drop the subscript $R$. With $\tau$ and $\delta$ as in~\eqref{eq:transition-parameter}, we have by the third inclusion in~\eqref{eq:nested-smoothing-1}, the relation~\eqref{eq:distance-attained}, and the triangle inequality that 
\[
 \widehat{\Psi}_{J}(B^{n-3, J}_{\delta\alpha} \times [0, \alpha) \times \partial\Omega_{J}) \subset B^{n-3, J}_{\delta\alpha} + (\Omega_{J} \cap \cA_{(2-\delta)\alpha}) \subset \cA_{2\alpha},
\]
which together with Lemma~\ref{lemm:preservation-by-flow} gives
\begin{equation}\label{eq:2-strata-nbd-contained-in-N}
C^{J}_{\delta\alpha, R}([0, \alpha) )  \subset  \cA_{2\alpha} \cap \cN_{R}.
\end{equation}
On the other hand, again using~\eqref{eq:distance-attained}, along with the last inclusion in~\eqref{eq:nested-smoothing-2}, we have
\[
\begin{split}
&\Psi_{J}([\mu_1\alpha, \alpha) \times \partial\Omega_{J} ) \cap \cA_{(\tau + 2\delta)\alpha}  = \emptyset,
\end{split}
\]
and consequently
\begin{equation}\label{eq:overlap-with-distance-neighborhood-2}
C^{J}_{\delta\alpha}([\mu_1\alpha, \alpha)) \cap \cA_{(\tau + \delta)\alpha} = \emptyset.
\end{equation}
Using instead the first inclusion in~\eqref{eq:nested-smoothing-2}, we obtain
\begin{equation}\label{eq:Psi-from-mu0-disjoint}
\Psi_{J}([\mu_{0}\alpha, \alpha) \times \partial\Omega_{J} ) \cap \cA_{\mu_{0}\alpha} = \emptyset,
\end{equation}
and thus
\begin{equation}\label{eq:overlap-with-distance-neighborhood-3}
C^{J}_{\frac{\mu_0\alpha}{2}}([\mu_0\alpha, \alpha)) \cap \cA_{\frac{\mu_{0}\alpha}{2}} = \emptyset.
\end{equation}

Recalling that the thresholds $\lambda_0, r_0, r_1$ and $r_2$ appearing in the arguments leading up to Proposition~\ref{prop:H-J-properties} all depend on the choice of $J \subset\{1, \cdots, n\}$ with $|J| = n-2$, we fix ${h} > 0$ such that 
\begin{equation}\label{eq:N'-height-threshold}
h < \frac{1}{4}\min\{\delta\alpha, \mu_0\alpha\}, \quad \text{and}\quad {h} < \min\{r_{2, J}\ |\ |J| = n-2 \}.
\end{equation}
Given $J$ as above, along with some $i \not\in J$, since $V^{n-1}_{J \cup \{i\}} \subset \cA_{\mu_0\alpha}$, we have by~\eqref{eq:Psi-from-mu0-disjoint} that
\[
V^{n-1}_{i} \cap \Psi_{J}([\mu_0\alpha, \frac{(1 + \mu_{5})\alpha}{2}]  \times (\partial\Omega_{J} \cap \cN)) = \emptyset.
\]
Thus, since $\partial\Omega_{J} \cap \cN$ is compact (Remark~\ref{rmk:portions}), upon decreasing ${h}$ if necessary, we can also assume that
\begin{equation}\label{eq:Phi-J-positive-distance}
\dist\big(V^{n-1}_{J'} ,\ \Psi_{J}([\mu_0\alpha, \frac{(1 + \mu_{5})\alpha}{2}] \times (\partial\Omega_{J} \cap \cN))\big) \geq 4h,
\end{equation}
whenever $J \neq J',\ |J| = |J'| = n-2$. In particular, 
\begin{equation}\label{eq:C-J-positive-distance}
\dist\big(C^{J'}_{h}([\mu_0\alpha, \frac{(1 + \mu_{5})\alpha}{2}]),\  C^{J}_{h}([\mu_0\alpha, \frac{(1 + \mu_{5})\alpha}{2}]) \big) \geq 2h.
\end{equation}

Now, given $R \in [3/8, 2]$ and $\ell \in [\mu_{1}\alpha, \mu_{5}\alpha]$, consider the domains
\begin{equation}\label{eq:N'-definition}
\begin{split}
\cN'_{R,\ell} : = (\cN_{R} \cap \cA_{\tau\alpha} ) \cup \bigcup_{|J| = n-2}C^{J}_{{h}, R}((\mu_0\alpha, \ell]),\\
\mathring{\cN}'_{R, \ell} : = (\mathring{\cN}_{R} \cap \cA_{\tau\alpha}) \cup \bigcup_{|J| = n-2}\mathring{C}^{J}_{{h}, R}((\mu_0\alpha, \ell)),
\end{split}
\end{equation}
which we abbreviate as $\cN'$ and $\mathring{\cN}'$, respectively, in the case $R = 2$ and $\ell = \mu_{5}\alpha$. Recalling from~\eqref{eq:cN-definition} that $\mathring{\cN}_{R}$ is an open set in $V^{n-1}$, we see that so is $\mathring{\cN}'_{R, \ell}$ with the help of~\eqref{eq:C-J-E-definition} and Lemma~\ref{lemm:preservation-by-flow}. Also, observe by~\eqref{eq:2-strata-nbd-contained-in-N} and its proof that 
\begin{equation}\label{eq:N'-contained-in-N}
\cN'_{R, \ell} \subset \cN_{R} \cap \cA_{2\alpha}, \quad  \mathring{\cN}'_{R, \ell} \subset \mathring{\cN}_{R} \cap \cA_{2\alpha}.
\end{equation}
We then define a map 
\[
G_{(1)} : \cN' \times S_0 \to M
\]
as follows:
\begin{enumerate}
\item[(1)] If $(y, p) \in (\cN \cap \cA_{\tau\alpha}) \times S_0$, we simply let
\begin{equation}\label{eq:G-1-definition-case-1}
G_{(1)}(y, p) = G(y, p).
\end{equation}
\vskip 1mm
\item[(2)] If $(y, p) \in C^{J}_{h}((\mu_0\alpha, \mu_5\alpha]) \times S_0$ for some $J \subset \{1, \cdots, n\}$ with length $n-2$, necessarily unique by~\eqref{eq:C-J-positive-distance}, we let $(x, s, z) = (\widehat{\Psi}_{J})^{-1}(y)$, and define
\begin{equation}\label{eq:G-1-definition-case-2}
G_{(1)}(y, p) = H_{J}(x, s, G(z, p)).
\end{equation}
\end{enumerate}
\begin{lemm}\label{lemm:N'-properties}
$\mathring{\cN}'$ is open in $V^{n-1}$. Also, we have $\cN' \setminus \mathring{\cN}' \subset \overline{V^{n-1}\setminus \cN'}$, and that $\cN' \subset \overline{\mathring{\cN}'}$.
\end{lemm}
\begin{proof}
The openness of $\mathring{\cN}'$ is already observed immediately after~\eqref{eq:N'-definition}. For the second conclusion, notice that 
\begin{equation}\label{eq:N'-properties-1}
\cN' \setminus \mathring{\cN}' \subset  [(\cN \setminus \mathring{\cN}) \cap \cA_{\tau\alpha} ] \cup \big(\cup_{|J| = 2}C^{J}_{h}(\{\mu_{5}\alpha\})\big)  \cup \big( \cup_{|J| = 2}\cD_{J} \big),
\end{equation}
where 
\[
\cD_{J} := \widehat{\Psi}_{J}(B_{h}^{n-3, J} \times (\mu_{0}\alpha, \mu_{5}\alpha) \times (\partial\Omega_{J} \cap (\cN\setminus \mathring{\cN}))). 
\]
With the help of~\eqref{eq:N-A-intersection} and~\eqref{eq:N'-contained-in-N}, we see that 
\begin{equation}\label{eq:N'-properties-inclusion-1}
(\cN \setminus \mathring{\cN}) \cap \cA_{\tau\alpha} \subset \overline{V^{n-1} \setminus\cN} \subset \overline{V^{n-1}\setminus\cN'}.
\end{equation}
Next, by~\eqref{eq:overlap-with-distance-neighborhood-2} and~\eqref{eq:C-J-positive-distance}, we have for all $|J| = 2$ that $C^{J}_{h}((\mu_{5}\alpha, (1+\mu_{5})\alpha/2]) \cap \cN' = \emptyset$, and hence 
\begin{equation}\label{eq:N'-properties-inclusion-2}
C^{J}_{h}(\{\mu_{5}\alpha\}) \subset \overline{V^{n-1}\setminus \cN'}.
\end{equation}
Finally, by Lemma~\ref{lemm:preservation-by-flow} we have $\cD_{J} \subset \cN \setminus \mathring{\cN}$. Using also~\eqref{eq:2-strata-nbd-contained-in-N}, we find that
\[
\begin{split}
\cD_{J} \subset (\cN \setminus \mathring{\cN}) \cap C^{J}_{h}((\mu_{0}\alpha, \mu_{5}\alpha)) \subset\ & (\cN \setminus \mathring{\cN}) \cap \cA_{2\alpha}\\
\subset\ & \overline{V^{n-1}\setminus \cN'},
\end{split}
\]
where the last inclusion follows in the same way as~\eqref{eq:N'-properties-inclusion-1}. Substituting the above along with~\eqref{eq:N'-properties-inclusion-1} and~\eqref{eq:N'-properties-inclusion-2} back into~\eqref{eq:N'-properties-1} gives the second assertion of the lemma. For the last conclusion, notice by~\eqref{eq:cN-definition} that $\cN \subset \overline{\mathring{\cN}}$, so that
\[
\cN \cap \cA_{\tau\alpha} \subset \overline{\mathring{\cN} \cap \cA_{\tau\alpha}}
\]
since $\cA_{\tau\alpha}$ is open. On the other hand, the expression~\eqref{eq:N-R-intersect-Omega} from Remark~\ref{rmk:portions}, along with the analogous expression involving $\mathring{\cN}$, shows that $\partial\Omega_{J} \cap \cN \subset \overline{\partial\Omega_{J} \cap \mathring{\cN}}$, and hence
\[
C^{J}_{h}((\mu_{0}\alpha, \mu_{5}\alpha]) \subset \overline{\mathring{C}_{h}^{J}((\mu_0\alpha, \mu_5\alpha))},
\]
thanks to~\eqref{eq:C-J-E-definition}. Combining the previous two displayed inclusions with~\eqref{eq:N'-definition} gives the result.
\end{proof}

\begin{prop}\label{prop:G-1-properties}
$G_{(1)}$ is well-defined Lipschitz map on $\cN' \times S_0$, and $y \mapsto G_{(1)}(y, \cdot)$ is continuous from $\mathring{\cN}'$ into $C^{1}(S_0; M)$. Moreover, the following hold.
\vskip 1mm
\begin{enumerate}
\item[(a)] $G_{(1)}$ restricts to a diffeomorphism on $\mathring{\cN}'_{7/8, \mu_{4}\alpha} \times S_0$. Also, we have
\begin{equation}\label{eq:G-1-somewhere-injective}
G_{(1)}(B^{n-1}_{\frac{\mu_0\alpha}{2}} \times S_0) \cap G_{(1)}((\cN' \setminus B^{n-1}_{\frac{\mu_0\alpha}{2}} ) \times S_0) = \emptyset.
\end{equation}
\vskip 1mm
\item[(b)] For all $y_0 \in \cN' \setminus \mathring{\cN}'$, we have $\cH^{1}(G_{(1)}(\{y_0\} \times S_0)) < \infty$, and that 
\begin{equation}\label{eq:G-1-area-boundary}
\Area(G_{(1)}(y, \cdot)) \to 0 \quad \text{as }y \to y_0 \text{ from within } \mathring{\cN}'.
\end{equation}
\vskip 1mm
\item[(c)] Given $(y, p) \in \cN' \times S_0$ and $i \in \{1, \cdots, n\}$, we have
\[
y \in V^{n-1}_{i} \quad \text{if and only if}\quad G_{(1)}(y, p) \in M_{i}.
\]
\end{enumerate}
\end{prop}
\begin{proof} 
Thanks to~\eqref{eq:C-J-positive-distance}, to see that $G_{(1)}$ is well-defined, we only have to ensure that~\eqref{eq:G-1-definition-case-1} agrees with~\eqref{eq:G-1-definition-case-2} when $y \in (\cN \cap \cA_{\tau\alpha}) \cap C^{J}_{h}((\mu_0\alpha, \mu_5\alpha])$, for some $J \subset \{1, \cdots, n\}$ with $|J| = n-2$. To that end, notice by~\eqref{eq:overlap-with-distance-neighborhood-2} that
\[
(\cN \cap \cA_{\tau\alpha}) \cap C^{J}_{h}\big( (\mu_0\alpha, \mu_5\alpha] \big) \subset \cN \cap C^{J}_{h}\big( (\mu_0\alpha, \mu_1\alpha) \big).
\]
Thus, given $y$ in the left-hand side, upon letting $(x, s, z) = (\widehat{\Psi}_{J})^{-1}(y)$, we must have $s \in (\mu_0\alpha, \mu_1\alpha)$, in which case, by the definition~\eqref{eq:H-J-definition} of $H_{J}$ and Proposition~\ref{prop:g-J-properties}(a), 
\[
\begin{split}
H_{J}(x, s, G(z, p)) =\ & g_J(x, s, G(z, p)) = G(y, p).
\end{split}
\]
This proves that $G_{(1)}$ is well-defined, and the argument also shows that
\begin{equation}\label{eq:G1-G-coincide}
G_{(1)} = G \quad\text{on }\big[(\cN \cap \cA_{\tau\alpha})  \cup \bigcup_{|J| = n-2}C^{J}_{h}((\mu_0\alpha, \mu_1\alpha))\big] \times S_0.
\end{equation}
In particular $G_{(1)}$ is Lipschitz on the set on the right-hand side. On the other hand, since $\partial\Omega_{J} \cap \cN$ is compact (Remark~\ref{rmk:portions}), we see that $(\widehat{\Psi}_{J})^{-1}$ maps $\overline{C^{J}_{h}((\mu_0\alpha, \mu_5\alpha])}$ diffeomorphically onto $\overline{B_{h}^{n-3, J}} \times [\mu_0\alpha, \mu_5\alpha] \times (\partial\Omega_{J} \cap \cN)$. The product of this last set with $S_0$ is then sent by the Lipschitz map
\[
(x, s, z, p) \mapsto (x, s, G(z, p))
\]
onto $\overline{B_{h}^{n-3, J}} \times [\mu_0\alpha, \mu_5\alpha] \times \partial \widetilde{M}_{J}$, on which $H_J$ is Lipschitz continuous. This means that $G_{(1)}$ is Lipschitz on each $C^{J}_{h}((\mu_0\alpha, \mu_5\alpha]) \times S_0$ as well. It follows that $G_{(1)}$ is Lipschitz on its entire domain, since, by~\eqref{eq:overlap-with-distance-neighborhood-2} and~\eqref{eq:C-J-positive-distance}, any pair of points in $\cN' \times S_0$ which are sufficiently close together but do not both lie in the region~\eqref{eq:G1-G-coincide} must belong to the same $C^{J}_{h}((\mu_0\alpha, \mu_5\alpha]) \times S_0$ for some $J$. 

Next, suppose $(y_{i})$ is a sequence converging to some $y_0$ in $\mathring{\cN}'$. In the case $y_0 \in \mathring{\cN} \cap \cA_{\tau\alpha}$, we have eventually that $y_i \in \mathring{\cN} \cap \cA_{\tau\alpha}$ as well, so that, by~\eqref{eq:G-1-definition-case-1} and Proposition~\ref{prop:1234-with-123}, 
\[
G_{(1)}(y_{i}, \cdot) = G(y_{i}, \cdot) \to G(y_0, \cdot) = G_{(1)}(y_0, \cdot),
\]
where the convergence takes place in $C^{1}(S_0; M)$. In the case $y_0 \in \mathring{C}_{h}^{J}((\mu_0\alpha, \mu_5\alpha))$, noting from~\eqref{eq:C-J-E-definition} and Lemma~\ref{lemm:preservation-by-flow} that 
\[
\mathring{C}_{h}^{J}((\mu_0\alpha, \mu_5\alpha)) = \big[B^{n-3, J}_{h} + \Psi_{J}((\mu_0\alpha, \mu_5\alpha) \times \partial\Omega_{J})\big] \cap \mathring{\cN},
\]
which is open in $V^{n-1}$, we can assume without loss of generality that $y_i \in \mathring{C}_{h}^{J}((\mu_0\alpha, \mu_5\alpha))$ for all $i$. Thus, it makes sense to define $(x_i, s_i, z_i) := (\widehat{\Psi}_{J})^{-1}(y_i)$, which yields a sequence in $B^{n-3, J}_{h} \times (\mu_0\alpha, \mu_5\alpha) \times (\partial\Omega_{J} \cap \mathring{\cN})$ that converges to $(x_0, s_0, z_0) := (\widehat{\Psi}_{J})^{-1}(y_0)$. We then get from Proposition~\ref{prop:1234-with-123}, and respectively Proposition~\ref{prop:H-J-properties}, that 
\[
G(z_i, \cdot) \to G(z_0, \cdot) \quad\text{in }C^{1}(S_0; M),
\]
and that 
\[
H_{J}(x_i, s_i, \cdot) \to H_{J}(x_0, s_0, \cdot)\quad\text{in }C^{1}(\partial\widetilde{M}_{J}; M).
\]
Since the maps $G(z_i, \cdot)$ all take values in $\partial\widetilde{M}_{J}$ thanks to~\eqref{eq:smoothing-contained-in-G-image}, we deduce from~\eqref{eq:G-1-definition-case-2} and the above two convergences that 
\[
G_{(1)}(y_i, \cdot) = H_{J}(x_i, s_i, G(z_i, \cdot)) \to H_{J}(x_0, s_0, G(z_0, \cdot)) = G_{(1)}(y_0, \cdot),
\]
in $C^{1}(S_0; M)$ as $i \to \infty$. This completes the proof that $y \mapsto G_{(1)}(y, \cdot)$ is continuous on $\mathring{\cN}'$ as a map into $C^1(S_0; M)$.

For the first claim in part (a), we use the following abbreviations:
\[
\cY_{1} = \mathring{\cN}'_{7/8, \mu_1\alpha},\quad \cY_{2} = \bigcup_{|J| = n-2} \mathring{C}^{J}_{h, 7/8}((\mu_0\alpha,\mu_{4}\alpha)).
\]
Notice that $\cY_{1}$ is open in $V^{n-1}$, as is each $\mathring{C}^{J}_{h, 7/8}((\mu_0\alpha,\mu_{4}\alpha))$. Also, clearly there holds
\[
\cY_{1} \cup \cY_{2} = \mathring{\cN}'_{7/8, \mu_{4}\alpha},
\]
whereas by~\eqref{eq:overlap-with-distance-neighborhood-2} and~\eqref{eq:C-J-positive-distance}, we have
\begin{equation}\label{eq:cY-intersection}
\cY_{1} \cap \cY_{2} = \bigcup_{|J| = n-2} \mathring{C}^{J}_{h, 7/8}((\mu_0\alpha, \mu_1\alpha)).
\end{equation}
Since $G_{(1)} = G$ on $\cY_{1} \times S_0$ by~\eqref{eq:G1-G-coincide}, and since $\cY_{1} \subset \mathring{\cN}_{7/8}$ by~\eqref{eq:N'-contained-in-N}, we deduce from Proposition~\ref{prop:1234-with-123}(a) that $G_{(1)}|_{\cY_{1} \times S_0}$ is a diffeomorphism. Analyzing the composition involved in~\eqref{eq:G-1-definition-case-2} as we did after~\eqref{eq:G1-G-coincide}, but this time restricting $(\widehat{\Psi}_{J})^{-1}$ to $\mathring{C}^{J}_{h, 7/8}((\mu_0\alpha, \mu_4\alpha))$, and noting instead that $G$ maps $(\partial\Omega_{J} \cap \mathring{\cN}_{7/8}) \times S_0$ diffeomorphically onto $\partial\widetilde{M}_{J} \cap G(\mathring{\cN}_{7/8} \times S_0)$ by Remark~\ref{rmk:portions}, we see from Proposition~\ref{prop:H-J-properties}(a) that $G_{(1)}$ restricts to a diffeomorphism on each $\mathring{C}^{J}_{h, 7/8}((\mu_0\alpha, \mu_{4}\alpha)) \times S_0$. Since these open sets have disjoint images under $G_{(1)}$ thanks to~\eqref{eq:H-J-image-property-1} and~\eqref{eq:U-mutually-disjoint}, we conclude that $G_{(1)}$ restricts to a diffeomorphism on $\cY_{2} \times S_0$ as well. Noting in addition that 
\[
\begin{split}
&G_{(1)}((\cY_1\setminus( \cY_1 \cap \cY_2)) \times S_0) \cap
G_{(1)}((\cY_2\setminus( \cY_1 \cap \cY_2)) \times S_0)\\
\subset\ & G((\cN \cap \cA_{\tau\alpha}) \times S_0) \cap \bigcup_{|J| = n-2} H_{J}(B^{n-3, J}_{h} \times [\mu_1\alpha, \mu_5\alpha] \times \partial\widetilde{M}_{J}) = \emptyset,
\end{split}
\]
where the equality at the end follows from Lemma~\ref{lemm:images-4-strata}(b), we conclude that $G_{(1)}$ is a diffeomorphism on $(\cY_1 \cup \cY_2) \times S_0 = \mathring{\cN}'_{7/8, \mu_{4}\alpha} \times S_0$. 

For the second conclusion of part (a), we first use the definition of $G_{(1)}$ to see that
\begin{equation}\label{eq:G-1-protected-nbhd-splitting}
\begin{split}
&G_{(1)}((\cN' \setminus \mathring{\cN}'_{7/8, \mu_{4}\alpha}) \times S_0) \\
&\subset G(((\cN \setminus \mathring{\cN}_{7/8}) \cap \cA_{\tau\alpha}) \times S_0) \cup\bigcup_{|J| = n-2}H_{J}(B_{h}^{n-3, J} \times (\mu_0\alpha, \mu_5\alpha] \times \partial\widetilde{M}_{J}).
\end{split}
\end{equation}
Next, recalling from~\eqref{eq:mu1-mu0-inequalities},~\eqref{eq:transition-parameter}, and~\eqref{eq:alpha-choice} that $\mu_0\alpha < \tau\alpha < \alpha < \frac{\rho}{16}$, we have
\begin{equation}\label{eq:G-1-protected-nbhd-inclusions}
\begin{split}
B^{n-1}_{\frac{\mu_0\alpha}{2}} \subset\ & \cN_{2/3} \cap \cA_{\frac{\mu_0\alpha}{2}}\\
\subset\ & \mathring{\cN}_{7/8}\cap \cA_{\tau\alpha} \subset \mathring{\cN}'_{7/8, \mu_{4}\alpha}.
\end{split}
\end{equation}
Using the first in the above string of inclusions, along with Proposition~\ref{prop:1234-with-123}(a) and respectively~\eqref{eq:H-J-image-property-2}, we see that $G(B^{n-1}_{\frac{\mu_0\alpha}{2}} \times S_0)$ is disjoint from each set on the second line of~\eqref{eq:G-1-protected-nbhd-splitting}. Noting also that $G_{(1)}(B^{n-1}_{\frac{\mu_0\alpha}{2}} \times S_0) = G(B^{n-1}_{\frac{\mu_0\alpha}{2}} \times S_0)$, we deduce that
\[
G_{(1)}(B^{n-1}_{\frac{\mu_0\alpha}{2}} \times S_0) \cap G_{(1)}((\cN' \setminus \mathring{\cN}'_{7/8, \mu_{4}\alpha}) \times S_0) = \emptyset.
\]
Since $B^{n-1}_{\frac{\mu_0\alpha}{2}} \subset  \mathring{\cN}'_{7/8, \mu_{4}\alpha}$, we get the second assertion of part (a) upon recalling that $G_{(1)}$ is injective on $\mathring{\cN}'_{7/8, \mu_{4}\alpha} \times S_0$.

For part (b), given $y_0 \in \cN' \setminus \mathring{\cN}'$, there are three possibilities: 
\vskip 1mm
\begin{enumerate}
\item[(i)] If $y_0 \in (\cN \setminus \mathring{\cN}) \cap \cA_{\tau\alpha}$, we have by Proposition~\ref{prop:1234-with-123}(b) that
\[
\cH^{1}(G_{(1)}(\{y_0\} \times S_0)) = \cH^{1}(G(\{y_0\} \times S_0)) < \infty.
\]
Moreover, if $(y_i)$ is a sequence in $\mathring{\cN}'$ converging to $y_0$, then by the second inclusion in~\eqref{eq:N'-contained-in-N}, along with the fact that $\cA_{\tau\alpha}$ is open, we have $y_{i} \in \mathring{\cN} \cap \cA_{\tau\alpha}$ for all sufficiently large $i$, so that
\[
G_{(1)}(y_{i}, \cdot) = G(y_{i}, \cdot).
\]
Since $(y_i)$ converges to a limit in $\cN \setminus \mathring{\cN}$, we conclude from~\eqref{eq:G-area-boundary} together with the above equality that~\eqref{eq:G-1-area-boundary} holds.
\vskip 1mm
\item[(ii)] If $y_0 = x_0 + \Psi_{J}(s_0, z_0)$ for some $|J| = n-2$ along with $x_0 \in B^{n-3, J}_{h}$, $s_0 \in (\mu_0\alpha, \mu_{5}\alpha)$, and $z_0 \in \partial\Omega_{J} \cap (\cN \setminus\mathring{\cN})$, then 
\[
G_{(1)}(\{y_0\} \times S_0) = H_{J}(\{x_0\} \times \{s_0\} \times G(\{z_0\} \times S_0)),
\]
and we get from Proposition~\ref{prop:H-J-properties}(c) that
\[
\cH^{1}(G_{(1)}(\{y_0\} \times S_0))  < \infty.
\]
Next take a sequence $(y_i)$ in $\mathring{\cN}'$ converging to $y_0$. With the help of~\eqref{eq:C-J-positive-distance}, up to taking a subsequence, we can assume either that $y_i \in \mathring{\cN} \cap \cA_{\tau\alpha}$ for all $i$, or that $y_i \in \mathring{C}^{J}_{h}((\mu_0\alpha, \mu_5\alpha))$ for all $i$. In the former case, we have by~\eqref{eq:G-1-definition-case-1} that $G_{(1)}(y_i, \cdot) = G(y_i,\cdot)$. Noting that $y_0 \in \cN \setminus \mathring{\cN}$ thanks to Lemma~\ref{lemm:preservation-by-flow}, we deduce from Proposition~\ref{prop:1234-with-123}(b) that~\eqref{eq:G-1-area-boundary} holds. In the later case, letting 
\[
(x_i, s_i, z_i) = (\widehat{\Psi}_{J})^{-1}(y_i), 
\]
we obtain a sequence in $B_{h}^{n-3, J} \times (\mu_0\alpha, \mu_{5}\alpha) \times (\partial\Omega_{J} \cap \mathring{\cN})$ converging to $(x_0, s_0, z_0)$, and from~\eqref{eq:G-1-definition-case-2} we have
\[
G_{(1)}(y_i, \cdot) = H_{J}(x_i, s_i, G(z_i, \cdot)).
\]
The assumption $z_0 \in \cN \setminus \mathring{\cN}$ and Proposition~\ref{prop:1234-with-123}(b) then yields $\Area(G(z_i, \cdot)) \to 0$. Since $G(z_i, \cdot)$ all map into $\partial\widetilde{M}_{J}$ by~\eqref{eq:smoothing-contained-in-G-image}, and since $H_{J}(x_i, s_i, \cdot) \to H_{J}(x, s, \cdot)$ in $C^{1}(\partial\widetilde{M}_{J}; M)$ by Proposition~\ref{prop:H-J-properties}, we are again led to~\eqref{eq:G-1-area-boundary}
\vskip 1mm
\item[(iii)] If $y_0 = x_0 + \Psi_{J}(\mu_{5}\alpha, z_0)$ for some $|J| = n-2$, $x_0 \in B^{n-3, J}_{h}$, and $z_0 \in \partial\Omega_{J} \cap \cN$, we see from Proposition~\ref{prop:H-J-properties}(b) that $\cH^{1}(G_{(1)}(\{y_0\} \times S_0)) < \infty$. In view of~\eqref{eq:C-J-positive-distance} and~\eqref{eq:overlap-with-distance-neighborhood-2}, any sequence $(y_i)$ in $\mathring{\cN}'$ converging to $y_0$ eventually lies in $\mathring{C}_{h}^{J}((\mu_0\alpha, \mu_5\alpha))$, and thus has the form $y_i = x_i + \Psi_J(s_i, z_i)$, where $s_i \in (\mu_0\alpha, \mu_5\alpha)$, $z_i\in \partial\Omega_{J} \cap \mathring{\cN}$, and
\[
(x_i, s_i, z_i) \to (x_0, (\mu_5\alpha)^{-}, z_0).
\]
By~\eqref{eq:G-1-definition-case-2} and the third case in~\eqref{eq:H-J-definition}, for all large enough $i$ we have
\[
G_{(1)}(y_i, \cdot) = f_{J}^{-1}(x_i, h_{J}(s_i, G(z_{i}, \cdot))),
\]
and hence, letting $L$ denote the $C^{1}$-norm of $f_J^{-1}$ on $\overline{B_{h}^{n-3, J}} \times h_{J}([\mu_2\alpha, \mu_5\alpha] \times \partial\widetilde{M}_{J})$, there holds
\[
\Area(G_{(1)}(y_i, \cdot)) \leq L^2\cdot \|\Lambda^2 dh_{J, s_i}\| \cdot \Area(G(z_i, \cdot)).
\]
Since $G(z_i, \cdot) \to G(z_0, \cdot)$ in $C^1(S_0; M)$ if $z_0 \in \mathring{\cN}$, and $\Area(G(z_i, \cdot)) \to 0$ if $z_0 \in \cN \setminus \mathring{\cN}$, we conclude from~\eqref{eq:h-J-shrink-area} that~\eqref{eq:G-1-area-boundary} holds.
\end{enumerate} 

Moving on to part (c), in the case where $y$ lies in $(\cN \cap \cA_{\tau\alpha}) \cup \bigcup_{|J| = n-2}C^{J}_{h}((\mu_0\alpha, \mu_1\alpha))$, we get the desired equivalence from~\eqref{eq:G1-G-coincide} and Proposition~\ref{prop:1234-with-123}(d). In the case where $y$ lies in $C^{J}_{h}([\mu_1\alpha, \mu_5\alpha])$ for some $|J| = n-2$, we have
\[
y = x + \Psi_{J}(s, z)
\]
for some $x \in B^{n-3, J}_{h}$, $s \in [\mu_1\alpha, \mu_5\alpha]$, and $z \in \partial\Omega_{J} \cap\cN$, and that
\begin{equation}\label{eq:H-J-properties-c-relation}
G_{(1)}(y, p) = H_J(x, s, G(z, p)) \in \cU_J,
\end{equation}
where the inclusion uses~\eqref{eq:H-J-image-property-1}. Since $s > 0$ and $z \in V^{n-1}_{J}$, we deduce from~\eqref{eq:gauss-image} that 
\[
\Psi_{J}(s, z) \in \mathring{V}^{n-1}_{J}.
\]
Now if $x = \0^{n-3}$, then $y = \Psi_{J}(s, z) \in \mathring{V}^{n-1}_{J}$, while by Proposition~\ref{prop:H-J-properties}(b) and~\eqref{eq:U-disjoint-4-strata}, we have
\[
G_{(1)}(y, p) \in \cU_{J} \cap \widetilde{M}_J \subset \mathring{M}_{J}.
\]
In particular, the asserted equivalence holds in this case since both ends reduce to $i \in J$. For the case $x \neq \0^{n-3}$, assume first that $y \in V^{n-1}_{i}$. Then since $\Psi_J(s, z) \in \mathring{V}^{n-1}_{J}$, we get from~\eqref{eq:projection-of-sectors-2} that $i \in J$, in which case~\eqref{eq:projection-of-sectors} gives $x \in V^{n-3, J}_{i}$, and we conclude that $G_{(1)}(y, p) \in M_{i}$ by~\eqref{eq:H-J-properties-c-relation} and Proposition~\ref{prop:H-J-properties}(d). Conversely, if $G_{(1)}(y, p)$ lies in $M_{i}$, then since it is in $\cU_J$ to begin with, we have by~\eqref{eq:U-disjoint-4-strata} that $i \in J$. Proposition~\ref{prop:H-J-properties}(d) then implies $x \in V^{n-3, J}_{i}$, which combines with~\eqref{eq:projection-of-sectors} to give $y \in V^{n-1}_{i}$. The proof is complete.
\end{proof}
\subsection{Patching across different strata (III)}\label{subsec:different-strata-III}
For this remaining step we specialize to $n = 4$. As in the previous section, we omit the subscripts $\rho$ and $\sigma$ and write $\cN_{\rho, \sigma, R}$ simply as $\cN_{R}$. To begin, we define, for each $I \subset \{1, 2, 3, 4\}$ with $|I| = 3$, 
\begin{equation}\label{eq:M_I-**-definition}
M_{I}^{**} := M_I^* \setminus h_{I}([0, 5/8] \times S_0) = h_I((5/8, 2] \times S_0),
\end{equation}
where $h_I$ is the map described right after~\eqref{eq:1234-parametrize-collar}, and the second equality uses~\eqref{eq:M_I*-expression} along with~\eqref{eq:h-disjoint-images}. Then, for instance, with $\alpha$ still being as in~\eqref{eq:alpha-choice}, we have by Lemma~\ref{lemm:1-strata-nbhd-facts}(a) that, for all $r \in (0, 4\alpha]$,
\[
G(((\cN\setminus \cN_{5/8}) \cap \cA_{r}) \times S_0) = \cup_{|I| = 3}\,f_I^{-1}(B^{2, I}_{r} \times M_I^{**}).
\]
In particular, the set on the left-hand side is open, as $M_I^{**}$ is open relative to $M_I^*$. Next, fix $j_0 \in \{1, 2, 3, 4\}$ and write its complement as $\{j_0\}^{c} = \{i_0,i_1,i_2\}$. Adopting indices modulo $3$ within $\{j_0\}^{c}$, we have for any $\lambda \in \{0, 1, 2\}$ that
\[
\{i_{\lambda}\}^{c} = \{j_0, i_{\lambda+1}, i_{\lambda+2}\}.
\]
Using Proposition~\ref{prop:1234-with-123}(d) along wit Lemma~\ref{lemm:1-strata-nbhd-facts}(a), recalling from~\eqref{eq:U-disjoint} that the image of $f_{\{j_0\}^{c}}^{-1}$ is disjoint from $M_{j_0}$, and then applying~\eqref{eq:f123-compatible}, we get
\begin{equation}\label{eq:product-structure-for-image-V}
G(((\cN\setminus \cN_{5/8}) \cap \cA_{r} \cap V^{3}_{j_0}) \times S_0)
= \cup_{\lambda \in \{0, 1, 2\}}f_{\{i_\lambda\}^c}^{-1}((B^{2, \{i_\lambda\}^c}_{r} \cap V^{2, \{i_{\lambda}\}^{c}}_{j_0} ) \times M_{\{i_\lambda\}^{c}}^{**}).
\end{equation}
Notice also that $\partial_2 V^{3}_{j_0} = \cup_{\lambda \in \{0, 1, 2\}} V^{3}_{\{i_{\lambda}\}^{c}}$, and that~\eqref{eq:distance-attained} becomes
\begin{equation}\label{eq:S2-neighborhood-specialized}
V^{3}_{j_0} \cap B_r\big(\partial_2 V^{3}_{j_0}\big) = \cup_{\lambda \in \{0, 1, 2\}}V^{3}_{j_0} \cap C^{\{i_{\lambda}\}^{c}}_{r}([0, \infty)) = V^{3}_{j_0} \cap \cA_{r}.
\end{equation}
In the context of Lemma~\ref{lemm:sectors-product-structure}, with $I$, $J$, and $i_0$ taken to be $\{j_0\}^{c}$, $\{j_0\}$, and $i_{\lambda}$, respectively, the spaces $W$ and $W'$ in the statement there reduce to $W = V^3$ and $W' = V^{2; \{i_\lambda\}^{c}}$, and the first conclusion of the lemma becomes
\[
V^{3} = V^{2, \{i_\lambda\}^{c}} \oplus \{\ba_{i_\lambda}\}.
\]
With $\alpha$ and $\mu_0$ being the constants from~\eqref{eq:alpha-choice} and~\eqref{eq:mu1-mu0-choices}, respectively, we set
\begin{equation}\label{eq:rho2-choice}
\rho_2 := \frac{\mu_0\alpha}{16}.
\end{equation}
Then, with the constants $\theta_{4, 3} \in (0, \frac{\pi}{2})$ and $C_0(4, 3) > 2$ being as in Lemma~\ref{lemm:foliation-from-flow}, and with $h$ given by~\eqref{eq:N'-height-threshold}, we choose $\alpha' > 0$ satisfying
\begin{equation}\label{eq:alpha'-threshold}
\alpha' < \min\big\{ \frac{\rho_2 \cdot \sin(\theta_{4,3})}{4\cdot (C_0(4, 3) + 1)},\ \frac{h}{4}\big\},
\end{equation}
and apply Proposition~\ref{prop:smoothing-flow} with the parameter ``$\alpha$'' replaced by $\alpha'$. With the help of~\eqref{eq:S2-neighborhood-specialized}, we see that the resulting domain $\Omega \subset V^{3}_{j_0}$ and vector field $\xi: \partial\Omega \to \mathring{V}^{3}_{j_0} \cap B_1$ satisfy
\begin{equation}\label{eq:Omega-agree-specialized}
\Omega \setminus \cA_{\rho_2} = V^{3}_{j_0}\setminus \cA_{\rho_2}, \quad \partial\Omega \setminus \cA_{\rho_2} = \partial V^{3}_{j_0}\setminus \cA_{\rho_2},
\end{equation}
and 
\begin{equation}\label{eq:w-lambda-definition}
\xi(y) = \frac{\ba_{i_{\lambda}; \{j_0, i_{\lambda}\}}}{|\ba_{i_{\lambda}; \{j_0, i_{\lambda}\}}|}= :\bw_{\lambda}, \quad\text{for all }y \in V^{3}_{\{j_0, i_{\lambda}\}}\setminus \cA_{\rho_2}.
\end{equation}
Recall also that there is some $\alpha_{1}' > \alpha'$ so that the map 
\begin{equation}\label{eq:Psi-specialized}
\Psi: (s, y) \mapsto y + s\xi(y)
\end{equation}
is a diffeomorphism from $(-\infty, \alpha_{1}') \times \partial\Omega$ onto an open set in $V^{3}$. In view of Proposition~\ref{prop:smoothing-product-structure} and the inequality 
\[
2\sqrt{\frac{4}{3}}\cdot\frac{C_0(4,3)\cdot \alpha'}{\sin(\theta_{4,3})} < \rho_2,
\]
we further obtain, for $\lambda \in \{0, 1, 2\}$, a smooth domain $\Omega'_{\lambda} \subset V^{2, \{i_{\lambda}\}^{c}}_{j_0}$ such that 
\begin{equation}\label{eq:product-structure-specialized}
\Omega \cap \big( V^{2, \{i_{\lambda}\}^{c}} + t\frac{\ba_{i_\lambda}}{|\ba_{i_\lambda}|} \big) = \Omega_{\lambda}' + t\frac{\ba_{i_\lambda}}{|\ba_{i_\lambda}|},\quad\text{for all }t \geq \rho_2,
\end{equation}
and that
\begin{equation}\label{eq:Omega-lambda-agree}
\begin{split}
\Omega'_{\lambda} \setminus B_{\rho_2} =\ & V^{2; \{i_{\lambda}\}^{c}}_{j_0} \setminus B_{\rho_2}, \\
\partial\Omega'_{\lambda} \setminus B_{\rho_2} =\ & \big( V^{2; \{i_{\lambda}\}^{c}}_{j_0, i_{\lambda+1}} \cup  V^{2; \{i_{\lambda}\}^{c}}_{j_0, i_{\lambda+2}}  \big) \setminus B_{\rho_2},
\end{split}
\end{equation}
where recall that $\partial\Omega_{\lambda}'$ stands for the relative boundary of $\Omega_{\lambda}'$ in $V^{2, \{i_\lambda\}^{c}}$. Moreover, there is a smooth vector field $\xi_{\lambda}': \partial\Omega'_{\lambda} \to V^{2, \{i_{\lambda}\}^{c}}$ such that
\begin{equation}\label{eq:xi-along-section-specialized}
\xi(y + t\frac{\ba_{i_\lambda}}{|\ba_{i_\lambda}|}) = \xi_{\lambda}'(y), \quad\text{for all }y \in \partial\Omega'_{\lambda},\ t \geq \rho_2.
\end{equation}
By~\eqref{eq:N-A-intersection}, the inclusion $\Omega \subset V^{3}_{j_0}$, and~\eqref{eq:cone-half-space-disjoint}, we have
\begin{equation}\label{eq:N-A-intersection-as-product}
\begin{split}
(\cN \setminus \cN_{5/8}) \cap \cA_{r} \cap \Omega = \ & \cup_{\lambda \in \{0, 1, 2\}} C^{\{i_\lambda\}^{c}}_{r}((5/8, 2]) \cap \Omega\\
=\ & \cup_{\lambda \in \{0, 1, 2\}} \big( (B_{r}^{2, \{i_{\lambda}\}^{c}} \cap \Omega_{\lambda}') + (5/8, 2] \cdot \frac{\ba_{i_{\lambda}}}{|\ba_{i_{\lambda}}|} \big),
\end{split}
\end{equation}
where the second equality uses~\eqref{eq:product-structure-specialized}. Using~\eqref{eq:N-A-closure-intersection-2} instead of~\eqref{eq:N-A-intersection}, we have
\[
(\cN \setminus \cN_{5/8}) \cap \overline{\cA_{r}} \cap \Omega 
=\cup_{\lambda \in \{0, 1, 2\}} \big( (\overline{B_{r}^{2, \{i_{\lambda}\}^{c}}} \cap \Omega_{\lambda}') + (5/8, 2] \cdot \frac{\ba_{i_{\lambda}}}{|\ba_{i_{\lambda}}|} \big).
\]
Combining these with~\eqref{eq:G-on-tube}, we see that
\begin{equation}\label{eq:product-structure-for-image-Omega}
G(((\cN\setminus \cN_{5/8}) \cap \cA_{r} \cap \Omega) \times S_0)
= \cup_{\lambda \in \{0, 1, 2\}}f_{\{i_\lambda\}^c}^{-1}((B^{2, \{i_\lambda\}^c}_{r} \cap \Omega'_{\lambda} ) \times M_{\{i_\lambda\}^{c}}^{**}),
\end{equation}
and that this continues to hold if we replace both $\cA_{r}$ and $B^{2, \{i_\lambda\}^c}_{r}$ by their respective closures. Also, recall from Proposition~\ref{prop:smoothing-product-structure} that~\eqref{eq:product-structure-specialized} holds with $\Omega$ and $\Omega_{\lambda}'$ replaced by $\partial\Omega$ and $\partial\Omega_{\lambda}'$. As a result, the relations~\eqref{eq:N-A-intersection-as-product} and~\eqref{eq:product-structure-for-image-Omega}, as well as their stated analogues involving closures, remain valid under the same modifications.

The following is an analogue of Lemma~\ref{lemm:preservation-by-flow}.
\begin{lemm}\label{lemm:N-preserved-by-flow}
Given $s \in [-\alpha', \alpha']$, $R \in [3/8, 2]$, and $y \in \partial\Omega \cap \cA_{4\alpha}$, there holds the following equivalence:
\[
y + s\xi(y) \in \cN_{R} \quad\text{if and only if}\quad y \in \cN_{R}.
\]
Furthermore, the above continues to hold with $\cN_{R}$ replaced by $\mathring{\cN}_{R}$ on both sides.
\end{lemm}
\begin{proof}
By~\eqref{eq:rho2-choice} and~\eqref{eq:alpha'-threshold}, and the fact that $\mu_0 \in (0, 1)$, we have in particular that $\alpha' < \alpha$, while from~\eqref{eq:alpha-choice} we have $B_{8\alpha} \subset \cN_{R}$. Thus, using also the inequalities $|\xi(y)| \leq 1$ and $|s| \leq \alpha'$, if either $y \not\in \cN_{R}$ or $y + s\xi(y) \not\in \cN_{R}$, then necessarily 
\[
y \in \partial\Omega \cap \cA_{4\alpha}\setminus B_{7\alpha} \subset V^{3}_{j_0} \cap \cA_{4\alpha}\setminus B_{7\alpha},
\]
in which case the second equality in~\eqref{eq:S2-neighborhood-specialized} shows that $y$ has the form
\[
y = y' + t\frac{\ba_{i_{\lambda}}}{|\ba_{i_{\lambda}}|},
\]
for some $\lambda \in \{0, 1, 2\}$, $y' \in B^{2, \{i_{\lambda}\}^{c}}_{4\alpha}$, and $t \geq 3\alpha > \rho_2$. We then get from~\eqref{eq:product-structure-specialized} and~\eqref{eq:xi-along-section-specialized} that $y' \in \partial\Omega_{\lambda}'$, and $\xi(y) = \xi_{\lambda}'(y') \in V^{2, \{i_{\lambda}\}^{c}}$. Thus $y$ and $y + s\xi(y)$ both lie in $C^{\{i_\lambda\}^{c}}_{5\alpha}(\{t\})$, so, by~\eqref{eq:tube-intersection}, no matter which end of the asserted equivalence we negate, we get $t> R$, in which case the other end fails also. This proves the result. The same argument gives the version with $\mathring{\cN}_{R}$ in place of $\cN_{R}$. We omit the details.
\end{proof}
To continue, we define
\begin{equation}\label{eq:5-strata-smoothing-definition}
\widetilde{M}_{j_0} = \big( M_{j_0} \setminus G((\cN \cap \overline{\cA_{\rho_2}})\times S_0 )\big) \cup G( (\cN \cap \cA_{3\rho_2} \cap  \Omega) \times S_0),
\end{equation}
Since $\Omega \subset V^{3}_{j_0}$, Proposition~\ref{prop:1234-with-123}(d) implies that $\widetilde{M}_{j_0} \subset M_{j_0}$. We claim that $\widetilde{M}_{j_0}$ is compact. To see this, we first use Proposition~\ref{prop:1234-with-123}(d) and Lemma~\ref{lemm:1-strata-nbhd-facts}(b) to get
\begin{equation}\label{eq:5-strata-smoothing-difference-1}
\begin{split}
M_{j_0} \setminus \widetilde{M}_{j_0} =\ & G((\cN \cap \overline{\cA_{\rho_2}} \cap V^{3}_{j_0}) \times S_0)\setminus G( ( \cN \cap \cA_{3\rho_2} \cap \Omega ) \times S_0 )\\
=\ & G((\cN \cap \overline{\cA_{\rho_2}} \cap V^{3}_{j_0}) \times S_0)\setminus G( ( \cN \cap \overline{\cA_{\rho_2}} \cap \Omega ) \times S_0 )\\
=\ & \big[G(((\cN\setminus \cN_{5/8}) \cap \overline{\cA_{\rho_2}} \cap V^{3}_{j_0}) \times S_0)\setminus G( ( \cN \cap \overline{\cA_{\rho_2}} \cap \Omega ) \times S_0 ) \big]\\
& \cup \big[ G((\mathring{\cN}_{2/3} \cap \overline{\cA_{\rho_2}} \cap V^{3}_{j_0}) \times S_0)\setminus G( ( \cN \cap \overline{\cA_{\rho_2}} \cap \Omega ) \times S_0 ) \big].
\end{split}
\end{equation}
For the set in the first pair of square brackets, we use Proposition~\ref{prop:1234-with-123}(a), followed by~\eqref{eq:product-structure-for-image-V} and~\eqref{eq:product-structure-for-image-Omega} to get that it is equal to
\[
\begin{split}
&G(((\cN\setminus \cN_{5/8}) \cap \overline{\cA_{\rho_2}} \cap V^{3}_{j_0}) \times S_0)\setminus G( ( (\cN \setminus \cN_{5/8}) \cap \overline{\cA_{\rho_2}} \cap \Omega ) \times S_0 )\\
=\ &\cup_{\lambda\in\{0, 1, 2\}} f_{\{i_\lambda\}^c}^{-1}((\overline{B_{\rho_2}} \cap V^{2, \{i_{\lambda}\}^{c}}_{j_0} \setminus \Omega_{\lambda}' ) \times M_{\{i_\lambda\}^{c}}^{**}) \\
=\ &\cup_{\lambda\in\{0, 1, 2\}} f_{\{i_\lambda\}^c}^{-1}((B_{\rho_2} \cap V^{2, \{i_{\lambda}\}^{c}}_{j_0} \setminus \Omega_{\lambda}' ) \times M_{\{i_\lambda\}^{c}}^{**}),
\end{split}
\]
where the second equality follows from~\eqref{eq:Omega-lambda-agree}. Again using Proposition~\ref{prop:1234-with-123}(a), but this time combining it with the injectivity of $G$ on $\mathring{\cN}_{7/8} \times S_0$, followed by~\eqref{eq:Omega-agree-specialized}, we get
\[
\begin{split}
G((\mathring{\cN}_{2/3} \cap \overline{\cA_{\rho_2}} \cap V^{3}_{j_0} ) \times S_0)\setminus G( ( \cN \cap \overline{\cA_{\rho_2}} \cap \Omega ) \times S_0 )= G((\mathring{\cN}_{2/3} \cap \cA_{\rho_2} \cap (V^{3}_{j_0} \setminus \Omega) ) \times S_0).
\end{split}
\]
Applying~\eqref{eq:f123-compatible} and, respectively, Proposition~\ref{prop:1234-with-123}(d) to the previous two observations, and substituting the results back into~\eqref{eq:5-strata-smoothing-difference-1}, we arrive at
\[
M_{j_0} \setminus \widetilde{M}_{j_0} = M_{j_0} \cap \cO,
\]
where
\[
\cO = G((\mathring{\cN}_{2/3} \cap\cA_{\rho_2} \setminus \Omega ) \times S_0 ) \cup \big( \cup_{\lambda \in \{0, 1, 2\}} f_{\{i_\lambda\}^{c}}^{-1}( (B_{\rho_2}^{2, \{i_\lambda\}^{c}} \setminus \Omega_{\lambda}') \times M_{\{i_\lambda\}^{c}}^{**} ) \big),
\]
which is open in $M$ since $\Omega$ and $\Omega_{\lambda}'$ are closed relative to $V^{3}$ and $V^{2, \{i_\lambda\}^{c}}$, respectively, while $M_{\{i_\lambda\}^{c}}^{**}$ is, by~\eqref{eq:M_I-**-definition}, open relative to $M_{\{i_\lambda\}^{c}}^*$. Thus we have shown that $\widetilde{M}_{j_0}$ is closed relative to $M_{j_0}$, and thus compact since the latter is.

\begin{prop}\label{prop:5-strata-smoothing}
$\widetilde{M}_{j_0}$ is a smooth domain in $M$. Moreover, given $\rho_2 \leq r \leq 4\alpha$, we have
\vskip 1mm
\begin{enumerate}
\item[(a)] $\widetilde{M}_{j_0} \setminus G((\cN \cap \overline{\cA_{r}}) \times S_0) = M_{j_0} \setminus G((\cN \cap \overline{\cA_{r}}) \times S_0)$.
\vskip 1mm
\item[(b)] $\widetilde{M}_{j_0} \cap G(( (\cN \setminus \cN_{5/8}) \cap \cA_{r} ) \times S_0)
= G(( (\cN \setminus \cN_{5/8}) \cap \cA_{r} \cap \Omega ) \times S_0)$. Consequently, for each $\lambda \in \{0, 1, 2\}$, we have $\widetilde{M}_{j_0} \cap  f_{\{i_{\lambda}\}^{c}}^{-1} (B_{r}^{2, \{i_{\lambda}\}^{c}} \times M_{\{i_\lambda\}^{c}}^{**} ) = f_{\{i_{\lambda}\}^{c}}^{-1}( (B_{r} \cap \Omega_{\lambda}') \times M_{\{i_\lambda\}^{c}}^{**})$.
\vskip 1mm
\item[(c)] $\widetilde{M}_{j_0} \cap G(( \mathring{\cN}_{2/3} \cap \cA_{r} ) \times S_0) = G(( \mathring{\cN}_{2/3} \cap \cA_{r} \cap \Omega ) \times S_0)$.
\end{enumerate}
Moreover, all three conclusions continue to hold with $\widetilde{M}_{j_0}$, $M_{j_0}$, $\Omega$, and $\Omega_{\lambda}'$ replaced, respectively, by $\partial\widetilde{M}_{j_0}$, $\partial M_{j_0}$, $\partial\Omega$, and $\partial\Omega_{\lambda}'$.
\end{prop}
\begin{proof}
We first prove (a), (b), and (c), and then explain how to deduce that $\widetilde{M}_{j_0}$ is a smooth domain. For part (a), in view of~\eqref{eq:5-strata-smoothing-definition}, the result clearly holds when $r \geq 3\rho_2$. In the case $\rho_2 \leq r < 3\rho_2$, we use Lemma~\ref{lemm:1-strata-nbhd-facts}(b) and~\eqref{eq:Omega-agree-specialized} to see that
\[
\begin{split}
&G( (\cN \cap \cA_{3\rho_2} \cap \Omega) \times S_0) \setminus G(( \cN \cap \overline{\cA_{r}} ) \times S_0) \\
=\ & G( (\cN \cap (\cA_{3\rho_2} \setminus \overline{\cA_{r}} ) \cap \Omega ) \times S_0) = G( (\cN \cap (\cA_{3\rho_2} \setminus \overline{\cA_{r}} )\cap V^{3}_{j_0} ) \times S_0)\\
=\ & M_{j_0} \cap  G( (\cN \cap (\cA_{3\rho_2} \setminus \overline{\cA_{r}})) \times S_0) \quad\text{(Proposition~\ref{prop:1234-with-123}(d))}\\
=\ & M_{j_0} \cap \big[ G( (\cN \cap \cA_{3\rho_2} ) \times S_0) \setminus G( (\cN \cap \overline{\cA_{r}} ) \times S_0)\big] \quad\text{(Lemma~\ref{lemm:1-strata-nbhd-facts}(b))}.
\end{split}
\]
From this and~\eqref{eq:5-strata-smoothing-definition}, and using again that $r \geq \rho_2$, we get the desired conclusion.

For part (b), we observe
\begin{equation}\label{eq:5-strata-smoothing-b-1}
\begin{split}
&M_{j_0} \cap G( ((\cN \setminus \cN_{5/8}) \cap \cA_{r}) \times S_0 ) \setminus G( (\cN \cap \overline{\cA_{\rho_2}}) \times S_0) \\
=\ & M_{j_0} \cap G( ((\cN \setminus \cN_{5/8}) \cap (\cA_{r} \setminus \overline{\cA_{\rho_2}} )) \times S_0 )  \quad \text{(Prop.~\ref{prop:1234-with-123}(a) and Lemm.~\ref{lemm:1-strata-nbhd-facts}(b))}\\
=\ & G( ((\cN \setminus \cN_{5/8}) \cap (\cA_{r} \setminus \overline{\cA_{\rho_2}} ) \cap V^{3}_{j_0}) \times S_0 ) \quad\text{(Prop.~\ref{prop:1234-with-123}(d))}\\
=\ & G( ((\cN \setminus \cN_{5/8}) \cap (\cA_{r} \setminus \overline{\cA_{\rho_2}} ) \cap \Omega) \times S_0 )  \quad\text{(by~\eqref{eq:Omega-agree-specialized})}.
\end{split}
\end{equation}
Similarly, 
\[
\begin{split}
&G( (\cN \cap \cA_{3\rho_2} \cap \Omega) \times S_0) \cap G( ((\cN\setminus \cN_{5/8}) \cap \cA_{r} ) \times S_0)\\
= \ &  G( ((\cN\setminus \cN_{5/8}) \cap \cA_{3\rho_2} \cap \Omega) \times S_0) \cap G( ((\cN\setminus \cN_{5/8}) \cap \cA_{r} ) \times S_0) \quad\text{(Prop.~\ref{prop:1234-with-123}(a))}\\
=\ & G( ((\cN\setminus \cN_{5/8}) \cap \cA_{\min\{r, 3\rho_2\}} \cap \Omega) \times S_0) \quad\text{(Lemm.~\ref{lemm:1-strata-nbhd-facts}(b))}.
\end{split}
\]
Combining this with~\eqref{eq:5-strata-smoothing-b-1} and~\eqref{eq:5-strata-smoothing-definition} gives the first equation in (b). The second follows from intersecting both sides with $f_{\{i_{\lambda}\}^{c}}^{-1}( B_{r}^{2, \{i_{\lambda}\}^{c}}  \times M_{\{i_\lambda\}^{c}}^{**} )$, and then using Lemma~\ref{lemm:1-strata-nbhd-facts}(a) along with~\eqref{eq:product-structure-for-image-Omega} and~\eqref{eq:U-disjoint}. The proof of (c) is exactly the same as (b), except that where Lemma~\ref{lemm:1-strata-nbhd-facts}(b) was used, we use the injectivity of $G$ on $\mathring{\cN}_{7/8} \times S_0$ instead. The details are omitted.

To see that $\widetilde{M}_{j_0}$ is a smooth domain, we fix any $s, r \in [\rho_{2}, 4\alpha]$ with $s < r$, and define
\begin{align*}
\cO_{1}: =\ & M \setminus G((\cN \cap \overline{\cA_{s}}) \times S_0),\\
\cO_{2, i} : =\ &  f_{\{i\}^{c}}^{-1}( B_{r}^{2, \{i\}^{c}}  \times M_{\{i\}^{c}}^{**} ) \quad\text{for each }i \in \{1, 2, 3, 4\},\\
\cO_{3} :=\ & G( (\mathring{\cN}_{2/3} \cap \cA_{r}) \times S_0).
\end{align*}
With the help of Lemma~\ref{lemm:1-strata-nbhd-facts}(a), we have 
\[
M \subset \cO_{1} \cup \cO_{2, 1} \cup \cdots \cup \cO_{2,4} \cup \cO_{3}.
\]
Also, it is clear that $\cO_{3}$ and each $\cO_{2, i}$ are open. By the compactness of $\cN \cap \overline{\cA_{s}}$, noted after~\eqref{eq:N-A-closure-intersection}, we see that $\cO_{1}$ is open as well. Lemma~\ref{lemm:1-strata-nbhd-facts}(c) together with part (a) above shows that $\widetilde{M}_{j_0} \cap \cO_{1}$ is a smooth domain in $\cO_{1}$. Using instead part (b) and the smoothness of $\Omega_{\lambda}'$, we see that $\widetilde{M}_{j_0} \cap \cO_{2, i_{\lambda}}$ is a smooth domain in $\cO_{2, i_{\lambda}}$ for $\lambda \in \{0, 1, 2\}$, while by~\eqref{eq:U-disjoint} we have $\widetilde{M}_{j_0} \cap \cO_{2, j_0} = \emptyset$. That $\widetilde{M}_{j_0} \cap \cO_{3}$ is a smooth domain in $\cO_{3}$ follows from part (c), the smoothness of $\Omega$, and the fact that $G|_{\mathring{\cN}_{7/8} \times S_0}$ is a diffeomorphism. Since $\cO_{1}$, $\cO_{2, 1}, \cdots, \cO_{2, 4}$, and $\cO_{3}$ form an open covering of $M$, we conclude from Lemma~\ref{lemm:submanifold-construction} that $\widetilde{M}_{j_0}$ is smooth, and that its boundary satisfies the indicated analogues of (a), (b), and (c). The proof is complete.
\end{proof}

Now let $\mu_0$, ..., $\mu_5$ , $\delta$, and $\tau$ be the positive constants chosen in the paragraph containing~\eqref{eq:transition-parameter}. In addition, for $\lambda \in \{0, 1, 2\}$, we let $g_{\{j_0, i_{\lambda}\}}$ and $H_{\{j_0,i_{\lambda}\}}$ be the maps from Proposition~\ref{prop:g-J-properties} and Proposition~\ref{prop:H-J-properties}, and let $f_{\{j_0, i_{\lambda}\}}$ and $h_{\{j_0,i_\lambda\}}$ be as in~\eqref{eq:chart-near-4-strata} and~\eqref{eq:collar-pushout-4-strata}. The subset~\eqref{eq:Omega-definition} of $V^{3}_{\{j_0, i_{\lambda}\}}$ and vector field appearing in~\eqref{eq:gauss-image} we henceforth denote by $\partial\Omega_{\{j_0, i_{\lambda}\}}$ and $\xi_{\{j_0, i_{\lambda}\}}$, which are used in Lemma~\ref{lemm:smoothed-strata} to define $\widetilde{M}_{\{j_0, i_{\lambda}\}}$. Likewise, we denote by $\Psi_{\{j_0, i_{\lambda}\}}$ the diffeomorphism in~\eqref{eq:3-strata-smoothing-Psi}. The next result supplements Proposition~\ref{prop:5-strata-smoothing}(a), or more precisely its counterpart for $\partial\widetilde{M}_{j_0}$. 

\begin{prop}\label{prop:4-strata-smoothing-decomp-far}
In the above notation, we have
\begin{equation}\label{eq:4-strata-smoothing-decomp-far}
\begin{split}
\partial \widetilde{M}_{j_0}  \setminus G(( \cN \cap \overline{A_{(\tau - 2\delta)\alpha}} ) \times S_0) \subset \ & \bigcup_{\lambda \in \{0, 1, 2\}} h_{\{j_0, i_{\lambda}\}} ((\mu_0\alpha, \mu_5\alpha] \times \partial\widetilde{M}_{\{j_0, \lambda\}})\\
\subset\ & \partial\widetilde{M}_{j_0}.
\end{split}
\end{equation}
\end{prop}
\begin{proof}
By the version of Proposition~\ref{prop:5-strata-smoothing}(a) for $\partial\widetilde{M}_{j_0}$, we have
\begin{equation}\label{eq:decomp-far-prelim}
\begin{split}
\partial\widetilde{M}_{j_0} \setminus G(( \cN \cap \overline{\cA_{(\tau - 2\delta)\alpha}} ) \times S_0)  =\ & \partial M_{j_0} \setminus G(( \cN \cap \overline{\cA_{(\tau - 2\delta)\alpha}} ) \times S_0) \\
=\ & \cup_{\lambda \in \{0, 1, 2\}} M_{\{j_0, i_{\lambda}\}} \setminus G(( \cN \cap \overline{\cA_{(\tau - 2\delta)\alpha}} ) \times S_0)\\
=\ & \cup_{\lambda \in \{0, 1, 2\}} \widetilde{M}_{\{j_0, i_{\lambda}\}} \setminus G(( \cN \cap \overline{\cA_{(\tau - 2\delta)\alpha}} ) \times S_0),
\end{split}
\end{equation}
where the last line follows from Proposition~\ref{prop:distance-neighborhood-in-smoothing}(b) and the fact that $\tau - 2\delta > \tau_0$ by~\eqref{eq:transition-parameter}. Recalling that 
\[
\widetilde{M}_{\{j_0, i_\lambda\}} = h_{\{j_0, i_{\lambda}\}}([0, \mu_5\alpha] \times \partial\widetilde{M}_{\{j_0, i_{\lambda}\}}),
\]
and noting, by property (4) below~\eqref{eq:collar-pushout-4-strata}, the first inequality in~\eqref{eq:transition-parameter}, and Proposition~\ref{prop:distance-neighborhood-in-smoothing}(a), that
\begin{equation}\label{eq:decomp-far-intermediate-inclusion-1}
h_{\{j_0, i_\lambda\}}([0, \mu_0\alpha] \times\partial\widetilde{M}_{\{j_0,i_\lambda\}}) \subset \widetilde{M}_{\{j_0, i_{\lambda}\}} \cap G( (\cN \cap \cA_{(\tau - 2\delta)\alpha}) \times S_0),
\end{equation}
we deduce from~\eqref{eq:decomp-far-prelim} the first inclusion in~\eqref{eq:4-strata-smoothing-decomp-far}. 

For the second inclusion, we split the interval $(\mu_0\alpha, \mu_5\alpha]$ into $(\mu_{0}\alpha, \mu_{1}\alpha) \cup [\mu_{1}\alpha, \mu_{5}\alpha]$. For the latter, observe that, by the first inclusion in Proposition~\ref{prop:distance-neighborhood-in-smoothing}(a) and~\eqref{eq:h-J-disjoint-images}, we have
\begin{equation}\label{eq:empty-intersection-for-decomp-far}
\begin{split}
&h_{\{j_0, i_\lambda\}}([\mu_1\alpha, \mu_5\alpha] \times\partial\widetilde{M}_{\{j_0,i_\lambda\}}) \cap G( (\cN \cap \cA_{(\tau+2\delta)\alpha}) \times S_0)\\
&\subset h_{\{j_0, i_\lambda\}}([\mu_1\alpha, \mu_5\alpha] \times\partial\widetilde{M}_{\{j_0,i_\lambda\}}) \cap h_{\{j_0, i_\lambda\}}([0, \mu_1\alpha) \times\partial\widetilde{M}_{\{j_0,i_\lambda\}}) = \emptyset.
\end{split}
\end{equation}
Here we have considered the slightly larger set $\cA_{(\tau+2\delta)\alpha}$ instead of $\overline{\cA_{(\tau - 2\delta)\alpha}}$ for later purposes. In any case, combining the above with~\eqref{eq:decomp-far-prelim} gives
\begin{equation}\label{eq:decomp-far-intermediate-inclusion-2}
\begin{split}
h_{\{j_0, i_\lambda\}}([\mu_1\alpha, \mu_5\alpha] \times\partial\widetilde{M}_{\{j_0,i_\lambda\}}) \subset\ & \widetilde{M}_{\{j_0, i_{\lambda}\}} \setminus G(( \cN \cap \overline{\cA_{(\tau - 2\delta)\alpha}} ) \times S_0) \\
\subset\ &  \partial\widetilde{M}_{j_0} \setminus G(( \cN \cap \overline{\cA_{(\tau - 2\delta)\alpha}} ) \times S_0).
\end{split}
\end{equation}
On the other hand, we have by~\eqref{eq:smoothing-contained-in-G-image}, Proposition~\ref{prop:g-J-properties}(a) and Lemma~\ref{lemm:preservation-by-flow} that
\begin{equation}\label{eq:4-strata-smoothing-inclusion-stage-1}
\begin{split}
&h_{\{j_0, i_{\lambda}\}}\big((\mu_0 \alpha, \mu_1\alpha) \times \partial\widetilde{M}_{\{j_0, i_{\lambda}\}}\big) \\ 
=\ & G( (\Psi_{\{j_0, i_{\lambda}\}}((\mu_0\alpha, \mu_1\alpha) \times (\cN \cap \partial\Omega_{\{j_0, i_{\lambda}\}}))) \times S_0)\\
=\ &  G( (\cN \cap \Psi_{\{j_0, i_{\lambda}\}}((\mu_0\alpha, \mu_1\alpha) \times \partial\Omega_{\{j_0, i_{\lambda}\}})) \times S_0).
\end{split}
\end{equation}
Using the inclusions~\eqref{eq:nested-smoothing-1} and~\eqref{eq:nested-smoothing-2} along with~\eqref{eq:distance-attained}, we see that
\[
\Psi_{\{j_0, i_{\lambda}\}}\big((\mu_0\alpha, \mu_1\alpha)\times \partial\Omega_{\{j_0, i_{\lambda}\}}\big) \subset \Omega_{\{j_0, i_\lambda\}} \cap \cA_{\alpha}\setminus \overline{\cA_{\mu_0\alpha}}.
\]
Recalling $\Omega_{\{j_0,i_{\lambda}\}} \subset V^{3}_{\{j_0, i_{\lambda}\}} \subset \partial V^{3}_{j_0}$, and then using~\eqref{eq:Omega-agree-specialized} and our choice~\eqref{eq:rho2-choice} of $\rho_2$, we get
\begin{equation}\label{eq:Psi-image-for-decomp-far}
\begin{split}
\Psi_{\{j_0, i_{\lambda}\}}\big((\mu_0\alpha, \mu_1\alpha)\times \partial\Omega_{\{j_0, i_{\lambda}\}}\big) \subset\ &  V^{3}_{\{j_0,i_{\lambda}\}} \cap (\cA_{\alpha}\setminus \overline{\cA_{\mu_0\alpha}})\\
\subset\ & \partial V^{3}_{j_0} \cap (\cA_{\alpha}\setminus \overline{\cA_{\mu_0\alpha}}) = \partial \Omega \cap (\cA_{\alpha}\setminus \overline{\cA_{\mu_0\alpha}}).
\end{split}
\end{equation}
Substituting this back into~\eqref{eq:4-strata-smoothing-inclusion-stage-1}, and using Lemma~\ref{lemm:1-strata-nbhd-facts}(b) and Proposition~\ref{prop:5-strata-smoothing}, we see that 
\[
\begin{split}
h_{\{j_0, i_{\lambda}\}}\big((\mu_0 \alpha, \mu_1\alpha) \times \partial\widetilde{M}_{\{j_0, i_{\lambda}\}}\big) \subset\ & G( ( \partial\Omega \cap (\cA_{\alpha}\setminus \overline{\cA_{\mu_0\alpha}}) \cap \cN )\times S_0)\\
= \ & G( ( \partial\Omega \cap \cA_{\alpha}\cap \cN )\times S_0)\setminus G( ( \partial\Omega \cap  \overline{\cA_{\mu_0\alpha}} \cap \cN )\times S_0)\\
\subset\ & \partial\widetilde{M}_{j_0} \cap G( ( \cA_{\alpha}\cap \cN )\times S_0).
\end{split}
\]
From this we obtain the second inclusion in~\eqref{eq:4-strata-smoothing-decomp-far} upon recalling~\eqref{eq:decomp-far-intermediate-inclusion-2}. The proof is complete.
\end{proof}
\begin{rmk}\label{rmk:Psi-image-rmk}
Repeating the argument leading to~\eqref{eq:Psi-image-for-decomp-far}, we see that
\begin{equation}\label{eq:Psi-image-rmk-subset}
\begin{split}
\Psi_{\{j_0, i_{\lambda}\}}\big((\mu_0\alpha, \mu_5\alpha]\times \partial\Omega_{\{j_0, i_{\lambda}\}}\big)\subset\ & V^{3}_{\{j_0, i_{\lambda}\}} \cap (\cA_{2\alpha}\setminus \overline{\cA_{\mu_0\alpha}})  \\
\subset\ & \partial \Omega \cap (\cA_{2\alpha}\setminus \overline{\cA_{\mu_0\alpha}}).
\end{split}
\end{equation}
Another consequence of~\eqref{eq:nested-smoothing-1},~\eqref{eq:nested-smoothing-2}, and~\eqref{eq:distance-attained} we note here is that
\begin{equation}\label{eq:Psi-rmk-supset-before-union}
\begin{split}
\Psi_{\{j_0, i_{\lambda}\}}( (\mu_0\alpha, \mu_1\alpha) \times \partial\Omega_{\{j_0, i_{\lambda}\}} ) \supset\ &
(\cA_{(\tau + 2\delta)\alpha} \setminus \cA_{(\tau - 2\delta)\alpha}) \cap \Omega_{\{j_0, i_{\lambda}\}}\\
=\ & (\cA_{(\tau + 2\delta)\alpha} \setminus \cA_{(\tau - 2\delta)\alpha}) \cap V^{3}_{\{j_0, i_{\lambda}\}}.
\end{split}
\end{equation}
The second line follows from Proposition~\ref{prop:distance-neighborhood-in-smoothing}(b), since $\tau - 2\delta > \tau_0$. Taking the union over $\lambda \in \{0, 1, 2\}$ and using~\eqref{eq:Omega-agree-specialized} leads to
\begin{equation}\label{eq:Psi-image-rmk-supset}
\begin{split}
\cup_{\lambda \in \{0, 1, 2\}}\Psi_{\{j_0, i_{\lambda}\}}( (\mu_0\alpha, \mu_1\alpha) \times \partial\Omega_{\{j_0, i_{\lambda}\}} ) \supset\ & (\cA_{(\tau + 2\delta)\alpha} \setminus \cA_{(\tau - 2\delta)\alpha}) \cap \partial V^{3}_{j_0}  \\
=\ & (\cA_{(\tau + 2\delta)\alpha} \setminus \cA_{(\tau - 2\delta)\alpha}) \cap \partial\Omega.
\end{split}
\end{equation}
\end{rmk}
\vskip 1em
To continue, we recall the definition made right before Lemma~\ref{lemm:M-J-*-properties}, and let
\begin{equation}\label{eq:M-J-double-star}
\begin{split}
M_{\{j_0, i_{\lambda}\}}^{**} :=\ & M_{\{j_0, i_{\lambda}\}}^{*} \setminus h_{\{j_0, i_{\lambda}\}}([0, \mu_2\alpha] \times \partial \widetilde{M}_{\{j_0, i_{\lambda}\}}) \\
:=\ & h_{\{j_0, i_{\lambda}\}}( (\mu_{2}\alpha, \mu_5\alpha] \times \partial\widetilde{M}_{\{j_0,i_\lambda\}} ),
\end{split}
\end{equation}
where the second equality uses~\eqref{eq:h-J-image-outside-collar} and~\eqref{eq:h-J-disjoint-images}. Then, from Propositions~\ref{prop:5-strata-smoothing} and~\ref{prop:4-strata-smoothing-decomp-far}, we obtain the following expressions for $\partial\widetilde{M}_{j_0}$: 
\begin{align}
\partial\widetilde{M}_{j_0} =\ & G(( \cN \cap \cA_{(\tau - \delta)\alpha} \cap \partial\Omega) \times S_0) \nonumber\\
&  \cup \big(\cup_{\lambda \in \{0, 1, 2\}} h_{\{j_0, i_{\lambda}\}}\big((\mu_0 \alpha, \mu_5\alpha] \times \partial\widetilde{M}_{\{j_0, i_{\lambda}\}}\big) \big) \label{eq:two-part-decomp}\\
=\ & G(( \mathring{\cN}_{2/3} \cap \cA_{(\tau - \delta)\alpha} \cap \partial\Omega) \times S_0)\nonumber\\
& \cup \big( \cup_{\lambda \in \{0, 1, 2\}} f_{\{i_\lambda\}^{c}}^{-1}( (B_{(\tau - \delta)\alpha} \times \partial\Omega_{\lambda}') \times M_{\{i_{\lambda}\}^{c}}^{**} )\big)\nonumber\\
& \cup \big(\cup_{\lambda \in \{0, 1, 2\}} h_{\{j_0, i_{\lambda}\}}\big((\mu_0 \alpha, \mu_4\alpha) \times \partial\widetilde{M}_{\{j_0, i_{\lambda}\}}\big) \big)\nonumber\\
& \cup \big( \cup_{\lambda \in \{0, 1,2\}} M_{\{j_0, i_{\lambda}\}}^{**} \big).\label{eq:four-part-decomp}
\end{align}
It follows from Proposition~\ref{prop:5-strata-smoothing} that $G(( \mathring{\cN}_{2/3} \cap \cA_{(\tau - \delta)\alpha} \cap \partial\Omega) \times S_0)$ and each $f_{\{i_\lambda\}^{c}}^{-1}( (B_{(\tau - \delta)\alpha} \times \partial\Omega_{\lambda}') \times M_{\{i_{\lambda}\}^{c}}^{**} )$ are open in $\partial\widetilde{M}_{j_0}$. That the same is true for the remaining pieces in~\eqref{eq:four-part-decomp} is the content of the next result.
\begin{lemm}\label{lemm:relative-openness}
For each $\lambda \in \{0, 1, 2\}$, the sets $M_{\{j_0, i_{\lambda}\}}^{**}$ and $h_{\{j_0, i_{\lambda}\}}((\mu_0 \alpha, \mu_4\alpha) \times \partial\widetilde{M}_{\{j_0, i_{\lambda}\}})$ are open relative to $\partial\widetilde{M}_{j_0}$.
\end{lemm}
\vskip 1mm
\begin{proof}
By the second inclusion in~\eqref{eq:4-strata-smoothing-decomp-far}, the restriction $h_{\{j_0, i_{\lambda}\}}|_{(\mu_0\alpha, \mu_4\alpha) \times \partial\widetilde{M}_{\{j_0, i_{\lambda}\}}}$ takes values in $\partial \widetilde{M}_{j_0}$, and thus is smooth as a map into the latter. Recalling Proposition~\ref{prop:g-J-properties}(c) and property (4) below~\eqref{eq:collar-pushout-4-strata}, we deduce that $h_{\{j_0, i_{\lambda}\}}$ is an injective local diffeomorphism from $(\mu_0\alpha, \mu_4\alpha) \times \partial\widetilde{M}_{\{j_0, i_{\lambda}\}}$ into $\partial\widetilde{M}_{j_0}$, so in particular its image is an open subset of $\partial\widetilde{M}_{j_0}$.

Turning to $M_{\{j_0, i_{\lambda}\}}^{**}$, from the first inclusion in~\eqref{eq:decomp-far-intermediate-inclusion-2}, followed by~\eqref{eq:decomp-far-intermediate-inclusion-1} and the expression~\eqref{eq:h-J-image-outside-collar}, we have
\begin{equation}\label{eq:M-J-**-M-J-*-relation}
\begin{split}
M_{\{j_0, i_{\lambda}\}}^{**} \subset\ & \widetilde{M}_{\{j_0, i_{\lambda}\}} \setminus G( (\cN \cap \overline{\cA_{(\tau-2\delta)\alpha}}) \times S_0)\\
=\ & M_{\{j_0, i_{\lambda}\}}^{*} \setminus G( (\cN \cap \overline{\cA_{(\tau-2\delta)\alpha}}) \times S_0),
\end{split}
\end{equation}
Noting from the first line of~\eqref{eq:M-J-double-star} that $M_{\{j_0, i_{\lambda}\}}^{**}$ is open in $M_{\{j_0, i_{\lambda}\}}^{*}$, we see from the above that it is also open in $\widetilde{M}_{\{j_0, i_{\lambda}\}} \setminus G( (\cN \cap \overline{\cA_{(\tau-2\delta)\alpha}}) \times S_0)$. This latter set, in turn, is open in $\partial\widetilde{M}_{j_0} \setminus G(( \cN \cap \overline{\cA_{(\tau - 2\delta)\alpha}} ) \times S_0)$ since, by the equality in~\eqref{eq:M-J-**-M-J-*-relation}, and Lemma~\ref{lemm:M-J-*-properties}, the three sets in the union on the last line of~\eqref{eq:decomp-far-prelim} have pairwise disjoint closures. Having thus shown that $M_{\{j_0, i_{\lambda}\}}^{**}$ is open in $\partial\widetilde{M}_{j_0} \setminus G(( \cN \cap \overline{\cA_{(\tau - 2\delta)\alpha}} ) \times S_0)$, we conclude the proof upon recalling that $G(( \cN \cap \overline{\cA_{(\tau - 2\delta)\alpha}} ) \times S_0)$ is compact by~\eqref{eq:N-A-closure-intersection}.\\
\end{proof}

Letting $\alpha'$ be the threshold in~\eqref{eq:alpha'-threshold}, our next task is to define a smooth map from $(-\alpha', \alpha') \times \partial\widetilde{M}_{j_0}$ into $M$. This will be done by first defining $g_{j_0}^{(1)}(s, \cdot), \cdots, g_{j_0}^{(4)}(s, \cdot)$ on the four regions on the right-hand side of~\eqref{eq:four-part-decomp}, and then showing that they patch up smoothly. 
\vskip 1mm
\begin{defi}\label{defi:g-j-definition}
Given $(s, q) \in (-\alpha', \alpha') \times \partial\widetilde{M}_{j_0}$, notice by~\eqref{eq:alpha'-threshold} and~\eqref{eq:N'-height-threshold} that 
\[
4\alpha' < h < \min\{\frac{\delta\alpha}{4}, \frac{\mu_0\alpha}{4}\},
\]
and consider the following four cases. 
\begin{enumerate}
\item[(1)] In the case $q \in G(( \mathring{\cN}_{2/3} \cap \cA_{(\tau - \delta)\alpha} \cap \partial\Omega) \times S_0)$, we write $(y, p) = \big(G|_{\mathring{\cN}_{7/8} \times S_0}\big)^{-1}(q)$. Then $y + s\xi(y) \in \mathring{\cN}_{2/3} \cap \cA_{\tau\alpha}$ by Lemma~\ref{lemm:N-preserved-by-flow} and the triangle inequality, and we define
\begin{equation}\label{eq:g-j-definition-case-1}
g_{j_0}^{(1)}(s, q) = G(y + s\xi(y), p).
\end{equation}
\vskip 1mm
\item[(2)] In the case $q \in f_{\{i_\lambda\}^{c}}^{-1}( (B_{(\tau - \delta)\alpha} \cap \partial\Omega_{\lambda}') \times M_{\{i_{\lambda}\}^{c}}^{**} )$ for some $\lambda \in \{0, 1, 2\}$, necessarily unique by~\eqref{eq:U-disjoint}, we write $(y', p) = f_{\{i_{\lambda}\}^{c}}(q)$. Then, since $\xi_{\lambda}'$ takes values in $V^{2, \{i_{\lambda}\}^{c}}$, we have $y' + s\xi_{\lambda}'(y') \in B^{2, \{i_{\lambda}\}^{c}}_{\tau\alpha}$, and can thus define
\begin{equation}\label{eq:g-j-definition-case-2}
g_{j_0}^{(2)}(s, q) = f_{\{i_{\lambda}\}^{c}}^{-1}(y' + s\xi_{\lambda}'(y'), p).
\end{equation}
\vskip 1mm
\item[(3)] If $q \in  h_{\{j_0, i_{\lambda}\}}\big((\mu_0 \alpha, \mu_4\alpha) \times \partial\widetilde{M}_{\{j_0, i_{\lambda}\}}\big)$ for some $\lambda \in \{0, 1, 2\}$, necessarily unique by~\eqref{eq:h-J-image-outside-collar} and~\eqref{eq:U-mutually-disjoint}, we write $(t, p) = \big(h_{\{j_0, i_{\lambda}\}}|_{(0, \mu_4\alpha) \times \partial\widetilde{M}_{\{j_0, i_{\lambda}\}}}\big)^{-1}(q)$, and define
\begin{equation}\label{eq:g-j-definition-case-3}
g_{j_0}^{(3)}(s, q) = H_{\{j_0, i_{\lambda}\}}(s\bw_{\lambda}, t, p),
\end{equation}
where $\bw_{\lambda}$ is defined by~\eqref{eq:w-lambda-definition}, and belongs to $V^{1, \{j_0, i_{\lambda}\}}$.
\vskip 1mm
\item[(4)] If $q \in M_{\{j_0, i_\lambda\}}^{**}$ for some $\lambda \in \{0, 1, 2\}$, then since $M_{\{j_0, i_\lambda\}}^{**} \subset M_{\{j_0, i_\lambda\}}^{*}$ and $\bw_{\lambda} \in V^{1, \{j_0, i_{\lambda}\}}$, it makes sense to define
\begin{equation}\label{eq:g-j-definition-case-4}
g_{j_0}^{(4)}(s, q) = f_{\{j_0, i_{\lambda}\}}^{-1}(s\bw_{\lambda}, q).
\end{equation}
\end{enumerate}
\end{defi}
\begin{lemm}\label{lemm:g-j-target}
For $l = 1, 2, 3, 4$, each $g_{j_0}^{(l)}$ is a diffeomorphism on its domain, and satisfies $g_{j_0}^{(l)}(0, q) = q$. Moreover $g_{j_0}^{(l)}(s, q) \in \widetilde{M}_{j_0}$ if and only if $s \geq 0$.
\end{lemm}
\begin{proof}
The conclusion regarding $g_{j_0}^{(l)}(0, q)$ follows from inspecting the definitions made above. Below we focus on the other two conclusions, and shall use directly the notation in each of the above cases without further explanation. 

For case (1), by the inclusion mentioned right before~\eqref{eq:g-j-definition-case-1} and the properties of the map~\eqref{eq:Psi-specialized}, it is routine to check that $g^{(1)}_{j_0}$ is a composition of diffeomorphisms. The said inclusion, along with Proposition~\ref{prop:5-strata-smoothing}(c), also shows that $g_{j_0}^{(1)}(s, q) \in \widetilde{M}_{j_0}$ if and only if 
\[
G(y + s\xi(y), p) \in G(( \mathring{\cN}_{2/3} \cap \cA_{\tau\alpha} \cap \Omega )\times S_0 ).
\]
By the injectivity of $G$ on $\mathring{\cN}_{2/3} \times S_0$, and the fact that $y + s\xi(y) \in \mathring{\cN}_{2/3} \cap \cA_{\tau\alpha}$ to start with, the above occurs exactly when $y + s\xi(y) \in \Omega$, which is equivalent to $s \geq 0$ by Remark~\ref{rmk:Psi-distance-bounds}.

For case (2), with the help of Proposition~\ref{prop:smoothing-product-structure}(b), we see that $g^{(2)}_{j_0}$ is a composition of diffeomorphisms. By Proposition~\ref{prop:5-strata-smoothing}(b) along with the injectivity of $f_{\{i_{\lambda}\}^{c}}^{-1}$, we see that $g_{j_0}^{(2)}(s, q) \in \widetilde{M}_{j_0}$ if and only if $y' + s\xi_{\lambda}'(y') \in \Omega_{\lambda}'$. By Proposition~\ref{prop:smoothing-product-structure}(b), this last inclusion holds exactly when $s \geq 0$.

For case (3), from Proposition~\ref{prop:H-J-properties}(a) we infer that $g_{j_0}^{(3)}$ is a diffeomorphism. Then, by~\eqref{eq:H-J-image-property-2}, we have
\[
g_{j_0}^{(3)}(s, q) \in M \setminus G( (\cN \cap \overline{\cA_{\mu_0\alpha/4}}) \times S_0),
\]
which together with Proposition~\ref{prop:5-strata-smoothing}(a) shows that, in this case, $g_{j_0}^{(3)}(s, q) \in \widetilde{M}_{j_0}$ is equivalent to $g_{j_0}^{(3)}(s, q) \in M_{j_0}$, which by Proposition~\ref{prop:H-J-properties}(d) happens exactly when $s\bw_{\lambda} \in V^{1, \{j_0, i_{\lambda}\}}_{j_0}$, namely, when $s \geq 0$.

For case (4), that $g_{j_0}^{(4)}$ is a diffeomorphism is clear, since $M_{\{j_0, i_\lambda\}}^{**}$ is open in $M_{\{j_0, i_\lambda\}}^{*}$. Next we observe by~\eqref{eq:U-disjoint-distance} that $g_{j_0}^{(4)}(s, q) \not\in  G((\cN \cap \overline{\cA_{\mu_0\alpha/4}})\times S_0)$, and argue as in case (3) to show that $g_{j_0}^{(4)}(s, q) \in \widetilde{M}_{j_0}$ exactly when $s \geq 0$. The details of this last step are omitted.
\end{proof}
For the next statement, we define
\begin{equation}\label{eq:cD-1-defi}
\cD_{1} := \big(\cA_{(\tau-\delta)\alpha} \cap \partial \Omega\big) \cup \big( \cup_{\lambda \in\{0, 1, 2\}}\Psi_{\{j_0, i_{\lambda}\}}((\mu_0\alpha, \mu_1\alpha) \times \partial\Omega_{\{j_0, i_{\lambda}\}}) \big).
\end{equation}
With the help of Lemma~\ref{lemm:preservation-by-flow} and~\eqref{eq:4-strata-smoothing-inclusion-stage-1}, we have
\begin{equation}\label{eq:cD-1-image}
\begin{split}
G((\cN \cap \cD_{1}) \times S_0) =\ & G( (\cN \cap \cA_{(\tau - \delta)\alpha} \cap \partial\Omega) \times S_0) \\
& \cup \big( \cup_{\lambda \in \{0, 1, 2\}}h_{\{j_0, i_\lambda\}}((\mu_0\alpha, \mu_1\alpha) \times \partial\widetilde{M}_{\{j_0, i_\lambda\}}) \big),
\end{split}
\end{equation}
Also, by~\eqref{eq:Psi-image-for-decomp-far} we have
\begin{equation}\label{eq:cD-inclusion}
\cD_{1} \subset \cA_{\alpha} \cap \partial \Omega. 
\end{equation}
Finally, using~\eqref{eq:Psi-image-rmk-supset}, we see that in fact
\begin{equation}\label{eq:cD-flexible}
\cD_{1} =  \big(\cA_{r} \cap \partial \Omega\big) \cup \big( \cup_{\lambda \in\{0, 1, 2\}}\Psi_{\{j_0, i_{\lambda}\}}((\mu_0\alpha, \mu_1\alpha) \times \partial\Omega_{\{j_0, i_{\lambda}\}}) \big),
\end{equation}
for all $r \in [(\tau-2\delta)\alpha, (\tau+2\delta)\alpha]$.
\begin{prop}\label{prop:g-j-patching}
$g_{j_0}^{(1)}, \cdots, g_{j_0}^{(4)}$ patch together to give a smooth local diffeomorphism 
\[
g_{j_0}:(-\alpha', \alpha') \times \partial\widetilde{M}_{j_0} \to M.
\]
Moreover, with $\cD_{1}$ as defined in~\eqref{eq:cD-1-defi}, we have
\begin{equation}\label{eq:g-j-expression-in-G}
g_{j_0}(s, G(y, p)) = G(s\xi(y) + y, p),
\end{equation}
for all $s \in (-\alpha', \alpha')$, $y \in \cD_{1} \cap \cN$, and $p \in S_0$. On the other hand, there holds
\begin{equation}\label{eq:g-j-expression-in-H}
g_{j_0}(s, h_{\{j_0, i_{\lambda}\}}(t, p)) = H_{\{j_{0}, i_{\lambda}\}}(s\bw_{\lambda}, t, p),
\end{equation}
for all $\lambda \in \{0, 1, 2\}$, $s \in (-\alpha', \alpha')$, $t \in (\mu_0\alpha, \mu_5\alpha]$, and $p \in \partial\widetilde{M}_{\{j_0, i_\lambda\}}$.
\end{prop}
\begin{proof}
We begin with a few observations. Given the following data:
\begin{equation}\label{eq:12-overlap-data}
\lambda \in \{0, 1, 2\},\quad y' \in \partial\Omega_{\lambda}' \cap B_{(\tau - \delta)\alpha},\quad t \in (5/8, 2],\quad p \in S_0,
\end{equation}
we have by~\eqref{eq:N-A-intersection-as-product} that
\begin{equation}\label{eq:case-1-case-2-prelim-y}
y := y' + t\frac{\ba_{i_{\lambda}}}{|\ba_{i_{\lambda}}|} \in (\cN \setminus \cN_{5/8}) \cap \cA_{(\tau - \delta)\alpha} \cap \partial\Omega,
\end{equation}
while from~\eqref{eq:xi-along-section-specialized} we obtain $\xi(y) = \xi_{\lambda}'(y')$, and hence 
\[
s\xi(y) + y = s\xi_{\lambda}'(y') + y' + t\frac{\ba_{i_{\lambda}}}{|\ba_{i_{\lambda}}|} \in C^{\{i_{\lambda}\}}_{\tau\alpha}((5/8, 2]) \subset (\cN \setminus \cN_{5/8}) \cap \cA_{\tau\alpha}.
\]
With the help of~\eqref{eq:G-on-tube}, we thus see that
\begin{align}
G(y, p) =\ & f_{\{i_{\lambda}\}^{c}}^{-1}(y', h_{\{i_{\lambda}\}^{c}}(t, p)),\label{eq:case-1-case-2-prelim-1}\\
G(s\xi(y) + y, p) =\ & f_{\{i_{\lambda}\}^{c}}^{-1}(s\xi_{\lambda}'(y') + y', h_{\{i_{\lambda}\}^{c}}(t, p)).\label{eq:case-1-case-2-prelim-2}
\end{align}
Next, suppose $q$ lies in the intersection of case (1) and case (2). Then, with the help of~\eqref{eq:M_I-**-definition}, we obtain $\lambda$, $y'$, $t$, and $p$ as in~\eqref{eq:12-overlap-data} such that
\[
q = f_{\{i_{\lambda}\}^{c}}^{-1}(y', h_{\{i_{\lambda}\}^{c}}(t, p)).
\]
With $y$ defined as in~\eqref{eq:case-1-case-2-prelim-y}, we have~\eqref{eq:case-1-case-2-prelim-1} and~\eqref{eq:case-1-case-2-prelim-2} available, so in particular $q = G(y, p)$. Since it falls into case (1) by assumption, we infer from Proposition~\ref{prop:1234-with-123}(a) and the inclusion in~\eqref{eq:case-1-case-2-prelim-y} that $y \in \mathring{\cN}_{2/3}$, so $(y, p) = (G|_{\mathring{\cN}_{7/8} \times S_0})^{-1}(q)$, and hence the left-hand side of~\eqref{eq:case-1-case-2-prelim-2} coincides $g_{j_0}^{(1)}(s, q)$. Since $g_{j_0}^{(2)}(s, q)$ is exactly the right-hand side of~\eqref{eq:case-1-case-2-prelim-2}, we get $g_{j_0}^{(1)}(s, q) = g_{j_0}^{(2)}(s, q)$. Recalling also that, by~\eqref{eq:product-structure-for-image-Omega}, the domains for $q$ in cases (1) and (2) union to $G( (\cN \cap \cA_{(\tau - \delta)\alpha} \cap \partial\Omega) \times S_0)$, we obtain a well-defined smooth map
\[
g_{j_0}^{(12)}:(-\alpha', \alpha') \times G( (\cN \cap \cA_{(\tau - \delta)\alpha} \cap \partial\Omega) \times S_0) \to M
\]
upon requiring that 
\[
g_{j_0}^{(12)}(s, q) =\left\{
\begin{array}{ll}
g_{j_0}^{(1)}(s, q), & \text{ if }q \in G( (\mathring{\cN}_{2/3} \cap \cA_{(\tau - \delta)\alpha} \cap \partial\Omega) \times S_0),\\
g_{j_0}^{(2)}(s, q), & \text{ if }q \in G( ((\cN\setminus \cN_{5/8}) \cap \cA_{(\tau - \delta)\alpha} \cap \partial\Omega) \times S_0).
\end{array}
\right.
\]
Notice that 
\begin{equation}\label{eq:g-j-12-expression}
g_{j_0}^{(12)}(s, G(y, p)) = G(s\xi(y) + y, p),
\end{equation}
for all $y \in \cN \cap \cA_{(\tau - \delta)\alpha} \cap \partial\Omega$ and $p \in S_0$. Indeed, in the case $y \in \mathring{\cN}_{2/3}$, this merely restates the definition~\eqref{eq:g-j-definition-case-1}, while if $y \in \cN \setminus\cN_{5/8}$, then by~\eqref{eq:N-A-intersection-as-product} we can express $y$ as in~\eqref{eq:case-1-case-2-prelim-y}, in which case we have by~\eqref{eq:case-1-case-2-prelim-1},~\eqref{eq:g-j-definition-case-2}, and~\eqref{eq:case-1-case-2-prelim-2} that
\[
\begin{split}
g_{j_0}^{(12)}(s, G(y, p)) =\ & g_{j_0}^{(2)}(s, f_{\{i_{\lambda}\}^{c}}^{-1}(y', h_{\{i_{\lambda}\}^{c}}(t, p)))\\
=\ & f_{\{i_{\lambda}\}^{c}}^{-1}(s\xi_{\lambda}'(y') + y', h_{\{i_{\lambda}\}^{c}}(t, p)) = G(s\xi(y) + y, p).
\end{split}
\]
This proves~\eqref{eq:g-j-12-expression}.

Next we turn to the overlap between cases (3) and (4). Notice that, given $\lambda \in \{0, 1, 2\}$, $t \in (\mu_2\alpha, \mu_{5}\alpha]$, and $p \in \partial\widetilde{M}_{\{j_0, i_{\lambda}\}}$, by the definition~\eqref{eq:H-J-definition} of $H_{\{j_0, i_{\lambda}\}}$, we have
\begin{equation}\label{eq:case-3-case-4-prelim}
H_{\{j_0, i_{\lambda}\}}(s\bw_{\lambda}, t, p) = f_{\{j_0, i_{\lambda}\}}^{-1}(s\bw_{\lambda}, h_{\{j_0, i_{\lambda}\}}(t, p)).
\end{equation}
Now suppose $q$ lies in the intersection of case (3) and case (4). Then by~\eqref{eq:h-J-image-outside-collar} and~\eqref{eq:U-mutually-disjoint} we get a common $\lambda \in \{0, 1, 2\}$ so that
\[
q \in  h_{\{j_0, i_{\lambda}\}}((\mu_0 \alpha, \mu_4\alpha) \times \partial\widetilde{M}_{\{j_0, i_{\lambda}\}}) \cap M_{\{j_0, i_{\lambda}\}}^{**},
\]
and~\eqref{eq:h-J-disjoint-images} implies that the pair $(t,p):=  \big(h_{\{j_0, i_{\lambda}\}}|_{(0, \mu_4\alpha)}\big)^{-1}(q)$ satisfies $t > \mu_2\alpha$, in which case we get from~\eqref{eq:case-3-case-4-prelim} that
\[
g_{j_0}^{(3)}(s, q) = g_{j_0}^{(4)}(s, q).
\]
As a result, we obtain a smooth map
\[
g_{j_0}^{(34)}:(-\alpha', \alpha') \times \cup_{\lambda \in \{0, 1, 2\}}h_{\{j_0, i_{\lambda}\}}((\mu_0 \alpha, \mu_5\alpha] \times \partial\widetilde{M}_{\{j_0, i_{\lambda}\}}) \to M
\] 
upon requiring for each $\lambda \in \{0, 1, 2\}$ that 
\[
g_{j_0}^{(34)}(s, \cdot) = \left\{
\begin{array}{ll}
g_{j_0}^{(3)}(s, \cdot), & \text{ on }h_{\{j_0, i_{\lambda}\}}((\mu_0 \alpha, \mu_4\alpha) \times \partial\widetilde{M}_{\{j_0, i_{\lambda}\}}),\\
g_{j_0}^{(4)}(s, \cdot), & \text{ on }M_{\{j_0, i_{\lambda}\}}^{**}.
\end{array}
\right.
\]
With the help of~\eqref{eq:case-3-case-4-prelim}, we have
\begin{equation}\label{eq:g-j-34-expression}
g_{j_0}^{(34)}(s, h_{\{j_0, i_{\lambda}\}}(t, p)) = H_{\{j_{0}, i_{\lambda}\}}(s\bw_{\lambda}, t, p),
\end{equation}
for all $\lambda \in \{0, 1, 2\}$, $t \in (\mu_0\alpha, \mu_5\alpha]$, and $p \in \partial \widetilde{M}_{\{j_0,i_\lambda\}}$.

As a preparation for patching together $g_{j_0}^{(12)}$ and $g_{j_0}^{(34)}$, observe that, given
\begin{equation}\label{eq:1234-overlap-data}
\lambda \in \{0, 1, 2\},\quad t \in (\mu_0\alpha, \mu_1\alpha),\quad z \in \cN \cap \partial\Omega_{\{j_0, i_{\lambda}\}},\quad  p \in S_0,
\end{equation}
by Lemma~\ref{lemm:preservation-by-flow} and the inclusion~\eqref{eq:Psi-image-for-decomp-far} derived in the proof of Proposition~\ref{prop:4-strata-smoothing-decomp-far}, we have
\begin{equation}\label{eq:case-3-y}
\begin{split}
y: = \Psi_{\{j_0,i_{\lambda}\}}(t, z) \in\ & \cN \cap (\cA_{\alpha}\setminus \overline{\cA_{\mu_0\alpha}}) \cap V^{3}_{\{j_0, i_{\lambda}\}} \\
\subset\ & \cN \cap (\cA_{\alpha}\setminus \overline{\cA_{\mu_0\alpha}}) \cap \partial\Omega,
\end{split}
\end{equation}
which along with~\eqref{eq:w-lambda-definition} implies
\begin{equation}\label{eq:case-3-xi}
\xi(y) = \bw_{\lambda}.
\end{equation}
Also, since $t \in (\mu_0\alpha, \mu_1\alpha)$, combining Proposition~\ref{prop:g-J-properties}(a) with the relationship between $h_{\{j_0, i_{\lambda}\}}$ and $g_{\{j_0, i_{\lambda}\}}$ stated below~\eqref{eq:collar-pushout-4-strata}, and respectively the definition~\eqref{eq:H-J-definition} of $H_{\{j_0, i_\lambda\}}$, we obtain the following two relations:
\begin{align}
h_{\{j_0, i_\lambda\}}(t, G(z, p)) =\ & G(y, p), \label{eq:12-34-prelim-1} \\
H_{\{j_0, i_{\lambda}\}}(s\bw_{\lambda}, t, G(z, p)) =\ & G(s\bw_{\lambda} + y, p). \label{eq:12-34-prelim-2}
\end{align}
Now suppose 
\[
q \in G( (\cN \cap \cA_{(\tau - \delta)\alpha} \cap \partial\Omega)\times S_0) \cap \big( \cup_{\lambda \in \{0, 1, 2\}} h_{\{j_0, i_\lambda\}}((\mu_0\alpha, \mu_5\alpha] \times \partial\widetilde{M}_{\{j_0, i_{\lambda}\}})\big).
\]
Then, by~\eqref{eq:empty-intersection-for-decomp-far} from the proof of Proposition~\ref{prop:4-strata-smoothing-decomp-far}, we have $q \in  h_{\{j_0, i_\lambda\}}((\mu_0\alpha, \mu_1\alpha) \times \partial\widetilde{M}_{\{j_0, i_{\lambda}\}})$ for some $\lambda \in \{0, 1, 2\}$. In particular, there are $t$, $z$ and $p$ as in~\eqref{eq:1234-overlap-data} such that
\[
q  = h_{\{j_0,i_\lambda\}}(t, G(z, p)) = G(y, p),
\]
where $y = \Psi_{\{j_0, i_{\lambda}\}}(t, z)$ as above, and the second equality is~\eqref{eq:12-34-prelim-1}. Since $q \in G( (\cN \cap \cA_{(\tau - \delta)\alpha} \cap \partial\Omega)\times S_0)$, we infer from~\eqref{eq:case-3-y} and Lemma~\ref{lemm:1-strata-nbhd-facts}(b) that $y \in \cN \cap \cA_{(\tau - \delta)\alpha} \cap \partial\Omega$. Thus, from~\eqref{eq:g-j-12-expression} and~\eqref{eq:case-3-xi}, followed by~\eqref{eq:12-34-prelim-2} and~\eqref{eq:g-j-34-expression}, we get
\[
\begin{split}
g_{j_0}^{(12)}(s, q) = G(s\bw_{\lambda} + y, p) = H_{\{j_0, i_\lambda\}}(s\bw_{\lambda}, t, G(z, p)) = g_{j_0}^{(34)}(s, q).
\end{split}
\]
In view of~\eqref{eq:two-part-decomp}, we obtain a smooth map $g_{j_0}: (-\alpha', \alpha') \times \partial\widetilde{M}_{j_0} \to M$ upon requiring that 
\[
g_{j_0}(s, q) = \left\{
\begin{array}{ll}
g_{j_0}^{(12)}(s, q), & \text{ if } q \in G( (\cN \cap \cA_{(\tau - \delta)\alpha} \cap \partial\Omega)\times S_0),\\
g_{j_0}^{(34)}(s, q), & \text{ if }q \in \cup_{\lambda \in \{0, 1, 2\}} h_{\{j_0, i_\lambda\}}((\mu_0\alpha, \mu_5\alpha] \times \partial\widetilde{M}_{\{j_0, i_{\lambda}\}}). 
\end{array}
\right.
\]
Since $g_{j_0}^{(1)}, \cdots, g_{j_0}^{(4)}$ are diffeomorphisms on their respective domains by Lemma~\ref{lemm:g-j-target}, and since each of these domains is open in $\partial\widetilde{M}_{j_0}$ by Lemma~\ref{lemm:relative-openness} and the comments prior to it, we conclude $g_{j_0}$ is a local diffeomorphism. 

For the remaining two conclusions of the proposition, notice that~\eqref{eq:g-j-expression-in-H} is just a restatement of~\eqref{eq:g-j-34-expression}. As for~\eqref{eq:g-j-expression-in-G}, observe first that, given $s$, $y$, and $p$ as in the statement, both sides of the equation makes sense. For the left-hand side, this is because by~\eqref{eq:cD-1-image} and~\eqref{eq:two-part-decomp}, we have $G(y, p) \in \partial \widetilde{M}_{j_0}$, while for the right-hand side, we use~\eqref{eq:cD-inclusion},  Lemma~\ref{lemm:N-preserved-by-flow}, and the triangle inequality to see that $y + s\xi(y) \in \cN \cap \cA_{2\alpha}$. Next, from~\eqref{eq:g-j-12-expression} we see that~\eqref{eq:g-j-expression-in-G} holds when $y \in \cN \cap \cA_{(\tau -\delta)\alpha} \cap \partial\Omega$. By Lemma~\ref{lemm:preservation-by-flow} and the definition of $\cD_{1}$, the remaining case is 
\[
y \in \cup_{\lambda \in \{0, 1, 2\}} \Psi_{\{j_0, i_\lambda\}}((\mu_0\alpha, \mu_1\alpha) \times (\partial\Omega_{\{j_0,i_\lambda\}} \cap \cN) ).
\]
In other words, there exist $\lambda$, $t$, and $z$ as in~\eqref{eq:1234-overlap-data}, such that $y$ can be expressed as in~\eqref{eq:case-3-y}. As a result we again have $\xi(y) = \bw_{\lambda}$. Using~\eqref{eq:12-34-prelim-1},~\eqref{eq:g-j-34-expression}, and~\eqref{eq:12-34-prelim-2}, we see that
\[
\begin{split}
g_{j_0}(s, G(y, p)) = \ & g_{j_0}^{(34)}(s, h_{\{j_0, i_\lambda\}}(t, G(z,p)))\\
=\ & H_{\{j_0, i_\lambda\}}(s\bw_{\lambda}, t, G(z, p)) = G(s\xi(y) + y, p).
\end{split}
\]
This proves~\eqref{eq:g-j-expression-in-G}, and we are done.
\end{proof}

\begin{prop}\label{prop:g-j-injective}
The map $g_{j_0}$ is injective on $(-\alpha', \alpha') \times \partial\widetilde{M}_{j_0}$, and thus a global diffeomorphism.
\end{prop}
\begin{proof}
Define 
\[
\begin{split}
E_1 :=\ & G((\cN \cap \cA_{(\tau - \delta)\alpha} \cap \partial\Omega) \times S_0)\\
& \cup \big( \cup_{\lambda \in \{0, 1, 2\}} h_{\{j_0, i_{\lambda}\}}((\mu_0\alpha, \mu_1\alpha) \times \partial\widetilde{M}_{\{j_0, i_{\lambda}\}})\big),\\
E_2:=\ & \big(\cup_{\lambda \in \{0, 1, 2\}}h_{\{j_0, i_\lambda\}}((\mu_0\alpha, \mu_{4}\alpha) \times \partial\widetilde{M}_{\{j_0, i_\lambda\}})\big) \cup \big(\cup_{\lambda \in \{0, 1, 2\}}M_{\{j_0, i_{\lambda}\}}^{**}\big).
\end{split}
\]
By~\eqref{eq:two-part-decomp}, we of course have $E_1 \cup E_2  = \partial\widetilde{M}_{j_0}$, while by~\eqref{eq:empty-intersection-for-decomp-far},~\eqref{eq:h-J-disjoint-images}, and~\eqref{eq:U-mutually-disjoint}, we have
\begin{equation}\label{eq:E1-intersect-E2}
E_1 \cap E_2 =  \cup_{\lambda \in \{0, 1, 2\}} h_{\{j_0, i_{\lambda}\}}((\mu_0\alpha, \mu_1\alpha) \times \partial\widetilde{M}_{\{j_0, i_{\lambda}\}}).
\end{equation}
Recall also that, by~\eqref{eq:cD-1-image},
\[
E_1 = G(( \cN \cap \cD_{1}) \times S_0).
\]
We divide the proof of Proposition~\ref{prop:g-j-injective} into the verification of a number of claims.
\begin{claim}\label{claim:g-j-injective-on-E1}
$g_{j_0}$ is injective on $(-\alpha', \alpha') \times E_1$.
\end{claim}
\begin{proof}
We further split $E_{1}$ as
\[
E_1 = G(( \cN_{5/8} \cap \cD_{1}) \times S_0) \cup G(( (\cN \setminus \cN_{5/8}) \cap \cD_{1}) \times S_0) =: E_{11} \cup E_{12}.
\]
Take a pair of points $(s_1, q_1), (s_2, q_2) \in (-\alpha', \alpha') \times E_{11}$ such that 
\[
g_{j_0}(s_1, q_1) = g_{j_0}(s_2, q_2).
\]
By~\eqref{eq:g-j-expression-in-G}, upon expressing each $q_k$ as $G(y_k, p_k)$, with $y_k \in \cN_{5/8}\cap \cD_{1}$ and $p_k \in S_0$, the above translates into
\[
G(s_1\xi(y_1) + y_1, p_1) = G(s_2\xi(y_2) + y_2, p_2).
\]
By~\eqref{eq:cD-inclusion}, both $y_{1}$ and $y_{2}$ lie in $\cN_{5/8} \cap \cA_{\alpha} \cap \partial\Omega$, so we infer from Lemma~\ref{lemm:N-preserved-by-flow} that
\[
s_1\xi(y_1) + y_1,\ s_2\xi(y_2) + y_2 \in \cN_{5/8}.
\]
Recalling that $G$ and the map $(s, y) \mapsto y + s\xi(y)$ are injective on $\mathring{\cN}_{7/8} \times S_0$ and, respectively, $(-\alpha', \alpha') \times \partial\Omega$, we deduce that $s_1 = s_2$, $y_1 = y_2$, and $p_1 = p_2$. It follows that $q_{1} = q_{2}$ as well. 

Next, take $(s, q), (\tilde{s}, \tilde{q}) \in (-\alpha', \alpha') \times E_{12}$ such that 
\[
g_{j_0}(s, q) = g_{j_0}(\tilde{s}, \tilde{q}).
\]
Writing $q = G(y, p)$ for some $y \in (\cN \setminus \cN_{5/8}) \cap \cD_{1}$ and $p \in S_0$, we obtain from~\eqref{eq:cD-inclusion} and~\eqref{eq:N-A-intersection-as-product} some $\lambda \in \{0, 1, 2\}$, $z \in B_{\alpha} \cap \partial\Omega_{\lambda}'$, and $t \in (5/8, 2]$ such that 
\[
y = z + t\frac{\ba_{i_\lambda}}{|\ba_{i_\lambda}|},\quad\quad \xi(y) = \xi_{\lambda}'(z).
\]
Consequently, with the help of~\eqref{eq:G-on-tube} and~\eqref{eq:g-j-expression-in-G}, we get
\begin{equation}\label{eq:case-2-3-injective-q}
q  = f_{\{i_{\lambda}\}^{c}}^{-1}(z, h_{\{i_\lambda\}^{c}}(t, p)),\quad\quad g_{j_0}(s, q) = f_{\{i_{\lambda}\}^{c}}^{-1}(z + s \xi'_{\lambda}(z), h_{\{i_\lambda\}^{c}}(t, p)).
\end{equation}
Similar expressions hold for $\widetilde{q}$ and $g_{j_0}(\tilde{s}, \tilde{q})$, with the corresponding ingredients decorated with tildes ($\tilde{\ }$). The assumption $g_{j_0}(s, q) = g_{j_0}(\tilde{s}, \tilde{q})$ then becomes
\[
f_{\{i_{\lambda}\}^{c}}^{-1}(z + s \xi'_{\lambda}(z), h_{\{i_\lambda\}^{c}}(t, p)) = f_{\{i_{\tilde{\lambda}}\}^{c}}^{-1}(\tilde{z} + \tilde{s}\cdot \xi'_{\tilde{\lambda}}(\tilde{z}), h_{\{i_{\tilde{\lambda}}\}^{c}} (\tilde{t}, \tilde{p})).
\]
The disjointness property~\eqref{eq:U-disjoint} forces $\lambda  = \tilde{\lambda}$, in which case we get
\begin{equation}\label{eq:case-2-3-injective-prelim}
h_{\{i_\lambda\}^{c}}(t, p) = h_{\{i_{\lambda}\}^{c}} (\tilde{t}, \tilde{p}),\quad
z + s \xi'_{\lambda}(z) = \tilde{z} + \tilde{s}\cdot \xi'_{\lambda}(\tilde{z}).
\end{equation}
Adding $\rho_2\frac{\ba_{i_{\lambda}}}{|\ba_{i_{\lambda}}|}$ to both sides of the second equation and using~\eqref{eq:xi-along-section-specialized} gives
\[
\Psi(s, z + \rho_2\frac{\ba_{i_{\lambda}}}{|\ba_{i_{\lambda}}|} ) = \Psi(\tilde{s}, \tilde{z} + \rho_2\frac{\ba_{i_{\lambda}}}{|\ba_{i_{\lambda}}|}),
\]
from which we get $s = \tilde{s}$ and $z = \tilde{z}$. Substituting the latter along with $\lambda = \tilde{\lambda}$ and the first equation in~\eqref{eq:case-2-3-injective-prelim} back into the expression~\eqref{eq:case-2-3-injective-q} and the analogous expression for $\tilde{q}$, we get $q = \tilde{q}$ also. 

Having shown that $g_{j_0}$ restricts to be injective on $(-\alpha', \alpha') \times E_{11}$ and $(-\alpha', \alpha') \times E_{12}$, we next notice that, by~\eqref{eq:g-j-expression-in-G}, along with~\eqref{eq:cD-inclusion} and Lemma~\ref{lemm:N-preserved-by-flow},
\[
g_{j_0}((-\alpha', \alpha') \times E_{11}) \subset G(\cN_{5/8} \times S_0),\quad\quad
g_{j_0}((-\alpha', \alpha') \times E_{12}) \subset G((\cN \setminus \cN_{5/8}) \times S_0),
\]
and thus Proposition~\ref{prop:1234-with-123}(a) shows that $(-\alpha', \alpha') \times E_{11}$ and $(-\alpha', \alpha') \times E_{12}$ have disjoint images under $g_{j_0}$. By what is said at the start of this paragraph, we deduce that in fact $g_{j_0}$ is injective on $(-\alpha', \alpha') \times E_{1}$. This proves the claim.
\end{proof}
\begin{claim}\label{claim:g-j-injective-on-E2}
$g_{j_0}$ is injective on $(-\alpha', \alpha') \times E_{2}$.
\end{claim}
\begin{proof}
Recalling~\eqref{eq:M-J-double-star}, have
\begin{equation}\label{eq:case-3-4-domain}
\begin{split}
E_2 =  \cup_{\lambda\in \{0 ,1, 2\}} h_{\{j_0, i_\lambda\}}((\mu_0\alpha, \mu_{5}\alpha] \times \partial\widetilde{M}_{\{j_0, i_\lambda\}}).
\end{split}
\end{equation}
By~\eqref{eq:g-j-expression-in-H}, as well as~\eqref{eq:H-J-image-property-1} and~\eqref{eq:U-mutually-disjoint}, the sets $(-\alpha', \alpha')\times h_{\{j_0, i_\lambda\}}((\mu_0\alpha, \mu_{5}\alpha] \times \partial\widetilde{M}_{\{j_0, i_\lambda\}})$ for $\lambda \in \{0, 1, 2\}$ have mutually disjoint images under $g_{j_0}$, so it suffices to show that $g_{j_0}$ is injective on each of them in order to prove the claim. Fix $\lambda \in \{0, 1, 2\}$ and take, for $ k \in \{1, 2\}$, some $(s_k, t_k, p_k) \in (-\alpha', \alpha') \times (\mu_0\alpha, \mu_5\alpha] \times \partial\widetilde{M}_{\{j_0, i_{\lambda}\}}$ such that
\begin{equation}\label{eq:case-3-4-injective-assumption}
H_{\{j_0, i_\lambda\}}(s_1\bw_{\lambda}, t_1, p_1) = H_{\{j_0, i_\lambda\}}(s_2\bw_{\lambda}, t_2, p_2),
\end{equation}
our goal is to deduce that
\begin{equation}\label{eq:case-3-4-injective-goal}
s_1 = s_2, \quad h_{\{j_0, i_{\lambda}\}}(t_1, p_1) = h_{\{j_0, i_{\lambda}\}}(t_2, p_2).
\end{equation}
To begin, by the second part of Proposition~\ref{prop:H-J-properties}(a), we infer from~\eqref{eq:case-3-4-injective-assumption} that $t_1, t_2$ are either both in $(\mu_0\alpha, \mu_2\alpha]$, or both in $(\mu_2\alpha, \mu_5\alpha]$. In the former case, the first part of Proposition~\ref{prop:H-J-properties}(a) gives $(s_1, t_1, p_1) = (s_2, t_2, p_2)$, and thus~\eqref{eq:case-3-4-injective-goal} holds as desired. In the case where $t_1, t_2 \in (\mu_2\alpha, \mu_5\alpha]$, the hypothesis~\eqref{eq:case-3-4-injective-assumption} reduces, by~\eqref{eq:H-J-definition}, to
\[
f_{\{j_0, i_\lambda\}}^{-1}(s_1\bw_{\lambda}, h_{\{j_0, i_{\lambda}\}}(t_1, p_1)) = f_{\{j_0, i_\lambda\}}^{-1}(s_2\bw_{\lambda}, h_{\{j_0, i_{\lambda}\}}(t_2, p_2)),
\]
and we get~\eqref{eq:case-3-4-injective-goal} by the injectivity of $f_{\{j_0, i_\lambda\}}^{-1}$. The proves the claim.
\end{proof}
\begin{claim}\label{claim:E1-diff-E2-disjoint-images}
The sets $(-\alpha', \alpha') \times \big(E_{1} \setminus (E_1 \cap E_2)\big)$ and $(-\alpha', \alpha') \times \big(E_{2} \setminus (E_1 \cap E_2)\big)$ have disjoint images under $g_{j_0}$.
\end{claim}
\begin{proof}
By~\eqref{eq:E1-intersect-E2} we have 
\[
E_1 \setminus (E_1 \cap E_2) \subset G((\cN \cap \cA_{(\tau -\delta)\alpha} \cap \partial\Omega) \times S_0),
\]
so that, by~\eqref{eq:g-j-expression-in-G} together with Lemma~\ref{lemm:N-preserved-by-flow} and the triangle inequality,
\begin{equation}\label{eq:g-j-image-of-E1-setminus}
g_{j_0}((-\alpha', \alpha') \times (E_1 \setminus (E_1 \cap E_2))) \subset G((\cN \cap \cA_{\tau\alpha}) \times S_0).
\end{equation}
On the other hand, by~\eqref{eq:case-3-4-domain} and~\eqref{eq:E1-intersect-E2}, we have
\[
E_{2} \setminus (E_1 \cap E_2) \subset \cup_{\lambda \in \{0, 1, 2\}}h_{\{j_0, i_\lambda\}}([\mu_1\alpha, \mu_5\alpha] \times \partial\widetilde{M}_{\{j_0, i_{\lambda}\}}),
\]
so that, by~\eqref{eq:g-j-expression-in-H},
\begin{equation}\label{eq:g-j-image-of-E2-setminus}
\begin{split}
&g_{j_0}((-\alpha', \alpha') \times (E_2 \setminus (E_1 \cap E_2))) \\
&\subset \cup_{\lambda \in \{0, 1, 2\}}H_{\{j_0, i_{\lambda}\}}(B^{1, \{j_0, i_{\lambda}\}}_{\alpha'} \times [\mu_1\alpha, \mu_5\alpha] \times \partial\widetilde{M}_{\{j_0, i_{\lambda}\}}).
\end{split}
\end{equation}
By~\eqref{eq:H-J-definition} and Lemma~\ref{lemm:images-4-strata}(b), the two sets on the right-hand side of~\eqref{eq:g-j-image-of-E2-setminus} and~\eqref{eq:g-j-image-of-E1-setminus} are disjoint from each other. This gives the desired conclusion.
\end{proof}
Recalling that $E_1 \cup E_2 = \partial\widetilde{M}_{j_0}$, we conclude from Claims~\ref{claim:g-j-injective-on-E1},~\ref{claim:g-j-injective-on-E2} and~\ref{claim:E1-diff-E2-disjoint-images} that $g_{j_0}$ is injective on $(-\alpha', \alpha') \times \partial\widetilde{M}_{j_0}$. The proof of Proposition~\ref{prop:g-j-injective} is complete.\\
\end{proof}

To continue, we define, for $s \leq 4\alpha'$, the distance neighborhood
\begin{equation}\label{eq:distance-neighborhood-2-strata}
\cA^{(2)}_{s} := B_{s}\big( \cup_{|J| = 2}V^{3}_{J} \big) = \cup_{|J| = 2} \big(B^{1, J}_{s} + V^{3}_{J} \big),
\end{equation}
where the second equality follows since the distance to $\cup_{|J| = 2}V^{3}_{J}$ from a point outside of it must be realized on the smooth part, namely $\cup_{|J| = 2}\mathring{V}^{3}_{J}$. It is also elementary to see that, given $j_0 \in \{1, 2, 3, 4\}$, writing $\{j_0\}^{c} = \{i_0, i_1, i_2\}$ as before, we have 
\begin{equation}\label{eq:V-j-disjoint}
V^{3}_{j_0} \cap (V^{1, \{i_{\lambda+1}, i_{\lambda+2}\}} + V^{3}_{\{i_{\lambda+1}, i_{\lambda+2}\}}) = V^{3}_{\{i_{\lambda}\}^{c}},
\end{equation}
for all $\lambda \in \{0, 1, 2\}$, and thus
\begin{equation}\label{eq:distance-2-strata-attained}
\begin{split}
V^{3}_{j_0} \cap \cA_{s}^{(2)} = \ &\cup_{\lambda\in \{0,1, 2\}} V^{3}_{j_0}\cap   \big( B^{1, \{j_0,i_\lambda\}}_{s} + V^{3}_{\{j_0, i_\lambda\}} \big)\\
=\ & V^{3}_{j_0} \cap B_{s}(\partial V^{3}_{j_0}).
\end{split}
\end{equation}
For later use, we derive here an analogue of~\eqref{eq:product-structure-for-image-Omega} based on Lemma~\ref{lemm:sectors-product-structure}(d). Specifically, given $r \in (0, 4\alpha]$ and $s \in (0, 4\alpha']$, by~\eqref{eq:S2-neighborhood-specialized} and~\eqref{eq:distance-2-strata-attained}, followed by~\eqref{eq:tube-intersection}, we have
\[
\begin{split}
(\cN \setminus \cN_{5/8}) \cap \cA_{r}\cap A_{s}^{(2)} \cap V^{3}_{j_0} =\ & \cup_{\lambda \in \{0, 1, 2\}} C^{\{i_{\lambda}\}^{c}}_{r}((5/8, 2]) \cap V^{3}_{j_0} \cap B_{s}(\partial V^{3}_{j_0}).
\end{split}
\]
Noting by~\eqref{eq:rho2-choice} and~\eqref{eq:alpha'-threshold} that $\rho_2 > 4\alpha' \sqrt{\frac{4}{3}} \geq s\sqrt{\frac{4}{3}}$, we deduce with the help of Lemma~\ref{lemm:sectors-product-structure}(d) that 
\[
\begin{split}
&C_{r}^{\{i_{\lambda}\}^{c}}((5/8, 2]) \cap \big( V^{3}_{j_0} \cap B_{s}(\partial V^{3}_{j_0}) \big)\\
& = \big( B_{r} \cap V^{2, \{i_{\lambda}\}^{c}}_{j_0} \cap B_{s}(\partial V^{2, \{i_\lambda\}^{c}}_{j_0})  \big) + (5/8, 2] \cdot \frac{\ba_{i_\lambda}}{|\ba_{i_\lambda}|}.
\end{split}
\]
Recalling the definition~\eqref{eq:G-on-tube} of $G$, we get
\begin{equation}\label{eq:product-structure-with-distance-nbhds}
\begin{split}
&G(((\cN \setminus \cN_{5/8}) \cap \cA_{r}\cap A_{s}^{(2)} \cap V^{3}_{j_0}) \times S_0)\\
=\ & \cup_{\lambda \in \{0, 1, 2\}} f_{\{i_\lambda\}^{c}}^{-1}((B_{r} \cap V^{2, \{i_\lambda\}^{c}}_{j_0} \cap B_{s}(\partial V^{2, \{i_\lambda\}^c}_{j_0})) \times M_{\{i_\lambda\}^{c}}^{**} ).
\end{split}
\end{equation}

Next, we recall from the previous section the domains $\cN'_{R, \ell}$ and $\mathring{\cN}'_{R, \ell}$, given by~\eqref{eq:N'-definition}, as well as the map
\[
G_{(1)}:\cN' \times S_0 \to M
\]
defined by~\eqref{eq:G-1-definition-case-1} and~\eqref{eq:G-1-definition-case-2}. 

\begin{lemm}\label{lemm:Psi-action-on-N'}
Fix $j_0 \in \{1, 2, 3, 4\}$, write $\{j_0\}^{c} =\{i_0, i_1, i_2\}$, and let $\Omega$, $\xi$, and $\Psi$ be the objects mentioned right after~\eqref{eq:alpha'-threshold}. Take also $R \in [3/8, 2]$ and $\ell \in [\mu_{1}\alpha, \mu_{5}\alpha]$.
\vskip 1mm
\begin{enumerate}
\item[(a)] $\partial\Omega \cap \cN'_{R, \ell}$ is a compact set. Moreover, we have 
\begin{equation}\label{eq:boundary-Omega-N'-intersection}
\begin{split}
\partial\Omega \cap \cN'_{R, \ell} =\ & \big( \partial\Omega \cap \cA_{r} \cap \cN_{R} \big) \\
&\cup \big( \cup_{\lambda \in \{0, 1, 2\}}\Psi_{\{j_0,i_{\lambda}\}}((\mu_0\alpha, \ell] \times (\partial\Omega_{\{j_0, i_{\lambda}\}} \cap \cN_{R})) \big),
\end{split}
\end{equation}
for any $r \in [(\tau - 2\delta)\alpha, (\tau + 2\delta)\alpha]$.
\vskip 1mm
\item[(b)] Given $(s, z) \in (-\alpha', \alpha') \times (\partial\Omega \cap \cN'_{R, \ell})$, if $z \in \partial\Omega \cap \cA_{(\tau - \delta)\alpha} \cap \cN_{R}$, then we have 
\[
\Psi(s, z) \in \cA_{\tau\alpha} \cap \cN_{R}.
\]
On the other hand, if $z \in \Psi_{\{j_0,i_{\lambda}\}}((\mu_0\alpha, \ell] \times (\partial\Omega_{\{j_0, i_{\lambda}\}} \cap \cN_{R}))$ for some $\lambda \in \{0, 1, 2\}$, then
\[
\Psi(s, z) = z + s\bw_{\lambda} \in C^{\{j_0,i_{\lambda}\}}_{\alpha', R}((\mu_0\alpha, \ell]).
\]
\vskip 1mm
\item[(c)] Suppose $s \in (-\alpha', \alpha')$ and $z \in \partial\Omega$. Then we have 
\[
\Psi(s, z) \in \cN'_{R, \ell} \quad\text{if and only if}\quad z \in \cN'_{R, \ell}.
\]
\end{enumerate}
\end{lemm}
\begin{proof}
For part (a), we first prove~\eqref{eq:boundary-Omega-N'-intersection}. Since $\ell \geq \mu_{1}\alpha$, we see from~\eqref{eq:Psi-image-rmk-supset} and Lemma~\ref{lemm:preservation-by-flow} that the right-hand side indeed does not depend on the choice of $r \in [(\tau -2\delta)\alpha, (\tau + 2\delta)\alpha]$. Thus we need only establish~\eqref{eq:boundary-Omega-N'-intersection} for $r = \tau\alpha$. To start, by the definition of $\cN'_{R, \ell}$, there holds
\begin{equation}\label{eq:boundary-Omega-N'-intersection-1}
\begin{split}
\partial \Omega \cap \cN'_{R, \ell} = \ &  \big( \partial\Omega \cap \cA_{\tau\alpha} \cap \cN_{R} \big) \cup \big(\cup_{|J| = 2} \partial\Omega \cap C^{J}_{h, R}((\mu_0\alpha, \ell])\big).
\end{split}
\end{equation}
Recalling that $h < \frac{\mu_0\alpha}{4}$, we see from~\eqref{eq:overlap-with-distance-neighborhood-3} and~\eqref{eq:Omega-agree-specialized} that
\[
\begin{split}
C^{J}_{h, R}( (\mu_0\alpha, \ell] ) \cap \partial\Omega=\ & C^{J}_{h, R}( (\mu_0\alpha, \ell] ) \cap \partial V^{3}_{j_0} \\
=\ & \cup_{\lambda \in \{0, 1, 2\}}C^{J}_{h, R}( (\mu_0\alpha, \ell] ) \cap V^{3}_{\{j_0, i_{\lambda}\}}.
\end{split}
\]
Taking the union over $|J| = 2$ and using~\eqref{eq:Phi-J-positive-distance}, we obtain
\[
\begin{split}
\big(\cup_{|J| = 2} \partial\Omega \cap C^{J}_{h, R}((\mu_0\alpha, \ell])\big) =\ & \cup_{\lambda \in \{0, 1, 2\}} C^{\{j_0, i_{\lambda}\}}_{h, R}( (\mu_0\alpha, \ell] ) \cap V^{3}_{\{j_0, i_{\lambda}\}}\\
=\ & \cup_{\lambda \in \{0, 1, 2\}}  \Psi_{\{j_0,i_{\lambda}\}}((\mu_0\alpha, \ell] \times (\partial\Omega_{\{j_0, i_{\lambda}\}} \cap \cN_{R})).
\end{split}
\]
Combining this with~\eqref{eq:boundary-Omega-N'-intersection-1} gives~\eqref{eq:boundary-Omega-N'-intersection}. To see that $\partial\Omega \cap \cN'_{R, \ell}$ is compact, we use Lemma~\ref{lemm:preservation-by-flow},~\eqref{eq:nested-smoothing-1}, and~\eqref{eq:nested-smoothing-2}, and argue as in the proof of~\eqref{eq:Psi-image-for-decomp-far} to see that, for each $\lambda \in \{0, 1, 2\}$,
\begin{equation}\label{eq:Psi-image-for-compact}
\begin{split}
\Psi_{\{j_0, i_{\lambda}\}}( \{\mu_0\alpha\} \times (\partial\Omega_{\{j_0, i_{\lambda}\}} \cap \cN_{R}) ) \subset\ & \cN_{R} \cap (\cA_{\tau\alpha}\setminus \cA_{\mu_{0}\alpha}) \cap V^{3}_{\{j_0, i_{\lambda}\}}\\
\subset\ & \cN_{R} \cap (\cA_{\tau\alpha}\setminus \cA_{\mu_{0}\alpha}) \cap \partial\Omega.
\end{split}
\end{equation}
In view of this and the admissible range of $r$ in~\eqref{eq:boundary-Omega-N'-intersection}, we deduce that
\begin{equation}\label{eq:boundary-Omega-N'-intersection-compact}
\begin{split}
\partial\Omega \cap \cN'_{R, \ell} =\ & \big( \partial\Omega \cap \overline{\cA_{\tau\alpha}}\cap \cN_{R} \big) \\
&\cup \big( \cup_{\lambda \in \{0, 1, 2\}}\Psi_{\{j_0,i_{\lambda}\}}([\mu_0\alpha, \ell] \times (\partial\Omega_{\{j_0, i_{\lambda}\}} \cap \cN_{R})) \big).
\end{split}
\end{equation}
Since both $\overline{\cA_{\tau\alpha}}\cap \cN_{R}$ and $\partial\Omega_{\{j_0, i_{\lambda}\}} \cap \cN_{R}$ are compact, thanks respectively to~\eqref{eq:N-A-closure-intersection} and Remark~\ref{rmk:portions}, we conclude that so is $\partial\Omega \cap \cN'_{R, \ell}$.

For part (b), the case $z \in \partial\Omega \cap \cA_{(\tau -\delta)\alpha}\cap\cN_{R}$ follows from Lemma~\ref{lemm:N-preserved-by-flow}, the triangle inequality, and our choice of $\alpha'$, while the other case uses~\eqref{eq:Psi-image-rmk-subset} and the property~\eqref{eq:w-lambda-definition} of $\xi$.

For part (c), the implication ``$\Leftarrow$'' follows from~\eqref{eq:boundary-Omega-N'-intersection} and part (b). Conversely, suppose by contradiction that $\Psi(s, z) \in \cN'_{R, \ell}$ but $z \not\in \cN'_{R, \ell}$. Since $\cN'_{R, \ell} \subset \cN_{R} \cap \cA_{2\alpha}$ by~\eqref{eq:N'-contained-in-N}, the triangle inequality implies that
\[
z \in \partial\Omega \cap \cA_{2\alpha + \alpha'},
\]
so Lemma~\ref{lemm:N-preserved-by-flow} can be applied, which together with our assumption on $\Psi(s, z)$ and $z$ gives
\begin{equation}\label{eq:N'-preserved-contradiction}
\begin{split}
z \in \partial\Omega \cap \cN_{R} \setminus \cN'_{R, \ell} \subset \ & \partial \Omega \cap \cN_{R} \setminus \cA_{\tau\alpha} =  \cN_{R} \cap \partial V^{3}_{j_0} \setminus \cA_{\tau\alpha},
\end{split}
\end{equation}
where the equality is~\eqref{eq:Omega-agree-specialized}. Now, in the case $\Psi(s, z) \in \cN_{R} \cap \cA_{\tau\alpha}$, the triangle inequality gives additionally that $z \in \cA_{(\tau + \delta)\alpha}$, which together with~\eqref{eq:N'-preserved-contradiction},~\eqref{eq:Psi-image-rmk-supset}, and Lemma~\ref{lemm:preservation-by-flow} yields some $\lambda \in \{0, 1, 2\}$ such that
\[
z \in \Psi_{\{j_0, i_{\lambda}\}}((\mu_0\alpha,\mu_1\alpha) \times (\partial\Omega_{\{j_0, i_\lambda\}} \cap \cN_{R})).
\]
The right-hand side being contained in $\cN'_{R, \ell}$ since $\ell \geq \mu_1\alpha$, we get a contradiction. Next, in the case where $\Psi(s, z) \in C^{J}_{h, R}((\mu_0\alpha, \ell])$ for some $J$, by~\eqref{eq:C-J-E-definition} we can write
\[
\Psi(s, z) = z' + x',
\]
for some $z' \in \Psi_{J}((\mu_0\alpha, \ell] \times (\partial\Omega_{J} \cap \cN_{R}))$ and $x' \in B^{1, J}_{h}$. On the other hand,~\eqref{eq:N'-preserved-contradiction} gives some $\lambda \in \{0, 1, 2\}$ such that $z \in \cN_{R} \cap V^{3}_{\{j_0, i_{\lambda}\}} \setminus \cA_{\tau\alpha}$, in which case~\eqref{eq:w-lambda-definition} implies
\[
\Psi(s, z) = s\bw_{\lambda} + z \in B^{1, \{j_0, i_{\lambda}\}}_{h} + V^{3}_{\{j_0, i_{\lambda}\}}.
\]
The inequality $|z - z'| < 2h$ and the property~\eqref{eq:Phi-J-positive-distance} then forces $J = \{j_0, i_{\lambda}\}$, in which case, upon comparing the two above decompositions for $\Psi(s, z)$, we get
\[
z = z' \in \cN'_{R, \ell},
\]
again a contradiction. The proof is complete.
\end{proof}
\begin{rmk}\label{rmk:N'-ring-version}
By the same argument as in the proof of part (a), we get 
\begin{equation}\label{eq:boundary-Omega-N'-ring-intersection}
\begin{split}
\partial\Omega \cap \mathring{\cN}'_{R, \ell} =\ & \big( \partial\Omega \cap \cA_{r} \cap \mathring{\cN}_{R} \big) \\
&\cup \big( \cup_{\lambda \in \{0, 1, 2\}}\Psi_{\{j_0,i_{\lambda}\}}((\mu_0\alpha, \ell) \times (\partial\Omega_{\{j_0, i_{\lambda}\}} \cap \mathring{\cN}_{R})) \big),
\end{split}
\end{equation}
for any $r \in [(\tau-2\delta)\alpha, (\tau + 2\delta)\alpha]$. Then, following the proof of parts (b) and (c), using the second inclusion in~\eqref{eq:N'-contained-in-N} instead of the first, and noting also that we may replace $\cN_{R}$ with $\mathring{\cN}_{R}$ in both Lemma~\ref{lemm:N-preserved-by-flow} and Lemma~\ref{lemm:preservation-by-flow}, we see that the equivalence in part (c) continues to hold with both occurrences of $\cN'_{R, \ell}$ replaced by $\mathring{\cN}'_{R, \ell}$.
\end{rmk}

\begin{prop}\label{prop:M-j-as-G-1-image}
Let $G_{(1)}:\cN' \times S_0 \to M$ be the map defined by~\eqref{eq:G-1-definition-case-1} and~\eqref{eq:G-1-definition-case-2}. Then, given $j_0 \in \{1, 2, 3, 4\}$, in the notation of Lemma~\ref{lemm:Psi-action-on-N'}, we have 
\begin{equation}\label{eq:M-j-as-G-1-image}
\partial\widetilde{M}_{j_0} = G_{(1)}(( \partial\Omega \cap \cN' ) \times S_0).
\end{equation}
\end{prop}
\begin{proof}
Tracing the definition of $G_{(1)}$, we see in particular that, for $y \in \partial\Omega \cap \cN'$, 
\begin{equation}\label{eq:G-1-definition-special-case}
G_{(1)}(y, p) =\left\{
\begin{array}{ll}
G(y, p),& \text{ if } y \in \cN \cap \cA_{\tau\alpha},\\
h_{\{j_0, i_{\lambda}\}}(s, G(z, p)),& \text{ if }y = \Psi_{\{j_0, i_\lambda\}}(s, z).
\end{array}
\right.
\end{equation}
Here $p$ is any element of $S_0$, and the second case involves no ambiguity since there is at most one way to write $y$ as $\Psi_{\{j_0, i_\lambda\}}(s, z)$ for $\lambda \in \{0, 1, 2\}$, $s \in (\mu_0\alpha, \mu_5\alpha]$, and $z \in \partial\Omega_{\{j_0, i_\lambda\}} \cap \cN$. Applying~\eqref{eq:boundary-Omega-N'-intersection} with $r = (\tau - \delta)\alpha$ and then using~\eqref{eq:G-1-definition-special-case} and~\eqref{eq:smoothing-contained-in-G-image}, we get
\[
\begin{split}
G_{(1)}((  \partial\Omega \cap  \cN') \times S_0)=\ & G((\partial \Omega \cap \cA_{(\tau -\delta)\alpha} \cap \cN) \times S_0) \\
& \cup \big( \cup_{\lambda \in \{0, 1, 2\}} h_{\{j_0, i_\lambda\}}((\mu_0\alpha, \mu_5\alpha] \times \partial\widetilde{M}_{\{j_0, i_\lambda\}}) \big),
\end{split}
\]
and we are done upon recalling~\eqref{eq:two-part-decomp}.
\end{proof}

\begin{lemm}\label{lemm:N'-intersect-2-strata-neighborhood}
Define, for each $J \subset \{1, 2, 3, 4\}$ with $|J| = 2$, 
\[
\cN_{J}': = \cN' \cap V^{3}_{J}.
\]
Then the following hold.
\vskip 1mm
\begin{enumerate}
\item[(a)] Each $\cN_{J}'$ is a compact set. Moreover, 
\begin{equation}\label{eq:N'-J-decomp}
\begin{split}
\cN'_{J} =\ & \big( \cN \cap  V^{3}_{J} \cap \cA_{\tau\alpha} \big)\cup  \Psi_{J}((\mu_0\alpha, \mu_5\alpha] \times (\partial\Omega_{J}\cap \cN))\\
=\ &\big( \cN \cap  V^{3}_{J} \cap \cA_{r} \big)\cup  \Psi_{J}((\mu_0\alpha, \mu_5\alpha] \times (\partial\Omega_{J}\cap \cN)),
\end{split}
\end{equation}
for any $r \in [(\tau-2\delta)\alpha, (\tau+2\delta)\alpha]$.
\vskip 1mm
\item[(b)] For $s \in (0, 4\alpha']$, we have
\[
\cN' \cap \cA^{(2)}_{s} = \cup_{|J| =2} \big(B^{1, J}_{s} + \cN'_{J} \big),
\]
and this continues to hold with $\cA_{s}^{(2)}$ and each $B^{1, J}_{s}$ replaced by their closures.
\end{enumerate}
\end{lemm}
\begin{proof}
For part (a), the first equality in~\eqref{eq:N'-J-decomp} follows straight from~\eqref{eq:Phi-J-positive-distance} and the definition~\eqref{eq:N'-definition} of $\cN'$, while the second equality is a consequence of~\eqref{eq:Psi-rmk-supset-before-union}. Next, by~\eqref{eq:Psi-image-for-compact} and the admissible range of $r$ in~\eqref{eq:N'-J-decomp}, we get
\[
\cN_{J}' = \big(\cN \cap \overline{\cA_{\tau\alpha}} \cap V^{3}_{J}\big) \cup \Psi_{J}([\mu_0\alpha, \mu_5\alpha] \times (\partial\Omega_{J} \cap \cN)),
\]
which shows that $\cN_{J}'$ is compact. See the observations made after~\eqref{eq:boundary-Omega-N'-intersection-compact}.

For part (b), since each $\cN'_{J}$ is a subset of $V^{3}_{J}$, it is clear that 
\[
\cup_{|J| =2} \big(B^{1, J}_{s} + \cN'_{J} \big) \subset \cA_{s}^{(2)}.
\]
Next, given $|J| = 2$, by the first equality in~\eqref{eq:distance-attained} along with~\eqref{eq:tube-intersection}, as well as the fact that $V^{1, J}$ is orthogonal to the plane containing $V^{3}_{J}$, it is not hard to see that 
\begin{equation}\label{eq:fattening-of-core}
B^{1, J}_{s} + \big(\cN \cap  V^{3}_{J} \cap \cA_{(\tau - \delta)\alpha}\big) \subset \cN \cap \cA_{\tau\alpha}.
\end{equation}
Using this in~\eqref{eq:N'-J-decomp} (with $r = \tau - \delta$), and recalling the definition~\eqref{eq:N'-definition} of $\cN'$, we obtain 
\[
\cup_{|J| =2} \big(B^{1, J}_{s} + \cN'_{J} \big)  \subset \cN',
\]
and we have proved the inclusion ``$\supset$'' in (b). 

Conversely, suppose $z \in \cN' \cap \cA^{(2)}_{s}$. Then there exist $J_0 \subset \{1, 2, 3, 4\}$ with $|J_0| = 2$, along with $y \in V^{3}_{J_0}$ and $x \in B^{1, J_0}_{s}$, such that 
\begin{equation}\label{eq:N'-A-intersection-z-decomp}
z = x + y.
\end{equation}
If $z \in  C^{J}_{h}((\mu_0\alpha, \mu_5\alpha])$ for some $|J| = 2$, then~\eqref{eq:Phi-J-positive-distance} forces $J=J_0$, and~\eqref{eq:N'-A-intersection-z-decomp} must coincide with the decomposition of $z$ coming from the definition of $C^{J_0}_{h}((\mu_0\alpha, \mu_5\alpha])$, as in the proof of Lemma~\ref{lemm:Psi-action-on-N'}(c). In particular, we get
\[
y \in \Psi_{J_0}((\mu_0\alpha, \mu_5\alpha] \times (\cN \cap \partial\Omega_{J_0})) \subset \cN'_{J_0},
\]
so that $z \in B^{1, J_0}_{s} + \cN'_{J_0}$. On the other hand, suppose $z \in \cN \cap \cA_{\tau\alpha}$. If in addition $z \in B_{(\tau + \delta)\alpha}$, then we get 
\[
y \in B_{(\tau + 2\delta)\alpha} \cap V^{3}_{J_0} \subset \cN \cap \cA_{(\tau + 2\delta)\alpha} \cap V^{3}_{J_0}\subset \cN_{J_0}',
\]
so again $z \in B^{1, J_0}_{s} + \cN_{J_0}'$. The remaining case is $z \in \cN \cap \cA_{\tau\alpha}\setminus B_{(\tau + \delta)\alpha}$. Here we can find $i \in\{1, 2, 3, 4\}$ such that
\begin{equation}\label{eq:z-decomp-2}
z = t\frac{\ba_{i}}{|\ba_{i}|}+ z_1,
\end{equation}
for some $z_1 \in B^{2, \{i\}^{c}}_{\tau\alpha}$ and $t \in [\delta \alpha, 2]$. The proof of Lemma~\ref{lemm:sectors-product-structure}(b) then gives
\[
\dist(z, V^{3}_{J}) \geq \delta\alpha, 
\]
for any $J$ with $|J| = 2$ that contains $i$. It follows that in~\eqref{eq:N'-A-intersection-z-decomp} we have $i \not\in J_0$, so that $\ba_{i} \cdot x = 0$ by~\eqref{eq:decomposition-wrt-lower-simplex}. Combining this with $\ba_{i} \cdot z_{1} = 0$ and $y \cdot x = 0$, and also rearranging~\eqref{eq:z-decomp-2} and~\eqref{eq:N'-A-intersection-z-decomp} to get
\[
z_1 = x + (y - t\frac{\ba_{i}}{|\ba_{i}|}),
\]
we see that $y  - t\frac{\ba_{i}}{|\ba_{i}|}$ is orthogonal to both $x$ and $\ba_{i}$. Recalling that $|z_1| < \tau\alpha$, and that $y \in V^{3}_{J_0}$ to begin with, we conclude
\[
y = t\frac{\ba_{i}}{|\ba_{i}|} + (y - t\frac{\ba_{i}}{|\ba_{i}|}) \in \cN \cap \cA_{\tau\alpha} \cap V^{3}_{J_0} \subset \cN'_{J_0},
\]
and again we get $z \in B^{1, J_0}_{s} + \cN'_{J_0}$. The very last assertion of the lemma can be proved by the exact same argument as above, and we omit the details.
\end{proof}

\begin{prop}\label{prop:smoothing-in-distance-nbhd}
Given $j_0 \in \{1, 2, 3, 4\}$, $s \in (0, \alpha']$, and $r \in (3\rho_2, 4\alpha]$, we have the following.
\vskip 1mm
\begin{enumerate}
\item[(a)] $\widetilde{M}_{j_0} \cap G((\cN \cap \cA_{r} \cap \cA_{s}^{(2)}) \times S_0) = G((\cN \cap \cA_{r} \cap \cA_{s}^{(2)} \cap {\Omega}) \times S_0)$.
\vskip 1mm
\item[(b)] $\widetilde{M}_{j_0} \cap G_{(1)}((\cN' \cap \cA_{s}^{(2)}) \times S_0) \subset g_{j_0}([0, s) \times \partial\widetilde{M}_{j_0})$.
\end{enumerate}
\end{prop}
\begin{proof}
For part (a), we first use Proposition~\ref{prop:5-strata-smoothing} to see that
\begin{equation}\label{eq:smoothing-in-distance-nbhd-a-1}
\begin{split}
&\widetilde{M}_{j_0} \cap G((\mathring{\cN}_{2/3} \cap \cA_{r} \cap \cA_{s}^{(2)}) \times S_0)\\
=\ &  G((\mathring{\cN}_{2/3} \cap \cA_{r} \cap {\Omega}) \times S_0) \cap G((\mathring{\cN}_{2/3} \cap \cA_{r} \cap \cA_{s}^{(2)}) \times S_0)\\
=\ & G((\mathring{\cN}_{2/3} \cap \cA_{r} \cap \cA_{s}^{(2)} \cap {\Omega}) \times S_0),
\end{split}
\end{equation}
where the second equality uses the injectivity of $G$ on $\mathring{\cN}_{7/8} \times S_0$. Again using Proposition~\ref{prop:5-strata-smoothing}, along with the inclusion $\widetilde{M}_{j_0} \subset M_{j_0}$ and Proposition~\ref{prop:1234-with-123}(d), we get
\begin{align*}
&\widetilde{M}_{j_0} \cap G(((\cN \setminus \cN_{5/8}) \cap \cA_{r} \cap \cA_{s}^{(2)}) \times S_0)\\
=\ &   G(((\cN \setminus \cN_{5/8}) \cap \cA_{r} \cap {\Omega}) \times S_0) \cap G(((\cN\setminus \cN_{5/8}) \cap \cA_{r} \cap\ \cA_{s}^{(2)} \cap V_{j_0}^{3}) \times S_0).
\end{align*}
By~\eqref{eq:product-structure-for-image-Omega} and~\eqref{eq:product-structure-with-distance-nbhds}, applied respectively to the two sets on the last line, as well as the inclusion $\Omega'_{\lambda} \subset V^{2, \{i_{\lambda}\}^{c}}_{j_0}$, we obtain the first of the following equalities:
\begin{align*}
& \widetilde{M}_{j_0} \cap G(((\cN \setminus \cN_{5/8}) \cap \cA_{r} \cap \cA_{s}^{(2)}) \times S_0)\\
=\ & \cup_{\lambda \in \{0, 1, 2\}}f_{\{i_{\lambda}\}^{c}}^{-1}( (B_{r} \cap \Omega_{\lambda}' \cap B_{s}(\partial V^{2, \{i_\lambda\}^c}_{j_0}) ) \times M_{\{i_\lambda\}^{c}}^{**} )\\
=\ & G( ( (\cN\setminus \cN_{5/8}) \cap \cA_{r} \cap {\Omega} \cap \cA_{s}^{(2)}   ) \times S_0 ),
\end{align*}
where for the second equality we argued as in the proof of~\eqref{eq:product-structure-with-distance-nbhds}, using in addition Proposition~\ref{prop:smoothing-product-structure}(a). We get part (a) from the above along with~\eqref{eq:smoothing-in-distance-nbhd-a-1}.

For part (b), since $\widetilde{M}_{j_0} \subset M_{j_0}$, we first have by Proposition~\ref{prop:G-1-properties}(c) that 
\[
\widetilde{M}_{j_0} \cap G_{(1)}((\cN' \cap \cA_{s}^{(2)}) \times S_0) = \widetilde{M}_{j_0} \cap G_{(1)}((\cN' \cap \cA_{s}^{(2)} \cap V^{3}_{j_0}) \times S_0).
\]
By Lemma~\ref{lemm:N'-intersect-2-strata-neighborhood} and~\eqref{eq:V-j-disjoint}, we get
\[
\begin{split}
\cN' \cap \cA_{s}^{(2)} \cap V^{3}_{j_0} =\ & \cup_{\lambda \in \{0, 1, 2\}} (B^{1, \{j_0, i_{\lambda}\}}_{s} + \cN'_{\{j_0, i_{\lambda}\}}) \cap V^{3}_{j_0}.
\end{split}
\]
For each $\lambda \in \{0, 1, 2\}$, using the expression~\eqref{eq:N'-J-decomp} for $\cN'_{\{j_0, i_\lambda\}}$ with $r = (\tau - 2\delta)\alpha$ and following the reasoning for~\eqref{eq:fattening-of-core}, we find that 
\[
(B^{1, \{j_0, i_{\lambda}\}}_{s} + \cN'_{\{j_0, i_{\lambda}\}}) \cap V^{3}_{j_0} \subset  \big( \cN \cap \cA_{(\tau - \delta)\alpha} \cap \cA_{s}^{(2)} \cap V^{3}_{j_0} \big) \cup\big(C^{\{j_0, i_\lambda\}}_{s}((\mu_0\alpha, \mu_5\alpha]) \cap V^{3}_{j_0}\big).
\]
Noting, with the help of~\eqref{eq:inward-pointing} and the definition~\eqref{eq:w-lambda-definition} of $\bw_{\lambda}$, that
\[
V^{3}_{j_0} \cap C^{\{j_0, i_\lambda\}}_{s}((\mu_0\alpha, \mu_5\alpha]) \subset [0, s) \cdot \bw_{\lambda} + \Psi_{\{j_0, i_{\lambda}\}}((\mu_0\alpha, \mu_5\alpha] \times (\partial \Omega_{\{j_0, i_{\lambda}\}} \cap \cN)),
\]
and recalling the definition of $G_{(1)}$, we arrive at
\begin{equation}\label{eq:2-strata-nbhd-image-key}
\begin{split}
&\widetilde{M}_{j_0} \cap G_{(1)}((\cN' \cap \cA_{s}^{(2)}) \times S_0) \\
\subset \ & \widetilde{M}_{j_0} \cap G( (\cN \cap \cA_{(\tau - \delta)\alpha} \cap \cA_{s}^{(2)})\times S_0)\\
&\cup \big( \cup_{\lambda \in \{0, 1, 2\}}H_{\{j_0,i_\lambda\}}([0, s)\cdot \bw_{\lambda} \times (\mu_0\alpha, \mu_5\alpha] \times \partial\widetilde{M}_{\{j_0, i_\lambda\}}) \big).
\end{split}
\end{equation}
To continue, note that, by~\eqref{eq:distance-2-strata-attained} and~\eqref{eq:distance-neighborhood-contained}, we have
\[
\begin{split}
\cN \cap \cA_{(\tau - \delta)\alpha} \cap \cA_{s}^{(2)} \cap \Omega \subset\ & \cN \cap \cA_{(\tau - \delta)\alpha} \cap \Psi([0, s) \times \partial\Omega) \\
\subset\ & \cN \cap \Psi([0, s) \times (\partial\Omega \cap \cA_{\tau\alpha})) \quad\text{(triangle inequality)}\\
=\ & \Psi([0, s) \times (\partial\Omega \cap \cA_{\tau\alpha} \cap \cN)) \quad\text{(Lemma~\ref{lemm:N-preserved-by-flow})}.
\end{split}
\]
From this and part (a), we obtain
\[
\begin{split}
\widetilde{M}_{j_0} \cap G( (\cN \cap \cA_{(\tau - \delta)\alpha} \cap \cA_{s}^{(2)} )\times S_0) =\ & G( (\cN \cap \cA_{(\tau - \delta)\alpha} \cap \cA_{s}^{(2)} \cap \Omega)\times S_0)\\
\subset\  & G(\Psi([0, s) \times (\partial\Omega \cap \cA_{\tau\alpha} \cap \cN)) \times S_0) \\
=\ & g_{j_0}([0, s) \times G( (\partial\Omega\cap \cA_{\tau\alpha} \cap \cN) \times S_0)) \\
&\text{(by~\eqref{eq:cD-flexible} and~\eqref{eq:g-j-expression-in-G})}\\
\subset\ & g_{j_0}([0, s) \times \partial\widetilde{M}_{j_0}) \quad\text{(Proposition~\ref{prop:5-strata-smoothing})}.
\end{split}
\]
Turning to the union over $\lambda$ at the end of~\eqref{eq:2-strata-nbhd-image-key}, we note by~\eqref{eq:g-j-expression-in-H} that
\[
\begin{split}
&H_{\{j_0,i_\lambda\}}([0, s)\cdot\bw_{\lambda} \times (\mu_0\alpha, \mu_5\alpha] \times \partial\widetilde{M}_{\{j_0, i_\lambda\}})\\
& = g_{j_0}([0, s) \times h_{\{j_0,i_\lambda\}}((\mu_0\alpha, \mu_5\alpha] \times \partial\widetilde{M}_{\{j_0, i_\lambda\}})) \subset g_{j_0}([0, s) \times \partial\widetilde{M}_{j_0}),
\end{split}
\]
the second line being a consequence of Proposition~\ref{prop:4-strata-smoothing-decomp-far}. Substituting the previous two observations back into~\eqref{eq:2-strata-nbhd-image-key} completes the proof of (b).
\end{proof}

In what follows, given $j_0 \in \{1, 2, 3, 4\}$, when we want to emphasize the $j_0$-dependence of the objects mentioned below~\eqref{eq:alpha'-threshold}, we denote them by $\Omega_{j_0}$, $\xi_{j_0}$, $\Psi_{j_0}$, and so on. By Remark~\ref{rmk:Psi-distance-bounds} and~\eqref{eq:distance-2-strata-attained}, in a manner similar to how we chose $\mu_0$, $\mu_1$, $\tau$, and $\delta$ in~\eqref{eq:mu1-mu0-choices} and~\eqref{eq:transition-parameter}, we obtain universal positive constants $\nu_{1}, \nu_{2}, \tau'$, and $\eta$ such that
\[
0 < \nu_{1} < \tau' - 2\eta < \tau' + 2\eta < \nu_{2} < 1,
\]
and that
\begin{align}
&\Omega_{j_0}  \cap \cA_{\nu_{1}\alpha'}^{(2)} \subset  \Psi_{j_0}([0, \nu_1\alpha') \times \partial\Omega_{j_0} ) \Subset \Omega_{j_0}  \cap \cA_{(\tau' - 2\eta)\alpha'}^{(2)},\label{eq:2-strata-neighborhood-squeezed-1}\\
&\Omega_{j_0} \cap \cA_{(\tau' + 2\eta)\alpha'}^{(2)} \subset \Psi_{j_0}([0, \nu_2\alpha') \times \partial\Omega_{j_0}) \Subset \Omega_{j_0}  \cap \cA_{\alpha'}^{(2)}, \label{eq:2-strata-neighborhood-squeezed-2}\\
&\Omega_{j_0}  \cap \cA_{\alpha'}^{(2)} \subset 
\Psi_{j_0}([0, \alpha') \times \partial\Omega_{j_0} ) \subset \Omega_{j_0}  \cap \cA_{2\alpha'}^{(2)}, \label{eq:2-strata-neighborhood-squeezed-3}\\
&\Omega_{j_0} \setminus \cA_{(\tau' - 2\eta)\alpha'}^{(2)} = V^{3}_{j_0}\setminus \cA_{(\tau' - 2\eta)\alpha'}^{(2)}.\label{eq:smoothing-agree-2-strata-nbhd}
\end{align}
We will often use~\eqref{eq:2-strata-neighborhood-squeezed-1} through~\eqref{eq:2-strata-neighborhood-squeezed-3} in conjunction with Lemma~\ref{lemm:Psi-action-on-N'}(c), without explicitly citing the latter. Fixing in addition $\nu_3 < \nu_{4}$ from $(\nu_2, 1)$, by the properties of $g_{j_0}$ established in Lemma~\ref{lemm:g-j-target} through Proposition~\ref{prop:g-j-injective}, we see that the restriction 
\[
g_{j_0}|_{[0, \nu_3\alpha') \times \partial\widetilde{M}_{j_0}}
\]
is an injective immersion into $\widetilde{M}_{j_0}$, and sends $(0, q)$ to $q$ for all $q \in \partial\widetilde{M}_{j_0}$. Thus, by Lemma~\ref{lemm:boundary-IVT}, the above parametrizes a collar neighborhood of $\partial\widetilde{M}_{j_0}$ in $\widetilde{M}_{j_0}$. Since the latter is a $1$-handlebody, we obtain, as before, a Lipschitz map 
\begin{equation}\label{eq:collar-5-strata}
h_{j_0}:[0, \nu_{4}\alpha'] \times \partial\widetilde{M}_{j_0} \to \widetilde{M}_{j_0} 
\end{equation}
with the following properties:
\vskip 1mm
\begin{enumerate}
\item $s \mapsto h_{j_0}(s, \cdot)$ is continuous from $[0, \nu_{4}\alpha')$ to $C^{1}(\partial\widetilde{M}_{j_0}; M)$.
\vskip 1mm
\item $\Gamma_{j_0}: = h_{j_0}(\{\nu_4\alpha'\} \times \partial\widetilde{M}_{j_0})$ has finite $\cH^{1}$-measure. Moreover, writing $h_{j_0, s}$ for $h_{j_0}(s, \cdot)$, and defining $\|\Lambda^2 dh_{j_0, s} \|$ as in~\eqref{eq:h-J-shrink-area}, we have 
\begin{equation}\label{eq:h-j0-shrink-area}
\|\Lambda^2 dh_{j_0, s}\| \to 0 \quad\text{as }s \to (\nu_{4}\alpha')^{-}.
\end{equation}
\vskip 1mm
\item $h_{j_0}([0, \nu_4\alpha']\times \partial\widetilde{M}_{j_0}) = \widetilde{M}_{j_0}$. Also, 
\begin{equation}\label{eq:h-j0-disjoint-images}
h_{j_0}(\{s\} \times \partial \widetilde{M}_{j_0}) \cap h_{j_0}(\{s'\} \times \partial \widetilde{M}_{j_0}) = \emptyset, \quad\text{whenever }s \neq s'.
\end{equation}
\vskip 1mm
\item $h_{j_0}(s, q) = g_{j_0}(s, q)$, for all $(s, q) \in [0, \nu_{3}\alpha'] \times \partial\widetilde{M}_{j_0}$. In particular $h_{j_0}|_{(0, \nu_3\alpha') \times \partial\widetilde{M}_{j_0}}$ is a diffeomorphism onto an open set in $\widetilde{M}_{j_0}\setminus \partial\widetilde{M}_{j_0}$. Note also that
\begin{equation}\label{eq:5-strata-smoothing-collar-deleted}
\begin{split}
M_{j_0}^{*} :=\ & \widetilde{M}_{j_0} \setminus g_{j_0}([0, \frac{\nu_1\alpha'}{2}] \times \partial\widetilde{M}_{j_0}) = h_{j_0}((\frac{\nu_1\alpha'}{2}, \nu_{4}\alpha'] \times \partial\widetilde{M}_{j_0}).
\end{split}
\end{equation}
\end{enumerate}
\begin{lemm}\label{lemm:M-j-star-disjoint}
Given $j_0 \in \{1, 2, 3, 4\}$, listing the elements of $\{j_0\}^{c}$ as $\{i_0, i_1, i_2\}$, we have 
\[
M_{j_0}^{*} \cap M_{i_\lambda} = \emptyset, \quad\text{for all }\lambda \in \{0, 1, 2\}.
\]
\end{lemm}
\begin{proof}
By Proposition~\ref{prop:smoothing-in-distance-nbhd}(b), we have
\begin{equation}\label{eq:M-j-star-disjoint-from-image}
M_{j_0}^{*} \subset \widetilde{M}_{j_0}\setminus G_{(1)}((\cN' \cap \cA^{(2)}_{\frac{\nu_1\alpha'}{2}}) \times S_0),
\end{equation}
and thus
\begin{equation}\label{eq:M-j-star-disjoint-start}
\begin{split}
M_{j_0}^{*} \cap M_{i_\lambda} \subset \ & M_{\{j_0, i_\lambda\}} \setminus G_{(1)}((\cN' \cap \cA^{(2)}_{\frac{\nu_1\alpha'}{2}} ) \times S_0)\\
=\ & M_{\{j_0, i_\lambda\}} \setminus G_{(1)}((\cN' \cap V^{3}_{\{j_0, i_\lambda\}}) \times S_0)
\end{split}
\end{equation}
where the last step uses Proposition~\ref{prop:G-1-properties}(c) and the inclusion $V^{3}_{\{j_0,i_\lambda\}} \subset \cA^{(2)}_{\frac{\nu_1\alpha'}{2}}$. From the expression~\eqref{eq:N'-J-decomp} for $\cN'\cap V^{3}_{\{j_0, i_\lambda\}}$ and the definition of $G_{(1)}$ in~\eqref{eq:G-1-definition-case-1} and~\eqref{eq:G-1-definition-case-2}, followed by Proposition~\ref{prop:1234-with-123}(d) and Proposition~\ref{prop:H-J-properties}(b), we get
\[
\begin{split}
G_{(1)}((\cN'  \cap V^{3}_{\{j_0, i_\lambda\}}) \times S_0) =\  & M_{\{j_0, i_\lambda\}} \cap G((\cN \cap \cA_{\tau\alpha}) \times S_0)\\
& \cup h_{\{j_0, i_\lambda\}}((\mu_0\alpha, \mu_5\alpha] \times \partial\widetilde{M}_{\{j_0, i_\lambda\}}).
\end{split}
\]
Combining this with~\eqref{eq:M-j-star-disjoint-start} leads to
\begin{equation}\label{eq:5-strata-smoothing-disjoint-key}
\begin{split}
M_{j_0}^{*} \cap M_{i_\lambda} \subset\ & \big(M_{\{j_0, i_\lambda\}} \setminus G((\cN \cap \cA_{\tau\alpha}) \times S_0)\big)\setminus h_{\{j_0, i_\lambda\}}((\mu_0\alpha, \mu_5\alpha] \times \partial\widetilde{M}_{\{j_0, i_\lambda\}}).
\end{split}
\end{equation}
Recalling that $\tau > \tau_0 + \mu_0\sin^2\theta > \max\{\tau_0, \mu_0\}$, and using parts (b) and (a) of Proposition~\ref{prop:distance-neighborhood-in-smoothing}, we see that
\[
\begin{split}
M_{\{j_0, i_\lambda\}} \setminus G((\cN \cap \cA_{\tau\alpha}) \times S_0) =\ & \widetilde{M}_{\{j_0, i_\lambda\}} \setminus G((\cN \cap \cA_{\tau\alpha}) \times S_0)\\
\subset\ & \widetilde{M}_{\{j_0, i_\lambda\}} \setminus h_{\{j_0, i_\lambda\}}([0, \mu_0\alpha] \times \partial\widetilde{M}_{\{j_0,i_\lambda\}}).
\end{split}
\]
Putting this back into~\eqref{eq:5-strata-smoothing-disjoint-key}, we get
\[
M_{j_0}^{*}\cap M_{i_\lambda} \subset \widetilde{M}_{\{j_0, i_\lambda\}} \setminus h_{\{j_0, i_\lambda\}}([0, \mu_5\alpha] \times \partial\widetilde{M}_{\{j_0,i_\lambda\}}) = \emptyset.
\]
The proof is complete.\\
\end{proof}

Now let
\begin{equation}\label{eq:N''-definition}
\begin{split}
\cN'' =\ & (\cN' \cap \cA_{\tau'\alpha'}^{(2)})  \cup \big(\cup_{j_0 =1}^4 \Psi_{j_0}((\nu_1\alpha', \nu_{4}\alpha'] \times (\partial\Omega_{j_0} \cap \cN') )\big),\\
\mathring{\cN}'' =\ & (\mathring{\cN}' \cap \cA_{\tau'\alpha'}^{(2)})  \cup \big(\cup_{j_0 =1}^4 \Psi_{j_0}((\nu_1\alpha', \nu_{4}\alpha') \times (\partial\Omega_{j_0} \cap \mathring{\cN}') )\big),\\
\mathring{\cD}'' = \ & (\mathring{\cN}'_{7/8, \mu_{4}\alpha}\cap \cA_{\tau'\alpha'}^{(2)}) \cup \big(\cup_{j_0 =1}^4 \Psi_{j_0}((\nu_1\alpha', \nu_{3}\alpha') \times (\partial\Omega_{j_0} \cap \mathring{\cN}'_{7/8, \mu_{4}\alpha}) ) \big).
\end{split}
\end{equation}
Recalling from the lines below~\eqref{eq:N'-definition} that the sets $\mathring{\cN}'_{R, \ell}$ are open in $V^{n-1}$, we see with the help of Lemma~\ref{lemm:Psi-action-on-N'}(c) that both $\mathring{\cN}''$ and $\mathring{\cD}''$ are open sets as well. Again by Lemma~\ref{lemm:Psi-action-on-N'}(c), along with~\eqref{eq:2-strata-neighborhood-squeezed-3}, we find that 
\begin{equation}\label{eq:N''-contained-in-N'}
\cN'' \subset \cN' \cap \cA_{2\alpha'}^{(2)}, \quad \mathring{\cN}'' \subset \mathring{\cN}' \cap \cA_{2\alpha'}^{(2)},
\end{equation}
and that
\begin{equation}\label{eq:D''-contained}
\mathring{\cD}'' \subset \mathring{\cN}'_{7/8, \mu_{4}\alpha} \cap \cA_{2\alpha'}^{(2)}.
\end{equation}
Also, since $\partial\Omega_{j_0} \subset V^{3}_{j_0}$, and since $\xi_{j_0}$ takes values in $\mathring{V}^{3}_{j_0}$, we have
\begin{equation}\label{eq:N''-limbs}
\Psi_{j_0}([\nu_1\alpha', \alpha') \times \partial\Omega_{j_0} ) \subset \mathring{V}_{j_0}^{3} \quad\text{for each }j_0 \in \{1, 2, 3, 4\},
\end{equation}
which implies, by the compactness of $\partial\Omega_{j_0} \cap \cN'$ coming from Lemma~\ref{lemm:Psi-action-on-N'}(a), that
\begin{equation}\label{eq:N''-limbs-separated}
\min_{j_0 \in \{1, 2, 3, 4\}} \dist\big( \Psi_{j_0}([\nu_1\alpha', \frac{(1+\nu_4)\alpha'}{2}] \times (\partial\Omega_{j_0} \cap \cN')),\ V^{3} \setminus \mathring{V}^{3}_{j_0}\big) > 0.
\end{equation}
With the above preparation, we then define $G_{(2)}:\cN''\times S_0 \to M$ as follows:
\vskip 1mm
\begin{itemize}
\item If $y \in \cN' \cap \cA_{\tau'\alpha'}^{(2)}$, we let 
\begin{equation}\label{eq:G-2-definition-case-1}
G_{(2)}(y, p) = G_{(1)}(y, p).
\end{equation}
\vskip 1mm
\item If $y\in \Psi_{j_0}((\nu_1\alpha', \nu_{4}\alpha'] \times (\partial\Omega_{j_0} \cap \cN') )$ for some $j_0 \in \{1, 2, 3, 4\}$, necessarily unique by~\eqref{eq:N''-limbs}, we write $(s, z) = \Psi_{j_0}^{-1}(y)$, recall from Proposition~\ref{prop:M-j-as-G-1-image} that $G_{(1)}(z, p) \in \partial\widetilde{M}_{j_0}$, and let 
\begin{equation}\label{eq:G-2-definition-case-2}
G_{(2)}(y, p) = h_{j_0}(s, G_{(1)}(z, p)).
\end{equation}
\end{itemize}

\begin{lemm}\label{lemm:cN''-compact}
The set $\cN''$ is compact, and $\mathring{\cN}''$ coincides with its interior. Moreover, $\cN'' \subset \overline{\mathring{\cN}''}$.
\end{lemm}
\begin{proof}
By~\eqref{eq:2-strata-neighborhood-squeezed-1},~\eqref{eq:2-strata-neighborhood-squeezed-2}, and~\eqref{eq:smoothing-agree-2-strata-nbhd}, we have for each $j_{0} \in \{1, 2, 3, 4\}$ that
\[
\begin{split}
\cN' \cap (\cA_{(\tau' + \eta)\alpha'}^{(2)} \setminus \cA_{(\tau' - \eta)\alpha'}^{(2)})) \cap V^{3}_{j_0} = \ & \cN' \cap (\cA_{(\tau' + \eta)\alpha'}^{(2)} \setminus \cA_{(\tau' - \eta)\alpha'}^{(2)})) \cap \Omega_{j_0}\\
\subset\ & \Psi_{j_0}( (\nu_{1}\alpha', \nu_{2}\alpha') \times (\partial\Omega_{j_0} \cap \cN') ).
\end{split}
\]
Since $V^{3} = \cup_{j_0 = 1}^{4}V^{3}_{j_0}$ by~\eqref{eq:sectors-cover}, we infer from the above that
\begin{equation}\label{eq:N''-flexible}
\begin{split}
\cN'' = (\cN' \cap \cA_{s}^{(2)})  \cup \big(\cup_{j_0 =1}^4 \Psi_{j_0}((\nu_1\alpha', \nu_{4}\alpha'] \times (\partial\Omega_{j_0} \cap \cN') )\big),
\end{split}
\end{equation}
for all $s \in [(\tau' - \eta)\alpha', (\tau' + \eta)\alpha']$. On the other hand, by~\eqref{eq:2-strata-neighborhood-squeezed-1}, we have
\[
\Psi_{j_0}(\{\nu_{1}\alpha'\} \times (\partial\Omega_{j_0} \cap \cN')) \subset  \cA_{(\tau' - \eta)\alpha'}^{(2)} \cap \cN'. 
\]
Combining this with~\eqref{eq:N''-flexible} gives
\[
\cN'' = (\cN' \cap \overline{\cA_{\tau'\alpha'}^{(2)}})  \cup \big(\cup_{j_0 =1}^4 \Psi_{j_0}([\nu_1\alpha', \nu_{4}\alpha'] \times (\partial\Omega_{j_0} \cap \cN') )\big).
\]
Since $\cN' \cap \overline{\cA_{\tau'\alpha'}^{(2)}}$ and $\partial\Omega_{j_0} \cap \cN'$ are compact, thanks respectively to Lemmas~\ref{lemm:N'-intersect-2-strata-neighborhood} and~\ref{lemm:Psi-action-on-N'}, we conclude that $\cN''$ is a compact set.

We have already seen in the comments after~\eqref{eq:N''-definition} that $\mathring{\cN}''$ is open. To prove that it coincides with the interior of $\cN''$, we first note that 
\begin{equation}\label{eq:N''-minus-interior}
\cN'' \setminus \mathring{\cN}''\subset (\cN' \setminus \mathring{\cN}') \cup \big( \cup_{j_0 = 1}^{4}\Psi_{j_0}(\{\nu_{4}\alpha'\} \times (\partial\Omega_{j_0} \cap \cN')) \big) \cup \big( \cup_{j_0 = 1}^{4}\cE_{j_0} \big),
\end{equation}
where 
\[
\cE_{j_0}: = \Psi_{j_0}((\nu_1\alpha', \nu_{4}\alpha') \times (\partial\Omega_{j_0} \cap (\cN'\setminus \mathring{\cN}')) ).
\]
By Lemma~\ref{lemm:Psi-action-on-N'}(c) and Remark~\ref{rmk:N'-ring-version}, we have $\cE_{j_0} \subset \cN' \setminus \mathring{\cN}'$, which together with Lemma~\ref{lemm:N'-properties} and~\eqref{eq:N''-contained-in-N'} gives
\[
(\cN' \setminus \mathring{\cN}') \cup \big( \cup_{j_0 = 1}^{4}\cE_{j_0} \big) \subset \overline{V^{3}\setminus \cN'} \subset \overline{V^{3} \setminus \cN''}.
\]
Next, by~\eqref{eq:2-strata-neighborhood-squeezed-2} and~\eqref{eq:N''-limbs}, we have $ \Psi_{j_0}((\nu_4\alpha', \alpha') \times (\partial\Omega_{j_0} \cap \cN') ) \cap \cN'' = \emptyset$, and thus
\[
\Psi_{j_0}(\{\nu_{4}\alpha'\} \times (\partial\Omega_{j_0} \cap \cN')) \subset \overline{V^{3} \setminus \cN''}.
\]
Putting the two previous displayed inclusions back into~\eqref{eq:N''-minus-interior}, we get 
\[
\cN'' \setminus \mathring{\cN}'' \subset \overline{V^{3} \setminus \cN''}, 
\]
which together with the openness of $\mathring{\cN}''$ shows that it is the interior of $\cN''$. 

For the last conclusion, recall from Lemma~\ref{lemm:N'-properties} that $\cN' \subset\overline{\mathring{\cN}'}$, and thus
\begin{equation}\label{eq:N''-properties-core-inclusion}
\cN' \cap \cA^{(2)}_{\tau'\alpha'} \subset \overline{\mathring{\cN}' \cap \cA^{(2)}_{\tau'\alpha'}}.
\end{equation}
Next we claim that 
\begin{equation}\label{eq:N''-properties-claim}
\partial\Omega_{j_0} \cap \cN' \subset \overline{\partial\Omega_{j_0} \cap \mathring{\cN}'}. 
\end{equation}
Indeed, given $\lambda \in \{0, 1, 2\}$, as pointed out in the proof of Lemma~\ref{lemm:N'-properties}, we have by Remark~\ref{rmk:portions} that $\partial\Omega_{\{j_0, i_{\lambda}\}} \cap\cN \subset \overline{\partial\Omega_{\{j_0, i_{\lambda}\}} \cap\mathring{\cN}}$, implying that 
\[
\Psi_{\{j_0, i_{\lambda}\}}((\mu_0\alpha, \mu_5\alpha] \times (\partial\Omega_{\{j_0, i_{\lambda}\}} \times \cN)) \subset \overline{\Psi_{\{j_0, i_{\lambda}\}}((\mu_0\alpha, \mu_5\alpha) \times (\partial\Omega_{\{j_0, i_{\lambda}\}} \times \mathring{\cN}))}.
\]
On the other hand, the argument leading to~\eqref{eq:product-structure-for-image-Omega} yields similar expressions for $(\cN \setminus \mathring{\cN}) \cap \cA_{\tau\alpha} \cap \partial\Omega_{j_0}$ and $(\mathring{\cN}\setminus \cN_{5/8}) \cap \cA_{\tau\alpha} \cap \partial\Omega_{j_0}$, from which we deduce that 
\[
\partial\Omega_{j_0} \cap \cA_{\tau\alpha} \cap \cN \subset \overline{\partial\Omega_{j_0} \cap \cA_{\tau\alpha} \cap \mathring{\cN}}.
\]
Combining by the two previous displayed inclusion, as well as the expressions~\eqref{eq:boundary-Omega-N'-intersection} and~\eqref{eq:boundary-Omega-N'-ring-intersection}, we get~\eqref{eq:N''-properties-claim} as claimed. We conclude the proof upon combining~\eqref{eq:N''-properties-claim} and~\eqref{eq:N''-properties-core-inclusion} with the definitions in~\eqref{eq:N''-definition}.
\end{proof}
\begin{prop}\label{prop:G-2-properties}
$G_{(2)}:\cN'' \times S_0 \to M$ is a well-defined Lipschitz map, and the assignment $y \mapsto G_{(2)}(y,\cdot)$ is continuous from $\mathring{\cN}''$ into $C^{1}(S_0; M)$. Moreover, the following hold:
\vskip 1mm
\begin{enumerate}
\item[(a)] $G_{(2)}$ restricts to a diffeomorphism from $\mathring{\cD}'' \times S_0$ onto an open set in $M$.
\vskip 1mm
\item[(b)] We have 
\[
G_{(2)}( B_{\frac{\nu_1\alpha'}{2}}^{3} \times S_0) \cap G_{(2)}( ( \cN'' \setminus B_{\frac{\nu_1\alpha'}{2}}^{3}) \times S_0) = \emptyset,
\]
and that $B_{\frac{\nu_1\alpha'}{2}}^{3} \subset \mathring{\cD}''$.
\vskip 1mm
\item[(c)] Given $y_0 \in \cN'' \setminus \mathring{\cN}''$, we have $\cH^{1}(G_{(2)}(\{y_0\} \times S_0)) < \infty$. Moreover, 
\begin{equation}\label{eq:G-2-area-boundary}
\Area(G_{(2)}(y, \cdot)) \to 0\quad  \text{as }y \to y_0 \text{ from within } \mathring{\cN}''.
\end{equation}
\vskip 1mm
\end{enumerate}
\end{prop}
\begin{proof}
To see that $G_{(2)}$ is well-defined, take $j_0 \in \{1, 2, 3, 4\}$ along with
\[
y \in (\cN' \cap \cA_{\tau'\alpha'}^{(2)}) \cap \Psi_{j_0}((\nu_1\alpha', \nu_{4}\alpha'] \times (\partial\Omega_{j_0} \cap \cN')).
\]
Since $\Psi_{j_0}([0, \alpha') \times \partial\Omega_{j_0}) \subset \Omega_{j_0}$, we have by~\eqref{eq:2-strata-neighborhood-squeezed-2} that
\begin{equation}\label{eq:beyond-nu2-disjoint}
\begin{split}
&\Psi_{j_0}([\nu_2\alpha', \alpha') \times \partial\Omega_{j_0} ) \cap  \cA_{(\tau' + 2\eta)\alpha'}^{(2)}\\
\subset \ &  \Psi_{j_0}([\nu_2\alpha',\alpha') \times \partial\Omega_{j_0}) \cap \Psi_{j_0}([0, \nu_2\alpha') \times \partial\Omega_{j_0} )  =\emptyset.
\end{split}
\end{equation}
Thus, upon letting $(s, z) = \Psi_{j_0}^{-1}(y)$, we must have $s \in (\nu_1\alpha', \nu_2\alpha')$, so that, by property (4) above Lemma~\ref{lemm:M-j-star-disjoint},
\begin{equation}\label{eq:G-2-overlap}
h_{j_0}(s, G_{(1)}(z, p)) = g_{j_0}(s, G_{(1)}(z,p)).
\end{equation}
In view of Lemma~\ref{lemm:Psi-action-on-N'}(a) there are two possibilities regarding $z \in \partial\Omega_{j_0} \cap \cN'$: 

If $z \in \partial\Omega_{j_0} \cap \cA_{(\tau-\delta)\alpha} \cap \cN$, then by~\eqref{eq:G-1-definition-special-case} along with~\eqref{eq:g-j-expression-in-G}, we have
\[
g_{j_0}(s, G_{(1)}(z, p)) = g_{j_0}(s, G(z, p)) = G(y, p),
\]
while from Lemma~\ref{lemm:Psi-action-on-N'}(b) we get $y  = \Psi_{j_0}(s, z) \in \Omega_{j_0} \cap \cA_{\tau\alpha}\cap \cN$, in which case~\eqref{eq:G-1-definition-special-case} gives
\[
G(y, p) = G_{(1)}(y, p).
\]
Combining the previous two observations with~\eqref{eq:G-2-overlap} gives $h_{j_0}(s, G_{(1)}(z, p)) = G_{(1)}(y, p)$ in this case. 

For the other case, where $z = \Psi_{\{j_0, i_{\lambda}\}}(t, z')$ for some $\lambda \in \{0, 1, 2\}$, $t \in (\mu_0\alpha, \mu_5\alpha]$, and $z' \in \partial\Omega_{\{j_0, i_{\lambda}\}} \cap \cN$, we have by~\eqref{eq:G-1-definition-special-case} and~\eqref{eq:g-j-expression-in-H} that
\begin{equation}\label{eq:G-2-well-defined-case2-1}
g_{j_0}(s, G_{(1)}(z, p)) = g_{j_0}(s, h_{\{j_0, i_{\lambda}\}} (t, G(z', p)) ) = H_{\{j_0, i_{\lambda}\}}(s\bw_{\lambda}, t, G(z', p)).
\end{equation}
On the other hand, Lemma~\ref{lemm:Psi-action-on-N'}(b) and the inequality $\alpha' < \frac{h}{4}$ gives 
\[
(\,y =\,)\,\Psi_{j_0}(s, z) = s\bw_{\lambda} + z \in \Omega_{j_0} \cap C^{\{j_0, i_{\lambda}\}}_{h}((\mu_0\alpha, \mu_5\alpha]),
\]
which together with $z = \Psi_{\{j_0, i_{\lambda}\}}(t, z')$ and~\eqref{eq:G-1-definition-case-2} implies
\[
H_{\{j_0, i_{\lambda}\}}(s\bw_{\lambda}, t, G(z', p)) = G_{(1)}(y, p).
\]
From this along with~\eqref{eq:G-2-well-defined-case2-1} and~\eqref{eq:G-2-overlap}, we get $h_{j_0}(s, G_{(1)}(z, p)) = G_{(1)}(y, p)$ in this case also. In view of~\eqref{eq:N''-limbs}, we have shown that $G_{(2)}$ is well-defined. Moreover, the argument in fact gives
\begin{equation}\label{eq:G-2-definition-1-domain}
\begin{split}
G_{(2)}(y, p) = G_{(1)}(y, p)\quad \text{provided}\quad
y \in\ & (\cN' \cap \cA^{(2)}_{\tau'\alpha'})\\
&\cup \big( \cup_{j_0= 1}^4 \Psi_{j_0}((\nu_1\alpha', \nu_2\alpha') \times (\partial\Omega_{j_0} \cap \cN') ) \big).
\end{split}
\end{equation} 
We can then argue as in the lines below~\eqref{eq:G1-G-coincide} in the proof of Proposition~\ref{prop:G-1-properties} to see that $G_{(2)}$ is Lipschitz on $\cN'' \times S_0$. We omit the details of this process, and only mention that we make use of the compactness of each $\partial\Omega_{j_0} \cap \cN'$ (Lemma~\ref{lemm:Psi-action-on-N'}) and Proposition~\ref{prop:M-j-as-G-1-image}, and also observe that, thanks to~\eqref{eq:N''-limbs-separated} and~\eqref{eq:beyond-nu2-disjoint}, given two points in $\cN'' \times S_0$ which are sufficiently close to each other, if at least one them lies outside of the region where~\eqref{eq:G-2-definition-1-domain} applies, then there is some $j_0$ such that $\Psi_{j_0}((\nu_1\alpha', \nu_4\alpha'] \times (\partial\Omega_{j_0} \cap \cN')) \times S_0$ contains both points.

Next, to prove that 
\[
y \to G_{(2)}(y, \cdot)
\]
is continuous on $\mathring{\cN}''$ as a map into $C^{1}(S_0; M)$, we first notice that $\mathring{\cN}' \cap \cA_{\tau'\alpha'}^{(2)}$ is open in $V^{3}$, and that, with the help of Lemma~\ref{lemm:Psi-action-on-N'} (and Remark~\ref{rmk:N'-ring-version}) so is each $\Psi_{j_0}((\nu_1\alpha', \nu_4\alpha') \times (\partial\Omega_{j_0} \cap \mathring{\cN}'))$. Now, by the definition~\eqref{eq:G-2-definition-case-1} and Proposition~\ref{prop:G-1-properties}, we see that the above assignment is continuous into $C^{1}(S_0; M)$ when $y \in \mathring{\cN}' \cap \cA_{\tau'\alpha'}^{(2)}$. Using the definition~\eqref{eq:G-2-definition-case-2} along with the continuity properties of $s \mapsto h_{j_0}(s, \cdot)$ and $z\mapsto G_{(1)}(z, \cdot)$, given respectively by item (1) below~\eqref{eq:collar-5-strata} and Proposition~\ref{prop:G-1-properties}, and recalling also~\eqref{eq:M-j-as-G-1-image}, we see that $y \mapsto G_{(2)}(y, \cdot)$ is continuous as a map into $C^{1}(S_0; M)$ on each $\Psi_{j_0}((\nu_1\alpha', \nu_4\alpha') \times (\partial\Omega_{j_0} \cap \mathring{\cN}'))$ as well. Thus we have shown that $G_{(2)}(y, \cdot)$ varies continuously in $C^{1}(S_0; M)$ as $y$ varies in $\mathring{\cN}''$, as claimed.

Towards proving part (a), we define
\[
\begin{split}
\cC_{1} =\ & (\mathring{\cN}'_{7/8,\mu_{4}\alpha} \cap \cA^{(2)}_{\tau'\alpha'}) \cup \big( \cup_{j_0 =1}^4 \Psi_{j_0}((\nu_1\alpha', \nu_2\alpha') \times (\partial\Omega_{j_0} \cap \mathring{\cN}'_{7/8, \mu_{4}\alpha}) ) \big)\\
\cC_{2} =\ &  \cup_{j_0=1}^4 \Psi_{j_0}((\nu_1\alpha', \nu_{3}\alpha') \times (\partial\Omega_{j_0} \cap \mathring{\cN}'_{7/8, \mu_{4}\alpha})).
\end{split}
\]
Note that both are open sets in $V^{3}$, and that $\cC_{1} \cup \cC_{2} = \mathring{\cD}''$. By~\eqref{eq:beyond-nu2-disjoint} and~\eqref{eq:N''-limbs}, we also have
\begin{equation}\label{eq:cC1-cC2-intersection}
\cC_1 \cap \cC_{2} = \cup_{j_0 =1}^4 \Psi_{j_0}((\nu_1\alpha', \nu_2\alpha') \times (\partial\Omega_{j_0} \cap \mathring{\cN}'_{7/8, \mu_{4}\alpha}) ).
\end{equation}
Since $G_{(2)}$ agrees with $G_{(1)}$ on $\cC_{1} \times S_0$ thanks to~\eqref{eq:G-2-definition-1-domain}, and since $\cC_{1} \subset \mathring{\cN}'_{7/8, \mu_{4}\alpha}$ by Lemma~\ref{lemm:Psi-action-on-N'}(c) and Remark~\ref{rmk:N'-ring-version}, we deduce from Proposition~\ref{prop:G-1-properties}(a) that 
\[
G_{(2)}|_{\cC_1 \times S_0}\quad \text{is a diffeomorphism.}
\]
On the other hand, by Proposition~\ref{prop:G-1-properties}(a) and Proposition~\ref{prop:M-j-as-G-1-image}, we see that $G_{(1)}$ maps $(\partial\Omega_{j_0} \cap \mathring{\cN}'_{7/8, \mu_{4}\alpha}) \times S_0$ diffeomorphically onto an open subset of $\partial\widetilde{M}_{j_0}$. Since $h_{j_0}|_{(\nu_1\alpha', \nu_{3}\alpha') \times \partial\widetilde{M}_{j_0}}$ is a diffeomorphism by property (4) listed below~\eqref{eq:collar-5-strata}, it follows that $G_{(2)}$ restricts to a diffeomorphism on each $\Psi_{j_0}((\nu_1\alpha', \nu_{3}\alpha') \times (\partial\Omega_{j_0} \cap \mathring{\cN}'_{7/8, \mu_{4}\alpha})) \times S_0$. Noting that these latter sets have mutually disjoint images under $G_{(2)}$ thanks to~\eqref{eq:5-strata-smoothing-collar-deleted}, Lemma~\ref{lemm:M-j-star-disjoint}, and the inclusions $M_{j_0}^* \subset \widetilde{M}_{j_0} \subset M_{j_0}$, we conclude that $G_{(2)}$ restricts to a diffeomorphism on $\cC_{2} \times S_0$ as well. 

To continue, note that by~\eqref{eq:cC1-cC2-intersection} and the definition of $G_{(2)}$, there holds
\[
\begin{split}
G_{(2)}( (\cC_1\setminus (\cC_1\cap\cC_2)) \times S_0) \subset\ & G_{(1)}( (\cN' \cap \cA_{\tau'\alpha'}^{(2)}) \times S_0),\\
G_{(2)}( (\cC_2\setminus (\cC_1\cap\cC_2)) \times S_0) \subset\ & \cup_{j_0 = 1}^{4} h_{j_0}([\nu_2\alpha', \nu_{3}\alpha')\times \partial\widetilde{M}_{j_0}) \subset \widetilde{M}_{j_0}.
\end{split}
\]
Hence we deduce from Proposition~\ref{prop:smoothing-in-distance-nbhd}(b) and property (4) below~\eqref{eq:collar-5-strata} that $(\cC_1\setminus (\cC_1\cap\cC_2)) \times S_0$ and $(\cC_2\setminus (\cC_1\cap\cC_2)) \times S_0$ have disjoint images under $G_{(2)}$. Taking into account what we proved in the previous paragraph, we conclude that $G_{(2)}$ is an injective local diffeomorphism on $(\cC_1\cup \cC_2)\times S_0 = \mathring{\cD}'' \times S_0$. This proves part (a).

For part (b), we first use the definition of $G_{(2)}$ to see that
\begin{equation}\label{eq:G-2-properties-b-initial-splitting}
\begin{split}
G_{(2)}(\cN'' \setminus \mathring{\cD}'') \times S_0) \subset\ & G_{(1)}((\cN' \setminus \mathring{\cN}'_{7/8, \mu_{4}\alpha})\times S_0) \\
&\cup \big(\cup_{j_0 = 1}^{4}h_{j_0}((\nu_1\alpha', \nu_4\alpha'] \times \partial\widetilde{M}_{j_0}) \big).
\end{split}
\end{equation}
Recalling from~\eqref{eq:mu1-mu0-inequalities} and~\eqref{eq:transition-parameter} that $\mu_0 < \tau < 1$, and using also~\eqref{eq:alpha-choice}, we get
\begin{equation}\label{eq:ball-mu0-alpha}
B_{\frac{\mu_0\alpha}{2}} \subset \mathring{\cN}_{7/8} \cap \cA_{\tau\alpha} \subset \mathring{\cN}'_{7/8, \mu_{4}\alpha},
\end{equation}
which together with Proposition~\ref{prop:G-1-properties}(a) implies
\begin{equation}\label{eq:G-2-properties-b-first-piece}
G_{(1)}((\cN' \setminus \mathring{\cN}'_{7/8, \mu_{4}\alpha})\times S_0)  \cap G_{(1)}(B_{\frac{\mu_0\alpha}{2}} \times S_0) = \emptyset. 
\end{equation}
On the other hand, by~\eqref{eq:5-strata-smoothing-collar-deleted} and~\eqref{eq:M-j-star-disjoint-from-image}, we have for each $j_0 \in \{1, 2, 3, 4\}$ that
\begin{equation}\label{eq:G-2-properties-b-second-piece}
h_{j_0}((\nu_1\alpha', \nu_4\alpha'] \times \partial\widetilde{M}_{j_0}) \big) \cap G_{(1)}((\cN' \cap \cA^{(2)}_{\frac{\nu_1\alpha'}{2}}) \times S_0) = \emptyset.
\end{equation}
Next, noting from~\eqref{eq:rho2-choice} and~\eqref{eq:alpha'-threshold} that $\nu_1\alpha' < \alpha' < \frac{\mu_0\alpha}{16}$, and using~\eqref{eq:ball-mu0-alpha}, we get
\begin{equation}\label{eq:G-2-properties-b-inclusion}
\begin{split}
B_{\frac{\nu_1\alpha'}{2}}  \subset  B_{\frac{\mu_0\alpha}{2}} \cap \cA^{(2)}_{\frac{\nu_1\alpha'}{2}} =\ &  B_{\frac{\mu_0\alpha}{2}} \cap \mathring{\cN}'_{7/8, \mu_{4}\alpha} \cap \cA^{(2)}_{\frac{\nu_1\alpha'}{2}}, 
\end{split}
\end{equation}
and hence
\[
G_{(1)}(B_{\frac{\nu_1\alpha'}{2}} \times S_0) \subset G_{(1)}(B_{\frac{\mu_0\alpha}{2}} \times S_0) \cap G_{(1)}((\cN' \cap \cA^{(2)}_{\frac{\nu_1\alpha'}{2}}) \times S_0),
\]
which together with~\eqref{eq:G-2-properties-b-second-piece},~\eqref{eq:G-2-properties-b-first-piece}, and~\eqref{eq:G-2-properties-b-initial-splitting} gives
\[
G_{(2)}(\cN'' \setminus \mathring{\cD}'') \times S_0) \cap G_{(1)}(B_{\frac{\nu_1\alpha'}{2}} \times S_0) =\emptyset.
\]
Upon using~\eqref{eq:G-2-properties-b-inclusion}, the inequality $\nu_1 < \tau'$, and the definition~\eqref{eq:G-2-definition-case-1} to see that 
\[
G_{(1)}(B_{\frac{\nu_1\alpha'}{2}} \times S_0) = G_{(2)}(B_{\frac{\nu_1\alpha'}{2}} \times S_0),
\]
and that $B_{\frac{\nu_1\alpha'}{2}} \subset \mathring{\cD}''$, while recalling also that $G_{(2)}$ is injective on $\mathring{\cD}'' \times S_0$, we get part (b).

Moving on to part (c), given $y \in \cN'' \setminus \mathring{\cN}''$, as in the proof of Proposition~\ref{prop:G-1-properties}, there are three cases: First, if $y_0 \in (\cN' \setminus \mathring{\cN}') \cap \cA^{(2)}_{\tau'\alpha'}$, then by~\eqref{eq:G-2-definition-case-1} and Proposition~\ref{prop:G-1-properties}(b) we have 
\[
\cH^{1}(G_{(2)}(\{y_0\} \times S_0)) = \cH^{1}( G_{(1)}(\{y_0\} \times S_0))< \infty.
\]
Now let $(y_i)$ be a sequence in $\mathring{\cN}''$ converging to $y_0$. Then by the second inclusion in~\eqref{eq:N''-contained-in-N'} and the openness of $\cA_{\tau'\alpha'}^{(2)}$, eventually all $y_i$ belong to $\mathring{\cN}' \cap \cA_{\tau'\alpha'}^{(2)}$. Since $(y_i)$ converges to a limit in $\cN' \setminus \mathring{\cN}'$, we conclude from~\eqref{eq:G-2-definition-case-1} and~\eqref{eq:G-1-area-boundary} that 
\[
\Area(G_{(2)}(y_{i}, \cdot)) = \Area(G_{(1)}(y_{i}, \cdot)) \to 0 \quad\text{as }i \to \infty,
\]
which finishes the proof in this case.

The second case is when $y_0 =\Psi_{j_0}(s_0, z_0)$ for some $j_0 \in \{1, 2, 3, 4\}$, $s_0 \in (\nu_1\alpha', \nu_4\alpha')$, and $z_0 \in \partial\Omega_{j_0} \cap (\cN' \setminus \mathring{\cN}')$. Here we have
\[
G_{(2)}(\{y_0\} \times S_0) = h_{j_0}(\{s_0\} \times G_{(1)}(\{z_0\} \times S_0)), 
\]
which has finite $\cH^{1}$-measure since $G_{(1)}(\{z_0\} \times S_0)$ does by Proposition~\ref{prop:G-1-properties}(b), and since $h_{j_0}(s_0, \cdot)$ is Lipschitz. Next suppose $(y_i)$ is a sequence converging to $y_0$ from within $\mathring{\cN}''$. Using the fact that $\mathring{V}^{3}_{j_0}$ is open, and is disjoint from $\mathring{V}^{3}_{j}$ if $j \neq j_0$, we see from~\eqref{eq:N''-limbs} and the definition of $\mathring{\cN}''$ that, up to taking a subsequence, either $y_i \in \mathring{\cN}' \cap \cA_{\tau'\alpha'}^{(2)}$ for all $i$, or there is a sequence of points $(s_i, z_i)$ in $(\nu_1\alpha', \nu_4\alpha') \times (\partial\Omega_{j_0} \cap \mathring{\cN}')$ converging to $(s_0, z_0)$ such that
\[
y_i = \Psi_{j_0}(s_i, z_i) \quad\text{for all }i.
\]
In the former situation, since $z_0 \in \cN' \setminus \mathring{\cN}'$, we obtain from Lemma~\ref{lemm:Psi-action-on-N'}(c) that $y_0 \in \cN' \setminus \mathring{\cN}'$ as well, and therefore, by~\eqref{eq:G-2-definition-case-1} and~\eqref{eq:G-1-area-boundary}, 
\[
\Area(G_{(2)}(y_i, \cdot)) = \Area(G_{(1)}(y_i, \cdot)) \to 0\quad\text{as }i \to \infty.
\]
In the latter situation, we use $z_0 \in \cN' \setminus \mathring{\cN}'$ and~\eqref{eq:G-1-area-boundary} again, this time to see that
\[
\Area(G_{(1)}(z_i, \cdot)) \to 0\quad \text{as }i \to \infty.
\]
Because $G_{(1)}(z_i, \cdot)$ all take values in $\partial\widetilde{M}_{j_0}$ by Proposition~\ref{prop:M-j-as-G-1-image}, and because $h_{j_0}(s_i, \cdot)$ converges to $h_{j_0}(s_0, \cdot)$ in $C^{1}$ on $\partial\widetilde{M}_{j_0}$, we conclude that
\[
\Area(G_{(2)}(y_i, \cdot)) = \Area(h_{j_0}(s_i, G_{(1)}(z_i, \cdot))) \to 0\quad\text{as }i \to \infty.
\]

Finally, the third case is when $y_0 = \Psi_{j_0}(\nu_{4}\alpha', z_0)$ for some $j_0 \in \{1, 2, 3, 4\}$ and $z_0 \in \partial\Omega_{j_0} \cap \cN'$. By property (2) mentioned after~\eqref{eq:collar-5-strata}, we have
\[
\cH^{1}(G_{2}(\{y_0\} \times S_0)) = \cH^{1}(h_{j_0}(\{\nu_4\alpha'\} \times G_{(1)}(\{z_0\} \times S_0))) < \infty.
\]
By~\eqref{eq:beyond-nu2-disjoint} and~\eqref{eq:N''-limbs}, any sequence $(y_i)$ in $\mathring{\cN}''$ that converges to $y_0$ eventually has the form $y_i = \Psi_{j_0}(s_i, z_i)$, where $(s_i, z_i) \in (\nu_1\alpha', \nu_4\alpha') \times (\partial\Omega_{j_0} \cap \mathring{\cN}')$, and
\[
s_i \to (\nu_4\alpha')^{-},\quad z_i \to  z_0, \quad\text{as }i \to \infty.
\]
Thus, for all large enough $i$ there holds
\[
\Area(G_{(2)}(y_i, \cdot)) = \Area(h_{j_0}(s_i, G_{(1)}(z_i, \cdot))) \leq \|\Lambda^2 dh_{j_0, s_i}\| \cdot \Area(G_{(1)}(z_i, \cdot)),
\]
and we arrive at~\eqref{eq:G-2-area-boundary} upon recalling~\eqref{eq:h-j0-shrink-area}, and noting from Proposition~\ref{prop:G-1-properties} that $\Area(G_{(1)}(z_i, \cdot))$ remains bounded whether or not $z_0 \in \mathring{\cN}'$. The proof is complete.
\end{proof}
\section{Convergence results for quasiconformal maps}\label{appendix:quasiconformal}
In this appendix, we first prove Lemmas~\ref{lemm:non-vanishing-Jacobian} and~\ref{lemm:qc-standard-properties}, and then move on to prove (a part of) Theorem~\ref{thm:qc-dependence}, along with Propositions~\ref{prop:qc-H-dependence},~\ref{prop:qc-smooth-dependence} and~\ref{prop:metric-inverse-convergence}. We begin by reviewing certain elements of the work of Ahlfors and Bers~\cite{Ahlfors-Bers} on the Beltrami equation. Given $p > 2$, following~\cite{Ahlfors-Bers} we define $\cB_p$ to be the space of continuous functions $f: \CC \to \CC$ such that $f(0) = 0$, that the distributional derivatives $f_{z}, f_{\overline{z}}$ lie in $L^{p}_{\loc}$, and that
\[
\|f\|_{\cB_p} : = [f]_{1 - \frac{2}{p}; \CC} + \|f_z\|_{p; \CC} + \|f_{\overline{z}}\|_{p; \CC} < \infty.
\]
Note that $(\cB_{p}, \|\cdot\|_{\cB_p})$ is a Banach space. Also, if $f$ is a continuous function on $\CC$ such that $f(0) = 0$ and $\nabla f \in L^p_{\loc}$, then $\zeta f \in \cB_p$ for any $\zeta \in C^\infty_{c}(\CC)$, thanks to Morrey's embedding~\cite[page 283, second displayed equation]{Evans}. 

\begin{lemm}[\cite{Ahlfors-Bers}, Lemmas 3 and 4]
\label{lemm:CZ}
For $f \in C^{\infty}_{c}(\CC)$, define
\[
\begin{split}
(Pf)(z) :=\ & \frac{1}{2\pi i} \int_{\CC} f(\zeta) \big( \frac{1}{\zeta - z} - \frac{1}{\zeta} \big) d\zeta \wedge d\overline{\zeta},\\
(Tf)(z) :=\ & \frac{1}{2\pi i} \lim_{\ep \to 0} \int_{\CC \setminus \bB_{\ep}(z)} \frac{f(\zeta)}{(\zeta - z)^2} d\zeta \wedge d\overline{\zeta}.
\end{split}
\]
\vskip 1mm
\begin{enumerate}
\item[(a)] $T$ extends to a bounded linear operator from $L^p(\CC)$ to itself for all $2 \leq p < \infty$. Moreover, the operator norms satisfy 
\[
\lim_{p \to 2}\|T\|_{p} = 1.
\]
\vskip 1mm
\item[(b)] $P$ extends to a bounded linear operator from $L^p(\CC)$ to $\cB_p$ for all $p > 2$. Also, in the sense of distributions we have
\[
(Pf)_{\overline{z}} = f,\quad (Pf)_z = Tf.
\] 
\end{enumerate}
\end{lemm}
\begin{proof}
This is essentially a restatement of~\cite[Lemmas 3 and 4]{Ahlfors-Bers}. See also \cite[Lemmas 4.20 and 4.21, and Propositions 4.22 and 4.23]{Imayoshi-Taniguchi}.
\end{proof}
\begin{lemm}[\cite{Ahlfors-Bers}, Theorem 1]
\label{lemm:inhomogeneous-Beltrami} Suppose $\mu \in L^{\infty}(\CC)$ and that $\|\mu\|_{\infty; \CC} < 1$. Taking any $p > 2$ such that 
\begin{equation}\label{eq:choice-of-p}
\|\mu\|_{\infty; \CC} \cdot \|T\|_p  < 1,
\end{equation}
which is possible by Lemma \ref{lemm:CZ}(a), we have the following.
\begin{enumerate}
\item[(a)] $w = P(w_{\overline{z}})$, for all $w \in \cB_p$.
\vskip 1mm
\item[(b)] Given $\sigma \in L^p(\CC)$, there exists a unique solution $w = w^{\mu, \sigma}$ in $\cB_{p}$ of the equation 
\[
w_{\overline{z}} = \mu w_{z} + \sigma.
\]
Moreover, this solution satisfies the estimate
\begin{equation}\label{eq:Beltrami-estimate}
\|w\|_{\cB_p} \leq C \|\sigma\|_{p; \CC},
\end{equation}
where $C$ depends only on $p$ and $\|\mu\|_{\infty; \CC}$.
\vskip 1mm
\item[(c)] Let $(\sigma_n)$ be a sequence converging to some $\sigma$ in $L^{p}(\CC)$, and $(\mu_n)$ a sequence in $L^{\infty}(\CC)$ that converges to $\mu$ almost everywhere, and satisfies $\|\mu_n\|_{\infty; \CC} \leq \|\mu\|_{\infty; \CC}$ for each $n$. Then 
\[
\|w^{\mu_n, \sigma_n} - w^{\mu, \sigma}\|_{\cB_{p}}  \to 0 \quad\text{as }n \to \infty.
\]
\end{enumerate}
\end{lemm}
\begin{proof}
For part (a), by Lemma~\ref{lemm:CZ}(b) and the definition of $\cB_{p}$, we see that $w - P(w_{\overline{z}})$ is holomorphic on $\CC$, vanishes at the origin, and satisfies for all $z \in \CC$ that
\[
|w(z) - P(w_{\overline{z}})(z)| \leq [w - P(w_{\overline{z}})]_{1 - \frac{2}{p}}\cdot |z|^{1 - \frac{2}{p}} \leq C\|w\|_{\cB_{p}}\cdot|z|^{1 - \frac{2}{p}}.
\]
Elementary complex analysis then shows that $w - P(w_{\overline{z}})$ is identically zero. For part (b), 
the existence and uniqueness statement is \cite[Theorem 1]{Ahlfors-Bers}, while the estimate~\eqref{eq:Beltrami-estimate} is stated in~\cite[Lemma 5]{Ahlfors-Bers}. Part (c) is~\cite[Lemma 6]{Ahlfors-Bers}. 
\end{proof}
\begin{lemm}[\cite{Ahlfors-Bers}, Section 3]
\label{lemm:normal-solution}
Let $\mu$ be a bounded measurable function on $\CC$ which satisfies $\|\mu\|_{\infty} < 1$ and has compact essential support. Then,  the equation
\begin{equation}\label{eq:normal-solution-eq}
F_{\overline{z}} = \mu F_{z}
\end{equation}
has a unique solution in $(C^{0} \cap W^{1, 2})_{\loc}(\CC)$ subject to the conditions that 
\begin{equation}\label{eq:normal-solution-integrability}
F_z - 1 \in L^{p}(\CC) \quad\text{for \emph{some} }p > 2 \text{ such that }\|\mu\|_{\infty} \cdot \|T\|_{p} < 1,
\end{equation}
and that $F(0) = 0$. This solution, denoted $F^{\mu}$, has the following additional properties.
\vskip 1mm
\begin{enumerate}
\item[(a)] $\nabla F^{\mu} \in L^{p}_{\loc}(\CC)$ for \emph{all} $p > 2 $ such that $\|\mu\|_{\infty} \cdot \|T\|_{p} < 1$.
\vskip 1mm
\item[(b)] $F^{\mu}$ is an orientation-preserving homeomorphism.
\vskip 1mm
\item[(c)] $F^{\mu}_{z}\neq 0$ almost everywhere. Also, both $F^{\mu}$ and $(F^{\mu})^{-1}$ take sets of measure zero to sets of measure zero. Consequently, they both take measurable sets to measurable sets.
\vskip 1mm
\item[(d)] Given $q \geq 2$, a domain $\Omega \subset \CC$, and a function $h \in (C^0 \cap W^{1, q})_{\loc}(F^{\mu}(\Omega))$, we have 
\[
\nabla (h \circ F^{\mu}) \in L^{\frac{pq}{p+q-2}}_{\loc}(\Omega),
\]
for all $p > 2$ such that $\|\mu\|_{\infty} \cdot \|T\|_{p} < 1$. Moreover the chain rule holds; that is
\[
(h\circ F^{\mu})_{z} = (h_{z} \circ F^{\mu}) F^{\mu}_{z} + (h_{\overline{z}} \circ F^{\mu}) (\overline{F^{\mu}})_{z},
\]
and similarly for the $\overline{z}$-derivative.
\vskip 1mm
\item[(e)] Given another compactly supported $\mu_{1} \in L^{\infty}(\CC)$ with $\|\mu_1\|_{\infty} < 1$, we have $F^{\mu_1} \circ (F^{\mu})^{-1} = F^{\lambda}$, where $\lambda$ is given by
\begin{equation}\label{eq:normal-solution-composition}
\lambda: = \Big( \frac{\mu_1 - \mu}{1 - \overline{\mu} \mu_1} \cdot \frac{F^{\mu}_{z}}{(\overline{F^{\mu}})_{\overline{z}}} \Big) \circ (F^{\mu})^{-1},
\end{equation}
and satisfies the estimate
\begin{equation}\label{eq:complex-dilatation-growth}
\|\lambda\|_{\infty} \leq \frac{\|\mu\|_{\infty} + \|\mu_{1}\|_{\infty}}{1 + \|\mu\|_{\infty}\cdot \|\mu_{1}\|_{\infty}} < 1.
\end{equation}
\vskip 1mm
\item[(f)] If in addition $\mu \in C^{1}$, then $F^{\mu}$ is an orientation-preserving $C^1$-diffeomorphism.
\end{enumerate}
\end{lemm}
\begin{proof}[Proof of uniqueness]
Suppose, for $a = 1, 2$, that $F_{a} \in (C^0 \cap W^{1, 2})_{\loc}(\CC)$ is a solution to~\eqref{eq:normal-solution-eq} such that $F_{a}(0) = 0$, and that 
\begin{equation}\label{eq:integrability-near-infty}
(F_{a})_{z} - 1 \in L^{p_{a}}(\CC),
\end{equation}
where $p_a > 2$ satisfies~\eqref{eq:choice-of-p}. By our assumption on $\mu$, each $F_a$ is holomorphic near $\infty$, and thus so is $(F_{a})_{z} - 1$. The mean-value property and~\eqref{eq:integrability-near-infty} then implies that 
\begin{equation}\label{eq:normal-solution-decay-1}
(F_{a})_{z}  - 1 \to 0 \quad\text{as }|z| \to \infty,
\end{equation}
and consequently $\lim_{|\zeta| \to \infty} \zeta\cdot ((F_{a})_{z}(\zeta) - 1)$ exists. From this, and assuming without loss of generality that $p_1 \leq p_2$, we infer that the difference
\[
f := F_{1} - F_{2},
\]
which is a solution to~\eqref{eq:normal-solution-eq}, satisfies $f_{z}, f_{\overline{z}} \in L^{p_1}(\CC)$. Since $f$ is already continuous and vanishes at the origin, we see from Morrey's embedding that $f \in \cB_{p_1}$, in which case the uniqueness part of Lemma~\ref{lemm:inhomogeneous-Beltrami}(b) forces $f \equiv 0$. This proves uniqueness.
\end{proof}
\begin{proof}[Proof of existence]
Fix any $p > 2$ satisfying~\eqref{eq:choice-of-p}, and notice that $\mu$ belongs to $L^{p}(\CC)$ since it is bounded and has compact support. We then let $w$ be the unique solution in $\cB_p$ to $w_{\overline{z}} = \mu w_z + \mu$ produced by Lemma~\ref{lemm:inhomogeneous-Beltrami}, and define
\begin{equation}\label{eq:normal-solution}
F(z) = z + w(z).
\end{equation}
Clearly $F$ is of class $C^0 \cap W^{1, 2}$ locally on $\CC$, vanishes at the origin, and solves the equation~\eqref{eq:normal-solution-eq}. Also, by construction we have $F_{z} - 1 = w_{z} \in L^{p}(\CC)$. This proves existence. For later use, we note from~\eqref{eq:normal-solution} and the definition of $\cB_{p}$ that 
\begin{equation}\label{eq:normal-solution-decay-2}
\lim_{|z| \to \infty} |z|^{-1}|F(z) - z| = 0
\end{equation}
\end{proof}
\begin{proof}[Proof of part (f)] 
This is~\cite[Lemma 7]{Ahlfors-Bers}. We sketch the idea for the convenience of the reader. Also, this part does not depend on any of the preceding five. We state it last only to make it clear that parts (a) through (e) do not require $\mu$ to be $C^{1}$.

Fixing an arbitrary $p > 2$ satisfying~\eqref{eq:choice-of-p}, the second paragraph of the proof of~\cite[Lemma 7]{Ahlfors-Bers} produces a $C^1$-solution $f$ to~\eqref{eq:normal-solution-eq} with $f_{z}$ having the form
\begin{equation}\label{eq:f-z-exp}
f_{z} = e^{\lambda},
\end{equation}
where $\lambda$ is an element of $\cB_{p}$ which, among other things, is holomorphic near $\infty$. The definition of $\|\cdot\|_{\cB_{p}}$ implies that $z^{-1}\lambda(z) \to 0$ as $|z| \to \infty$, and thus there exists $a_0 \in \CC$ such that $\lim_{|z| \to \infty} z \cdot (\lambda(z) - a_0)$ exists, from which we deduce that 
\[
e^{-a_0}f_{z} - 1 \in L^{p}(\CC),
\]
and hence $e^{-a_0}(f - f(0))$ coincides with $F^{\mu}$. This together with~\eqref{eq:f-z-exp} and the equation~\eqref{eq:normal-solution-eq} shows that $F^{\mu}$ is $C^{1}$, and moreover that
\begin{equation}\label{eq:normal-solution-positive-Jacobian}
|F^{\mu}_{z}|^2 - |F^{\mu}_{\overline{z}}|^2 \geq (1 - \|\mu\|_{\infty}^2)|F^{\mu}_{z}|^2 > 0 \quad\text{everywhere}.
\end{equation}
Since $F^{\mu}$ is also proper by~\eqref{eq:normal-solution-decay-2}, we conclude that it is a covering map from $\CC$ onto itself, and hence a global diffeomorphism since $\CC$ is simply-connected.
\end{proof}
\begin{proof}[Proof of remaining parts]
Part (a) follows from the proof of existence, and the uniqueness statement. 

Part (b) is~\cite[Lemma 8]{Ahlfors-Bers}. We describe the framework of the proof to make an observation. To start, again fix $p > 2$ such that~\eqref{eq:choice-of-p} holds. By mollifying $\mu$, we obtain a sequence $(\mu_n)$ of $C^{1}$-functions that satisfies $\|\mu_n\|_{\infty} \leq \|\mu\|_{\infty}$ and approximates $\mu$ in the sense indicated at the start of~\cite[Section 3.3]{Ahlfors-Bers}. Then by~\eqref{eq:normal-solution} and Lemma~\ref{lemm:inhomogeneous-Beltrami}(c), we see that $F^{\mu_n}  - F^{\mu} \in \cB_{p}$, and that
\begin{equation}\label{eq:Fn-F-in-Bp}
\lim_{n \to \infty}\|F^{\mu_n} - F^{\mu}\|_{\cB_{p}} = 0.
\end{equation}
On the other hand, part (f) guarantees that $(F^{\mu_n})^{-1}$ exists for each $n$. We then observe that, by the characterization of $(F^{\mu_n})^{-1}$ given right after~\cite[page 392, equation (25)]{Ahlfors-Bers}, along with the estimate on the line that follows, the sequence of inverses is equicontinuous and uniformly bounded on compact subsets of $\CC$. Taking~\eqref{eq:Fn-F-in-Bp} into account, we conclude that $(F^{\mu_n})^{-1}$ converges in $C^{0}_{\loc}(\CC)$, and that the limit inverts $F^{\mu}$. In particular $F^{\mu}$ is a homeomorphism. That it preserves orientation follows from part (f) and the local uniform convergence implied by~\eqref{eq:Fn-F-in-Bp}.

Next, part (c) is contained in Lemma 9 of~\cite{Ahlfors-Bers}, the proof of which depends on the convergence of $(F^{\mu_n})^{-1}$ to $(F^{\mu})^{-1}$ noted above. Part (d) is Lemma 10. 

For part (e), it follows from part (c) that $\lambda$ is measurable, and vanishes almost everywhere outside of a compact set. Also, letting
\[
k:= \|\mu\|_{\infty}, \quad k_1 := \|\mu_{1}\|_{\infty},
\]
then by a direct computation, we have
\[
1 - |\lambda \circ F^{\mu}|^2 = 1 - \big| \frac{\mu_1 - \mu}{1 - \overline{\mu}\mu_1} \big|^2 \geq \frac{(1 - k_1^2)(1-k^2)}{(1 + k_1k)^2},
\]
from which~\eqref{eq:complex-dilatation-growth} follows easily. Next, by Lemma~\ref{lemm:CZ}(a), there is $p_0 > 2$ such that
\begin{equation}\label{eq:choice-of-p0}
\frac{k+k_1}{1 + k k_1} \cdot \|T\|_{p} < 1\quad\text{for all }p \in (2, p_0].
\end{equation}
As $\frac{k+k_1}{1 + k k_1} > \max\{k, k_1\}$, we have by part (a) that both $\nabla F^{\lambda}$ and $\nabla F^{\mu}$ lie in $L^{p}_{\loc}$ for all $p \in (2, p_0]$. Applying part (d) with $h = F^{\lambda}$ and $q = p = p_0$ yields some $p_1 \in (2, p_0)$ such that 
\begin{equation}\label{eq:composition-locally-Lp}
\nabla (F^{\lambda} \circ F^{\mu}) \in L^{p_1}_{\loc}(\CC).
\end{equation}
Substituting the equation~\eqref{eq:normal-solution-eq}, as well as the version for $F^{\lambda}$, into the formulas produced by the chain rule, and recalling the definition of $\lambda$, we find that
\[
(F^{\lambda} \circ F^{\mu})_{\overline{z}} = \mu_1 \cdot (F^{\lambda} \circ F^{\mu})_{z}.
\]
Noting that~\eqref{eq:normal-solution-decay-2} and~\eqref{eq:normal-solution-decay-1} apply to both $F^{\mu}$ and $F^{\lambda}$, we see, again using the formula for $(F^{\lambda} \circ F^{\mu})_{z}$ from the chain rule, that 
\[
(F^{\lambda} \circ F^{\mu})_{z} \to 1\quad\text{as }|z| \to \infty.
\]
Since $\mu_1$ has compact support, and hence $F^{\lambda} \circ F^{\mu}$ is holomorphic near $\infty$, this latter asymptotic behavior implies that $(F^{\lambda} \circ F^{\mu})_{z} = 1 + O(|z|^{-1})$ as $|z| \to \infty$, which combines with~\eqref{eq:composition-locally-Lp} to give 
\[
(F^{\lambda} \circ F^{\mu})_{z}  - 1 \in L^{p_1}(\CC).
\]
Recalling from~\eqref{eq:choice-of-p0} that $k_1 \|T\|_{p_1} < 1$, we conclude that $F^{\lambda} \circ F^{\mu}$, which clearly vanishes at the origin, must coincide with $F^{\mu_1}$. This finishes part (e).
\end{proof}
\begin{rmk}[\cite{Ahlfors-Bers}, page 387]
\label{rmk:choice-of-p}
Below we assume that for each $k \in [0, 1)$, a certain $p = p(k) > 2$ is chosen so that $k \|T\|_{p} < 1$. This is how the exponent $p(k)$ in the statement of Lemma \ref{lemm:qc-standard-properties}, Theorem~\ref{thm:qc-existence}, and Theorem~\ref{thm:qc-existence-H}, among other places, is to be understood.
\end{rmk}
We now describe the proof of Lemma~\ref{lemm:non-vanishing-Jacobian} and Lemma~\ref{lemm:qc-standard-properties}. 
\begin{proof}[Proof of Lemma~\ref{lemm:non-vanishing-Jacobian}]
Define
\[
\mu = \left\{
\begin{array}{ll}
f_{\overline{z}}/f_{z},& \text{ if }f_z \neq 0,\\
0,& \text{ if }f_z = 0,
\end{array}
\right.
\]
so that $\|\mu\|_{\infty; \CC} \leq \frac{K-1}{K + 1} < 1$, and that~\eqref{eq:Beltrami-eq} holds. Next, given a disk $B \Subset \Omega$, we let $\zeta \in C^{\infty}_{c}(\Omega)$ be a cut-off function that equals $1$ on $\overline{B}$, and write $F$ for the map $F^{\zeta \mu}$ produced by Lemma~\ref{lemm:normal-solution}. By the characterization of $F^{-1}$ in Lemma~\ref{lemm:normal-solution}(e), followed by part (d) there, we see that 
\[
\varphi: = f \circ F^{-1}:F(\Omega) \to f(\Omega)
\]
is of class $C^{0} \cap W^{1, 2}$ locally on $F(\Omega)$. Thus we may apply Lemma~\ref{lemm:normal-solution}(d) again, this time to compute the derivatives of $f = \varphi \circ F$ on $\Omega$. From the resulting formulas along with our choice of $\zeta$, as well as Lemma~\ref{lemm:normal-solution}(c), we see that $\varphi_{\overline{z}} = 0$ almost everywhere on $F(B)$. Thus $\varphi: F(B) \to f(B)$, being a holomorphic bijection, is a conformal map in the classical sense. Moreover, we have
\begin{equation}\label{eq:nabla-f-chain-rule}
f_{z} = (\varphi_{z} \circ F) F_{z}, \quad f_{\overline{z}} = (\varphi_{z} \circ F) F_{\overline{z}}, 
\end{equation}
almost everywhere on $B$. From Lemma~\ref{lemm:normal-solution}(c), and the fact that $\varphi_{z}$ is non-vanishing on $F(B)$, we see that $f_{z} \neq 0$ almost everywhere on $B$. This finishes the proof.\\
\end{proof}

\begin{proof}[Proof of Lemma~\ref{lemm:qc-standard-properties}]
As in the previous proof, it suffices to establish the desired conclusions on disks in $\Omega$. Thus, we again take any closed disk $\overline{B} \subset \Omega$, and let $\zeta \in C^{\infty}_{c}(\Omega)$ be the same cut-off function as before. Writing $\mu$ for $\mu_{f}$, and letting $F = F^{\zeta \mu}$, we again have that 
\[
\varphi: = f \circ F^{-1}
\]
is a biholomorphic map from $F(B)$ onto $f(B)$, and that~\eqref{eq:nabla-f-chain-rule} holds. Part (a) then follows from its counterpart in Lemma~\ref{lemm:normal-solution}, since $\|\zeta \mu\|_{\infty} \leq k$. Next, as a consequence of the expression $f = \varphi \circ F$, along with Lemma~\ref{lemm:normal-solution}(c) and the smoothness of $\varphi$ on $F(B)$, we see that $|f(E \cap B)| = 0$ whenever $E \subset \Omega$ has measure zero. Likewise, writing $f^{-1} = F^{-1} \circ \varphi^{-1}$ shows that $f^{-1}(E)\cap B$ has measure zero whenever $E \subset f(\Omega)$ does. The disk $\overline{B} \subset \Omega$ being arbitrary, this proves (b).

For part (c), by Lemma~\ref{lemm:normal-solution}(e) we see that $F^{-1}$ is of class $C^0 \cap W^{1, 2}$ locally on $\CC$. Consequently, the chain rule in part (d) there can be applied to $F^{-1}\circ F$, giving us 
\[
(F^{-1})_{z} \circ F = \frac{1}{(1 - |\mu|^2)F_{z}}.
\]
Also, since $\varphi^{-1}:f(B) \to F(B)$ is a conformal map, we see in particular that $f^{-1} = F^{-1} \circ \varphi^{-1}$ has distributional derivatives in $L^{2}_{\loc}$ on $f(B)$, and that these can be computed with the chain rule. Doing so, and using the above expression for $(F^{-1})_{z}$ and the holomorphicity of $\varphi^{-1}$, we have almost everywhere on $B$ that
\[
\begin{split}
(f^{-1})_{z} \circ f = [((F^{-1})_{z} \circ \varphi^{-1}) \cdot (\varphi^{-1})_{z}] \circ f =\ & \frac{1}{(1 - |\mu|^2)F_{z} } \cdot \frac{1}{(\varphi_{z}\circ F)}\\
=\ &  \frac{1}{(1 - |\mu|^2)f_{z}},
\end{split}
\]
where for the last step we used~\eqref{eq:nabla-f-chain-rule}. This gives the first identity in~\eqref{eq:f-inverse-derivatives}, and the second can be derived similarly. Since $|\mu| \leq k < 1$ by assumption, we conclude from~\eqref{eq:f-inverse-derivatives} that $f^{-1}$ is quasiconformal. The asserted expression for $\mu_{f^{-1}}$ is also an immediate consequence.

For part (d), we introduce the abbreviations 
\[
k_1: = \|\mu_{g}\|_{\infty; \Omega}, \quad\quad h: = g \circ f^{-1},
\]
and also let $G: = F^{\zeta \mu_{g}}$, with $\zeta$ being the cut-off function from the beginning. Then, as before, $\psi: = g \circ G^{-1}$ is a biholomorphic map from $G(B)$ onto $g(B)$. From the expression
\[
h = \psi \circ G \circ F^{-1} \circ \varphi^{-1},
\]
along with Lemma~\ref{lemm:normal-solution}(e) applied to $G \circ F^{-1}$, we see that $h$ is of class $C^0 \cap W^{1, 2}$ locally on $f(B)$, and that its derivatives are related to those of $G\circ F^{-1}$ by the chain rule. Carrying out the computation, and taking into account also the formula~\eqref{eq:normal-solution-composition}, our choice of $\zeta$, and the estimate~\eqref{eq:complex-dilatation-growth}, we see that, almost everywhere on $f(B)$, there holds
\[
|h_{\overline{z}}| \leq \frac{k + k_1}{1 + k k_1}|h_{z}|,
\]
and hence $h$ is quasiconformal. Recalling also~\eqref{eq:nabla-f-chain-rule}, we obtain 
\[
\frac{h_{\overline{z}}}{h_{z}} = \Big( \frac{\mu_g - \mu_f}{1 - \overline{\mu_{f}}\mu_{g}} \Big)\circ f^{-1} \cdot \frac{F_{z}\circ f^{-1}}{\overline{F}_{\overline{z}} \circ f^{-1}} \cdot \frac{\varphi_{z} \circ \varphi^{-1}}{\overline{\varphi}_{\overline{z}} \circ \varphi^{-1}} =\Big( \frac{\mu_g - \mu_f}{1 - \overline{\mu_{f}}\mu_{g}} \Big)\circ f^{-1} \cdot \frac{f_{z} \circ f^{-1}}{\overline{f}_{\overline{z}} \circ f^{-1}},
\]
which is the desired formula~\eqref{eq:composition-Beltrami-coefficient}. This proves part (d).

Finally, if $\mu$ is $C^1$, then so is $\zeta \mu$, and we deduce from~\eqref{eq:nabla-f-chain-rule} and Lemma~\ref{lemm:normal-solution}(f) that $\nabla f$ is continuous on $B$, and that $f_z$ is non-vanishing. Letting $Jf$ denote the Jacobian determinant of $f$ regarded as a map between subsets of $\RR^2$, we have
\[
Jf = |f_{z}|^2 - |f_{\overline{z}}|^2 \geq (1 - k^2)|f_{z}|^2 > 0,\quad \text{everywhere on }B.
\]
Since $\overline{B} \subset \Omega$ is arbitrary, and since $f:\Omega \to f(\Omega)$ is a homeomorphism to start with, this gives the desired conclusion. The proof is complete.
\end{proof}
\begin{rmk}\label{rmk:mu-conformal-converse}
The method of proof for the previous lemma, combined with Lemma~\ref{lemm:normal-solution}(b), shows that any $\mu$-conformal homeomorphism in the sense of~\cite{Ahlfors-Bers} automatically preserves orientation, and are thus quasiconformal maps according to Definition~\ref{defi:quasiconformal-maps}.
\end{rmk}
\subsection{Proofs of Theorem~\ref{thm:qc-dependence} and Proposition~\ref{prop:qc-H-dependence}}\label{subsec:beltrami-proofs}
Before giving the promised proofs, we recall the following local $C^0\cap W^{1, p}$-estimate for the quasiconformal maps given by Theorem \ref{thm:qc-existence}.
\begin{prop}[\cite{Ahlfors-Bers}, equations (35) and (36)]
\label{prop:qc-estimates}
Given $\mu \in L^{\infty}(\CC)$ such that $k: = \|\mu\|_{\infty} < 1$, let $p = p(k)$ and $f^{\mu}:\CC \to \CC$ be as in Theorem~\ref{thm:qc-existence}. Then, for all $R > 0$ we have
\[
\|f^{\mu}\|_{\infty; \bB_R} + [f^\mu]_{1 - \frac{2}{p}; \bB_R} + \|\nabla f^\mu\|_{p; \bB_R}  \leq C,
\]
where $C$ depends only on $k$ and $R$.
\end{prop}
\begin{proof}[Proof of Theorem~\ref{thm:qc-dependence}]
For brevity we write $f_n = f^{\mu_n}$ and $f = f^{\mu}$. Assume by contradiction that for some $N > 0$, there exist $\alpha > 0$ and a subsequence of $f_n$, which we do not relabel, such that 
\begin{equation}\label{eq:qc-dependence-contra}
\|\nabla f_n - \nabla f\|_{p; \bB_N} +  \|f_n - f\|_{\infty; \bB_N} \geq \alpha,\quad \text{for all }n.
\end{equation}
Letting 
\[
\nu_n = -\Big(\mu_n \cdot \frac{(f_{n})_z}{(\overline{f_{n}})_{\overline{z}}}\Big) \circ f_{n}^{-1},
\]
we see from Lemma~\ref{lemm:qc-standard-properties}(c) and the uniqueness part of Theorem~\ref{thm:qc-existence} that
\[
h_n : = f_n^{-1} = f^{\nu_n}.
\]
The estimates in Proposition~\ref{prop:qc-estimates} are then applicable to both $f_n$ and $h_n$. Thus, up to taking a subsequence, we may assume that $f_n$ and $h_n$ converge in $C^0_{\loc}(\CC)$ to $\widetilde{f}$ and $\widetilde{h}$, respectively, and that $\nabla f_{n}$ converge weakly to $\nabla \widetilde{f}$ in $L^{p}_{\loc}(\CC)$. 

We next show that $\widetilde{f}$ coincides with $f$. To that end, note that $f_n \circ h_n$ converges uniformly locally to $\widetilde{f} \circ \widetilde{h}$, so the latter coincides with the identity map. For the same reason, $\widetilde{h} \circ \widetilde{f} = \id$, and thus $\widetilde{f}$ is a homeomorphism. On the other hand, since $\nabla f_n \to \nabla \widetilde{f}$ weakly in $L^p_{\loc}$ while $\mu_n \to \mu$ almost everywhere, with the help of Egorov's theorem we have
\[
\mu_n\cdot f_{n, z} \to \mu \cdot \widetilde{f}_z \text{ as distributions,}
\]
from which it follows that
\[
(\widetilde{f})_{\overline{z}} = \mu \cdot \widetilde{f}_z \text{ almost everywhere on } \CC.
\]
Thus we conclude that $\widetilde{f}$ is quasiconformal, having already seen that it is a homeomorphism with distributional derivatives in $L^2_{\loc}$. To check the normalization condition \eqref{eq:qc-normalization}, note that since $f_n(0) = 0$ and $f_n(1) = 1$ for all $n$, we have $\widetilde{f}(0) = 0$ and $\widetilde{f}(1) = 1$. Moreover, for all $R > 0$, Proposition~\ref{prop:qc-estimates} gives
\[
h_n(\bB_R) \subset \bB_{C_{k, R}} \text{ for all }n,
\]
which passes to the limit to give the same bound with $h_n$ replaced by $\widetilde{h}$, and hence $|\widetilde{f}(z)| \to \infty$ as $|z| \to \infty$. The uniqueness part of Theorem~\ref{thm:qc-existence} now forces $\widetilde{f} = f$, as asserted above. In particular,
\begin{equation}\label{eq:contradict-uniform}
\lim_{n \to \infty}\|f_n - f\|_{\infty; \bB_N} = 0.
\end{equation}
To obtain a contradiction to~\eqref{eq:qc-dependence-contra}, it remains to improve the weak convergence of $\nabla f_n$ to strong convergence. For that, we note
\[
(f_n - f)_{\overline{z}} = \mu_n f_{n, z} - \mu f_z = (\mu_n - \mu)f_{z} + \mu_n(f_{n} - f)_{z}.
\]
Thus, letting $\zeta$ be a cut-off function that equals $1$ on $\bB_N$ and vanishes outside of $\bB_{2N}$, we have
\begin{equation}\label{eq:difference-cutoff-Beltrami}
[\zeta(f_n - f)]_{\overline{z}} = \mu_n [\zeta (f_n - f)]_z + \sigma_n,
\end{equation}
where
\[
\sigma_n = (\zeta_{\overline{z}} - \mu_n \zeta_z)(f_n - f) + \zeta (\mu_n - \mu) f_{z}.
\]
Since $\zeta(f_n - f) \in \cB_p$, the a priori estimate in Lemma~\ref{lemm:inhomogeneous-Beltrami}(b) applied to~\eqref{eq:difference-cutoff-Beltrami} gives
\[
\| \zeta(f_n - f) \|_{\cB_p} \leq C_{k, p}\|\sigma_n\|_{p; \CC}.
\]
Noting that $\lim_{n \to \infty}\|\sigma_n\|_{p; \CC} =0$ thanks to~\eqref{eq:contradict-uniform} along with the dominated convergence theorem and the fact that $f_{z} \in L^p(\bB_{2N})$, we deduce from the definition of the $\|\cdot\|_{\cB_p}$-norm and our choice of $\zeta$ that
\[
\lim_{n \to \infty}\|\nabla f_n - \nabla f\|_{p; \bB_N}=0,
\]
which together with~\eqref{eq:contradict-uniform} yields a contradiction to~\eqref{eq:qc-dependence-contra}. The proof is complete.\\
 \end{proof}
Next we turn to Proposition~\ref{prop:qc-H-dependence}. We first prove a preliminary $W^{2, p}$-regularity result with estimates.
\begin{lemm}\label{lemm:inverse-W2p}
Given $\mu \in L^{\infty}(\CC)$ with $k := \|\mu\|_{\infty; \CC} < 1$, let $f^{\mu}$ and $p = p(k)$ be as in Theorem~\ref{thm:qc-existence} and write $h^{\mu} = (f^{\mu})^{-1}$. Suppose further that for some convex domain $\Omega \subset \CC$ we have $\mu \in C^1_{\loc}(\Omega)$. Then $h^{\mu}\in W^{2, p}_{\loc}(f(\Omega))$. Moreover, for any disk $\overline{B} \subset f(\Omega)$ and cut-off function $\zeta \in C^{\infty}_{c}(B)$ the following inequality holds:
\begin{equation}\label{eq:inverse-W2p-bound}
\|\zeta \nabla h^{\mu}\|_{1, p} \leq C(\|\zeta\|_{\infty} + \|\nabla\zeta\|_{\infty})\|\nabla h^{\mu}\|_{p; B} + C\|D\mu\|_{\infty; h^{\mu}(\overline{B})} \|\nabla h^{\mu}\|_{p; B} \|\zeta\nabla h^{\mu}\|_{\infty},
\end{equation}
where both constants $C$ on the right-hand side depend only on $k$ and $p$.
\end{lemm}
\begin{proof}
We write $f$ and $h$ for $f^{\mu}$ and $h^\mu$, respectively. By Lemma~\ref{lemm:qc-standard-properties}(e), we have $h \in C^{1}_{\loc}(f(\Omega))$. To obtain $W^{2, p}$-regularity, we employ a standard difference quotient argument. Specifically, for $t \in \RR\setminus \{0\}$, we define the following horizontal shift and difference quotient operators:
\[
g_{[t]}(z) = g(z + t),\ \ \ D_{[t]}g = \frac{g_{[t]} - g}{t}.
\]
With the help of the equation \eqref{eq:conjugate-Beltrami}, we see that for all $t \neq 0$ with $|t|$ sufficiently small, we have on the set $\{z \in f(\Omega)\ |\ \dist(z, f(\partial \Omega)) > |t|\}$ that 
\begin{equation}\label{eq:diff-q-Beltrami}
(D_{[t]}h)_{\overline{z}} = -D_{[t]}(\mu\circ h) \cdot (\overline{h}_{\overline{z}})_{[t]} - (\mu \circ h) \cdot (D_{[t]}\overline{h})_{\overline{z}}.
\end{equation}
Next, given $z_0 \in f(\Omega)$, we choose $r > 0$ such that $\overline{\bB_{4r}(z_0)} \subset f(\Omega)$, and define 
\[
K_{z_0,r} = \text{ the convex hull of }h(\overline{\bB_{4r}(z_0)}),
\]
which is a compact subset of $\Omega$ since the latter is convex. Next we let $\zeta \in C^{\infty}_c(\bB_{2r}(z_0))$ be a cut-off function that equals $1$ on $\overline{\bB_r(z_0)}$. It follows from~\eqref{eq:diff-q-Beltrami} that on $\bB_{2r}(z_0)$, there holds whenever $0 < |t| < r$ that
\begin{equation}\label{eq:diff-q-Beltrami-cutoff}
\begin{split}
(\zeta D_{[t]}h)_{\overline{z}} =\ & \Big( \zeta_{\overline{z}} D_{[t]}h - \zeta D_{[t]}(\mu\circ h) \cdot (\overline{h}_{\overline{z}})_{[t]} + \zeta_{\overline{z}} (\mu \circ h) \cdot D_{[t]}\overline{h} \Big) - (\mu\circ h) \cdot (\zeta D_{[t]}\overline{h})_{\overline{z}}\\
=: &\   \sigma(t) - (\mu\circ h) \cdot (\zeta D_{[t]}\overline{h})_{\overline{z}}.
\end{split}
\end{equation}
We now estimate the $L^p$-norm of each of the three terms comprising $\sigma(t)$. First, again with $0 < |t| < r$, we have
\begin{equation}\label{eq:sigma-t-estimate-1}
\| \zeta_{\overline{z}} D_{[t]}h \|_{p; \CC} \leq C\|\nabla \zeta\|_{\infty} \|D_{[t]}h\|_{p; \bB_{2r}(z_0)} \leq C\|\nabla\zeta\|_\infty \|\nabla h\|_{p; \bB_{3r}(z_0)},
\end{equation}
where the last inequality is a standard fact about difference quotients (see \cite[Section 5.8.2, Theorem 3(i)]{Evans}). Likewise, for the third term in the definition of $\sigma(t)$, 
\begin{equation}\label{eq:sigma-t-estimate-2}
\| \zeta_{\overline{z}} (\mu \circ h) \cdot D_{[t]}\overline{h} \|_{p; \CC} \leq Ck\|\nabla\zeta\|_{\infty} \|\nabla h\|_{p; \bB_{3r}(z_0)}.
\end{equation}
Finally, for the second term in the expression defining $\sigma(t)$ in \eqref{eq:diff-q-Beltrami-cutoff}, we begin by noting that, for $0 < |t| < r$,
\begin{equation}\label{eq:sigma-t-estimate-3-1}
\begin{split}
\| \zeta D_{[t]}(\mu \circ h) \cdot (\overline{h}_{\overline{z}})_{[t]} \|_{p; \CC} \leq\ & C\| \zeta D_{[t]}(\mu \circ h) \|_{\infty; \CC} \|\nabla h\|_{p; \bB_{3r}(z_0)}.
\end{split}
\end{equation}
By the fundamental theorem of calculus, 
\[
\zeta D_{[t]}(\mu \circ h) = \Big(\int_{0}^{1}[D\mu \circ (sh_{[t]} + (1-s)h)]ds \Big)\cdot \zeta D_{[t]}h.
\]
Since $\zeta$ vanishes outside of $\bB_{2r}(z_0)$ and since $sh(z + t) + (1-s)h(z) \in K_{z_0, r}$ for all $z \in \bB_{2r}(z_0)$, $|t| < r$ and $s \in [0, 1]$, we deduce that
\[
\begin{split}
\| \zeta D_{[t]}(\mu \circ h) \|_{\infty; \CC} \leq \ &\|D\mu\|_{\infty; K_{z_0, r}} \|\zeta D_{[t]}h\|_{\infty; \CC}.
\end{split}
\]
Combining this with~\eqref{eq:sigma-t-estimate-3-1},~\eqref{eq:sigma-t-estimate-2} and~\eqref{eq:sigma-t-estimate-1} shows that 
\begin{equation}\label{eq:sigma-t-estimate}
\|\sigma(t)\|_{p; \CC} \leq C\|\nabla\zeta\|_{\infty} \|\nabla h\|_{p; \bB_{3r}(z_0)} + C\|D\mu\|_{\infty; K_{z_0, r}}  \|\zeta D_{[t]}h\|_{\infty; \CC}\|\nabla h\|_{p; \bB_{3r}(z_0)}.
\end{equation}
Since $\mu$ and $h$ are of class $C^1$ on $\Omega$ and $f(\Omega)$, respectively, we conclude that
\begin{equation}\label{eq:sigma-t-estimate-coarse}
\|\sigma(t)\|_{p; \CC} \leq C, \text{ with $C$ independent of $t$}.
\end{equation}
Noting that $w: = \zeta D_{[t]}h - (\zeta D_{[t]}h)(0)$ lies in the space $\cB_{p}$, from Lemma \ref{lemm:inhomogeneous-Beltrami}(a), Lemma \ref{lemm:CZ}(b), and the equation~\eqref{eq:diff-q-Beltrami-cutoff}, we see that
\[
(\zeta D_{[t]}h)_{z} = -T\big((\mu\circ h) \cdot (\zeta D_{[t]}\overline{h})_{\overline{z}}\big) + T(\sigma(t)).
\]
Since $(\zeta D_{[t]}\overline{h})_{\overline{z}}$ is just the complex conjugate of $(\zeta D_{[t]}h)_{z}$, upon taking the $L^p$-norm on both sides in the above equation and recalling that $k\|T\|_p < 1$, we get
\begin{equation}\label{eq:x-diff-quotient-estimate}
\|(\zeta D_{[t]}h)_{z} \|_{p; \CC} \leq C_{k, p} \|\sigma(t)\|_{p; \CC},
\end{equation}
which together with~\eqref{eq:sigma-t-estimate-coarse} implies that $\|(\zeta D_{[t]}h)_{z} \|_{p; \CC}$ is bounded independently of $t$. Going back to equation~\eqref{eq:diff-q-Beltrami-cutoff}, we infer that the same is true of $\|(\zeta D_{[t]}h)_{\overline{z}} \|_{p; \CC}$. Recalling our choice of $\zeta$, we have thus obtained a bound on $\|D_{[t]}\nabla h\|_{p; \bB_{r}(z_0)}$ that is independent of $t$, which implies that the distributional derivative $(\nabla h)_x$ exists and lies in $L^p(\bB_{r}(z_0))$. Repeating this argument with vertical instead of horizontal difference quotients shows that the same holds for $(\nabla h)_{y}$.

Given a disk $\overline{B} \subset f(\Omega)$ and $\zeta \in C^{\infty}_{c}(B)$ as in the statement, now that we know $h \in C^1 \cap W^{2, p}$ locally on $f(\Omega)$, we may differentiate the equation satisfied by $h$ in Lemma~\ref{lemm:conjugate-Beltrami} to get the following analogue of~\eqref{eq:diff-q-Beltrami-cutoff}:
\begin{equation}\label{eq:diff-Beltrami-cutoff}
(\zeta h_{x})_{\overline{z}} = \zeta_{\overline{z}}h_{x} - (D\mu \circ h)(\zeta h_{x} ) \overline{h}_{\overline{z}} + (\mu \circ h) \zeta_{\overline{z}} \cdot \overline{h}_{x} - (\mu \circ h) (\zeta\overline{h}_{x})_{\overline{z}}.
\end{equation}
By the argument leading to~\eqref{eq:x-diff-quotient-estimate}, using, among other things, the fact that $\zeta h_{x} - (\zeta h_{x})(0)$ belongs to $\cB_{p}$ and that $k \|T\|_{p} < 1$, we get
\[
\begin{split}
\|(\zeta h_{x})_{z}\|_{p} \leq\ & C_{k, p}\|\zeta_{\overline{z}}h_{x} - (D\mu \circ h)\cdot (\zeta h_{x}) \cdot \overline{h}_{\overline{z}} + (\mu \circ h)\zeta_{\overline{z}}\cdot \overline{h}_{x}\|_{p}\\
\leq\ & C_{k, p}\big(  \|\nabla\zeta\|_{\infty}\|h_{x}\|_{p; B} +  \|D\mu\|_{\infty; h(\overline{B})} \|\nabla h\|_{p; B} \|\zeta h_{x}\|_{\infty} \big).
\end{split}
\]
Combining this with~\eqref{eq:diff-Beltrami-cutoff} yields an $L^{p}$-estimate on $(\zeta h_{x})_{\overline{z}}$, and hence
\[
\|\zeta h_{x}\|_{1,p} \leq C_{k, p} (\|\zeta\|_{\infty} + \|\nabla \zeta\|_{\infty}) \|h_{x}\|_{p; B} + C_{k, p} \|D\mu\|_{\infty; h(\overline{B})} \|\nabla h\|_{p; B} \|\zeta h_{x}\|_{\infty}.
\]
A similar estimate holds with $h_y$ in place of $h_x$, and we obtain~\eqref{eq:inverse-W2p-bound}.\\
\end{proof}
We can now give the proof of Proposition~\ref{prop:qc-H-dependence}.
\begin{proof}[Proof of Proposition~\ref{prop:qc-H-dependence}]
Define $\widehat{\mu}_n$ and $\widehat{\mu}$ in terms of $\mu_n$ and $\mu$, respectively, by the formula~\eqref{eq:qc-reflection}. By our assumptions on $(\mu_n)$, we have that $\|\widehat{\mu}_n\|_{\infty; \CC} \leq k$ and that  $\widehat{\mu}_n \to \widehat{\mu}$ almost everywhere on $\CC$. Recalling the notation of Theorems \ref{thm:qc-existence} and \ref{thm:qc-existence-H}, in what follows, we write 
\[
f_n = f^{\widehat{\mu}_n},\ \ f = f^{\widehat{\mu}},\quad \ h_n = (f_n)^{-1},\ \ h = f^{-1},
\]
and also set
\[
w_n = w^{\mu_n},\ \ w = w^{\mu}, \quad \ v_n = (w_{n})^{-1},\ \ v = w^{-1}.
\]
Recall that $w_n = f_n|_{\HH}$ and $w = f|_{\HH}$. In particular, for all $z \in \HH$, we have
\[
h_n(w_n(z)) = h_n(f_n(z)) = z = v_n(w_n(z)),
\]
which implies that $h_n = v_n$ on $\HH$. A similar argument shows that $h|_{\HH} = v$.

As with Theorem~\ref{thm:qc-dependence}, we demonstrate the asserted convergence~\eqref{eq:qc-H-dependence} by a contradiction argument. Assume, then, that there exist a disk $\overline{\bB_{2r}(z_0)} \subset \HH$ and some $\alpha > 0$ such that along a subsequence we have
\begin{equation}\label{eq:inverse-dependence-contra}
\|v_n - v\|_{\infty; \bB_{r}(z_0)} + \|\nabla v_n - \nabla v\|_{p; \bB_{r}(z_0)} \geq \alpha \text{ for all }n.
\end{equation}
Since $\|\widehat{\mu}_n\|_{\infty; \CC} \leq k$, arguing as in the proof of Theorem~\ref{thm:qc-dependence}, we may again apply Proposition \ref{prop:qc-estimates} to both $(f_n)$ and $(h_n)$, which gives the uniform bound
\begin{equation}\label{eq:inverses-uniform-bound}
\|h_n\|_{\infty; \bB_R} + \|\nabla h_n\|_{p; \bB_R} \leq C_{k, R} \quad \text{for all }n \in \NN \text{ and }R > 0,
\end{equation}
and also produces a subsequence of $(f_n, h_n)$ that converges uniformly locally on $\CC$ to a pair of homeomorphisms $(\widetilde{f}, \widetilde{h})$ that are inverses of each other. As $\widehat{\mu}_n \to \widehat{\mu}$ almost everywhere on $\CC$, we see from the conclusion of Theorem~\ref{thm:qc-dependence} that $\widetilde{f} = f$, and hence $\widetilde{h} = h$. Since $h_n|_{\HH} = v_n$ and $h|_{\HH} = v$, we deduce that
\begin{equation}\label{eq:inverse-uniform-contradict}
v_n \to v \text{ in }C^0_{\loc}(\HH).
\end{equation}
The next step is to estimate $\|\nabla v_n - \nabla v\|_{p; \bB_{r}(z_0)}$. By Lemma~\ref{lemm:conjugate-Beltrami} and a straightforward computation, we find that
\begin{equation}\label{eq:inverse-Beltrami-difference}
(h_n - h)_{\overline{z}} = -(\widehat{\mu}_n\circ h_n - \widehat{\mu}\circ h)(\overline{h_n})_{\overline{z}} - (\widehat{\mu}\circ h)(\overline{h_n} - \overline{h})_{\overline{z}}.
\end{equation}
Introducing, similar to the proof of Theorem~\ref{thm:qc-dependence}, a cut-off function $\zeta$ that equals $1$ on $\bB_{r}(z_0)$ and vanishes outside $\bB_{2r}(z_0)$, we deduce from the above equation that
\begin{equation}\label{eq:inverse-cutoff-difference}
\begin{split}
[\zeta(v_n - v)]_{\overline{z}} = - (\mu\circ v)[\zeta(\overline{v_n} - \overline{v})]_{\overline{z}} + \tau_n,
\end{split}
\end{equation}
where
\[
\tau_n =  \zeta_{\overline{z}} (v_n - v)  - \zeta (\mu_n \circ v_n - \mu\circ v)(\overline{v_n})_{\overline{z}} + \zeta_{\overline{z}} (\mu\circ v)(\overline{v_n} - \overline{v}).
\]
Since our assumptions imply in particular that $\mu_n \to \mu$ in $C^0_{\loc}(\HH)$, we deduce from~\eqref{eq:inverse-uniform-contradict} that $\mu_n\circ v_n \to \mu\circ v$ uniformly locally on $\HH$, which together with~\eqref{eq:inverses-uniform-bound} gives
\[
\|\zeta(\mu_n \circ v_n - \mu\circ v)(\overline{v_n})_{\overline{z}}\|_{p; \CC} \to 0.
\]
Using again~\eqref{eq:inverse-uniform-contradict}, we obtain
\begin{equation}\label{eq:quasi-linear-Beltrami-error}
\lim_{n \to \infty}\|\tau_n\|_{p; \CC} = 0.
\end{equation}
Now, since the function $\zeta (v_n - v)$, extended to be $0$ outside $\bB_{2r}(z_0)$, belongs to $\cB_{p}$, by Lemma~\ref{lemm:inhomogeneous-Beltrami}(a), equation~\eqref{eq:inverse-cutoff-difference}, and Lemma~\ref{lemm:CZ}(b), we have
\[
[\zeta(v_n - v)]_{z} = -T((\mu\circ v)[\zeta(\overline{v_n} - \overline{v})]_{\overline{z}}) + T(\tau_n).
\]
Taking the $L^p$-norm gives
\begin{equation}\label{eq:inverse-cutoff-estimate}
\begin{split}
\| [\zeta(v_n - v)]_{z} \|_{p; \CC} \leq\ & \|T\|_p \cdot k \cdot \| [\zeta(v_n - v)]_{z} \|_{p; \CC} + \|T\|_p \|\tau_n\|_{p; \CC}.
\end{split}
\end{equation}
Since $k \|T\|_p < 1$, we see from \eqref{eq:inverse-cutoff-estimate} and \eqref{eq:quasi-linear-Beltrami-error} that
\[
\lim_{n \to \infty}\| [\zeta(v_n - v)]_{z} \|_{p; \CC} = 0,
\]
which together with~\eqref{eq:inverse-cutoff-difference} and~\eqref{eq:quasi-linear-Beltrami-error} shows that 
\[
\lim_{n \to \infty}\| [\zeta(v_n - v)]_{\overline{z}} \|_{p; \CC} = 0.
\]
By our choice of $\zeta$, these previous two convergences and~\eqref{eq:inverse-uniform-contradict} yields a contradiction to~\eqref{eq:inverse-dependence-contra}. This proves (a).

For part (b), note first that $\nabla v_n$ and $\nabla v$ are indeed $C^0 \cap W^{1, p}$-maps on $\HH$, thanks to Lemma~\ref{lemm:qc-standard-properties}(e), Lemma~\ref{lemm:inverse-W2p}, and the assumption that $\mu$ and $\mu_n$ are $C^1$ on $\HH$. Also, since $W^{2,p}$ embeds compactly into $C^{1}$ on compact sets in dimension $2$, and since we have already shown in part (a) that $v_n \to v$ in $C^0 \cap W^{1, p}$ locally on $\HH$, to prove part (b) it suffices to establish the second conclusion, namely \eqref{eq:uniform-W2p-on-compact}, which we now turn to. Given any $z_0 = (x_0, y_0) \in \HH$ and any $r_0 < \frac{y_0}{8}$, we take a cut-off function $\zeta$ that equals $1$ on $\bB_{r_0}(z_0)$ and vanishes outside $\bB_{2r_0}(z_0)$, and apply the estimate~\eqref{eq:inverse-W2p-bound} to get
\begin{equation}\label{eq:inverse-W2p-applied}
\| \zeta \nabla v_n \|_{1, p; \CC}\leq C_{k, p, r_0} \|\nabla v_n\|_{p; \bB_{\frac{y_0}{2}}(z_0)} + C_{k, p}\|D\mu_n\|_{\infty; K_n} \|\nabla v_n\|_{p; \bB_{\frac{y_0}{2}}(z_0)} \|\zeta \nabla v_n\|_{\infty; \CC},
\end{equation}
where $K_n = v_n(\overline{\bB_{\frac{y_0}{2}}(z_0)})$. Since, as noted earlier, $v_n \to v$ in $C^0 \cap W^{1, p}$ locally on $\HH$, we get some constant $A > 0$ and some fixed compact set $K' \subset \HH$ such that for all $n$ we have
\[
\|\nabla v_n\|_{p; \bB_{\frac{y_0}{2}}(z_0)} \leq A, \quad K_n \subset K'.
\]
Since $\mu_n \to \mu$ in $C^1_{\loc}(\HH)$, we further get some $B > 0$ such that
\[
\|D\mu_n\|_{\infty; K'} \leq B, \text{ for all }n.
\]
Putting these bounds back into \eqref{eq:inverse-W2p-applied} and also estimating the very last term by
\[
\|\zeta \nabla v_n\|_{\infty; \bB_{2r_0}(z_0)} \leq (2r_0)^{1 - \frac{2}{p}}\cdot [\zeta \nabla v_n]_{1 - \frac{2}{p}; \bB_{2r_0}(z_0)} \leq C_{p}\cdot r_0^{1 - \frac{2}{p}} \|\nabla (\zeta \nabla v_n)\|_{p; \CC},
\]
we arrive at
\[
\|\zeta \nabla v_n\|_{1, p; \CC} \leq C_{k,p, r_0} A + C_{k, p}AB r_0^{1 - \frac{2}{p}}\|\zeta \nabla v_n\|_{1, p; \CC}.
\]
Choosing $r_0$ to be sufficiently small depending only on $k, p, A, B$, it follows that
\begin{equation}\label{eq:uniform-W2p-on-small-balls}
\|\nabla v_n\|_{1, p; \bB_{r_0}(z_0)} \leq C(k, p, A, B).
\end{equation}
That is, each $z_0 \in \HH$ possesses a neighborhood on which the $W^{1, p}$-norms of $\nabla v_n$ are bounded independently of $n$. A standard covering argument then gives \eqref{eq:uniform-W2p-on-compact}. As remarked earlier, this implies that $v_n \to v$ in $C^{1}_{\loc}(\HH)$, and we are done with part (b).

For part (c), we only sketch the argument to avoid repetition. Applying to $f_n$ a difference quotient argument on compact subsets of $\HH$ in a way that is similar to the proof of Lemma~\ref{lemm:inverse-W2p}, we get that $w_n \in W^{2, p}_{\loc}(\HH)$, along with the following estimate analogous to \eqref{eq:inverse-W2p-bound}: for any disk $\overline{B} \subset \HH$ and $\zeta \in C^{\infty}_{c}(B)$, there holds
\[
\|\zeta \nabla w_n\|_{1, p} \leq C_{k, p}(\|\zeta\|_{\infty} + \|\nabla \zeta\|_{\infty}) \|\nabla w_n\|_{p; B} + C\|D \mu_n\|_{\infty; B} \|\zeta \nabla w_n\|_{p; B}.
\]
Combining this estimate and the $C^0 \cap W^{1, p}$-convergence of $w_n$ to $w$ on compact subsets of $\HH$, the latter a consequence of Theorem~\ref{thm:qc-dependence}, we conclude that $w_n$ converges to $w$ in $C^{1}_{\loc}(\HH)$, as asserted.
\end{proof}
\subsection{Proofs of Propositions~\ref{prop:qc-smooth-dependence} and~\ref{prop:metric-inverse-convergence}}\label{subsec:inverse-beltrami-proofs}
\begin{proof}[Proof of Proposition~\ref{prop:qc-smooth-dependence}]
As in the proof of Proposition \ref{prop:qc-H-dependence}, we reflect $\mu_n$ and $\mu$ across $\partial \HH$ as indicated in~\eqref{eq:qc-reflection} to obtain $\widehat{\mu}_n$ and $\widehat{\mu}$. Then again $w_n = f^{\widehat{\mu}_n}|_{\HH}$ and $w = f^{\widehat{\mu}}|_{\HH}$. Since $\widehat{\mu}_n$ converges to $\widehat{\mu}$ almost everywhere on $\CC$, and since $\|\widehat{\mu}_{n}\|_{\infty} \leq k$ for all $n$, we may invoke Theorem~\ref{thm:qc-dependence} to infer that 
\begin{equation}\label{eq:qc-harmonic-C0-convergence}
w_n \to w \text{ in } C^0 \cap W^{1, p} \text{ on compact subsets of }\HH,
\end{equation}
where $p = p(k)$ is as in Theorem \ref{thm:qc-existence} and in particular satisfies $k \|T\|_{p} < 1$. To obtain higher order estimates on $w_n$, take any disk such that $\overline{\bB_{2r}(z_0)} \subset \HH$ and let $\zeta \in C^\infty_c(\bB_{2r}(z_0))$ be a cut-off function that is identically $1$ on $\bB_{r}(z_0)$. Since $(w_n)_{\overline{z}} = \mu_n (w_n)_{z}$, we have
\[
\begin{split}
(\zeta w_n)_{\overline{z}} =\ & \zeta_{\overline{z}} w_n + \zeta \mu_n (w_n)_{z} = \mu_n(\zeta w_n)_{z} + (\zeta_{\overline{z}} - \mu_n\zeta_{z})w_n.
\end{split}
\]
Taking also an integer $s \geq 2$ and applying to both sides an arbitrary differential operator of the form $D^\alpha : = (\paop{x})^{\alpha_1} (\paop{y})^{\alpha_2}$ with $\alpha_1 + \alpha_2 = s - 1$, we obtain
\[
\begin{split}
[D^\alpha(\zeta w_n)]_{\overline{z}} =\ & \mu_n (D^\alpha (\zeta w_n))_{z} + [D^\alpha, \mu_n](\zeta w_n)_{z} + D^\alpha((\zeta_{\overline{z}} - \mu_n\zeta_{z})w_n),
\end{split}
\]
where $[D^{\alpha}, \mu_n]$ denotes the commutator of $D^{\alpha}$ and multiplication by $\mu_n$. Fixing $n$ and introducing the shorthands
\[
\varphi = D^\alpha(\zeta w_n) - D^\alpha(\zeta w_n)(0),\quad F = [D^\alpha, \mu_n](\zeta w_n)_{z} + D^\alpha((\zeta_{\overline{z}} - \mu_n\zeta_{z})w_n),
\]
then $\varphi \in \cB_p$ and $F \in L^p(\CC)$, and furthermore $\varphi$ satisfies
\[
\varphi_{\overline{z}} = \mu_n \varphi_{z} + F.
\]
From this together with Lemma \ref{lemm:inhomogeneous-Beltrami}(a), Lemma \ref{lemm:CZ}(b) and the fact that $k\|T\|_{p} < 1$, we get
\begin{equation}\label{eq:estimate-with-w-F}
\|[D^\alpha(\zeta w_n)]_{z}\|_{p; \CC} + \|  [D^\alpha(\zeta w_n)]_{\overline{z}}\|_{p; \CC}  =  \|\varphi_z\|_{p; \CC} + \|\varphi_{\overline{z}}\|_{p; \CC} \leq C_{p, k} \|F\|_{p; \CC}.
\end{equation}
Since $F$ involves only derivatives of $w_n$ up to order $s-1$, while $\mu_n$ and all its derivatives are converging uniformly on $\bB_{2r}(z_0)$, we deduce from the estimate~\eqref{eq:estimate-with-w-F} that
\[
\|w_n\|_{s, p; \bB_{r}(z_0)} \leq C\|w_n\|_{s - 1, p; \bB_{2r}(z_0)}, \text{ where }C \text{ is independent of }n.
\]
Thus, starting from~\eqref{eq:qc-harmonic-C0-convergence} and arguing inductively, we see that 
\[
\sup_{n} \|w_n\|_{s, p; K} < \infty, \quad \text{ for all }s \geq 2 \text{ and compact set }K \subset \HH,
\]
which suffices to upgrade~\eqref{eq:qc-harmonic-C0-convergence} to smooth convergence on compact subsets of $\HH$.\\
\end{proof}

As preparation for the proof of Proposition~\ref{prop:metric-inverse-convergence}, we note the following facts:
\vskip 1mm
\begin{enumerate}
\item[(i)] In the notation of Proposition~\ref{prop:qc-smooth-dependence}, the maps $v_n := w_n^{-1}$ are smooth, since $w_n$ are assumed to be diffeomorphisms. Moreover, Lemma~\ref{lemm:conjugate-Beltrami} gives
\begin{equation}\label{eq:metric-inverse-Beltrami}
(v_n)_{\overline{z}} = -(\mu_n \circ v_n)(\overline{v_n})_{\overline{z}}.
\end{equation}
\vskip 1mm
\item[(ii)] Given a differential operator $D^\alpha = (\paop{x})^{\alpha_1} (\paop{y})^{\alpha_2}$, with $m: = \alpha_1 + \alpha_2 \geq 2$, we have, as a consequence of the chain rule, the follow pointwise bound:
\begin{equation}\label{eq:multi-chain-rule}
\begin{split}
&\Big|D^\alpha(\mu_n \circ v_n) - [(D\mu_n) \circ v_n] D^\alpha v_n \Big|\\
 \leq\ &   C_{m}\sum_{j = 2}^{m} \big|(D^j\mu_n)\circ v_n\big| \Big( \sum_{i_1 + \cdots + i_j = m} \big|D^{i_1}v_n\big| \cdots \big| D^{i_j}v_n \big| \Big),
\end{split}
\end{equation}
where $i_1, \cdots, i_j$ are positive integers. Since $j \geq 2$, there holds the additional constraint that $i_1, \cdots, i_j \leq m-1$.
\end{enumerate}
\begin{proof}[Proof of Proposition~\ref{prop:metric-inverse-convergence}]
As in the previous proof we let $p = p(k) > 2$ be the exponent from Theorem \ref{thm:qc-existence}, so that $k\|T\|_{p} < 1$. Under our current hypotheses, we get from Proposition~\ref{prop:qc-H-dependence}(b) that
\begin{equation}\label{eq:metric-inverse-uniform}
v_n : = w_n^{-1} \longrightarrow w^{-1} =: v \text{ in }C^1  \text{ locally on }\HH,
\end{equation}
and that $(\nabla^2 v_n)$ is a bounded sequence in $L^p$ on each compact subset of $\HH$. Consequently, given a disk $ \bB_{4r_0}(z_0) \subset \HH$, with $r_0 < 1$, there exist a constant $A$ and a fixed compact set $K$ such that, for all $n$,
\begin{equation}\label{eq:prelim-bounds}
\|v_n\|_{1, \infty; \bB_{2r_0}(z_0)} + \|D^2v_n\|_{p;  \bB_{2r_0}(z_0)} \leq A,\quad \text{and} \quad v_n(\bB_{2r_0}(z_0)) \subset K.
\end{equation}
The smooth convergence of $\mu_n$ in turn yields constants $B_0 \leq B_1 \leq  \cdots$ such that for all $m \geq 0$ there holds
\begin{equation}\label{eq:mu-uniform-bounds}
\|\mu_n\|_{m, \infty; K} \leq B_m, \text{ for all }n.
\end{equation}
Below, to save space, we drop the subscript $n$ from $v_n$ and $\mu_n$. Taking a cut-off function $\zeta \in C^{\infty}_{c}(\bB_{2r_0}(z_0))$ that equals $1$ on $\bB_{r_0}(z_0)$, we compute using~\eqref{eq:metric-inverse-Beltrami} that
\[
(\zeta v)_{\overline{z}} = -(\mu \circ v) (\zeta \overline{v})_{\overline{z}} + \zeta_{\overline{z}}v + (\mu \circ v)\zeta_{\overline{z}}\overline{v} .
\]
Fixing any $s \geq 3$ and applying $D^\alpha : = (\paop{x})^{\alpha_1} (\paop{y})^{\alpha_2}$ to both sides as in the previous proof, where $\alpha_1 + \alpha_2 = s - 1$, we get
\[
\begin{split}
[D^\alpha(\zeta v)]_{\overline{z}} =\ & - (\mu \circ v) [D^\alpha (\zeta \overline{v})]_{\overline{z}} - [D^\alpha, \mu \circ v](\zeta \overline{v})_{\overline{z}} + D^{\alpha}\Big( \zeta_{\overline{z}}v + (\mu \circ v)\zeta_{\overline{z}}\overline{v} \Big)\\
=\ &  - (\mu \circ v) [D^\alpha (\zeta \overline{v})]_{\overline{z}}  -\Big( (\zeta\overline{v})_{\overline{z}} D^\alpha(\mu \circ v) - \zeta_{\overline{z}}\overline{v} D^\alpha (\mu \circ v)\Big) \\
& - \Big( [D^\alpha, \mu \circ v] (\zeta \overline{v})_{\overline{z}} - (\zeta \overline{v})_{\overline{z}} D^\alpha(\mu \circ v)  \Big) + [D^\alpha, \zeta_{\overline{z}}\overline{v}](\mu \circ v) + D^\alpha(\zeta_{\overline{z}}v).
\end{split}
\]
Thus, letting $\varphi = D^\alpha(\zeta v) - D^\alpha(\zeta v)(0)$ and also introducing
\begin{equation}\label{eq:F1-F2}
\begin{split}
F_1 = \ & -\big( (\zeta\overline{v})_{\overline{z}} D^\alpha(\mu \circ v) - \zeta_{\overline{z}}\overline{v} D^\alpha (\mu \circ v)\big) \\
=\ & -\big(\zeta \overline{v}_{\overline{z}}\big)D^\alpha(\mu \circ v), \\
F_2 =\ & D^\alpha(\zeta_{\overline{z}}v) - \big(  \underbrace{[D^\alpha, \mu \circ v] (\zeta \overline{v})_{\overline{z}} - (\zeta \overline{v})_{\overline{z}} D^\alpha(\mu \circ v)}_{F_{21}}  \big)  + \underbrace{[D^\alpha, \zeta_{\overline{z}}\overline{v}](\mu \circ v)}_{F_{22}},
\end{split}
\end{equation}
we see that $\varphi \in \cB_{p}$, and that
\begin{equation}\label{eq:inverse-derivative-Beltrami}
\varphi_{\overline{z}} = -(\mu \circ v) \overline{\varphi}_{\overline{z}} + F_1 + F_2.
\end{equation}
Applying Lemma \ref{lemm:inhomogeneous-Beltrami}(a) and Lemma \ref{lemm:CZ}(a) as we did in the proof of Proposition~\ref{prop:qc-H-dependence} right after equation~\eqref{eq:inverse-cutoff-difference}, we get
\[
\varphi_{z} = -T((\mu \circ v)\overline{\varphi}_{\overline{z}} + F_1 + F_2).
\]
Taking the $L^p$-norm, recalling that $\|\mu \circ v\|_{\infty; \HH}\cdot \|T\|_{p} \leq k \|T\|_{p} < 1$, and using \eqref{eq:inverse-derivative-Beltrami} again, we have
\begin{equation}\label{eq:inverse-prelim-Lp}
\|[D^\alpha(\zeta v)]_{z}\|_{p; \CC} + \|[D^\alpha(\zeta v)]_{\overline{z}}\|_{p; \CC} \leq C_{k, p}\big( \|F_1\|_{p; \CC} + \|F_2\|_{p; \CC} \big).
\end{equation}
We next estimate the $L^p$-norm of $F_1$ and $F_2$. Using the second expression for $F_1$ in \eqref{eq:F1-F2}, rewriting it as 
\[
\begin{split}
F_1 =\ &  -  (\zeta \overline{v}_{\overline{z}}) [(D\mu) \circ v] D^\alpha v  - (\zeta \overline{v}_{\overline{z}}) \Big( D^\alpha (\mu \circ v) - [(D\mu) \circ v] D^\alpha v \Big),
\end{split}
\]
and then recalling~\eqref{eq:multi-chain-rule}, we have
\begin{equation}\label{eq:inverse-F1-estimate}
\begin{split}
\|F_1\|_{p; \CC} \leq \ & \|Dv\|_{\infty; \bB_{2r_0}(z_0)} \cdot \big\| |(D\mu) \circ v||D^{s-1}v| \big\|_{p; \bB_{2r_0}(z_0)} \\
& + C_{s}  \|v\|_{1, p; \bB_{2r_0}(z_0)}  \cdot \sum_{j = 2}^{s-1} \sum_{i_1 + \cdots + i_j = s-1} \big\||(D^j \mu) \circ v||D^{i_1}v| \cdots |D^{i_j}v|  \big\|_{\infty; \bB_{2r_0}(z_0)}.
\end{split}
\end{equation}
As noted after \eqref{eq:multi-chain-rule}, each $i_\lambda$ in the summation on the second line of \eqref{eq:inverse-F1-estimate} is at most $s - 2$, so we have by Sobolev embedding that
\[
\|D^{i_\lambda}v\|_{\infty; \bB_{2r_0}(z_0)} \leq C_{p, r_0}\|D^{i_\lambda}v\|_{1, p; \bB_{2r_0}(z_0)} \leq C_{p, r_0}\|v\|_{s-1, p; \bB_{2r_0}(z_0)}.
\]
Substituting this back into~\eqref{eq:inverse-F1-estimate}, and combining the result with~\eqref{eq:prelim-bounds} and \eqref{eq:mu-uniform-bounds}, we obtain
\begin{equation}\label{eq:F1-estimate}
\begin{split}
\|F_1\|_{p; \CC} \leq\ & A \cdot B_1 \|v\|_{s- 1, p; \bB_{2r_0}(z_0)} + C(s, p, r_0, B_{s-1}, \|v\|_{s-1, p; \bB_{2r_0}(z_0)}).
\end{split}
\end{equation}
Turning to $F_2$, we note that
\[
\begin{split}
|F_{21}| =\ & \big|D^{\alpha}\big( (\mu \circ v) (\zeta \overline{v})_{\overline{z}} \big) - (\mu \circ v) D^{\alpha}(\zeta \overline{v})_{\overline{z}} - (\zeta \overline{v})_{\overline{z}} D^\alpha(\mu \circ v)\big| \\
\leq\ & C_{s}\sum_{j = 1}^{s-2} |D^{j}(\mu \circ v)|\cdot |D^{s-1-j}(\zeta\overline{v})_{\overline{z}}|,
\end{split}
\]
and that a similar estimate holds for $F_{22}$. Consequently, the following holds pointwise on $\bB_{2r}(z_0)$:
\[
\begin{split}
|F_2| \leq\ & |D^{s-1}(\zeta_{\overline{z}}v)| + C_{s}\sum_{j = 1}^{s-2} \|D^j(\mu \circ v)\|_{\infty; \bB_{2r_0}(z_0)}\cdot | D^{s - 1 - j}(\zeta \overline{v})_{\overline{z}}|\\
&+  C_{s}\sum_{j = 0}^{s - 2} \|D^j (\mu \circ v)\|_{\infty; \bB_{2r_0}(z_0)}\cdot |D^{s - 1 - j}(\zeta_{\overline{z}}\overline{v})|.
\end{split}
\]
Imitating what we did with $F_{1}$, we use~\eqref{eq:multi-chain-rule} and the Sobolev embedding $W^{1, p} \hookrightarrow C^0$ to estimate the derivatives of $\mu \circ v$, and bring in the bounds \eqref{eq:prelim-bounds} and \eqref{eq:mu-uniform-bounds}. The resulting estimate is
\[
\|F_2\|_{p; \CC} \leq C(s, p, r_0, B_{s-2}, \|v\|_{s-1, p; \bB_{2r_0}(z_0)}).
\]
Combining this with~\eqref{eq:F1-estimate}, and then recalling~\eqref{eq:inverse-prelim-Lp}, we obtain
\begin{equation}\label{eq:F1-F2-estimate}
\|v_n\|_{s,p;\bB_{r_0}(z_0)} \leq  C(p, k, s, r_0, A, B_{s-1}, \|v_n\|_{s-1, p; \bB_{2r_0}(z_0)}).
\end{equation}
We can now use~\eqref{eq:F1-F2-estimate} inductively, starting with $s = 3$ and the local $W^{2, p}$-boundedness of $(v_n)$, to get that for all $s \geq 3$ and compact subset $K \subset \HH$, we have $\sup_{n}\|v_n\|_{s, p; K} < \infty$. Recalling~\eqref{eq:metric-inverse-uniform}, we conclude that $v_n$ converges to $v$ in $C^{\infty}_{\loc}(\HH)$. 
\end{proof}
\section{Basic estimate for iterated harmonic replacements}\label{appendix:estimates-iterated-replacement}
In this appendix we prove Proposition \ref{prop:iterated-replacement-estimates} and Corollary~\ref{coro:iterated-iterated-replacement}.
\begin{proof}[Proof of Proposition \ref{prop:iterated-replacement-estimates}]
As in the proof of Proposition \ref{prop:iterated-replacement}, we let $(\tau, u) = \Psi(\sigma, v)$ and use freely the notation from the diagram \eqref{eq:Phi-diagram}. To start, by~\eqref{eq:iterated-replacement-smallness} and Lemma~\ref{lemm:energy-after-replacement} we have 
\begin{equation}\label{eq:second-replacement-possible}
\sum_{B \in \fB_2} \int_{B}|\nabla v_1|^2_{\sigma} \vol_{\sigma} < \frac{2\ep_0}{3},
\end{equation}
so that the harmonic replacements involved in~\eqref{eq:iterated-replacement-estimate-1} and~\eqref{eq:iterated-replacement-estimate-2} make sense. To prove~\eqref{eq:iterated-replacement-estimate-1}, we label the members of $\fB_{1}$ and $\fB_{2}$ by writing $\fB_{1} = \{B_{1, \alpha}\}_{\alpha \in A}$ and $\fB_{2} = \{B_{2, i}\}_{i \in I}$, and then define the following subsets of $I$:
\[
I_{+, \alpha} = \{i \in I \ |\ \frac{1}{2}B_{2, i} \subset \Int B_{1, \alpha}\}, 
\]
\[
I_{-} = \{i \in I\ |\ \frac{1}{2}B_{2, i}\not\subset \Int B_{1, \alpha} \text{ for any }\alpha \in A\}.
\]
As $\{B_{1, \alpha}\}_{\alpha\in A}$ is a disjoint collection, we obtain a partition of $I$ by writing
\[
I = \big(\cup_{\alpha \in A}I_{+,\alpha}\big) \cup I_{-}.
\]
Given $\alpha \in A$, notice that $v_1\big|_{B_{1, \alpha}}$ minimizes $E(\sigma, \cdot)$ among maps from $B_{1, \alpha}$ into $M$ that agree with $v$ on the boundary. Thus, using as a competitor the map which agrees with $v$ on $B_{1, \alpha} \setminus \big(\cup_{i \in I_{+, \alpha}}\frac{1}{2}B_{2, i}\big)$, and with $\cR(\sigma, v, \frac{1}{2}\fB_{2})|_{\frac{1}{2}B_{2, i}}$ on $\frac{1}{2}B_{2, i}$ for each $i \in I_{+, \alpha}$, we infer that
\[
\int_{B_{1, \alpha}} |\nabla v|_{\sigma}^2 - |\nabla v_1|_{\sigma}^2 \vol_{\sigma} \geq \sum_{i \in I_{+, \alpha}} \int_{\frac{1}{2}B_{2, i}} |\nabla v|_{\sigma}^2 - |\nabla \cR(\sigma, v, \frac{1}{2}\fB_{2})|^2_{\sigma} \vol_{\sigma}.
\]
Summing over $\alpha \in A$ gives
\begin{equation}\label{eq:iterated-estimate-1-+}
E(\sigma, v) - E(\sigma, v_1) \geq \frac{1}{2}\sum_{\alpha \in A}\sum_{i \in I_{+, \alpha}}\int_{\frac{1}{2}B_{2, i}} |\nabla v|_{\sigma}^2 - |\nabla \cR(\sigma, v, \frac{1}{2}\fB_{2})|^2_{\sigma} \vol_{\sigma}.
\end{equation}
In view of the form of the desired estimate \eqref{eq:iterated-replacement-estimate-1}, we note that the left-hand side of \eqref{eq:iterated-estimate-1-+} can be estimate from above by
\begin{equation}\label{eq:iterated-estimate-1-+-lhs}
\begin{split}
E(\sigma, v) - E(\sigma, v_1) =\ & \Big(\frac{1}{2} \sum_{\alpha \in A}\int_{B_{1, \alpha}} |\nabla v|_{\sigma}^2 - |\nabla v_1|_{\sigma}^2 \vol_{\sigma} \Big)^{\frac{1}{2}} \cdot \big[ E(\sigma, v) - E(\sigma, v_1)\big]^{\frac{1}{2}}\\ \leq \ &\ep_0^{\frac{1}{2}} \cdot \big[ E(\sigma, v) - E(\sigma, v_1)\big]^{\frac{1}{2}}.
\end{split}
\end{equation}
Turning to $I_{-}$, we claim that 
\begin{equation}\label{eq:boundary-non-inclusion}
\partial(\mu B_{2, i}) \not\subset \cup_{\alpha \in A}\Int B_{1, \alpha},\quad \text{for all }i \in I_{-} \text{ and } \mu \in [\frac{1}{2}, 1].
\end{equation}
To see this, suppose by contradiction that inclusion holds for some choice of $i, \mu$ as above, and let $\widetilde{B}_{2, i}$ be a lift of $B_{2, i}$ with respect to $p$. Then
\[
\partial (\mu \widetilde{B}_{2, i}) \subset p^{-1}\big( \partial(\mu B_{2, i}) \big) \subset \cup_{\alpha \in A}\big(p^{-1}(\Int B_{1, \alpha})\big).
\]
Note that $\{p^{-1}(\Int B_{1, \alpha})\}_{\alpha \in A}$ is a mutually disjoint collection, and each member is further a disjoint union of open disks, namely the interiors of the lifts of $B_{1, \alpha}$. Since $\partial (\mu \widetilde{B}_{2, i})$ is connected, there exist $\alpha \in A$ and a lift $\widetilde{B}_{1, \alpha}$ of $B_{1, \alpha}$ with respect to $p$ such that 
\[
\partial (\mu \widetilde{B}_{2, i})  \subset \Int \widetilde{B}_{1, \alpha}.
\]
Consequently $\Int \widetilde{B}_{1, \alpha}$ contains the entire disk $\mu \widetilde{B}_{2, i}$, and hence $\frac{1}{2}B_{2,i} \subset \mu B_{2, i} \subset \Int B_{1, \alpha}$, which contradicts the fact that $i \in I_{-}$. Therefore~\eqref{eq:boundary-non-inclusion} holds, from which we deduce that 
\begin{equation}\label{eq:boundary-agreement}
v = v_1 \text{ somewhere on }\partial(\mu B_{2, i}), \quad \text{for all } i \in I_{-} \text{ and }\mu \in [\frac{1}{2}, 1].
\end{equation}
Thanks to~\eqref{eq:iterated-replacement-smallness},~\eqref{eq:second-replacement-possible}, and~\eqref{eq:boundary-agreement}, we may apply Proposition \ref{prop:iterated-replacement} with $\fB$ being the collection
\[
\fB_{2, -} : = \{B_{2, i}\ |\ i \in I_{-}\},
\]
$v_1$ being the current one, and with $v_2$ being just $v$. Doing so, and noting that 
\[
\cR(\sigma, v, \frac{1}{2}\fB_{2, -})|_{\frac{1}{2}B_{2, i}} = \cR(\sigma, v, \frac{1}{2}\fB_{2})|_{\frac{1}{2}B_{2, i}},\quad \text{for all }i \in I_{-},
\]
and that
\[
E(\sigma, v_1) - E(\sigma, \cR(\sigma, v_1, \fB_{2, -})) \leq E(\sigma, v_1) - E(\sigma, \cR(\sigma, v_1, \fB_{2})),
\]
we arrive at
\begin{equation}\label{eq:iterated-estimate-1--}
\begin{split}
&\frac{1}{2}\sum_{i \in I_{-}} \int_{\frac{1}{2}B_{2, i}} |\nabla v|_{\sigma}^2 - |\nabla \cR(\sigma, v, \frac{1}{2}\fB_{2})|_{\sigma}^2 \vol_{\sigma}\\
\leq\ & E(\sigma, v_1) - E(\sigma, \cR(\sigma, v_1, \fB_{2})) + \frac{1}{\kappa}\Big( \sum_{i \in I_{-}} \int_{B_{2, i}} |\nabla v - \nabla v_1|_{\sigma}^2 \vol_{\sigma} \Big)^{\frac{1}{2}}.
\end{split}
\end{equation}
To estimate the last term, we use Remark \ref{rmk:v-replacement}(2) to see that
\begin{equation}\label{eq:energy-difference}
\begin{split}
\sum_{i \in I_{-}}\int_{B_{2, i}} |\nabla v - \nabla v_1|_{\sigma}^2 \vol_{\sigma}  \leq \int_{S}|\nabla v - \nabla v_1|_{\sigma}^2  \vol_{\sigma} \leq 4(E(\sigma, v) - E(\sigma, v_1)).
\end{split}
\end{equation}
Putting this back into~\eqref{eq:iterated-estimate-1--}, summing the result with~\eqref{eq:iterated-estimate-1-+}, and recalling \eqref{eq:iterated-estimate-1-+-lhs}, we get \eqref{eq:iterated-replacement-estimate-1}.

The proof of~\eqref{eq:iterated-replacement-estimate-2} is similar, and hence we only give an outline. For any $\alpha \in A$ and any $i \in I_{+, \alpha}$, since $v_1|_{B_{1, \alpha}}$ is already a harmonic map with respect to $\sigma$, with $E(\sigma, \cdot)$-energy at most $\ep_0$, we have that 
\[
v_1(x) = \cR(\sigma, v_1, \frac{1}{2}\fB_2)(x)\quad  \text{ for all } x \in \frac{1}{2}B_{2, i}.
\]
Therefore, in place of~\eqref{eq:iterated-estimate-1-+}, we simply have
\begin{equation}\label{eq:iterated-estimate-2-+}
\sum_{\alpha \in A}\sum_{i \in I_{+, \alpha}}\int_{\frac{1}{2}B_{2, i}} |\nabla v_1|_{\sigma}^2 - |\nabla \cR(\sigma, v_1, \frac{1}{2}\fB_{2})|^2_{\sigma} \vol_{\sigma} = 0.
\end{equation}
To obtain the counterpart of~\eqref{eq:iterated-estimate-1--} we apply Proposition \ref{prop:iterated-replacement} again with $\fB=\fB_{2, -}$, but instead with $v_1$ there chosen to be $v$, and with $v_2$ being the current $v_1$. Combining the resulting inequality with~\eqref{eq:energy-difference} gives
\[
\begin{split}
&\frac{1}{2}\sum_{i \in I_{-}} \int_{\frac{1}{2}B_{2, i}} |\nabla v_1|_{\sigma}^2 - |\nabla \cR(\sigma, v_1, \frac{1}{2}\fB_2)|_{\sigma}^2 \vol_{\sigma}\\
\leq\ & E(\sigma, v) - E(\sigma, \cR(\sigma, v, \fB_{2})) + \frac{C}{\kappa}\big[ 
E(\sigma, v) - E(\sigma, v_1)\big]^{\frac{1}{2}}.
\end{split}
\]
Summing this with~\eqref{eq:iterated-estimate-2-+} gives~\eqref{eq:iterated-replacement-estimate-2}. The proof is complete.
\end{proof}
Next we prove Corollary~\ref{coro:iterated-iterated-replacement}.
\begin{proof}[Proof of Corollary~\ref{coro:iterated-iterated-replacement}]
Again, the maps appearing in the asserted estimates all make sense thanks to Remark~\ref{rmk:many-replacements}. For convenience, we define 
\[
v_0 = v,\quad v_{j} = \cR(\sigma, v, \fB_{1}, \cdots, \fB_{j}),\quad \text{for }j = 1, \cdots, L-1.
\]
To prove~\eqref{eq:many-iteration-1}, take any $j \in \{1, \cdots, L-1\}$. By~\eqref{eq:iterated-iterated-smallness} and Lemma \ref{lemm:energy-after-replacement}, we see that the hypothesis \eqref{eq:iterated-replacement-smallness} of Proposition \ref{prop:iterated-replacement-estimates} is fulfilled when $v, \fB_1$, and $\fB_2$ are taken to be $v_{L - j - 1}, \fB_{L - j}$, and $2^{-j + 1}\fB_{L}$, respectively. With this choice, we obtain from the estimate \eqref{eq:iterated-replacement-estimate-1} that
\[
\begin{split}
&E(\sigma, v_{L - j}) - E(\sigma, \cR(\sigma, v_{L  - j}, 2^{-j + 1}\fB_{L}))\\
\geq\ & E(\sigma, v_{L - j - 1}) - E(\sigma, \cR(\sigma, v_{L - j - 1}, 2^{-j}\fB_L)) - \frac{1}{\kappa}\big[ E(\sigma, v_{L - j - 1}) - E(\sigma, v_{L - j}) \big]^{\frac{1}{2}}\\
\geq\ & E(\sigma, v_{L - j - 1}) - E(\sigma, \cR(\sigma, v_{L - j - 1}, 2^{-j}\fB_L)) - \frac{1}{\kappa}\big[ E(\sigma, v) - E(\sigma, v_{L - 1}) \big]^{\frac{1}{2}}.
\end{split}
\]
Summing from $j = 1$ to $j = L-1$ leads to
\begin{equation}\label{eq:L-iteration-proof-1}
\begin{split}
E(\sigma, v_{L-1}) - E(\sigma, v_{L})\geq\ & E(\sigma, v) - E(\sigma, \cR(\sigma, v, 2^{-L + 1}\fB_{L})) \\
&- \frac{L-1}{\kappa}\big[ E(\sigma, v) - E(\sigma, v_{L-1}) \big]^{\frac{1}{2}},
\end{split}
\end{equation}
which gives the first inequality in~\eqref{eq:many-iteration-1}. Using the fact that $v =v_{L}$ outside of the union of all the disks in $\cup_{j = 1}^L\fB_j$, we estimate
\[
\begin{split}
E(\sigma, v) - E(\sigma,v_{L})
=\ & \frac{1}{2}\int_{\cup_{j = 1}^L \cup_{B \in \fB_{j}}B}|\nabla v|_{\sigma}^2 - |\nabla v_{L} |_{\sigma}^2 \vol_{\sigma} \\
\leq\ & \frac{1}{2}\sum_{j= 1}^{L}\sum_{B \in \fB_{j}} \int_{B}|\nabla v|_{\sigma}^2 \vol_{\sigma} \leq \frac{L \ep_0}{2 \cdot 3^{L-1}} < \ep_0.
\end{split}
\]
Thus,
\[
E(\sigma, v_{L-1}) - E(\sigma, v_{L}) \leq E(\sigma, v) - E(\sigma, v_{L}) \leq \ep_0^{\frac{1}{2}}\cdot \big[ E(\sigma, v) - E(\sigma, v_{L})  \big]^{\frac{1}{2}}.
\]
Combining this with \eqref{eq:L-iteration-proof-1} gives the second inequality in \eqref{eq:many-iteration-1}. The proof of~\eqref{eq:many-iteration-2} is similar. Namely, given $j = 1, \cdots, L-1$ as above, we apply instead~\eqref{eq:iterated-replacement-estimate-2}, with $v, \fB_1$, and $\fB_2$ taken to be $v_{L - j  -1}, \fB_{L - j}$, and $2^{-L + j + 1}\fB_L$, to get 
\[
\begin{split}
&E(\sigma, v_{L - j}) - E(\sigma, \cR(\sigma, v_{L-j}, 2^{-L + j}\fB_{L}))\\
\leq\ & E(\sigma, v_{L - j - 1}) - E(\sigma, \cR(\sigma, v_{L - j - 1}, 2^{-L + j + 1}\fB_{L})) + \frac{1}{\kappa}\big[ E(\sigma, v_{L - j - 1}) - E(\sigma, v_{L - j}) \big]^{\frac{1}{2}}.
\end{split}
\]
Summing from $j = 1$ to $j = L-1$ gives~\eqref{eq:many-iteration-2}.
\end{proof}

\section{Semi-continuity of the maximal energy drop}\label{appendix:semi-continuity}
Here we give the proof of Proposition \ref{prop:e-semicontinuous}. 
\begin{proof}[Proof of Proposition~\ref{prop:e-semicontinuous}]
We adopt the same notation as in the diagrams~\eqref{eq:Phi-diagram} and~\eqref{eq:n-Phi-diagram} from Section~\ref{subsec:mappings}, without repeating the relevant definitions. A few remarks to set the stage: 
\vskip 1mm
\begin{enumerate}
\item[(i)] Since $\sigma_n \to \sigma$ smoothly by assumption, we have by Proposition \ref{prop:metric-qc-harmonic} that $\widetilde{f}_n \to \widetilde{f}$ and $(\widetilde{f}_n)^{-1} \to (\widetilde{f})^{-1}$, smoothly locally on $\HH$. As a result, letting $k_n = \widetilde{f} \circ (\widetilde{f}_{n})^{-1}$, we see that 
\begin{equation}\label{eq:semi-conti-kn-converge}
k_n,\ k_n^{-1} \to \id,\quad \text{in }C^{\infty}_{\loc}(\HH).
\end{equation}
\vskip 1mm
\item[(ii)] To continue, recall that $\Gamma_0$ is the group of deck transformations of $p_0$, while 
\begin{equation}\label{eq:semi-conti-groups}
\Gamma_{\tau_n} := \widetilde{f}_n \cdot \Gamma_0 \cdot \widetilde{f}_n^{-1},\quad \Gamma_{\tau} := \widetilde{f} \cdot \Gamma_0 \cdot \widetilde{f}^{-1},
\end{equation}
are respectively the deck transformations of $p_n$ and $p$. 
\vskip 1mm
\item[(iii)] By Proposition \ref{prop:fundamental-domain}, there exists a compact set $F_0 \subset \HH$ with the property that it intersects every orbit of $\Gamma_0$ on $\HH$, and it follows that the compact sets
\[
F_n = \widetilde{f}_n(F_0),\quad F = \widetilde{f}(F_0),
\]
enjoy a similar property, with $\Gamma_0$ replaced by $\Gamma_{\tau_n}$ and $\Gamma_{\tau}$, respectively. By the local smooth convergence of $\widetilde{f}_n$ to $\widetilde{f}$, we may find another compact set $K \subset \HH$ such that
\begin{equation}\label{eq:fixed-compact-set}
F \subset K,\quad \text{and}\quad F_n \subset K\quad \text{for all }n.    
\end{equation}
\end{enumerate}
With these in mind we begin the actual proof. Note that since $v$ is not a harmonic map with respect to $\sigma$, we have by Theorem~\ref{thm:convexity} that $e(\sigma, v, \ep, \lambda) > 0$; see also Remark~\ref{rmk:v-replacement}(2). We shall prove that 
\begin{equation}\label{eq:e-delta-semicontinuous}
\limsup_{n \to \infty} e(\sigma_n, v_n, \frac{\ep}{3}, \frac{\lambda}{4}) \leq (1 + \delta)e(\sigma, v, \ep,  \lambda), \text{ for all }\delta > 0.
\end{equation}
Suppose towards a contradiction that~\eqref{eq:e-delta-semicontinuous} fails for some $\delta > 0$. Then, up to taking a subsequence of $(\sigma_n, v_n)$, we get 
\[
e(\sigma_n, v_n, \frac{\ep}{3}, \frac{\lambda}{4}) > (1 + \delta) e(\sigma, v, \ep, \lambda),\quad \text{for all }n,
\]
which means for all $n$ we can find a finite disjoint collection $\fB_n = \{B_{n, i}\}_{i \in I_n}$ of geodesic disks in $(S, \sigma_n)$, each with radius at most $\rho_0$, so that, upon also choosing a lift $\widetilde{B}_{n, i} = \overline{B_{g_{\hyp}}(z_{n, i}, r_{n, i})}$ for each $B_{n, i}$ with respect to $p_n$, we have
\begin{equation}\label{eq:n-energy-bound}
\sum_{i}\int_{B_{n, i}}|\nabla v_n|_{\sigma_n}^2 \vol_{\sigma_n} = \sum_{i} \int_{\widetilde{B}_{n, i}} |\nabla u_n|^2 < \frac{\ep}{3},
\end{equation}
\begin{equation}\label{eq:energy-drop-lower}
\begin{split}
(1 + \delta) e(\sigma, v, \ep, \lambda) < \ & 
E(\sigma_n, v_n) - E(\sigma_n, \cR(\sigma_n, v_n, \frac{\lambda}{4}\fB_{n}))\\
=\ & \frac{1}{2}\sum_{i} \int_{\frac{\lambda}{4}\widetilde{B}_{n, i}} |\nabla u_n|^2 - |\nabla \cR(u_n, \frac{\lambda}{4}\widetilde{B}_{n, i})|^2.
\end{split}
\end{equation}
Notice that $r_{n, i} \leq \rho_0$ for each $n, i$. Furthermore, replacing $\widetilde{B}_{n, i}$ with $\gamma(\widetilde{B}_{n, i})$ for a suitable $\gamma \in \Gamma_{\tau_n}$ if needed, we can assume without loss of generality that 
\[
z_{n, i} \in F_n,\quad (\text{ so that }k_n(z_{n, i}) \in F\ ).
\]
With $K$ as in \eqref{eq:fixed-compact-set}, we choose $R > 0$ such that
\begin{equation}\label{eq:fixed-compact-2}
\{z \in \HH\ |\ d_{\hyp}(z, K) \leq 4\rho_0 \} \subset \overline{B_{g_{\hyp}}((0, 1), R)} =: K', \quad \text{ for all }n,
\end{equation}
where $d_{\hyp}$ denotes distance with respect to $g_{\hyp}$. 

\begin{claim}\label{claim:semi-conti-inclusion}
Provided $n$ is large enough, we have for all $i \in I_n$ and $r \in (0, r_{n, i}]$ that
\[
k_n\big(B_{g_{\hyp}}(z_{n, i}, r) \big) \subset B_{g_{\hyp}}(k_n(z_{n, i}), 2r) \subset k_n\big(B_{g_{\hyp}}(z_{n, i}, 4r) \big).
\]
\end{claim}
\begin{proof}[Proof of Claim]
With $K'$ as in \eqref{eq:fixed-compact-2}, there exists $\Lambda = \Lambda_{K'} >0$ such that 
\[
\Lambda^{-1} \cdot g_{\euc} \leq g_{\hyp} \leq \Lambda \cdot g_{\euc}, \quad \text{ on }K'.
\]
Also, by \eqref{eq:semi-conti-kn-converge} we see that both $(k_n)^* g_{\hyp}$ and $(k_n^{-1})^*g_{\hyp}$ converge to $g_{\hyp}$ uniformly on $K'$. From these observations, and the fact that $K'$ is geodesically convex with respect to $g_{\hyp}$, it follows that there exists $N \in \NN$ so that provided $n \geq N$, we have for all $z, w \in K'$ that
\begin{equation}\label{eq:semi-conti-distance-ratio}
d_{\hyp}(k_n(z), k_n(w))\leq 2 d_{\hyp}(z, w),\quad d_{\hyp}(k_n^{-1}(z), k_n^{-1}(w))\leq 2 d_{\hyp}(z, w).
\end{equation}
Fix any $n \geq N$ and $i \in I_n$. Given $r \in (0, r_{n, i}]$ and $z \in B_{g_{\hyp}}(k_n(z_{n, i}), 2r)$, since $k_n(z_{n, i}) \in F$ and $r_{n, i} \leq \rho_0$, we see from \eqref{eq:fixed-compact-set} and \eqref{eq:fixed-compact-2} that $z$ and $k_n(z_{n, i})$ both lie in $K'$, and by the second estimate in \eqref{eq:semi-conti-distance-ratio} we conclude that $k_n^{-1}(z) \in B_{g_{\hyp}}(z_{n, i}, 4r)$. This proves the second asserted inclusion. The first inclusion is proved similarly and we omit the details.
\end{proof}

To continue, define $\widehat{B}_{n, i} = \overline{B_{g_{\hyp}}(k_n(z_{n, i}), \frac{r_{n, i}}{2})}$. By Claim \ref{claim:semi-conti-inclusion}, as soon as $n$ is large enough, we have for all $i \in I_n$ that
\begin{equation}\label{eq:ball-inclusion}
k_n\big( \frac{\lambda}{4}\widetilde{B}_{n, i} \big) \subset \lambda \widehat{B}_{n, i},\quad   \widehat{B}_{n,i} \subset k_n(\widetilde{B}_{n, i}).
\end{equation}
Therefore, since $\{\widetilde{B}_{n, i}\}_{i \in I_n}$ is $\tau_n$-admissible in the sense of Definition \ref{defi:admissible-disks-upstairs}, and since $\Gamma_{\tau} = k_n \cdot \Gamma_{\tau_n} \cdot (k_n)^{-1}$, eventually $\{\widehat{B}_{n, i}\}_{i \in I_n}$ is $\tau$-admissible. Thus
\[
\fB'_n := \{p(\widehat{B}_{n, i}) \}_{i \in I_n}
\]
is a disjoint collection of geodesic disks in $(S, \sigma)$, each with radius at most $\frac{\rho_0}{2}$. 
\begin{claim}\label{claim:semi-conti-remaining-step-1}
For all large enough $n$, we have
\begin{equation}\label{eq:energy-drop-upper}
\begin{split}
e(\sigma, v, \ep, \lambda) \geq \frac{1}{2}\sum_{i \in I_n} \int_{\lambda\widehat{B}_{n, i}} |\nabla u|^2 - |\nabla \cR(u, \lambda\widehat{B}_{n, i})|^2.
\end{split}
\end{equation}
\end{claim}
\begin{proof}[Proof of Claim]
By our assumption and Proposition \ref{prop:coming-down}, we see that $u_n$ converges to $u$ in the $C^0 \cap W^{1, 2}$-topology  locally on $\HH$. Combining this with \eqref{eq:semi-conti-kn-converge}, and arguing as in the proof of Proposition \ref{prop:coming-down}, we get
\begin{equation}\label{eq:un-kn-converge}
u_n \circ k_n^{-1} \to u\quad \text{in } C^0 \cap W^{1, 2} \text{ locally on }\HH.
\end{equation}
Now, by Claim \ref{claim:semi-conti-inclusion}, and the inclusions~\eqref{eq:fixed-compact-set} and~\eqref{eq:fixed-compact-2}, we have for all sufficiently large $n$ that
\begin{equation}\label{eq:disks-contained-in-K'}
\cup_{i \in I_n} k_n(\widetilde{B}_{n, i}) \subset K', \quad \cup_{i \in I_{n}} \widetilde{B}_{n, i} \subset K'.
\end{equation}
Therefore \eqref{eq:un-kn-converge} implies that
\begin{equation}\label{eq:n-energy-compare-1}
\begin{split}
\lim_{n \to \infty}  \int_{\cup_i k_n(\widetilde{B}_{n, i})}\Big| |\nabla (u_n \circ k_n^{-1})|^2 -  |\nabla u|^2 \Big| = 0.
\end{split}
\end{equation}
Also, by the chain rule and the area formula, and using the convergence \eqref{eq:semi-conti-kn-converge} again, it can be shown that
\[
\begin{split}
\lim_{n \to \infty}\Big|\int_{\cup_ik_n(\widetilde{B}_{n, i})} |\nabla (u_n \circ k_n^{-1})|^2 - \int_{\cup_i \widetilde{B}_{n, i}} |\nabla u_n|^2 \Big| = 0,
\end{split}
\]
which together with~\eqref{eq:n-energy-bound} implies that eventually
\[
\sum_{i}\int_{k_n(\widetilde{B}_{n, i})}|\nabla(u_n \circ k_n^{-1})|^2 < \frac{2\ep}{3}.
\]
Combining this with~\eqref{eq:n-energy-compare-1} and~\eqref{eq:ball-inclusion}, we get for large enough $n$ that
\begin{equation}\label{eq:energy-bound-limit}
\sum_{B \in \fB_n'} \int_{B}|\nabla v|_{\sigma}^2 \vol_{\sigma} =    \sum_{i} \int_{\widehat{B}_{n, i}} |\nabla u|^2 \leq \int_{\cup_{i}k_n(\widetilde{B}_{n, i})} |\nabla u|^2 < \ep.
\end{equation}
The desired conclusion follows from the definition of $e(\sigma, v, \ep, \lambda)$, along with the fact that
\[
\frac{1}{2}\sum_{i \in I_n} \int_{\lambda\widehat{B}_{n, i}} |\nabla u|^2 - |\nabla \cR(u, \lambda\widehat{B}_{n, i})|^2= E(\sigma, v) - E(\sigma, \cR(\sigma, v, \lambda\fB_n')).
\]
\end{proof}
\begin{claim}\label{claim:semi-conti-remaining-step-2}
We have 
\begin{equation}\label{eq:n-energy-compare-2}
\begin{split}
e(\sigma, v, \ep, \lambda) \geq\ & \limsup_{n \to \infty}\frac{1}{2}\sum_{i\in I_n}  \int_{\lambda\widehat{B}_{n, i}} |\nabla (u_n \circ k_n^{-1})|^2 - |\nabla \cR(u_n\circ k_n^{-1}, \lambda\widehat{B}_{n, i})|^2.
\end{split}
\end{equation}
\end{claim}
\begin{proof}[Proof of Claim]
By the inclusions~\eqref{eq:ball-inclusion} and the triangle inequality we have
\begin{equation}\label{eq:drop-difference}
\begin{split}
&\Big|\sum_{i} \int_{\lambda\widehat{B}_{n, i}} |\nabla (u_n \circ k_n^{-1})|^2 - |\nabla \cR(u_n\circ k_n^{-1}, \lambda\widehat{B}_{n, i})|^2 -  \sum_{i} \int_{\lambda\widehat{B}_{n, i}} |\nabla u|^2 - |\nabla \cR(u, \lambda\widehat{B}_{n, i})|^2\Big|\\
\leq\ & \int_{\cup_{i}k_n(\widetilde{B}_{n, i})} \Big| |\nabla (u_n \circ k_n^{-1})|^2 - |\nabla u|^2 \Big|\\
&+ \sum_{i}\Big|\int_{\lambda \widehat{B}_{n,i}} |\nabla \cR(u_n \circ k_n^{-1}, \lambda\widehat{B}_{n,i})|^2 - \int_{\lambda\widehat{B}_{n, i}}|\nabla \cR(u, \lambda\widehat{B}_{n,i})|^2 \Big|.
\end{split}
\end{equation}
To bound the last line, we use the estimate \eqref{eq:replaced-energy} on each $\lambda \widehat{B}_{n, i}$ to get
\[
\begin{split}
&\sum_{i \in I_n}\Big|\int_{\lambda \widehat{B}_{n,i}} |\nabla \cR(u_n \circ k_n^{-1}, \lambda\widehat{B}_{n,i})|^2 - \int_{\lambda\widehat{B}_{n, i}}|\nabla \cR(u, \lambda\widehat{B}_{n,i})|^2 \Big|\\
\leq\ & C\sum_{i \in I_{n}}\|u - u_n\circ k_n^{-1}\|_{\infty; \lambda\widehat{B}_{n, i}} \Big( \int_{\lambda\widehat{B}_{n, i}} |\nabla (u_n\circ k_n^{-1})|^2 + |\nabla u|^2 \Big)\\
&+ C\sum_{i \in I_n}\Big(\int_{\lambda\widehat{B}_{n, i}} |\nabla (u_n \circ k_n^{-1}) - \nabla u|^2 \Big)^{\frac{1}{2}}\Big( \int_{\lambda\widehat{B}_{n, i}} |\nabla (u_n\circ k_n^{-1})|^2 + |\nabla u|^2 \Big)^{\frac{1}{2}}\\
\leq\ & C\|u - u_n \circ k_n^{-1}\|_{\infty; \cup_{i}\widehat{B}_{n, i}}\Big( \int_{\cup_{i}\widehat{B}_{n, i}} |\nabla (u_n\circ k_n^{-1})|^2 + |\nabla u|^2\Big)\\
& + C\Big(\int_{\cup_{i}\widehat{B}_{n, i}} |\nabla (u_n\circ k_n^{-1}) - \nabla u|^2 \Big)^{\frac{1}{2}}\Big( \int_{\cup_{i}\widehat{B}_{n, i}} |\nabla (u_n \circ k_n^{-1})|^2 + |\nabla u|^2 \Big)^{\frac{1}{2}},
\end{split}
\]
where we applied the Cauchy--Schwarz inequality to get the second inequality. Noting that $\cup_{i}\widehat{B}_{n, i} \subset K'$ by~\eqref{eq:ball-inclusion} and~\eqref{eq:disks-contained-in-K'}, and using the convergence \eqref{eq:un-kn-converge}, we deduce from the above estimate that 
\begin{equation}\label{eq:replacement-energy-difference-lim}
\sum_{i \in I_n}\Big|\int_{\lambda \widehat{B}_{n,i}} |\nabla \cR(u_n \circ k_n^{-1}, \lambda\widehat{B}_{n,i})|^2 - \int_{\lambda\widehat{B}_{n, i}}|\nabla \cR(u, \lambda\widehat{B}_{n,i})|^2 \Big| \to 0 \text{ as }n \to \infty.
\end{equation}
Combining this and~\eqref{eq:n-energy-compare-1} with~\eqref{eq:drop-difference}, while also using Claim \ref{claim:semi-conti-remaining-step-1}, we arrive at \eqref{eq:n-energy-compare-2} as asserted.
\end{proof}
\begin{claim}\label{claim:semi-conti-remaining-step-3}
We have
\[
\begin{split}
&\limsup_{n \to \infty} \sum_{i \in I_n}\int_{\frac{\lambda}{4}\widetilde{B}_{n, i}} |\nabla u_n|^2 - |\nabla \cR(u_n, \frac{\lambda}{4}\widetilde{B}_{n, i})|^2\\
\leq\ & \limsup_{n \to \infty}\sum_{i\in I_n}  \int_{\lambda\widehat{B}_{n, i}} |\nabla (u_n \circ k_n^{-1})|^2 - |\nabla \cR(u_n\circ k_n^{-1}, \lambda\widehat{B}_{n, i})|^2.
\end{split}
\]
\end{claim}
\begin{proof}[Proof of Claim]
Since $\cR(u_n, \frac{\lambda}{4}\widetilde{B}_{n, i})$ agrees with $u_n$ outside of $\frac{\lambda}{4}\widetilde{B}_{n, i}$, the latter contained in $k_n^{-1}(\lambda\widehat{B}_{n, i})$ due to~\eqref{eq:ball-inclusion}, we infer that
\begin{equation}\label{eq:un-kn-agreement}
\cR(u_n, \frac{\lambda}{4}\widetilde{B}_{n, i}) \circ k_n^{-1}(z) = u_n \circ k_n^{-1}(z)\quad \text{for all  }z \not\in \lambda \widehat{B}_{n, i}.
\end{equation}
Thus, using the energy minimizing property of $\cR(u_n\circ k_n^{-1}, \lambda\widehat{B}_{n, i})$, we have
\begin{equation}\label{eq:n-energy-compare-3}
\begin{split}
&  \int_{\lambda \widehat{B}_{n,i}}  |\nabla (u_n \circ k_n^{-1})|^2 - |\nabla \cR(u_n\circ k_n^{-1}, \lambda\widehat{B}_{n, i})|^2\\
\geq\ &\int_{\lambda \widehat{B}_{n, i}} |\nabla (u_n \circ k_n^{-1})|^2 - |\nabla [\cR(u_n, \frac{\lambda}{4}\widetilde{B}_{n, i}) \circ k_n^{-1}]|^2\\
=\ & \int_{k_{n}^{-1}(\lambda \widehat{B}_{n, i})} |\nabla u_n|_{k_n^*g_{\euc}}^2 - |\nabla \cR(u_n, \frac{\lambda}{4}\widetilde{B}_{n, i})|_{k_n^*g_{\euc}}^2 \vol_{k_n^*g_{\euc}}\\
=\ & \int_{\frac{\lambda}{4}\widetilde{B}_{n, i}} |\nabla u_n|_{k_n^*g_{\euc}}^2 - |\nabla \cR(u_n, \frac{\lambda}{4}\widetilde{B}_{n, i})|_{k_n^*g_{\euc}}^2 \vol_{k_n^*g_{\euc}},
\end{split}
\end{equation}
where the last equality follows from \eqref{eq:ball-inclusion} and the fact that $u_n$ and $\cR(u_n, \frac{\lambda}{4}\widetilde{B}_{n, i})$ agree on $k_n^{-1}(\lambda \widehat{B}_{n, i}) \setminus \frac{\lambda}{4}\widetilde{B}_{n, i}$. 
With $K'$ being as in~\eqref{eq:fixed-compact-2}, and noting that 
\[
g_n:= k_n^*g_{\euc} \to g_{\euc}\quad \text{uniformly on }K',
\]
we can find a sequence $\alpha_n \to 0$ such that for all sufficiently large $n$, and any $z \in K'$ and $\xi \in \RR^2$, there holds
\[
(1 - \alpha_n) |\xi|^2_{g_{\euc}} \leq g_n^{ij}(z)\xi_i \xi_j\sqrt{\det g_n(z)} \leq (1 + \alpha_n) |\xi|^2_{g_{\euc}}.
\]
From this we get
\begin{equation}\label{eq:step-3-convergence-1}
\begin{split}
&\Big|\sum_{i \in I_n} \int_{\frac{\lambda}{4}\widetilde{B}_{n, i}} |\nabla \cR(u_n, \frac{\lambda}{4}\widetilde{B}_{n, i})|_{g_n}^2 \vol_{g_n} - \sum_{i \in I_n}\int_{\frac{\lambda}{4}\widetilde{B}_{n, i}} |\nabla \cR(u_n, \frac{\lambda}{4}\widetilde{B}_{n, i})|_{g_{\euc}}^2 \vol_{g_{\euc}}\Big|\\
\leq\ & \alpha_n \cdot \sum_{i \in I_n}  \int_{\frac{\lambda}{4}\widetilde{B}_{n, i}} |\nabla \cR(u_n, \frac{\lambda}{4}\widetilde{B}_{n, i})|_{g_{\euc}}^2  \vol_{g_{\euc}}\\
\leq\ & \alpha_n \cdot \sum_{i \in I_n}  \int_{\frac{\lambda}{4}\widetilde{B}_{n, i}} |\nabla u_n|_{g_{\euc}}^2  \vol_{g_{\euc}} \leq \alpha_n \cdot \int_{K'} |\nabla u_n|^2 \to 0 \text{ as }n \to \infty,
\end{split}
\end{equation}
where in the last step we used \eqref{eq:fixed-compact-2}. Similarly, 
\begin{equation}\label{eq:step-3-convergence-2}
\Big|\sum_{i \in I_n} \int_{\frac{\lambda}{4}\widetilde{B}_{n, i}} |\nabla u_n|_{g_n}^2 \vol_{g_n} - \sum_{i \in I_n}\int_{\frac{\lambda}{4}\widetilde{B}_{n, i}} |\nabla u_n|_{g_{\euc}}^2 \vol_{g_{\euc}}\Big| \to 0 \text{ as }n \to \infty.
\end{equation}
We finish the proof of the claim upon summing both sides of \eqref{eq:n-energy-compare-3} over $i \in I_n$, taking the limsup as $n \to \infty$, and using \eqref{eq:step-3-convergence-1} and \eqref{eq:step-3-convergence-2}.
\end{proof}
Combining Claims \ref{claim:semi-conti-remaining-step-3} and \ref{claim:semi-conti-remaining-step-2} with \eqref{eq:energy-drop-lower} gives
\[
\begin{split}
e(\sigma, v, \ep, \lambda) \geq (1 + \delta) e(\sigma, v, \ep, \lambda),
\end{split}
\]
which is a contradiction since $e(\sigma, v, \ep, \lambda) > 0$. The proof of Proposition \ref{prop:e-semicontinuous} is complete.\\
\end{proof}
\section{Operations on min-max sequences}\label{appendix:operations}
In this appendix we prove Lemma~\ref{lemm:smooth-sweepout} and Lemma~\ref{lemm:cut-and-paste}. For convenience, given $r > 0$, we define
\[
\cN_{r}(M) = \{y \in \RR^{q}\ |\ \dist(y, M) < r\}.
\]
Recall that $d_0$ is chosen so that $\overline{\cN_{4d_0}(M)}$ is contained in a tubular neighborhood of $M$. The nearest-point projection $\Pi$ onto $M$ then has all its derivatives bounded on $\overline{\cN_{4d_0}(M)}$. In particular, the assignment $v \mapsto \Pi \circ v$ is continuous from $(C^0 \cap W^{1, 2})(S; \cN_{2d_0}(M))$ into $(C^0 \cap W^{1, 2})(S; M)$. Indeed, whenever $w, v \in (C^{0} \cap W^{1, 2})(S; \cN_{2d_0}(M))$ are such that $
\|v - w\|_{C^0} < d_0$, noting from the triangle inequality that 
\[
tv(x) + (1 - t)w(x) \subset \cN_{3d_0}(M)\quad\text{for all }t \in [0, 1],\ x \in S,
\]
we may apply the fundamental theorem of calculus to bound quantities such as $|\Pi(v) - \Pi(w)|$ and $|(d\Pi)_{v} - (d\Pi)_{w}|$. Doing so results in the estimate
\begin{equation}\label{eq:Pi-uniform-continuity}
\begin{split}
\|\Pi(v) - \Pi(w)\|_{C^0 \cap W^{1, 2}} \leq\ & C_{\Pi} \cdot \|v - w\|_{C^0}\cdot (1 + \|\nabla v\|_{L^2})\\
& + C_{\Pi} \cdot \|\nabla v - \nabla w\|_{L^2}.
\end{split}
\end{equation}
A similar consideration involving second derivatives shows that $v \mapsto \Pi \circ v$ also defines a continuous map from $C^2(S; \cN_{2d_0}(M))$ into $C^{2}(S; M)$.

\begin{proof}[Proof of Lemma~\ref{lemm:smooth-sweepout}]
Let $\{\varphi_{\alpha}: \bB_{3} \subset \RR^2 \to S\}_{\alpha = 1, \cdots, L}$ be a collection of charts such that 
\[
{S} \subset \cup_{\alpha = 1}^{L}\varphi_{\alpha}(\bB_1),
\]
and suppose $\{\zeta_{\alpha}\}_{\alpha = 1, \cdots, L}$ is a partition of unity on $S$ subordinate to the covering. Also, take $\{\eta_{\ep} = \frac{1}{\ep^2}\eta(\frac{\cdot}{\ep})\}_{\ep > 0}$ to be a family of standard mollifiers on $\RR^2$, with $\supp(\eta) \subset \bB_{1}$. Given $u \in C^0 \cap W^{1,2}({S}; \RR^q)$, $\alpha \in \{1, \cdots, L\}$, and $\ep \in (0, 1)$, we define
\[
\widetilde{u}_{\alpha} : = u \circ \varphi_{\alpha} \in (C^0 \cap W^{1, 2})(\bB_3; \RR^q),
\]
\[
\widetilde{u}_{\alpha, \ep}: = \eta_{\ep} * \widetilde{u}_{\alpha} \in C^{\infty}(\bB_2;\RR^q).
\]
The following facts are standard. First, as $\ep \to 0$ we have
\begin{equation}\label{eq:mollified-C0-W12-converge}
\widetilde{u}_{\alpha, \ep} \to \widetilde{u}_{\alpha} \text{ in }C^0 \cap W^{1, 2}(\bB_2; \RR^q).
\end{equation}
Next, we have the bounds
\begin{align*}
\|\widetilde{u}_{\alpha, \ep}\|_{C^3(\bB_2)} \leq\ & C \ep^{-3}\|\widetilde{u}_{\alpha}\|_{C^0(\bB_3)}, \\ 
\|\widetilde{u}_{\alpha, \ep}\|_{C^0 \cap W^{1, 2}(\bB_2)} \leq\ & \|\widetilde{u}_{\alpha}\|_{C^0 \cap W^{1, 2}(\bB_3)},
\end{align*}
where the constant $C$ depends only on the mollifier $\eta$. These then imply the following continuity properties. Suppose $\ep_i \to \ep \in [0, 1)$, and $u_i \to u$ in $(C^0 \cap W^{1, 2})({S}; \RR^q)$, then
\[
(\widetilde{u}_i)_{\alpha, \ep_i} \to \widetilde{u}_{\alpha, \ep} \text{ in }C^2(\bB_2; \RR^q), \text{ if }\ep > 0;
\]
\[
(\widetilde{u}_i)_{\alpha, \ep_i} \to \widetilde{u}_{\alpha} \text{ in }(C^0 \cap W^{1, 2})(\bB_2; \RR^q), \text{ if }\ep = 0.
\]
Therefore, upon defining $\mathscr{S}:[0, 1) \times (C^0 \cap W^{1, 2})({S}; \RR^q)  \longrightarrow (C^0 \cap W^{1, 2})({S}; \RR^q)$ by
\[
\mathscr{S}(\ep, u): = \left\{
\begin{array}{ll}
\sum_{\alpha = 1}^{L} \zeta_{\alpha} \cdot\big(\widetilde{u}_{\alpha, \ep} \circ (\varphi_{\alpha}|_{\bB_2})^{-1}\big), & \text{ if }\ep \in (0, 1) ,\\
u, & \text{ if }\ep = 0,
\end{array}
\right.
\]
we see that $\mathscr{S}$ is continuous, and that the restriction $\mathscr{S}|_{(0, 1) \times (C^0 \cap W^{1, 2})({S}; \RR^q)}$ is continuous as a map into $C^2({S}; \RR^q)$.

Now suppose $\bv \in \cP$, $\delta \in (0, \frac{1}{4})$, and $\mu \in (0, d_0)$, as in the statement of Lemma~\ref{lemm:smooth-sweepout}. We let $f: I^m \to [0, 1]$ be a continuous function such that
\[
f(t) = 1 \text{ on }I_{\delta}^m,\quad f(t) = 0 \text{ outside of }I_{2\delta/3}^m,
\]
choose any positive sequence $\ep_i \to 0$, and define, for $(s,t) \in [0, 1] \times \Int(I^m)$,
\[
\widehat{\bh}_i(s, t) =\left\{
\begin{array}{ll}
\mathscr{S}(s \ep_i \cdot f(t), \bv(t)), & \text{ if } t \in  I_{\delta/2}^m,\\
\bv(t), & \text{ if } t \in \Int(I^m) \setminus I_{\delta/2}^m.
\end{array}
\right.
\]
Since $\bv$ is continuous and $I_{\delta/2}^m$ is compact, we infer from the continuity of $\mathscr{S}$ that
\begin{equation}\label{eq:choice-of-ep-bar}
\sup_{(s, t) \in [0, 1] \times I_{\delta/2}^m}\|\widehat{\bh}_i(s, t) - \bv(t)\|_{C^0 \cap W^{1, 2}} \to 0 \text{ as }i \to \infty.
\end{equation}
Combining this with the fact that $\widehat{\bh}_i(s, t) = \bv(t)$ when $t \not\in I_{2\delta/3}^m$, we see that $\widehat{\bh}_{i}(s, t)$ maps $S$ into $\cN_{2d_0}(M)$ for all $(s, t) \in [0, 1] \times \Int(I^m)$ and sufficiently large $i$, which allows us to define
\[
\bh_i(s, t) = \Pi \circ \widehat{\bh}_i(s, t), \text{ for all }(s, t) \in [0,1] \times \Int(I^m).
\]
With the help of the continuity properties of $\mathscr{S}$, and the remarks made at the start of this appendix, each $\bh_i$ is continuous as a map from $[0, 1] \times \Int(I^m)$ to $(C^0 \cap W^{1, 2})({S}; M)$, and restricts to a continuous map from $(0, 1] \times I_{\delta}^m$ into $C^2({S}; M)$. Furthermore,  we have
\begin{equation}\label{eq:hi-boundary}
\bh_i(s, t) = \bv(t), \text{ whenever } t \not\in I_{2\delta/3}^m \text{ or }s = 0.
\end{equation}
Therefore, in view of Definition \ref{defi:relation-on-P}, we see that letting 
\[
\bw_{i}(t) = 
\left\{
\begin{array}{ll}
\bh_{i}(1, t),& \text{ if }t \in \Int(I^m),  \\
\bv(t),& \text{ if }t \in \partial I^m,
\end{array}
\right.
\]
produces an element of $[\bv]$ satisfying conclusions (a) and (b). Next, thanks to~\eqref{eq:choice-of-ep-bar} and the estimate~\eqref{eq:Pi-uniform-continuity}, we see that 
\[
\sup_{t\in  I_{\delta/2}^m}\|\bw_{i}(t) - \bv(t)\|_{C^0 \cap W^{1,2}} \to 0 \text{ as } i \to \infty.
\]
We complete the proof upon taking $\bw = \bw_{i}$ with a sufficiently large $i$ so that conclusion (c) is fulfilled as well.
\end{proof}
\vskip 1em

\begin{proof}[Proof of Lemma~\ref{lemm:cut-and-paste}]
Take a non-negative increasing $\zeta \in C^{\infty}(\RR; [0, 1])$ such that 
\[
\zeta(t) = 0, \text{ if }t \leq 3\tau,\ \ \zeta(t) = 1, \text{ if }t \geq 1 - 3\tau.
\]
Also, by suitably mollifying $t \mapsto \max\{1, t\}$, we obtain a smooth, increasing function $\varphi: \RR \to \RR$ such that 
\[
\varphi(t) = 1, \text{ if }t \leq 1 - \tau,\ \ \varphi(t) = t, \text{ if } t \geq 1 + \tau.
\]
Below we first define an operator 
\begin{equation}\label{eq:T-domains}
\widehat{T}: (0, 1] \times C^2(\bB_2; \RR^q) \to C^2(\bB_2; \RR^q),
\end{equation}
out of which the desired $T$ is constructed. Given $u \in C^2(\bB_2; \RR^q)$ and $\rho \in (0, 1]$, we let
\begin{equation}\label{eq:T-rho-definition}
\widehat{T}(\rho, u)(x)=  \left\{
\begin{array}{ll}
u(0), & \text{ if } |x| \leq 2\tau \rho,\\
\zeta(\frac{|x|}{\rho}) u(\rho\frac{x}{|x|}) + [1 - \zeta(\frac{|x|}{\rho}) ] u(0), & \text{ if } 2\tau\rho \leq |x| \leq (1 - 2\tau)\rho,\\
u(\rho\varphi(\frac{|x|}{\rho})\frac{x}{|x|}), & \text{ if }(1- 2\tau)\rho \leq |x|.
\end{array}
\right.
\end{equation}
Note that, first of all, $\widehat{T}$ has the following scaling property:
\begin{equation}\label{eq:T-scaling}
\widehat{T}(\rho, u)(x) = \widehat{T}(1, u(\rho\ \cdot))(\frac{x}{\rho}), \text{ for }x \in \bB_{2\rho}.
\end{equation}

Secondly, the pieces in the definition~\eqref{eq:T-rho-definition} join smoothly, so in particular $\widehat{T}(\rho, u)$ does lie in $C^2(\bB_2;\RR^q)$. For later use, we observe more specifically that
\begin{equation}\label{eq:T-regions}
\widehat{T}(\rho, u)(x) = 
\left\{
\begin{array}{ll}
u(0), & \text{ if }|x| \leq 3\tau\rho, \\
u(\rho\frac{x}{|x|}), & \text{ if } (1 - 3\tau)\rho \leq |x| \leq (1 - \tau)\rho,\\
u(x), & \text{ if } (1 + \tau)\rho \leq |x|.
\end{array}
\right.
\end{equation}

Thirdly, it is straightforward to work out from~\eqref{eq:T-rho-definition} that
\begin{equation}\label{eq:T-rho-estimates}
\|\widehat{T}(\rho, u)\|_{C^1(\bB_2)} \leq C\|u\|_{C^1(\bB_2)},\quad \|D^2 \widehat{T}(\rho, u)\|_{C^0(\bB_2)} \leq C\rho^{-1}\|u\|_{C^2(\bB_2)},
\end{equation}
where $C$ does not depend on $\rho$. These bounds in turn imply that $\widehat{T}$ is continuous with respect to the domain and codomain in~\eqref{eq:T-domains}. Indeed, suppose $\rho_i \to \rho \in (0, 1]$ and $u_i \to u$ in $C^2(\bB_2;\RR^q)$. Then we eventually have
\[
\begin{split}
\|\widehat{T}(\rho_i, u_i) - \widehat{T}(\rho, u)\|_{C^2(\bB_2)} \leq\ & 
\|\widehat{T}(\rho_i, u_i) - \widehat{T}(\rho_i, u)\|_{C^2(\bB_2)} + \|\widehat{T}(\rho_i, u) - \widehat{T}(\rho, u)\|_{C^2(\bB_2)} \\
\leq\ & C\rho^{-1}\|u - u_i\|_{C^2(\bB_2)} + \|\widehat{T}(\rho_i, u) - \widehat{T}(\rho, u)\|_{C^2(\bB_{(1 + 2\tau)\rho})}.
\end{split}
\]
Here in getting the second inequality we used~\eqref{eq:T-rho-estimates} and the linearity of $\widehat{T}$ in the second variable, as well as the last case in~\eqref{eq:T-regions}. By the $C^2$-convergence $u_i \to u$ together with~\eqref{eq:T-scaling}, we conclude that $\widehat{T}(\rho_i, u_i) \to \widehat{T}(\rho, u)$ in $C^2$, proving the asserted continuity.

Lastly, since $\varphi$ is increasing, we have 
\[
\rho\varphi(\frac{r}{\rho}) \leq (1 + \tau)\rho\quad\text{whenever } r \leq (1 + \tau)\rho.
\]
Combining this with~\eqref{eq:T-rho-definition} and the last case in~\eqref{eq:T-regions}, we infer that 
\begin{equation}\label{eq:T-u-difference-estimate}
\begin{split}
\|\widehat{T}(\rho, u) - u\|_{\infty; \bB_2} =\ & \|\widehat{T}(\rho, u) - u\|_{\infty; \bB_{(1 + \tau)\rho}} \leq \osc_{\bB_{(1 + \tau)\rho}}u\leq 4\rho \|Du\|_{\infty; \bB_2}.
\end{split}
\end{equation}

Now, given $K > 0$, consider the space
\[
\cX_{K} = \{u \in C^2(\bB_2; M)\ |\ \|Du\|_{\infty; \bB_2} \leq K\}.
\]
By~\eqref{eq:T-u-difference-estimate}, if we require
\[
\overline{\rho} \leq \min\big\{\frac{d_0}{4(K+ 1)}, 1\big\}, 
\]
then it makes sense to define, for all $\rho \in (0,\overline{\rho}]$ and $u \in \cX_{K}$,
\begin{equation}\label{T-nohat-definition}
T(\rho, u) = \Pi \circ \widehat{T}(\rho, u).
\end{equation}
Observe that $T:(0, \overline{\rho}] \times \cX_{K} \to C^2(\bB_2; M)$ is a continuous map. Also, conclusion (a) follows straight from the first and third cases in~\eqref{eq:T-regions}.

To see property (b), we begin by estimating the mapping area of $\widehat{T}(\rho, u)$. To that end we employ polar coordinates and notice from the second case in~\eqref{eq:T-regions} that
\begin{equation}\label{eq:area-bound-0}
\widehat{T}(\rho, u)_{r} \wedge \widehat{T}(\rho, u)_{\theta} = 0, \text{ when }(1-3\tau)\rho \leq r \leq (1 - \tau)\rho.
\end{equation}
On the other hand, since $\varphi' \geq 0$, by the change of variables $s = \rho\varphi(\frac{r}{\rho})$ we have
\[
\int_{(1-\tau)\rho}^{(1 + \tau)\rho}\int_{0}^{2\pi} |\widehat{T}(\rho, u)_{r} \wedge r^{-1}\widehat{T}(\rho, u)_{\theta}| \, r d\theta dr = \int_{\rho}^{(1 + \tau)\rho}\int_{0}^{2\pi} |u_r \wedge u_\theta| \, d\theta dr.
\]
Combining this with~\eqref{eq:area-bound-0}, and recalling that $\widehat{T}(\rho, u)= u$ outside of $\bB_{(1 +\tau)\rho}$, we obtain
\[
\begin{split}
A(\widehat{T}(\rho, u); \bB_2 \setminus \bB_{(1 - 3\tau)\rho}) =\ & 
A(u; \bB_2 \setminus \bB_{(1 +\tau)\rho}) + A(\widehat{T}(\rho, u); \bB_{(1 + \tau)\rho} \setminus \bB_{(1 - \tau)\rho}) \leq A(u).
\end{split}
\]
Now since $u \in \cX_{K}$ maps into $M$, from~\eqref{eq:T-rho-definition} and~\eqref{eq:T-regions} we see that $\widehat{T}(\rho, u) = T(\rho, u)$ on $\bB_2 \setminus \bB_{(1 - 3\tau)\rho}$, and hence the above estimate immediately becomes
\begin{equation}\label{eq:T-area-estimate-1}
A(T(\rho, u); \bB_2 \setminus \bB_{(1 - 3\tau)\rho}) \leq A(u).
\end{equation}
On the other hand, by a direct computation using~\eqref{eq:T-rho-definition}, we have on $\bB_{(1 - 3\tau)\rho}$ that 
\[
\begin{split}
|T(\rho, u)_{r} \wedge T(\rho, u)_{\theta}| \leq \ & \|d\Pi\|_{\infty; \cV}^2|\widehat{T}(\rho, u)_{r}| |\widehat{T}(\rho, u)_{\theta}| \leq  \|d\Pi\|_{\infty; \cV}^2 \cdot (K\|\zeta'\|_{\infty}) \cdot (\rho K),
\end{split}
\]
which together with~\eqref{eq:T-area-estimate-1} gives
\begin{equation}\label{eq:T-area-estimate-final}
A(T(\rho, u)) \leq A(u) + C \|d\Pi\|_{\infty; \cV}^2 K^2 \rho^2,
\end{equation}
where $C$ is a universal constant, and we get the estimate in conclusion (b) upon further decreasing $\overline{\rho}$ depending only on $K, \|d\Pi\|_{\infty; \cV}$ and $\mu$.

To prove the continuity asserted in (c), suppose $(\rho_i)$ is a sequence in $(0, \overline{\rho}]$ converging to $0$, and that $u_i \to u$ in $\cX_{K}$. With the help of~\eqref{eq:T-rho-estimates} and~\eqref{eq:T-u-difference-estimate}, we see that
\[
\begin{split}
\|\widehat{T}(\rho_i, u_i) - u\|_{C^0(\bB_2)}\leq\ & \|\widehat{T}(\rho_i, u_i) - \widehat{T}(\rho_i, u)\|_{C^0(\bB_2)} + \|\widehat{T}(\rho_i, u) - u\|_{C^0(\bB_2)}\\
\leq\ & C\|u - u_i\|_{C^1(\bB_2)} + \osc_{\bB_{(1 + \tau)\rho_i}}u \longrightarrow 0 \text{ as }i \to \infty.
\end{split}
\]
On the other hand, by~\eqref{eq:T-rho-estimates} and~\eqref{eq:T-regions} we have
\[
\begin{split}
\|D \widehat{T}(\rho_i, u_i) - D u\|_{L^{2}(\bB_2)} \leq\ & \|D \widehat{T}(\rho_i, u_i) - D \widehat{T}(\rho_i,u)\|_{L^2(\bB_2)} + \|D \widehat{T}(\rho_i, u) -D  u\|_{L^2(\bB_{(1 + \tau)\rho_i})}\\
 \leq\ & C\|D \widehat{T}(\rho_i, u_i) - D \widehat{T}(\rho_i,u)\|_{L^\infty(\bB_2)} + \|D \widehat{T}(\rho_i, u)\|_{L^2(\bB_{(1 + \tau)\rho_i})} \\
 &+ \|D  u\|_{L^2(\bB_{(1 + \tau)\rho_i})}\\
\leq\ & C\|u - u_i \|_{C^1(\bB_2)} + C\|u\|_{C^1(\bB_2)} \cdot \rho_i \longrightarrow 0 \text{ as }i \to \infty.
\end{split}
\]
As in the end of the previous proof, the two above convergences along with~\eqref{eq:Pi-uniform-continuity} implies that 
\[
T(\rho_i, u_i) \to u  \text{ in }(C^0 \cap W^{1, 2})(\bB_2; M),
\]
and we are done.
\end{proof}
\section{Energy of good min-max sequence on long cylinders}\label{appendix:long-cylinders}
In this appendix we give the proofs of Propositions~\ref{prop:neck} and~\ref{prop:collar-bubble}. The parameters $L \geq 1$ and $\ep_2 > 0$ appearing below are those chosen before the statement of Proposition~\ref{prop:neck}.
\begin{proof}[Proof of Proposition~\ref{prop:neck}]
Throughout the proof, by properties (p1), (p2) and (p3) we mean those listed in Definition~\ref{defi:well-prepared}. Let $k_n$ be the positive integer such that $k_n L \leq 2T_n < (k_n + 1)L$, and define $L_n = \frac{2T_n}{k_n}$. Notice that
\begin{equation}\label{eq:length-control}
L \leq L_n < \frac{k_n + 1}{k_n}L \to L \text{ as }n \to \infty.
\end{equation}
To establish~\eqref{eq:no-energy-on-neck}, it suffices by~\eqref{eq:length-control} and~\eqref{eq:no-energy-on-finite-part} to prove that
\begin{equation}\label{eq:no-energy-with-gn}
\lim_{n \to \infty}E(g_n, u_n; \bC_{-T_n + L_n, T_n -L_n}) = 0.
\end{equation}
We will first prove~\eqref{eq:no-energy-with-gn} under the extra hypothesis that 
\begin{equation}\label{eq:small-globally}
E(g_n, u_n; \bC_{-T_n, T_n}) \leq \ep_2,\quad \text{for all large enough }n.
\end{equation}
Assume, towards a contradiction, that for some $\alpha > 0$, and along some subsequence, there holds
\begin{equation}\label{eq:leftover-on-cylinder}
E(g_n, u_n; \bC_{-T_n + L_n, T_n -L_n}) \geq \alpha,\quad \text{for all }n.
\end{equation}
We divide the subsequent argument into a number of claims. 
\begin{claim}\label{claim:compactness-property}
Define $R_n := \sup_{|t| \leq T_n} (f_n(t))^{-1}$. For all sufficiently large $n$, given any finite, disjoint collection $\fB = \{B_i\}_{i \in I}$ of closed geodesic disks in $(\bC_{-T_n + L_n, T_n - L_n}, g_n)$, each with with radius at most $R_n$, the following hold.
\vskip 1mm
\begin{enumerate}
\item[(a)] With $\psi_n(\fB): = \{  \psi_n(B_i)\}_{i \in I}$, the harmonic replacement $\cR(\sigma_n, v_n, \psi_n(\fB))$ makes sense.
\vskip 1mm
\item[(b)] Let $\mu_n: = \sqrt{\frac{\delta_n}{\alpha}}$, where $\delta_n$ is given by Proposition~\ref{prop:good-min-max-seq}, and also define
\begin{equation}\label{eq:replacement-on-cylinder}
\cR(g_n, u_n, 2^{-2N-2}\fB) : = \cR(\sigma_n, v_n, 2^{-2N-2}\psi_n(\fB)) \circ \psi_n.
\end{equation}
Then we have
\[
\int_{\bC_{-T_n + L_n, T_n - L_n}} |\nabla u_n - \nabla \cR(g_n, u_n, 2^{-2N-2}\fB)|^2_{g_n} \vol_{g_n}\leq \mu_n^2 E(g_n, u_n; \bC_{-T_n + L_n, T_n - L_n}).
\]
\end{enumerate}
\end{claim}
\begin{proof}[Proof of Claim]
By (p1), we have $\lambda_n R_n \leq \rho_0$ for all large enough $n$. Fix such an $n$ and let $\fB$ be as in the statement. Then since $\psi_n$ maps $(\bC_{-T_n, T_n}, \lambda_n^2 g_n)$ isometrically into $(S, \sigma_n)$, we see that $\psi_n(\fB)$ is a finite, disjoint collection of closed geodesic disks in the latter, each with radius at most $\rho_0$. Moreover, by~\eqref{eq:small-globally} and our choice of $\ep_2$, we find that
\begin{equation}\label{eq:compactness-property-smallness}
\begin{split}
\sum_{i \in I} \int_{\psi_n(B_i)}|\nabla v_n|_{\sigma_n}^2\vol_{\sigma_n} = \sum_{i \in I} \int_{B_i} |\nabla u_n|^2_{g_n}\vol_{g_n} \leq\ & 2E(g_n, u_n; \bC_{-T_n + L_n, T_n - L_n})\\
\leq\ & 2\ep_2 < \frac{\ep_1}{3^{N + 2}}.
\end{split}
\end{equation}
In particular, harmonic replacement can be applied to $(\sigma_n, v_n)$ on the collection $\psi_n(\fB)$, as asserted in part (a). For part (b), we first observe that, for each $i \in I$, we have 
\begin{equation}\label{eq:iso-concentric}
2^{-2N-2}\psi_n(B_i) = \psi_n(2^{-2N-2}B_i),
\end{equation}
where on the right-hand side concentric scaling does not depend on the choice between $g_n$ or $\lambda_n^2 g_n$. From this we infer that $\cR(g_n, u_n; 2^{-2N-2}\fB) = u_n$ outside of $\bigcup_{i\in I}2^{-2N-2}B_i$, and thus
\[
\begin{split}
&\int_{\bC_{-T_n + L_n, T_n - L_n}} |\nabla u_n - \nabla \cR(g_n, u_n, 2^{-2N-2}\fB)|^2_{g_n} \vol_{g_n}\\
=\ & \sum_{i \in I} \int_{2^{-2N-2}B_i} |\nabla u_n - \nabla \cR(g_n, u_n, 2^{-2N-2}\fB)|^2_{g_n} \vol_{g_n}\\
=\ & \sum_{i \in I}\int_{2^{-2N-2}\psi_n(B_i)} |\nabla v_n - \nabla \cR(\sigma_n, v_n, 2^{-2N-2}\psi_n(\fB))|_{\sigma_n}^2 \vol_{\sigma_n} \leq \delta_n,
\end{split}
\]
where in passing to the last line we used~\eqref{eq:iso-concentric} again, and for the inequality at the end we used~\eqref{eq:compactness-property-smallness} and Proposition~\ref{prop:good-min-max-seq}(c). Combining this with~\eqref{eq:leftover-on-cylinder} and the definition of $\mu_n$, we get the asserted estimate in conclusion (b).
\end{proof}

For all large enough $n$ such that Claim \ref{claim:compactness-property} applies, we consider, for $m \in \{4, \cdots, k_n -4\}$, the sub-cylinders
\[
C_{n,m} = \bC_{T_n - (m + 3)L_n, T_n - (m-3)L_n}, \quad C_{n, m}' = \bC_{T_n - (m+1)L_n, T_n - (m-1)L_n},
\]
and define
\[
\beta_{n, m} = \left\{
\begin{array}{ll}
0, & \text{ if }E(g_n, u_n; C_{n, m}) = 0,\\ 
\frac{\sup_{\fB}\int_{C_{n, m}}|\nabla u_n - \nabla \cR(g_n, u_n, 2^{-2N-2}\fB)|^2_{g_n} \vol_{g_n}}{E(g_n, u_n; C_{n, m})}, & \text{ if }E(g_n, u_n; C_{n, m}) > 0,
\end{array}
\right.
\]
where in the second case, the supremum is taken over all finite, disjoint collection $\fB$ of geodesic disks in $(C_{n, m}, g_n)$ with radius at most $R_n$. 
\begin{claim}\label{claim:bad-cylinders}
Letting 
\[
I_n := \{4 \leq m \leq k_n - 4\ |\ \beta_{n, m} > \mu_n \},
\]
we have for all sufficiently large $n$ that 
\begin{equation}\label{eq:bad-annuli}
\sum_{m \in I_n} \int_{C_{n, m}}|\nabla u_n|^2\leq o(1) +  20\mu_n\int_{\bC_{-T_n + 2L_n, T_n - 2L_n}} |\nabla u_n|^2.
\end{equation}
\end{claim}
\begin{proof}[Proof of Claim]
For each $m \in I_n$, there exists by definition a finite, disjoint collection $\fB_{n, m}$ of geodesic disks in $(C_{n, m}, g_n)$, each with radius at most $R_n$, such that
\begin{equation}\label{eq:bad-cylinder-lower-bound}
\mu_n \cdot E(g_n, u_n; C_{n, m}) \leq \int_{C_{n, m}} |\nabla u_n - \nabla \cR(g_n, u_n, 2^{-2N-2}\fB_{n, m})|_{g_n}^2 \vol_{g_n}.
\end{equation}
Note that, given $l < m$, we have $C_{n, l} \cap C_{n, m} \neq \emptyset$ only if $l \geq m - 6$. Therefore, it is possible to partition $\{C_{n, m}\}_{m \in I_n}$ into at most $10$ subcollections consisting of mutually disjoint cylinders. (See the labeling scheme explained in~\cite[pages 34 to 35]{Evans-Gariepy}.) If $\{C_{n, m}\}_{m\in I_{n}'}$ denotes one such subcollection, then with the help of the estimate from Claim~\ref{claim:compactness-property}(b) applied with $\fB = \cup_{m \in I_{n}'}\fB_{n, m}$, we see that
\[
\sum_{m \in I_{n}'} \mu_n E(g_n, u_n; C_{n, m}) \leq \mu_n^2 E(g_n, u_n; \bC_{-T_n + L_n, T_n - L_n}).
\]
Summing over the at most $10$ subcollections then gives
\begin{equation}\label{eq:bad-annuli-small}
\sum_{m \in I_{n}} E(g_n, u_n; C_{n, m}) \leq 10 \mu_n E(g_n, u_n; \bC_{-T_n + L_n, T_n - L_n}).
\end{equation}
Using (p2) and~\eqref{eq:no-energy-on-finite-part}, we infer from~\eqref{eq:bad-annuli-small} that eventually
\[
\begin{split}
\sum_{m \in I_n} \int_{C_{n, m}}\frac{|\nabla u_n|^2}{2} \leq\ & (1 + \alpha_n)^2 \cdot 10\mu_n E(g_n, u_n; \bC_{-T_n + L_n, T_n - L_n})\\
\leq \ & o(1) + (1 + \alpha_n)^4 \cdot 10\mu_n E(g_{\pro}, u_n; \bC_{-T_n + 2L_n, T_n - 2L_n})\\
\leq\ & o(1) +  20\mu_n\int_{\bC_{-T_n + 2L_n, T_n - 2L_n}} \frac{|\nabla u_n|^2}{2}.
\end{split}
\]
This proves~\eqref{eq:bad-annuli}.
\end{proof}
\begin{claim}\label{claim:good-cylinders}
For all sufficiently large $n$, we have
\begin{equation}\label{eq:good-annuli}
\int_{C_{n, m}'} |(u_n)_{\theta}|^2 \leq \frac{1}{30}\int_{C_{n, m}} |\nabla u_n|^2, \quad \text{for all }m \not\in I_n.
\end{equation}
\end{claim}
\begin{proof}[Proof of Claim]
Suppose towards a contradiction that, up to taking a subsequence, there exists for each $n$ some $m_n \in \{4, \cdots, k_n - 4\}\setminus I_n$ so that
\begin{equation}\label{eq:good-annuli-2}
\int_{C_{n, m_n}'} |(u_n)_\theta|^2 > \frac{1}{30}\int_{C_{n, m_n}} |\nabla u_n|^2.
\end{equation}
To shorten notation we let $t_n = T_n - m_{n}L_n$, and define
\[
\widetilde{u}_n := u_n \circ \boldsymbol{\tau}_{t_n} = v_n \circ (\psi_n \circ \boldsymbol{\tau}_{t_n}),\quad\quad \widetilde{g}_n := (f_n(t_n))^2 \cdot \boldsymbol{\tau}_{t_n}^*g_n = \big(\frac{\lambda_n}{f_n(t_n)} \big)^{-2} (\psi_n \circ \boldsymbol{\tau}_{t_n})^*\sigma_n.
\]
Noting that $t_n \in [-T_n + 4L_n, T_n - 4L_n]$, we have by (p3) and~\eqref{eq:length-control} that
\begin{equation}\label{eq:smooth-convergence-to-conformal-flat}
\widetilde{g}_n \to g \quad \text{smoothly on }\bC_{-3L, 3L}, 
\end{equation}
where $g$ is a metric conformal to $g_{\pro}$. Now, given any geodesic disk $B = \overline{B_{g}(x, t)}$ in $(\bC_{-3L, 3L}, g)$ with radius $t \leq 1$, eventually $B_{\widetilde{g}_n}(x, \frac{2t}{3})$ is a geodesic disk with respect to $\widetilde{g}_n$, and hence 
\[
B_n: = B_{g_n}(\tau_{t_n}(x), \frac{2t}{3f_n(t_n)}) = \tau_{t_n}\big(B_{\widetilde{g}_n}(x, \frac{2t}{3})\big)
\]
is a geodesic disk in $(C_{n, m_n}, g_n)$, and has radius not exceeding $\frac{2R_n}{3}$. Claim~\ref{claim:compactness-property} then allows us to define
\[
h_n: = \cR(g_n, u_n, 2^{-2N-2}B_n) \circ \tau_{t_n},
\]
which yields an energy-minimizing map on $B_{\widetilde{g}_n}(x, \frac{t}{3 \cdot 2^{2N +1}})$ with respect to $\widetilde{g}_n$ that agrees with $\widetilde{u}_n$ on the boundary. Following the argument in the proof of Lemma~\ref{lemm:ep-reg-almost-harmonic} leading up to the estimates~\eqref{eq:L2-convergence},~\eqref{eq:hn-small-energy-g}, and~\eqref{eq:W12-close-g}, we obtain for all sufficiently large $n$ that
\begin{equation}\label{eq:hn-small-energy-cylinder}
\begin{split}
\int_{B_{g}(x, \frac{t}{2^{2N + 3}})}|\nabla h_n|_{g}^2 \vol_{g} \leq 4\int_{B_{\widetilde{g}_n}(x, \frac{t}{3 \cdot 2^{2 N + 1}})}|\nabla \widetilde{u}_n|_{\widetilde{g}_n}^2 \vol_{\widetilde{g}_n} \leq\ & 8E(\widetilde{g}_n, \widetilde{u}_n; \bC_{-3L_n, 3L_n})\\
\leq\ & 8\ep_2 < \ep_{\reg},
\end{split}
\end{equation}
and that
\begin{equation}\label{eq:W12-close-cylinder}
\begin{split}
\int_{B_{g}(x, \frac{t}{2^{2N + 3}})} |h_n - \widetilde{u}_n|^2 +  |\nabla h_n - \nabla \widetilde{u}_n|_{g}^2 \vol_{g} \leq\ & C\int_{B_{\widetilde{g}_n}(x, \frac{t}{3 \cdot 2^{2N + 1}})} |\nabla h_n - \nabla\widetilde{u}_n|^2_{\widetilde{g}_n} \vol_{\widetilde{g}_n}\\
=\ & C\int_{2^{-2N-2}B_n} |\nabla (h_n \circ\tau_{t_n}^{-1}) - \nabla u_n|^2_{g_n} \vol_{g_n}\\
\leq\ &  C\mu_n \cdot E(\widetilde{g}_n, \widetilde{u}_n; \bC_{-3L_n, 3L_n}).
\end{split}
\end{equation}
In~\eqref{eq:hn-small-energy-cylinder}, the last two inequalities  follows from~\eqref{eq:small-globally} and our choice of $\ep_2$ in~\eqref{eq:ep2-choice}. In~\eqref{eq:W12-close-cylinder}, the constants $C$ are independent of $n$, and the last inequality follows from the fact that $m_n \not\in I_n$.

To continue, in the case where 
\begin{equation}\label{eq:cylinder-case-1}
\limsup_{n \to \infty} E(\widetilde{g}_n, \widetilde{u}_n; \bC_{-3L_n, 3L_n}) > 0, 
\end{equation}
we use the above estimates together with Theorem~\ref{thm:ep-regularity} to extract a subsequence of $(\widetilde{u}_n)$, not relabeled, that converges strongly in $W^{1,2}_{\loc}$ on $(-3L, 3L) \times S^1$ to a harmonic map $u$ with respect to $g$. By property (p2) and~\eqref{eq:small-globally} we have
\begin{equation}\label{eq:small-energy-case-1}
\begin{split}
\int_{\bC_{-3L, 3L}} \frac{|\nabla u|^2}{2} \leq \liminf_{n \to \infty}\int_{\bC_{-3L_n, 3L_n}} \frac{|\nabla \widetilde{u}_n|^2}{2} \leq\ & \liminf_{n \to \infty} (1 + \alpha_n)^{2}E(g_n, u_n; C_{n, m_n})\\
\leq\ & \ep_2 < \ep_{\text{cyl}}.
\end{split}
\end{equation}
On the other hand, from~\eqref{eq:good-annuli-2} and~\eqref{eq:length-control}, we get
\begin{equation}\label{eq:uneven-case-1}
\begin{split}
\int_{\bC_{-L, L}} |u_{\theta}|^2=\ & \lim_{n \to \infty}\int_{\bC_{-L_n, L_n}} |(\widetilde{u}_n)_{\theta}|^2 \\
\geq\ & \limsup_{n \to \infty}\frac{1}{30}\int_{\bC_{-3L_n, 3L_n}} |\nabla\widetilde{u}_n|^2\geq\frac{1}{30}\int_{\bC_{-2L, 2L}} |\nabla u|^2.
\end{split}
\end{equation}
Estimating the penultimate term above using property (p2) and~\eqref{eq:cylinder-case-1} instead, we have
\begin{equation}\label{eq:lower-bound-case-1}
\begin{split}
\int_{\bC_{-L, L}} |u_{\theta}|^2 \geq \ & \limsup_{n \to \infty}\frac{(1 + \alpha_n)^{-2}}{15}E(\widetilde{g}_n, \widetilde{u}_n; \bC_{-3L_n, 3L_n}) > 0.
\end{split}
\end{equation}
Now observe $u$ is also a harmonic map with respect to $g_{\pro}$, since the latter is conformal to $g$. Thus, from~\eqref{eq:small-energy-case-1} along with~\cite[Proposition B.1]{Colding-Minicozzi08b}, we have~\eqref{eq:CM-estimate} with $u$ in place of $v$, which together with~\eqref{eq:uneven-case-1} contradicts~\eqref{eq:lower-bound-case-1}.

It remains to treat the case
\begin{equation}\label{eq:cylinder-case-2}
\limsup_{n \to \infty} E(\widetilde{g}_n, \widetilde{u}_n; \bC_{-3L_n, 3L_n})  = 0.
\end{equation}
For brevity, we let $E_n: = E(\widetilde{g}_n, \widetilde{u}_n; \bC_{-3L_n, 3L_n})$, which is positive by the strict inequality~\eqref{eq:good-annuli-2}. Following the proof of~\cite[Lemma B.20]{Colding-Minicozzi08b}, we define
\[
v_n = \frac{\widetilde{u}_n - \fint_{\bC_{-3L_n, 3L_n}}\widetilde{u}_n}{\sqrt{E_n}},
\]
where the average is taken with respect to the product metric $g_{\pro}$. Since $L_n \to L$ by~\eqref{eq:length-control}, the Poincar\'e inequality gives a constant $A$ independent of $n$ such that
\begin{equation}\label{eq:vn-W12}
\int_{\bC_{-3L_n, 3L_n}} |v_n|^2 + |\nabla v_n|^2 \leq \frac{A \int_{\bC_{-3L_n, 3L_n}}|\nabla \widetilde{u}_n|^2}{E_n} \leq 2A\cdot (1 + \alpha_n)^2.
\end{equation}
Moreover, for each geodesic disk $\overline{B_g(x, r)}$ in $(\bC_{-3L, 3L}, g)$ with $r \leq 1$, letting $h_n$ be the harmonic map appearing in~\eqref{eq:hn-small-energy-cylinder} and~\eqref{eq:W12-close-cylinder}, we define
\[
H_n = \frac{h_n - \fint_{\bC_{-3L_n, 3L_n}}\widetilde{u}_n}{\sqrt{E_n}},
\]
and observe by~\eqref{eq:W12-close-cylinder} that
\begin{equation}\label{eq:kn-vn-W12}
\lim_{n \to \infty}\int_{B_g(x, \frac{r}{2^{2N + 3}})} |H_n - v_n|^2 + |\nabla H_n - \nabla v_n|^2 = 0,
\end{equation}
which together with~\eqref{eq:vn-W12} shows that $(H_n)$ is bounded in $W^{1,2}$ on $B_g(x, \frac{r}{2^{2N + 3}})$. In addition, recall from~\eqref{eq:hn-small-energy-cylinder} that
\[
\int_{B_{g}(x, \frac{r}{2^{2N + 3}})}|\nabla h_n|_{g}^2 \vol_{g} \leq 8E_n < \ep_{\reg}.
\]
Applying Theorem~\ref{thm:ep-regularity} to $h_n$ in isothermal charts centered around points in $B_{g}(x, \frac{r}{2^{2N + 4}})$ as we did in the proof of Lemma~\ref{lemm:ep-reg-almost-harmonic}, but this time dividing through by $E_n$ in the gradient estimate~\eqref{eq:ep-reg-estimate}, and using in addition the $W^{1, 2}$-bound noted after~\eqref{eq:kn-vn-W12}, we obtain after some routine steps a subsequence of $(H_n)$ that converges in $C^{1}$ on $B_g(x, \frac{r}{2^{2N+4}})$ to a limit which, because $E_n \to 0$, is a harmonic function (with respect to $g$ on the domain) into the flat $\RR^{q}$. Returning to~\eqref{eq:kn-vn-W12}, and performing another covering argument, we get a subsequence of $(v_n)$ that converges strongly in $W^{1, 2}_{\loc}$ on $(-3L, 3L) \times S^1$ to a harmonic function
\[
v: ((-3L, 3L) \times S^1, g ) \to (\RR^{q}, g_{\euc}).
\]
With this convergence of $(v_n)$ at hand, by the arguments leading to the estimates~\eqref{eq:uneven-case-1} and~\eqref{eq:lower-bound-case-1}, we get
\[
\begin{split}
\int_{\bC_{-L, L}} |v_{\theta}|^2 =\lim_{n \to \infty}\int_{\bC_{-L_n, L_n}} |(v_n)_{\theta}|^2\geq\ & \limsup_{n \to \infty} \frac{\frac{1}{30}\int_{\bC_{-3L_n, 3L_n}}|\nabla \widetilde{u}_n|^2}{E_n}\\
\geq\ & \limsup_{n \to \infty}\frac{(1 + \alpha_n)^{-2}}{15} = \frac{1}{15} > 0,
\end{split}
\]
and that
\[
\int_{\bC_{-L, L}} |v_{\theta}|^2 \geq \frac{1}{30}\int_{\bC_{-2L, 2L}} |\nabla v|^2.
\]
Since $g$ is conformal to $g_{\pro}$, the inequality~\eqref{eq:CM-estimate} holds by our choice of $L$, and a contradiction results as in the previous case. Thus eventually we must have~\eqref{eq:good-annuli}. This ends the proof of Claim~\ref{claim:good-cylinders}.
\end{proof}

\begin{claim}\label{claim:finished-under-extra-hypo}
Proposition~\ref{prop:neck} holds assuming~\eqref{eq:small-globally}.
\end{claim}
\begin{proof}[Proof of Claim]
Given what we have proved thus far, we will derive a contradiction with~\eqref{eq:leftover-on-cylinder}. To begin, for all sufficiently large $n$, by Claims~\ref{claim:bad-cylinders} and~\ref{claim:good-cylinders}, we have available both~\eqref{eq:bad-annuli} and~\eqref{eq:good-annuli}. Summing the latter over $m \not\in I_n$, bounding the overlap as suggested in~\cite[page 2574]{Colding-Minicozzi08b}, and using~\eqref{eq:no-energy-on-finite-part}, we get
\[
\sum_{m \not\in I_n} \int_{C_{n, m}'} |(u_n)_{\theta}|^2 \leq \frac{7}{30}\int_{\bC_{-T_n+ L_n, T_n - L_n}} |\nabla u_n|^2\leq o(1) + \frac{7}{30}\int_{\bC_{-T_n+ 2L_n, T_n - 2L_n}} |\nabla u_n|^2,
\]
which when added to~\eqref{eq:bad-annuli} gives
\[
\int_{\bC_{-T_n + 3L_n, T_n - 3L_n}} |(u_n)_{\theta}|^2 \leq o(1) +  (\frac{7}{30}+ 20\mu_n)\int_{\bC_{-T_n + 2L_n, T_n - 2L_n}} |\nabla u_n|^2.
\]
Using again~\eqref{eq:no-energy-on-finite-part}, and noting that $\mu_n \to 0$ as $n \to \infty$, we arrive at
\begin{equation}\label{eq:uneven-distribution}
\int_{\bC_{-T_n + 2L_n, T_n -2L_n}} |(u_n)_{\theta}|^2 \leq o(1) + \frac{1}{3}\int_{\bC_{-T_n + 2L_n, T_n-2L_n}} |(u_n)_{t}|^2 + |(u_n)_{\theta}|^2,
\end{equation}
for all sufficiently large $n$. Rearranging this and combining it with the elementary inequality
\[
\frac{|u_{t}|^2 + |u_{\theta}|^2}{2} - \sqrt{|u_{t}|^2 |u_{\theta}|^2 - \bangle{u_{t}, u_{\theta}}^2} \geq \frac{(|u_{t}| - |u_{\theta}|)^2}{2},
\]
we have
\begin{equation}\label{eq:uneven-consequence}
\begin{split}
&\frac{1}{3}\int_{\bC_{-T_n + 2L_n, T_n -2L_n}} |\nabla u_n|^2 - o(1)\leq \int_{\bC_{-T + 2L_n, T_n -2L_n}}|(u_n)_t|^2 - |(u_n)_{\theta}|^2 \\
\leq\ & \Big( \int_{\bC_{-T_{n} + 2L_{n}, T_{n} -2L_{n}}} (|(u_{n})_t| - |(u_{n})_\theta|)^2  \Big)^{\frac{1}{2}}\Big(2\int_{\bC_{-T_{n} + 2L_{n}, T_{n} -2L_{n}}} |(u_{n})_t|^2 + |(u_{n})_\theta|^2 \Big)^{\frac{1}{2}}\\
\leq\ & 2^{\frac{3}{2}}\Big( E(g_{\pro}, u_{n}; \bC_{-T_{n}, T_{n}}) - A(u_{n}; \bC_{-T_{n}, T_{n}}) \Big)^{\frac{1}{2}}\cdot \big( E(g_{\pro}, u_{n}; \bC_{-T_{n}, T_{n}}) \big)^{\frac{1}{2}}.
\end{split}
\end{equation}
To estimate the last line, note that by property (p2) and~\eqref{eq:small-globally}, we have
\begin{equation}\label{eq:small-globally-euc}
E(g_{\pro}, u_n; \bC_{-T_n, T_n}) \leq  (1 + \alpha_n)^2\cdot E(g_n, u_n; \bC_{-T_n, T_n}) \leq (1 + \alpha_n)^2\cdot \ep_2.
\end{equation}
By the first inequality in~\eqref{eq:small-globally-euc}, and using~\eqref{eq:small-globally} again, we have
\begin{equation}\label{eq:energy-comparable}
\begin{split}
&E(g_{\pro}, u_n; \bC_{-T_n, T_n}) - A(u_n; \bC_{-T_n, T_n})\\
\leq\ & E(g_n, u_n; \bC_{-T_n, T_n}) - A(u_n; \bC_{-T_n, T_n}) + [(1 + \alpha_n)^2 - 1]\cdot \ep_2\\
=\ & E(\sigma_{n}, v_n; \psi_{n}(\bC_{-T_n, T_n})) - A(v_n; \psi_n(\bC_{-T_n, T_n})) + [(1 + \alpha_n)^2 - 1] \cdot \ep_2\\
\leq\ & E(\sigma_n, v_n) - A(v_n) + [(1 + \alpha_n)^2 - 1]\cdot \ep_2 \longrightarrow 0, \quad \text{as }n \to \infty, 
\end{split}
\end{equation}
the convergence at the end being a consequence of Proposition~\ref{prop:good-min-max-seq}(b). Substituting~\eqref{eq:energy-comparable} and~\eqref{eq:small-globally-euc} back into~\eqref{eq:uneven-consequence} yields
\[
\lim_{n \to \infty}\int_{\bC_{-T_n + 2L_n, T_n -2L_n}} |\nabla u_n|^2 =0,
\]
which together with property (p2) and the assumption~\eqref{eq:no-energy-on-finite-part} gives 
\[
E(g_n, u_n; \bC_{-T_n + L_n, T_n-L_n}) \to 0,
\]
contradicting~\eqref{eq:leftover-on-cylinder}, and we are done.
\end{proof}

\begin{claim}\label{claim:small-globally}
Under the original assumptions of Proposition~\ref{prop:neck}, the condition~\eqref{eq:small-globally} must hold.
\end{claim}
\begin{proof}[Proof of Claim]
Assume by contradiction that, along a subsequence, there holds 
\[
E(g_{n}, u_n; \bC_{-T_n, T_n}) > \ep_2\quad \text{for all }n.
\]
For each $n$, consider the function $\eta_n:[0, T_n] \to \RR$ defined by
\[
\begin{split}
\eta_n(T) :=\ & \sup_{-T_n + T \leq t \leq T_n - T} E(g_{n}, u_n; \bC_{t-T, t+T})\\
=\ & \sup\big\{ E(g_n, u_n; \bC_{a, b})\ \big|\ [a, b] \subset [-T_n, T_n], \ b - a = 2T \big\}.
\end{split}
\]
Observe that $\lim_{T \to 0}\eta_n(T) = 0$, and that
\begin{equation}\label{eq:eta-n-subadditive}
\eta_n(T) \leq \eta_n(T') \leq \eta_n(T) + \eta_n(T' - T),\quad\text{whenever }T \leq T'
\end{equation}
which together shows that $\eta_n$ is continuous. Since $\eta_n(T_n) > \ep_2$ by assumption, while $\lim_{n \to \infty}\eta_n(T) = 0$ for each fixed $T > 0$ by~\eqref{eq:no-energy-on-finite-part}, a standard argument yields a further subsequence $(n_k)$, along with sequences $T_k' \to \infty$ and $t_k \in [-T_{n_k} + T_k', T_{n_k} - T_k']$, such that 
\begin{equation}\label{eq:energy-persist}
\ep_2 = \eta_{n_k}(T_k') =  E(g_{n_k}, u_{n_k}; \bC_{t_{k} -T_k', t_k + T_k'}).
\end{equation}
It is straightforward to verify that the sequence 
\[
\{\bC_{-T_k',T_k'},\ \psi_{n_k} \circ \boldsymbol{\tau}_{t_k},\ \lambda_{n_k},\ f_{n_k}(t_k + \cdot)\}
\]
is well-prepared in the sense of Definition~\ref{defi:well-prepared}, where in property (p2) we use $\sigma_{n_k}$ in place of $\sigma_n$. Moreover, letting
\[
\widehat{g}_k := \lambda_{n_k}^{-2} (\psi_{n_k} \circ \boldsymbol{\tau}_{t_k})^*\sigma_{n_k} = \boldsymbol{\tau}_{t_k}^*g_{n_k},\quad \widehat{u}_{k} := v_{n_k} \circ (\psi_{n_k} \circ \boldsymbol{\tau}_{t_k}) = u_{n_k} \circ \boldsymbol{\tau}_{t_k},
\]
and noting that~\eqref{eq:no-energy-on-finite-part} as well as~\eqref{eq:small-globally} both hold with $n$ replaced by $k$, $T_n$ replaced by $T_k'$, and with $g_n$ and $u_n$ replaced by $\widehat{g}_k$ and $\widehat{u}_k$, it follows from Claim~\ref{claim:finished-under-extra-hypo} that
\[
0 =  \lim_{k \to \infty}E(\widehat{g}_k, \widehat{u}_k ; \bC_{-T_k', T_k'}) = \lim_{k \to \infty}E(g_{n_k}, u_{n_k}; \bC_{t_k-T_k', t_k +T_k'}),
\]
which contradicts~\eqref{eq:energy-persist}.
\end{proof}
In view of Claim~\ref{claim:finished-under-extra-hypo} and Claim~\ref{claim:small-globally}, the proof of Proposition~\ref{prop:neck} is complete.\\
\end{proof}

\begin{proof}[Proof of Proposition~\ref{prop:collar-bubble}]
As in the previous proof, by properties (p1), (p2) and (p3) we mean the ones from Definition~\ref{defi:well-prepared}. We first notice that by Proposition~\ref{prop:good-min-max-seq}(a)(b) there holds
\begin{equation}\label{eq:cylinder-uniform-energy-bound}
\limsup_{n \to \infty} E(g_{n}, u_n; \bC_{-T_n, T_n}) \leq W.
\end{equation}
Next, towards a contradiction, suppose that no matter what $Q$ is, along any subsequence of $(n)$ and for any sequences of intervals $I_{n, i}' \subset I_{n, i} \subset [-T_n, T_n]$ ($i = 1, \cdots, Q$) satisfying conclusions (a), (b), (c) and (d), we have that (e) fails.

As a first step, note that~\eqref{eq:there-is-energy} gives some $\eta > 0$, a subsequence of $(g_n, u_n)$, which we do not relabel, and, for each $n$, some $t_{n, 1} \in [-T_n + 1, T_n - 1]$, such that
\begin{equation}\label{eq:energy-somewhere}
E(g_{n}, u_n; \bC_{t_{n, 1} - 1, t_{n, 1} + 1}) \geq \eta.
\end{equation}
By~\eqref{eq:iv-weakened} we must have
\begin{equation}\label{eq:concentration-far-from-endpoint}
S_{n, 1} : = \min\{T_n - t_{n, 1}, t_{n, 1} + T_n\} \to \infty \text{ as }n \to \infty.
\end{equation}
Now define
\[
w_{n, 1} = u_{n} \circ \boldsymbol{\tau}_{t_{n, 1}}, \quad h_{n, 1} = (f_{n}(t_{n, 1}))^2\cdot \boldsymbol{\tau}_{t_{n, 1}}^* g_{n}, \quad \text{on }\bC_{-S_{n, 1}, S_{n, 1}}.
\]
By~\eqref{eq:concentration-far-from-endpoint} and property (p3), we obtain a metric $g$ on $\RR \times S^1$, conformal to $g_{\pro}$, such that $h_{n, 1}$ converge smoothly to $g$ on $\bC_{-T, T}$, for all $T > 0$. Note also that $w_{n, 1}$ and $h_{n, 1}$ can alternatively be expressed as 
\[
w_{n, 1} = v_{n} \circ (\psi_{n} \circ \boldsymbol{\tau}_{t_{n, 1}}), \quad h_{n, 1} = \Big(\frac{\lambda_n}{f_n(t_{n, 1})}\Big)^{-2} (\psi_n \circ \boldsymbol{\tau}_{t_{n, 1}})^*\sigma_n,
\]
and that $\lim_{n \to \infty}\frac{\lambda_n}{f_n(t_{n, 1})} = 0$ by property (p1). Thus we find ourselves in a position to invoke Proposition~\ref{prop:domain-convergence}, with 
\[
\Omega = \RR \times S^1, \quad \Omega_n = \Int(\bC_{-S_{n, 1}, S_{n, 1}}),\quad \varphi_n = \psi_n \circ \boldsymbol{\tau}_{t_{n, 1}},
\]
and with ``$\lambda_n$'' chosen to be $\frac{\lambda_n}{f_n(t_{n, 1})}$. As a result, and recalling that $g$ is conformal to $g_{\pro}$, we obtain a finite-energy, weakly conformal harmonic map $w_0: (\RR \times S^1, g_{\pro}) \to M$, along with a finite set $\cA \subset \RR \times S^1$, such that 
\begin{equation}\label{eq:convergence-to-mu1}
\frac{1}{2}|\nabla w_{n, 1}|_{h_{n, 1}}^2 \vol_{h_{n, 1}} \to \frac{1}{2}|\nabla w_0|^2_{g_{\pro}} \vol_{g_{\pro}} + \sum_{x \in \cA}m_x'\delta_x =: \mu_1,
\end{equation}
where each $m_x' > \ep_2$ by~\eqref{eq:ep2-choice}, and that
\begin{equation}\label{eq:convergence-to-w0}
w_{n, 1} \to w_0\quad \text{in }W^{1, 2}_{\loc}((\RR \times S^1) \setminus \cA).
\end{equation}
Since $(\RR \times S^1, g_{\pro})$ is conformally equivalent to $(S^2 \setminus \{\text{poles}\}, g_{S^2})$ and since $w_0$ has finite energy, the removable singularity theorem (\cite{Sacks-Uhlenbeck81}) and the fact that $\ep_2 < \ep_{\text{gap}}$ implies that 
\[
E(g_{\pro}, w_0) > \ep_2,\quad \text{provided }w_0 \text{ is not a constant.}
\]
Now, by~\eqref{eq:cylinder-uniform-energy-bound} and the convergence~\eqref{eq:convergence-to-mu1} we infer that
\begin{equation}\label{eq:finite-mass}
\|\mu_1\|: = \mu_1(\RR \times S^1) \leq W < \infty.
\end{equation}
On the other hand, by~\eqref{eq:energy-somewhere} we have
\[
E(h_{n, 1}, w_{n, 1}; \bC_{-1, 1}) = E(g_n, u_n; \bC_{t_{n, 1} - 1, t_{n, 1} + 1}) \geq \eta, \quad \text{for all }n.
\]
From this lower bound and the definition of $\mu_1$ in~\eqref{eq:convergence-to-mu1}, we conclude that either $w_0$ is non-constant or $\cA \neq \emptyset$, and consequently
\begin{equation}\label{eq:energy-gap}
\|\mu_1\| = E(g_{\pro}, w_0) + \sum_{x \in \cA} m'_x > \ep_2.
\end{equation}
We then choose a subsequence $(h_{n_k, 1}, w_{n_k, 1})$ such that $S_{n_k, 1} > k^3$ and that
\[
\big|E(h_{n_k, 1}, w_{n_k, 1}; \bC_{-k, k}) - \mu_1(\bC_{-k, k})\big| < \frac{1}{k},
\]
\[
\big|E(h_{n_k, 1}, w_{n_k, 1}; \bC_{-k^3, k^3} \setminus \bC_{-k, k}) - \mu_1(\bC_{-k^3,k^3} \setminus \bC_{-k, k})\big| < \frac{1}{k}.
\]
From these we deduce:
\begin{align}
\big| E(h_{n_k, 1}, w_{n_k, 1}; \bC_{-k, k}) - \|\mu_1\| \big| = o(1), \label{eq:close-to-mu1}\\
E(h_{n_k, 1}, w_{n_k, 1}, \bC_{-k^3, k^3} \setminus \bC_{-k, k}) = o(1).\label{eq:transition-region}
\end{align}
To continue, we write $(g_{n}, u_n)$ for $(g_{n_k}, u_{n_k})$ and denote $t_{n_k, 1}, S_{n_k, 1}$ simply by $t_{n, 1}, S_{n, 1}$, respectively. Also, we define
\begin{equation}\label{eq:first-intervals}
J_{n, 1} = [t_{n, 1} - n^3, t_{n, 1} + n^3],\ \ I_{n, 1} = [t_{n, 1} - n^2, t_{n, 1} + n^2],\ \ I_{n, 1}' = [t_{n, 1} - n, t_{n, 1} + n].
\end{equation}
Then by~\eqref{eq:transition-region} we have
\begin{equation}\label{eq:small-transition-middle}
\limsup_{n \to \infty} E(g_{n}, u_n; (I_{n, 1} \setminus I_{n, 1}') \times S^1) \leq \limsup_{n \to \infty} E(g_{n}, u_n; (J_{n, 1} \setminus I_{n, 1}') \times S^1) = 0,
\end{equation}
so that (c) holds, while~\eqref{eq:close-to-mu1} implies
\begin{equation}\label{eq:branch-captured}
\begin{split}
\big|E(g_{n}, u_{n}; I_{n, 1}' \times S^1) - \|\mu_1\|\big| = o(1),
\end{split}
\end{equation}
which together with~\eqref{eq:convergence-to-mu1},~\eqref{eq:finite-mass} and~\eqref{eq:energy-gap} yields (d). Also, $\diam I_{n, 1}'$ and $\dist(I_{n, 1}', \partial I_{n, 1})$ both tend to infinity as $n \to \infty$, and, since $S_{n, 1} > n^3$, we have
\[
\dist(I_{n, 1}, \{-T_n, T_n\}) \geq n^3 - n^2 \to \infty.
\]
This proves (a). There being only one sequence of intervals $(I_{n, 1})$ at this stage, conclusion (b) holds vacuously. Finally, thanks to~\eqref{eq:cylinder-uniform-energy-bound},~\eqref{eq:energy-gap} and~\eqref{eq:branch-captured}, we have 
\begin{equation}\label{eq:definite-drop}
\begin{split}
\limsup_{n \to \infty}E(g_{n}, u_n; ([-T_n, T_n] \setminus I_{n, 1}')\times S^1) \leq W - \ep_2.
\end{split}
\end{equation}
Having verified conclusions (a), (b), (c) and (d) for the intervals $I_{n, 1}'$ and $I_{n, 1}$, we infer from the assumption made at the beginning of the proof (with $Q = 1$) that, up to taking a subsequence, there exist $\eta > 0$ and, for each $n$, some $t_{n, 2} \in [-T_n+1, T_n-1]$, such that $[t_{n, 2} - 1, t_{n, 2} + 1] \cap I_{n, 1} = \emptyset$ and that 
\[
E(g_{n}, u_n; \bC_{t_{n, 2} - 1, t_{n, 2} + 1}) \geq \eta.
\]
By~\eqref{eq:iv-weakened} and~\eqref{eq:small-transition-middle}, we must have
\begin{equation}\label{eq:concentration-far-step-2}
S_{n, 2}: =\dist(t_{n, 2}, I_{n, 1} \cup \{-T_n, T_n\}) \to \infty \text{ as }n \to \infty.
\end{equation}
Defining, as before, 
\[
w_{n, 2} = u_n \circ \boldsymbol{\tau}_{t_{n, 2}}, \quad h_{n, 2} = f_n(t_{n, 2})^2 \cdot \boldsymbol{\tau}_{t_{n, 2}}^*g_n,\quad \text{on }\bC_{-S_{n,2}, S_{n, 2}},
\]
we see from~\eqref{eq:concentration-far-step-2} and property (p3) that $h_{n,2}$ converges smoothly locally on $\RR \times S^1$ to a limiting metric which is conformal to $g_{\pro}$. Again noting the expressions $w_{n, 2} = v_n \circ (\psi_n \circ \boldsymbol{\tau}_{t_{n, 2}})$ and
\[
h_{n, 2} = \big( \frac{\lambda_n}{f_n(t_{n, 2})} \big)^{-2}\cdot(\psi_n \circ \boldsymbol{\tau}_{t_{n, 2}})^*\sigma_n,
\]
with $\lim_{n \to \infty}\frac{\lambda_n}{f_n(t_{n, 2})} = 0$, the latter being a consequence of (p1), we may repeat the previous argument to extract a further subsequence of $(g_{n}, u_{n}), t_{n, 2}$ and $S_{n, 2}$, which we do not relabel, such that $S_{n, 2} > n^3$, and that, defining $J_{n, 2}$, $I_{n, 2}$ and $I_{n, 2}'$ by~\eqref{eq:first-intervals} with $t_{n, 1}$ replaced by $t_{n, 2}$, we have 
\vskip 1mm
\begin{enumerate}
\item[(1)] $\lim_{n \to \infty}E(g_{n}, u_n; (J_{n, 2}\setminus I_{n, 2}') \times S^1) = 0$,
\vskip 1mm
\item[(2)] $\boldsymbol{\tau}_{t_{n, 2}}^*\big(\frac{1}{2}|\nabla u_n|_{g_n}^2 \vol_{g_{n}}\big)$  weak$^*$-converge on compact subsets on $\RR \times S^1$ to a limiting measure $\mu_2$ with $\|\mu_2\| \in [\ep_2, W - \ep_2]$. Moreover, 
\[
\big| E(g_{n}, u_n; I_{n, 2}' \times S^1) - \|\mu_2\|  \big| = o(1),
\]
\vskip 1mm
\item[(3)] $\min\{\diam I_{n, 2}', \dist(I_{n, 2}', \partial I_{n, 2})\} \geq \min\{2n, n^2 - n\} \to \infty$,
\vskip 1mm
\item[(4)] $\min\{\dist(I_{n, 2}, \{-T_n, T_n\}), \dist(I_{n, 2}, I_{n, 1})\} \geq S_{n, 2} - n^2 \to \infty$ as $n \to \infty$,
\vskip 1mm
\item[(5)] $\limsup_{n \to \infty} E(g_{n}, u_n; ([-T_n, T_n] \setminus (I_{n, 1}' \cup I_{n, 2}'
))\times S^1) \leq W - 2\ep_2$.
\end{enumerate}
In particular, the intervals $I_{n, 1}', I_{n, 1}, I_{n, 2}', I_{n, 2}$ fulfill conclusions (a), (b), (c) and (d). Our initial assumption, this time with $Q = 2$, allows us to iterate the construction further. Continuing in this fashion, we inductively obtain, for each $Q \in \NN$, a subsequence of $(g_n, u_n)$, along with intervals $I_{n, i}', I_{n, i}$ for $i = 1, \cdots, Q$, satisfying among other things that
\[
\limsup_{n \to \infty}E(g_{n}, u_n; ([-T_n, T_n] \setminus \cup_{i = 1}^Q I_{n, i}') \times S^1) \leq W - Q \cdot \ep_2.
\]
This however implies that $Q\cdot \ep_2 \leq W$ for all $Q \in \NN$, which is a contradiction.
\end{proof}

\bibliographystyle{amsplain}
\bibliography{main}
\end{document}